\documentstyle[makeidx,12pt, amscd]{amsbook}
\input amssymb.sty
\font\tenfrakt=eufm10
\font\sevfrakt=eufm7
\font\fivfrakt=eufm5
\newfam\fraktfam
\textfont\fraktfam=\tenfrakt
\scriptfont\fraktfam=\sevfrakt
\scriptscriptfont\fraktfam=\fivfrakt
\def\frakt{\fam\fraktfam}

\newtheorem{thm}{Theorem}[section]
\newtheorem{lem}{Lemma}[section]
\newtheorem{prop}{Proposition}[section]
\newtheorem{cor}{Corollary}[section]
\theoremstyle{remark}
\newtheorem{rem}{Remark}[section]
\newtheorem{clm}{Claim}[section]
\theoremstyle{definition}
\newtheorem{defn}{Definition}[section]
\newtheorem{exam}{Example}[section]

\def\Ad{\mathop{\operatorname {Ad}}\nolimits}
\def\ad{\mathop{\operatorname {ad}}\nolimits}
\def\Aff{\mathop{\operatorname{Aff}}\nolimits}
\def\aff{\mathop{\operatorname{aff}}\nolimits}
\def\Aut{\mathop{\operatorname {Aut}}\nolimits}
\def\arcctg{\mathop{\operatorname{arcctg}}\nolimits}
\def\Card{\mathop{\operatorname {Card}}\nolimits}

\def\cent{\mathop{\operatorname {cent}}\nolimits}
\def\Char{\mathop{\operatorname {Char}}\nolimits}
\def\codim{\mathop{\operatorname {codim}}\nolimits}
\def\ctg{\mathop{\operatorname{ctg}}\nolimits}
\def\Coh{\mathop{\operatorname {Coh}}\nolimits}
\def\Coker{\mathop{\operatorname{Coker}}\nolimits}

\def\curv{\mathop{\operatorname{Curv}}\nolimits}
\def\Curv{\mathop{\operatorname{Curv}}\nolimits}
\def\Der{\mathop{\operatorname {Der}}\nolimits}

\def\End{\mathop{\operatorname {End}}\nolimits}
\def\Ext{\mathop{\operatorname{{\cal E}xt}}\nolimits}
\def\Fred{\mathop{\operatorname {Fred}}\nolimits}
\def\GL{\mathop{\operatorname {GL}}\nolimits}
\def\Ham{\mathop{\operatorname {Ham}}\nolimits}
\def\Hom{\mathop{\operatorname {Hom}}\nolimits}
\def\SL{\mathop{\operatorname{SL}}\nolimits}
\def\Id{\mathop{\operatorname {Id}}\nolimits}
\def\Ind{\mathop{\operatorname{Ind}}\nolimits}
\def\ind{\mathop{\operatorname{Ind}}\nolimits}
\def\inv{\mathop{\operatorname{inv}}\nolimits}
\def\Ker{\mathop{\operatorname{Ker}}\nolimits}
\def\Lie{\mathop{\operatorname {Lie}}\nolimits}
\def\ln{\mathop{\operatorname{ln}}\nolimits}
\def\MD{\mathop{\operatorname {MD}}\nolimits}
\def\Mat{\mathop{\operatorname {Mat}}\nolimits}

\def\Mp{\mathop{\operatorname {Mp}}\nolimits}
\def\mult{\mathop{\operatorname {mult}}\nolimits}
\def\Nil-Rad{\mathop{\operatorname {Nil-Rad}}\nolimits}

\def\pr{\mathop{\operatorname {pr}}\nolimits}
\def\Pfaffian{\mathop{\operatorname {Pfaff}}\nolimits}
\def\Rad{\mathop{\operatorname {Rad}}\nolimits}
\def\rank{\mathop{\operatorname {rank}}\nolimits}

\def\semi-linear{\mathop{\operatorname {semi-linear}}\nolimits}
\def\sgn{\mathop{\operatorname {sgn}}\nolimits}
\def\supp{\mathop{\operatorname {supp}}\nolimits}
\def\Sh{\mathop{\operatorname {Sh}}\nolimits}
\def\SL{\mathop{\operatorname {SL}}\nolimits}
\def\Spin{\mathop{\operatorname {Spin}}\nolimits}
\def\Sp{\mathop{\operatorname {Sp}}\nolimits}
\def\sgn{\mathop{\operatorname{sgn}}\nolimits}
\def\supp{\mathop{\operatorname{Supp}}\nolimits}

\def\Tr{\mathop{\operatorname{Tr}}\nolimits}
\def\tr{\mathop{\operatorname{Tr}}\nolimits}
\def\USp{\mathop{\operatorname{USp}}\nolimits}
\def\Vect{\mathop{\operatorname {Vect}}\nolimits}

\def\ker{\mathop{\operatorname{Ker}}\nolimits}

\def\CC{\Bbb C}
\def\HH{\Bbb H}
\def\NN{\Bbb N}
\def\RR{\Bbb R}

\def\ZZ{\Bbb Z}
\def\SS{\Bbb S}

\makeindex
\begin{document}
\title {Non Commutative Geometry Methods \\ for Group C*-algebras} 
\author {Do Ngoc Diep}
\address{Institute of Mathematics, National Center for Natural 
Science and Technology, P. O. Box 631, Bo Ho 10.000, Hanoi, Vietnam}
\email{dndiep@@thevinh.ncst.ac.vn }
\thanks{{\bf Acknowledgment.} This work was partially supported by
the International Centre for Theoretical Physics, Trieste, Italy, 
partially by the Alexander von Humboldt Foundation of Germany 
and partially by the National Fundamental Science Foundation of Vietnam}
\dedicatory{ To my father }
\abstract{This book is intended  to provide a quick introduction to the 
subject. The exposition is scheduled in the sequence, as possible for more 
understanding for beginners. The author exposed a K-theoretic approach 
to study group C*-algebras: started in the elementary part, with one 
example of 
description of the structure of C*-algebra of the group of affine 
transformations of the real straight line, continued then for some special 
classes 
of solvable and nilpotent Lie groups. In the second advanced part, he 
introduced the main tools of the theory. In particular, the conception of 
multidimensional geometric quantization and the index of group 
C*-algebras were created and developed. }  \endabstract

\maketitle
\vfill\eject

\centerline{\Large\bf Preface}
\vskip 1truecm
The main idea to create this book is to provide a quick introduction to
the subject from the lowest level of beginners up to the actual research
work level of working researchers. Certainly, we can not collect all
together. So that we must do some choice. Our main choice is based on
experiences with our research group on "K-theory, Harmonic Analysis and
Mathematical Physics". From the teaching experience with postgraduate
students in our group, we decide to expose the material in such a sequence
that, any time the reader can stop over to take out the corresponding
research problem. This means also that in the introduction, which is
reserved to the experts, we introduce the main ideas. Perhaps these
readers need not to read in detail the rest of the book. Other readers are
the beginners. They need to start from some very concrete examples to
illustrate the main ideas. So they should read from the part one. They
could stop at the end of this part to do by themselves some research for
other classes of examples. The third possibility is reserved to the
readers who need to known only the results, what are new in the general
theory. They could read the introduction chapter and then go directly to
the part II of the advanced theory. And finally, for the fourth class of
readers who want known all in the subject. They should read the chapters
in the same sequence as scheduled in the book. We hope that many of
readers could take contribution in perfecting the theory to solve the
problem which is very important from the point of view of theory and of
applications in physics and in mathematics. 

The author is very much indebted to his teacher, Prof. Dr. A. A.
Kirillov, for introducing him to the subject. Prof. Dr. G. G. Kasparov
discussed with him many times on the subject during 1977.
After returned from Moscow
University (in 1977) he was happy to meet Prof. Dr. P. Cartier, who
has totally supported him to develop research in this domain and
arranged for him a one-year visit (1983) at IHES Bures s/Yvette France,
where
he could work with Prof. Dr. A. Connes. Discussions with Prof. Dr. A.
Connes were not so much, but by which the author has believed some
close relation of the subject with Non Commutative Geometry.
Prof. Dr. A. Bak and Prof. Dr.
J. Cuntz provided him a possibility to work in Germany as a Humboldt
fellow (1991-1993). By the way, he could work in the nice conditions of
Alexander von Humboldt Foundation and German
universities Bielefeld and Heidelberg.
Prof. Dr. K. H. Hofmann arranged for him some visits to Technische
Hochschule Darmstadt as visiting C3 professor and completely supported
him in research. During this time and especially
during the time when the author prepared 3 invited Lectures for annual session
(Summer 1993) of the ``Seminar Sophus Lie", the idea about a book on the subject
was appeared.
The International Centre for Theoretical Physics has agreed to provide
him a scientific stay (August - October 1996) to support his project to
complete writing this book.

It is a great pleasure for the author to express the   deep and sincere
thanks to them and the
institutions for these important, effective and productive help and support.

{\parindent=0.75\hsize
Author         }

\vfill\eject
\tableofcontents

\chapter{Introduction}

The main problem we are interested in is how to characterize groups 
and their group algebras. This problem is well solvable with the 
complex representation theory of finite or compact groups. 
It seems to be quite difficult for locally compact groups. We focus
our attention, in this survey only on the topological method 
of characterizing the group C*-algebras.

\section{The Scope and an Example}
\subsection{The Problem}

Let us first of all consider a finite group $G$, $\vert G \vert \leq \infty$.
It is easy to see that the group $G$ can be included in some (co-)algebras, 
more precisely some Hopf bialgebra, for example in its complex group 
(Hopf bi-)algebra, $G \hookrightarrow {\Bbb C}[G]$, $$g 
\in G \mapsto \sum_{g' \in G} \delta_g(g') g'\in {\Bbb C}[G],$$   
which consists of the formal linear combinations of form 
$\sum_{g \in G} c_g g$,
and where $\delta_g$ is the Kronecker symbol corresponding to $g$. 
It is well-known that the group representation theory of $G$ is equivalent 
to the algebra representation theory of ${\Bbb C}[G]$. The last one is 
practically more flexible to describe. With each representation $\pi$ of 
this group one considers the Fourier-Gel'fand 
transforms of the algebra ${\Bbb C}[G]$
 to the matrix algebra, corresponding to representation $\pi$,
$$\sum_{g \in G} c_g g \mapsto \sum_{g\in G} c_g \pi(g).$$
Let us denote $\hat{G}$ the dual of $G$, i.e. the set of equivalence classes
of irreducible representations of $G$. We normally identify it with some
set of representatives of equivalence classes.
It is well-known that :

\begin{enumerate}
\item[a)]
The set $\hat{G}$ is finite, i.e. there is only finite number of nonequivalent
irreducible representations, say $\pi_1, \dots, \pi_N,$
\item[b)]
Each irreducible representation is finite dimensional, say of dimension $n_i,
i= 1,2,\dots,N,$
\item [c)]
This algebra ${\Bbb C}[G]$, by using the Fourier-Gel'fand transformation,
is isomorphic to the finite Cartesian product
of matrix algebras
$${\Bbb C}[G] \cong \prod_{i=1}^N \Mat_{n_i}({\Bbb C}).$$
\end{enumerate}

This means that the structure of the group algebra ${\Bbb C}[G]$ and therefore
of the group $G$ is well defined, if we could for $G$ do :

\begin{enumerate}
\item[1)]
a good construction of all the irreducible representations $\pi_i, i=1,\dots ,n$
of $G$ and
\item[2)]
the Fourier-Gel'fand transformation, realizing the above cited isomorphism. 
\end{enumerate}

{\it The problem is to extend this machinery to infinite, say locally compact,
topological groups. }

Let us from now on, consider a locally compact group $G$ and consider some
appropriate group algebras. The group algebras ${\Bbb C}[G]$ for $G$ as an 
abstract group is not enough to define the structure of $G$. We must find a
more effective group algebra. 

It is well known that for any locally compact group $G$ one must in place of 
the general linear representations consider the unitary ones. It is related 
with the fact that in general case one must consider also the infinite 
dimensional representations, which are not every time completely reducible. 
The unitary representations however are completely irreducible.
In case of a locally compact group $G$ 
there is a natural left-(right-)invariant Haar measure $dg$. The space
$L^2(G):= L^2(G,dg)$ of the square-integrable functions plays an important 
role in harmonic analysis. If the group is of type I, $L^2(G)$ admits a 
spectral decomposition with respect to
 the left and right regular representations into a sum of the direct sum 
(the so called discrete series) and/or the direct integral (the continuous 
series) of irreducible unitary representations. The space $L^1(G)=L^1(G,dg)$ 
of the  functions with integrable module plays a crucial role. 
With the well-defined convolution product, 
$$\varphi,\psi \in L^1(G) \mapsto \varphi * \psi \in L^1(G);$$
$$(\varphi * \psi)(x) := \int_G \varphi(y)\psi(y^{-1}x)dy
$$ it becomes a Banach algebra.
There is also a well-defined Fourier-Gel'fand transformation on $L^1(G)$,
$$\varphi \in L^1(G,g) \mapsto \hat{\varphi},$$
$$\hat{\varphi}(\pi) := \pi(\varphi) = \int_G \pi(x)\varphi(x)dx.$$
There is a one-to-one correspondence between the (irreducible) unitary 
representations of $G$ and the  non-degenerate (irreducible) *-representations 
of the involutive Banach algebra $L^1(G)$. The general theorems of the spectral theory of the 
representations of $G$ are then proved with the help of an appropriate 
translation in to the corresponding theory of $L^1(G)$, for which one can use more 
tools from functional analysis and topology. One can also define on $L^1(G)$ 
an involution $\varphi \mapsto \varphi^*$,
$$\varphi^*(g) := \overline{\varphi(g^{-1})}.$$ 
However the norm of the involutive Banach algebra $L^1(G)$ is not regular, 
i.e. in general
$$\Vert a^*a\Vert_{L^1(G)} \ne \Vert a\Vert^2_{L^1(G)}\;.$$
It is therefore more useful to consider the corresponding regular norm 
$\Vert.\Vert_{C^*(G)}$,
$$\Vert \varphi  \Vert_{C^*(G)} := \sup_{\pi \in \hat{G}} \Vert \pi(\varphi)
\Vert$$ and its completion $C^*(G)$. The spectral theory of
unitary representations
of $G$ is equivalent to the spectral theory of non-degenerate *-representations
 of the C*-algebra $C^*(G)$. The general theorems of harmonic analysis say that
the structure of $G$ can be completely definite by the structure of $C^*(G)$.
One poses therefore the problem of description of 
the structure of the C*-algebras of locally compact groups. This means that 
we must answer to the questions:

\begin{enumerate}
\item[1)] How to realize the irreducible unitary representations of the 
locally compact group $G$.
\item[2)] How to describe the images of the Fourier-Gel'fand transformation
and in particular, of the inclusion of $C^*(G)$ into some ``continuous'' 
product of the algebras ${\cal L}({\cal H}_\pi)$, $\pi\in \hat{G}$ of bounded
operators in the separable Hilbert space ${\cal H}_\pi$ of representation 
$\pi$.
\end{enumerate}

To see that this is a good setting the problem for finite group to the 
locally compact groups, let us consider these questions for the compact 
group case. Consider for the moment a compact group $G$. 
For compact group, all representations are unitarizable, 
i.e. are equivalent to some unitary ones. It is well-known also that :

\begin{enumerate}
\item[a)]
The family of irreducible representations is not more than countable.
\item[b)]
Each irreducible representation is finite dimensional, say 
$n_i, i= 1,\dots, \infty$ and there is some good realizations of these 
representations, say in tensor spaces, or last time, in cohomologies.
\item[c)]
The Fourier-Gel'fand transformation gives us an isomorphism
$$C^*(G) \cong {\prod_{i=1}^\infty}' \Mat_{n_i}({\Bbb C}),$$
where the prime in the product means the subset of "continuous vanishing at 
infinity" elements. \end{enumerate}

This means that in compact group case the group C*-algebra  plays the 
same role as the group algebra of finite groups.

Let us now return to the general case of locally compact group. 

{\it The main problem is how to describe the group algebra in general
, and in particular the C*-algebra $C^*(G)$. }

In general the problem of describing the structure of
C*-algebras of non-compact groups rests open up-date. 
This review outlines only the well-known cases, where is a nice 
interaction of the methods from noncommutative geometry, say  Orbit Method, 
category ${\cal O}$, KK-theory, deformation quantization, cyclic theories,...
We restrict mainly onto the case of Lie groups.

\subsubsection{Analytic Method}
The first nontrivial example is the group $\SL_2({\Bbb C})$. Its C*-algebra
was studied by J.M.G. Fell in 1961 in \cite{fell}. He described exactly the 
Fourier-Gel'fand transforms of $C^*(G)$ as some C*-algebra of sections 
of a continuous field of operator algebras over the dual. Many other 
mathematicians attempted to generalize his beautiful but complicate 
analytic result to other groups. Nevertheless, until the moment the only 
groups, the structure of whose C*-algebras were explicitly described are:
the Abelian or compact groups and a few semi-simple Lie groups, say 
$\SL_2({\Bbb R})$ and its universal covering $\widetilde{\SL_2({\Bbb R})}$, the de Sitter
group $\Spin(4,1)$ and recently a family $G(p,q,\alpha)$ of two step solvable
Lie groups (see [De], [F], [M], [KM], [BM] and \cite{wang1}). A fair amount
is known about the C*-algebras of nilpotent Lie groups (see [P]), including
the Heisenberg groups. The C*-algebra of the Euclidean motion group were
studied by Evans \cite{evans}. Also P. Green \cite{green} proposed some 
another  analytic method for studying the C*-algebras of several solvable
Lie groups. The result are given very slowly and spectacularly.

One need therefore to develop another method, say to obtain some topological
invariants, which will be described in the rest of this paper.

\subsubsection{K-Theory Approach}

The very useful K-functor for our approach is the operator KK-functor of 
G. G. Kasparov \cite{kasparov1}, generalizing the BDF K-functor \cite{bdf1},
which characterizes the isomorphic classes of short exact sequences of 
C*-algebras. 

We are trying to decompose our C*-algebras into some towers of ideals 
and step-by-step define the associated extensions by KK-functors or 
their generalizations. The resulting invariants form just our index. 
This idea was proposed in \cite{diep1} and developed in \cite{diep2} 
for a large class of type I C*-algebras. Hence, there are two general problems:
\begin{enumerate}
\item[(1)] Find out the C*-algebras which can be characterized by the 
well-known K-functors, say by the operator K-functors.
\item[(2)] Generalize the theory of K-functors in such a way that they are applicable for a large class of C*-algebras. 
\end{enumerate}

Concerning the first problem, we propose \cite{diep8} a general 
construction and some reduction procedure of the K-theory invariant 
$\Ind\, C^*(G)$ of group C*-algebras. Using the orbit method 
\cite{kirillov},
\cite{diep4} - \cite{diep7}, we reduces $Index~C^*(G)$ to a family of 
Connes' 
foliation C*-algebras  indices $Index~C^*(V_{2n_i},{\cal F}_{2n_i})$, 
see \cite{connes1}-\cite{connes2}, by a family of KK-theory invariants. 
Using some generalization of the Kasparov type condition (treated by G.G. 
Kasparov in the nilpotent Lie group case \cite{kasparov2}), we reduces 
every \newline $Index C^*(V_{2n_i},{\cal F}_{2n_i})$ to a family of 
KK-theory invariants of the same type valuated in KK(X,Y) type groups. 
The last ones are in some sense computable by using the cup-cap product 
realizing the Fredholm operator indices.  

To demonstrate the idea, we consider the C*-algebra of the group of 
affine transformations of the real straight line, 
but first of all we need some new K-functor tool. 
It is described in the next two subsections.

\subsection{BDF K-Homology functor}

Let us recall in this subsection the well-known BDF K-functor ${\cal E}xt$.
The main reference is \cite{bdf1}.
Denote by $C(X)$ the C*-algebra of continuous complex-valued functions
over a fixed metrizable compact $X$, ${\cal H}$ a fixed separable Hilbert
space over complex numbers, ${\cal L}({\cal H})$ and ${\cal K}({\cal H})$
the C*-algebras of bounded and respectively, compact linear operators in
${\cal H}$. An extension of C*-algebras means a short exact sequence of
C*-algebras and *-homomorphisms of special type
$$0 \longrightarrow {\cal K}({\cal H}) \longrightarrow {\cal E} 
\longrightarrow C(X)\longrightarrow 0.$$ Two extensions are by definition
equivalent iff there exists an isomorphism
$ \psi : {\cal E_1} \longrightarrow {\cal E_2}$  and
its restriction
$\psi\vert_{{\cal K}({\cal H}_1)} : {\cal K}({\cal H}_1) 
\longrightarrow {\cal K}({\cal H}_2)$ such that the following diagram
is commutative
$$ 
\begin{array}{ccccccccc}
0 & \longrightarrow & {\cal K}({\cal H}_1 ) & \longrightarrow & {\cal E}_1 & \longrightarrow & C(X) & \longrightarrow & 0\\
 &    & \Big\downarrow\vcenter{%
     \rlap{$\psi\vert_.$}} & & \Big\downarrow\vcenter{\rlap{$\psi$}}
  &  & \Big\Vert &  &  \\
0 & \longrightarrow & {\cal K}({\cal H}_2) & \longrightarrow  & {\cal E}_2 
& \longrightarrow & C(X) &\longrightarrow & 0
\end{array}
$$

There is a canonical universal
extension of C*-algebras
$$0 \longrightarrow {\cal K}({\cal H}) \longrightarrow {\cal L}({\cal H})
\longrightarrow {\cal A}({\cal H}) \longrightarrow 0,$$
the quotient algebra ${\cal A}({\cal H}) \cong {\cal L}({\cal H})/{\cal K}
({\cal H})$ is well-known as the Calkin algebra. By the construction of
fiber product, there is one-to-one correspondence between the extensions of
type $$0\longrightarrow {\cal K}({\cal H}) \longrightarrow {\cal E} 
\longrightarrow C(X) \longrightarrow 0$$ and the unital monomorphisms of type
$$\varphi : C(X) \hookrightarrow {\cal A}({\cal H}).$$ 
Thus we can identify the 
extensions with the inclusions of $C(X)$ into ${\cal A}({\cal H})$. Because
\cite{kirillov} all separable Hilbert spaces are isomorphic and
the automorphisms of ${\cal K}({\cal H})$ are inner and 
$$\Aut{\cal K}({\cal H}) \cong {\cal P}{\cal U}({\cal H}),$$
the projective unitary group, where ${\cal U}({\cal H})$
denotes the unitary operator group,
we can identify the equivalences classes of extensions with the unitary 
conjugacy classes of unital inclusions of $C(X)$ into the Calkin algebra: 
Two extensions $\tau_1$ and $\tau_2$ are equivalent 
iff there exists a unitary operator 
$U : {\cal H}_1 \longrightarrow {\cal H}_2$, 
such that $\tau_2 = \alpha_U \circ \tau_1 $, 
where by definition 
$\alpha_U : {\cal A}({\cal H}_1) \longrightarrow {\cal A}({\cal H}_2)$ 
is the isomorphism obtained from the inner isomorphism 
$$U.(-).U^{-1} : {\cal L}({\cal H}_1) \longrightarrow {\cal L}({\cal H}_2 ).$$
Extension $\tau : C(X) \hookrightarrow {\cal A}({\cal H})$ 
is called trivial iff there exists a unital inclusion 
$\sigma : C(X) \hookrightarrow  {\cal L}({\cal H})$ 
such that $\tau = \pi \circ \sigma, $ 
where $\pi : {\cal L}({\cal H}) \longrightarrow {\cal A}({\cal H}) = 
{\cal L}({\cal H})/{\cal K}({\cal H})$ is the canonical quotient map.  
This inclusion $\tau$ corresponds to the split short exact sequence.
The sum of two extensions 
$\tau_i : C(X) \hookrightarrow {\cal A}_i, i= 1,2$ 
is defined as the extension 
$$\tau_1 \oplus \tau_2 : C(X) \hookrightarrow {\cal A}({\cal H}_1) \oplus
{\cal A}({\cal H}_2) \hookrightarrow 
{\cal A}({\cal H}_1 \oplus 
{\cal H}_2).$$ 
This definition is compatible also with the equivalence classes of extensions. 
In \cite{bdf1} the authors proved that:

\begin{enumerate}
\item[1)] 
The equivalence class of trivial extension is the identity element  
with respect to this sum.
\item[2)] 
For every metrizable compact $X$, 
the set ${\cal E}xt_1(X)$ of the equivalence classes of extensions is 
an Abelian group. One defines the higher groups by 
${\cal E}xt_{1-n}(X) := {\cal E}xt_1({\Bbb S}^{n}\wedge X)$, $n=0,1,2,\dots,$
\item[3)] 
${\cal E}xt_*$ is a generalized K-homology. In particular, the group 
${\cal E}xt_1(X)$ is dependent only of the homotopy type of $X$ and 
there is a homomorphism 
$$ Y_\infty : {\cal E}xt_1(X) \longrightarrow 
\Hom_{\Bbb Z}(K^{-1}(X), {\Bbb Z})$$ 
which will be an isomorphism if $X \subset {\Bbb R}^3$. 
\end{enumerate}

This K-homology is well developed and fruitfully applicable. 
It has many application in operator theory and in our problem of 
characterizing the group C*-algebras.
Let us demonstrate this in the first example of the group of affine 
transformations of the real straight line.

\subsection{Topological Invariant Index }

Let us in this subsection denote by $G$ the group of all affine transformations of the real straight line.

\begin{thm} Every irreducible unitary representation of group $G$ is
unitarilly equivalent to one of the following mutually nonequivalent
representations: \begin{enumerate} \item[a)] the representation $S$,
realized in the space $L^2({\Bbb R}^*, \frac{dx}{ \vert x \vert})$, where
${\Bbb R}^* := {\Bbb R} \setminus (0)$, and acting in according with the
formula $$(S_gf)(x) = e^{\sqrt{-1}bx} f(ax),\text{ where } g =
\begin{pmatrix} \alpha & b \\ 0 & 1\end{pmatrix}.$$ \item[b)] the
representation $U^\varepsilon_\lambda$, realized in ${\Bbb C}^1$ and given
by the formula $$U_\lambda^\varepsilon (g) = \vert \alpha
\vert^{\sqrt{-1}\lambda}.(\sgn\alpha)^\varepsilon, \text{ where }\lambda
\in {\Bbb R}; \varepsilon = 0,1.$$ \end{enumerate} \end{thm} \begin{pf}
See \cite{gelfandnaimark}. \end{pf}

This list of all the irreducible unitary representations gives the
corresponding list of all the irreducible non-degenerate unitary
*-representations of the group C*-algebra $C^*(G)$. In \cite{diep1} it
was proved that

\begin{thm} The group C*-algebra with formally adjoined unity
$C^*(G)^\sim$ can be included in a short exact sequence of C*-algebras and
*-homomorphisms $$ 0 \longrightarrow {\cal K} \longrightarrow C^*(G)^\sim
\longrightarrow C({\Bbb S}^1 \vee {\Bbb S}^1) \longrightarrow 0,$$ i.e.
the C*-algebra $C^*(G)^\sim$, following the BDF theory, is defined by an
element, called the index and denoted by $Index\,C^*(G)^\sim$, of the
groups ${\cal E}xt({\Bbb S}^1 \vee {\Bbb S}^1) \cong {\Bbb Z} \oplus {\Bbb
Z}$. \end{thm} \begin{pf} See \cite{diep1}. \end{pf}

The infinite dimensional representation $S$ realizes the inclusion said above.
 Since $${\cal E}xt({\Bbb S}^1 \vee {\Bbb S}^1) \cong \Hom_{\Bbb 
Z}(\pi^1({\Bbb
 S}^1\vee {\Bbb S}^1, {\Bbb Z})$$ it realized by 
a homomorphism from $\pi^1({\Bbb S}^1 
\vee {\Bbb S}^1)$ to ${\Bbb C}^*$. Since the isomorphism $$Y_\infty : 
{\cal E}xt({\Bbb S}^1 \vee {\Bbb S}^1)\cong \Hom_{\Bbb Z}(\pi^1({\Bbb 
S}^1 \vee 
{\Bbb S}^1), {\Bbb Z})$$ is obtained by means of computing the indices
and because the general type of elements of $\pi^1({\Bbb S}^1 \vee 
{\Bbb S}^1)$ is $g_{k,l} = [g_{0,1}]^k[g_{1,0}]^l$, $k, l \in {\Bbb Z}$, 
we have $$\Ind\,(g_{k,l}) = k.\Ind\, T(g_{1,0}) + l.\Ind\,T(g_{0,1}),$$ 
where $T$ is the the *-isomorphism corresponding to $S$. It is enough 
therefore to compute the pair of two indices $\Ind\, T(g_{1,0})$ and 
$\Ind\, T(g_{0,1})$. The last ones are directly computed by the indices 
of the corresponding Fredholm operators. 

\begin{thm}
$$Index\, C^*(G) = (1,1) \in {\cal E}xt({\Bbb S}^1 \vee {\Bbb S}^1) 
\cong {\Bbb Z} \oplus {\Bbb Z}.$$
\end{thm}
\begin{pf} See \cite{diep1}. 
\end{pf}

Let us now go to the general situation. To do this we must introduce
also some preparation about, first of all, the construction of
irreducible unitary representations, we mean the orbit method, then
a method of decomposing the C*-algebra into a tower of extensions and
lastly compute the index with the help of the general KK-theory.

\section{Multidimensional Orbit Methods}

Let us in this section consider the problem of realization of
irreducible unitary representations of Lie groups. There are two
versions of the orbit method; one is the multidimensional quantization,
the other is the infinitesimal orbit method, related with the so
called category ${\cal O}$.

\subsection{Multidimensional Quantization}

The orbit method can be constructed from the point of view of the theory 
of holomorphly induced representations and also from the point of view
of the ideas of quantization from physics.

\subsubsection{Construction of Partially Invariant Holomorphically 
Induced Representations}

Let us consider now a connected and simply connected Lie group $G$ 
with Lie algebra ${\frakt g} := \Lie(G)$. Denote by ${\frakt g}_{\Bbb C}
$ the complexification of ${\frakt g}$. The complex conjugation in 
the Lie algebra will be also denoted by an over-line sign. 
Consider the dual space ${\frakt g}^*$ to the Lie algebra ${\frakt g}$. 
The group $G$ acts on itself by the inner automorphisms
$$A(g):= g.(.).g^{-1} : G \longrightarrow G,$$ for each $g\in G$, 
conserving the identity element $e$ as some fixed point. It follows 
therefore that the associated adjoint action $A(g)_* $ maps 
${\frakt g} = T_eG$ into itself and the co-adjoint action 
$K(g) := A(g^{-1})^* $ 
maps the dual space 
${\frakt g}^*$ 
into itself. The orbit space 
${\cal O}(G):= {\frakt g}^*/G$ 
is in general a bad topological space, namely non- Hausdorff. Consider 
one orbit 
$\Omega\in{\cal O}(G)$
and an element 
$F \in {\frakt g}^*$ in it. The stabilizer is denote by 
$G_F$, 
its connected component by $(G_F)_0$ and its Lie algebra by 
${\frakt g}_F := \Lie(G_F)$. It is well-known that
$$
\begin{array}{ccc}
G_F & \hookrightarrow & G\\
    &                 & \Big\downarrow\\
    &                  & \Omega_F
    \end{array}
$$
is a principal bundle with the structural group $G_F$. 
Let us fix some {\it connection in this principal bundle, } 
\index{connection on 
principal bundle} i.e. some {\it trivialization } \index{trivialization}
of this bundle, see \cite{sulankewintgen}.  
We want to construct representations in some cohomology spaces 
with coefficients in the sheaf of sections of some vector bundle 
associated with this principal bundle. 
It is well know \cite{sulankewintgen} that every vector bundle 
is an induced one with respect to some representation of the structural group
in the typical fiber.
It is natural to fix some unitary representation 
$\tilde{\sigma}$ 
of $G_F$ such that its kernel contains $(G_F)_0$, the character
$\chi_F$ of the connected component of stabilizer
$$\chi_F(\exp{X}) := \exp{(2\pi\sqrt{-1}\langle F,X\rangle )}$$
and therefore the differential 
$D(\tilde{\sigma}\chi_F)=\tilde{\rho}$ 
is some representation  of the Lie algebra ${\frakt g}_F$. 
We suppose that the representation $D(\tilde{\rho}\chi_F)$ was extended to 
the complexification $({\frakt g}_F)_{\Bbb C}$. 
The whole space of all sections seems to be so large for the construction 
of irreducible unitary representations. 
One consider the invariant subspaces with the help of some so 
called polarizations.

\begin{defn}
We say that a triple 
$({\frakt p}, \rho,\sigma_0)$ is some 
$(\tilde{\sigma},F)$-{\it polarization, } \index{polarization}
iff :
\begin{enumerate}
\item[a)]
${\frakt p}$ is some subalgebra of the complex Lie algebra 
$({\frakt g})_{\Bbb C}$, containing ${\frakt g}_F$.
\item[b)]
The subalgebra ${\frakt p}$ is invariant under the action 
of all the operators of type $Ad_{{\frakt g}_{\Bbb C}}x, x\in G_F$.
\item[c)]
The vector space ${\frakt p} + \overline{\frakt p}$ 
is the complexification of some real subalgebra 
${\frakt m} = ({\frakt p} + \overline{\frakt p}) \cap {\frakt g}$. 
\item[d)]
All the subgroups $M_0$, $H_0$, $M$, $H$ are closed. where by definition 
$M_0$ (resp., $H_0$) is the connected subgroup of 
$G$ with the Lie algebra ${\frakt m}$ 
(resp., ${\frakt h} := {\frakt p} \cap {\frakt g}$) and $M:= G_F.M_0$,
$H:= G_F.H_0$.
\item[e)]
$\sigma_0$ is an irreducible representation of $H_0$ in some Hilbert space 
$V$ such that : 
1. the restriction $\sigma|_{G_F \cap H_0}$ 
is some multiple of the restriction 
$\chi_F.\tilde{\sigma}|_{G_F \cap H_0}$, 
where by definition $\chi_F(\exp X) := \exp{(2\pi\sqrt{-1}\langle F,X\rangle )}$;
2. under the action of $G_F$ on the dual $\hat{H}_0$, 
the point $\sigma_0$ is fixed.
\item[f)]
$\rho$ is some representation of the complex Lie algebra 
${\frakt p}$ in $V$ , 
which satisfies the E. Nelson conditions for $H_0$ and 
$\rho|_{\frakt h} = D\sigma_0$.
\end{enumerate}
\end{defn}

Let us recall that R. Blattner introduced the notion of mixed manifold 
of type $(k,l)$, 
see for example \cite{kirillov}. We consider the fiber bundle, 
the base of which is some type $(k,l)$ mixed manifold and the fibers of 
which are smooth $m$ dimensional manifold. 
We say that this fiber bundle is some {\it mixed manifold of type } 
\index{manifold!mixed -} $(k,l,m)$.

\begin{thm}
Let us keep all the introduced above notation of $\Omega_F$, $\tilde{\sigma}$,
$G_F$,etc. and 
let us denote $\chi_F$ the character of the group $G_F$ such that 
$D\chi_F = 2\pi\sqrt{-1}F|_{{\frakt g}_F}$. Then :
\begin{enumerate}
\item[1)]
On the K-orbit $\Omega_F$ there exists a structure of some mixed manifold 
of type $(k,l,m)$, where $$k = \dim G - \dim M,$$
$$l = {1 \over 2}(\dim M - \dim H),$$
$$m = \dim H - \dim G_F.$$
\item[2)]
There exists some irreducible unitary representation 
$\sigma$ of the group $H$ such that its restriction $\sigma|_{G_F}$ 
is some multiple of the representation 
$\chi_F.\tilde{\sigma}$ and $\rho|_{\frakt h} = D\sigma$.
\item[3)]
On the $G$-fiber bundle ${\cal E}_{\sigma|_{G_F}} = G \times_{G_F} V$ 
associated with the representation $\sigma|_{G_F}$, there exists a structure
of a partially invariant and partially holomorphic Hilbert vector 
$G$-bundle ${\cal E}_{\sigma,\rho}$ 
such that the natural representation of $G$ on the space of (partially 
invariant and partially holomorphic) sections is equivalent to the 
representation by right translations of $G$ in the space 
$C^\infty (G; {\frakt p}, \rho, F, \sigma_0)$ of $V$-valued 
$C^\infty$-functions on $G$ satisfying the equations 
$$f(hx) = \sigma(h)f(x), \forall h \in H, \forall x\in G,$$
$$L_Xf + \rho(X)f = 0, \forall X\in \overline{\frakt p},$$
where $L_X$ denotes the Lie derivative along the vector field 
$\xi_X$ on $G$, corresponding to $X$.
\end{enumerate}
\end{thm}

\begin{pf} The first assertion is clear. 
The second one can be deduced from the remark that the formula
$$(x,h) \mapsto I_{V'} \otimes \chi_F.\tilde{\sigma})(x).\sigma_0(h)$$
defines an irreducible representation of the direct product 
$G_F \times H_0$ which is trivial on the kernel of the surjection 
$$G_F \times H_0 \longrightarrow G_F.H_0.$$ 
This point is essential in the sense that with the assumption about fixed 
point property of $\sigma_0$ we can ignore the Mackey obstacle, 
appeared when we take the representations which are multiple of 
some representations $\tilde{\sigma}\chi_F$ at the restriction to 
some normal subgroup. 
M. Duflo \cite{duflo1} considered two-fold covering to avoid this obstacle.
See \cite{diep5}.
\end{pf}

One can than apply the construction of unitarization 
$\overline{\cal E}_{\sigma,\rho}$ to obtain the corresponding 
unitary representation,
which is noted by $\Ind(G;{\frakt p},F,\rho,\sigma_0)$. 
One can define also the representations in cohomologies with coefficients 
in this sheaf of partially invariant and  partially holomorphic sections, 
which will be noted by $(L^2-\Coh)\Ind(G;{\frakt p},F,\rho,\sigma_0)$. 

\begin{rem} 
One introduces some order in the set of all 
$(\tilde{\sigma},F)$- polarizations 
$$({\frakt p},\rho,\sigma_0) \leq ({\frakt p}',\rho',\sigma'_0) 
\Longleftrightarrow {\frakt p} \subseteq {\frakt p}', \sigma'_0|_{H_0} 
\simeq \sigma_0, \rho'|_{\frakt p} \simeq \rho.$$ 
To have some irreducible representation, 
one must take the maximal polarizations in this construction. 
It is interesting that this representations are coincided with the 
representations appeared from the geometric quantization.
\end{rem}

\subsubsection{Multidimensional Geometric Quantization} 

Let us now consider the general conception of multidimensional 
geometric quantization. 
Consider a symplectic manifold $(M,\omega)$, 
i.e. a smooth manifold equipped with a non-degenerate closed skew-symmetric 
differential 2-form $\omega$. The vector space $C^\infty (M,\omega)$, 
with respect to the Poisson brackets  
$$f_1, f_2 \in C^\infty \mapsto \{ f_1,f_2\} \in C^\infty(M,\omega)$$
become an infinite dimensional Lie algebra. 

\begin{defn} 
A {\it procedure of quantization }  \index{procedure of quantization}
is a correspondence associating to each {\it classical quantity } 
\index{classical quantity} $f\in C^\infty(M)$ a {\it quantum quantity } 
\index{quantum quantity}
$Q(f) \in {\cal L}({\cal H})$, i.e. a continuous, perhaps unbounded,
normal operator, which is auto-adjoint if $f$ is a real-valued function,
in some Hilbert space ${\cal H}$, such that
$$Q(\{f_1,f_2\}) = {i \over \hbar}[Q(f_1),Q(f_2)],$$
$$Q(1) = Id_{\cal H},$$
where $\hbar := h/2\pi$ is the normalized Planck constant, 
and $h$ is the unnormalized Planck's constant.
\end{defn}

Let us denote by ${\cal E}$ a fiber bundle into Hilbert spaces, 
$\Gamma$ a fixed connection conserving the Hilbert structure 
on the fibers; in other words,
If $\gamma$ is a curve connecting two points $x$ and $x'$, the parallel 
transport along the way $\gamma$ provides an scalar preserving isomorphism 
from the fiber ${\cal E}_x$ onto the fiber ${\cal E}_{x'}$. In this case we can
define the corresponding covariant derivative $\nabla_\xi$, $\xi \in Vect(M):=
Der\, C^\infty(M)$
in the space of smooth sections. One considers the invariant Hilbert space 
$L^2({\cal E}_{\rho,\sigma})$, which
is the completion of the space $\Gamma({\cal E}_{\rho,\sigma})$
of square-integrable partially invariant and partially holomorphic sections.

Suppose from now on that $M$ is a homogeneous $G$-space.
Choose a trivialization $\Gamma$ of the principal bundle 
$G_x \rightarrowtail G \twoheadrightarrow M$, 
where $G_x$ is the stabilizer of the fixed point $x$ on $M$.
Let us denote by $L_\xi$ the Lie derivation corresponding to the vector field
$\xi\in Vect(M)$. Let us denote by $\beta\in \Omega^1(M)$ the form of affine
connection on ${\cal E}$, corresponding to the
connection $\Gamma$ on the principal bundle. It is more comfortable to
consider the normalized connection form
$\alpha(\xi) = {\hbar \over \sqrt{-1}}\beta(\xi)$,
the values of which are anti-auto-adjoint operators on fibers. One has therefore
$$\nabla_\xi = L_\xi + {\sqrt{-1}\over \hbar}\alpha(\xi),$$
see for example \cite{sulankewintgen} for the finite dimensional case.

For each function $f \in C^\infty(M)$ one denotes $\xi_f$ the corresponding 
Hamiltonian vector field, i.e.
$$i(\xi_f)\omega + df = 0.$$

\begin{defn}
We define the geometrically quantized operator $Q(f)$ as
$$Q(f) := f + {\hbar\over \sqrt{-1}}\nabla_{\xi_f} = f + {\hbar\over\sqrt{-1}}
L_{\xi_f} +\alpha(\xi_f).$$
\end{defn}

\begin{thm}
The following three conditions are equivalent.
\begin{enumerate}
\item[1)] 
$$\xi\alpha(\eta)-\eta\alpha(\xi)-\alpha([\xi,\eta])+{\sqrt{-1}\over\hbar}
[\alpha(\xi),\alpha(\eta)] = -\omega(\xi,\eta).Id; \forall \xi,\eta.$$
\item[2)]
The curvature of the affine connection $\nabla$ is equal to $-{\sqrt{-1}\over
\hbar}\omega(\xi,\eta).Id$, i.e.
$$[\nabla_\xi , \nabla_\eta] - \nabla_{[\xi,\eta]} = -{\sqrt{-1}\over\hbar}
\omega(\xi,\eta).Id; \forall \xi, \eta.$$
\item[3)]
The correspondence $f\mapsto Q(f)$ is a quantization procedure.
\end{enumerate}
\end{thm}
\begin{pf} See \cite{diep6}.
\end{pf}

Suppose that the Lie group   $G$ act  on $M$ by the symplectomorphisms. 
Then each element  $X$ of the Lie algebra ${\frakt g}$ corresponds to 
one-parameter
subgroup $\exp{(tX)}$ in $G$, which acts on $M$. Let us denote by $\xi_X$
the corresponding strictly Hamiltonian vector field. Let us denote also $L_X$
the Lie derivation along this vector field. We have 
$$[L_X,L_Y] = L_{[X,Y]},$$ and $$L_Xf = \{f_X,f\}.$$

Suppose that $f_X$ depends linearly on $X$. One has then a 2-cocycle of
the action
$$c(X,Y) := \{f_X,f_Y\} - f_{[X,Y]}.$$

\begin{defn}
We say that the action of $G$ on $M$ is {\it flat } \index{flat action}
iff this 2-cocycle is trivial.
\end{defn}

In this case we obtain from the quantization procedure a representation 
$\wedge$ of the Lie algebra ${\frakt g}$
by the anti-auto-adjoint operators 
$$X \mapsto {\sqrt{-1}\over\hbar}Q(f_X)$$
and also a representation of ${\frakt g}$ by the functions
$$X \mapsto f_X.$$ If the E. Nelson conditions are satisfied, we have a unitary
representation of the universal covering of the group $G$.

\begin{thm}
The Lie derivative of the partially invariant and holomorphically induced
representation $\Ind(G;{\frakt p},F,\rho,\sigma_0)$ of a connected Lie 
group $G$
is just the representation obtained from the procedure of multidimensional
geometric quantization, corresponding to a fixed connection $\nabla$ of 
the partially invariant partially holomorphic induced unitarized bundle
$\overline{\cal E}_{\sigma,\rho}$, i.e.
$$\Lie_X(\Ind(G;{\frakt p},F,\rho,\sigma_0)) = {\sqrt{-1}\over\hbar}Q(f_X).$$
\end{thm}
\begin{pf} See \cite{diep4}. 
\end{pf}

\begin{rem}
The multidimensional version of the orbit method was developed independently
by the author in language of multidimensional quantization \cite{diep4} - 
\cite{diep7}
and by M. Duflo 
\cite{duflo1}, see also, \cite{kirillov} in the language
of Mackey method of small subgroups. The result
show that for the most connected Lie groups the construction gives us at least
a quantity of irreducible unitary representations, enough to decompose the
regular representations of $G$ in $L^2(G)$, i.e. enough to prove the
 Plancher\`el formula \cite{duflo1}.
\end{rem} 

\begin{rem} There are some reductions of this multidimensional 
quantization procedure to the radical or nil-radical of stabilizer
 of type $G_F$, see \cite{dongvui}, and lifting them to $U(1)$-coverings
\cite{vui1} - \cite{vui3}, \cite{dong1} - \cite{dong2}.
\end{rem}

\begin{rem} In \cite{diep10} the author proposed some method for 
common quantization for foliations, the fibers of which are the
K-orbits, and its relation with the integral Fourier operators.
\end{rem}

\subsection{Category ${\cal O}$ and globalization of Harish-Chandra modules }

The construction of irreducible unitary representations \newline
$(L^2-\Coh) \Ind (G;{ \frakt p}, F, \rho, \sigma_0)$ in the Hilbert space 
$L^2({\cal E}_{\rho,\sigma}) \cong L^2(G; {\frakt p},F,\rho,\sigma_0)$
can be in restricted case considered as some globalization of some so called
$({\frakt g}, K)$-module, i.e. ${L^2({\cal E}_{\rho, \sigma_0})}_{(K)}$ 
itself is some $({\frakt g},K)$-module,
where $K$ is some maximal compact subgroup of $G$. 
It is therefore interesting to consider these
$({\frakt g},K)$-module as some infinitesimal version of the orbit method.
Let us see this in this subsection. 

  \subsubsection{Admissible representations} 
Let us in this sub-subsection recall some result about the 
Borell-Weil-Bott-Kostant theorem and the construction of admissible 
representations of finite dimensional semi-simple Lie groups as $({\frakt g},K)$-modules, \cite{wolf},
\cite{milicic2}.

If $G$ is a compact connected Lie group, and $F\in {\frakt g}^*$ 
is a well-regular integral functional on its Lie algebra, then the stabilizer is a 
maximal torus $T$, (If the Harish-Chandra criterion for existence of 
discrete series holds, it is a compact Cartan subgroup.) and a choice of 
positive root system $\Phi^+ = \Phi^+({\frakt g},{\frakt t})$
defines a $G$-invariant complex manifold structure on $G/T$ in such a way that
$\sum_{\alpha\in\Phi^+}{\frakt g}_\alpha$
represents the holomorphic tangent space.
The character $\chi_\lambda, \lambda:= {\sqrt{-1}\over\hbar}F$ 
can be extended to a character of the stabilizer $G_F = T$, 
if the orbit is as usually supposed to be integral.
Let us denote in this case the induced bundle 
${\cal E}_{\rho,\sigma}$ simply by ${\Bbb E}_\lambda$
as in \cite{wolf}. 
It is the associated homogeneous holomorphic hermitian line bundle.
One writes ${\cal O}({\Bbb E}_\lambda)\longrightarrow G/T$ 
for the sheaf of germs of holomorphic sections of 
${\Bbb E}_\lambda \longrightarrow G/T$. The group $G$ acts every where,
including the cohomologies $H^q(G/T; {\cal O}({\Bbb E}_\lambda))$.
One denotes by $\rho  := {1\over 2}\sum_{\alpha\in\Phi^+} \alpha$
the half-sum of positive roots.
We cite from \cite{wolf} the Borel-Weil-Bott-Kostant theorem

\begin{thm}
If $\lambda + \rho$     is singular then every 
$H^q(G/T;{\cal O}({\Bbb E}_\lambda))$ is trivial. 
If $\lambda+\rho$ is regular, let $w$ denote the unique element
such that $$\langle w(\lambda+\rho),\alpha\rangle  > 0, \forall \alpha\in\Phi^+$$
and let $\ell(w)$ denote its length as a word in the simple root reflections.
Then 
\begin{enumerate}
\item[i)]
$H^q(G/T;{\cal O}({\Bbb E}_\lambda)) = 0$ for all $q\ne \ell(w)$, and
\item[ii)]
the action of $G$ in $H^q(G/T;{\cal O}({\Bbb E}_\lambda))$ is the 
representation with highest weight $w(\lambda + \rho) - \rho$.
\end{enumerate}
\end{thm}

This result was then extended for realizing {\it the discrete 
series representations
 of general semi-simple Lie groups. } \index{representations! discrete 
series -} It is well known that one can induce
from these discrete series representations of reductive part of parabolic
subgroups to obtain the {\it tempered admissible 
representations } \index{representations! tempered admissible - } of $G$. It 
was then
remarked that the representation of $G$ in $L^2({\cal E}_{\rho,\sigma})$
can be considered as the {\it globalization} \index{globalization} of some $({\frakt 
g},K)$-module , namely, $L^2({\cal E}_{\rho,\sigma})_{(K)}$.  

Lastly the {\it tempered admissible representations } 
\index{representations!tempered admisssible -}  are described in the cohomologies
corresponding to ${\cal D}$-modules \cite{milicic2}.

\subsubsection{Discrete Series for loop groups}
Let us now consider the loop groups associated with compact Lie groups. 
With the help of the Zuckermann's derived functor, we can construct the 
infinitesimal version of the ``discrete series'' for loop groups, \cite{diep9}
The algebraic realization of these representations are described in 
\cite{diep9} as a version of the Borel-Weil-Bott-Kostant theorem.
It is very interesting to develop a theory of ``tempered representations''
for loop groups.

\section{KK-theory Invariant $Index C^*(G)$}

 \subsection{About KK-Functors}
 
We now recall some essential points of the Kasparov's setting of the 
KK-theory. It is an analogy of the Brown-Douglas-Fillmore theory, 
but settled for the general case.

\subsubsection{Definitions}
 The main reference for this sub-subsection is \cite{jensenthomsen}.
 Let $A$, $B$, $E$ to be the C*-algebras, ${\cal K}$ the ideal of compact 
operators in some fixed separable Hilbert space. Let us consider the 
extensions of type 
 $$0 \longrightarrow B \otimes {\cal K} \longrightarrow E \longrightarrow A 
\longrightarrow 0.$$ Two extensions are said to be {\it equivalent } 
\index{extensions!equivalent - } iff there is 
some isomorphism $\psi :E \longrightarrow E'$  such that it induces the 
identity isomorphisms on the ideal $B \otimes {\cal K}$ and on the quotient 
$A$, i.e. the following diagram is commutative
$$
\begin{array}{ccccccccc}
 0 & \longrightarrow & B\otimes {\cal K} & \longrightarrow & E &
                                        \longrightarrow & A & \longrightarrow &    
                                                                 0\\
  &  & \Big\Vert &  & \Big\downarrow\vcenter{%
                       \rlap{$\psi$}}  &  & \Big\Vert & & \\
 0 & \longrightarrow &  B \otimes {\cal K} & \longrightarrow & E' 
& \longrightarrow & A & \longrightarrow & 0
\end{array}
$$
The extension is called {\it trivial } \index{extension!trivial}
if the the exact sequence can be lifted.
Also due to well-known result of R.C. Busby, we can identify each extension
with some *-homomorphism from $A$ to the algebra of exterior multiplicators
of $B \otimes {\cal K}$, $\tau : A \longrightarrow {\cal O}
(B \otimes {\cal K})$.
{\it The sum } \index{extensions!sum of - }
of two extensions $$\tau_i : A \longrightarrow {\cal O}
(B \otimes {\cal K})$$
can be therefore defined as the extension
$$\tau_1 \oplus \tau_2 : A \longrightarrow {\cal O}(B \otimes {\cal K})
\oplus {\cal O}(B \otimes {\cal K}) \hookrightarrow {\cal O}(B \otimes
{\cal K}) \otimes M_2 \cong {\cal O}(B \otimes {\cal K}),$$
where $M_2$ is the full algebra of $2\times 2$-matrices over the complexes
numbers. 
Two extensions $\tau_i, i=1,2$ are {\it  stably equivalent } 
\index{extensions!stably equivalent - }
if there exist two trivial extensions $\sigma_1$ and $\sigma_2$ such that 
the sums $\tau_i + \sigma_i$, $i=1,2$ are equivalent. 

Kasparov \cite{kasparov1} proved that:
\begin{enumerate}
\item[i)] 
when $A$ is a nuclear separable and
$B$ has at list an approximative unity, the set $Ext(A,B)$ of the 
stably equivalent classes of extensions is an Abelian group. 
\item[ii)]
$K^*(.)=Ext_*(A,.)$ is a K-homology theory and $K_*=Ext_*(.,B)$ 
is the algebraic K-theory of C*-algebras.
\item[iii)]
There is a natural realization of $KK^{*,*}$ as some K-bi-functor and its
direct relation with the $Ext_*$-groups see \cite{kasparov2}.
\end{enumerate}

J. Rosenberg and C. Schochet \cite{rosenbergschochet} proved the K\"unneth
formula for these groups, i.e. there is some homomorphism
$$Y: Ext_i(A,B) \longrightarrow \oplus_{j (mod \,2)}
\Hom_{\Bbb Z}(K_{i+j}(A),K_{i+j+1(B)}).$$
Let us see this in more detail in the next sub-subsection.

\subsubsection{Relation with K-groups of C*-algebras}

The most important for us is the relation of the theory with K-groups
of C*-algebras. Let $A$ be an algebra with unity. By definition,
$K_0(A)$ is the Grothendieck group of the semi-group of the stably
equivalent classes of projective $A$-modules of finite type.  When $A$
has no unity element, one considers the algebra $A^+$ with the
formally adjoint unity and defines the K-group as $$K_*(A) :=
\Ker \{K_*(A^+)\longrightarrow K_*({\Bbb C})={\Bbb Z}\}.$$ This
definition is compatible with the above defined K-groups also for
algebras with unity element. For $A=C(X)$, there is a natural
isomorphism between these K-groups with the corresponding topological
groups $K^*(X)$, see for example \cite{karoubi}. One defines the
higher groups $K_n(A)$ as $$K_n(A) := K_0(A \otimes C_0({\Bbb R}^n)),
\forall n \geq 0.$$ The Bott theorem says that $K_0(A) \cong K_2(A)$.
The Connes-Kasparov theorem says that for any connected and
simply-connected solvable Lie group $G$,

$$
K_0(C^*(G)) = \begin{cases} 
                  {\Bbb Z} & \text{if $\dim\, G$ is even,}\cr
		  0 & \text{if others, }
	       \end{cases}
$$

$$
K_1(C^*(G)) = \begin{cases}
                  0 & \text{if $\dim\, G$ is even,}\cr
		  {\Bbb Z} & \text{ if otheres.}
		  \end{cases}
$$

For each extension 
$$0 \longrightarrow J \longrightarrow E \longrightarrow A \longrightarrow 0,$$
there is a six-term exact sequence of K-groups
$$
\begin{array}{ccccccc}
    &    & K_0(E) & \longrightarrow & K_0(A) &      &    \\
    & \nearrow &   &    &    & \searrow &   \\
 K_0(J)  &   &    &    &    &    & K_1(J)\\
   & \nwarrow &    &    &    & \swarrow &    \\
    &   & K_1(A) & \longrightarrow & K_0(E) &   &   
 \end{array}
 $$
 
 Let us consider the case $J= B \otimes {\cal K}$, There is an isomorphism 
between $K_*(J)$ and $K_*(B)$. The group $K_*(A)$ consists of the formal
differences of equivalence classes of projectors in $A\otimes {\cal K}$.
One obtain therefore the well-known exact sequence
$$
\begin{array}{ccccccc}
 
   &      & K_0(E) & \longrightarrow & K_0(A) &   &   \\
   & \nearrow &   &   &     & \searrow\partial_0 &  \\
K_0(B) & &   &   &   &   & K_1(B)\\
   & \partial_1\nwarrow &   &   &   & \swarrow &   \\
   &   & K_1(A) & \longleftarrow & K_1(E) &   &  
\end{array}
$$
It is therefore clear that each element of $\Ext(A,B)$ induces a pair of 
homomorphisms $(\partial_0,\partial_1)$ of K-groups, and one has a homomorphism
$$\gamma : \Ext_i(A,B) \longrightarrow \oplus_{j \in {\Bbb Z}/(2)} 
\Hom_{\Bbb Z} (K_{i+j}(A),K_{i+j+1}(B)),$$
associating to each extension a pair of connecting homomorphisms
$(\partial_0, \partial_1)$.

J. Rosenberg and S. Schochet \cite{rosenbergschochet} have proved the
following exact sequence $$0 \longrightarrow \sum_{i\in {\Bbb Z}/(2)}
\Ext^1_{\Bbb Z}(K_{i+j}(A),K_{i+j+1}(B))
\longrightarrow \Ext_i(A,B) \longrightarrow $$
$$ \longrightarrow\sum_{j\in{\Bbb Z}/(2)} 
\Hom_{\Bbb Z} (K_{i+j}(A),K_{i+j+1}(B)) \longrightarrow 0.$$

\subsection{Construction and reduction of the K-Theory Invariant $Index\; C^*(G)$}

We review in this section a construction for obtaining the short exact sequence
of C*-algebras.

\subsubsection{Measurable foliations} 
In this section we propose a canonical method for constructing the
measurable \cite{connes1} foliations, consisting of the adjoint orbits
of fixed dimension, and therefore their C*-algebras. The last ones
are included in group C*-algebras or their quotients.

Let us denote by G a connected and simply connected Lie group,
${\frakt g} = \Lie (G)$ its Lie algebra, ${\frakt g}^* = \Hom_{\Bbb
R}({\frakt g}, {\Bbb R})$ the dual vector space,${\cal O} = {\cal O}%
(G)$ the space of all the co-adjoint orbits of G in ${\frakt g}^*$. This space is a disjoint union of subspaces of co-adjoint orbits of  fixed dimension, i.e. 
$${\cal O}%
= \amalg_{0 \leq 2n \leq \dim G}{\cal O}%
_{2n} ,$$
$${\cal O}%
_{2n} := \{ \Omega \in {\cal O}%
; \dim \Omega = 2n \} .$$
We define  
$$V_{2n} := \cup_{\dim \Omega = 2n}\Omega .$$
Then it is easy to see that $V_{2n}$ is the set of points of a fixed rank of the Poisson structure bilinear function
$$\{X,Y\}(F) = \langle F,[X,Y]\rangle  ,$$
 suppose it is a foliation, at least for $V_{2n}$, with $2n = max$ . 

First, we shall show that the foliation $V_{2n}$
can be obtained by the associated action of ${\Bbb R}^{2n}$ on $V_{2n}$ via 
2n times repeated action of ${\Bbb R}$ . 

Indeed, fixing any basis $X_1,X_2,\dots ,X_{2n}$ of the tangent space
${\frakt g} /{\frakt g}_F$ of $\Omega$ at the point $F \in \Omega$ ,
we can define an action ${\Bbb R}^{2n} \curvearrowright V_{2n}$ as
$$({\Bbb R} \curvearrowright ({\Bbb R} \curvearrowright (\dots {\Bbb
R} \curvearrowright V_{2n})))$$ by $$(t_1,t_2,\dots,t_{2n})
\longmapsto \exp(t_1X_1)\dots\exp(t_{2n}X_{2n}) F .$$ Thus we have the
Hamiltonian vector fields $$\xi_k := {d \over dt} |_{t=0}\exp(t_kX_k)F
, k = 1,2,\dots,2n$$ and the linear span $$F_{2n} =
\{\xi_1,\xi_2,\dots,\xi_{2n}\}$$ provides a tangent distribution.

\begin{thm} 
$(V_{2n},F_{2n})$ is a measurable foliation. 
\end{thm}
\begin{pf} See \cite{diep8}. \end{pf}

\begin{cor}
 The Connes C*-algebra $C^*(V_{2n},F_{2n}) , o \leq {2n} \leq
\dim G$ are well defined.
\end{cor}

\subsubsection{Reduction of $Index C^*(G)$ to $Index C^*(V_{2n},F_{2n})$ }

  Now we assume that the orbit method (see[Ki],[D4]-[D6]) gives us a  complete list of irreducible representations
of $G$ , $$\pi_{\Omega_F,\sigma} = Ind(G,\Omega_F,\sigma ,{\frakt p}), \sigma \in {\cal X}%
_G(F) ,$$
the finite set of Duflo's data.

    Suppose that
$${\cal O} = \cup_{i=1}^k{\cal O}%
_{2n_i}$$ 
is the decomposition of the orbit space on a stratification of orbits of dimensions $2n_i$, where $n_1 > n_2>  \dots > n_k > 0$

We include $C^*(V_{2n_1},F_{2n_1})$ into $C^*(G)$. It is well known that the Connes C*-algebra of foliation can be included in the algebra of pseudo-differential operators of degree 0 as an ideal. This algebra of pseudo-differential operators of degree 0

is

 included in C*(G).

We define $$J_1 = {\bigcap_{\Omega_F \in {\cal O}%
(G) \setminus {\cal O}%
_{2n_1}}} \Ker  \pi_{\Omega_F ,\sigma},$$
and $$ {A_1} = {C^*(G)/J_1}.$$
Then $$C^*(G)/C^*(V_{2n_1},F_{2n_1}) \cong A_1$$
and we have 

$$\vbox{\halign{ #\quad &#\quad &#\quad &#\hfill\cr
$ 0 \rightarrow $ & $J_1 \rightarrow$ & $C^*(G) \rightarrow$ & $A_1 \rightarrow 0$ \cr
          \hfill      & $\enskip\downarrow$ & $\quad \downarrow Id$ & $\quad\downarrow$   \cr
$0 \rightarrow$ & $C^*(V_{2n_1},F_{2n_1}) \rightarrow$ & $C^*(G) \rightarrow$ & $C^*(G)/C^*(V_{2n_1},F_{2n_1}) \rightarrow 0$ \cr}}$$

Hence $J_1 \simeq C^*(V_{2n_1},F_{2n_1})$ and we have 
$$O \rightarrow C^*(V_{2n_1},F_{2n_1}) \rightarrow C^*(G) \rightarrow A_1 \rightarrow 0 .$$   
Repeating the procedure in replacing $$C^*(G),C^*(V_{2n_1},F_{2n_1}),A_1,J_1$$ by $$A_1,C^*(V_{2n_1},F_{2n_1}),A_2,J_2,$$ we have 
$$0 \rightarrow C^*(V_{2n_2},F_{2n_2}) \rightarrow A_1 \rightarrow A_2 \rightarrow 0$$
etc .... 

So we obtain the following result.

\begin{thm}
 The group C*-algebra C*(G)  can be included in a finite sequence of extensions
$$ 0 \rightarrow C^*(V_{2n_1},F_{2n_1}) \rightarrow C^*(G) \rightarrow 
A_1 \rightarrow 0\leqno{(\gamma_1):}$$
$$ \quad 0 \rightarrow C^*(V_{2n_2},F_{2n_2}) \rightarrow A_1 
\rightarrow A_2 \rightarrow 0,\leqno{(\gamma_2):}$$
$$\ldots\dots\dots \dots \dots$$
$$ 0 \rightarrow C^*(V_{2n_k},F_{2n_k}) \rightarrow A_{k-1} \rightarrow A_k \rightarrow 0,\leqno{(\gamma_k):}$$
where $\widehat{A_k} \simeq Char(G)$
\end{thm}

\begin{cor}
$Index C^*(G)$ is reduced to the system $Index C^*(V_{2n_i},$ $ 
F_{2n_i}), i = 1,2 ,\dots , k$  by  the  invariants 
$$[\gamma_i] \in KK(A_i,C^*(V_{2n_i},F_{2n_i})) , i= 1,2,\dots ,k.$$
\end{cor}

\begin{rem}
Ideally, all these invariants $[\gamma_i]$ could be computed step-by-step from $[\gamma_k]$ to $[\gamma_1]$.
\end{rem}

\subsubsection{Reduction of  $Index C^*(V_{2n_i},F_{2n_i})$ to the 
computable extension indices valuated in topological KK-groups of pairs of spaces}

Let us consider $C^*(V_{2n_i},F_{2n_i})$ for a fixed i. We introduce the 
following assumptions which were considered by Kasparov 
in nilpotent cases [K2]:

$(A_1)$ \quad 
There exists $k \in {\Bbb Z} , 0 < k \leq 2n_i$ such that the foliation 
$$V_{gen} := V_{2n_i} \setminus (\Lie \Gamma)^\perp$$
has its C*- algebra 
$$C^*(V_{gen},F|_{V_{gen}}) \cong C({\cal O}%
_{gen}^\sim) \otimes {\cal K}%
(H) ,$$
where $$\Gamma := {\Bbb R}^k \hookrightarrow {\Bbb R}^{2n_i} 
\hookrightarrow G ,$$  
$$\Lie \Gamma = {\Bbb R}^k \hookrightarrow {\frakt g}/{\frakt g}_{F_i} ,
(\Lie \Gamma)^\perp \subset {\frakt g}^* \cap V_{2n_i}.$$

 \begin{exam} If $V_{gen}$ is a principal bundle, or the space ${\cal O}%
_{gen} = V_{gen}/G$ is a Hausdorff space, then $C^*(V_{gen},F|_{V_{gen}}) \simeq C({\cal O}%
_{gen}^\sim) \otimes {\cal K}%
(H)$ 
\end{exam}

 It is easy to see that  if the condition $(A_1)$ holds,
$C^*(V_{2n_i},F_{2n_i})$ is an extension of 
$C^*(V_{2n_i} \setminus V_{gen},F_{2n_i}|_.)$ by $C({\cal O}%
_{gen}^\sim) \otimes {\cal K}%
(H)$, where ${\cal O}_{gen}^\sim = \{ \pi_{\Omega_F,\sigma}; 
\Omega_F \in {\cal O}%
_{gen},\sigma \in {\cal X}%
_G(F)\}$, described by the multidimensional orbit method from the
previous section.
If $k = 2n_i , ({\Bbb R}^{2n_i})^\perp = \{ O \}$, 
$V_{2n_i} = V_{gen}$, we have $$C^*(V_{2n_i},F_{2n_i}) \simeq C({\cal O}%
_{2n_i}^\sim \otimes {\cal K}%
(H)) .$$
If $k = k_1 < 2n_i$ , then ${\Bbb R}^{2n_i-k_1}$ acts on $V_{2n_i} 
\setminus V_{gen}$ and we suppose that  a similar assumption $(A_2)$ holds

$(A_2)$ There exists $k_2,0 < k_2 \leq 2n_i-k_1$ such that 
$$(V_{2n_i} \setminus V_{gen})_{gen} := (V_{2n_i} \setminus V_{gen}) \setminus 
({\Bbb R}^{k_2})^\perp$$
has its C*-algebra 
$$C^*((V_{2n_i} \setminus V_{gen})_{gen},F_{2n_i}|.) \simeq C(({\cal O}%
_{2n_i} \setminus {\cal O}%
_{gen})_{gen})^\sim \otimes {\cal K}
(H) .$$
As above, if $k_2 = 2n_i - k_1$ , $C^*(V_{2n_i} \setminus V_{gen},F_{2n_i}|.) 
\simeq C(({\cal O}%
_{2n_i} \setminus {\cal O}%
_{gen})_{gen}^\sim) \otimes {\cal K}%
(H)$. In other case we repeat the procedure and go to assumption $(A_3)$, etc....

The procedure must be finished after a finite number of steps, say in m -th step,$$C^*((\dots(V_{2n_i} \setminus V_{gen}) \setminus (V_{2n_i} \setminus V_{gen})_{gen} \setminus \dots ,F_{2n_i}|_.) \simeq C((\dots ({\cal O}%
_{2n_i} \setminus {\cal O}%
_{gen}) \setminus \dots )) \otimes {\cal K}%
(H) .$$
Thus we have the following result.

\begin{thm}
If all the arising assumptions $(A_1),(A_2),\dots$ hold, the C*-algebra $C^*(V_{2n_i},F_{2n_i})$ can be included in a finite sequence of extensions
$$ 0 \rightarrow C({\cal O}%
_{gen}^\sim) \otimes {\cal K}%
(H) \rightarrow C^*(V_{2n_i},F_{2n_i}) \rightarrow C^*(V_{2n_i} \setminus V_{gen},F_{2n_i}|_.) \rightarrow 0$$
$$ 0 \rightarrow C(({\cal O}%
_{2n_i} \setminus {\cal O}%
{gen})_{gen}^\sim) \otimes {\cal K}%
(H) \rightarrow C^*(V_{2n_i} \setminus V_{gen},F_{2n_i}) \rightarrow C^*(\dots)
\rightarrow 0$$
$$\ldots \ldots \ldots \ldots \ldots$$
$$ 0 \rightarrow C((\dots({\cal O}%
_{2n_i} \setminus {\cal O}%
_{gen}) \setminus ({\cal O}%
_{2n_i} \setminus {\cal O}%
_{gen}))_{gen}\dots{^\sim}) \otimes {\cal K}%
(H) \rightarrow $$ $$ \rightarrow C^*(\dots) \rightarrow C^*(\dots) \otimes 
{\cal K}%
(H) \rightarrow 0.$$
\end{thm}

\subsubsection{General remarks concerning computation of Index C*(G)}

We see that the general computation procedure of Index C*(G) is reduced to the case of short exact sequences of type 
$$ 0 \rightarrow C(Y) \otimes {\cal K}%
(H) \rightarrow {\cal E}%
\rightarrow C(X) \otimes {\cal K}%
(H) \rightarrow 0, \leqno{(\gamma)}$$
and the index is
$$ [\gamma] = Index {\cal E}%
\in KK(X,Y) .$$
The group $KK_i(X,Y)$ can be mapped onto 
$$\oplus_{j \in {\Bbb Z}/(2)}\Hom_{\Bbb Z}(K^{i+j}(X),K^{i+j+1}(Y))$$
with kernel 
$$\oplus_{j \in {\Bbb Z}/(2)}\Ext_{\Bbb Z}^1(K^{i+j}(X),K^{i+j+1}(Y))$$
by the well known cap-product, see [K2]. So $[\gamma] = (\delta_0,\delta_1)$
$$\delta_0 \in \Hom_{\Bbb Z}(K^0(X),K^1(Y)) = \Ext_0(X) \wedge K^1(Y)$$
$$\delta_1 \in \Hom_{\Bbb Z}(K^1(X),K^0(Y)) = \Ext_1(X) \wedge K^0(Y).$$
Suppose $e_1,e_2,\dots,e_n \in \pi^1(X)$ to be generators and $\phi_1,\phi_2,
\dots,\phi_n \in {\cal E}$%
\quad the corresponding Fredholm operators, $T_1,T_2,\dots,T_n$ the 
Fredholm operators, representing the generators of $K^1(Y) = Index 
[Y,\Fred ]$ . We have therefore $$[\delta_0] = \sum_jc_{ij}\enskip Index T_j ,
$$ where
$$\delta_0 = (c_{ij}) \in \Mat_{\rank K^0(X) \times rank K^1(Y)}({\Bbb Z}) .$$
In the same way $\delta_1$ can be computed.

\section{Deformation Quantization and Cyclic Theories}

Let us finish this  survey with some indication about some relations
of the problem with some new developments. Recall that the group
algebra of finite or compact groups are in fact some Hopf bi-algebras.
One deforms this Hopf bialgebra structure to obtain the corresponding quantum
groups. Our problem is therefore closely related with the interesting 
problem to describe these quantum groups. 
One of the method is deformation quantization which is closely related
with orbit method. The others which are closely related with KK-theory are
the periodic cyclic (co-)homologies. We finish this survey by indication
the subjects, the author is working with.

\subsection{Star-Products and Star-Representations}
See \cite{gutt} and the references there.

\subsection{Periodic Cyclic Homology}
See \cite{cuntz}, \cite{cuntzquillen1} - \cite{cuntzquillen3}.

\subsection{Chern Characters}
See for example \cite{cuntz}, \cite{puschnigg}.

\section{Bibliographical Remarks}

The material exposed in this introduction was the subject of the author
for the talks at Seminar on Representation Theory and Gelfand's Seminar on
Functional Analysis in Moscow University (1975), Centre International des
Rencontres Math\'ematiques in Luminy-Marseille (1983), Seminar on Group
Representation \`a l'Universit\'e de Paris VI et Paris VII (1983),
l'Universit\'e de Lyon (1983), International Banach Center (1988), Seminar
on C*-algebras at the Iniversit\"at Heidelberg, the Humboldt Universit\"at
zu Berlin, Seminar `` Sophus Lie" at the Technische Hochschule Darmstadt
(1993) and many conference talks.  The idea to create this book appeared
during preparation of lectures at Seminar ``Sophus Lie".

\part{Elementary Theory: An Overview \\ Based on Examples}
There are already fifteen years ago from the first moment of attacking
the C*-algebra structure problem by using the K-functors. A general
method is not yet constructed, but there is a lot of accommulated
results, confirming the usefulness of the topological K-theory
invariant $Index~C^*(G)$. We think therefore reasonable to propose an
overview based on experimented examples of theses research works. In
this survey we intend to give a concrete reflection of the current
research . It is an introduction to a writing book on the subject.

 {\bf 1. Group C*-algebra Structure } Let $G$ be a locally compact group.
There exists a naturally normalized left invariant Haar measure $dg$ .The
space $L^2(G,dg)$ of the square-integrable module functions plays an
important role, say in Harmonic Analysis studies , first of all the
spectral decomposition of the regular representation of $G$ in
$L^2(G,dg)$ into a direct integral ( or a sum ) of irreducible unitary
$G$-module. The space $L^1(G,dg)$ of functions with integrable module
plays a crucial role . Following the non commutative Fourier-Gel'fand
transformation 
$$\begin{array}{rl} \varphi \in L^1(G,dg) \mapsto &\hat{\varphi},\\
  \hat{\varphi}(\pi) = \pi(\varphi) &:= \int_G\;\pi(x) \varphi(x)~dx 
\end{array}$$
We have a one-to-one correspondence between the (irreducible) unitary
$G$-modules and the (irreducible) non-degenerate *-modules over $L^1(G)$. 
So the general theorems of the spectral theory for $G$-modules can be
translated and proved in the corresponding theory for $L^1(G)$ , which is
more analytical then for which are applicable the strong results of
Functional Analysis , say the Hahn-Banach theorem, the Banach principles
of linear functional analysis. However, as involutive algebra, $L^1(G)$
has its non-regular norm, i.e. in general $$\parallel
a^*a\parallel_{L^1(G)} \ne \parallel a\parallel^2_{L^1(G)} \; ,$$ it is
more useful to consider the corresponding regular norm $\parallel .
\parallel_{C^*(G)}$ , $$\parallel \varphi \parallel_{C^*(G)} :=
\sup_{\pi\in\hat{G}} \parallel \pi(\varphi)\parallel$$ and take its
completion $C^*(G)$. Ideally, the spectral theory for unitary $G$-modules
is equivalent to the same one for the $C^*(G)$-modules. The last theory is
closely related with subjects of functional analysis and its applications
in the physical field theories and statistical mechanics. 

So, what is the structure of $C^*(G)$ for a given $G$?

In general, the problem rests open up-to-date! This review outlines only
that for the concrete examples the problem requires the tools of various
nature from topology and analysis, namely the K-functors. 

It is useful to add to $C^*(G)$ the formal jointed unit element if there
is no such one. 

{\bf 2. The analytic method } The first nontrivial example
$C^*(\SL_2({\Bbb C}))$ was done by J.M.G. Fell in 1961, but until the
moment the only groups, the structure of whose $C^*$-algebra was
explicitly described are the Abelian or compact groups and a few
semi-simple Lie groups :  $\SL_2({\Bbb R})$, $\widetilde{\SL_2({\Bbb R})}$
i.e. the universal covering group, $Spin(4,1)$ and recently, a family
$G(p,q,\alpha)$ of two-step solvable Lie groups (see
\cite{delaroche},\cite{fell},\cite{milicic1},\cite{kraljevicmilicic},
\cite{boyermartin},
and \cite{wang1}).A fair amount is known about the C*-algebras of nilpotent
Lie groups (see \cite{perdrizet}), including the Heisenberg groups. The 
C*-algebra
of the Euclidean motion group were studied by Evans \cite{evans}. P. Green
also proposed another analytic method for studying the C*-algebras of
several solvable Lie groups, see \cite{green}.  So it is very interesting
to characterize the group C*-algebras by topological invariants by another
nature tools, say by K-functors. Such an idea was suggested by the author
in \cite{diep1},\cite{diep2}.

{\bf 3. K-theory invariant Index~C*(G)} The very useful K-functor for our
approach is the operator KK-functor of G. G. Kasparov \cite{kasparov1},
generalizing the BDF K-functor \cite{bdf1}, which characterizes the
isomorphic classes of short exact sequences of C*-algebras . 

We are trying to decompose our C*-algebras into some towers of ideals and
step-by-step define the associated extensions by KK-functors or their
generalizations. The resulting invariants form just our index . This idea
was proposed in \cite{diep1} and develop ed in \cite{diep2} for a large 
class of type I C*-algebras. Hence, there are two general problems:
\begin{itemize} \item Find out the C*-algebras which can be characterized
by the well-known K-functors, say by the operator K-theory functors. \item
Generalize the theory of K-functors in such a way that they are applicable
for a large class of C*-algebras.  \end{itemize}

{\bf 4. Construction and reduction of Index~C*(G)} Concerning the first
problem, we propose \cite{diep8} a general construction and some reduction
procedure of the K-theory invariant Index~C*(G) of group C*-algebras.
Using the orbit method \cite{kirillov},\cite{diep4}-\cite{diep7}, we
reduces $Index~C^*(G)$ to a family of Connes' foliation C*-algebras
indices $Index~C^*(V_{2n_i},{\cal F}_{2n_i})$ , see 
\cite{connes1}-\cite{connes2},
by a family of KK-theory invariants. Using some generalization of the
Kasparov type condition (treated by G.G. Kasparov in the nilpotent Lie
group case \cite{kasparov2}), we reduces every $Index~C^*(V_{2n_i},{\cal
F}_{2n_i})$ to a family of KK-theory invariants of the same type valuated
in KK(X,Y) type groups. The last ones are in some sense computable by
using the cup-cap product realizing the Fredholm operator indices. 

Following this procedure we now describe the obtained experimental results
concerning the structure of the C*-algebras of concrete groups. 

{\bf 5. Case-by-case examples} We divide the examples into three classes
following the complexity of computing indices.  {\it a. Absorbing
extensions} \index{extension!absorbing} The first promising example was the group 
$\Aff{\Bbb R}$ 
of the affine transformations of the real straight line. Its group C*-algebra can
be included in the short exact sequence $$0 \rightarrow {\cal K}(H)
\rightarrow C^*(\Aff{\Bbb R}) \rightarrow C({\Bbb S}^1 \vee {\Bbb S}^1)
\rightarrow 0\;,$$ where ${\cal K}(H)$ denotes the ideal of compact
operators in a separable Hilbert space. Because the extension is
absorbing, the structure of $C^*(\Aff{\Bbb R})$ is just defined by the
element $$Index~C^*(\Aff{\Bbb R}) = (1,1) \in {\Bbb Z} \oplus {\Bbb Z}
\cong {\cal E}xt({\Bbb S}^1 \vee {\Bbb S}^1) \cong KK({\Bbb S}^1 \vee
{\Bbb S}^1, pt)\;.$$ Here we need only the BDF K-functor , which is the
source and the inspiring particular case of the KK-functors (see
\cite{diep1}). 

For the connected and simply connected group $\Aff_0{\Bbb R}$ of affine
transformations of the straight line the analogous results hold. Its
C*-algebra can be included in the sort exact sequence $$ 0 \rightarrow
{\cal K}(H) \oplus {\cal K}(H) \rightarrow C^*(\Aff_0{\Bbb R}) \rightarrow
C({\Bbb S}^1) \rightarrow 0 \;.$$ This absorbing extension, and hence the
isomorphic class of $C^*(\Aff_0{\Bbb R})$ can be characterized by the
topological invariant $$Index~C^*(\Aff_0{\Bbb R}) = (1,1) \in KK({\Bbb S}^1
,\{pt\} \cup \{pt\}) \;.$$

The C*-algebra of the group $\Aff{\Bbb C}$ of the affine transformations of
the complex straight line is included in the short exact sequence $$ 0
\rightarrow {\cal K}(H) \rightarrow C^*(\Aff{\Bbb C}) \rightarrow C(X)
\rightarrow \;,$$ where $X$ is the one-point compactification of the so
called "Hawaiian necklace"  $$\{ z \in {\Bbb C} ; \vert z - 2^{-n}\vert =
2^{-n} , n = 1,2,... \}\;.$$ Hence , the isomorphism class of
$C^*(\Aff{\Bbb C})$ is characterized by the topological invariant
$$Index~C^*(\Aff{\Bbb C}) = (-1,-1,...) \in KK(X,pt) \cong {\Bbb Z} \oplus
{\Bbb Z} \oplus \dots ,$$ see \cite{rosenberg1} . 

Our method is applied also to the group $\Aff~K$ of affine transformations
of any non-discrete totally disconnected locally compact field K. Its
C*-algebra can be included in the short exact sequence $$0\rightarrow
{\cal K}(H) \rightarrow C^*(\Aff~K) \rightarrow C(X) \rightarrow 0 \;,$$
now $X$ is the one-point compactification of ${\Bbb S}^1 \times \hat{H}$,
where $H$ is the multiplicative group of the elements of $K$ with absolute
values 1, $\hat{H}$ is its dual isomorphic to a countably infinite and
discrete set. The C*-algebra $C^*(\Aff K)$ is characterized by the
topological invariant $$Index~C^*(\Aff K) = (\dots,1,1,\dots) \in KK(X,pt)
\cong {\cal E}xt(X) \cong {\Bbb Z}^{\hat{H}}\;,$$ see \cite{rosenberg1}. 

In all these examples the index $Index~C^*(G)$ take values in the BDF
K-groups ${\cal E}xt(X)$ , which are isomorphic to $\Hom_{\Bbb
Z}(\pi^1(X),{\Bbb Z})$, given by the Fredholm index map, where $$\pi^1(X)
= [X,{\Bbb S}^1]$$ is the cohomotopy group. The group $\Hom_{\Bbb
Z}(\pi^1(X),{\Bbb Z})$ is product of countably many copies of ${\Bbb Z}$,
one for each generator of $\pi^1(X)$. The index is then given by a
sequence of integers , namely the Fredholm indices of the image under the
infinite-dimensional representations of a sequence of elements from
$C^*(G)$ mapping into the generators of $\pi^1(X)$. These sequences of
elements from the group C*-algebra $C^*(G)$ and their Fredholm indices are
defined and calculated firstly by the author for $\Aff{\Bbb R}$ in
\cite{diep1}, then by J. Rosenberg for $\Aff{\Bbb C}$ and $\Aff K$ in 
\cite{rosenberg1}. 

The C*-algebra of the universal covering $\widetilde{\Aff{\Bbb C}} = {\Bbb
C} \ltimes {\Bbb C}$ of the group $\Aff{\Bbb C}= ({\Bbb C}\setminus 0)
\ltimes {\Bbb C}$ can be included in the short exact sequence $$0
\rightarrow C({\Bbb S}^1) \otimes {\cal K}(H)
 \rightarrow C^*(\widetilde{\Aff{\Bbb C}}) \rightarrow C({\Bbb S}^2)
\rightarrow 0$$ and the structure of $C^*(\widetilde{\Aff{\Bbb C}})$ is
uniquely defined by the index $$Index C^*(\widetilde{\Aff{\Bbb C}}) = 1 \in
KK({\Bbb S}^2,{\Bbb S}^1) = {\cal E}xt(C({\Bbb S}^2),C({\Bbb S}^1)) =
{\Bbb Z}\;.$$

{\it b. Non-absorbing extensions} 

The first example which requires essentially the KK-theory , but not
enough the ${\cal E}xt$-functor BDF K-theory arisen as a class of two-step
solvable Lie groups ${\Bbb R} \ltimes_{\alpha>0}{\Bbb R}^m$, considered by
J. Rosenberg \cite{rosenberg1} and then $G(p,q) = {\Bbb R} \ltimes {\Bbb 
R}^{p+q}$
considered by X. Wang \cite{wang1}, including the Heisenberg group ${\Bbb
H}_{2n+1}$ considered before by G. G. Kasparov 
\cite{kasparov1},\cite{kasparov2}. But in
these examples we obtain the non-absorbing exact sequences and then the
topological invariant $Index C^*(G)$ does not define the isomorphic classes
of C*-algebras. Let us describe the results in more detailed form. 

The C*-algebras of the , say {\it elliptic semi-simple } product 
\index{product!elliptic semi-simple -} ${\Bbb R}
\ltimes {\Bbb R}^m$ ($m$ is any positive integer) , where the action of
${\Bbb R}$ on ${\Bbb R}^m$ have roots with the real part of the same sign
, are isomorphic each to another and can be included in the short exact
sequences of type
 $$ 0 \rightarrow C({\Bbb S}^{m-1},{\cal K}(H)) \rightarrow C^*({\Bbb R}
\ltimes {\Bbb R}^m) \rightarrow C({\Bbb S}^1) \rightarrow 0 \;.$$ Hence we
have , but not describe by, the topological invariant $$Index~C^*({\Bbb R}
\ltimes {\Bbb R}^m) \in KK({\Bbb
 S}^{m-1},{\Bbb S}^1) $$ of the group C*-algebras.

The C*-algebra of Heisenberg group ${\Bbb H}_{2n+1}$ can be included in 
the short , but not absorbing , exact sequence $$0 \rightarrow C({\Bbb S}^1 
\vee {\Bbb S}^1,{\cal K}(H)) \rightarrow C^*({\Bbb H}_{2n+1} 
\rightarrow C({\Bbb S}^{2n}) \rightarrow 0 \;$$

Hence , the index $$Index~C^*({\Bbb H}_{2n+1}) = (1,(-1)^n) \in KK({\Bbb
S}^{2n},{\Bbb S}^1 \vee {\Bbb S}^1) \cong {\Bbb Z} \oplus {\Bbb Z} \;.$$
This result has been generalized to the general case of connected and
simply connected Lie groups. The group C*-algebras in these cases can be
included in the exact sequences of type $$0 \rightarrow C(X) \otimes {\cal
K}(H) \rightarrow C^*(G) \rightarrow C^*(G/R) \rightarrow 0$$ if the union
of the co-adjoint G-orbits of maximal dimension is the complement to the
annihilator of $Lie\,\Gamma$ in ${\frakt g}^*$, $\gamma = R$, where the
set of all the co-adjoint orbits of maximal dimension is denoted by $X$,
see \cite{kasparov2}. So, the C*-algebra $C^*(G)$ admits the topological
invariant $Index~C^*(G)= (1,(-1)^{n/2})$, where $n$ is the half-dimension
of the co-adjoint orbits of maximal dimension. 

{\it c. Non-absorbing extension, associated with towers of C*-ideals}
\index{extensions!non-absorbing -}
There are only a finite number of non-isomorphic C*-algebras between the
group C*-algebras of the 3-dimensional real solvable Lie groups. These
C*-algebras are easily characterized all but the subclass
$G_{3,2}(-\alpha)$, $\alpha >0$, the C*-algebras of which are isomorphic
one-to-other and are included into the short exact (but non absorbing)
sequences $$0 \rightarrow C({\Bbb S}^1 \vee {\Bbb S}^1 \vee {\Bbb S}^1
\vee {\Bbb S}^1) \otimes {\cal K}(H) \rightarrow C^*(G_{3,2}(-\alpha))
\rightarrow A^1 \rightarrow 0\leqno{(\gamma_1)}\;,$$ $$0 \rightarrow {\Bbb
C}^4 \otimes {\cal K}(H) \rightarrow A^1 \rightarrow C({\Bbb S}^1)
\rightarrow 0.$$ Hence, $Index~C^*(G_{3,2}(-\alpha)) = ([\gamma_1],
[\gamma_2])$ ; $$[\gamma_1] = \pmatrix 1 & 0 & 0 & 1 \\ 1 & 1 & 0 & 0 \\ 0
& 1 & 1 & 0 \endpmatrix \in KK(A^1,{\Bbb S}^1 \vee {\Bbb S}^1 \vee {\Bbb
S}^1 \vee {\Bbb S}^1) \cong \Hom_{\Bbb Z}({\Bbb Z}^3,{\Bbb Z}^4)$$ and
$$[\gamma_2] = (1,1,-1,-1) \in KK({\Bbb S}^1, 4\,pt) \cong {\Bbb Z}^4
\;.$$ It is very interesting to consider the following class $MD$ (resp.,
$\overline{MD}$) of connected and simply connected solvable Lie groups,
all the co-adjoint orbit of which have maximal dimension or zero (resp.
equal to the dimension of the groups or zero). It is easy to see 
\cite{sonviet}
that the only non commutative groups in the class $\overline{MD}$ are the
groups of affine transformations of the real or complex straight line. 

All the C*-algebras of the groups in class $\overline{MD}$ are
characterized up to isomorphism by the topological invariant
$Index~C^*(G)$ : $$Index~C^*(\Aff_0{\Bbb R}) = (1,1) \in KK({\Bbb S}^1,
2\,pt) \cong {\Bbb Z} \oplus {\Bbb Z} \;,$$ $$Index
~C^*(\widetilde{\Aff{\Bbb C}}) = 1 \in KK({\Bbb S}^2, {\Bbb S}^1) \cong
{\Bbb Z} \;.$$

The $MD$-groups are not yet classified, but the subclass $\MD_4$ of
four-dimensional solvable $MD$-groups is completely listed, including the
real diamond group ${\Bbb R} \ltimes {\Bbb H}_3$ (see
\cite{viet1},\cite{viet2}); There are 13 concrete connected and simply 
connected
$\MD_4$-groups or series of such ones with exact description of commutator
relations . All the C*-algebras of these $\MD_4$-groups can be described by
the direct analytic method, but only three cases of groups $G$, the Lie
algebra of which is $$Lie\;G \cong \langle T,X,Y,Z\rangle , [X,Y] = Z$$ and
\begin{itemize}
\item $adT = \pmatrix
 \cos\varphi & \sin\varphi & 0 \\ -\sin\varphi & \cos\varphi & 0 \\ 0 & 0
& \lambda \endpmatrix$, the Euclidean motion group $G_{\varphi,\lambda}$
\item $adT = \pmatrix
 0 & 1 & 0 \\ -1 & 0 & 0 \\ 0 & 0 & 0 
 \endpmatrix$, the Harmonic Oscillator ${\Bbb R} \ltimes_J{\Bbb H}_3$ \item
$adT = \pmatrix -1 & 0 & 0 \\ 0 & 1 & 0 \\ 0 & 0 & 0 \endpmatrix$, the
Real diamond group ${\Bbb R} \ltimes {\Bbb H}_3$ 
\end{itemize}

The C*-algebra $C^*(G_{\varphi,\lambda})$ can be included in the short
exact sequences $$0 \rightarrow C^*(V_{\varphi,\lambda},{\cal F})
\rightarrow C^*(G_{\varphi,\lambda}) \rightarrow C({\Bbb S}^1) \rightarrow
0\leqno{(\gamma_1)}$$ $$0 \rightarrow C({\Bbb S}^2 \vee {\Bbb S}^2)
\otimes {\cal K}(H) \rightarrow C^*(V_{\varphi,\lambda},{\cal F})
\rightarrow C({\Bbb S}^1) \otimes {\cal K}(H) \rightarrow
0\leqno{(\gamma_2)}$$ and is characterized by the index
$Index~C^*(G_{\varphi,\lambda}) \cong ([\gamma_1], [\gamma_2])$,
$$[\gamma_1] \in KK(C({\Bbb S}^1),C^*(V_{\varphi,\lambda},{\cal F})) \;,$$
$$[\gamma_2] = (1,1) \in KK({\Bbb S}^1,{\Bbb S}^2 \vee {\Bbb S}^2) \cong
{\Bbb Z} \oplus {\Bbb Z}\;.$$

The C*-algebra $C^*({\Bbb R} \ltimes_J {\Bbb H}_3)$ can be included in the
short exact sequences $$ 0 \rightarrow C^*(V_{{\Bbb R} \ltimes_J {\Bbb
H}_3}, {\cal F}) \rightarrow C^*({\Bbb R} \ltimes_J {\Bbb H}_3)
\rightarrow C({\Bbb S}^1) \otimes {\cal K}(H) \rightarrow
0\leqno{(\gamma_1)}$$ $$0 \rightarrow C(({\Bbb R}^\times \times {\Bbb
R})_{cpt} \otimes {\cal K}(H) \rightarrow C^*({\Bbb R} \ltimes_J {\Bbb
H}_3) \rightarrow C({\Bbb S}^1) \otimes {\cal K}(H) \rightarrow
0\leqno{(\gamma_2)}$$ and has the same index $Index~C^*({\Bbb R} \ltimes_J
{\Bbb H}_3) = ([\gamma_1],[\gamma_2])$, $$[\gamma_1] \in KK(C({\Bbb
S}^1),C^*(V_{{\Bbb R} \ltimes_J {\Bbb H}_3},{\cal F})$$ $$[\gamma_2] =
(1,1) \in KK({\Bbb S}^1,{\Bbb S}^2) \cong {\Bbb Z} \oplus {\Bbb Z} \;.$$

The C*-algebra of the real diamond group ${\Bbb R} \ltimes {\Bbb H}_3$ can
be included in three exact sequences $$0 \rightarrow C({\Bbb S}^2 \vee
{\Bbb S}^2) \otimes {\cal K}(H) \rightarrow C^*({\Bbb R} \ltimes {\Bbb
H}_3) \rightarrow A^1 \rightarrow 0\leqno{(\gamma_1)}$$ $$0 \rightarrow
C({\Bbb S}^2\vee {\Bbb S}^2 \vee {\Bbb S}^2 \vee {\Bbb S}^2) \otimes {\cal
K}(H) \rightarrow A^1 \rightarrow A^2 \rightarrow 0\leqno{(\gamma_2)}$$
$$0 \rightarrow {\Bbb C}^4 \otimes {\cal K} \rightarrow A^2 \rightarrow C(
{\Bbb S}^1) \rightarrow 0\leqno{(\gamma_3)}$$ and $$Index~C^*({\Bbb R}
\ltimes {\Bbb H}_3) = ([\gamma_1], [\gamma_2],[\gamma_3])\;,$$ where
$$[\gamma_1] = (1,1) \in KK(A^1,C({\Bbb S}^2 \vee {\Bbb S}^2) \cong
\Hom_{\Bbb Z}({\Bbb Z},{\Bbb Z}^2) \cong {\Bbb Z}^2,$$ $$[\gamma_2] =
\pmatrix -1 & 0 & 0 & 1 \\ 1 & -1 & 0 & 0 \\ 0 & 1 & -1 & 0 \\ 0 & 0 & 1 &
-1 \endpmatrix \in KK(A^2, C({\Bbb S}^2 \vee{\Bbb S}^2 \vee {\Bbb S}^2
\vee {\Bbb S}^2)) \cong \Hom_{\Bbb Z}({\Bbb Z}^4,{\Bbb Z}^4),$$
$$[\gamma_3] = (1,1,-1,-1) \in KK(C({\Bbb S}^1),{\Bbb C}^4) \cong
\Hom_{\Bbb Z}({\Bbb Z},{\Bbb Z}^4).$$ So the description of $\MD_4$-group
C*-algebras is achieved. One can hope to describe the C*-algebras of the
whole class $MD$ by the same method. The question rests open update. 

Finally, the C*-algebras of the hyperbolic semi-direct product $${\Bbb R}
\ltimes {\Bbb R}^{p+q} = G(p,q,\alpha)$$ with $p$ negative roots
$-\alpha_1,\dots,-\alpha_p$ and $q$ positive roots
$\alpha_{p+1},\dots,\alpha_{p+q} $ can be included in two short exact
sequences $$0 \rightarrow C^*(U,F)^\sim \rightarrow C^*(G(p,q,\alpha))
\rightarrow C({\Bbb S}^1) \rightarrow 0\;,$$ where $U$ is the one-point
compactification of ${\Bbb R}^{p+q} \setminus (0)$ and $$0 \rightarrow
C^*(U_1,F)^\sim \rightarrow C^*(U,F)\rightarrow \rightarrow C({\Bbb S}^p
\vee {\Bbb S}^q) \rtimes {\Bbb R} \rightarrow 0,$$ where $U_1 = {\Bbb
R}^{p+q} \setminus ({\Bbb R}^p \vee {\Bbb R}^q)$ . So the C*-algebra
$C^*(G(p,q,\alpha))$ admits the topological invariant
$Index~C^*(G(p,q,\alpha))$ valuated in the Kasparov groups $KK(C({\Bbb
S}^1),$ $C^*(U,F))$ and $KK(C({\Bbb S}^p \vee {\Bbb S}^q) \rtimes {\Bbb
R})$. 

In all these examples, the invariant $Index~C^*(G)$ is a sequence of type
$(\delta_0,\delta_1)$-homomorphisms in the six-term exact sequences of
$K$-groups. They are in general expressed by using intersection
cup-cap-products . Only for the examples in subsection a. one can use
analytic method of calculating the Fredholm index. So it is easy to see
that we must develop the theory of K-functors admitting some intersection
products. 

In the next four chapters we shall expose these cases in detail.
We hope that after reading this part the beginner can already work in these
problems more or less productively.

\chapter{Classification of $\overline{\MD}$-Groups}
\section{Definitions} First of all we recall the notion of $K$-action. Let
us denote by $G$ a connected and simply connected Lie group, ${\frakt g}=
T_eG$ its Lie algebra as the tangent space at the neutral element $e$. It
is easy to see that to each element $g\in G$ one can associate a map
$$A(g) : G \to G$$ by the conjugacy, in fixing the identity element $e\in
G$. Therefore, the corresponding tangent map $A(g)_* : {\frakt g} \to
{\frakt g}$ $$X\in {\frakt g} \mapsto {d\over dt}|_{t=0} g\exp(tX)g^{-1}
\in {\frakt g}.$$ It is easy to see that this really defines an action,
denoted as usually by $\Ad$ of group $G$ in its Lie algebra ${\frakt g}$.
One defines therefore an so called co-adjoint action of group $G$ in the dual 
vector space
${\frakt g}^*$ by the formula $$\langle K(g)F, X\rangle := \langle F,
\Ad(g^{-1})X \rangle ,$$
for all $F\in {\frakt g}^*$, $X\in {\frakt g}$ and $g\in G$.
It is easy to check that this defines a real action of $G$ on $\frakt g^*$.

\begin{defn}
The orbits of this action are called the {\it co-adjoint orbits} 
\index{orbit!co-adjoint -} or {\it $K$-orbit}. \end{defn}

As an easy consequence, one deduce that the dual space ${\frakt g}^*$ is 
decomposed into a disconnected sum of the K-orbit.

\begin{defn} We say that a real Lie algebra ${\frakt g}$ is in the class 
$\overline{\MD}$ if every K-orbit is of dimension, equal 0 or 
$\dim{\frakt g}$. \end{defn}

This means that the structure of the orbit space of such a Lie algebra 
must be rather simple: There is only two strata of orbits and the union 
of the maximal dimension is dense in the whole space 
${\frakt g}^*$. The other evident consequence is the 
fact that the dimension $\dim{\frakt g}= \dim G$ must be an even 
number.

Recall that in each K-orbit there is a natural differential form, 
associated with the bilinear 
form $$B_F(X,Y) := \langle F, [X,Y]\rangle , \forall X,Y \in {\frakt g}.$$
It is not hard to verify the following assertion 

\begin{clm}
The kernel of this bilinear form is just the Lie algebra ${\frakt g}_F := 
\Lie G_F$ of the stabilizer $G_F$ of the point $F\in {\frakt g}^*$ under 
the co-adjoint action.
 \end{clm}
\begin{pf}
For a connected Lie group, every element can be obtained as some 
product of elements from 
the image of the exponential map $$\exp : {\frakt g} \to G.$$ We can 
therefore restrict to the case of element of type $\exp(X), X\in {\frakt g}$
. For the elements of this kind it is enough to remember a formula from 
the Lie theory
$$\Ad(\exp X) = \exp(\ad X).$$ 
\end{pf}

This means that the form $B_F(.,.)$ is invariant under the action of 
the stabilizer group $G_F$ of the $K$-orbit passing through $F$. One can 
therefore translate this form to other points in order to have a 
differential form on the $K$-orbit $\Omega_F = G/G_F$, passing through $F$. 

\begin{defn}
The corresponding symplectic form is called Kirillov form on K-orbits.
\end{defn}

\section{$\overline{\MD}$-Criteria}

Let us denote by $\frakt g^1 := [\frakt g,\frakt g]$ the commutator of the 
Lie algebra.

\begin{prop}
If $F\in {\frakt g}^* \setminus ({\frakt g}^1)^*$, i.e. if $F$ is a 
functional on 
${\frakt g}$ which is non-vanishing on ${\frakt g}^1$, then the K-orbit 
$\Omega_F$
passing through $F$, is of maximal dimension, $\dim \Omega_F = \dim G$. 
Moreover ${\frakt g}^1$ is commutative.\end{prop}

\begin{pf} We prove this proposition by contradiction argument.
Assume that $F$ is some functional on ${\frakt g}$ which is non-vanishing 
on the commutator ${\frakt g}^1 = [{\frakt g},{\frakt g}]$ but $\dim 
\Omega_F \ne \dim{\frakt g}$, i.e. $\dim\Omega_F = 0.$ This means 
that $$\dim G_F = \dim G - \dim \Omega_F = \dim {\frakt g}.$$
Hence, ${\frakt g}_F = {\frakt g}$, i.e. $\ker B_F = {\frakt g}$. This 
contradicts the assumptions. Hence $\dim\Omega_F = \dim{\frakt g}= \dim G$.

Now we prove that ${\frakt g}^1$ is commutative. Denote by ${\frakt g}^2$ 
the second derived ideal, ${\frakt g}^2 = [{\frakt g}^1,{\frakt g}^1]$. 
Because ${\frakt g}$ is solvable, $\dim{\frakt g}^2 < \dim{\frakt g}^1$. 
Hence there exists a nonzero functional $F \in ({\frakt g}^1)^*$, 
vanishing on ${\frakt g}^2$. This means that ${\frakt g}^2 \subset 
{\frakt g}_F = 0$. What mean that ${\frakt g}^1$ is commutative.
\end{pf}

\begin{prop}{\bf ($\overline{\MD}$-criterion)} Lie algebra ${\frakt g}$ is
of class $\overline{\MD}$ if and only if $$\ad_X({\frakt g}) = [X,{\frakt
g}] = {\frakt g}^1, \forall 0\ne X\in {\frakt g}.$$
\end{prop}
\begin{pf}
Consider ${\frakt g}$ of class $\overline{\MD}$. Following the previous 
proposition, for every element $F\in {\frakt g}^* \setminus ({\frakt g}^1 
)^*$ , we have $\dim \Omega_F = \dim {\frakt g}$ and ${\frakt g}_F = 
\ker B_F = 0$\end{pf}. Suppose that there exists some element $X\in 
{\frakt g}$, $X \ne 0$ and $[X,{\frakt g}] \ne {\frakt g}^1$. There 
exists hence an element $F\in {\frakt g}^*$ vanishes in ${\frakt g}^1$ 
and doesn't vanish on $[X,{\frakt g}]$. This means that $0 \ne X \in \ker 
B_F= 0.$ This contradiction proves  that $[X,{\frakt g}] = {\frakt g}^1$.

Conversely, Suppose that for every $X\ne 0$ $[X,{\frakt g}] = {\frakt 
g}^1$. If $F\in {\frakt g}$ vanishes on ${\frakt g}^1$, we have $\ker B_F 
= {\frakt g}$. In this case, ${\frakt g}_F = {\frakt g},$ $\dim G_F = 
\dim {\frakt g}$, $\dim \Omega_F = 0$.

If $F^*\in {\frakt g}^*$ doesn't vanish on ${\frakt g}^1$, we have $[X, 
{\frakt g}] = {\frakt g}^1, \forall X\ne 0$, $\ker B_F =0$ and hence $\dim 
\Omega_F = \dim{\frakt g}$.
 
\section{Classification Theorem}

Let us denote $\ad^1_Y := \ad_Y|_{{\frakt g}^1}$.

\begin{lem}
If ${\frakt g}$ is of class $\overline{\MD}$, the operators of type 
$\ad^1_Y, Y\in {\frakt g}$ are pairwise commuting 
\end{lem}
\begin{pf}
We have the well-known Jacobi identity
$$[X,[Y,Z]] + [Y,[Z,X]] + [Z,[X,Y]] = 0, \forall X,Y,Z \in {\frakt g}$$. 
In particular, if Z is in the commutative derived ideal ${\frakt 
g}^1$, then $$[Z, [X,Y]] = 0$$ and we have
$$(\ad_X\circ \ad_Y - \ad_Y \circ \ad_X )Z \equiv 0, \forall X,Y\in 
{\frakt g}, \forall Z\in {\frakt g}^1.$$
\end{pf}

Recall that the Lie algebra of affine transformations of the real 
straight line is described as follows.
The Lie group $\Aff{\Bbb R}$ of affine transformations of the real 
straight line is  the group of affine transformations of type 
$$x\in {\Bbb R} \mapsto ax+b,$$ for some parameters $a,b\in {\Bbb R}$ and 
$a\ne 0$. For this reason some time one refers this group by 
"ax+b"-group. It is easy to prove that 
$$\Aff{\Bbb R} \cong \left\{ \pmatrix a & b \\ 1 & 0\endpmatrix \quad 
\vert \quad a,b \in {\Bbb R}, a \ne 0\right\}.$$  

It is easy to see that its Lie algebra $\aff{\Bbb R} = 
\Lie\Aff{\Bbb R}$ is
$$\aff{\Bbb R} = \left\{ \pmatrix \alpha & \beta \\ 0 & 0 \endpmatrix \quad 
\vert \quad \alpha,\beta\in {\Bbb R} \right\}.$$ By an easy direct 
calculation 
it is easy also to see that the Lie algebra $\aff{\Bbb R}$ is generated 
by two generators with the only nontrivial Lie brackets 
\index{brackets!Lie -} $ [X,Y] = Y$, i.e. 
$$\aff{\Bbb R} = \{aX + bY \quad \vert \quad [X,Y]=Y; a,b \in {\Bbb R}\}.$$

Now let us consider the group $\Aff{\Bbb C}$ of complex affine 
transformations of the complex straight line. The most easy method is to 
consider $X,Y$ as complex generators,  $X = X_1 +iX_2$ and $Y = Y_1 +iY_2$. 
Then from the relation $[X,Y]= Y$ we get
$$[X_1 +iX_2,Y_1+iY_2] = Y_1 + iY_2.$$
This means that
$$[X_1,Y_1] - [X_2,Y_2] +i([X_1,Y_2] + [X_2,Y_1]) = Y_1 +iY_2.$$
This means that the Lie algebra $\aff{\Bbb C}$ of affine transformations of 
the complex straight line is a real 4-dimensional Lie algebra, having 
4 generators with the only nonzero Lie brackets \index{brackets!Lie -}
$$[X_1,Y_1] -[X_2,Y_2] = Y_1, \quad [X_2,Y_1]+[X_1,Y_2] = Y_2.$$
Later in the proof of the theorem of classification what follows we shall 
choose another basis noted also by the same letters to have more clear Lie 
brackets of 
this Lie algebra, $$\begin{array}{ll} [X_1, Y_1] = Y_1, & [X_1, Y_2] = Y_2\cr
				      [X_2,Y_1] = Y_2, & [X_2, Y_2] = - Y_1\end{array}$$ 	

\begin{thm}
Up to isomorphism every algebra of class $\overline{\MD}$ is one of the 
following:
\begin{itemize}
\item Commutative Lie algebra,
\item Lie algebra $\aff{\Bbb R}$ of affine transformations of the real 
straight line,
\item Lie algebra $\aff{\Bbb C}$ of affine transformations of the complex 
straight line. \end{itemize}
\end{thm}
\begin{pf}
{\bf Step 1.}
\begin{clm}
There are only two possibilities: $\dim{\frakt g}^1 = $ either 1 or 2
\end{clm}
Consider the representation $\ad^1$ of Lie algebra ${\frakt g}$ in 
${\frakt g}^1$. The criterion $\overline{\MD}$ is just the irreducibility 
of this representation. In other hand, the operators $\ad^1_Y$ commute 
one with others. There is therefore an invariant complex line $D$ in the 
complexification ${\frakt g}_{\Bbb C}^1 := {\frakt g}^1 \otimes_{\Bbb R} 
{\Bbb C} $ for all operators $\ad^1_Y,Y\in{\frakt g}$. There are two 
cases : 
\begin{itemize}
\item {\it $D$ coincide with its complex conjugation $\overline{D}$}.  In 
this case $D = \delta \otimes_{\Bbb R} {\Bbb C}$, where $\delta \subset 
{\frakt g}^1$ is a real line invariant under all the operators $ad^1_Y$, 
$Y\in {\frakt g}$.

\item $D \ne \overline{D}$. In this case we have
$D \otimes \overline{D} = \delta \otimes_{\Bbb R} {\Bbb C}$, where 
$\delta$ is some real 2-dimensional plane, invariant under all the operators
$\ad^1_Y, Y\in {\frakt g}$. 
\end{itemize}
In both the cases we have $\delta = {\frakt g}^1$, in virtue of the 
irreducibility of the representation $\ad^1$. 

We do some remarks:
\begin{rem}
In the first case the existence of one dimensional real line, invariant 
under all the operators $\ad^1_Y$, $Y\in{\frakt g}$ deduces that 
$\dim{\frakt g}^1 = 1$.
\end{rem} 
\begin{rem}
If $[Z_1,{\frakt g}^1] = [Z_2,{\frakt g}^1] \equiv 0$ for some 
$Z_1,Z_2\in{\frakt g}$, then $[Z_1,Z_2] = 0$. 
\end{rem}
Really, following Jacobi identities, we have
$$[[Z_1,Z_2],{\frakt g}] + [[Z_2,{\frakt g}],Z_1] + [[{\frakt 
g},Z_1],Z_2] \equiv 0.$$ The last two summands are 0, then 
$$[[Z_1,Z_2],{\frakt g}] \equiv 0.$$ Following the 
$\overline{\MD}$-criterion we have $[Z_1,Z_2] = 0$.

{\bf Step 2.}

{\it Case 1: $\dim{\frakt g}^1 = 1$}.

Choose some $Y\ne 0$ in ${\frakt g}^1$. Following the 
$\overline{\MD}$-criterion $\ad_Y : {\frakt g} \to {\frakt g}^1$ is an 
surjection then there exists some $X \in {\frakt g}$ such that $\ad_Y(X) 
= - Y$, i.e. $[X,Y]=Y$.
We show that $\ker\ad_Y = {\frakt g}^1$. Following the remark 3.2, 
$\ker\ad_Y$ is commutative. If $\ker\ad_Y \ne{\frakt g}^1$, then there 
exists $X_1 \in \ker\ad_Y$, but $X_1$ is not in ${\frakt g}^1$. Because, 
$[X_1,X] \in {\frakt g}^1$, we have
$$[X_1, X] = \lambda Y,$$
for some $\lambda\in {\Bbb R}$. Then, $[X_1 +\lambda Y, {\frakt g}] 
\equiv 0$. From the $\overline{\MD}$-criterion we deduces that 
$X_1 +\lambda Y = 0 $, what  contradicts to the assumption that 
$\ker\ad_Y \ne 0$. Thus The Lie algebra has two generators $X,Y$ with the 
only nonzero relation $[X,Y]=Y$, i.e. ${\frakt g} \cong \aff{\Bbb R}$.

{\it Case 2: $\dim{\frakt g}^1 = 2$}

Suppose that ${\frakt g} = {\frakt g}^1 \otimes L,$ for some $L$ with $\dim 
L > 2$. Consider the operator $\ad_X : L \to {\frakt g}^1, \forall 
X\in{\frakt g}^1.$ Because $\dim L > \dim {\frakt g}^1$, there exists 
$Y\in L,Y\ne 0$ such that  
$\ad_XY = 0$. We shall prove that $[Y,{\frakt g}^1] = 0$.
Indeed, if $[Y,{\frakt g}^1] \ne 0$, $X$ must be the unique 0-vector of 
the operator $\ad^1_Y$. Because the operators $\ad^1_Y$, commute one with 
another, ${\Bbb R}X$ become an invariant space for all 
$\ad^1_T,T\in{\frakt g}^1$. Following Remark 3.1, $\dim{\frakt g}^1 =1$. 
This contradiction prove that $[Y,{\frakt g}^1] \equiv 0$. If the 
co-dimension of ${\Bbb R}Y$ in $L$ stills bigger than 2, there exists a 
vector $Z$ in its the complement, such that $[Z,{\frakt g}^1] \equiv  0$. 
We continue this procedure until the case where $\dim L = 2$. In this 
case, $L = L_1 \otimes L_2$, such that $\dim L_1 = 2$ and $[L_2,{\frakt 
g}^1] \equiv 0$. Following Remark 3.2, $L_2$  must be commutative. If 
$L_2\ne 0$, then $\dim L_2 \geq 2$, because the K-orbit have even dimension.

Consider the surjective map 
$$\ad_Y : {\frakt g}= {\frakt g}^1 \oplus L_2 \oplus L_1  \to {\frakt 
g}^1.$$ if $Y\in {\frakt g}^1 \oplus L_2$, we have $[Y, {\frakt g}^1 
\oplus L_2] = 0.$ Therefore, the operator $\ad_Y : L_1 \to {\frakt g}^1 $
is a surjection for all $Y\in {\frakt g}^1 \oplus L_2  \setminus \{0\}$. 
Because $\dim L_1 = \dim {\frakt g}^1 = 2$, the operator $\ad_Y$ is an 
isomorphism. We have a linear map 
$${\frakt g}^1 \oplus L_2 \setminus \{0\} \to ISO(L_1,{\frakt g}^1).$$
Moreover, following $\overline{\MD}$-criterion, this is an injective map. 
This is a contradiction.
 Thus, $L_2 = 0$, ${\frakt g} = {\frakt g}^1 \oplus L$, where $\dim 
{\frakt g}^1 = \dim L = 2$. Suppose that $Y_1,Y_2$ provide a basis of 
${\frakt g}^1$. Because, $\ad_{Y_i} : L \to {\frakt g}^1$ is surjective, 
there exists some $X_1\in L$ such that $[X_1,Y_1] = Y_1$, i.e. 
$\ad^1_{X_1}Y_1 = Y_1.$ Following Remark 3.1, $Y_1$ can not be the unique 
eigen vector up to a scalar multiple with eigenvalue 1 of the operator 
$\ad_{X_1}^1$. Thus $\ad^1_{X_1} = Id.$

Consider epimorphism $\ad_{X_1} : {\frakt g} \to {\frakt g}^1.$ We have 
$\dim \ker\ad_{X_1} = 2$, thus there exists $X_2'$ such that $[X_1,X_2'] = 
0$ and $X_1,X_2',Y_1,Y_2$ form a basis of ${\frakt g}$.
Because $\ad^1_{X_1} = Id, $ $\ad^1_{X_2'}$ can not have an eigenvector, 
following $\overline{\MD}$-criterion. Change basis in ${|frakt g}^1$ we have
$$\ad_{X_2'}^1 - \pmatrix a & -b\\ b & a \endpmatrix,$$ for some $a,b \in 
{\Bbb R}$, $b\ne 0$. Choose $$X_2 = {1\over b} (X_2' - aX_1),$$ we have
$$\ad^1_{X_2} = \pmatrix 0 & -1 \\ 1 & 0 \endpmatrix .$$
Thus Lie algebra ${\frakt g}$ is isomorphic to to an algebra generated by 
 a basis $X_1,X_2,Y_1,Y_2$ with the following only nonzero brackets
$$\begin{array}{ll} [X_1,Y_1] = Y_1, & [X_1, Y_2] =  Y_2 \cr
                [X_2, Y_1] = Y_2,  & [X_2, Y_2] = -Y_1\cr \end{array}$$
This means that Lie algebra ${\frakt g}$ is isomorphic to $\aff{\Bbb C}$.   

\end{pf}
\begin{rem}
This theorem let us to restrict our consideration for the class 
$\overline{\MD}$ practically to Lie algebras of affine 
transformations of the real or complex straight lines and the 
corresponding Lie groups
\begin{itemize}
\item the Lie group $\Aff{\Bbb R}$ of all the transformations of the real 
straight line,
\item the universal covering $\widetilde{\Aff{\Bbb R}}$, which is just 
isomorphic to the connected component $\Aff_0{\Bbb R}$ of the identity 
element, 
\item the Lie group of complex affine transformations of the complex 
straight line,
\item and its universal covering $\widetilde{\Aff{\Bbb C}}$ 
which is isomorphic to the connected component of identity 
$\Aff_0{\Bbb R}$
\end{itemize} 
\end{rem}

\section{Bibliographical Remarks}
The main idea to classify of this class of Lie algebra $\overline{\MD}$ 
belongs to the author of this book. He posed this problem for his 
postgraduate student H. H. Viet and a young colleague V. M. Son and 
solved together with them. The 
solution of this question was published in J. Operator Theory \cite{sonviet}.
Professor P. Cartier, during his scientific trip in Hanoi, helped them to 
shortcut the long proof.

\chapter  {The Structure of C*-Algebras of $\overline{\MD}$-Groups}  
\section{The C*-Algebra of $\Aff{\Bbb R}$}

\subsection{Statement of Theorems}
 
Let us denote by $G$ the group of affine transformations of the real 
straight line. We want to study the structure of this group and its 
C*-algebra $C^*(\Aff\Bbb R)$. To do this, we need first known its dual 
object, i.e. all its unitary representations up to unitary equivalence.
 
\begin{thm}{see \cite{gelfandnaimark}}
Each irreducible unitary representation of the group of affine 
transformations of the real straight line, up to unitary equivalence 
belongs to the following list of nonequivalent irreducible unitary 
representations:
\begin{enumerate} 
\item[a)] the representation  $S$, realized in th space 
$L^2({\Bbb R}^*, \frac{dx}{ \vert x \vert})$, where  ${\Bbb R}^* := 
{\Bbb R} \setminus (0)$, and the action is given by the formula
$$(S_gf)(x) = e^{\sqrt{-1}bx} f(ax),\quad\mbox{ where } g =
\pmatrix \alpha & b \\  0 & 1 \endpmatrix $$
\item[b)]
the representations $U^\varepsilon_\lambda$, realized in
${\Bbb C}^1$ and is given by the formula 
$$U_\lambda^\varepsilon (g) = \vert \alpha \vert^{\sqrt{-1}\lambda}.
(\sgn\alpha)^\varepsilon, \quad\mbox{ where }\lambda \in {\Bbb R}; 
\varepsilon = 0,1.$$

\end{enumerate}
\end{thm}

\begin{pf} see \cite{gelfandnaimark} and section 1.2 below.
\end{pf}

This list of irreducible unitary representations give us also the 
corresponding list 
of irreducible non-degenerate *-representations of the corresponding group 
C*-algebra 
$C^*(G)$. The next deal is to study the structure of this group C*-algebra. In 
\cite{diep1}
the author proved that 
this group C*-algebra $C^*(G)$ can be considered as some extension of the 
C*-algebra of continuous functions on a compact by  the elementary C*-algebra of 
compact operators in a separable Hilbert space.

\begin{thm}
The C*-algebra with a formally jointed unity element 
$C^*(G)^\sim$ can be included in a short exact sequence of C*-algebras and 
*-homomorphisms
$$ 0 \longrightarrow {\cal K} \longrightarrow C^*(G)^\sim \longrightarrow 
C({\Bbb S}^1 \vee {\Bbb S}^1) \longrightarrow 0,$$ 
This means that the C*-algebra $C^*(G)^\sim$, following BDF theory, 
is uniquely determined, up to isomorphisms class,by an element, said to be  its
{\rm the index } and is denoted by  $Index\,C^*(G)$.
\end{thm}
\begin{pf}
\cite{diep1} and section 1.3 below
\end{pf}

The infinite irreducible unitary representation $S$ realizes the indicated inclusion. 
Because  $${\cal E}xt({\Bbb S}^1 \vee {\Bbb S}^1) \cong \Hom_{\Bbb 
Z}(\pi^1({\Bbb
 S}^1\vee {\Bbb S}^1),$$ it provides a homomorphism from
 $\pi^1({\Bbb S}^1 
\vee {\Bbb S}^1)$ to ${\Bbb C}^*$. In virtue of the fact that the homomorphism 
$$Y_\infty : {\cal E}xt({\Bbb S}^1 \vee {\Bbb S}^1)\cong 
\Hom_{\Bbb Z}(\pi^1({\Bbb S}^1 \vee {\Bbb S}^1), {\Bbb Z}^1)$$ 
is obtained by computing of indices Fredholm operators and the general type of 
elements of $\pi^1({\Bbb S}^1 \vee 
{\Bbb S}^1)$ is $g_{k,l} = [g_{0,1}]^k[g_{1,0}]^l$, $k, l \in {\Bbb Z}$, 
we have $$\Ind\,(g_{k,l}) = k.\Ind\; T(g_{1,0}) + l.\Ind\;T(g_{0,1}),$$ 
where $T$ is a *-isomorphism, corresponding to  $S$. Hence, we need only to 
compute a pair of indices  
 $\Ind\, T(g_{1,0})$ and
$\Ind\, T(g_{0,1})$. These indices can be computed directly by the methods of 
computing the indices of Fredholm operators. We obtained the following exact 
computed results. 
\begin{thm}
$$Index\, C^*(G) = (1,1) \in {\cal E}xt({\Bbb S}^1 \vee {\Bbb S}^1) 
\cong {\Bbb Z} \oplus {\Bbb Z}.$$
\end{thm}
\begin{pf}
\cite{diep1} and section 1.4 below.\end{pf}

These results give us an interesting  
{\it improvement for a question
of {\sc A.A. Kirillov}} raised in
{\sc the I.M. Gelfand's Seminar} at the Moscow University in 1974
the the extension of the C*-algebra $C^*(G)$ of the group of affine 
transformations of the real straight line is split. 
This example does also spirit to apply the K-functors to some more general group 
C*-algebras. This was also a reason conducting
{\sc G.G. Kasparov} to create  his beautiful general  
KK-functors.

Returning to the general question, we need first of all some preparation
on construction of all the irreducible unitary representations of groups.
In the advanced theory, we shall do this by the construction of the
multidimensional orbit method. Secondly, we shall propose some methods to
decompose the group C*-algebras into some repeated short exact sequences
of C*-algebras and finally we shall propose some methods to reduces the
indices to some computable cases by using the general
Atiyah-Hirzebruch-Singer index formula for elliptic operators on topological 
spaces.  This program shall be done in the second part of this book.

\subsection{Proof of Theorem 1.1}
Theorem 1.1, was well-known as the first result about infinite 
dimensional unitary representations, after what the fruitful theory was 
well developed. In this first pioneer work. 
{\sc I.  M. Gel'fand}  and {\sc  M. A.  Naimark} have proved the theorem, 
using only the classical analysis and functional analysis. Now the day, 
we have enough strong tools to prove it shorter and more clear, for 
example the method of small subgroups of
{\sc G. W. Mackey}\cite{mackey},   or the orbit method, said in this book.   

It is easy to see that $\Aff  {\Bbb R}= {\Bbb R}^* . \ltimes  {\Bbb R}$ is 
the well-known semi-direct product of 
${\Bbb R}^*:= {\Bbb R} \setminus \{0\}$ and ${\Bbb R}$.   
we use therefore the Mackey method of small subgroups for the 
commutative normal subgroup  $N  = {\Bbb R}$.
The dual object $\hat{N}\approx  {\Bbb R}$  of  $N  \approx {\Bbb  R}$ 
consists of characters $\chi_\lambda$, $\lambda \in {\Bbb R}$, defined by 
the formula  $$\chi_\lambda(n) = \exp( i\lambda n), \forall n\in N,$$
for each fixed $\lambda$.

The group $\Aff {\Bbb R}$ acts on the dual object  
$\hat{N} \approx {\Bbb R}$ following the formula
$$(g.\chi_\lambda)(n)   :=   \chi_\lambda(gng^{-1})  =  \chi_\lambda(an)  =
\chi_\lambda(n),$$
where $g = (a,b)\in \Aff {\Bbb R}$. We have therefore
$$g.\chi_\lambda = \chi_{a\lambda}.$$
Hence in $\hat{N} \approx {\Bbb R}$ there are only two orbits  
$$\{ 0 \} \quad \mbox{ and } \quad {\Bbb R} \setminus \{ 0 \}.$$
Following the Mackey theory of induction from small subgroups,   the orbit 
$\{ 0 \}$
has the stabilizer subgroup, coincided with the whole group $\Aff {\Bbb  R}$. 
The characters of this stabilizer are in a $1-1$ correspondence with    
elements of the dual object
$$\widehat{{\Bbb R}^\times_+}=\widehat{{\Bbb R} \setminus \{ 0\}},$$
i.e. the unitary representations  
$U^\varepsilon_\lambda,   \forall  \varepsilon = 0,1;
\lambda \in {\Bbb  R}$. On the orbit  ${\Bbb R}  \setminus \{0\}$ choose 
a fixed point $\chi_1$,  i.e. with $\lambda = 1$  .  The stabilizer of the
fixed point
$\chi_1$ is $G_{\chi_1}= N = {\Bbb R}$. Hence, following the Mackey theory,
the orbit ${\Bbb  R}\setminus \{0\}$  corresponds to a unique 
representation $\Ind^G_N\chi_1 $, following the induction rule
$$(\Ind^G_N(a,b))f(x)= \exp(ibx) f(xa),$$
realized in the space $L^2(X,   d\mu_X(x))$  with  $X  = N\setminus G
\approx {\Bbb R}^* \setminus \{0\}$ and  $d\mu_X(.)= {dx \over \vert 
x\vert}$ as the quasi-invariant measure on the corresponding orbit. 
It is easy to see that this representation is exactly the representation 
$S$ in the statement of the theorem. Following the Mackey theory this list
is the complete list of all irreducible unitary representations of the 
group. The theorem is therefore proved. 

\subsection{Proof of Theorem 1.2}

As usually,  we add to the group C*-algebras the formally jointed unity 
element. It is easy to see that the representation $S$ is an exact 
representation, i.e. 
$\Ker S = \{0\}$. Because our group is a type   I group, i.e.  the 
group C*-algebra is  GCR, in  $C^*(G)$ there is an ideal, isomorphic to 
the ideal 
${\cal K}(H)$ of compact operators, where $H \cong L^2({\Bbb R}  ^*, 
{dx\over |x|})$.

\begin{lem}
Consider $\varphi\in L^1(\Aff{\Bbb R},   {dadb\over |a|})$. The 
conditions what follow are equivalent 
\begin{enumerate}
\item[1)] $$\varphi\in   \bigcap_{\varepsilon,\lambda}   \Ker
U^\varepsilon_\lambda,$$
\item[2)] $$\int_{-\infty}^{+\infty} \varphi(a,b) db =  0,  \mbox{ a. e. 
w.  r. t. the measure } {da\over |a|}.$$
\end{enumerate} 
\end{lem} 
\begin{pf}
Suppose $\psi(a)   :=  \int_{-\infty}^{+\infty}  \varphi(a,b)  db$.  The 
function
$\psi$ shall be decomposed into the sum of the even part $\psi_1$  and the  
odd part $\psi_2$. We have
$$\begin{array}{ll}
\varphi\in   \bigcap_{\varepsilon,\lambda}   \Ker
U^\varepsilon_\lambda   \Leftrightarrow   \int_{0<|a|<+\infty}
\int_{-\infty<b<+\infty}  |a|^{i\lambda}(\sgn  a)^\varepsilon  \psi(a)  {da
\over |a|} &= 0,\forall \lambda, \varepsilon\cr
  \int_{0<|a|<+\infty}   |a|^{i\lambda}(\sgn   a)^\varepsilon   \psi(a)
{da\over |a|} &= 0,\forall \lambda,\varepsilon .\cr    \end{array}$$
Because  $\psi_1$ is even and  $\psi_2$ is odd 
the above assertion is equivalent to the two condition what follow:
$$\cases \int_{0<|a|<+\infty}  |a|^{i\lambda}\psi_2(a) {da\over  |a|} =0 &,
\forall \lambda \in {\Bbb R}\\
\int_{0<|a|<+\infty}  |a|^{i\lambda}\psi_1(a) \sgn  a {da\over  |a|} = 0 &,
\forall \lambda\in {\Bbb R}.\endcases$$
Following the parity of $\psi_2$ and the exactness of the Fourier 
transformation, we have the equivalent conditions
$$\psi_1(a) = 0, \mbox{ a. e. w. r. p. t. } {da\over 
|a|},$$
$$\psi_2(a) = 0, \mbox{ a. e. w. r. t. } {da\over 
|a|},$$
and therefore we have an equivalent condition $$\psi(a) = 0, \mbox{ a. e. 
w. r. t. } 
{da\over |a|}.$$ 
\end{pf}

\begin{lem} Assume $\varphi\in L^1(\Aff{\Bbb 
R},{dadb\over|a|})$,satisfies the condition 
$$\int_{-\infty}^{+\infty}\varphi(a,b) db = 0, \mbox{ a. e. on } {\Bbb
R}\setminus \{0\} \mbox{w. r. t.} {da\over |a|}. $$ Then $S(\varphi)$
is a compact operator,   where $S$ in th single irreducible unitary 
representation of the group $\Aff{\Bbb R}$.  
\end{lem}
\begin{pf}
In \cite{gelfandnaimark} was proved separately. But we can deduce this
Lemma also from our compactness criteria in the advanced part.   
To avoid some repetition we omit this proof here. 
\end{pf}

\begin{lem}
The infinite dimensional representation $S$ provides n isomorphism of 
C*-algebras
$$\CD \bigcap_{\varepsilon,\lambda}  \Ker U^\varepsilon_\lambda  @>\cong >>
{\cal K}(H).\endCD$$ 
\end{lem}
\begin{pf}
1. First of all we must show that the representation $S$  is an exact 
representation, i.e. in  $C^*(G)$,  $\Ker S = 0$. Indeed,  Following the 
results of {\sc
J. M. Fell}[47],  every representation  $U^\varepsilon_\lambda$  is weakly
contained in the representation $S$ in the topology of the dual object 
$\widehat{\Aff{\Bbb R}}$.
This means that if   $S(\varphi)   =   0$,   then
$U^\varepsilon_\lambda(\varphi) = 0,   \forall \varepsilon = 0,1;   \forall
\lambda\in {\Bbb R}$. Hence,
$$\sup_{\pi\in\hat{G}} \vert \pi(\varphi)\Vert = 0,   \mbox{  i.e. 
} \Vert \varphi \Vert_{C^*(G)} = 0.$$
We have from here, $\varphi =0$.

2. Following the previous lemmas we have an inclusion 
$$S(\bigcap_{\varepsilon,\lambda}   \Ker   U^\varepsilon_\lambda)
v\hookrightarrow {\cal K}(H).$$ Because
${\cal K}(H)$ is an elementary  C*-algebra,   see 
{\sc J. Dixmier}[37],  one rests to show that 
$$S(\bigcap_{\varepsilon,\lambda}   \Ker
U^\varepsilon_\lambda)$$ is a two sided closed ideal in ${\cal K}(H)$.

Assume that $K   \in   S(\bigcap_{\varepsilon,\lambda}   \Ker
U^\varepsilon_\lambda)$ and $A \in  {\cal K}(H)$ is an arbitrary element. 
Because the group is of type I,
following the  Dixmier-Glimme-Sakai theorem,  the image 
$S(C^*(G))$ contain at least one compact operator. There exists therefore 
an element $\varphi\in C^*(G)$ such that $$A = S(\varphi).$$
Following the assumption, there exists 
$$K = S(\varphi_1),   \mbox{ where } \varphi_1 \in 
\bigcap_{\varepsilon, \lambda} \Ker U^\varepsilon_\lambda.$$
We have in this case
$$A.K = S(\varphi)S(\varphi_1) = S(\varphi \* \varphi_1),$$
$$K.A = S(\varphi_1)S(\varphi) = S(\varphi_1 \* \varphi).$$

Because $$U^\varepsilon_\lambda(\varphi   \*   \varphi_1)   =
U^\varepsilon_\lambda(\varphi)   U^\varepsilon_\lambda(\varphi_1)   =   0,
\forall \varepsilon = 0,1;\forall \lambda \in {\Bbb R},$$
we have
$$\varphi \* \varphi_1   \in   \bigcap_{\varepsilon,\lambda}   \Ker
U^\varepsilon_\lambda,$$
hence,
$$A.K =S(\varphi\*\varphi_1)   \in   S(\bigcap_{\varepsilon,\lambda}   \Ker
U^\varepsilon_\lambda).$$y analogy, we also have
$$K.A \in S(\bigcap_{\varepsilon, \lambda} \Ker U^\varepsilon_\lambda).$$
Thus, $$S(\bigcap_{\varepsilon,\lambda}\Ker U^\varepsilon_\lambda)
= {\cal K}(H).$$
\end{pf}

\subsubsection{End of Proof of Theorem 1.2}

Because the representation $S$  is exact then we can identify th ideal
$\bigcap_{\varepsilon,\lambda} \Ker U^\varepsilon_\lambda$ 
with its image 
$$S(\bigcap_{\varepsilon,\lambda}  \Ker  U^\varepsilon_\lambda)\cong  {\cal
K}(H). $$ We show that the quotient C*-algebra  
$C^*(G)/{\cal K}(H)$ is commutative. 

Indeed, because the representations $U^\varepsilon_\lambda$ are all of
dimension 1, then $$U^\varepsilon_\lambda(\varphi\*\psi - \psi\*\varphi) =
U^\varepsilon_\lambda(
 \varphi)U^\varepsilon_\lambda(\psi)   -   U^\varepsilon_\lambda(\psi)
U^\varepsilon_\lambda( \varphi) = 0,   \forall \varepsilon = 0,1;   \forall
\lambda \in {\Bbb R}.$$

One deduce that $$\varphi\*\psi - \psi\*\varphi) \in \bigcap_{\varepsilon,
\lambda} \Ker U^\varepsilon_\lambda,$$ i.e.  the quotient C*-algebra
$C^*(G)/{\cal K}(H)$ is a commutative C*-algebra. One shows that the
maximal ideals of the quotient C*-algebra are just the ideals obtained
from One deduce that $$\varphi\*\psi  -  \psi\*\varphi)  \in \bigcap_{\varepsilon,
\lambda} \Ker U^\varepsilon_\lambda,$$
i.e.  the quotient C*-algebra $C^*(G)/{\cal  K}(H)$ 
is a commutative C*-algebra. One shows that the maximal ideals 
of the quotient C*-algebra are just the ideals obtained from  
the corresponding ideals, which contain the ideal ${\cal K}(H)$.
It is easy to see from here that $$C^*(\Aff{\Bbb R})/{\cal K}(H) \cong {\Bbb 
C}({\Bbb S}^1 \vee {\Bbb S}^1),$$ in other words, we have the exact 
sequence in the theorem 1.2:
$$\CD 0 @>>> {\cal  K}(H)  @>>>  \widetilde{C^*(\Aff{\Bbb  R})}  @>>> {\Bbb
C}({\Bbb S}^1 \vee {\Bbb S}^1)@>>> 0.\endCD$$

\subsection{Proof of Theorem 1.3}

It is well-known that $$\Ext({\Bbb S}^1 \vee {\Bbb S}^1) = \Hom_{\Bbb Z}(
K^{-1}({\Bbb S}^1 \vee {\Bbb S}^1),{\Bbb Z} ).$$
Because of
$$ K^{-1}({\Bbb S}^1 \vee {\Bbb S}^1) \cong {\Bbb Z} \oplus {\Bbb Z},$$
we have
$$\Hom_{\Bbb Z}(K^{-1}({\Bbb S}^1 \vee {\Bbb Z}^1),{\Bbb Z}) \cong {\Bbb Z}
\oplus {\Bbb Z}$$
and so ,
$$\Ext({\Bbb S}^1 \vee {\Bbb S}^1) \cong {\Bbb Z} \oplus {\Bbb Z}.$$

We must prove that the  C*-algebra of the group of affine transformations 
of the real straight line corresponds to the pair $(1,1)$ in ${\Bbb Z} \oplus
{\Bbb Z}$. The proof is rather long and we must divide it in a sequence 
of steps. 

First of all, we must find out the generators of the group 
$\Ext({\Bbb S}^1  \vee {\Bbb  S}^1) \cong  {\Bbb Z}  \oplus {\Bbb  Z}$. 
Remark that we added to the C*-algebra the formal unity element is 
equivalent 
to consideration the one point compactification of its dual object.   
Certainly that the one point compactification of a pair of 
parallel real straight lines is just    
${\Bbb S}^1 \vee {\Bbb  S}^1$.  This means that we can enumerate the 
points of 
${\Bbb S}^1 \vee  {\Bbb  S}^1$  by the pairs of numbers $(\varepsilon,   
\lambda), \forall \varepsilon = 0,1; \lambda \in {\Bbb R}$, i.e.
$$X = {\Bbb S}^1 \vee {\Bbb  S}^1 = \{(\lambda,   \varepsilon)  \;   :   \;
\lambda \in {\Bbb R}; \varepsilon = 0,1\} \cup \{ \infty \}.$$

From definition,   $$\pi^1({\Bbb S}^1 \vee {\Bbb S}^1) := [ {\Bbb S}^1 \vee
{\Bbb S}^1 , {\Bbb C}^*] \cong {\Bbb Z} \oplus {\Bbb Z},$$
where ${\Bbb C}^* := {\Bbb  C} \setminus \{ 0\}$ and $[.,.]$  denote 
the homotopy class of maps. More precisely, with each pair of integers
 $(k,l) \in {\Bbb Z}
\oplus {\Bbb Z}$ we can choose a map $g_{k,l}$ as follows,
$$g_{k,l}(\lambda, \varepsilon) : = \cases \exp[k.i.2.\arcctg \lambda/2] &,
\mbox{ if      } \varepsilon = 0, \lambda \in {\Bbb R}\\
  \exp [k.i.2\arcctg \lambda/2] &,   \mbox{ if } \varepsilon =1,   
\lambda \in {\Bbb R}\\
  1 &, \mbox{ if } \lambda = \infty .\endcases$$ Then the homotopy 
classes $[g_{k,l}]$ provide a group which is equivalent with  
${\Bbb Z}  \oplus {\Bbb Z}$
with two generators $[g_{0,1}]$ and $[g_{1,0}]$,
$$[g_{k,l}] = [g_{0,1}]^k.[g_{1,0}]^l.$$

Recall that for an arbitrary compact  $X$,   we have an exact sequence
$$\CD 0 @>>> Ext_{\Bbb Z}(K^0(X),{\Bbb Z}) @>>>  \Ext(X) @>Y_\infty>> \Hom(
K^{-1}(X),{\Bbb Z}) @>>> 0\endCD$$
Moreover, if $X \subset  {\Bbb C}$ and  $\dim X \leq  1$,   then 
the homomorphism
$Y_\infty$ is an isomorphism,  see {\sc L. G. Brown,  R. G. Douglas and 
P. A. Fillmore}[7]. Because $X\subset {\Bbb C}$, 
$$K^{-1}(X) = [X, \GL_1({\Bbb C})] = [X, {\Bbb C}^*] = \pi^1(X).$$
Assume $\tau$ to be an arbitrary extension, i.e.
$\tau$ is an inclusion  $$\tau  :   {\Bbb C}(X) \hookrightarrow {\cal
A}(H).$$ The the invertible elements $g\in {\Bbb C}(X)^*$ must be mapped 
into invertible elements $\tau(g)$  in  Calkin algebra ${\cal A}(H)$. 
Hence, $\tau(g)$ is defined by an unique  Fredholm operator,  up to 
compact operator perturbations. Moreover, the index 
$\Ind\tau(g)$  depends only on the homotopy class of map
$g\in {\Bbb C}(X)^*$.

Then to our extension $\tau$ corresponds a homomorphism   $\Ind$
in $\Hom_{\Bbb Z}(\pi^1(X), {\Bbb Z})$,   with $X = {\Bbb S}^1 \vee {\Bbb
S}^1$. We have
$$\Ind\tau([g_{k,l}]) = k.\Ind\tau([g_{0,1}]) + l.\Ind([g_{1,0}]) .$$

The irreducible infinite dimensional unitary representation $S$ of
the group $\Aff{\Bbb R}$ gives us an extension of type  
$$\CD 0 @>>> {\cal   K}(H)   @>>>   \widetilde{C^*(\Aff{\Bbb   R})}
@>U^\varepsilon_\lambda>>  {\Bbb  C}({\Bbb  S}^1  \vee  {\Bbb  S}^1) @>>> 0
\endCD$$.  following BDF theory,  to him corresponds an unique 
homomorphism of type  
$$\tau :  {\Bbb C}({\Bbb S}^1 \vee {\Bbb S}^1) \hookrightarrow {\cal A}(H).
$$

Remark that   if $\varphi\in  C^*(\Aff{\Bbb  R})$,  and if 
$U^\varepsilon_\lambda(\varphi) = g(\lambda,  \varepsilon),   \varepsilon =
0,1; \lambda\in  {\Bbb  R}$,   with  $g\in  {\Bbb  C}( {\Bbb S}^1 \vee {\Bbb
S}^1)^*$, then $S(\varphi)$ is a  Fredholm operator and
$$\Ind S(\varphi) = \Ind \tau(g).$$

We conclude that we should compete our study of the structure of the  
group C*-algebra following the program:

{\bf Program of computing indices }
\begin{enumerate}
\item[a.] Find out the functions  $\varphi_1,   \varphi_2 \in C^*(\Aff{\Bbb 
R})$ such that
$$U^\varepsilon_\lambda(\varphi_1)   =   g_{0,1}(\varepsilon,\lambda),
\forall \varepsilon = 0,1; \forall \lambda \in {\Bbb R},$$
$$U^\varepsilon_\lambda(\varphi_2)   =   g_{1,0}(\varepsilon,\lambda),
\forall \varepsilon = 0,1; \forall \lambda \in {\Bbb R},$$
\item[b.] Compute the indices of  Fredholm operators
$$\Ind S(\varphi_1) = \Ind\tau([g_{0,1}]),$$
$$\Ind S(\varphi_2) = \Ind\tau([g_{1,0}]),$$

Then,the topological invariant of our C*-algebra is just
$$Index C^*(\Aff{\Bbb R}) = (\Ind S(\varphi_1),\Ind S(\varphi_2)).$$
\end{enumerate}

\begin{lem}
$$-2 \int_{-\infty}^{+\infty}\exp \{ -2  |a| + i\lambda  a\} da =  {-8\over
\lambda^2 +4}\leqno{1.}$$
$$2 \int_{-\infty}^{+\infty} \exp\{-2|a|  +  i\lambda a\}\sgn  a  da  =  4i
{\lambda \over \lambda^2 +4}. \leqno{2.}$$ 
\end{lem}
\begin{pf}
We do the exact computation by integrating in part.

1. Put $${\cal J} = \int_0^\infty e^{-2a }\cos{\lambda a} da,$$ then
$${\cal J} = {1\over 2} - {\lambda^2\over 4}{\cal J}.$$ Hence,
$${\cal J} = {2\over \lambda^2 + 4}.$$
We have
$$\begin{array}{ll} -2 \int_{-\infty}^{+\infty}\exp{ \{ -2  |a| + 
i\lambda  a\}} da &=
-2 \int_{-\infty}^{+\infty}   e^{-2|a|}   \cos\lambda   a   da   -   2i
\int_{-\infty}^{+\infty}  e^{-1|a|}\sin \lambda a da\cr  
 &= -4 \int_0^{+infty} e^{-2a} \cos \lambda a da \cr
 &= -4 {\cal J}\cr
 &= {-8\over  \lambda^2 +4} \end{array}$$

2. By the same type computation.
\end{pf}          

\begin{lem}
$$\exp [i2\arcctg  ({\lambda\over  2})]  -1  =  {-8\over \lambda^2 + 4} +4i
{\lambda\over \lambda^2 +4}.$$ 
\end{lem}
\begin{pf}
We do some trigonometric transform.
Pose $\alpha = 2\arcctg{\lambda\over 2}.$  Then 
$${\lambda\over 2} = \ctg {\alpha\over 2} = t.$$
We have
$$\begin{array}{rl} \exp \{i 2\arcctg{\lambda\over 2} \} &= e^{i\alpha} 
= \cos\alpha + i\sin\alpha \cr
&= {t^2 -1 \over t^2 +1} + i{2t \over t^2 +1} \cr
&= {({\lambda\over 2})^2 -1 \over ({\lambda\over 2})^2 +1}+ i
{2({\lambda\over 2})\over ({\lambda\over 2})^2 +1}\cr
&= {\lambda^2 -4\over \lambda^2 +4}  +  i4{\lambda\over\lambda^2 + 4} 
.\end{array}  $$ \end{pf} 

\begin{lem}
Suppose $$\psi_i(a) = \cases 0 &, \mbox{ if } 1 < |a| < +\infty ,\\
                    -2a^2  (\sgn a)^{i-1}  &,   \mbox{ if  } 0  < |a| \leq
1,\endcases$$
for $i=1,2$. Then

$$\begin{array}{ll} {1.}\qquad &\int_{0<|a|<+\infty} 
|a|^{i\lambda}\psi_1(a)  {da\over |a|} = 
               \exp \{i2.\arcctg{\lambda\over 2} \} -1,\cr
2.\qquad &\int_{0<|a|<+\infty} |a|^{i\lambda}\psi_1(a) \sgn a {da\over |a|}
\equiv 0, \cr
3.\qquad &\int_{0<|a|<+\infty} |a|^{i\lambda}\psi_2(a) {da\over |a|} \equiv
0,\cr 
4.\qquad &\int_{0<|a|<+\infty} |a|^{i\lambda}\psi_2(a)\sgn a  {da\over |a|}
= 
               \exp \{i2.\arcctg{\lambda\over 2} \} -1.\end{array}$$  
\end{lem}
\begin{pf}
$$\begin{array}{rl} 1.\qquad &\int_{0<|a|<+\infty} |a|^{i\lambda}\psi_1(a) {da\over
|a|} \cr 
  &= -2 \int_{0<|a|\leq 1} |a|^{i\lambda +2} {da\over |a|} \cr
  &= -4 \int_0^1 a^{i\lambda +2} {da\over a} \cr
  &= -4 \int_0^1 e^{2\ln a + i\lambda \ln a} d\ln a \cr
  &= -4 \int_{-\infty}^0 e^{2a' + i\lambda a'} da'\cr
  &= \int_{-\infty}^0  \{-2e^{2a'}  +  2\sgn  a'.e^{2a'}\}e^{i\lambda a'}da'\cr 
  &= {-8\over \lambda^2 +4} + 4i{\lambda\over \lambda^2 +4}\cr
  &=   \exp \{i2.\arcctg{\lambda\over 2} \} -1.\end{array}$$
By the same way the other integrals are computed. 
\end{pf}

\begin{lem}
Assume that
$$\tilde{\varphi}_i(a,b)   =   \psi_i(a){1\over\sqrt{2\pi}}\exp(-{b^2
\over 2}),\quad \varphi_i = \tilde{\varphi}_i + 1, i = 1,2,$$
where $1$ is the formal identity element of  
$C^*(\Aff{\Bbb C}^\sim$. Then
$$U^\varepsilon_\lambda(\varphi_1) = g_{1,0}(\lambda,\varepsilon),  \forall
\lambda\in {\Bbb R}; \varepsilon = 0,1,$$
$$U^\varepsilon_\lambda(\varphi_2) = g_{0,1}(\lambda,\varepsilon),  \forall
\lambda\in {\Bbb R}; \varepsilon = 0,1.$$ 
\end{lem}
\begin{pf} 
Using the previous lemmas, we have
$$\begin{array}{ll}U^\varepsilon_\lambda\tilde{\varphi}_i)   &=
\iint_{0<|a|<+\infty\atop -\infty < b< +\infty}   |a|^{i\lambda}(\sgn
a)^\varepsilon\tilde{\varphi}_i(a \,b) {dadb\over |a|}\cr 
 &={1\over \sqrt{2\pi}}\int_{-\infty}^{+\infty} \exp(-{b^2\over 2})db
\int_{0<|a|<+\infty}  |a|^{i\lambda}(\sgn  a)^\varepsilon\psi_i(a) {da\over
|a|}. \end{array}$$ 
\end{pf}

We have thus found out the functions $\varphi_i,   i=1,2$  on the 
group, satisfying all the necessary conditions. following the program of 
computing indices, one rests to compute the indices of Fredholm operators  
$\Ind S(\varphi_i), i=1,2$. Firstly,   we write out the explicit action 
formula of the operators  $S(\varphi_i), i=1,2$.

From definition, we have
$$\begin{array}{rl} (S(\tilde{\varphi}_i)f)(x) 
&= -\sqrt{2\over\pi} \int_{-\infty}^{+\infty} \exp(-ibx -{b^2\over 2})db
   \int_{0<|a|<1} f(xa)a^2(\sgn a)^{i-1}{da\over|a|}\cr
&= -2\exp (-{x^2\over 2})\int_{-1}^1f(xa)|a|(\sgn a)^{i-1} da .\end{array}$$
Hence, we deduce the exact action formulas
$$[S(\varphi_i)f](x)   =   f(x)   -   \exp(-{x^2\over   2})   \int_{-1}^1
f(xa)|a|(\sgn a)^{i-1} da.$$

Thus, in order to compute the indices of these Fredholm operators,   
we consider the differential equations 

$$f(x) - \exp(-{x^2\over 2}) \int_{_1}^1 f(xa)|a| da = 0 \eqno{(1)}$$
$$f(x) - \exp(-{x^2\over 2}) \int_{_1}^1 f(xa)a da = 0 \eqno{(2)}$$

\begin{lem}
Each solution of the equation (1) must be an even function and each solution 
of the equation (2) must be an odd function, on the symmetric domain 
${\Bbb R}\setminus \{0\}$, if exist. 
\end{lem}
\begin{pf}
Suppose $f$ to be a solution of the equation (1). Because the domain 
is symmetric, we can decompose it into the sum  $f = f_1 + f_2$ 
of its even part  $$f_1 = {f + 
\check{f}\over 2}$$ and its odd part
$$f_2 = {f - \check{f}\over 2},$$
where
$$\check{f}(x) := f(-x), \forall x \in {\Bbb R}.$$
Then, for a fixed $x$,
$$\int_{-1}^1 f_2(xa)|a| da = 0.$$
Then following the equation (1),
$$\begin{array}{rl}f(x) &= \exp(-{x^2\over 2}) \int_{-1}^1 f(xa)|a| da\cr
  &= \exp(-{x^2\over 2})\int_{-1}^1 f_1(xa)|a|da\cr
  &= \exp(-{x^2\over 2})\int_{-1}^1 f_1(-xa) |a| da\cr
  &= f(-x). \end{array}$$

The second part is proved by the same way.
\end{pf}

From this lemma we can reduce our study of solutions of equations (1) and 
(2) on the domain  $0 < x < +\infty$,   the if necessary extend them 
following symmetry on the whole domain ${\Bbb R}^* = {\Bbb  R} \setminus 
\{0\}$.

\begin{lem}
In the Hilbert space $L^2({\Bbb R}^*, {dx\over |x|})$, we have 
$$\dim \Ker S(\varphi_i) = 1, \quad i= 1,2.$$ 
\end{lem}
\begin{pf}
This lemma is proved by reducing to some differential equation and then 
estimate the asymptotic behavior of solutions. 
For positive values of $x$, the equations (1) and (2) have the same 
form $$f(x) - \exp(-{x^2\over 2}) \int_0^1 f(xa)ada = 0.$$
It is the same as 
$$f(x) = {4\exp(-{x^2\over   2})\over   x^2}\int_0^x\xi   f(\xi)d\xi .
\leqno{(3)}$$ 
Put $$F(x) := \int_0^x \xi f(\xi)d\xi ,$$
we have an differential equation for  $F(x)$
$$F^\prime(x) - {4 \exp(-{x^2\over 2})\over x}F(x) = 0 \leqno{(4)}$$
Certainly, here we consider the generalized Sobolev derivatives.
For some fixed value $x_0$ in the domain  $0<x_0<+\infty$ 
we define a unique solution
$F(x)$ of the Cauchy problem,  and the function $f$ can be computed as 
$$f(x) = {F^\prime(x) \over x}.$$

The rest id to decide, whether the function $f$ belongs to the Hilbert space
$L^2({\Bbb R}^*,   {dx\over|x|})$. To do this, we study its asymptotic 
behavior when $x \to 0$ and when $x \to \infty$.

So, assume $F(x)$ to be a solution of the differential equation (4) and  
$$f(x) = {F^\prime(x)\over x},$$ then
$$x^2 f(x) - 4\exp(-{x^2\over 2}) \int_0^x f(\xi)d\xi = 0,$$
or the same
$$\begin{array}{rl} x^2 \exp ({x^2\over 2}) f(x) &= 4\int_0^x f(\xi)\xi 
d\xi\cr
           [x^2\exp({x^2\over 2})f(x)]^\prime &= 4 f(x)x \cr
     {[x^2\exp({x^2\over  2})f(x)]^\prime \over  x^2\exp({x^2\over 2})f(x)}
&= 4{\exp(-{x^2\over 2}) \over x}\cr
\ln | x^2\exp({x^2\over  2})f(x)|  &=  4\int_a^x  \exp(-{t^2\over  2})  {dt
\over t} + c\cr
           &\sim {\cases c_1  &, \mbox{ if } x \to \infty ,\\
                        4\ln x +  c_2 &,  \mbox{ if  } x \to  0,\endcases}
\end{array}$$
where $0<a<+\infty$, $c_1, c_2, c$ are some constants.
We obtain thus the asymptotic behavior of solutions as follows
$$x^2 \exp({x^2\over 2}) f(x) \sim  \cases e^{c_1} &,   \mbox{ if }  x \to
\infty,\\ x^4e^{C_2} &, \mbox{ if } x \to 0,\endcases$$
$$f(x) \sim \cases x^{-2}\exp(-{x^2\over 2}) &,  \mbox{ if } x \to +\infty
,\\ x^2 &, \mbox{ if } x \to 0 .\endcases$$
We conclude that the solutions are square-integrable with respect to the 
measure ${dx\over |x|}$ on ${\Bbb R}^* = {\Bbb R}\setminus \{0\}$.
\end{pf}

\begin{lem}
The image of $L^2({\Bbb  R}^*,{dx\over|x|})$ under the maps
$S(\varphi_i)$ are dense in itself; i.e.
$$\dim \Coker S(\varphi_i) = \{0\}, i=1,2.$$ 
\end{lem}
\begin{pf}
1. {\it Assume $g(x)$ is an even function,  then the solutions of the 
equation 
$$[S(\varphi_1)f](x) = f(x)  - 2\exp(-{x^2\over 2})\int_{-1}^1  f(xa)|a| da
= g(x)$$ are also even functions.}

Indeed, we decompose the function $f$ into the sum of its even part $f_1$ 
and odd part
$f_2$, we have
$$\int_{-1}^1 f_2(xa)|a| da = 0$$
and by this reason,
$$\begin{array}{rl} f(x) &= 2\exp({x^2\over 2}) \int_{-1}^1 f_1(xa)|a| da 
+ g(x)\cr &=2\exp({(-x)^2\over 2}) \int_{-1}^1 f_1(-xa)|a| da + g(-x)\cr
&= f(-x).\end{array}$$

2. An analogous assertion is valid for $\varphi_2$:

{\it Assume $g(x)$ to be an odd function,   then the solutions of the 
equation  
$$[S(\varphi_2)f](x) = f(x) - 2\exp(-{x^2\over 2})  \int_{-1}^1 f_1(xa)a
da = g(x)$$ are also odd functions .}

3. Assume $g(x)$ is an arbitrary function with compact support on 
${\Bbb R}^*= {\Bbb R}\setminus \{0\}$, in other words, 
there exists some number $N$, big enough such that
$$g(x) = 0, \quad \forall x;  |x| \leq {1\over N} \mbox{ or } |x| \geq N.$$
WE shall show that the equations
$$[S(\varphi_1)f](x) = g(x) \mbox{ and }\leqno{(1')}$$
$$[S(\varphi_2)f](x) = g(x).\leqno{(2')}$$
always have solutions in $L^2({\Bbb R}^*,{dx\over|x|})$

Assume $g = g_1  + g_2$  is the decomposition of $g$ into the sum of its 
even and odd parts. Then because $g_2$ is odd,
$$\int_{-1}^1 g_2(xa)|a| da = 0.$$
Hence,
$$g_2(x) - 2\exp(-{x^2 \over 2}) \int_{-1}^1 g_2(xa) |a| da = g_2(x).$$
We try to find the solutions of $(1')$ in form 
$$f = \tilde{f} + g_2,$$
where
$$\tilde{f}(x) - 2  \exp(-{x^2\over 2}) \int_{-1}^1 \tilde{f}(xa)  |a| da =
g_1(x).$$
Following 1., $\tilde{f}$ must be an odd function. We have
$$\tilde{f}(x) - 4\exp(-{x^2\over  2}) \int_0^1 \tilde{f}(xa) ada  = g_1(x)
= 0, \mbox{ if } |x| \leq {1\over n} \mbox{ or } |x| \geq N.$$
Thus outside the interval $[{1\over N}, N]$,   $\tilde{f}(x)$  must be an 
even solution of the equation
(1) and have the following asymptotic behavior
$$\tilde{f}(x) \sim \cases {1\over x^2}\exp(-{x^2\over 2}) &,  \mbox{ when }
x \to \infty,\\ x^2 &, \mbox{ when } x \to 0.\endcases$$

Because $g$ has compact support,   $\check{g}$ has also compact support 
and therefore the even and odd parts  $g_1 =
{g + \check{g}\over  2}$  and  $g_2  =  {g  - \check{g}\over 2 }$ have 
also the compact supports. Hence, 
$$f(x) = \tilde{f} + g_2 \in L^2({\Bbb R}^*, {dx\over|x|}).$$

By analogy, we can always solve the equations 
$$[S(\varphi_2)f](x)  =  g(x),$$  where  $g(x)$  has compact support.  
\end{pf}
The proof of the theorem 1.3 is therefore also accomplished.

\section{The Structure of $C^*(\widetilde{\Aff{\Bbb C}})$}

\begin{rem}
The exponential map $\exp : {\Bbb C} \to {\Bbb C}^* := {\Bbb C} \setminus 
\{ 0\}$, giving by  $z \mapsto e^z$  is just the covering map, and 
therefore $\widetilde{\Bbb C}^* \cong {\Bbb C}$. As a consequence, one 
deduces that $$\widetilde{\Aff{\Bbb C}}\cong {\Bbb C}\ltimes {\Bbb C} 
\cong \{ (z,w) \quad\vert \quad z,w\in {\Bbb C}\}$$
 with the following multiplication law 
$$(z,w).(z',w') := (z+z', w+e^zw').$$
\end{rem}

\begin{thm}
Up to unitary equivalence, every irreducible unitary representation of 
the universal covering $\widetilde{\Aff{\Bbb C}}$ of the group of 
affine transformations of the complex straight line is unitarily 
equivalent to one of the following list of its one-to-another 
nonequivalent irreducible unitary representations:
\begin{itemize}
\item The infinite dimensional irreducible representations $T_\alpha, 
\alpha\in {\Bbb S}^1$, realized in the Hilbert space $L^2({\Bbb 
R}\times{\Bbb S}^1)$ by the formula
$$T_\alpha(z,w)f(x) := \exp(i(\Re(we^x)+2\pi \alpha[{\Im(x+z)\over 
2\pi}]))f(x\oplus z),$$
where $(z,w)\in \widetilde{\Aff{\Bbb C}}$, $x\in {\Bbb R} \times {\Bbb 
S}^1 = {\Bbb C}/\{2i\pi {\Bbb Z}\}$, $f\in L^2({\Bbb R}\times {\Bbb 
S}^1)$ and $$x \oplus z := \Re(x+z) + i2\pi \{{\Im(x+z)\over 2\pi}\}.$$
\item
The unitary characters of the group, i.e. the 1-dimensional unitary 
representations $U\lambda, \lambda\in {\Bbb C}$, acting in ${\Bbb C}$ 
following the formula
$$U_\lambda(z,w) := e^{i\Re(z\overline{\lambda})}, \forall (z,w)\in 
\widetilde{\Aff{\Bbb C}}.$$
\end{itemize}
\end{thm}
\begin{pf}
In the group $\widetilde{\Aff{\Bbb C}}$ consider the normal subgroup 
$$N := \{ (0,a) ; a \in {\Bbb C} \}.$$
We have 
$$\begin{array}{rl}
(z,w).(0,a)(z,w)^{-1} &= (z,w)(0,a)(-z,-we^-z)\\
		      &= (z,w+e^za)(-z,-we^-z)\\
		      &= (0,w +e^z a - e^zwe^{-z})\\
		      &=(0,e^za)
\end{array}$$
This means that $N$ is a commutative normal subgroup. Its dual object 
$\hat{N}$ consists of the characters $\chi_\lambda, \lambda\in {\Bbb C}$,
$$\chi_\lambda(0,a) := e^\Re(a\overline{\lambda}).$$ The covering group 
$\widetilde{\Aff{\Bbb C}}$ acts on $\hat{N}$ following the formula
$$\begin{array}{rl}
((z,w)\chi_\lambda)(0,a) &= \chi_\lambda((z,w)(0,a)(z,w)^{-1})\\
			 &= \chi_\lambda(0,e^za)\\
			 &= \exp(i\Re(e^za\overline{\lambda}))\\
			 &= \exp(i\Re(a\overline{e^x\lambda}))\\
			 &= \chi_{\exp(\bar{z})\lambda}(0,a).
\end{array}$$ 

This means that under the action of $\widetilde{\Aff{\Bbb 
C}}$ on $\hat{N}$, there are only two orbits $\{0\}$ and $\hat{N} 
\setminus \{0\}$. Following the Mackey theory of induction from small 
subgroups, we have:
\begin{enumerate}
\item 
The one dimensional unitary representations (i.e. the unitary 
characters), corresponding to the K-orbit $\{ 0\}$, extended from the 
trivial representation of the normal subgroup $N$, being of type
$$U_\lambda(z,w) = e^{i\Re(z\overline{\lambda})}, \forall \lambda\in 
{\Bbb C}, (z,w)\in \widetilde{\Aff {\Bbb C}}.$$
\item
The infinite dimensional unitary representations $T_\alpha := 
\ind^G_{G_{\chi_1}}S_\alpha, \alpha\in {\Bbb S}^1$, acting on the Hilbert 
space $L^2({\Bbb R}\times {\Bbb S}^1)$, following the formulas
$$(T_\alpha(z,w)f)(x) = \exp(i(\Re(we^x) + 2\pi\alpha[{\Im(x+z)\over 
2\pi}])) f(x\oplus z),$$ where $(z,x)\in \widetilde{\Aff{\Bbb C}},$ $f\in 
L^2({\Bbb R}\times {\Bbb S}^1)$, $x\in {\Bbb R}\times {\Bbb S}^1 \approx 
{\Bbb C}/\{i2\pi {\Bbb Z}\}$, and $$x\oplus z := \Re(x+z) + 2\pi i \{ 
{\Im(x+z)\over 2\pi}\}.$$
\end{enumerate} 

Indeed, fix a point, say $\chi_1$ on the K-orbit ${\Bbb C} \setminus 
\{0\}$. The irreducible unitary representations of its stabilizer 
$$G_{\chi_1} = \{ (i2\pi n,h) \quad\vert\quad n\in {\Bbb Z}, h\in {\Bbb C} 
\},$$ which are multiples of the unitary character $\chi_1$, are the  
irreducible unitary representations of form
$$S_\alpha \in \widehat{G_{\chi_1}}, \alpha\in {\Bbb S}^1,$$
$$\begin{array}{rl}
S_\alpha(i2\pi n,h) & = \exp(i(\Re h.\bar{1} + \sin\alpha))\\
		    & = \exp(i(\Re h + 2\pi n\alpha))
\end{array}$$
Thus the induced representations $T_\alpha + \ind_{G_{\chi_1}}^G 
S_\alpha$ should be realized in the Hilbert space 
$$L^2(\widetilde{\Aff{\Bbb C}}/G_{\chi_1}) = L^2({\Bbb C}.{\Bbb C} /{\Bbb 
C}.\{i2\pi {\Bbb Z}\}) = L^2({\Bbb C}/\{ i2\pi {\Bbb C}) \cong L^2({\Bbb 
R} \times {\Bbb S}^1).$$
\end{pf}

\begin{prop}
The infinite dimensional induced unitary representations $T_\alpha = 
\ind_{G_{\chi_1}} S_\alpha$, $\alpha\in {\Bbb S}^1$ are *-homomorphism 
from the ideal $I := \bigcap_{\lambda\in{\Bbb C}} \ker U_\lambda$ into 
the ideal of compact operators ${\cal K}(H)$ in a separable Hilbert space.
\end{prop}
\begin{pf}
We are in the same situation as in the case of the group $\Aff{\Bbb R}$ 
of real affine transformations of the real straight line ${\Bbb R}$.
First of all we must prove two lemmas:
\begin{lem}
For each $\varphi\in L^1(\widetilde{\Aff{\Bbb C}})$, the conditions what 
follow are equivalent:
\begin{itemize}
\item $\varphi\in I = \bigcap_{\lambda\in{\Bbb C}} \ker U_\lambda.$
\item $\iint_{\Bbb C} \varphi(z,w) dw = 0, \mbox{ a. e. w. w. r. t. 
measure } dz,$
 where $(z,w) \in \widetilde{\Aff{\Bbb C}}$.
\end{itemize}
\end{lem}
Indeed, consider the function 
$$\psi(z) = \int_{\Bbb C} \varphi(z,w) dw.$$ Then,
$$\begin{array}{rl}
\varphi\in I &\Leftrightarrow \iint_{{\Bbb C}^2 }
e^{i\Re(z\overline{\lambda})} \varphi(z,w)dzdw = 0\\
 	& \Leftrightarrow \int_{\Bbb C} 
e^{i\Re(z\overline{\lambda})}(\int_{\Bbb C} \varphi(z,w)dw)dz=0\\
	& \Leftrightarrow \int_{\Bbb C} e^{i\Re(z\overline{\lambda})}\psi(z)dz 
= 0,\forall \lambda\\
	& \Leftrightarrow \widehat{\psi(\lambda)} = 0,\forall \lambda
\end{array}$$
where $\widehat{\psi}$ is the Fourier-Laplace image of $\psi$.
Because the Fourier-Laplace transformation is exact, $\psi(z) = 0$ in 
$L^1({\Bbb C},dz)$, i.e. $\psi$ vanishes almost everywhere. The lemma is 
therefore proved.

\begin{lem}
For each degenerate $\varphi(z,w) := \psi(z)\chi(w)$ in 
$L^1(\widetilde{\Aff{\Bbb C}})$, where 
$$\psi\in L^1({\Bbb C},dz) \cap L^2({\Bbb C},dz), \Vert \psi\Vert_{L^2} \ne 
0,$$ $$\chi\in L^1({\Bbb C},dw), \int_{\Bbb C} \chi(w) dw = 0,$$
the operator $T_\alpha(\varphi), \alpha\in {\Bbb S}^1$ are compact.
\end{lem} 
Really, for each function $f\in L^2({\Bbb R} \times {\Bbb S}^1,dx)$, we have
$$\begin{array}{rl}
(T_\alpha(\varphi)f)(x) &= \iint e^{i(\Re(we^x) + 
2\pi\alpha[{\Im(x+z)\over 2\pi}])}f(x\oplus z) \varphi(z,w) dzdw \\
	&= \int_{\Bbb C} e^{i\Re(we^x)}\chi(w)(\int_{\Bbb C} e^{i2\pi 
\alpha[{\Im(x+z)\over 2\pi}]})dw\\
	&= \widetilde{\chi}(e^x)\int_{\Bbb C} e^{i2\pi \alpha[{\Im(x+z)\over 
2\pi}]}f(x\oplus z) \psi(z) dz,
\end{array}$$
where $\widetilde{\chi}$ is the Fourier-Laplace image of $\chi$. It is 
well-known that $\widetilde{\chi}$
vanishes at the infinity. Then $\widetilde{\chi(e^x)} \to 0$ when $\Re(x) 
\to +\infty$ and $$\widetilde{\chi(e^x)} \to \widetilde{\chi}(0)= 
\int_{\Bbb C} \chi(w) dw = 0,$$ when $\Re(x) \to -\infty$. 

This means that  for each $\varepsilon > 0$, there exists a number $N$ 
such that for each $x$ in the domain $\vert x\vert > N$ we have 
estimation $\vert \widetilde{\chi}(e^x)\vert < {\varepsilon \over A},$ where
$$A := \sqrt{\int_{\Bbb C} \vert \psi(z)\vert^2 dz} = \Vert 
\psi\Vert_{L^2} \ne 0.$$ Consider the continuous function $\theta_N(x)$ 
satisfying
$$\theta_N(x) := \cases 0 & \mbox{ if } \vert x\vert > N+1,\\
			1 & \mbox{ if } \vert x\vert \leq N. \endcases$$
We have $$\widetilde{\chi}(e^x)(1-\theta_N(x)) = 0, \forall x; \vert 
x\vert \leq N,$$
$$\vert \widetilde{\chi}(e^x)(1-\theta_N(x))\vert < {\varepsilon\over 
A},\forall 
x.$$ Hence, $$\sup\vert \widetilde{\chi}(e^x)(1-\theta_N(x))\vert \leq 
{\varepsilon\over A}.$$ Consider the operator
$$(A_N(\varphi)f)(x) := \theta_N(x)\widetilde{\chi}(e^x)\int_{\Bbb C} 
e^{2\pi i\alpha[{\Im(x+z)\over 2\pi}]}f(x\oplus z) \psi(z)dz.$$ We prove 
that $A_N$ converge in norm to $T_\alpha(\varphi)$, when $N \to \infty$. 
Because the subspace of continuous function with compact support is dense 
in $L^1\cap L^2$, we can choose a sequence of continuous functions with 
compact support $\{ \psi_n \}$, approximating $\psi$. Put $\varphi_n 
:= \psi_n . \chi$. Consider the operator $A_N(\varphi_n)$ as above. We 
have $$A_N(\varphi) \to T_\alpha(\varphi_n) \mbox{ when } N \to \infty.$$ 
It is not hard to see that it is a uniform convergence on $N$. Thus we 
can change the order of limits and have
$$\lim_{N\to\infty} \lim_{n\to\infty} A_N(\varphi_n) = 
\lim_{n\to\infty}\lim_{N\to\infty} A_N(\varphi_n).$$ This means that 
$$\lim_{N\to\infty} A_N(\varphi) = \lim_{n\to\infty}T_\alpha(\varphi_n) = 
T_\alpha(\varphi).$$ We need therefore only to prove that $A_N(\varphi) 
\to T_\alpha(\varphi)$, $N \to \infty$ for continuous $\varphi$ with 
compact support. In that case we have
$$\begin{array}{rl}
\Vert (T_\alpha(\varphi)-A_N)f\Vert^2_{L^2} &\leq 
\sup\vert\widetilde{\chi}(1-\theta_N)\vert^2.\int_{{\Bbb R}\times{\Bbb 
S}^1} \vert e^{2\pi i\alpha[.]}f(z)\psi(z)dz\vert^2 dx\\
	& \leq \int_{\supp\; \psi} dz. 
\sup_x\vert\widetilde{\chi}(1-\theta_N)\vert^2 \int_{{\Bbb R}\times{\Bbb 
S}^1}\vert f\vert^2 \vert\psi\vert^2 dzdx
\end{array}$$
Remark that it is easy to see that
$$\vert\int hdt\vert^2 \leq \int_{\supp h} dt. \int\vert h\vert^2 dt.$$
Put $C := \int_{\supp h}$, we have
$$\begin{array}{rl}
\Vert(T_\alpha(\varphi)-A_N)f\Vert_{L^2}^2 &\leq C{\varepsilon^2\over A^2} 
\int_{\Bbb C}\vert\psi\vert^2(\int_{{\Bbb R}\times{\Bbb 
S}^1}\vert f\vert^2dx)dz\\	&\leq C{\varepsilon^2\over A^2} \int_{{\Bbb 
R}_+\times{\Bbb S}^1} \vert f\vert^2 dx \int\vert\psi\vert^2 dz\\
	& \leq C\varepsilon^2 \Vert f \Vert^2.
\end{array}$$
Thus we have $$\Vert T_\alpha(\varphi) - A_N)f\Vert \leq \sqrt{C} 
\varepsilon\Vert f\Vert$$ and therefore
$$\Vert T_\alpha(\varphi) - A_N\Vert \leq \sqrt{C} \varepsilon.$$
Choose a sequence $\varepsilon_n \to 0$, we have $$\Vert T_\alpha(\varphi) 
- A_{N_n}\Vert \to 0.$$

It is enough therefore to prove that the operators $A_N$ are compact. 
In fact, we have
$$(A_Nf)(x) = \theta_N(x) \widetilde{\chi}(e^x) \int_{\Bbb C} e^{2\pi 
i\alpha[{\Im(x+z)\over 2\pi}]}f(x\oplus z) \psi(z) dz.$$ Put $\xi :=
\Re(x+z)$ and $\eta := \Im(x+z)$, we have
$$\begin{array}{rl}
(A_Nf)(x) &= \theta_N(x)\widetilde{\chi}(e^x)\int_{\Bbb C}e^{2\pi 
i\alpha[{\eta\over 2\pi}]}f(\xi +2\pi i\{{\eta\over 
2\pi}\})\psi(\xi+i\eta-x)d\xi d\eta\\	&= \int K(x,\xi,\eta)f(\xi + 2\pi 
i\{{\eta\over 2\pi}\})d\xi d\eta ,
\end{array}$$
where $$K(x,\xi,eta) := \theta_N(x) \widetilde{\chi}(e^x) e^{2\pi 
i[{\eta\over 2\pi}]}\psi(\xi+i\eta-x).$$
We show now that 
$$\iint_{{\Bbb C}^2} \vert K(x,\xi,\eta)\vert^2 dxd\xi d\eta < +\infty.$$
In fact, 
$$\iint_{{\Bbb C}^2} \vert K\vert^2 d\xi d\eta = \int_{\Bbb C} \vert 
\theta_N\widetilde{\chi}\vert^2 \vert \psi\vert^2 d\xi d\eta = 
\vert\theta_N\widetilde{\chi} \vert^2 \Vert \psi\Vert^2.$$
Recall that $\theta_N(x) = 0$ for all $x$, $\vert x\vert \geq N+1$. Thus 
$\vert\theta_N\widetilde{\chi} \vert^2$ is a continuous function with 
compact support. Thus we have
$$\iint_{{\Bbb C}^2} \vert K(x,\xi,\eta)\vert^2 d\xi d\eta dx = 
\int_{\Bbb C} \vert \theta_N \widetilde{\chi}\vert^2 dx 
\Vert\psi\Vert^2_{L^2} < +\infty,$$
what means that $A_N$ is a compact operator. The lemma is proved

Now we are going to complete the proof of our proposition. Consider an
element $\varphi\in L^1({\Bbb C}\ltimes dzdw)$ and $\varphi\in I$, we
shall prove that $T_\alpha(\varphi)$ is a compact operator.  Because
$\varphi$ is of class $L^1$ then it can be approximated by finite linear
combinations of functions on separate variables $$\varphi_n(z,w) =
\sum_{k=1}^{N_n} \psi_k(z)\chi_k^*(w),$$ where $\psi_k\in L^1({\Bbb C},
dz)$, $\chi_k^*\in L^1({\Bbb C},dw)$.  Because $L^1 \cap L^2$ is dense in
$L^1$, we can assume that $\psi_k$ are in $L^1 \cap L^2$. Remove, when
necessary, the summands with $\Vert \psi_k\Vert_{L^2} = 0$, we can assume
that $\psi, \chi^*$ satisfy all the conditions but perhaps the condition
$$\alpha_k := \int_{\Bbb C} \chi^*(w)d \ne 0.$$ Following Lemma 3.2,
$$\varphi\in I \Longleftrightarrow \int_{\Bbb C} \varphi(z,w)dw = 0,
\mbox{ a. e. w. r. t. } dz.$$ We have the following a.e. inegalities
$$\begin{array}{rl} \vert\sum_{k=1}^{N_n} \alpha_k\psi_k \vert &=
\vert\sum_{k=1}^{N_n} \psi_k \int_{\Bbb C} \chi^*(w)dw - \int_{\Bbb C}
\varphi(z,w)dw \vert\\
	&\leq \int_{\Bbb C} \vert \sum_{k=1}^{N_n} \psi(k)\chi^*(w)
-\varphi(z,w) \vert dw \end{array}$$ Integrating both parts on the
variable $z$, we have $$\int_{\CC} \vert \sum_{k=1}^{N_n}
\alpha_k\psi_k\vert dz \leq \iint_{{\CC}^2} \vert\sum_{k=1}^{N_n}
\psi_k(z)\chi^*_k(w) - \varphi(z,w) \vert dzdw = \Vert \varphi_n^* -
\varphi\Vert_{L^1} \to 0 (n\to \infty).$$ 
Choose a function $\tau(w) \in
L^1({\CC}, dw )$ such that $$\tau(w) \geq 0, \forall w, \quad \int_{\CC}
\tau(w) dw = 1.$$ Consider the new approximation $$\varphi_n :=
\sum_{k=1}^{N_n} \psi_k (\chi^* - \alpha_k\tau).$$ We have $$\vert
\varphi_n -\varphi\vert \leq \vert \sum \psi_k\chi^*_k - \varphi\vert +
\vert \sum\alpha_k\psi_k\tau\vert = \vert \varphi^*- \varphi\vert +
\vert\sum\alpha_k\psi_k\vert \tau$$ and $$\begin{array}{rl} \Vert
\varphi_n - \varphi\Vert_{L^1} &\leq \iint_{{\CC}^2}\vert \varphi^*
-\varphi\vert dzdw + \int_{\CC} \vert\sum\alpha_k\psi_k \vert dz \int_{\CC}
\tau(w) dw\\
	&= \Vert\varphi^*_n - \varphi\Vert_{L^1} + 
\int\vert\sum\alpha_k\psi_k\vert dz.
\end{array}$$ 
Hence, $$\Vert \varphi - \varphi_n \Vert \to 0 \quad (n\to \infty).$$
\end{pf}

For each element $\varphi\in C^*(\widetilde{\Aff\RR})$ we have a natural 
map $$\Phi(\varphi) : {\Bbb S}^1 \to B(H),$$
$$\Phi(\varphi)(\alpha = T_\alpha(\varphi).$$

\begin{lem}
The map $\Phi$ is just a *-homomorphism from 
$C^*(\widetilde{\Aff(\CC)})^\sim$ into the C*-algebra $C({\SS}^1, B(H))$ of 
continuous functions on ${\SS}$ with values in the algebra $B(H)$ of 
bounded operators in a separable Hilbert space $H$.
\end{lem}
\begin{pf}
Because the functions with compact support are dense in 
$C^*(\widetilde{\Aff\CC})$ it is enough to prove the lemma for the case 
when $\varphi$ has compact support. Denote $c= mes(\supp \varphi)$. We 
have $$\begin{array}{rl}
\Vert \Phi(\varphi)(\alpha) - \Phi(\varphi)(\beta) \Vert &= \sup_{\Vert 
f\Vert \leq 1} \Vert (\Phi(\varphi)(\alpha) - \Phi(\varphi)(\beta) )f 
\Vert_{L^2} \\	 
\Vert(X (\Phi(\varphi)(\alpha) - \Phi(\varphi)(\beta))f \Vert^2_{L^2} &= 
\Vert (T_\alpha(\varphi) - T_{\beta}(\varphi))f\Vert^2_{L^2}= 
\end{array}$$
$$\begin{array}{rl} 
	&= \int_X \vert \iint_{{\CC}^2} e^{i\Re(we^x)} f(x\oplus z) 
\varphi(z,w) (e^{2\pi i\alpha[{\Im(x+z)\over 2\pi}]}- e^{2\pi 
i\beta}[{\Im(x+z)\over 2\pi}])dzdw \vert^2 dx\\
	&\leq \int_X(\iint_{\supp\varphi} \vert \varphi\vert^2 dzdw 
\iint_{\supp\varphi}\vert f\vert^2 \vert e^{2\pi 
i\alpha[{\Im(x+z)\over 2\pi}]}- e^{2\pi i\beta[{\Im(x+z)\over 
2\pi}]}\vert^2dzdw)dx \\
	&= \iint_{\supp\varphi} \vert \varphi\vert^2 dzdw . 
\int_X(\iint_{\supp\varphi} \vert f\vert^2 \vert e^{2\pi 
i\alpha[{\Im(x+z)\over 2\pi}]} - e^{2\pi i\beta[{\Im(x+z)\over 2\pi}]} 
\vert^2 dzdw)dx\\
	&\leq \iint_{\supp\varphi} \vert \varphi \vert^2 dzdw . \int_X 
\iint_{\supp\varphi} \vert f\vert^2 dzdw dx\times \\
	&\times \sup_{{\Im(x) \in 
{\SS}^1}\atop{\Im(z)\in \supp\varphi}}\vert e^{2\pi 
i\alpha[{\Im(x+Z)\over 2\pi}]}- e^{2\pi i\beta[{\Im(x+z)\over 2\pi}]}\vert^2.
\end{array}$$

Changing variable, we have
$$\int_X\iint_{\supp\varphi} \vert f(x\oplus z)\vert^2 dzdwdx = 
\iint_{\supp\varphi} (\int_X\vert f(x) \vert^2 dy) dzdw = c\Vert f 
\Vert^2_{L^2}.$$
Hence, we have
$$ (\Phi(\varphi)(\alpha) - \Phi(\varphi)(\beta))f\Vert^2_{L^2} \leq C. 
\iint_{\sup\varphi} \vert \varphi \vert^2 dwdz . \Vert f\Vert^2_{L^2} 
\sup_{{\Im(x)\in {\SS}^1} \atop {\Im(z)\in \supp\varphi}}\vert e^{2\pi i[.]}-
e^{2\pi i\beta[.]}\vert,$$
$$ \Vert \Phi(\varphi)(\alpha) - \Phi(\varphi)(\beta) \Vert^2_{L^2} \leq
C. \iint_{\supp\varphi} \vert \varphi\vert^2 dwdz . 
\sup_{{\Im(x)\in {\SS}^1}\atop {\Im(z)\in \supp\varphi}}\vert e^{2\pi
i\alpha[.]} - e^{2\pi i\beta[.]} \vert^2.$$
Because $\supp\varphi$ and ${\SS}^1$ are compact, we have
$$\sup_{{\Im(x)\in {\SS}^1}\atop {\Im(z)\in \supp\varphi}} \vert e^{2\pi
i\alpha[.]} - e^{2\pi i\beta[.]}\vert \to 0 \mbox{ when } \vert \alpha - 
\beta \vert \to 0.$$
This means that 
$$\Vert \Phi(\varphi)(\alpha) - \Phi(\varphi)(\beta) \Vert \to 0 \mbox{ 
when } \vert \alpha - \beta \vert \to 0.$$
\end{pf}

\begin{prop}
$$\Phi : I= \bigcap_{\lambda\in \CC} \ker U_\lambda \to C({\Bbb S}^1, {\cal 
K}(H))$$ is a *-isomorphism of C*-algebras.
\end{prop}
\begin{pf}
{\it a. Injectivity.}
If $\Phi(\varphi) = 0$ in $C({\Bbb S}^1, B(H))$. This means that 
$\Phi(\varphi)(\alpha) = T_\alpha(\varphi) = 0, \forall \alpha.$ Thus 
$$T_\alpha(\varphi)f \equiv 0, \forall \alpha \mbox { in } L^2(X,dx),$$ 
where $X \cong {\Bbb R} \times {\Bbb S}^1 = \{ x = t + i\theta ; t\in 
{\Bbb R}, \theta\in {\Bbb S}^1\}.$ 
Choose an orthonormal basis  $\{f_n \}$ of form
$$f_n(t,\theta) = g_n(t).h_n(\theta),$$ where  $\{ g_n\}$ and $\{ h_n\}$ 
are the corresponding orthonormal basis in $L^2({\Bbb R})$ and $L^2({\Bbb 
S}^1)$, respectively. We have therefore
$$\begin{array}{rl}
f_n(x\otimes z) &= g_n(\Re(x+z)) h_n(2\pi\{{\Im(x+z)\over 2\pi }\})\\
                &= g_n(t+z_1)h_n(2\pi\{{\theta+z_2\over 2\pi}\}),
\end{array}$$
for $z = z_1 +iz_2$. From the condition
$$T_\alpha(\varphi)f_n)(x) = \iint_{{\Bbb C}^2} e^{i(\Re(we^x) + 2\pi
\alpha[{\theta+z_2\over 2\pi}])}\varphi(z,w)f_n(x\oplus z)dzdw = 0,$$
 for a.e. w.r.t.  $dx.$  This means that
$$\int_{\Bbb R} g_n(t+z_1) (\iiint_{{\Bbb R}^3} e{i\Re(we^x)}
\varphi(z_1,w)h_n(2\pi\{{z_2\over 2\pi}\})e^{2\pi\alpha i[{z_2 +
\theta\over 2\pi}]})dz_1 = 0,$$ for a.e. w.r.t. $ dx$.
From here one deduces that
$$\int_{\Bbb R^1}(\iint_{\Bbb C^1} e^{\i\Re(we^x)}\varphi(z_1 , w)dw)
h_n(2\pi\{{z_2 + \theta\over 2\pi }\})e^{2\pi i\alpha[{z_2+\theta\over
2\pi}]}dz_2 = 0, \mbox{ a.e. w.r.t. } dz_1.$$
Put $z_2 + \theta = \xi$, we have 
$$\int_{\Bbb R} \widetilde{\varphi}(z_1e^x) h_n(2\pi\{{\xi\over 2\pi}\}) 
e^{2\pi i\alpha[{\xi\over 2\pi}]}d_\xi = 0, $$ a.e. w.r.t. $ dw.$ 
This means that 
$$\sum_{- \infty}^{+\infty} a_k e^{2\pi i\alpha.k} = 0, \forall 
\alpha\in{\Bbb S}^1,$$ where $$a_k := \int_{2\pi k}^{2\pi (k+1)} 
\widetilde{\varphi}(z, e^x) h_n(2\pi\{{\xi\over 2\pi}\})).$$
Because the Fourier-Laplace transformation is exact, we conclude that $a_k 
= 0$ for all $k$, i.e.
$$\int_{2\pi k}^{2\pi(k+1)} \widetilde{\varphi}(z,e^x) 
h_n(2\pi\{{\xi\over 2\pi }\})d\xi = 0, \forall k.$$ Because $h_n$ is a 
orthonormal basis in $L^2({\Bbb S}^1)$, we conclude that 
$$\widetilde{\varphi}(z,e^x) = 0,\mbox{ a.e. w.r.t. }dz_1d\xi = 
dz.$$ This means that 
$$\iint_{\Bbb C} e^{i\Re(we^x)} \varphi(z,w) dw = 0, \mbox{ a.e. w.r.t. } 
dz_1d\xi,$$ for all $x\in M$. 
For a sequence of points $\{x_n\}$, such that $\Re(x_n) \to -\infty
(m\to \infty)$ we have
$$\iint_{\Bbb C} \varphi(z,w)dw = \widetilde{\varphi}(z,0) = \lim_{m\to 
\infty} \iint_{\Bbb C} e^{i\Re(we^{x_n})} \varphi(z,w) dw = 0.$$ Following 
lemma 3.1, this is equivalent to the condition
$$U_\lambda(\varphi) = 0,\forall \lambda$$.

{\it b. Surjectivity.} It is easy to see that $\Phi(I)$ is a 
C*-subalgebra in $C({\Bbb S}^1, {\cal K}(H))$. But $T_\alpha(I) \subset 
{\cal K}(H)$, for all $\alpha\in {\Bbb S}^1$. Following J. Dixmier 
\cite{dixmier} for every $\alpha_1 \ne \alpha_2$ in $\Bbb S^1$ and for 
every $\xi_1, \xi_2$ in ${\cal K}(H)$, there exists $\varphi \in I$ such 
that $T_{\alpha_1}(\varphi) = \xi_1$ and $T_{\alpha_2}(\varphi) = \xi_2$. 
This means also that $\Phi(I) = C({\Bbb S}^1, {\cal K}(H)).$
\end{pf}

\begin{thm}
The C*-algebra $C^*(\widetilde{\Aff{\Bbb C}})^\sim$ can be included in a 
short exact sequence of C*-algebras and *-homomorphisms
$$0 \to C({\Bbb S}^1, {\cal K}(H)) \to C^*(\widetilde{\Aff{\Bbb C}}) 
\to C_0({\Bbb R}^2) \to 0.$$
\end{thm}
\begin{pf}
We remark first that {\it the quotient C*-algebra 
$C^*(\widetilde{\Aff{\CC}})/I$ is commutative, where by definition, 
$I = \bigcap_{\lambda\in{\CC}} \ker U_\lambda.$} 
Indeed, for every $\varphi_1$, $\varphi_2$ in 
$C^*(\widetilde{\Aff{\CC}})$, we have
$$U_\lambda(\varphi_1 * \varphi_2 - \varphi_2 * \varphi_1) = 
U_\lambda(\varphi_1) U_\lambda( \varphi_2) - U\lambda( 
\varphi_2)U_\lambda( \varphi_1) \equiv 0, \forall \lambda .$$ This means 
that $$\varphi_1 * \varphi_2 - \varphi_2 * \varphi_1 \in I.$$ Following 
J. Dixmier \cite{dixmier}{3.2.1} we see that the dual object of the 
quotient algebra consists of
$$\widehat(C^*(\widetilde{\Aff{\CC}}) /I) = \{ U_\lambda ; \lambda 
\in {\CC} \}\approx {\RR}^2.$$ Thus the Fourier-Gel'fand transformation 
gives us an isomorphism of C*-algebras 
$$C^*(\widetilde{\Aff{\CC}})/I \cong C_0({\RR}^2).$$ We have therefore an 
exact sequence 
$$\CD 0 @>>> C({\SS}^1,{\cal K}(H)) @>>> C^*(\widetilde{\Aff{\CC}}) @>>> 
C_0({\RR}^2) @>>> 0. \endCD$$

 \end{pf}

\begin{thm}
The structure of the C*-algebra $C^*(\widetilde{\Aff{\Bbb C}})^\sim$ is 
, up to isomorphic class, exactly defined by the invariant $index 
C^*(\widetilde{\Aff{\Bbb C}})^\sim = 1$
in the Kasparov group $$\Ext(C_0({\Bbb R}^2), C({\Bbb S}^1)) \cong {\Bbb 
Z}.$$ 
\end{thm}
\begin{pf}
It is easy to compute the K-groups as following
$$\begin{array}{rl}
K_0C_0({\RR}^2) &= \ker(K_0({\SS}^2 \to {\ZZ}) = \ker(K^0({\SS}^2) \to 
{\ZZ})\\
	&= \ker({\ZZ} \oplus {\ZZ} \to {\ZZ}) = {\ZZ},\\
K_0C({\SS}^1) &= K^0({\SS}^1) = {\ZZ}\\
K_1C_0({\RR}^2) &= K_0((C_0({\RR}^2 \otimes C_0({\RR})) = \ker( 
K_0C({\SS}^3) \to {\ZZ})\\
	&= \ker({\ZZ} \to {\ZZ}) = 0\\K_1C({\SS}^1) &= K_0(C({\SS}^1 \otimes 
C_0({\RR})) = \ker (K_0(C({\SS}^2) \to {\ZZ})\\
	&= \ker({\ZZ} \oplus {\ZZ} \to {\ZZ}) = {\ZZ}.
\end{array}$$
Following the universal coefficient formula, we have
$$\begin{array}{rl}
\Ext(C_0({\RR}^2,C({\SS}^1)) &\cong \Hom(K_0C_0({\RR}^2,K_1C({\SS}^1)) 
\oplus \Hom(K_1C_0({\RR}^2),K_0C({\SS}^1)) \\
	&\cong \Hom({\ZZ},{\ZZ} \cong {\ZZ}.
\end{array}$$
We have therefore a 6-terms exact sequence
$$\CD {\ZZ} @>>> K_0C^*(\widetilde{\Aff{\CC}}) @>>> {\ZZ} \\
@AA\delta_1 A @.                                   @. @VV\delta_0V\\
0 @<<< K_1C^*(\widetilde{\Aff{\CC}}) @<<< {\ZZ} \endCD$$
This means that $\delta_1$ must be $\delta_1=0$. Following the theorem of 
A. Connes \cite{connes3}, we have $K_1C^*(\widetilde{\Aff{\CC}}) = 0$ 
Thus $\delta_0$ is an epimorphism and hence is an isomorphism. This means 
that the element $Index C^*(\widetilde{\Aff{\CC}}) = 1$ in the KK-group 
$$\Ext(C_0({\RR}^2), C({\SS}^1)) \cong {\ZZ}.$$  
\end{pf}

\section{Bibliographical Remarks}
The results of sections 1 and 2 were firstly created by the author of 
this book \cite{diep1}, by using the BDF homological K-functor $\Ext$. It 
was then easily generated by J. Rosenberg 
\cite{rosenberg1}in the cases of the complex affine transformations group of 
the 
complex straight line and of  totally discrete local fields. Nevertheless
the case of the universal covering $\widetilde{\Aff{\Bbb C}}$ of the 
group of complex affine transformations of the complex straight line was 
difficult  and requested essentially to use KK-theory. It was done in the 
Ph.D. dissertation of Ho Huu Viet\cite{sonviet}, who was a Ph. D. student of the
author.

\chapter {Classification of $\MD_4$-Groups} 
%
\section{Real Diamond Group and Semi-direct Products ${\RR}\ltimes {\Bbb H}_3$}
Let us consider a connected and simply connected Lie group $G$, and its 
Lie algebra ${\frakt g} := \Lie G$. The action $\Ad : G \to \Aut{\frakt g}$ 
is defined as $${d\over dt}|_{t=0} A(\exp(tX)) : {\frakt g} = T_eG \to 
{\frakt g}= T_eG,$$ defined by the formula 
$$\Ad_X(Y) := ({d \over dt}|_{t=0} A(\exp(tX)))Y,$$ which is computed 
certainly as $${d\over ds}|_{s=0} \exp(tX)\exp(sY)\exp(-tX) \in 
{\frakt g}=  T_eG, $$ 
where $c(s):= \exp(sY)$ is just the curve passing through $e$ and with 
tangent $Y$ at $e$, i.e. the unique local solution of the system
$$\left\{ \begin{array}{rl}
\dot c(s) &= Y(c(s)),\\
c(0) &= e 
\end{array}\right.$$
We define the dual action $K:= \Ad^*(\inv)$ of $G$ on the dual 
vector space ${\frakt g}^*$ of its Lie algebra ${\frakt g}$,
$$\langle K(g)F, X\rangle := \langle F, \Ad(g^{-1}) \rangle ,$$ for all 
$X\in {\frakt g}$ and $F \in {\frakt g}^*$, where $\langle .,.\rangle$ is 
the pairing between ${\frakt g}$ and ${\frakt g}^*$. It is easy to see 
that \begin{itemize}
\item The stabilizer $G_F$ of a fixed point $F \in \Omega_F := G.F$ is a 
closed subgroup and the tangent space to the orbit at this point can be 
identified with the quotient-space ${\frakt g}/{\frakt g}_F$, where 
${\frakt g}_F := \Lie G_F$ is the Lie algebra of the stabilizer $G_F$
\item The Kirillov form $\langle F,[X,Y]\rangle$ define a non-degenerate 
bilinear form $B_F$ on the tangent space $T_F\Omega$
\item The right translations define an action of $G$ on the orbit 
$\Omega_F$ and transpose the Kirillov bilinear form $B_F$ providing a 
right-invariant differential form $\omega_F$ on the orbit $\Omega_F$ 
passing through $F$.
\item And finally, each K-orbit is of even dimension.
\end{itemize} 

\begin{defn}
We say that a solvable Lie group $G$ belongs to the class $\MD$ (resp., 
$\overline{\MD}$) iff every its K-orbit has dimension 0 or maximal 
(resp., equal to its dimension $\dim G$).
A Lie algebra is of class $\MD$ (resp., $\overline{\MD}$) iff its 
corresponding Lie group is of the same class. 
\end{defn}

\begin{prop}
If ${\frakt g}$ is a $\MD$-algebra. Then the second derived ideal 
$[[{\frakt g},{\frakt g}],[{\frakt g},{\frakt g}]]$ is commutative.
\end{prop}
\begin{pf}
This proposition was proved in chapter 2.
\end{pf}

\begin{exam}
The so called {\it real diamond Lie algebra} \index{Lie algebra! real 
diamond - }  is the 4-dimensional 
solvable Lie algebra ${\frakt g}$ with a basis $X,Y,Z,T$ satisfying the 
following commutation relations
$$[X,Y] = Z, [T,X] = -X, [T,Y] = Y,$$
$$[Z,X] = [Z,Y] = [T,Z] = 0.$$  
\end{exam}
This example shows that our real diamond Lie algebra ${\RR}\ltimes 
{\frakt h}_3$ is 
just the extension of the 1-dimensional Lie algebra ${\RR} T$ by the 
Heisenberg algebra ${\frakt h}_3$ with basis $X,Y,Z$, where the action of 
$T$ on ${\frakt h}_3$ is defined by the matrix
$$\pmatrix 
-1 & 0 & 0\\
 0 & 1 & 0\\
 0 & 0 & 0 
\endpmatrix$$
It is easy to verify that this real diamond Lie algebra belongs to the 
class $\MD_4$.

\begin{exam}
Let us consider the Lie algebra ${\frakt g} := {\RR}\ltimes_J{\frakt h}_3$ 
with a basis $X,Y,Z,T$ and 
the action of Lie algebra ${\RR}T$ on the 3-dimensional Heisenberg ideal 
${\frakt h}_3$ defined by the matrix 
$$
\ad_T := \pmatrix  0  & -1 & 0 \\
		   1  &  0 & 0 \\
		   0  &  0 & 0 \\
	\endpmatrix$$ 
It is easy to verify that it belongs to the class $\MD_4$

\end{exam}
In the next section, we shall give a classification of all 
Lie $\MD_4$-algebras.

\section{Classification Theorem}

\begin{prop}
Assume ${\frakt g}$ is a $\MD_4$-algebra with generators $T,X,Y,Z$. 
Denote by ${\frakt g}^1 = [{\frakt g},{\frakt g}]$ the first derived 
ideal and ${\RR}^n$ the commutative subalgebra in ${\frakt g}$.
\begin{enumerate}
\item[I] If ${\frakt g}$ is decomposable, i.e. cam be written as a
direct product of two ideals, then $${\frakt g} = {\RR}^n \oplus 
\widetilde{\frakt g},$$ for some $n = 1,2,3,4$ and some indecomposable 
ideal $\widetilde{\frakt g}$.
\item[II] If ${\frakt g}$ is indecomposable, it is is isomorphic to one 
of the Lie algebras what follow:
{\begin{enumerate}
\item[1] 
${\frakt g}^1 = {\RR}Z, [X,Y] = aZ, [X,Z] = xZ, [X,T] =cZ$
$$[Y,Z] = yZ, [Y,T] = bZ, [Z,T] = tZ,$$ for some $a,b,c,x,y,t \in {\RR}$ 
non-vanishing all together and $at+bx+cy = 0$ 
\item[2]
${\frakt g}^1 = {\RR}Y + {\RR}Z$, $[T,X] = 0$ and
        \begin{enumerate}
	\item[2.1] $\ad_T \in \Aut{\frakt g}^1 \cong \GL_2(\RR),$ $\ad_X = 
	\alpha \ad_T$, $(\alpha\in {\RR})$

	\item[2.2] ${\frakt g} = \aff{\CC} = \Lie(\Aff{\CC})$, i.e. $\ad_X 
	= \pmatrix 0 & 1\\  -1 & 0 \endpmatrix$ and $\ad_T = \pmatrix 1 & 0 \\
	0 & 1 \endpmatrix$ 
	\end{enumerate}
\item[3]
${\frakt g}^1 = {\RR}X + {\RR}Y + {\RR}Z\cong {\RR}^3$, $\ad_T \in 
\Aut_{\RR}{\frakt g}^1 \cong \GL_3(\RR)$
\item[4]
${\frakt g}^1 = {\RR}X + {\RR}Y + {\RR}Z \cong {\frakt h}_3$, the 
3-dimensional Heisenberg Lie algebra
and $$\ad_T = \pmatrix a_{11} & a_{12} & 0 \\
		       a_{21} & a_{22} & 0 \\ 
		       a_{31} & a_{32} & 0 ,\endpmatrix \in \End_{\RR}{\frakt 
g}^1, a_{11}^2 + a_{12}.a_{21} \ne 0 $$
 \end{enumerate}
}
\end{enumerate} 
\end{prop}
\begin{pf}
We proved that ${\frakt g}^2$ is commutative, say isomorphic to ${\RR}^n$, 
for $n=1,2,3,4$. 

\end{pf}

\begin{thm}[\bf Classification of $\MD_4$-algebras]
We keep the same notation as in the previous proposition.
\begin{enumerate}
\item[I] 
 If ${\frakt g}$ is decomposable, i.e. cam be written as a
direct product of two ideals, then $${\frakt g} = {\RR}^n \oplus 
\widetilde{\frakt g},$$ for some $n = 1,2,3,4$ and some indecomposable 
ideal $\widetilde{\frakt g}$.
\item[II]
Suppose that ${\frakt g}$ {\rm indecomposable }%
into a direct product of two proper ideals. Then ${\frakt g}$ is of
class $\MD_4$ if and only if it is generated by the generators $T,X,Y,Z$
with {\rm the only non-trivial }%
commutation relations in one
of the following cases :

	\begin{enumerate}
        \item[\bf 1.] ${\frakt g} = {\RR}Z \cong {\RR}$ and
		\begin{enumerate}
                \item[1.1] $$[T,X] = Z ,\eqno{ ({\frakt g}_{4,1,1})}$$
                \item[1.2] $$[T,Z] = Z , \eqno{({\frakt g}_{4,1.2})}$$
		\end{enumerate}
        \item[\bf 2.] ${\frakt g}^1 = {\RR}Y + {\RR}Z \cong {\RR}^2$  and
		\begin{enumerate}
                \item[2.1] $$[T,X] = \lambda X , [T,Y] = Y ; \lambda \in {\Bbb R}^* = {\Bbb R}
\setminus (0), \eqno{({\frakt g}_{4,2,1(\lambda)})}$$
                \item[2.2] $$[T,X] = X , [T,Y] = X + Y,\eqno{({\frakt 
g}_{4,2,2})}$$
                \item[2.3] $$\ad_T = \begin{pmatrix}\cos{\varphi} & 
\sin{\varphi} & 0 \cr
                 -\sin{\varphi} & \cos{\varphi} & 0 \cr
                  0            &       0        & 0 \cr
\end{pmatrix} , \eqno{({\frakt 
g}_{4,2,3(\varphi)})}$$
                \item[2.4] $$\ad_T = \begin{pmatrix} 1 & 0 & 0 \cr
                  0 & 1 & 0 \cr
                  0 & 0 & 0 \cr
\end{pmatrix} , \ad_X = \begin{pmatrix} 0 & 1 & 0 \cr
                   -1 & 0 & 0 \cr
                    0 & 0 & 0 \cr\end{pmatrix} ,\eqno{(
{\frakt g}_{4,2,4} = Lie(\Aff{\Bbb C}))}$$
		\end{enumerate}
        \item[\bf 3.] ${\frakt g}^1 = {\RR}X + {\RR}Y + {\RR}Z$, commutative 
derived
        ideal
		\begin{enumerate}
                \item[3.1]          $$\ad_T = \begin{pmatrix} \lambda_1 &    0 &      0 \cr
                  0 & \lambda_2 & 0 \cr
                  0 & 0 & 1 \cr\end{pmatrix} , \lambda_1,\lambda_2 \in {\Bbb R}^* 
, \eqno{({\frakt g}_{4,3,1(\lambda_1,\lambda_2)})}$$
                \item[3.2] $$\ad_T = \begin{pmatrix} \lambda & 1 & 0 \cr
                  0 & \lambda & 0 \cr
                  0 & 0 & 1       \cr\end{pmatrix} , \lambda \in {\Bbb R}^* ,
\eqno{({\frakt g}_{4,3,2(\lambda)})}$$
                \item[3.3] $$\ad_T = \begin{pmatrix} 1 & 1 & 0 \cr
                  0 & 1 & 1 \cr
                  0 & 0 & 1 \cr\end{pmatrix}, \eqno{({\frakt g}_{4,3,3})}$$
                \item[3.4] $$\ad_T = \begin{pmatrix} \cos{\varphi} & \sin{\varphi} & 0 \cr
                 -\sin{\varphi} & \cos{\varphi} & 0 \cr
                 0 & 0 & \lambda \cr\end{pmatrix} , \lambda \in {\Bbb R}^* , \varphi \in 
(0,\pi) , \eqno{({\frakt g}_{4,3,4(\lambda)})}$$
		\end{enumerate}
        \item[\bf 4.] ${\frakt g}^1 = {\RR}X + {\RR}Y + {\RR}Z \cong {\frakt 
h}_3$, the 3-dimensional Heisenberg Lie algebra
and
        	\begin{enumerate}
                \item[4.1] $$\ad_T = \begin{pmatrix} 0 & 1 & 0 \cr
                  -1 & 0 & 0 \cr
                  0 & 0 & 0 \cr\end{pmatrix}, [X,Y] = Z ,\eqno{({\frakt 
g}_{4,4,1} = Lie({\Bbb R} \ltimes_J {\Bbb H}_3))}$$
                \item[4.2] $$\ad_T = \begin{pmatrix} -1 & 0 & 0 \cr
                  0 & 1 & 0 \cr
                  0 & 0 & 0 \cr\end{pmatrix} , [X,Y] = Z , \eqno{({\frakt 
g}_{4,4,2} = Lie({\Bbb R} \ltimes {\Bbb H}_3))}$$
( in this case the group is called the real diamond group ).
		\end{enumerate}
	\end{enumerate}
\end{enumerate}
\end{thm}
\begin{pf}

I.

The part I of the theorem is proved in the previous proposition. Let us 
prove the second part. 

II.1. 

If the Lie algebra ${\frakt g}$ is in the case 
II.1 of Proposition 2.1 then ${\frakt g} = {\RR}X + {\RR}Y + {\RR}Z + 
{\RR}T$, $ {\frakt g}^1 = {\RR}Z \cong {\RR}$ and we have the 
commutation relations $$[X,Y] = aZ, \quad [Y,T] = bZ, \quad [T,X] = cZ,$$
$$[X,Z] = xZ, \quad [Y,Z] = yZ, \quad [T,Z] = tZ,$$ where $a,b,c,x,y$ are 
real numbers not vanishing all together, and $$at + bx + cy = 0.$$

Suppose that $t\ne 0$. Change the basis $\{ X,Y,Z,T\}$ to the new one 
$\{ X',Y',=Z',T'\}$, where 
$$\begin{array}{rl}
X' &:= X - {c\over t}Z - {x\over t}T\\
Y' &:= Y + {b\over t}Z - {y\over t}T\\
Z' &:= Z\\
T' &:= {1\over t}T,
\end{array}$$
then $[T',Z'] = Z', [T',X'] =[T',Y']= [X',Y'] = [X',Z'] = [Y',Z'] = 0$. 
This means that ${\frakt g} \cong {\frakt g}_{4,1,2}$. By a similar 
argument, one deduces that if $x\ne 0$ or $y\ne 0$ the Lie algebra 
${\frakt g}$ is of the same class ${\frakt g}_{4,1,2}$.

Suppose now that $x=y=t=0$ and $c \ne 0$. Changing the basis $X,Y,Z,T$ by 
the new one $X',Y',Z',T'$, 
$$\begin{array}{rl}
X' &:= X,\\
Y' &:= Y + {b\over c}X + {a\over c}T,\\
Z' &:= Z,\\
T' &:= {1\over c}T
\end{array}$$
we have 
$$[T',X'] = Z', [T',Y'] = [T',Z'] = [X',Y'] = [X',Z']= [Y',Z'] = 0,$$ 
i.e. ${\frakt g} \cong {\frakt g}_{4,1,1}$. By similar argument, we 
prove also that if $a\ne 0$ or $b\ne 0$ then ${\frakt g}\cong {\frakt 
g}_{4,1,1}$. 

II.2.  ${\frakt g}^1 \cong {\RR}Y + {\RR}Z$ and $[T,X] =0$:

If the Lie algebra is in the case II.2 but not in the case II.2.4 of the 
Lie algebra of affine transformations of the complex straight line 
${\frakt g}_{4,2,4} = \aff{\CC}$, defined by the only nontrivial commutation 
relations $$\ad_X = \pmatrix 0 & 1\\ -1 & 0 \endpmatrix \mbox{ and } 
\ad_T = \pmatrix 1 & 0\\ 0 & 1\endpmatrix $$ Then $\ad_T \in \End_{\RR} 
{\frakt g}^1 \cong \GL_2({\RR})$, and $\ad_X = \alpha \ad_T,$ for some 
$\alpha\in {\RR}$.

Do change the basis $X,Y,Z,T$ to the new basis $X'= X-\alpha T, 
Y'=Y,Z'=Z, T' = T$ we obtain $\ad_{X'}=0$. Therefore we can, and do,  assume 
from the beginning that $\ad_X=0$. The later classify depends therefore 
on the canonical Jordan form of $\ad_T \in \End({\frakt g}^1)\cong 
\Mat_2({\RR})$. We have therefore three sub-cases
\begin{enumerate}
\item[II.2.1] $\ad_T$ is a non-degenerate diagonalizable matrix 
equivalent to $\pmatrix \lambda_1 & 0\\ 0 & \lambda_2 \endpmatrix,$ for 
$\lambda_1,\lambda_2 \in {\RR}$ and $\lambda_1 .\lambda_2 \ne 0$, i.e. 
$\lambda_1, \lambda_2 \ne 0$. Change the basis $X,Y,Z,T$ to the basis 
$X,Y,Z,T' = {1\over\lambda_2}T$, we have $$\ad_{T'} = \pmatrix \lambda & 
0\\ 0 & 1\endpmatrix$$ where $\lambda = {\lambda_1\over \lambda_2}$. 
This means that ${\frakt g} \cong {\frakt g}_{4,2,1(\lambda)}$, for some 
$\lambda \in {\RR}^*$.  

\item[II.2.2]
The matrix $\ad_T$ can not be diagonalizable, but has only real 
eigenvalues. In this case, it has an eigenvalue $\lambda$ of multiplicity 2.
Change $T$ by $T' = {1\over \lambda}T$ we have
$$\ad_{T'} = \pmatrix 1 & 1 \\ 0 & 1 \endpmatrix.$$ This means that 
${\frakt g} \cong {\frakt g}_{4,2,2}$.

\item[II.2.3]
$\ad_T$ is diagonalizable but has complex eigenvalue. In this case it 
should has the second eigenvalue also complex and is conjugate to the 
first. This means that $\det(\ad_T)> 0$. Change $T$ by $T' = {1\over 
\sqrt{\det(\ad_T)}}T$, we have $\det\ad_{T'} =1$. In this case the real 
Jordan form of $\ad_{T'}$ is just $$\ad_{T'} = \pmatrix \cos \varphi & 
\sin\varphi \\ -\sin\varphi & \cos\varphi\endpmatrix.$$ This means that 
${\frakt g} \cong {\frakt g}_{4,2,3}$. 

\end{enumerate}

II.3.
Let us consider now the case when ${\frakt g} \cong {\RR}X + {\RR}Y + 
{\RR}Z + {\RR}T$ and ${\frakt g}^1 \cong {\RR}X + {\RR}Y + {\RR}Z \cong 
{\RR}^3$ is a commutative ideal and $\ad_T \in \Aut_{\RR}({\frakt g}^1) 
\cong \GL_3({\RR})$. We have also 4 cases for the Jordan form of the 
matrix $\ad_T$. 
\begin{enumerate}
\item[II.3.1]
$$\ad_T \cong \pmatrix \tilde{\lambda_1} & 0 & 0\\ 0 & \tilde{\lambda_2} 
& 0\\ 0 & 0 & \tilde{\lambda_3} \endpmatrix, \tilde{\lambda_1}, 
\tilde{\lambda_2},\tilde{ \lambda_3} \in {\RR}^*.$$

In this case, we change $T$ by $T' = {1\over \tilde{\lambda_3}}T$, we have
$$\ad_{T'} \cong \pmatrix \lambda_1 & 0 & 0\\ 0 & \lambda_2 & 0\\ 0 & 0 
& 1\endpmatrix ,$$ with $\lambda_1 := {\tilde{\lambda_1}\over 
\tilde{\lambda_3}}, \lambda_2 := {\tilde{\lambda_2}\over 
\tilde{\lambda_3}}.$ This means that ${\frakt g} \cong {\frakt 
g}_{4,3,1(\lambda_1,\lambda_2)}, \lambda_1, \lambda_2 \in {\RR}^*.$

\item[II.3.2]
$$\ad_T \cong \pmatrix \lambda_1 & 1 & 0\\ 0 &\lambda_1 & 0\\ 0 & 0
\lambda_2\endpmatrix, \lambda_1, \lambda_2 \in {\RR}^*$$

In this case, change $T$ by $T' = {1\over \lambda_2}T$, we have 
$$\ad_{T'} \cong \pmatrix \lambda & 1 & 0\\ 0 & \lambda & 0\\ 0 & 0 & 
1\endpmatrix,$$ with $\lambda = {\lambda_1 \over \lambda_2}$. This means 
that ${\frakt g} \cong {\frakt g}_{4,3,2(\lambda)}, \lambda\in {\RR}^*$.

\item[II.3.3]
$$\ad_T \cong \pmatrix \lambda & 1 & 0\\ 0 & \lambda & 1\\ 0 & 0 & \lambda
\endpmatrix, \lambda \in {\RR}^*$$

In this case, we change $T$ by $T' = {1\over \lambda}T$ and have 
$$\ad_{T'} \cong \pmatrix 1 & 1 & 0\\ 0 & 1 & 0\\ 0 & 0 & 1 \endpmatrix 
$$ this means that ${\frakt g}\cong {\frakt g}_{4,3,3}$.

\item[II.3.4]
$$\ad_T \cong \pmatrix a & b & 0\\ -b & a & 0 \\ 0 & 0 & \tilde{\lambda} 
\endpmatrix, \tilde{\lambda}\in {\RR}^*, b >0$$

In this case, we change $T$ by $T' = {1\over \sqrt{a^2 + b^2}}T$ and we 
have $$\ad_{T'} \cong \pmatrix \cos\varphi & \sin\varphi & 0\\ 
-\sin\varphi & \cos\varphi & 0 \\ 0 & 0 & \lambda \endpmatrix,$$ for 
$\lambda = {\tilde{\lambda}\over \sqrt{a^2 + b^2}}, \varphi = \arccos{a 
\over \sqrt{a^2+b^2}}\in (0,\pi).$ This means that ${\frakt g} \cong 
{\frakt g}_{4,3,4(\lambda,\varphi)}, \lambda\in {\RR}^*, \varphi\in (0,\pi)$.
\end{enumerate}

II.4. If ${\frakt g} \cong {\RR}X + {\RR}Y + {\RR}Z + {\RR}T,$ and 
${\frakt g}^1 \cong {\frakt h}_3$ with the action of $T$ as $$\ad_T = 
\pmatrix a_{11} & a_{12} & 0 \\ a_{21} & -a_{11} & 0 \\ a_{31} & a_{32} & 
0 \endpmatrix \in \End_{\RR}{\frakt g}^1 \cong \Mat_3({\RR}),$$ such that 
$$a_{11}^2 + a_{12}a_{21} \ne 0.$$

Let us denote by $a,b$ the solution of the Cramer system
$$\left\{ \begin{array}{rl} a_{11}x + a_{21}y &= a_{31} \\ a_{12}x - 
a_{11}y &= a_{32}\end{array}\right.$$
Because $a_{11}^2 + a_{12}a_{21} \ne 0$ this system has a unique 
solution. Change the basis $X,Y,Z,T$ to the basis $X',Y',Z',T'$
$$\begin{array}{rl}
X' &= X + aZ\\
Y' &= Y + bZ\\
Z' &= Z\\
T' &= T
\end{array}$$
then we also have the same commutation relation of the Heisenberg algebra 
${\frakt h}_3$ 
$$[X',Y'] = Z', [X',Z'] = [Y',Z'] = 0$$ and in the basis $X',Y',Z'$ of 
${\frakt g}^1$ $$\ad_T = \pmatrix a_{11} & a_{12} & 0\\ a_{21} & -a_{11} 
& 0\\ 0 & 0 & 0 \endpmatrix .$$ We can do this from the beginning and thus 
we have two cases:
\begin{enumerate}
\item[II.4.1]
$a_{11}^2 + a_{12}a_{21} < 0$.
This deduces also that $a_{12}a_{21} \ne 0$. Pose $\lambda = 
\sqrt{-(a_{11}^2 + a_{12}a_{21})}.$ Change the basis $X,Y,Z,T$ by the new 
basis $X',Y',Z',T'$ 
$$\begin{array}{rl}
X' &:= -\lambda X\\
Y' &:= a_{11}X + a_{21}Y\\
Z' &:= -a_{21}Z \\
T' &:= {1\over \lambda}T
\end{array}$$
It is easy to see that $${\RR}X' + {\RR}Y' + {\RR}Z' \cong {\frakt h}_3$$ 
and $\ad_{T'}$ acts on ${\frakt h}_3$ by the matrix 
$$\ad_{T'} = \pmatrix 0 & 1 & 0\\ -1 & 0 & 0\\ 0 & 0 & 0 \endpmatrix.$$ 
This means that in this case ${\frakt g} \cong {\frakt g}_{4,4,1}$.

\item[II.4.2]
$a_{11}^2 + a_{12}a_{21} > 0$. Put $\lambda := \sqrt{a_{11}^2 + 
a_{12}a_{21}}. $ 

If $\lambda \ne a_{11}$, we do the following base change
$$\begin{array}{rl}
X' &= (a_{11}-\lambda)X + a_{21}Y,\\
Y' &= -a_{12}X +(a_{11} - \lambda)Y,\\
Z' &= 2\lambda(\lambda-a_{11})Z,\\
T' &= {1\over \lambda}T.
\end{array}$$
It is easy to check that $[X',Y'] = Z', [X',Z'] = [Y',Z'] = 0$, ${\RR}X' 
+ {\RR}Y' + {\RR}Z' \cong {\frakt h}_3$ and the action of $T'$ on ${\frakt 
h}_3$ is defined by the matrix $$\ad_{T'} =\pmatrix -1 & 0 & 0\\ 0 & 1 & 0\\ 
0 & 
0 & 0\endpmatrix.$$ 
This means that ${\frakt g} \cong {\frakt g}_{4,4,2}$.

If $\lambda = a_{11}$ then we use the following basis change
$$\begin{array}{rl}
X' &= a_{12}X - 2a_{11}Y',\\
Y' &= 2a_{11} + a_{21}Y,\\
Z' &= 4a_{11}^2 Z,\\
T' &= {1\over \lambda}T.\end{array}$$ 
We have the same result as in the 
previous case. Thus ${\frakt g} \cong {\frakt g}_{4,4,2}$.
\end{enumerate}
\end{pf}

\begin{rem}
Following part I of Proposition 2.1 and Theorem 2.1, if ${\frakt g}$ is 
an $\MD_4$-algebra, then it should be ${\RR}^4$ or decomposed as 
${\frakt g} = {\RR}^n \oplus \widetilde{\frakt g}$, for $1 \leq n \leq 3$ 
and hence ${\frakt g}^* = {\RR}_n \oplus \widetilde{\frakt g}^*$ and the 
correspondence simply connected Lie group $G$ should be decomposed as 
${\RR}^n \times \widetilde{G}$. The K-action of $G$ in ${\frakt g}^*$ 
should be decomposed into the product of the trivial action of $G$ on 
${\RR}_n$ and the K-action of $\widetilde{G}$ on $\widetilde{\frakt 
g}^*$. Thus $\widetilde{G}$ should be a $\MD_{4-n}$-group ($1\leq n \leq 
3$). So that we consider only the $\MD_4$-algebra which are either ${\RR}^4$
or indecomposable.
\end{rem}

\begin{rem}
We shall denote the corresponding Lie $\MD_4$-groups by the capitals with 
the same indices as their Lie algebras, e.g. $G_{4,2,1},....$
\end{rem}

\section{Description of the Co-adjoint Orbits}

\subsection{Some remarks about the co-adjoint representation}
Recall first of all some background of co-adjoint representations. For a 
connected and simply connected Lie group $G$ with its Lie algebra 
${\frakt g}$, the formula of the co-adjoint action of $G$ on ${\frakt 
g}^*$ is given by 
$$\langle K(g)F, X \rangle = \langle F, \Ad(g^{-1})\rangle, $$ for each 
$F\in {\frakt g}^*, g \in G$ and each $X\in {\frakt g}$. The co-adjoint 
orbit of $G$ in ${\frakt g}^*$ passing through $F$ is denoted by 
$$\Omega_F := K(G)F := \{ K(g)F ; g\in G \}.$$ 

Recall that the exponential map 
$$\exp : {\frakt g} \to G$$ defines a local diffeomorphism of some 
neighborhood of $0\in {\frakt g}$ to some neighborhood of the identity 
element $e\in G$ which is commutes with the adjoint representations
$$\ad : {\frakt g} \to \End_{\RR}{\frakt g} $$ and
$$\Ad : G \to \Aut_{\RR}({\frakt g})$$ in the sense that the following 
diagram is commutative
$$\CD
 G @>\Ad >> \Aut_{\RR}({\frakt g})\\
@A \exp AA      @AA\exp A\\
{\frakt g} @>\ad >> \End_{\RR}({\frakt g})
\endCD$$

Recall the Lie group $G$ is called to be exponential, if and only if the 
exponential map $$\exp : {\frakt g} \to G$$ is a diffeomorphism.

\begin{prop}
Let $G$ be a connected and simply connected solvable (finite 
dimensional, real) Lie group with Lie algebra ${\frakt g}$. Then the 
following assertions are equivalent:
\begin{enumerate}
\item The exponential map $\exp {\frakt g} \to G$ is a global 
diffeomorphism.
\item
For all $X\in {\frakt g}$, the operators $\ad_X$ have no purely imaginary 
eigenvalues.
\end{enumerate}
\end{prop}
\begin{pf}
See {\sc M. Saito} \cite{saito} and {\sc N. Bourbaki} \cite{bourbaki1}
\end{pf}

\begin{cor}
All $\MD_4$-groups, except for the groups $G_{4,2,3(\pi/2)}$, $G_{4,2,4}$  
$= \widetilde{\Aff{\CC}}$, $G_{4,3,4(\lambda,\pi/2)}, (\lambda\in 
{\RR}^*)$ and $G_{4,4,1}$ are exponential.
\end{cor}
\begin{pf}
Consider a general element $U = aX + bY + cZ + dT \in {\frakt g}$. 
\begin{enumerate}
\item[1.1] If ${\frakt g} = {\frakt g}_{4,1,1}$ then in the basis $X,Y,Z,T$,
$$\ad_U = \pmatrix 0 & 0 & 0 & 0\\ 0 & 0 & 0 & 0 \\ d & 0 & 0 & -a\\ 0 & 
0 & 0 & 0 \endpmatrix$$ and $\ad_U$ has only real eigenvalue 0 of 
multiplicity 4.
\item[1.2] If ${\frakt g} = {\frakt g}_{4,1,2}$ then in the basis $X,Y,Z,T$,
$$\ad_U = \pmatrix 0 & 0 & 0 & 0\\ 0 & 0 & 0 & 0\\ 0 & 0 & d & -c\\ 0 & 0 
& 0 & 0 \endpmatrix$$ and hence $\ad_U$ has only real eigenvalues 0 of 
multiplicity 3 and $d$.
\item[2.1] If ${\frakt g} = {\frakt g}_{4,2,1(\lambda)}, (\lambda \in 
{\RR}^*)$, then in the basis $X,Y,Z,T$,
$$\ad_U = \pmatrix 0 & 0 & 0 & 0\\ 0 & d\lambda & 0 & -\lambda b\\ 0 & 0 & 
d & -c\\ 0 & 0 & 0 & 0 \endpmatrix$$ and hence $\ad_U$ has only the real 
eigenvalues 0 of multiplicity 2 and $d$ and $d\lambda$.
\item[2.2] If ${\frakt g} = {\frakt g}_{4,2,2}$, then in the basis 
$X,Y,Z,T$, $$\ad_U = \pmatrix 0 & 0 & 0 & 0\\ 0 & d & d & -(b+c)\\ 0 & 0 
& d & -c\\ 0 & 0 & 0 & 0 \endpmatrix$$ and then $\ad_U$ has only the real 
eigenvalues 0 of multiplicity 2 and d of multiplicity 2. 
\item[2.3] If ${\frakt g} = {\frakt g}_{4,2,3}$, then in the basis 
$X,Y,Z,T$, 
$$\ad_U = \pmatrix 0 & 0 & 0 & 0 \\ 0 & d\cos\varphi & d\sin\varphi & 
-b\cos\varphi - c\sin\varphi\\ 0 & -d\sin\varphi & d\cos\varphi & 
b\sin\varphi -\cos\varphi\\ 0 & 0 & 0 & 0\endpmatrix$$, and hence $\ad_U$ 
has only the real eigenvalues 0 of multiplicity 2 and $de^{\pm 
i\varphi}$. The last number should not be purely imaginary iff 
$\varphi\in (0,\pi)$ and $\varphi \ne \pi/2$.

\item[3.1] If ${\frakt g} = {\frakt 
g}_{4,3,1(\lambda_1,\lambda_2)}(\lambda_1,\lambda_2\in {\RR}^*)$ then in 
the basis $X,Y,Z,T$,
$$\ad_U = \pmatrix d\lambda_1 & d & 0 & -\lambda_1a\\ 0 & d\lambda_2 & 0 
& -b\lambda_2\\ 0 & 0 & d & -c\\ 0 & 0 & 0 & 0 \endpmatrix$$ and hence 
$\ad_U$ has only the real eigenvalues $0, d, d\lambda_1, d\lambda_2$. 

\item[3.2] If ${\frakt g} = {\frakt g}_{4,3,2(\lambda)}(\lambda\in 
{\RR}^*)$, then in the basis $X,Y,Z,T$, $$\ad_U = \pmatrix d\lambda & d & 
0 & -a\lambda -b\\ 0 & d\lambda & 0 & -b\lambda\\ 0 & 0 & d & -c \\ 0 & 0 
& 0 & 0 \endpmatrix$$ and hence $\ad_U$ has only the real eigenvalues $0, 
d, d\lambda$ of multiplicity 2.

\item[3.3] If ${\frakt g} = {\frakt g}_{4,3,3}$ then in the basis 
$X,Y,Z,T$, $$\ad_U = \pmatrix d & d & 0 & -a-b\\ 0 & d & d & -b-c\\ 0 & 0 
& d & -c\\ 0 & 0 & 0 & 0 \endpmatrix $$ and hence $\ad_U$ has only the 
real eigenvalues 0 and d of multiplicity 3.

\item[3.4] If ${\frakt g} = {\frakt g}_{4,3,4(\lambda,\varphi)}$, then in 
the basis $X,Y,Z,T$, $$\ad_U = \pmatrix d\cos\varphi & d\sin\varphi & 0 & 
-a\cos\varphi -b\sin\varphi\\ -d\sin\varphi & d\cos\varphi & 0 & 
a\sin\varphi - b\cos\varphi\\ 0 & 0 & d\lambda & -c\lambda \\ 0 & 0 & 0 & 
0 \endpmatrix$$ and hence $\ad_U$ has no purely imaginary eigenvalues $0, 
d\lambda, de^{\pm i\varphi}, \forall \varphi\in (0,\pi) , \varphi \ne 
\pi/2$. 

\item[4.2]
If ${\frakt g} = {\frakt g}_{4,4,2}$, then in the basis $X,Y,Z,Y$, $$\ad_U 
= \pmatrix -d & 0 & 0 & a\\ 0 & d & 0 & -b\\ -b & a & 0 & 0\\ 0 & 0 & 0 & 
0 \endpmatrix$$ and hence $\ad_U$ has only the real eigenvalues 0 of 
multiplicity 2 and $\pm d$.

\end{enumerate}
\end{pf}

\subsection{Description of co-adjoint orbits}

We introduce the following convention. Fix a Lie $\MD_4$-algebra ${\frakt 
g}$ with the standard basis $X,Y,Z,T$ as in the Theorem 2.1. It is 
isomorphic to ${\RR}^4$ as vector spaces. The coordinates in this standard 
basis is denote by $(a,b,c,d)$. We identify its dual vector space 
${\frakt g}^*$ with ${\RR}_4$ with the help of the dual basis 
$X^*,Y^*,Z^*,T^*$ and with the local coordinates as 
$(\alpha,\beta,\gamma,\delta)$. Thus the general form of an element of 
${\frakt g}$ is $U = aX + bY + cZ + dT$ and the general form of an 
element of ${\frakt g}^*$ is $F= \alpha X^* + \beta Y^* + \gamma Z^* + 
\delta T^*$. We denote the co-adjoint orbit passing through $F\in {\frakt 
g}^*$ by $\Omega_F$. 

\begin{thm}[\bf The Picture of Co-adjoint Orbits]
\begin{enumerate}
\item[1.1] Case $G = G_{4,1,1}$.
	\begin{enumerate}
	\item[i.] Each point $F$ with the coordinate $\gamma = 0$ is a 
0-dimensional co-adjoint orbit $\Omega_F = 
\Omega_{(\alpha,\beta,0,\delta)}$
	\item[ii.] There a family of one 2-dimensional co-adjoint orbit 
$$\Omega_F =\Omega_{\beta,\gamma\ne 0}= \{(x,\beta,\gamma,t) ; x,t \in 
{\RR}\}.$$ 
	\end{enumerate}

\item[1.2] Case $G = G_{4,1,2}$.
	\begin{enumerate}
	\item[i.] Every point $F=\alpha X^* + \beta Y^* + \delta T^*$, 
with the coordinate $\gamma = 0$ is a 0-dimensional co-adjoint orbit.
	\item[ii.] The subset $\gamma \ne 0$ decomposes into a family of 
2-dimensional co-adjoint orbits 
$$\Omega_F = \Omega_{\alpha,\beta} = \{ (\alpha,\beta, z,t) ; z,t \in 
{\RR} \gamma z > 0\} ,$$ which are all half-planes, parameterized by the 
coordinates $\alpha, \beta \in {\RR}$.
	\end{enumerate}
\item[2.1-2.2] Case $G = G_{4,2,1(\lambda)},\lambda\in{\RR}^*$ or 
$G_{4,2,2}$. 
	\begin{enumerate}
	\item[i.] Each point on the plane $\beta = \gamma = 0$ is a 
0-dimensional co-adjoint orbit, $$\Omega_F = \Omega_{\alpha,0,0,\delta}.$$
	\item[ii.] The open set $\beta^2 + \gamma^2 \ne 0 $ is decomposed 
into the union of 2-dimensional cylinders of form
$$\Omega_F = \left\{ \begin{array}{ll} \{ 
(\alpha,\beta e^{s\lambda},\gamma e^s,t); s,t \in {\RR} \} & \mbox{ in 
case } G = G_{4,2,1(\lambda)}\\ \{(\alpha, \beta e^s, \beta se^s + 
\gamma e^s,t); s,t \in {\RR}  \} &\mbox{ in case } G= G_{4,2,2}    
\end{array} \right.$$
	\end{enumerate}

\item[2.3] Case $G = G_{4,2,3(\varphi)}(\varphi\in(0,\pi))$. We identify 
${\frakt g}^*= {\frakt g}_{4,2,3(\varphi)}$ with ${\RR} \times {\CC} 
\times {\RR}$ with coordinates $(\alpha, \beta+i\gamma, \delta)$: 
	\begin{enumerate}
	\item[i.] Each point $(\alpha, 0,\delta)$ is a 0-dimensional 
co-adjoint orbit $\Omega_{\alpha, 0 +i0, \delta}$

	\item[ii.] The open set $\beta + i\gamma \ne 0$ is decomposed into a 
family of disjoint co-adjoint orbit 
$$\Omega_F = \{ (\alpha, (\beta + i\gamma)e^{se^{i\varphi}}, t) ; s,t \in 
{\RR} \}$$ which are also cylinders.
	\end{enumerate}

\item[2.4] Case $G = G_{4,2,4} = \widetilde{\Aff{\CC}}$.
	\begin{enumerate}
	\item[i.] Each point $(\alpha, 0, 0,\delta)$ is a 0-dimensional 
co-adjoint orbit $\Omega_{\alpha, 0,0,\delta}$.

	\item[i..] The open set $\beta^2 + \gamma^2 \ne 0$ is the single 
4-dimensional co-adjoint orbit 
$$\Omega_F = \Omega_{\beta^2 + \gamma^2 \ne 0}.$$

	\end{enumerate}
\item[3.1-3.3] Case $G$ is one of the groups 
$G_{4,3,1(\lambda_1,\lambda_2)}(\lambda_1,\lambda_2 \in {\RR}^*)$.
$G_{4,3,2(\lambda)} (\lambda\in {\RR}^*)$ or $G_{4,3,3}$.
	\begin{enumerate}
	\item[i.] Each point $F= \delta T^*$ on the line $\alpha = \beta 
= \gamma = 0$ is a 0-dimensional co-adjoint orbit.
	\item[ii.] The open set defined by the condition $\alpha^2 + 
\beta^2 + \gamma^2\ne 0$ is decomposed into a family of 
2-dimensional cylinders 
$$\Omega_F = \left\{ \begin{array}{l} \{(\alpha e^{s\lambda_1}, \beta 
e^{s\lambda_2}, \gamma e^s, t) ; s,t \in {\RR} \}  \quad\quad  \mbox{ in 
case } 
G=G_{4,3,1(\lambda_1,\lambda_2)}, \lambda_1, \lambda_2 \in {\RR}^*\\
\{(\alpha e^{s\lambda}, \alpha se^{s\lambda} + \beta e{s\lambda}, 
\gamma e^s, t) ; s,t \in {\RR} \} \quad\quad \mbox{ in case } G = 
G_{4,3,2(\lambda)}, \lambda\in {\RR}^*\\ 
\{(\alpha e^s, \alpha se^s + \beta e^s, {1\over 2}\alpha s^2 e^s + 
\beta se^s + \gamma e^s, t); s,t\in {\RR} \}\quad \mbox{ in case } G = 
G_{4,3,3}      \end{array}\right.$$
	\end{enumerate}

\item[3.4] Case $G = G_{4,3,4(\lambda,\varphi)}$. We identify ${\frakt 
g}^*_{4,3,4(\lambda,\varphi)}$ with ${\CC} \times {\RR}^2$ and $F = 
(\alpha,\beta,\gamma,\delta)$ with $(\alpha + i\beta , \gamma, \delta)$ 
for $\lambda \in {\RR}^*$ and $\varphi\in (0,\pi)$. 
	\begin{enumerate}
	\item[i.]
	Each point of the line defined by the condition $\alpha = \gamma 
=\beta =0$ is a 0-dimensional co-adjoint orbit.

	\item[ii.]
	The open set $\vert \alpha + i\beta\vert^2 + \gamma^2 \ne 0$ is 
decomposed into an union of co-adjoint orbits which are just cylinders 
$$\Omega_F = \{((\alpha + i\beta\vert^2)e^{se^{i\varphi}}, \gamma 
e^{s\lambda}, t); s,t \in {\RR}\}$$

	\end{enumerate}
\item[4.1] Case $G = G_{4,4,1} = {\RR}\ltimes_J {\frakt h}_3$. 
	\begin{enumerate}
	\item[i.]
	Each point of the line defined by the conditions $ \alpha= \beta 
= \gamma = 0$ is a 0-dimensional orbit $\Omega_F = 
\Omega_{(0,0,0,\delta)}=\{ (0,0,0,\delta)$. 

	\item[ii.] The open set $\gamma \ne 0$ is decomposed into a union 
of 2-dimensional co-adjoint orbits $$\Omega_F = \{ (x,y,\gamma, t); 
x^2+y^2-2\gamma 
t = \alpha^2+\beta^2-2\gamma\delta \},$$ which are just rotation 
paraboloids, and

	\item[iii.] the set $\alpha^2+\beta^2 \ne 0$, $\gamma =0$ is a 
union of 2-dimensional co-adjoint orbits, which are just cylinders.
$$\Omega_F = \{(x,y,0,t); x^2+y^2 = \alpha^2+\beta^2 \}.$$
	\end{enumerate}

\item[4.2] Case $G = G_{4,4,2} = {\RR} \ltimes {\HH}_3$, the real diamond 
group. 
	\begin{enumerate}
	\item[i.]
	Each point of the line $\alpha=\beta=\gamma=0$ is a 0-dimensional 
co-adjoint orbit $\Omega_F = \Omega_{(0,0,0,\delta)}$
	\item[ii.]
	The set $\alpha \ne 0, \beta= 0 = \gamma$ is union of two 
2-dimensional co-adjoint orbits, which are just the half-planes $$\Omega_F 
= \{ (x,0,0,t); x,t\in {\RR}, \alpha x > 0\}.$$

	\item[iii.] The set $\alpha =0= \gamma$, $\beta\ne 0$ is union of 
two 2-dimensional co-adjoint orbits, which are just the half-planes $$\Omega_F 
\{(0,y,0,t); y,t\in {\RR}    \}.$$
	\item[iv.] The open set $\gamma \ne 0$ is decomposed into a 
family of 2-dimensional co-adjoint orbits, which are just the hyperbolic 
paraboloids $$\Omega_F = \{(x,y,\gamma,t); x,y,t\in{\RR}, xy-\alpha\beta 
= \gamma(t-\delta)  \}.$$
	\end{enumerate}
\end{enumerate}
\end{thm}
\begin{pf}
Mainly, the theorem is proved by a direct computation. We have in general 
$$ \Omega_F = \{ K(g)F ; g\in G \}.$$  Because of Corollary 
3.1, we have 
$$\Omega_F  = \{K(\exp(U)F; U \in {\frakt g}  \}.$$
Recall that $$\langle K(\exp(U)F,X\rangle = \langle F, 
\exp(\ad_U)X\rangle.$$
Thus we consider a general element $$K(\exp(U))F = xX^* + yY^* + zZ^* 
+tT^* = (x,y,z,t)\in {\RR}_4,$$ where $$\begin{array}{lll}
x &= \langle K(\exp(U)F,X\rangle &= \langle F, \exp(\ad_U)X\rangle,\\ 
y &= \langle K(\exp(U)F,Y\rangle &= \langle F, \exp(\ad_U)Y\rangle,\\
z &= \langle K(\exp(U)F,Z\rangle &= \langle F, \exp(\ad_U)Z\rangle,\\
t &= \langle K(\exp(U)F,T\rangle &= \langle F, \exp(\ad_U)T\rangle .
\end{array}$$

By a direct computation, for $U = aX+bY+cZ+dT$ we have:

1.1. $G = G_{4,1,1}$, 

$$\ad_U = \pmatrix 0 & 0 & 0 & 0\\ 0 & 0 & 0 
& 0\\ d & 0 & 0 & -a \\ 0 & 0 & 0 & 0 \endpmatrix \mbox{  and  } 
\exp(\ad_U) = \pmatrix 1 & 0 & 0 & 0\\ 0 & 1 & 0 & 0\\ d & 0 & 1 & -a\\ 0 
& 0 & 0 & 1 \endpmatrix.$$ 
This means that $$\begin{array}{rl}
x &= \alpha + \gamma d,\\
y &= \beta,\\
z &= \gamma,\\
t &= -\gamma a + \delta
\end{array}$$
and therefore if $\gamma = 0$, each point is unchanged, and is therefore 
a 0-dimensional co-adjoint orbit. In other words, the 0-dimensional orbits 
are parameterized by the points $(\alpha,\beta,\gamma\ne 0,\delta) \in 
{\RR}^2 \times {\RR}^* \times {\RR}$. If $\gamma \ne 0$, the coordinates $x = 
\alpha + \gamma d$ and $t = -\gamma a + \delta$ run over two coordinates 
lines, while the coordinates $y=\beta$ and $z=\gamma \ne 0$ are fixed. 
Thus we have a family of 2-dimensional co-adjoint orbits, parameterized by 
the points $(\beta,\gamma) \in {\RR} \times {\RR}^*$.

1.2. $G = G_{4,1,2}$. 

$$\ad_U = \pmatrix 0 & 0 & 0 & 0\\ 0 & 0 & 0 & 0\\ 0 & 0 & d & -c\\ 0 & 
0 & 0 & 0 \endpmatrix \mbox{  and  } \exp(\ad_U) = \pmatrix 1 & 0 & 0 & 0 
\\ 0 & 1 & 0 & 0\\ 0 & 0 & e^d & -c\sum_{n=1}^\infty {d^{n-1} \over n!}\\ 
0 & 0 & 0 & 1 \endpmatrix.$$
This means that 
$$\begin{array}{rl}
x &= \alpha,\\ 
y &= \beta,\\
z &= \gamma e^d,\\
t &= -\gamma c\sum_{n=1}^\infty {d^{n-1}\over n!} + \delta
\end{array}$$ 
and therefore the point $F$ should be unchanged if $\delta 
= 0$; otherwise, $\gamma \ne 0$, the coordinate $t$ run over a line and 
the coordinate $z$ run over a half-line. This means that the closet set 
$\gamma=0$ is decomposed into 0-dimensional co-adjoint orbits. The open 
set $\gamma \ne 0$ decomposed into a family of 2-dimensional co-adjoint 
orbits, which are just half-planes
$$\Omega_F = \{(\alpha,\beta,z,t); z,t\in {\RR}, z.\gamma > 0 \},$$
parameterized by $(\alpha,\beta,\sgn(\gamma)) \in {\RR}^2 \times \{ \pm \}.$

2.1. $G = G_{4,2,1(\lambda)} (\lambda \in {\RR}^*)$. 

$$\ad_U = \pmatrix 0 & 0 & 0 & 0\\ 0 & d\lambda & 0 & -\lambda b\\ 0 & 0 
& d & -c\\ 0 & 0 & 0 & 0 \endpmatrix \mbox{  and  } \exp(\ad_U) = 
\pmatrix 1 & 0 & 0 & 0\\ 0 & e^{d\lambda} & 0 & -b\sum_{n=1}^\infty 
{\lambda^n d^{n-1}\over n!}\\ 0 & 0 & e^d & -c\sum_{n=1}^\infty 
{d^{n-1}\over n!}\\ 0 & 0 & 0 & 1 \endpmatrix .$$ This means that 
$$\begin{array}{rl} 
x &= \alpha,\\
y &= \beta e^{d\lambda},\\
z &= \gamma e^d,\\
t &= \sum_{n=1}^\infty (-\beta b\lambda^n) - \gamma c)({d^{n-1}\over n!}) 
+ \delta \end{array}$$ and therefore the point $F$ should be unchanged if 
$\beta = \gamma = 0$ and if one of them is nonzero then the coordinates 
$(x,y,z,t)$ cover a 2-dimensional cylinder $(\alpha, \beta e^{s\lambda}, 
\gamma e^s, t) ; s,t \in {\RR}$. Thus the open set $\beta^2 + \gamma^2 \ne 
0$ is decomposed into a family of 2-dimensional co-adjoint orbits, which 
are cylinders of form
$$\Omega_F = \{(\alpha,\beta e^{d\lambda}, z,t) ; z,t\in {\RR} \},$$ 
parameterized by ${\RR}\times {\SS}^1$. 

2.2. $G = G_{4,2,2}$. 

$$\ad_U = \pmatrix 0 & 0 & 0 & 0\\ 0 & d & d & -b-c\\ 0 & 0 & d & -c\\ 0 
& 0 & 0 & 0 \endpmatrix \mbox{  and  } \exp(\ad_U) = \pmatrix 1 & 0 & 0 & 0\\
0 & e^d & de^d & -b\sum_{n=1}^\infty {d^{n-1} \over n!} - c e^d\\ 0 & 0 & 
e^d & -c\sum_{n=1}^\infty {d^{n-1} \over n!}\\ 0 & 0 & 0 & 1 \endpmatrix.$$
This means that 
$$\begin{array}{rl}
x &= \alpha,\\
y &= \beta e^d,\\
z &= \beta de^d + \gamma e^d,\\
t &= -\beta b\sum_{n=1}^\infty {d^{n-1}\over n!} - \beta ce^d - \gamma 
c\sum_{n=1}^\infty {d^{n-1}\over n!} + \delta\end{array}$$ and therefore the 
point 
$F$ is unchanged if $\beta = \gamma = 0$; otherwise the coordinates 
$(x,y,z,t)$ cover a 2-dimensional cylinder 
$$\{(\alpha, \beta e^s, \beta se^s + \gamma e^s, t); s,t \in {\RR} \}.$$

2.3. $G= G_{4,2,3(\varphi)} (\varphi \in (0,\pi))$.

$$\ad_U = \pmatrix 0 & 0 & 0 & 0\\ 0 & d\cos\varphi & d\sin\varphi & 
-b\cos\varphi-c\sin\varphi\\ 0 & -d\sin\varphi & d cos\varphi & 
b\sin\varphi -c\cos\varphi\\ 0 & 0 & 0 & 0 \endpmatrix $$ and 
$\exp(\ad_U) = $ $$\pmatrix 1 & 0 & 0 & 0 \\ 0 & 
e^{d\cos\varphi}\cos(d\sin\varphi) & e^{d\cos\varphi}\sin(d\cos\varphi) & 
-\sum_{n=1}^\infty {d^{n-1}\over n!}(b\cos(n\varphi) + c\sin(n\varphi))\\ 
0 & -e^{d\cos\varphi}\sin(d\sin\varphi) & 
e^{d\cos\varphi}\cos(d\sin\varphi) & -c\sum_{n=1}^\infty {d^{n-1}\over 
n!} \cos(n\varphi)\\ 0 & 0 & 0 & 1 \endpmatrix $$.
This means that 
$$\begin{array}{rl}
x &= \alpha,\\
y &= \beta e^{d\cos\varphi} \cos(d\sin\varphi) - \gamma e^{d\cos\varphi} 
\sin(d\sin\varphi) \\ 
z &= \beta e^{d\cos\varphi} \sin(d\sin\varphi) + \gamma e^{d\cos\varphi} 
\cos(d\sin\varphi)\\ 
t &= -\beta\sum_{n=1}^\infty {d^{n-1}\over n!}(b\cos(n\varphi) + 
c\sin(n\varphi)) - \gamma c\sum_{n=1}^\infty {d^{n-1}\over n!} 
\cos(n\varphi) + \delta 
\end{array}. $$
We identify ${\frakt g}^*$ with ${\RR}\times {\CC} \times {\RR}$, in 
writing the coordinates $\alpha, \beta,\gamma,\delta$ as $(\alpha, 
\beta+i\gamma, \delta)$, we have
$$\begin{array}{rl}
x &= \alpha,\\
y+iz &= e^{de^{i\varphi}}(\beta + i\gamma),\\
t &= -\beta\sum_{n=1}^\infty {d^{n-1}\over n!} (b\cos(n\varphi) + 
c\sin(n\varphi)) - \gamma c\sum_{n=1}^\infty {d^{n-1}\over n!} 
\cos(n\varphi) + \delta
\end{array}$$

This means that every point $(\alpha, 0 ,0 ,\delta)$ on the line $\beta + 
i\gamma = 0$ is unchanged 
under the co-adjoint action and provides a 0-dimensional co-adjoint orbit; 
otherwise, the open set $\beta + i\gamma \ne 0$ is decomposed into a 
family of cylinders,which are just the 2-dimensional co-adjoint orbits.

2.4. $G=G_{4,2,4} = \widetilde{\Aff{\CC}}$.

$$\ad_U = \pmatrix 0 & 0 & 0 & 0\\ -c & d & a & -b\\ b & -a & d & -c\\ 0 
& 0 & 0 & 0 \endpmatrix \mbox{ and } \exp(\ad_U) = \pmatrix 1 & 0 & 0 & 
0\\ L & M & N & P\\ -P & -N & M & L\\ 0 & 0 & 0 & 0 & 1 \endpmatrix,$$
where if $a^2 + d^2 \ne 0$ then $$\begin{array}{rl}
L &:= {1\over a^2+d^2}\left[ 
(ab+cd)(1-e^{{a^2+d^2+d\over\sqrt{a^2+d^2}}}\cos({a\over \sqrt{a^2+d^2}}) 
+(bd-ac)e^{{a^2+d^2+d\over\sqrt{a^2+d^2}}}\sin({a\over\sqrt{a^2+d^2}})\right],\\
M &:= e^{{a^2+d^2+d\over\sqrt{a^2+d^2}}} \cos({a\over\sqrt{a^2+d^2}}),\\
N &:= e^{{a^2+d^2+d\over\sqrt{a^2+d^2}}}\sin({a\over\sqrt{a^2+d^2}}),\\
P &:= {1\over a^2+d^2}\left[(ac-bd)(e^{{a^2+d^2+d\over\sqrt{a^2+d^2}}} 
\cos({a\over\sqrt{a^2+d^2}}) -1) 
-(ab+cd)e^{{a^2+d^2+d\over\sqrt{a^2+d^2}}}\sin ({a\over\sqrt{a^2+d^2}}) 
\right] \end{array}$$
and if $a=d=0$,
$L := -c, M := 1, N := 0, P := -b$.
This means that 
$$\begin{array}{rl}
x &= \alpha + \beta L -\gamma P;\\
y &= \beta M - \gamma N'\\
z &= \beta N + \gamma M,\\
t &= \beta P + \gamma L + \delta
\end{array}$$
and hence every point $F$ of the coordinate plane $\{(\alpha,0,0,\delta), 
\alpha,\delta\in {\RR}\}$ is a 0-dimensional co-adjoint orbits, otherwise, 
if $\beta^2 + \gamma^2 \ne 0$, this open set is just the single 
4-dimensional co-adjoint orbits.

3.1. $G=G_{4,3,1(\lambda_1,\lambda_2)}(\lambda_1,\lambda_2 \in {\RR}^*)$

$$\ad_U = \pmatrix d\lambda_1 & 0 & 0 & -a\lambda_1\\ 0 & d\lambda_2 & 0 
& -b\lambda_2\\ 0 & 0 & d & -c\\ 0 & 0 & 0 & 0 \endpmatrix \mbox{ ; } 
\exp(\ad_U) = \pmatrix e^{d\lambda_1} & 0 & 0 & -a\sum_{n=1}^\infty 
{(d\lambda_1)^{n-1}\over n!}\\ 0 & e^{d\lambda_2} & 0 & 
-b\sum_{n=1}^\infty {(d\lambda_2)^{n-1}\over n!}\\ 0 & 0 & e^d & 
-c\sum_{n=1}^\infty {d^{n-1}\over n!}\\ 0 & 0 & 0 & 1\endpmatrix.$$ This 
means that 
$$\begin{array}{rl}
x &= \alpha e^{d\lambda_1},\\
y &= \beta e^{d\lambda_2},\\
z &= \gamma e^d,\\
t &= -\alpha a\sum_{n=1}^\infty {(d\lambda_1)^{n-1}\over n!} - \beta 
b\sum_{n=1}^\infty {(d\lambda_2)^{n-1}\over n!} - \gamma 
c\sum_{n=1}^\infty {d^{n-1}\over n!} + \delta
\end{array}$$
and hence each point of the line $\alpha = \beta = \gamma = 0$ is a 
0-dimensional co-adjoint orbit; otherwise, the open set $\alpha^2 = 
\beta^2 + \gamma^2 \ne 0$ is decomposed into a family of 2-dimensional 
co-adjoint orbits,which are just the cylinders
$$\Omega_F = \{(\alpha e^{s\lambda_1}, \beta e^{s\lambda_2} \gamma e^s, 
t) ; s,t \in {\RR}    \}.$$

3.2. $G = G_{4,3,2(\lambda)} (\lambda \in {\RR}^*)$. 

$$\ad_U = \pmatrix d\lambda & d & 0 & -\lambda a - b\\ 0 & d\lambda & 0 & 
-b\lambda \\ 0 & 0 & d & -c\\ 0 & 0 & 0 & 0 \endpmatrix ;$$ 
and $$ \exp(\ad_U) = 
\pmatrix e^{d\lambda} & d e^{d\lambda} & 0 & -a\sum_{n=1}^\infty 
{(d\lambda)^{n-1}\over n!} - be^{d\lambda}\\ 0 & e^{d\lambda} & 0 & 
-b\sum_{n=1}^\infty {(d\lambda)^{n-1}\over n!}\\ 0 & 0 & e^d 
-c\sum_{n=1}^\infty {d^{n-1}\over n!}\\ 0 & 0 & 0 & 1 \endpmatrix.$$ This 
means that $$\begin{array}{rl}
x &= \alpha e^{d\lambda},\\
y &= \alpha de^{d\lambda} + \beta e^{d\lambda},\\
z &= \gamma e^{d\lambda},\\
t &= -\alpha a\sum_{n=1}^\infty {(d\lambda)^{n-1}\over n!} - \alpha 
be^{d\lambda} - \beta b\sum_{n=1}^\infty {(d\lambda)^{n-1}\over n!} - 
\gamma c\sum_{n=1}^\infty {d^{n-1}\over n!} + \delta
\end{array}$$
and hence every point of the line $\alpha=\beta=\gamma=0$ is just a 
0dimensional co-adjoint orbit $$\Omega_{(0,0,0,\delta)} = \{(0,0,0,\delta) 
\};$$ otherwise, the open set $\alpha^2 = \beta^2 + \gamma^2 \ne 0$ is 
decomposed into a family of 2-dimensional co-adjoint orbits which are just 
the cylinders
$$\Omega_F = \{(\alpha e^{s\lambda}, \alpha se^{s\lambda} + \beta 
e^{s\lambda}, \gamma e^{s\lambda}, t) ; s,t \in {\RR}   \}.$$

3.3. $G = G_{4,3,3}$.

$$\ad_U = \pmatrix d & d & 0 & -a-b\\ 0 & d & d & -b-c\\ 0 & 0 & d & -c\\ 
0 & 0 & 0 & 0 \endpmatrix ;$$ 
$$ \exp(\ad_U) = \pmatrix e^d & de^d & {1\over 
2} d^2e^d & -\alpha\sum_{n=1}^\infty {d^{n-1}\over n!} -(b+{1\over 
2}cd)e^d\\ 0 & e^d & de^d & -b\sum_{n=1}^\infty {d^{n-1}\over n!} - 
ce^d\\ 0 & 0 & e^d & -c\sum_{n=1}^\infty {d^{n-1}\over n!}\\ 0 & 0 & 0 & 
1 \endpmatrix.$$ This means that $$\begin{array}{rl}
x &= \alpha e^d,\\
y &= \alpha de^d +\beta e^d,\\
z &= {1\over 2}\alpha d^2 e^d + \beta de^d + \gamma e^d,\\
t &= -\alpha\left(a\sum_{n=1}^\infty {d^{n-1}\over n!} (b+{1\over2}cd)e^d    
\right) - \beta \left( b\sum_{n=1}^\infty {d^{n-1}\over n!} + ce^d\right) 
- \gamma c\sum_{n=1}^\infty {d^{n-1}\over n!} + \delta\end{array}$$ 
and hence every 
point of the line $\alpha =\beta=\gamma = 0$ is a 0-dimensional co-adjoint 
orbit $$\Omega_{(0,0,0,\delta)}= \{(0,0,0,\delta); \delta = \mbox{fixed} 
\};$$ otherwise, the open set $\alpha^2+\beta^2 + \gamma^2 \ne 0$ is 
decomposed into a family of 2-dimensional co-adjoint orbits, which are 
just the cylinders
$$\Omega_F = \{(\alpha e^s, \alpha se^s + \beta e^s, {1\over 2}\alpha s^2 
e^s + \beta se^s + \gamma e^s, t); s,t\in {\RR}   \}.$$

3.4. $G=G_{4,3,4(\lambda,\varphi)}(\lambda \in {\RR}^*, \varphi\in (0,\pi))$.

$$\ad_U = \pmatrix d\cos\varphi & d\sin\varphi 7 0 & 
-a\cos\varphi-b\sin\varphi\\ -d\sin\varphi & d\cos\varphi & 0 & 
a\sin\varphi -b\cos\varphi\\ 0 & 0 & d\lambda & -c\lambda\\ 0 & 0 & 0 & 0 
\endpmatrix$$ and $\exp(\ad_U) = $
$$\pmatrix e^{d\cos\varphi}\cos(d\sin\varphi) & 
e^{d\cos\varphi}\sin(d\sin\varphi) & 0 & \sum_{n=1}^\infty {d^{n-1}\over 
n!} (a\cos(n\varphi) + b\sin(n\varphi))\\ 
-e^{d\cos\varphi}\sin(d\sin\varphi) & e^{d\cos\varphi}\cos(d\sin\varphi) 
& 0 & \sum_{n=1}^\infty{d^{n-1}\over n!} (a\sin(n\varphi) 
-b\cos(n\varphi)\\ 0 & 0 & e^{d\lambda} & -c\sum_{n=1}^\infty 
{d^{n-1}\over n!}\\ 0 & 0 & 0 & 1 \endpmatrix.$$ This means that 
$$\begin{array}{rl}
x &= \alpha e^{d\cos\varphi}\cos(d\sin\varphi) - \beta 
e^{d\cos\varphi}\sin(d\sin\varphi),\\
y &= \alpha e^{d\cos\varphi} \sin(d\sin\varphi) + \beta 
e^{d\cos\varphi}\cos(d\sin\varphi),\\
z &= \gamma e^{d\lambda},\\
t &= -\alpha\sum_{n=1}^\infty {d^{n-1}\over n!} 
(a\cos(n\varphi)+b\sin(n\varphi)) +\beta\sum_{n=1}^\infty {d^{n-1}\over 
n!} (a\sin(n\varphi)- \\
  &\phantom{= } -b\cos(n\varphi)) - \gamma \sum_{n=1}^\infty 
{d^{n-1}\over n!} + \delta 
\end{array}.$$ We identify ${\frakt g}_{4,3,4(\lambda,\varphi)}^*$ with 
${\CC}\times {\RR}^2$ by identifying $(\alpha,\beta,\gamma,\delta)$ with 
$(\alpha+i\beta, \gamma,\delta)$, then we can rewrite
$$x+iy = (\alpha+i\beta)e^{de^{i\varphi}},$$ and hence each point of the 
line $\alpha + i\beta = 0, \gamma = 0$ is just a 0-dimensional 
co-adjoint orbit $$\Omega_{(0,0,0,\delta)} = \{(0,0,0,\delta); \delta= 
\mbox{ fixed } \},$$ otherwise the open set $\vert \alpha 
+i\beta\vert^2+\gamma^2 \ne 0$ is decomposed into a family of 
2-dimensional co-adjoint orbits, which are just the cylinders of form
$$\Omega_F = \{(\alpha+i\beta)e^{se^{i\varphi}},\gamma e^{s\lambda},t); 
s,t\in {\RR}  \}.$$

4.1. $G= G_{4,4,1}$.

$$\ad_U = \pmatrix 0 & d & 0 & -b\\ -d & 0 & 0 & 0\\ 0 & 0 & 0 & 0 \\ 0 & 
0 & 0 & 0 \endpmatrix$$  and 
$$\exp(\ad_U)= \pmatrix \cos d & \sin d & 0 & b\sum_{n=0}^\infty (-1)^n 
{d^{2n} \over (2n+1)!}- a\sum_{n=1}^\infty (-1)^n{d^{2n-1}\over (2n)!}\\ 
-\sin d & \cos d & 0 & a\sum_{n=0}^\infty (-1)^n{d^{2n}\over (2n+1)!}- 
b\sum_{n=1}^\infty (-1)^n {d^{2n-1}\over (2n)!}\\a_{31} & a_{32} & a_{33} 
& a_{34}\\ 0 & 0 & 0 & 1 \endpmatrix ,$$ 
where $$\begin{array}{rl}
a_{31} &:=
-b\sum_{n=0}^\infty(-1)^n{d^{2n}\over(2n+1)!} + 
a\sum_{n=1}^\infty(-1)^n{d^{2n-1}\over (2n)!}\\
a_{32} &:=  
a\sum_{n=0}^\infty(-1)^n{d^{2n}\over(2n+1)!}+b\sum_{n=1}^\infty 
(-1)^n{d^{2n-1}\over (2n)!} \\
a_{33} &:= 1\\ 
a_{34} &:=
-(a^2+b^2)\sum_{n=1}^\infty(-1)^n{d^{2n-2}\over (2n)!}
\end{array} .$$ 
This means that $$\begin{array}{rl}
x &= \alpha\cos d - \beta\sin d + \gamma\left[a\sum_{n=1}^\infty 
(-1)^n{d^{2n-1}\over (2n)!} - b\sum_{n=0}^\infty(-1)^n{d^{2n}\over 
(2n+1)!}   \right]\\ 
y &= \alpha\sin d + \beta\cos d + 
\gamma\left[a\sum_{n=0}^\infty(-1)^n{d^{2n}\over (2n+1)!} + 
b\sum_{n=0}^\infty (-1)^n{d^{2n-1}\over (2n)!} \right]\\
z &= \gamma,\\
t &= -\alpha\left[b\sum_{n=0}^infty (-1)^n {d^{2n}\over (2n+1)!} + 
a\sum_{n=1}^\infty (-1)^n{d^{2n+1}\over (2n)!} \right] + 
   \beta [  
a\sum_{n=0}^\infty (-1)^n \times\\  &  \times {d^{2n}\over (2n+1)!}    
 - b\sum_{n=1}^\infty (-1)^n {d^{2n-1}\over (2n)!} ] - 
\gamma(a^2+b^2)\sum_{n=1}^\infty (-1)^n{d^{2n-2}\over (2n)!} +\delta
\end{array}$$
and hence every point of the line $\alpha=\beta=\gamma =0$ is just a 
0-dimensional co-adjoint orbit $$\Omega_{(0,0,0,\delta)} = 
\{(0,0,0,\delta); \delta =\mbox{ fixed } \},$$ otherwise the set 
$\alpha^2 + \beta^2 \ne 0$, $\gamma = 0$ is decomposed into a family of 
2-dimensional co-adjoint orbits, which are just the rotation cylinders
$$\Omega_F = \{ (x,y,0,t); y,y,t\in {\RR} , x^2+y^2 = \alpha^2 + 
\beta^2\}$$ and the open set $\gamma \ne 0$ is decomposed into a family 
of 2-dimensional co-adjoint orbits, which are just the elliptic paraboloids
$$\Omega_F = \{(x,y,\gamma,t); x,y,t \in {\RR}, x^2+y^2 -2\gamma t = 
\alpha^2 + \beta^2 - 2\gamma\delta \}.$$

4.2. $G = G_{4,4,2} = {\RR} \ltimes {\frakt h}_3$, the real diamond group.

$$\ad_U = \pmatrix -d & 0 & 0 & a\\ 0 & d & 0 -b\\ -b & a & 0 & 0\\ 0 & 0 
& 0 & 0\endpmatrix$$ and 
$$\exp(\ad_U) = \pmatrix e^{-d} & 0 & 0 & a\sum_{n=1}^\infty 
(-1)^n{d^{n-1}\over n!}\\ 0 & e^d & 0 & -b\sum_{n=1}^\infty {d^{n-1}\over 
n!}\\ b\sum_{n=1}^\infty (-1)^n{d^{n-1}\over n!} & a\sum_{n=1}^\infty 
{d^{n-1}\over n!} & 1 & -ab\sum_{n=2}^\infty ((-1)^n +1){d^{n-2}\over 
n!}\\ 0 & 0 & 0 & 1 \endpmatrix.$$ This means that 
$$\begin{array}{rl}
x &= \alpha e^{-d} + \gamma b\sum_{n=1}^\infty (-1)^n {d^{n-1}\over n!},\\
y &= \beta e^d + \gamma a\sum_{n=1}^\infty {d^{n-1}\over n!},\\
z &= \gamma,\\
t &= \alpha a\sum_{n=1}^\infty(-1)^n{d^{n-1}\over n!} - \beta 
b\sum_{n=1}^\infty {d^{n-1}\over n!} - \gamma ab\sum_{n=2}^\infty 
((-1)^n+1){d^{n-2}\over n!} + \delta\end{array}$$ and hence every point of 
the line 
$\alpha=\beta=\gamma=0$ is a 0-dimensional co-adjoint orbit $$\Omega_F = 
\{(0,0,0,\delta); \delta = \mbox{ fixed } \},$$ otherwise the set 
$\alpha\ne 0$, $\beta = \gamma = 0$ is decomposed into two 
2-dimensional co-adjoint orbits, which are just two coordinates half-planes
$$\Omega_F = \{(x,0,0,t); x,t\in {\RR}, \alpha x > 0 \},$$ the set $\beta 
\ne 0$, $\alpha=\gamma = 0$ is decomposed into two co-adjoint orbits, 
which are just two coordinate half-planes 
$$\Omega_F = \{(0,y,0,t); y,t\in{\RR}, \beta y > 0  \},$$ the set 
$\alpha\beta \ne 0$, $\gamma = 0$ is decomposed into a family of 
2-dimensional co-adjoint orbits, which are just hyperbolic cylinders
$$\Omega_F = \{(x,y,0,t); x,y,t\in{\RR}, \alpha x> 0, \beta y > 0, xy = 
\alpha\beta \}$$ and finally, the open set $\gamma\ne 0$ is decomposed 
into a family of 2-dimensional co-adjoint orbits, which are just parabolic 
hyperboloids
$$\Omega_F = \{(x,y,\gamma,t); x,y,t\in{\RR}, xy-\alpha\beta = 
\gamma(t-\delta) \}.$$ 

\end{pf}

\section{Measurable MD4-Foliation} 
In this section, we study the C*-algebras associated with co-adjoint 
orbits. in many cases this should give us adequate informations about th 
structure of group C*-algebras. We recall first of all the background 
from Connes theory of measurable foliations. We show later that the 
generic co-adjoint orbits provide connes measurable foliations and finally 
we consider topological classification of these foliations.

\subsection{Measurable foliations after A. Connes}
Let us recall that a {\it integrable tangent distribution}\index{tangent 
distribution!integrable} is by 
definition a smooth sub-fibration ${\cal F}$ of the tangent bundle $TV$ of a 
smooth manifold $V$, such that each point $x\in V$ can be included in a 
smooth sub-manifold $W \subseteq V$, which recognizes ${\cal F}= TW$ as 
its tangent bundle, i.e. the fiber ${\cal F}_x$ at $x$ is coincided with 
the tangent space $T_xW$, for all $x\in W$. The manifold $W$ in this case 
is called the {\it integral manifold} \index{manifold!integral - } of 
${\cal F}$. Recall also the Frobenius criteria of integrability
\begin{prop}[A. Connes \cite{connes1}]
The following conditions are equivalent:
\begin{enumerate}
\item[i.] The tangent distribution ${\cal F}$ is integrable on $V$.
\item[ii.] For all $x\in V$, there exists an open sub-manifold $U$ in $V$, 
containing $x$, and a submersion $p : U \to {\RR}^q (q \codim {\cal F} := 
\dim V - \dim {\cal F})$ such that ${\cal F}_y = \ker(p_*)_y, \forall y\in 
W$.
\item[iii.] $C^\infty({\cal F}) = \{ s\in C^\infty(TV); s_x \in {\cal 
F}_x, \forall x\in V \}$ is a Lie subalgebra of smooth vector fields.
\item[iv.] The ideal $I({\cal F})$ of differential forms, vanishing on 
${\cal F}$ is stable under the operation of exterior differentiation $d$.
\end{enumerate}
\end{prop}
This deduces in particular that all the 1-dimensional tangent 
distributions are integrable. 

A manifold $V$, equipped with an integrable tangent distribution ${\cal 
F}$ is called a 
{\it foliation} \index{foliation} or a {\it foliated manifold}, 
\index{manifold!foliated -} denoted by $(V,{\cal F})$. 
Each maximal connected integral sub-manifold $L$ of ${\cal F}$ is called a 
{\it leaf} of the foliation \index{foliation!leaf of - } $(V,{\cal F})$. 
It is 
reasonable to recall (see for example A. Connes \cite{connes1}) that:
\begin{enumerate}
\item[i.] The family of all leaves of a foliation ${\cal F}$ on $V$ form a 
partition of $V$.
\item[ii.] For each point $x\in V$, there exists a coordinate 
neighborhood $\{U,(x^1,\dots,x^n)\}$, $n = \dim V$, such that if a fiber 
$L$ intersects with this neighborhood, $L \cap U \ne \empty$, each 
connected component of $L \cap U$, which is called a {\it plaque}, 
\index{foliation!plaque of - } is given by the following equations
$$x^{k+1} = c^1, \dots , x^n = c^{n-k}, n = \dim{\cal F},$$ where $c^1, 
\dots, c^{n-k}$ are some constants, depending on the plaque. 
\end{enumerate}
The atlas of 
this kind coordinate charts is called an {\it atlas of foliated manifold. 
} \index{atlas of foliated manifold} Conversely, one can use these last two 
properties to defines a larger 
class of foliations (see I. Tamura \cite{tamura}[Ch. 4, pp. 121--126].): 
one consider a family ${\cal C}$ of sub-manifolds, satisfying two 
conditions and for each fiber $L\in {\cal C}$ there exists a single 
integrable distribution ${\cal F}$ such that $L$ is its maximal 
connected integral sub-manifold of ${\cal F}$. Then we have also a 
foliations.  Remark that locally, all foliations of a fixed dimension 
have the same local structure. But globally, they are quite different, 
e.g. compactness of leaves, existence of dense leaves, .... The space of 
leaves $V/{\frakt F}$.

\begin{defn} Two foliations $(V_1,{\cal f}_1)$ and $(V_2,{\cal F}_2)$ are 
called {\it topologically equivalent}, \index{foliations!topologically 
equivalent -} if there exists a 
leaf-wise homeomorphism $h : (V_1,{\cal F}_1) \to (V_2,{\cal F}_2)$.

Recall that a sub-manifold $N$ of the foliated manifold $(V,{\cal F})$ is 
called {\it transversal} \index{sub-manifold!transversal - } if at each 
point $p\in V$, we can split $T_pV$ 
as the direct sum $T_pV = {\cal F}_p \oplus TpN$, for all $p \in N$. 
Certainly that in this case $\dim N = codim{\cal F}$. In a small 
neighborhood $(U,(x^1,\dots, x^n))$ there is a 1-1 correspondence between the 
plaques and the point of $U\cap N$. If for a Borel set $B$, the set $U\cap 
B$ is countable, the transversal set $B$ is called {\it Borel 
transversal}. \index{Borel transversal}
 It was proved in \cite{connes1} that there exists an injection $\psi : B 
\to N$, where $N$ is a transversal sub-manifold such that $\psi(x) \in 
L_x,$ the leaf containing $x$.

\end{defn}

\begin{defn}
A $\sigma$-additive measure $B \mapsto \Lambda(B)$ from the set of Borel 
transversals to the set $[0,+\infty]$ is called a {\it transversal 
measure} \index{measure! transversal -} if:
\begin{enumerate}
\item[$(\Lambda_1)$]{\bf (Borel equivalence)}. the measure is invariant 
w.r.t. 
Borel bijections $\psi : B_1 \to B_2,$ $\Lambda(B_1) = \Lambda(B_2).$
\item[$(\Lambda_2)$] $\Lambda(K) < +\infty$ if $K$ is a compact subset of 
a transversal sub-manifold.
\end{enumerate}
A foliation $(V,{\cal F})$, equipped with a transversal measure is called 
a {\it measurable foliation}. \index{foliation! measurable -}
\end{defn}

Let us recall finally some relation between the transversal measure and 
the ordinary measure. 

For oriented foliations, we can deduce a more clear relationship between 
the transversal measures and the ordinary measure on foliated manifolds.
 Choose an orientation of ${\cal F}$. Then the fiber bundle $\wedge^k 
{\cal F} (k=\dim{\cal F})$ is decomposed into two parts $(\wedge^k{\cal 
F})^+$ and $(\wedge^k{\cal F})^-$ by the zero section. Fix a $k$-vector 
field $X\in C^\infty(\wedge^k{\cal F})^+$ and a measure $\mu$ on the 
foliated manifold $(V,{\cal F})$. If $U$ is some local coordinate cart of 
the foliation, then $U$ can be identified with the direct product 
$N\times\pi$ of some transversal sub-manifold $N$ and a typical fiber 
$\pi$. Thus the restriction $\mu_U$ is separated into the product of
measures $\mu_N$ on $N$ and $\mu_\pi$ on $\pi$. Denote $\mu_X$ the 
measure along leaves defined by the volume element $X$. The measure $\mu$ 
is $X$-invariant iff $\mu_X$ and $\mu\pi$ are proportional for all 
coordinate cart of foliation. Two pairs $(X,\mu)$ and $(Y,\nu)$ are said 
to be equivalent iff there exists a function $\varphi\in C^\infty(V)$ 
such that $Y=\varphi X$ and $\mu = \varphi \nu$.  
\begin{prop}[A. Connes \cite{connes1}] If $(V,{\cal F})$ is an oriented 
foliated manifold, there is a bijective correspondence between 
equivalent classes of pairs $(X,\mu)$ and the transversal measures.
\end{prop}
The transversal measure, corresponding to the pair $(X,\mu)$ is just 
given by the formula
$$\Lambda(B) = \int_N \Card (B\cap N) d\mu_N(\pi),$$ for all transversal 
Borel set $B$ in $U$, and is continued to other Borel set by 
$\sigma$-additivity. 

We conclude that in order to describe some foliation as a measurable 
foliation, we need to pick out a suitable pair $(X,\mu)$.

\subsection{Measurable $\MD_4$-foliations}

\begin{thm}
If $G$ is a undecomposable connected and simply connected $\MD_4$-group 
and ${\cal F}_G$ is the foliation, formed by all the orbits of maximal 
dimension, $V_G := \bigcup_{\Omega\in {\cal F}_G} \Omega$. Then 
$(V_G,{\cal F}_G)$ is a measurable foliation, called the associated 
$\MD_4$-foliation. 
\end{thm}
\begin{pf}
We prove the theorem in two steps. 
\begin{enumerate}
\item[Step 1.]
Find out the integrable tangent distribution, also denoted by ${\cal 
F}_G$ on $V_G$, 
having co-adjoint orbits as maximal connected integral sub-manifolds.

\item[Step 2.] 
Equip to each $(V_G,{\cal F}_G)$ a transversal measure.

\end{enumerate}

For the first step, we find out the differential system $S_G$ defining 
our distributions. 
Following is the list of differential systems by which we choose:

Case 1.1. $$\left\{ \begin{array}{rl} X_1(x,y,z,t) &= (x,0,0,0)\\ 
X_2(x,y,z,t) &= 0,0,0,-z) \end{array}   \right. \leqno{{\cal S}_{1,1}:}$$ 
on the manifold $V_{G_{4,1,1}} = {\RR}^2 \times {\RR}^* \times {\RR}.$

Case 1.2. $$\left\{ \begin{array}{rl} X_1(x,y,z,t) &= (0,0,z,0)\\ 
X_2(x,y,z,t) &= (0, 0,0,-z) \end{array} \right. \leqno{{\cal S}_{1,2}:}$$
on the manifold $V_{G_{4,1,2}} = {\RR}^2 \times {\RR}^* \times {\RR}.$

Case 2.1. $$\left\{ \begin{array}{rl} X_1(x,y,z,t) &= (0,\lambda y, z, 0)\\
X_2(x,y,z,t) &= (-\lambda y,0,0,0)\\ X_3(x,y,z,t) &= (-z,0,0,0)   
\end{array} \right. \leqno{{\cal S}_{4,2,1(\lambda)}:}$$
on the manifold $V_{G_{4,2,1}} = {\RR} \times ({\RR}^*)^2  \times {\RR}$.

Case 2.2. $$\left\{ \begin{array}{rl} X_1(x,y,z,t) &= (0,y,y+z,0)\\ 
X_2(x,y,z,t) &= (-y,0,0,0)\\ X_3(x,y,z,t) &= (-(y+z),0,0,0) \end{array}   
\right. \leqno{{\cal S}_{2,2}:}$$ on the manifold $V_{G_{4,2,2}} = 
{\RR}\times ({\RR}^*)^2 \times {\RR}$.

Case 2.3. $$\left\{ \begin{array}{rl} X_1(x,y+iz,t) &= 
(0,(y+iz)e^{i\varphi},0)\\ X_2(x,y+iz,t) &= -y\cos\varphi + 
z\sin\varphi,0,0)\\ X_3(x,y+iz,t) &= (-y\sin\varphi-z\cos\varphi, 0,0)   
\end{array}  \right.\leqno{{\cal S}_{2,3(\varphi)}:}$$ on the manifold 
$V_{G_{4,2,3}} = {\RR} \times {\CC}^* \times{\RR}$.

Case 2.4. $$\left\{ \begin{array}{rl} X_1(x,y,z,t) &= (0,0,0,1)\\ 
X_2(x,y,z,t) &= (1,0,0,0)\\ X_3(x,y,z,t) &= (0,y,z,0)\\ X_4(x,y,z,t) &= 
(0,-z,y,0) \end{array} \right. \leqno{{\cal S}_{2,4}:}$$ on 
the manifold $V_{G_{4,2,4}} = {\RR} \times ({\RR}^2)^* \times {\RR}$.

Case 3.1. $$\left\{ \begin{array}{rl} X_1(x,y,z,t) &= (\lambda_1 x, 
\lambda_2 y, z, 0)\\ X_2(x,y,z,t) &= 0,0,0, -\lambda_1 x)\\ X_3(x,y,z,t) 
&= (0,0,0,-\lambda_2 y),\\ X_4(x,y,z,t)  &= (0,0,0, -z) \end{array} 
\right. \leqno{{\cal S}_{3,1(\lambda_1,\lambda_2)} (\lambda_1, \lambda_2 
\in {\RR}^*):}$$ on the manifold 
$V_{G_{4,3,1(\lambda_1,\lambda_2)}} = ({\RR}^3)^* \times {\RR}$.

Case 3.2. $$\left\{ \begin{array}{rl} X_1(x,y,z,t) &= \lambda x, x+ 
\lambda y, z, 0)\\ X_2(x,y,z,t) &= (0,0,0, -\lambda x),\\ X_3(x,y,z,t) &= 
(0,0,0,-x -\lambda y),\\ X_4(x,y,z,t) &= (0,0,0,-z) \end{array} \right. 
\leqno{{\cal S}_{3,2(\lambda)} (\lambda\in {\RR}^*):}$$ on the manifold 
$(V_{G_{4,3,2(\lambda)}} = {\RR}^3)^* \times {\RR}$.

Case 3.3. $$\left\{ \begin{array}{rl} X_1(x,y,z,t) &= (x,x+y, y+z, 0),\\ 
X_2(x,y,z,t) &= (0,0,0,-x), \\ X_3(x,y,z,t) &= (0,0,0, -x-y),\\ 
X_4(x,y,z,t) &= (0,0,0,-y-z) \end{array} \right. \leqno{{\cal 
S}_{3,3}:}$$ on the manifold $V_{G_{4,3,3}}= ({\RR}^3)^* \times {\RR}$.

Case 3.4. $$\left\{ \begin{array}{rl} X_1(x+iy,z,t) &= 
((x+iy)e^{i\varphi}, \lambda z, 0)\\ X_2(x+iy,z,t) &= (0,0,-x\cos\varphi 
+ y\sin\varphi)\\ X_3(x+iy,z,t) &= 0,0,-x\sin\varphi - y\cos\varphi)\\ 
X_4(x+iy,z,t) &= (0,0,-\lambda z) \end{array} \right. \leqno{{\cal 
S}_{3,4(\lambda,\varphi)} (\lambda \in {\RR}^*, \varphi\in (0,\pi)) :}$$ 
on the manifold $V_{G_{4,3,4(\lambda,\varphi)}} = ({\CC}\times {\RR})^* 
\times {\RR}$.

Case 4.1. $$\left\{ \begin{array}{rl} X_1(x,y,z,t) &= (-y,x,0,0)\\ 
X_2(x,y,z,t) &= (0,z,0,y)\\ X_3(x,y,z,t) &= (-z, 0,0,-y) \end{array} 
\right. \leqno{{\cal S}_{4,1}:}$$ on th manifold $V_{G_{4,4,1}} = 
({\RR}^3)^* \times {\RR}$.

Case 4.2. $$\left\{ \begin{array}{rl} X_1(x,y,z,t) &= (-x,y,0,0)\\ 
X_2(x,y,z,t) &= (0,z,0,x)\\ X_3(x,y,z,t) &= (-z,0,0,-y) \end{array} 
\right. \leqno{{\cal S}_{4,4,2}:}$$ on the manifold $V_{G_{4,4,2}} = 
({\RR}^3)^* \times {\RR}$.

It is easy to verify that:
\begin{itemize} 
\item All the indicated differential systems are of 
rank 2, but the system ${\cal S}_{2,4}$ of rank 4. 
\item Each co-adjoint orbit $\Omega$ from ${\cal F}_G$ is a maximal 
connected integral manifold of the tangent distribution, generated by the 
corresponding system ${\cal S}_G$. 
\end{itemize}
Thus we have foliation $(V_G,{\cal 
F}_G)$ for each undecomposable connected and simply connected 
$\MD_4$-group. 

To realize the second step, we show that our foliations are orientable in 
finding out a non-vanishing multi-vector of maximal degree $X_G \in 
C^\infty(\wedge^{\dim {\cal F}_G} {\cal F}_G)$ for each case of $G$:
$$\begin{array}{rl}
X_{1,1} &:= X_1 \wedge X_2,\\
X_{1,2} &:= X_1 \wedge X_2, \\
X_{2,1(\lambda\in {\RR}^*)} &:= X_1 \wedge X_2 + X_1 \wedge X_3,\\
X_{2,2} &:= X_1 \wedge X_2 + X_1 \wedge X_3,\\
X_{2,3(\varphi\in (0,\pi))} &:= X_1 \wedge X_2 + X_1 \wedge X_3,\\
X_{2,4} &:= X_1 \wedge X_2 \wedge X_3 \wedge X_4,\\
X_{3,1(\lambda_1,\lambda_2 \in {\RR}^*)} &:= X_1 \wedge X_2 + X_1 \wedge 
X_3 + X_1 \wedge X_4,\\
X_{3,2(\lambda\in {\RR}^*)} &:= X_1 \wedge X_2 + X_1 \wedge X_3 + X_1 
\wedge X_4,\\ 
X_{3,3} &:= X_1 \wedge X_2 + X_1 \wedge X_3 + X_1 \wedge X_4,\\
X_{3,4(\lambda\in{\RR}^*,\varphi\in(0,\pi))} &:= X_1 \wedge X_2 + X_1 
\wedge X_3 + X_1 \wedge X_4,\\ 
X_{4,1} &:= X_1 \wedge X_2 + X_1 \wedge X_3 + X_2 \wedge X_3,\\
X_{4,2} &:= X_1 \wedge X_2 + X_1 \wedge X_3 + X_2 \wedge X_3
\end{array}$$
It is easy also to verify that these multi-vectors are invariant with 
respect to the Lebesgue measure $\mu$. It is just equivalent to its 
invariance w. r. t. the co-adjoint representation of $G$ in ${\frakt 
g}^*$. This means that $(X_G, \mu)$ is an invariant pair.

\end{pf}

\subsection{Topological classification of $\MD_4$-foliations}

Let us recall that two foliations have the same topological type iff 
there exists a leaf-wise homeomorphism between them.

\begin{thm}[\bf Topological Classification]

1.There are exactly 9 topological type of foliations: ${\cal F}_1 - {\cal 
F}_9$ 
$$ (V_{G_{4,1,1}}, {\cal F}_{1,1})\leqno{({\cal F}_1)}$$
$$(V_{G_{4,1,2}},{\cal F}_{1,2})\leqno{({\cal F}_2)}$$
$$\begin{array}{rl}(V_{G_{4,2,1(\lambda)}},{\cal F}_{2,1(\lambda)})(\lambda 
\in {\RR}^*) &\cong (V_{G_{4,2,1(\lambda=1)}},{\cal F}_{2,1(\lambda 
=1)} \\  &\cong  (V_{G_{4,2,2}},{\cal F}_{2,2}) \end{array} \leqno{({\cal 
F}_3)}$$ 
$$ (V_{G_{4,2,3(\varphi)}},{\cal F}_{2,3(\varphi)})(\varphi\in 
(0,\pi))\cong (V_{G_{4,2,3(\pi/2)}},{\cal F}_{2,3(\pi/2)}\leqno{({\cal 
F}_4)}$$ 
$$(V_{G_{4,2,4}},{\cal F}_{2,4})\leqno{({\cal F}_5)}$$
$$ \left\{\begin{array}{rl} (V_{G_{4,3,1(\lambda_1,\lambda_2)}},{\cal 
F}_{3,1(\lambda_1,\lambda_2)} ) &\cong (V_{G_{4,3,1(1,1)}},{\cal 
F}_{4,3,1(1,1)})\\ \cong (V_{G_{4,3,2(\lambda)}},{\cal F}_{3,2(\lambda)}) 
& \cong (V_{G_{4,3,2(1)}}, {\cal F}_{3,2(1)}) \\  &\cong 
(V_{G_{4,3,3}},{\cal F}_{3,3}) \end{array}\right. \leqno{({\cal F}_6)}$$

$$ (V_{G_{4,3,4(\lambda,\varphi)}},{\cal 
F}_{3,4(\lambda,\varphi)}) \cong (V_{G_{4,3,4(1,\pi/2)}},{\cal 
F}_{3,4(1,\pi/2)}) \leqno{({\cal F}_7)}$$
$$(V_{G_{4,4,1}},{\cal F}_{4,1}) \leqno{({\cal F}_8)}$$
$$(V_{G_{4,4,2}},{\cal F}_{4,2}) \leqno{({\cal F}_9)}$$

2. The $\MD_4$-foliations of type ${\cal F}_1 -{\cal F}_6$ are given by 
fibration (trivial with connected fibers), over the bases ${\RR} \times 
{\RR}^*, {\RR}^2 \cup {\RR}^2, {\RR} \times {\SS}^1 , {\RR}_+ \times 
{\RR}, \{ pt \}, {\SS}^2$, resp. , where $\{ pt \}$ is a one-point set.

3. The $\MD$-foliation of type ${\cal F}_7, {\cal F}_8, {\cal F}_9$ are 
given by continuous actions of commutative Lie group ${\RR}^2$ on 
foliated manifolds $({\CC} \times {\RR})^*, ({\RR}^3)^* \times {\RR}, 
({\RR}^3)^* \times {\RR},$ respectively .

\end{thm}
\begin{pf} 

1.

Consider the maps $$h_{2,1(\lambda)} : V_{G_{4,2,1(\lambda)}} 
\approx {\RR} \times ({\RR}^2)^* \times {\RR} \to V_{G_{4,2,1(1)}} = {\RR} 
\times ({\RR}^2)^* \times {\RR},$$ defined by
$$h_{2,1(\lambda)}(x,y,z,t) := (x, \sgn(y) \vert y\vert^{1\over\lambda}, 
z, t),$$ and $$h_{2,2} : V_{G_{4,2,2}}\approx {\RR}\times ({\RR}^2)^* \times 
{\RR} \to V_{G_{4,2,1(1)}} 
\approx {\RR} \times ({\RR}^2)^* \times {\RR},$$ defined by
$$h_{2,2}(x,y,z,t) = \left\{ \begin{array}{ll} (x,y,z-y\ln |y|,t) &\mbox{ 
if } y \ne 0\\ (x,0,z,t) &\mbox { otherwise }, y = 0. \end{array} \right.$$
It is easy to see that these maps are just the leaf-wise homeomorphisms. 
By the same way we construct the maps which realize the leaf-wise 
homeomorphisms: 
$$h_{2,3(\varphi)} : V_{G_{4,2,3(\varphi)}} \approx {\RR} \times {\CC}^* 
\times {\RR} \to V_{G_{4,2,3(\pi/2)}} \approx {\RR}\times {\CC}\times 
{\RR},$$ defined by the formula 
$$h_{2,3(\varphi)}(x,re^{i\theta},t) = (x, e^{(\ln r + i\theta) i 
e^{-i\varphi}},t),$$
$$h_{3,1(\lambda_1,\lambda_2)} : V_{G_{4,3,1(\lambda_1,\lambda_2)}} 
\approx ({\RR}^3)^* \times {\RR} \to V_{G_{3,1(1,1)}} \approx ({\RR}^3)^* 
\times {\RR},$$ defined by the formula 
$$h_{3,1(\lambda_1,\lambda_2)}(x,y,z,t) = (\sgn(x) |x|^{\lambda_1}, 
\sgn(y)|y|^{\lambda_2}, z,t),$$ 
$$h_{3,2(\lambda)} : V_{G_{4,3,2(\lambda)}} \approx ({\RR}^3)^* 
\times 
{\RR} \to V_{G_{4,3,2(1)}} = ({\RR}^3)^* \times {\RR},$$ defined by the 
formula $$h_{3,2)(\lambda)}(x,y,z,t) = (\tilde{x}, \tilde{y}, \tilde{z}, 
\tilde{t}),$$ where $$\begin{array}{ll}
\tilde{x} &= \sgn(x) |x|^{1\over\lambda},\\
\tilde{y} &= \left\{ \begin{array}{ll} \sgn(y-{1\over\lambda} x\ln|x|).|y - 
x\ln |x||^{1\over \lambda} & \mbox{ if } x \ne 0,\\  \sgn(y) |y|^{1\over 
\lambda} & \mbox{ if } x = 0\end{array} \right. \\
\tilde{z} &= z,\\
\tilde{t} &= t, \end{array}$$
$$h_{3,3}(x,y,z,t) = (\tilde{x}, \tilde{y}, \tilde{z}, \tilde{t}),$$ 
where $$\begin{array}{rl}
\tilde{x} &= x,\\
\tilde{y} &= {\left\{ \begin{array}{ll} y -x\ln |x|, &\mbox { if } x\ne 0, 
\\ y &\mbox{ if } x = 0,\end{array}\right. }\\
\tilde{z} &= {\cases z - {1\over 2}y\ln|x| - {1\over 2}(y-x\ln|x|)\ln|y - 
x\ln|x||, &\mbox { if } x \ne 0, y \ne x\ln|x|,\\  
z - {1\over 2} y\ln|x|, &\mbox{ if } x \ne 0, y = x\ln|x|,\\
z &\mbox{ if } x = 0,\endcases} \\
\tilde{t} &= t \end{array}$$
$$h_{3,4(\lambda,\varphi)} : V_{G_{4,3,4(\lambda,\varphi)}} \approx 
({\CC} \times {\RR})^* \times {\RR} \to V_{G_{4,3,4(1,\pi/2)}} \approx 
({\CC}\times {\RR})^* \times {\RR},$$ defined by the formula 
$$h_{3,4(\lambda,\varphi)}(re^{i\theta},z,t) = (e^{\ln r + 
i\theta)ie^{-i\theta}}, \sgn(z)|z|^{1\over \lambda}, t).$$

2. 

From Theorem 3.1 on the structure of co-adjoint orbits, it is easy to see 
that the type ${\cal F}_1 - {\cal F}_9$ are non topologically equivalent 
and that the foliations of type ${\cal F}_1, {\cal F}_2$ and $ {\cal F}_5$ 
are trivial fibration over the bases ${\RR} \times {\RR}^*, {\RR}^2 
\cup {\RR}^2$ and $\{ pt \}$, respectively. 
It is easy to see that the foliations of type ${\cal F}_3, {\cal F}_4, 
{\cal F}_6$ are fibrations defined by the submersions
$$p_{2,1(1)} : V_{G_{4,2,1(1)}} \approx {\RR} \times ({\RR}^2)^* \times 
{\RR} \cong {\RR} \times {\SS}^1 \times {\RR}_+ \times {\RR} \to {\RR} 
\times {\SS}^1,$$ the projection on the first two components,
$$p_{2,3(\pi/2)} : V_{G_{4,2,3(\pi/2)}} \approx {\RR} \times {\CC}^* \times 
{\RR} \to {\RR}_+ \times {\RR},$$ the projection on the first and radian 
part of the second component,
$$p_{3,1(1,1)} : V_{G_{4,3,1(1,1)}} \approx ({\RR}^3)^* \times {\RR} 
\approx {\SS}^2 \times {\RR}_+ \times {\RR} \to {\SS}^2,$$ the projection 
on the first component. 

3.  

Let us defines the actions of ${\RR}^2$ on foliated manifolds 
$V_{G_{4,3,4}} \approx V_{G_{4,4,1}} \approx V_{G_{4,4,2}}$ as following
$$\rho_{3,4} : {\RR}^2 \times V_{G_{4,3,4}}\approx {\RR}^2 \times ({\CC} 
\times {\RR})^* \times {\RR} \to V_{G_{4,3,4}},$$ 
$$\rho_{3,4}((r,s),(x+iy,z,t) = (x+iy)e^{is}, ze^s, t+r);$$
$$\rho_{4,1} : V_{G_{4,4,1}} \approx ({\RR}^3)^* \times {\RR} \to 
V_{G_{4,4,1}},$$ $$\rho_{4,1}((r,s),(x,y,z,t)) = (\tilde{x}, \tilde{y}, 
\tilde{z},\tilde{t}),$$ where 
$$\begin{array}{rl}
\tilde{x} &= x\cos r - y\sin r - sz,\\
\tilde{y} &= x\sin r + y\cos r -sz,\\
\tilde{z} &= z,\\
\tilde{t} &= t -s(x+y) \cos r + s(y-x) \sin r + s^2 z;\end{array}$$
$$\rho_{4,2} : {\RR}^2 \times V_{G_{4,4,2}} \approx {\RR}^2 \times 
({\RR}^3)^* \times {\RR} \to V_{G_{4,4,2}},$$
$$\rho_{4,2}((r,s),(x,y,z,t)) = (\tilde{x}, \tilde{y}, \tilde{z}, 
\tilde{t}),$$ where $$\begin{array}{rl}
\tilde{x} &= e^{-s}(x+r){yz\over x^2 + y^2+ z^2},\\
\tilde{y} &= e^{-s}(y+r){xz\over x^2+y^2+z^2},\\
\tilde{z} &= z,\\
\tilde{t} &= t + r{x^2+y^2\over x^2+y^2+z^2} + r^2{xyz\over 
x^2+y^2+z^2}.\end{array}$$
It is easy to check that these actions give us just the foliations of 
type ${\cal F}_7, {\cal F}_8, {\cal F}_9$.  

\end{pf}

\section{Bibliographical Remarks} 

Dao van Tra had first proved the proposition 2.1 by using the old method
of classification of low dimension Lie algebras. The classification Theorem
2.1, Theorem 3.1 about the picture of K-orbits and the Theorem about the
measurability of foliations of generic K-orbits belong to Le Anh Vu, who
was a Ph.D. student under author's supervision. The problem of
classification of Lie $\MD_n$-algebras, for $n \geq 5$ rests open
up-to-date.

\chapter{The Structure of C*-Algebras of $\MD_4$-Foliations}
%
\section{C*-Algebras of Measurable Foliations}
We recall in this section the well-known A. Connes' theory of C*-algebras 
of measurable foliations.

\subsection{Holonomy group of foliation}
Let us consider a foliated manifold $(V,{\cal F})$, $\dim V = n, \dim 
{\cal F} = q$. Holonomy groupoid $H$ is a manifold (not necessarily 
Hausdorff) of dimension $\dim H = \dim V + \dim {\cal F}$ (see e.g. 
\cite{winkelkemper}, \cite{connes1}, \cite{torpe}): An element $\gamma$ 
of $H$ is just done by a pair of source $s(\gamma)\in V$ and rival 
$r(\gamma)$ in the same leaf of $(V,{\cal F})$ as $s(\gamma)$, together 
with an homotopy equivalence class $[\gamma]$ of path connecting the same 
source and rival, lying in the same leaf. 
The product of two elements $(s(\gamma),r(\gamma),[\gamma])$ and 
$(s(\delta), r(\delta),[\delta])$ can be defined iff $r(\gamma) = 
s(\delta)$ and in this case $\delta\circ\gamma$ is defined as usually,
$$(s(\gamma),r(\gamma),[\gamma]) \circ (s(\delta) = r(\gamma), 
r(\delta),[\delta]) := (s(\gamma),r(\delta), [\delta\circ\gamma]).$$ 
The topology of $H$ is defined by the system of local chart (atlas) 
$$\Gamma := \{([\gamma],s(\gamma)=x, r(\gamma)=y); x\in U, y\in L_x  
\}.$$ The first component should be discrete, for fixed $x$ in a 
neighborhood $U$ and $y$ in the same leaf as $x$. Thus $\dim H = \dim V + 
\dim {\cal F}$. It is known:
\begin{prop}[\cite{torpe},Prop.2.1] The groupoid $H$ of the foliation 
$(V,{\cal F})$ iff for each pair of point $x$, $y$ in a leaf $L$, such 
that for each pair of smooth paths $\gamma_1,\gamma_2$ in $L$, connecting 
$x$ and $y$ the holonomy maps $h(\gamma_1) = h(\gamma)_2)$ if they are 
coincided in a small enough open neighborhood, the closure of which 
contains $x$. 
\end{prop}
 
\subsection{Half-density bundle}

Let us consider an oriented  $k$-dimensional foliation ${\cal F}$ in a 
$n$-dimensional smooth manifold $V$. For each point $x\in V$, consider 
the set of so called half-densities 
$$\Omega_x^{1/2} := \{ \rho: \wedge^k {\cal F}_x \to {\CC}; \rho(\lambda 
v) = \vert\lambda\vert^{1/2}\rho(v), \forall v\in \wedge^k{\cal F}_x, 
\forall \lambda \in {\RR}   \}.$$  
It is easy to see that $\Omega_x^{1/2}$ is a complex 1-dimensional space 
and $\bigcup_{x\in V} \Omega_x^{1/2}$ is a complex 1-dimensional fiber 
bundle over $V$, called the {\it bundle of half-densities}. 
\index{bundle of half-densities} Because we 
assumed that ${\cal F}$ is oriented, the bundle of half-densities is 
trivial. Fix a trivialization, we have $\bigcup_{x\in V} \Omega_x^{1/2} 
\cong V \times {\CC}$. For each element $\gamma\in H$, pose 
$\Omega_\gamma := \Omega_x^{1/2} \otimes \Omega_y^{1/2},$ where $x= 
s(\gamma)$ and $y = r(\gamma)$. Thus $\Omega_\gamma^{1/2}$ is also a 
complex 1-dimensional vector space ${\CC}$. 

If $H$ is Hausdorff, we define
$$C^\infty_c(H,\Omega^{1/2}) = \{f : \gamma\in H \to \Omega^{1/2}_\gamma; 
f \mbox{ is smooth and with compact support }  \},$$ 
each element of $C^\infty_c(H,\Omega^{1/2})$ is called a {\it 
smooth half-density with compact support}. \index{smooth half-densities 
with compact support} If $H$ is non Hausdorff, we 
define $C^\infty_c(H,\Omega^{1/2})$ as the set of finite combinations of 
smooth functions of form $\varphi \circ \chi$, where $\chi : \Gamma \to 
{\RR}^{n+k}$ is a local coordinate chart of $H$, and $\varphi \in 
C^\infty_c({\RR}^{n+k},\Omega^{1/2})$ such that $\supp\varphi\in 
\chi(\Gamma)$.

In virtue of trivialization $$\CD \nu : \bigcup_{x\in V}\Omega^{1/2}_x 
@>\cong>> V \times {\CC},\endCD$$ we can identify smooth half-density with 
compact support $f \in C^\infty_c(H,\Omega^{1/2})$ with ${\CC}$-valued 
smooth function with  compact support $f(\nu .s \otimes r .s)$ on $H$.

\subsection{C*-algebras of measurable foliations}

Let us now define \cite{connes1} the C*-algebras of measurable 
foliations.  Define a convolution product on $C^\infty_c(H,\Omega^{1/2})$ as
$$(f * g)(\gamma) = \int_{\gamma_1\circ \gamma_2 = \gamma} 
f(\gamma_1)g(\gamma_2), \forall f,g\in C^\infty_c(H,\Omega^{1/2}), 
\forall \gamma\in H.$$ Define involution $f \mapsto f^*$, $f \in 
C^\infty_c(H,\Omega^{1/2})$ by the formula 
$$f^*(\gamma) := \overline{f(\gamma^{-1})}.$$ It is easy to see that with 
this convolution product and involution $C^\infty_c(H,\Omega^{1/2})$ is 
an involutive algebra.

For each $x\in V$, define $H_x := \{ \gamma\in H; s(\gamma) = x\}$. There 
is a natural representation $\pi_x$ of $C^\infty_c(H,\Omega^{1/2})$ in 
$L^2(H_x,\Omega^{1/2})$, given by the following formula
$$(\pi_x(f)\eta)(\gamma) = \int_{\gamma_1\circ\gamma_2= \gamma} 
f(\gamma_1)\eta(\gamma_2), \forall f\in C^\infty_c(H.\Omega^{1/2}), 
\forall \eta\in L^2(H_x,\Omega^{1/2}),\forall \gamma\in H.$$

\begin{defn}
The C*-algebra of foliation $C^*(V,{\cal F})$ is defined as C*-hull of 
*-algebra $C^\infty_c(H,\Omega^{1/2})$, with respect to the norm of operators
$$\Vert f\Vert = \sup_{x\in V} \Vert \pi_x(f)\Vert, \forall f\in 
C^\infty_c(H,\Omega^{1/2}).$$
\end{defn}

Recall that a C*-algebra $A$ is called {\it stable} \index{C*-algebra! 
stable -} iff $A \cong A \otimes {\cal 
K}(H)$, where ${\cal K}(H)$ denotes the ideal of compact operators in a 
separable Hilbert space $H$. 

\begin{prop}[\cite{torpe},Prop. 2.1.4] If the foliations $(V,{\cal F})$ 
and $(V',{\cal F}')$ are of one topological type, then the corresponding 
C*-algebras are isomorphic $$C^*(V,{\cal F}) \cong C^*(V',{\cal F}').$$
\end{prop}

Let us recall that cross-product $A\rtimes_\rho G$ of a group $G$ with a 
C*-algebra $A$ is just the 
C*-hull of the involutive algebra of functions with compact support on 
$G$ with values in $A$, endowed with the convolution product
$$(a*b)(g) := \int_G a(g_1)\rho_{g_1}b(gg_1^{-1}) dg_1, \forall a,b \in 
L^1(G,A),\forall g \in G$$
and the involution $$a^*(g) := \rho_g(a(g^{-1})^*), \forall a\in 
L^1(G,A).$$ where $dg$ is the right-invariant Haar measure on $G$ and 
$\rho : G \to \Aut A$ is a representation of $G$, i.e. the map $G \times A 
\to A,$ $(g, x) \mapsto \rho(g)x$ is continuous in norm.

It is easy to see that if $f: A \to A'$ is a $G$-equivariant 
C*-homomorphism, then $f$ induces also a C*-morphism of crossed products
$$f_\# : A \rtimes_\rho G \to A'\rtimes_{\rho'} G',$$
$$f_\#(a)(g) := f(a(g)), \forall a\in L^1(G,A), \forall g\in G.$$

Let us list now some basic properties of C*-algebras of foliations.

\begin{prop}[\cite{connes3},Lem I.1]
If the sequence $$\CD 0 @>>> J @>>> A @>>> B @>>> 0\endCD$$ is an  
$G$-equivariant short exact sequence, then the corresponding sequence of 
crossed products $$\CD 0 @>>> J \rtimes G @>>> A\rtimes G @>>> B\rtimes G 
@>>> 0\endCD$$ is also exact. If the first sequence is split then the 
same is the second. 
\end{prop}

\begin{prop}[\cite{connes1},\S5]
Assume that the foliation $(V,{\cal F})$ is given by an action of Lie 
group $G$ on the manifold $V$, such that the holonomy group $H$ of 
$(V,{\cal F})$ is exactly of form $H = V \times G$. Then $C^*(V,{\cal F}) 
= C_0(V) \rtimes G$, where $C_0(V)$ is the algebra of continuous functions 
vanishing at infinity
\end{prop}

\begin{prop}[\cite{connes1},\S5]
If $V'$ is some open set of a foliated manifold $(V,{\cal F})$ and ${\cal 
F}' := {\cal F}|_{V'}$. Then the holonomy groupoid $H'$ of foliation 
$(V',{\cal F}')$ is an open subset of the holonomy groupoid $H$ of 
$(V,{\cal F})$. Moreover, the inclusion $C^\infty_c(H',\Omega^{1/2}) 
\hookrightarrow C^\infty_c(H,\Omega^{1/2})$ can be extended to an 
*-homomorphism $\imath : C^*(V',{\cal F}') \to C^*(V,{\cal F})$.
\end{prop}

\begin{prop}[\cite{connes2},\S5]
Assume that foliation $(V,{\cal F})$ is given by a fiber bundle $p : V 
\to M$ with connected fibers. Then the holonomy groupoid $H$ is just the 
sub-manifold $\{(x,y) \in V \times V; p(x) = p(y)\}$ of $V \times V$ and 
$C^*(V,{\cal F})\cong C_0(M) \otimes {\cal K}(L^2(\mbox{typical fiber}))$.
\end{prop}

The open sub-manifold $V'$ of foliated manifold $(V,{\cal F})$ is said to 
be {\it satured}, \index{satured foliated manifold} iff it contains the 
whole 
leaves, which intersect (i.e. have a non empty intersection) with $V'$. 

\begin{prop}[\cite{torpe},\S2.2]If $V'$ is a satured open sub-manifold of 
the foliated manifold $(V,{\cal F})$ and ${\cal F}' := {\cal F}|_{V'}$ 
then $C^*(V',{\cal F}')$ is an ideal in $C^*(V,{\cal F})$.
\end{prop}
 
Remark that in this case, the groupoid $H\setminus H'$ is closed in $H$ 
but in general, it is different from the groupoid of the foliation 
$V\setminus V',{\cal F}|_{V \setminus V'})$. Nevertheless, we can also 
define the representation $\pi_x (x\in V\setminus V')$ of the *-algebra 
$C^\infty_c(H\setminus H',\Omega^{1/2})$ in $L^2(H_x \setminus 
H'_x,\Omega^{1/2})$. The C*-hull of it is a C*-algebra, denoted by $C^*(V 
\setminus V', {\cal F}|_{V\setminus V'})$.

The inclusion $H \setminus H' \hookrightarrow H$ gives us a 
*-homomorphism $$\mu' : C^\infty)_c(H, \Omega^{1/2}) \to C^\infty_c(H 
\setminus H', \Omega^{1/2}),$$ which can be extended to a *-epimorphism
$$\CD C^*(V,{\cal F}) @>\mu >> C^*(V\setminus V', {\cal F}|_{V\setminus V'}) 
@>>> 0.\endCD$$ We have thus a short sequence
$$\CD 0 @>>> C^*(V',{\cal F}') @>\imath>> C^(V,{\cal F}) @>\mu >> 
C^(V\setminus V',{\cal F}|_{V\setminus V'}) @>>> 0, \endCD$$  which is 
exact at all terms, but perhaps,is not in the middle term $C^*(V,{\cal 
F})$. 

\begin{prop}[\cite{torpe},Lemma 2.2.1]
If the foliation $(V,{\cal F})$ is given by an action of an amenable Lie 
group $G$, such that $H\setminus H' = (V\setminus V') \times G$ then 
the previous sequence is exact.
\end{prop}
 
\subsection{Connes-Thom isomorphism}
\begin{prop}[\bf Connes-Thom Isomorphism, \cite{connes3}, Thm, IV.2]
Assume that the commutative group ${\RR}^n$ acts continuously on a
C*-algebra $A$ by an action by a continuous representation $\rho$. Then
there is a natural isomorphism
$$\CD \varphi^j_\rho : K_j(A) @>\cong>> K_{j+n(\mod 2)}(A \rtimes_\rho
{\RR}^n , \endCD$$ where $j \in {\ZZ}/2{\ZZ}$.
\end{prop}

We return to the case of extensions of type $$\CD 0 @>>> J @>>> E @>>> A
@>>> 0. \endCD$$ Assume that the commutative group ${\RR}^n$ acts
continuously on each C*-algebras of the sequence and the *-homomorphisms
are equivariant. We have then also an exact sequence of the
corresponding cross-product
$$\CD 0 @>>> J\rtimes G @>>> E\rtimes G @>>> A \rtimes G @>>>0 \endCD.$$
For each of these two exact sequence, we have 6-term exact sequences of
K-groups, in one hand and in other hand we have Connes-Thom isomorphisms
between the corresponding terms. This should be very useful in computing
connecting homomorphisms in many cases.

It is easy to see that if the foliation $(V,{\cal F})$ is given by an
action $\rho$ of some commutative Lie group ${\RR}^n$ in such a way that
its holonomy groupoid is $H = V \times {\RR}^n$, then $C^*(V,{\cal F})
\cong C_0(V) \rtimes_\rho {\RR}^n$ and hence
$$K_j(C_(V)) \cong K_{j+n\mod 2}(C^*(V,{\cal F}) = K^{j+n \mod
2}(V/F).$$


\section{The C*-Algebras of Measurable $\MD_4$-Foliations  }
 
We complete studying the structure of C*-algebras of $\MD_4$-groups in 
this section in the same way as in the previous sections. 
\subsection{C*-algebras of $\MD_4$-foliations of bundle type}
Following the theorem on topological type of $\MD_4$-foliations, the
foliations of type ${\cal F}_1 - {\cal F}_6$ are of bundle type with
connected fibers. Then as an easy corollary of the above result, we have
\begin{prop}
\begin{enumerate}
\item $$C^*(V_{G_{4,1,1}},{\cal F}_{1,1}) \cong C_0({\RR} \times {\RR}^*)
\otimes {\cal K}.$$
$$C^*(V_{G_{4,1,2}},{\cal F}_{1,2}) \cong C_0({\RR}^2 \cup {\RR}^2)
\otimes {\cal K} \cong (C_0({\RR}^2) \oplus C_0({\RR}^2)) \otimes {\cal
K}.$$
\item $$C^*(V_{G_{4,2,1(\lambda)}},{\cal F}_{2,1(\lambda)}) \cong
C^*(V_{G_{4,2,2}},{\cal F}_{2,2}) \cong C_0({\RR} \times {\SS}^1)
\otimes {\cal K}, \forall \lambda\in {\RR}.$$
$$C^*(V_{G_{,2,3(\varphi)}}, {\cal F}_{2,3(\varphi)}) \cong C_0({\RR}_+
\times {\RR}) \otimes {\cal K}, \forall \varphi\in (0,\pi),$$
$$C^*(V_{G_{4,2,4}}, {\cal F}{2,4}) \cong {\CC} \otimes {\cal K} \cong
{\cal K}.$$
\item $$C^*(V_{G_{4,3,1(\lambda_1,\lambda_)}}, {\cal
F}_{3,1(\lambda_1,\lambda_2)}) \cong C^*(V_{G_{4,3,2(\lambda)}}, {\cal
F}_{3,2(\lambda)}) $$ $$\cong C^*(V_{G_{4,3,3}}, {\cal F}_{3,3}) \cong
C({\SS}^2) \otimes {\cal K}, \forall \lambda, \lambda_1, \lambda_2 \in
{\RR}^*.$$
\end{enumerate}
\end{prop}

\subsection{C*-algebras of $\MD_4$-foliations of crossed product type} 

Let us consider the C*-algebras $C^*(V_{G{4,3,4(\lambda,\varphi)}},{\cal
F}_{3,4(\lambda,\varphi)}) \cong C^*(V_{G_{4,3,4(1,\pi/2)}},{\cal
F}_{3,4(1,\pi/2)})$.
\begin{thm}

\begin{enumerate}
\item
The C*-algebra $C^*(V_{G_{4,3,4(1,\pi/2)}},{\cal F}_{3,4(1,\pi/2)})$ can
be included in the extension
$$\CD 0 @>>> J_3 @>>> C^*(V_{G_{4,3,4(1,\pi/2)}}, {\cal F}_{3,4(1,\pi/2)}) @>>> B_3 @>>> 0,\endCD
\leqno{(\gamma_1)}$$  where $J_3 \cong C_0({\RR}^2 \times {\RR}^* \times
{\RR}) \rtimes_{\rho_{3,4}} {\RR}^2 \cong C_0({\RR}^2 \cup {\RR}^2)
\otimes {\cal K},$ $B_3 \cong C_0(({\RR}^2)^* \times {\RR})
\rtimes_{\rho_{3,4}} {\RR}^2 \cong C_0({\RR}_+) \otimes {\cal K}.$
\item 
The C*-algebra $C^*(V_{G_{4,4,1}},{\cal F}_{4,1})$ can be included
in the extension
$$\CD 0 @>>> J_{4,1} @>>> C^*(V_{G_{4,4,1}}, {\cal F}) @>>> B_{4,1} @>>>
0  ,\endCD, \leqno{(\gamma_2)}$$ where $J_{4,1} \cong C_0({\RR}^2
\times {\RR}^* \times {\RR}) \rtimes_{\rho_{4,1}} {\RR}^2 \cong
C_0({\RR}^* \times {\RR}) \otimes {\cal K},$ $B_{4,1} \cong
C_0(({\RR}^2)^* \times {\RR}) \rtimes_{\rho{4,1}} {\RR}^2 \cong
C_0({RR}_+ ) \otimes {\cal K}.$
\item 
The C*-algebra $C^*(V_{G_{4,4,2}}, {\cal F}_{4,2})$ of the diamond
$\MD_4$-foliation  can be included in the following two repeated exact
sequences
$$\CD 0 @>>> J_{4,2,1} @>>> C^*(V_{G_{4,4,2}}, {\cal F}_{4,2}) @>>>
B_{4,2,1} @>>> 0   ,\endCD \leqno{(\gamma_3)}$$
$$\CD 0 @>>> J_{4,2,2} @>>> B_{4,2,1} @>>> B_{4,2,2} @>>> 0  ,\endCD
\leqno{(\gamma_4)}$$ where
$$\begin{array}{rl}
C^*(V_{G_{4,4,2}},{\cal F}_{4,2}) &\cong C_0(V_{G_{4,4,2}}) \rtimes
{\RR}^2 \cong C_0(({\RR}^3)^* \times {\RR}) \rtimes_{\rho_{4,2}}
{\RR}^2,\\
J_{4,2,1} &\cong C_0({\RR}^2 \times {\RR}^* \times {\RR})
\rtimes_{\rho_{4,2}} {\RR}^2 \cong C_0({\RR}^* \times {\RR}) \otimes
{\cal K},\\
B_{4,2,1} &\cong C_0(({\RR}^2)^* \times {\RR})
\rtimes_{\rho_{4,2}} {\RR}^2, \\
J_{4,2,2} &\cong C_0(({\RR}^*)^2 \times {\RR}) \rtimes_{\rho_4,2}
{\RR}^2 \cong C_0({\RR}^* \cup {\RR}^*) \otimes {\cal K},\\
B_{4,2,2} &\cong C_0(({\RR}^* \cup {\RR}^*) \times {\RR})
\rtimes_{\rho_{4,2}} {\RR}^2 \cong {\CC}^4 \otimes {\cal K}.
\end{array}$$
\end{enumerate}
\end{thm}
\begin{pf}
Let us recall that for all the $\MD_4$-foliations of type 
${\cal F}_7, {\cal F}_8, {\cal F}_9$ we have
$$V_{G_{4,3,4(1,\pi/2)}} \cong V_{G_{4,4,1}} \cong V_{G_{4,2,2}} \cong 
({\RR}^3)^* \times {\RR}.$$ Let us consider therefore the following two 
subsets:
$$V = \{(x,y,z,t) \in ({\RR}^3)^* \times {\RR} ; z \ne 0 \} \cong {\RR}^2 
\times {\RR}^* \times {\RR},$$
$$W = \{ ({\RR}^3)^* \times {\RR} \backslash V = ({\RR}^2)^* 
\times \{0\} \times {\RR} \cong ({\RR}^2)^* \times {\RR}.$$ From the 
theorem on the geometric picture of co-adjoint orbits it is easy to see 
that $V$ is an satured open set and w. r. t. the foliations ${\cal F}_7, 
{\cal F}_8, {\cal F}_9$. Denote the restrictions of theses 3 foliation on 
$V$ by $(V,{\cal F}_{3,4(1,\pi/2)}), (V, {\cal F}_{4,1})$, $(V,{\cal 
F}_{4,2}$. By the same way, it is easy to see that $W$ is a satured close 
set (the complement of $V$ in the foliations 
$(V_{G_{4,3,4(1,\pi/2)}}, {\cal F}_7)$, $(V_{G_{4,4,1}},{\cal F}_8)$, 
$(V_{G_{4,4,2}}, {\cal F}_9)$. We denote the restrictions of these 
foliations to $W$ also by $(W, {\cal F}_7)$, $(W,{\cal F}_8)$ and 
$(W,{\cal F}_9)$, respectively. It is easy to see that these foliations
$(V,{\cal F}_7)$, ... , $(W,{\cal F}_9)$ are fiber bundles: 
$$\begin{array}{rl}
p_{3,4} &: V \approx {\RR}^2 \times {\RR}^* \times {\RR}\to {\RR}^2 
\times \{1,-1\} \approx {\RR}^2 \cup {\RR}^2,\\ 
p_{4,1} &: V \approx {\RR}^2 \times {\RR}^* \times {\RR} \to {\RR}^* 
\times {\RR},\\ 
p_{4,2} &: V \approx {\RR}^2 \times {\RR}^* \times {\RR} \to {\RR}^* \times 
{\RR},\\ 
q_{3,4} &: W \approx ({\RR}^2)^* \times {\RR} \approx {\RR}_+ \times 
{\SS}^1 \to {\RR}_+, \\ q_{4,1} &: W \approx ({\RR}^2)^* \times {\RR} 
\approx {\RR}_+ \times {\SS}^1 \times {\RR} \to  {\RR}_+, \end{array}$$
where by definition, $$\begin{array}{rl}
p_{3,4}(x,y,z,t) &:= (x,y,\sgn z),\\
p_{4,1}(x,y,z,t) &:= (z, t- {x^2 + y^2 \over z}),\\
p_{4,2}(x,y,z,t) &:= (z, {t-xy\over z}),\\
q_{3,4}(r, \theta, t) &:= q_{4,1}(r, \theta,t) =\\
	&:= r, \forall (r,\theta,t)\in {\RR}_+ \times {\SS}^1 \times {\RR}.
\end{array}$$.
Following Proposition 1.6, we have
$$\begin{array}{rl}
J_3 &= C^*(V,{\cal F}_{3,4}) \cong C_0({\RR} \times \{\pm\} \otimes {\cal 
K},\\
B_3 &= C^*(W, {\cal F}_{3,4}) \cong C_0({\RR}_+) \otimes {\cal K},\\
J_{4,1} &= C^*(V,{\cal F}_{4,1}) \cong C_0({\RR}^* \times {\RR}) \otimes 
{\cal K},\\
B_{4,1} &= C^*(W,{\cal F}_{4,1}) \cong C_0({\RR}_+) \otimes {\cal K},\\
J_{4,2,1} &= C^*(V, {\cal F}_{4,2}) \cong C_0({\RR}^* \times {\RR}) 
\otimes {\cal K}.
\end{array}$$
Let us recall that the foliations $(V_{G_{3,4(1,\pi/2)}},{\cal 
F}_{3,4(1,\pi/2)})$, $(V_{G_{4,1}},{\cal F}_{4,1})$ and $(V_{G_{4,2}}, 
{\cal F}_{4,2})$ are given by continuous actions $\rho_{3,4}, \rho_{4,1}, 
\rho_{4,2}$ of ${\RR}^2$. These actions conserve $V$ and $W$ in any 
cases, them we can think about $V$ and $W$ in each case as some foliation 
raised from the continuous actions $\rho_{3,4}, \rho_{4,1}, \rho_{4,2}$. 
They have trivial holonomy of fibers. From Prop. 1.4, we have:
$$\begin{array}{rl}
J_3 &= C^*(V, {\cal F}_{3,4(1,\pi/2)}) \cong C^*(V) 
\rtimes_{\rho{3,4}} {\RR}^2 \cong 
C_0({\RR}^2 \times {\RR}^* \times {\RR}) \rtimes_{\rho_{3,4}} {\RR}^2,\\
B_3 &= C^*(W, {\cal F}_{3,4(1,\pi/2)}) \cong C_0(W) \rtimes_{\rho_{3,4}} 
{\RR}^2 \cong C_0(({\RR}^2)^* \times {\RR}) \rtimes_{\rho_{3,4}} {\RR}^2,\\
J_{4,1} &= C^*(V, {\cal F}_{4,1}) \cong C_0(V) \rtimes_{\rho_{4,1}} 
{\RR}^2 \cong C_0({\RR}^2 \times {\RR}^* \times {\RR}) 
\rtimes_{\rho_{4,1}} {\RR}^2,\\
B_{4,1} &= C^*(W, {\cal F}_{4,1}) \cong C_0(W) \rtimes_{\rho_{4,1}} 
{\RR}^2 \cong C_0(({\RR}^2)^* \times {\RR}) \rtimes_{\rho_{4,1}} {\RR}^2,\\
J_{4,2,1} &= C^*(V, {\cal F}_{4,2}) \cong C_0(V) \rtimes_{\rho_{4,2}} 
{\RR}^2 \cong C_0({\RR}^2 \times {\RR}^* \times {\RR}) 
\rtimes_{\rho_{4,2}} {\RR}^2,\\
B_{4,2,1} &= C^*(W,{\cal F}_{4,2}) \cong C_0(W) \rtimes_{\rho_{4,2}} 
{\RR}^2 \cong C_0(({\RR}^2)^* \times {\RR}) \rtimes_{\rho_{4,2}} {\RR}^2.
\end{array}$$
By the same way, we see that the foliation $(V_{G_{3,4(1,\pi/2)}},{\cal 
F}_{3,4(1,\pi/2)})$, $(V_{G_{4,4,1}},{\cal F}_{4,1})$ and 
$(V_{G_{4,4,2}}, {\cal F}_{4,2})$ satisfy the conditions of Props. 1.4 - 
1.6., and we also have
$$\begin{array}{rl}
C^*(V_{G_{4,3,4(1,\pi/2)}} {\cal F}_{3,4(1,\pi/2)}) &\cong 
C_0(V_{G_{4,3,4(1,\pi/2)}}) \rtimes_{\rho_{3,4}} {\RR}^2 \cong 
VC_0(({\RR}^3)^* \times {\RR}) \rtimes_{\rho_{3,4}} {\RR}^2,\\
C^*(V_{G_{4,4,1}}, {\cal F}_{4,1}) &\cong C_0(V_{G_{4,4,1}}) \rtimes 
{\RR}^2 \cong C_0(({\RR}^3)^* \times {\RR}) \rtimes_{\rho{4,1}} {\RR}^2,\\
C^*(V_{G_{4,4,2}},{\cal F}_{4,2}) &\cong C_0(V_{G_{4,4,2}}) 
\rtimes_{\rho_{4,2}} {\RR}^2 \cong C_0(({\RR}^3)^* \times {\RR}) 
\rtimes_{\rho_{4,2}} {\RR}^2.
\end{array}$$
We have therefore the extensions $\gamma_1$, $\gamma_2$ and $\gamma_3$. 
It rests to show existence of the extension $\gamma_4$.
We consider the following sub-manifolds:
$$W_1 := \{(x,y,0,t) \in W ; xy \ne 0   \}\cong ({\RR}^*)^2 
\times {\RR},$$ and 
$$ W_2 = W \backslash W_1 = ({\RR}^* \times \{0\} \times {\RR}) 
\cup (\{0\} \times {\RR}^* \times {\RR} \cong ({\RR}^* \cup 
{\RR}^*) \times {\RR}.$$ It is easy to check that $W_1$, $W_2$ is a satured 
open, resp. closed in the foliated manifold $(W,{\cal F}_{4,2})$ 
and its restrictions to $W_1$ and $W_2$ are the foliations 
$(W_1,{\cal F}_{4,2})$ and $(W_2,{\cal F}_{34,2})$ formed under 
the action $\rho_{4,2}$ of ${\RR}^2$ on $W_1,W_2$, resp. These 
foliations are indeed fiber bundles
$$p_{4,2,2} : W_1 \approx ({\RR}^*)^2 \times {\RR}\to {\RR}^* \cup {\RR}^* 
\leqno{(W_1,{\cal F}_{4,2}):}$$
$$\left\{\begin{array}{l}
q_{4,2,2} : W_2 \approx ({\RR}^* \times \{0\} 
\times {\RR}) \cup (\{0\} 
\times {\RR}^* \times {\RR}) \approx\\ \approx ({\RR}^* \cup {\RR}^*) \times 
{\RR} \to 
\{(1,0),(0,1), (-1,0),(0,-1)\},\end{array}\right. \leqno{(W_2,{\cal 
F}_{4,2}):}$$           defined by the formulae
$$p_{4,2,2}(x,y,t) := (\sgn(x), \sgn(y), \sqrt{|xy|}), \forall (x,y,t) 
\in W_1,$$
$$q_{4,2,2}(x,y,t) := (\sgn(x),\sgn(y)), \forall (x,y,t) \in W_2.$$ It is 
easy to check the conditions of propositions 1.4 -1.5. Thus, we have the 
extension 
$$\CD 0 @>>> J_{4,2,2} @>>> B_{4,2,1} @>>> 
B_{4,2,2} @>>> 0,\endCD\leqno{(\gamma_4)}$$ where
$$\begin{array}{rl}
J_{4,2,2} &= C^*(W_1,{\cal F}_{4,2}) \cong C_0(W_1) \rtimes_{\rho_{4,2}} 
{\RR}^2 \cong C_0({\RR}^*)^2 \times {\RR}) \rtimes_{\rho{4,2}} {\RR}^2 \\
	&\cong C_0({\RR}^* \cup {\RR}^*) \otimes {\cal K},\\
B_{4,2,1} &= C^*(W,{\cal F}_{4,2}) \cong C_0(W) \rtimes_{\rho_{4,2}} 
{\RR}^2 \cong C_0(({\RR}^2)^* \times {\RR}) \rtimes_{\rho_{4,2}} {\RR}^2,\\
B_{4,2,2} &= C^*(W_2,{\cal F}_{4,2}) \cong C_0(W_2) \rtimes_{\rho_{4,2}} 
{\RR}^2 \cong ({\RR}^* \cup {\RR}^*) \times {\RR}) \rtimes_{\rho_{4,2}} 
{\RR}^2\\   &\cong C_0((-1,0),(1,0),(0,-1),(0,1)\} \otimes 
{\cal K} \cong {\CC}^4 \otimes {\cal K}.
\end{array}$$
\end{pf}

Use the Connes-Thom isomorphisms and Bott periodic property of K-groups 
and the fact that our extensions are all appeared from some action of 
${\RR}^2$, it is easy to see that:

\begin{rem}
1.

The six-term exact sequences, associated with the extensions $\gamma_1, 
\gamma_2,\gamma_3$ are term-wise equivalent to the following six-term exact 
sequence
$$\CD K_1(I) @>>> K_1(C_0(({\RR}^3)^*)) @>>> K_1(A)\\
@A\delta_0 AA @.     @VV\delta_1 V\\
K_0(A) @<<< K_0(C_0(({\RR}^3)^*)) @<<< K_0(I) \endCD$$ where $I := 
C_0({\RR}^2 \times {\RR}^*), A := C_0(({\RR}^2)^*).$

 2. The six-term exact sequence, associated with the extension $\gamma_4$ 
is term-wise equivalent to the following six-term exact sequence
$$\CD K_0(I_1) @>>> K_0(C({\SS}^1)) @>>> K_0({\CC}^4)\\
@A \delta_1 AA @.  @V\delta_0VV\\
K_1({\CC}^4) @<<< K_1(C({\SS}^1) @<<< K_1(I_1) \endCD$$ where $I_1 := 
C_0((0,\pi/2) \cup (\pi/2,\pi) \cup (\pi,3\pi/2) \cup (3\pi/2, 2\pi))$, 
${\SS}^1 \approx [0,2\pi)$ and ${\CC}^4 \cong C_0(\{0,\pi/2,\pi, 3\pi/2, 
2\pi\}).$ 
\end{rem}

\begin{rem}
$$\begin{array}{rlrl}
K_0({\CC}) &\cong {\ZZ}, & K_1({\CC}) & = 0,\\
K_0(C_0({\RR})) & = 0, & K_1(C_0({\RR})) &\cong {\ZZ},\\
K_0(C_0({\RR}^2)) &\cong {\ZZ}, & K_1(C_0({\RR}^2)) & = 0,\\
K_0(C_0({\RR}^3)) &= 0, & K_1(C_0({\RR}^3)) &\cong {\ZZ},\\
K_0(C({\SS}^1)) &\cong {\ZZ}, & K_1(C({\SS}^1)) &\cong {\ZZ}.
\end{array}$$
\end{rem}

\begin{rem}1.
Consider the function $u : {\RR} \to {\SS}^1$, defined by the formula
$$u(t) = \exp(2\pi i{t\over \sqrt{1+t^2}}), \forall t\in {\RR}$$ and 
$u_\pm := u|_{{\RR}_\pm}$. Then the homotopy class $[u_\pm]$ are just the 
generators of the K-groups
$$\begin{array}{rl}
K^{-1}({\RR}_+) &\cong K_1(C_0({\RR}_+)) \cong {\ZZ},\\
K^{-1}({\RR}_-) &\cong K_1(C_0({\RR}_-)) \cong {\ZZ}.
\end{array}$$

2. The homotopy class of the constant function $1: {\SS}^1 \to {\SS}^1$ 
having only the value 1 at every point, is the generator of $K^0({\SS}^1) 
\cong K_0(C({\SS}^1)) \cong {\ZZ}$ and the homotopy class of the identity 
function $\Id : {\SS}^1 \to {\SS}^1$ is the canonical generator of the 
group $K^{-1}({\SS}^1) \cong K_1(C({\SS}^1)) \cong {\ZZ}$.
\end{rem}

\begin{lem}
\begin{enumerate}
\item $K_0(I) = 0$, $K_1(I) = {\ZZ}^2$ and is generate by two elements 
$[b] \otimes [u_+]$ and $[b] \otimes [u_-]$.
\item $K_0(C_0(({\RR}^3)^*)) = 0$, $K_1(C_0(({\RR}^3)^*)) \cong {\ZZ}^2.$
\item $K_0(A) \cong {\ZZ}$ and is generated by the element $[\Id] \otimes 
[u_+]$. $K_1(A) \cong {\ZZ}$ and is generated by the element $[1] \otimes 
[u_+]$.
\end{enumerate}
\end{lem}
\begin{pf}
1.

Because ${\RR}^2 \times {\RR}^* = ({\RR}^2 \times {\RR}_+) \cup ({\RR}^2 
\times {\RR}_-),$ we have 
$$I = C_0({\RR}^2 \times {\RR}^*) \cong C_0(({\RR}^2) \otimes 
C_0({\RR}_+) \oplus C_0({\RR}^2) \otimes C_0({\RR}_-).$$ Thus, we have 
$$K_j(I) \cong K_j(C_0({\RR}^2) \otimes C_0({\RR}_+)) \oplus 
K_j(C_0({\RR}^2) \otimes {\RR}_-); j = 0,1$$  
and hence,
$K_0(I) = 0$. From the Bott periodicity and [\cite{connes2}, Corollary 
VI.3], it is easy to see that $[b] \otimes [u_+]$ (resp., $[b] \otimes 
[u_-]$) is just the generator of $K_1(C_0({\RR}^2)\otimes C_0({\RR}_+))$ 
(resp. $K_1(C_0({\RR}^2) \otimes C_0({\RR}_-))$). Hence, $K_1(I) \cong 
{\ZZ}^2$ is generated by two elements $[b]\otimes [u_+]$ and $[b] \otimes 
[u_-]$. 

2.

Consider the natural exact sequence 
$$\CD 0 @>>> C_0(({\RR}^3)^*) @>\sigma>> C_0({\RR}^3) @>\kappa>> {\CC} 
@>>> 0,\endCD$$ where $\kappa$ is the restriction of function to the 
point $\{0\}$. From the associate six-term exact sequence and Remark 3.2, 
we have $$K_0(C_0({\RR}^3)^*) = 0, K_1(C_0({\RR}^3)^*) \cong {\ZZ}^2.$$  

3.

Recall that $A=C_0(({\RR}^2)^*),$ where $({\RR}^2)^* := {\RR}^2 \{0\} 
\approx {\SS}^1 \times {\RR}_+$. Thus we have $A \cong C_0({\SS}^1 \times 
{\RR}_+) \cong C_0({\SS}^1) \otimes C_0({\RR}_+)$ and $$K_0(A) \cong 
K_0(C({\RR}^1) \otimes C_0({\RR}_+)) \cong K_1(C_0({\SS}^1)) \cong {\ZZ}.$$
Also in virtue of the Bott periodicity and [\cite{connes3}, Cor. VI.3], 
this group $K_0(A)$ is generated by $[\Id]\otimes [u_+]$. By analogy, 
$$K_1(A) \cong K_1(C_0({\SS}^1) \otimes C_0({\RR}_+)) \cong 
K_0(C_({\SS}^1)) \cong {\ZZ}$$ is generated by element $[1] \otimes 
[u_+]$. 
\end{pf}

\begin{cor}
The six-term exact sequences, associated with the extensions 
$\gamma_1,\gamma_2,\gamma_3$ can be identified with the following exact 
sequence
$$\CD {\ZZ}^2 @>>> {\ZZ}^2 @>>> {\ZZ}\\
@A\delta_0 AA @. @VV\delta_1=0V\\
{\ZZ} @<<< 0 @<<< 0 \endCD$$
\end{cor}

Let us consider two matrix functions: The constant $(2\times 
2)$-matrix-valued function
$$\pmatrix 1 & 0\\ 0 & 0 \endpmatrix : ({\RR}^2)^* := {\RR}^2 \backslash 
\{(0,0)\} \to \Mat_2({\CC})$$ and the function
$$p : ({\RR}^2)^* \approx {\SS}^1 \times {\RR}_+ \to \Mat_2({\CC}),$$ 
given by the formula
$$p(e^{i\varphi},r) := {1\over 2} \pmatrix 1-\cos(r\pi)
&e^{1\varphi}\sin(r\pi)\\ e^{-i\varphi}\sin(r\pi) &1+\cos(r\pi)
\endpmatrix, \forall (e^{i\varphi},r) \in {\SS}^1 \times {\RR}_+ \approx
({\RR}^2)^*.$$

\begin{lem}
1. For each $(e^{i\varphi},r) \in {\SS}^1 \times {\RR}_+ \cong
({\RR}^2)^*$, the matrix $p(e^{1\varphi},r)$ is an idempotent.

2. $[p]-\left[ \pmatrix 1 & 0\\ 0 & 0 \endpmatrix \right]$ is equal to
the generator $[\Id] \otimes [u_+]$ of $K_0(A) = K_0(C_0(({\RR}^2)^*))
\cong {\ZZ}.$
\end{lem}
\begin{pf}
The first assertion is proved by a direct computation. The second
assertion is deduced also from a direct computation in using
[\cite{connes3}, Lemma .2].
\end{pf}

\begin{rem}
It is easy to see that
$$\begin{array}{rl}
K_j(I_1) &= K_j(C(0,\pi/2) \cup (\pi/2, \pi) \cup (\pi, 3\pi/2) \cup
(3\pi/2, 2\pi))\\     &K_j(C_0(0,\pi/2)) \oplus K_j(C_0(\pi/2, \pi)) \oplus
K_j(C_0(\pi, 3\pi/2)) \oplus \\ &\oplus K_j(C_0(\pi/2, 2\pi)), \forall
j=0,1.\end{array}$$ Thus $$K_0(I_1) = 0 \mbox{ and } K_1(I_1) \cong
{\ZZ}^4.$$

Certainly, $$K_j(\{ 0,\pi/2,\pi,3\pi/2,2\pi\}) \cong K_1({\CC}^4) \cong
\cases {\ZZ}^4 &\mbox{ if } j = 0\\ 0 &\mbox{ if } j =1. \endcases$$
\end{rem}

\begin{rem}
Consider the function $f: [0,2\pi] \to {\SS}^1$, given by
$$f(\varphi) := e^{4i\varphi}, \forall \varphi\in [0,2\pi].$$ Denote by
$u_1$, $u_2$, $u_3$, $u_4$ the restrictions of $f$ to $(0,\pi/2),$
$(\pi/2,\pi)$, $(\pi,3\pi/2)$, $(3\pi/2,\pi)$, resp., we see that all of
them have rotation number 1 and hence they are just the generators of
$K_1(C_0(0,\pi))$, $K_1(C_0(\pi/2, \pi))$, $K_1(C_0(\pi, 3\pi/2))$,
$K_1(C_0(3\pi/2,2\pi))$, resp. .
\end{rem}

\begin{rem}
Define
$$\left. \begin{array}{rl} v_1 &= (1,0,0,0),\\ v_2 &= (0,1,0,00), \\ v_3 &=
(0,0,1,0),\\ v_4 &= (0,0,0,1)\end{array} \right\} \in
C(\{0,\pi/2,\pi,3\pi/2,2\pi\}) \cong {\CC}^4.$$ Then it is easy to see
that their homotopy classes are just the generators of $$K_0(C(\{
0,\pi/2,\pi,3\pi/2,2\pi\})) \cong K_0({\CC}^4) \cong {\ZZ}^4.$$
\end{rem}

\begin{cor}
The six-term exact sequence, associated with the extension $\gamma_4$
can be identified with the following six-term exact sequence
$$\CD 0 @>>> {\ZZ} @>>> {\ZZ}^4\\
@A\delta_1=0 AA @. @VV\delta_0V\\
0 @<<< {\ZZ} @<<< {\ZZ}^4 \endCD$$
\end{cor}

\begin{rem}[\cite{taylor},\S5]
Let us denote $\widetilde{C_0({\RR}^2)} = C({\SS}^2)$ the algebra of
compact support continuous functions on ${\RR}^2$ with the formal
adjoint unity element $1$,
$$Q(\widetilde{C_0({\RR}^2)}) := \{a\in C_0({\RR}^2); \exp(2\pi i a) = 1  \},$$
and
$$Q_n(\widetilde{C_0({\RR}^2)}) := \{a\in \Mat_n(C_0({\RR}^2));
\exp(2\pi i a) = 1_n  \}, n\in {\NN}$$. We have a trace map
$$\Tr : \bigcup_n Q_n(\widetilde{C_0({\RR}^2)}) \to
Q(\widetilde{C_0({\RR}^2)}).$$ Each element $[f] \in K_1(C_0({\RR}^2
\times {\RR}_\pm))$ has a representative, as a function
$$ f : {\RR}_\pm \to \bigcup_n Q_n(\widetilde{C_0({\RR}^2)}),$$ such
that $$\lim_{t\to 0} f(t) = \lim_{t\to\pm\infty} f(t) $$ and
$$[f] = w_f.([b] \otimes [u_\pm],$$ where $w_f$ is the winding number of
$f$, defined by the formula
$$w_f := {1\over 2\pi i} \int_{{\RR}_\pm} \Tr(f'(t).f(t)^{-1}) dt.$$
\end{rem}

\begin{thm}
The isomorphic class of C*-algebras of $\MD_4$-foliations of type ${\cal
F}_7$, ${\cal F}_8$, ${cal F}_9$ are defined exactly by the following
KK-theory invariants
\begin{enumerate}
\item
$index C^*(V_{G_{4,3,4(1,\pi/)}},{\cal F}_{3,4(1,\pi/2)} = [\gamma_1] =
(1,1)$ in the KK-group $$\Ext(C_0({\RR}_+) \otimes {\cal K},
C_0({\RR}^2\cup {\RR}^2) \otimes {\cal K}) \cong
\Hom_{\ZZ}({\ZZ},{\ZZ}^2) \cong {\ZZ}^2.$$
\item
$Index C^*(V_{G_{4,4,1}},{\cal F})_{4,1}) = [\gamma_2] = (1,1)$ in the
KK-group $$\Ext(C_0({\RR}_+)\otimes {\cal K}, C_0({\RR}^* \times {\RR})
\otimes {\cal K} \cong \Hom_{\ZZ}({\ZZ}, {\ZZ}^2) \cong {\ZZ}^2.$$
\item
$Index C^*(V_{G_{4,4,2}},{\cal F}_{4,2}) = ([\gamma_3],[\gamma_4]) ,$
where $[\gamma_3] = (1,1)$ in the KK-group
$$\Ext(C_0(({\RR}^2)^* \times {\RR}) \rtimes_{\rho_{4,2}} {\RR}^2,
C_0({\RR}^* \times {\RR}) \otimes {\cal K} \cong
\Hom_{\ZZ}({\ZZ},{\ZZ}^2) \cong {\ZZ}^2,$$
$$[\gamma_4] = \pmatrix -1 & 1 & 0 & 0\\ 0 & -1 & 1 & 0\\ 0 & 0 & -1 &
1\\ 1 & 0 & 0 & -1 \endpmatrix$$ in the KK-group
$$\Ext({\CC}^4 \otimes {\cal K}, C_0({\RR}^* \cup {\RR}^*) \otimes {\cal
K}) \cong \Hom_{\ZZ}({\ZZ}^4, {\ZZ}^4) \cong \Mat_4({\ZZ}).$$
\end{enumerate}
\end{thm}
\begin{pf}
Following the general conception of index, it is easy to see that that
isomorphic class of $C^*(G)$ is defined by the KK-invariant $Index
C^*(G)= [\gamma_1], [\gamma_2]$, or $([\gamma_3],[\gamma_4])$. The
classes $[\gamma_1]$, $[\gamma_2]$, $[\gamma_3]$ in virtue of Remark 2.1
can be identified with the connecting homomorphisms
$$\delta_0\in \Hom_{\ZZ}(K_0(A),K_1(I)) =
\Hom_{\ZZ}(K_0(C_0(({\RR}^2)*),$$
$$ K_1(C_0({\RR}^2 \times {\RR}^*) \cong
\Hom_{\ZZ}({\ZZ},{\ZZ}^2),$$
and $$\delta_1 = 0\in \Hom_{\ZZ}(K_1(A),K_0(I)) \cong
\Hom_{\ZZ}({\ZZ},0) = 0.$$
By the same reason, in virtue of Remark ....., the extension $\gamma_4$ is
characterized by the connecting homomorphisms
$$\delta_0 \in \Hom(K_0({\CC}^4),K_1(I_1)) \cong \Hom_{\ZZ}({ZZ}^4,
{\ZZ}^4) \cong \Mat_4({\ZZ}).$$
It rests therefore to prove the connecting homomorphisms for two cases
$$\delta_0\in \Hom_{\ZZ}(K_0(A),K_1(I))$$
and
$$\delta_0 \in \Hom(K_0({\CC}^4),K_1(I_1)).$$

1. Computation of $\delta_0\in \Hom_{\ZZ}(K_0(A),K_1(I))$

As it was said, $K_0(A) = K_0(C_0(({\RR}^2)^*)) \cong {\ZZ}$ is
generated by the class $[p] - \left[\pmatrix 1 & 0\\ 0 & 0 \endpmatrix
\right]$, we need only to compute $$\delta_0([p] - \left[\pmatrix 1 & 0\\ 0 & 0 \endpmatrix
\right]) = \delta_0([p]) - \delta_0(\left[\pmatrix 1 & 0\\ 0 & 0 \endpmatrix
\right]).$$  Recall that following J. Taylor (\cite{taylor}, p. 170)
for each idempotent $f \in \Mat_n(C_0(\widetilde{({\RR}^2)^*)}) =
\Mat_n(\widetilde{A})$, i.e. $[f] \in K_0(\widetilde{A})$, the value
$\delta_0([f])$ is given by
$$\delta_0([f]) = [\exp(2\pi i \tilde{f})],$$ where $\tilde{f}\in
\Mat_n(\widetilde{C_0(({\RR}^3)^*)})$ such that  $\nu{\tilde{f}} = f$,
i.e. the restriction of $\widetilde{f}$ on
$\Mat_n(\widetilde{C_0(({\RR}^2)^*)})$ is $f$. For $f = \pmatrix 1 & 0
\\ 0 & 0 \endpmatrix,$ we can choose $$\tilde{f}(x,y,z) = f.1_{\RR} =
\pmatrix 1 & 0\\ 0 & 0 \endpmatrix,$$ where
$1_{\RR}$ is the function with a single value 1 on ${\RR}$.
$$\delta_0\left(\left[\pmatrix 1 & 0 \\ 0 & 0 \endpmatrix
\right]\right) = [\exp(2\pi i\tilde{f}] \in K_1(I) = K_1(C_0({\RR}^2
\times {\RR}_+)) \oplus K_1(C_0({\RR}^2 \times {\RR}_-)) \cong
{\ZZ}^2.$$ It is easy to see that
$$[\exp(2\pi i\tilde{f})] = ([\exp(2\pi i\tilde{f}_+)],[\exp(2\pi
i\tilde{f}_-)]),$$ where $\tilde{f}_\pm := f.1_{{\RR}_\pm}$.  It is easy
to see, by using the winding number formula, that $$\delta_0([\pmatrix 1
& 0 \\ 0 & 0 \endpmatrix ]) = (0,0).$$	To compute $\delta([p])$, recall
that $p : ({\RR}^2)^* \to \Mat_2({\CC})$. Choose $$\tilde{p} :
({\RR}^3)^* \to \Mat_2({\CC}),$$ as follows
$$\tilde{p}(x,y,z) := {z\over\sqrt{x^2+y^2}} p(x,y), \forall (x,y,z) \in
({\RR}^3)^*.$$ Denote the restrictions of $\tilde{p}$ on ${\RR}^2 \times
{\RR}_\pm$ by $\tilde{p}_\pm$. We have
$$\delta_0([p]) = [\exp(2\pi i \tilde{p})] = ([\exp(2\pi i
\tilde{p}_+)],$$ $$[\exp(2\pi i\tilde{p}_+)] ) \in K_1(C_0({\RR}^2 \times
{\RR}_+)) \oplus K_1(C_0({\RR}^2 \times {\RR}_-)).$$ Following the
winding formula, we have
$$[\exp(2\pi i\tilde{p}_\pm)] = w_\pm.([b] \otimes [u_\pm]),$$ where
$$\begin{array}{rl} w_\pm &:= {1\over 2\pi i} \int_{{\RR}_\pm} \Tr({d\over
dz}(exp(2\pi i\tilde{p}_\pm)\exp(-2\pi i\tilde{p}_\pm)dz \\
&= \int_{{\RR}_\pm} {z \over \sqrt{x^2+y^2}} dz = 1.\end{array}$$ Thus we 
have $$\delta_([p]) = ([b] \otimes [u_+],[b] \otimes [u_-]),$$
$$\delta_([p]-\left[\pmatrix 1 & 0 \\ 0 & 0 \endpmatrix\right]) = ([b]
\otimes [u_+],[b] \otimes [u_-]),$$  This means that $$\delta = (1,1) \in
\Hom_{\ZZ}({\ZZ},{\ZZ}).$$

2. Computation of $\delta_0 \in \Hom(K_0({\CC}^4),K_1(I_1)).$

We constructed the generators $[v_1], [v_2], [v_3], [v_4]$ of the K-groups 
$$K_0(\{0,\pi/2,\pi,3\pi/2,2\pi\}) \cong K_0({\CC}^4) \cong {\ZZ}^4 $$ 
and the generators $[u_1], [u_2],[u_3],[u_4]$ of the K-groups
$$K_1(\{0,\pi/2,\pi,3\pi/2,2\pi\}) \cong K_1({\CC}^4) \cong {\ZZ}^4.$$
Choose a continuous function $\ell_1 \in C([0,2\pi])$, such that $\ell_1(0) 
= \ell_1(2\pi) = 1$ ,$\ell_1|_{[\pi/2,3\pi/2]} \equiv 0$ and linear 
outside this interval, and define 
$\tilde{v}_1(e^{i\varphi}) := \ell_1(\varphi),$ then the restriction of 
$\tilde{v}_1$ on $\{0,\pi/2,\pi,3\pi/2,2\pi \}$ is just $v_1$. We have, 
following the definition of $\delta_0$,
$\delta_0([v_1]) = [e^{2\pi i\tilde{v_1}}] \in K_1(I_1).$
Recall that
$$\begin{array}{rl}
K_1(I_1) &= K_1(C_0((0,\pi/2) \cup (\pi/2,\pi)\cup (\pi,3\pi/2) \cup 
(3\pi/2,2\pi)))\\   &\cong K_1(C_0(0,\pi/2)) \oplus K_1(C_0(\pi/2,\pi)) 
\oplus K_1(C_0(\pi, 3\pi/2))\oplus  \\  &\oplus K_1(C_0(3\pi/2, 2\pi))\\
  &\cong{\ZZ} \oplus {\ZZ} \oplus {\ZZ} \oplus {\ZZ} = {\ZZ}^4.\end{array}$$
Denote $\tilde{v}_{11}, \tilde{v}_{12}, \tilde{v}_{13}, \tilde{v}_{14}$ 
the restrictions of $\tilde{v}_1$ on 
$(0,\pi/2)$,$(\pi/2,\pi)$,$(\pi,3\pi/2)$,$(3\pi/2,2\pi)$, respectively. Then,
$$[e^{2\pi i\tilde{v}_1}] = ([e^{2\pi i\tilde{v}_{11}}], [e^{2\pi 
i\tilde{v}_{12}}], [e{2\pi i\tilde{v}_{13}}], [e^{2\pi 
i\tilde{v}_{14}}]).$$ From the definition, it is easy to see that 
$\tilde{v}_{12} = \tilde{v}_{13} = 0$ and hence $$[e^{2\pi 
i\tilde{v}_{12}}] = [e^{2\pi i\tilde{v}_{13}}] = 0.$$ 
Following the winding formula we have  \cite{taylor}
$$[e^{2\pi i\tilde{v}_{11}}] = w_1[u_1], [e^{2\pi i\tilde{v}_{14}}] =
w_4[u_4],$$ where $w_1$ and $w_4$ are the winding numbers of
$\tilde{v}_{1}$ and $\tilde{v}_{14}$, respectively,
$$w_1 = {i\over 2\pi i}\int_0^{\pi/2} {d\over d\varphi}(e^{2\pi
i\tilde{v}_{11}}) e^{-2\pi i\tilde{v}_{11}} d\varphi = -1,$$
$$w_4 = {i\over 2\pi i}\int_0^{\pi/2} {d\over d\varphi}(e^{2\pi
i\tilde{v}_{14}}) e^{-2\pi i\tilde{v}_{14}} d\varphi = +1.$$  Thus,
$$\delta_0([v_1]) = (-[u_1],0,0,[u_1]).$$
By analogy, we have also
$$\begin{array}{rl}
\delta_0([v_2]) &= ([u_1], -[u_2],0,0),\\
\delta_0([v_3]) &= (0, [u_2], -[u_3],0),\\
\delta_0([v_4]) &= (0,0,[u_3], -[u_4]). \end{array}$$ This means that
$$\delta_0 = \pmatrix -1 & 1 & 0 & 0\\ 0 & -1 & 1 & 0\\ 0 & 0 & -1 & 1
\\ 1 & 0 & 0 -1 \endpmatrix  \in \Hom_{\ZZ}({\ZZ}^4, {\ZZ}^4) \cong
\Mat_4({\ZZ}).$$
\end{pf}

\section{Bibliographical Remarks}
The main idea of studying the structure of C*-algebras of $\MD_4$-group
was proposed by the author. The new idea of this chapter is to study the
C*-algebras of measurable foliations, consisting of generic co-adjoint
orbits. It is just the intersection of 1-dimensional irreducible
representations. This idea was then developed for the general situation
by th author in \cite{diep11}.
The results of this chapter was done in Ph. D. Dissertation (1990) of Le Anh Vu
under supervision of the author.


\part {Advanced Theory: Multidimensional Quantization and Index of 
Group C*-Algebras} 
\chapter {Multidimensional Quantization}
\section{Induced Representations. Mackey Method of small subgroups}

We begin our exposition with a survey of key concepts and results in the
so called Mackey method of small subgroups. This method will play an
essential role  in the development of the orbit method.

Throughout this section $G$ denotes a locally compact group and all actions
of $G$ are continuous. 

\subsection{Criterion of inductibility}

\begin{defn} Let $M$ be a right $G$-space and $V$ a left $G$-space. Let $M \times
V$ denote the Cartesian product of $M$ and $V$ with the product topology. Let $\sim$
denote the equivalence relation on $M \times V$ given by $$(m,v) \sim (m',v') \iff
\exists g \in G$$ such that $m' = mg^{-1}$, $v' = gv$. Define $$M \times_G V:= ( M
\times V )/\sim $$ with the quotient topology. $M \times_G V$ is called the {\it
fibered product of M and V over G. } It will be usually denoted by ${\cal E}_V$.
\end{defn}

\begin{rem}

(1) The right $G$-space $M$ can be considered as a left $G$-space in the
obvious way via the action $m \mapsto g.m:= mg^{-1}$. Giving the product
$M \times V$ the diagonal action of $G$, we can consider the left
$G$-space $M \times V$, and the orbits of $G$ therein. Two points
$(m,v)$ and $(m',v')$ are on the same $G$-orbit, if and only if there
exists an element $g \in G$, such that $m' = mg^{-1},v' = gv $. The
space of all orbits of $G$ in $M \times V$ is the space $M \times_G V$
defined above. 

(2) If the quotient map $M \longrightarrow M/G$ is a fibration ( moreover
a {\it principal bundle } ), \index{bundle!principal -} then $$ V 
\rightarrowtail {\cal E} _V
\twoheadrightarrow M/G$$ is the so called {\it vector bundle, associated
with the $G$-action on V }. \index{bundle!vector -,associated} It is 
easy to see that there exists a
bijection between the sections of this bundle and the $V$-valued functions
on $M$ satisfying the following equations $$f(mg) = g^{-1}f(m),$$ for all
$m \in M$ and $g \in G$.

(3) We shall often apply the construction of {\it product over G } 
\index{product!over G} to the
left spaces. For this, it is enough to remark that if $M$ is a left
$G$-space, then it is also a right $G$-space by the action $g^{-1}M, g
\in G$. In this case, ${\cal E} _V = ( M \times V)/\sim$ and, by
definition, $(m,v) \sim (m',v')$ iff $m' = gm, v' = gv$. The sections
$s \in \Gamma({\cal E} _V)$ of the associate bundle $V \rightarrowtail
{\cal E} _V \twoheadrightarrow G/M$ can be identified with the $V$-valued
functions on $M$, satisfying the equations $$f_s(gm) = gf_s(m), \forall
g \in G, \forall m \in M \quad.$$

(4) We shall use the construction to the case of induced representations.
\end{rem}

\begin{defn} Let $G$ be a locally compact group,
$H$ a closed subgroup and $(\sigma,V)$ a unitary representation of $H$.
The (left) $G$-action on the space $\Gamma({\cal E} _V)$ of sections of
the associate $G$-bundle $$V \rightarrowtail {\cal E} _V
\twoheadrightarrow H \setminus G$$ will be called the {\it representation
of G induced } \index{representation!induced} from the unitary 
representation 
$(\sigma,V)$ of $H$. It is denoted by $\Ind^G_H(\sigma,V)$.  \end{defn}

\begin{rem}

(1) It is easy to see that $G$ is a left $H$-space, and that $H
\rightarrowtail G \twoheadrightarrow H \setminus G$ is a fibration. We
shall always suppose that there is a Borel section $s: H \setminus G
\rightarrow G $, $ x \mapsto s_x$ such that the decomposition $g = h.s_x$
is unique. It is easy to see that for every $x \in X = H \setminus G $
and every $g \in G$, there exists a unique $h = h(g,x)$ such that $$s_x.g
= h(g,x).s_{xg}$$ and $$\sigma(h(g_1,x)) \sigma(h(g_2,xg_1)) =
\sigma(h(g_1g_2,x)),$$ for all $x \in X$ and $ g_1, g_2 \in G$. We have
therefore a $\sigma(H)$-valued 1-cocycle $\sigma(h(.,.))$.

(2) The 1-cocycle $\sigma(h(.,.))$ acts on the fiber $V$ and $G$ acts on
$H \setminus G$ by the right translations. Together, we have therefore a
left $G$-action on the sections of the associate bundle ${\cal E} _V$,$$
(g.s)(x) = \sigma(h(g,x))s(xg), \forall g \in G, \forall x \in X \quad
.$$

(3) We always suppose that $H \rightarrowtail G \twoheadrightarrow H
\setminus G$ is a locally trivial principal $H$-bundle and fix a {\it
trivialization } $\Gamma$ on it. It is easy to see that on the associated
bundle ${\cal E} _V$, there is a natural associated ( affine ) connection
$\nabla$.

(4) It is also easy to see that on the right $G$-space $X:= H \setminus
G$, there exists a unique, up to scalar factor, quasi - invariant
measure $ d\mu = d\mu_s$, depending on the Borel section $s$. Let
$B(.,.)$ be the 1-cocycle defined by $$B(g,x):= \left[{d\mu(xg) \over
d\mu(x)}\right]^{-1/2}\sigma(h(g,x)) \quad.$$ Then the corresponding induced
representation $\Ind^G_HB(.,.)$ is unitary iff $B(.,.)$ is a unitary
operator - valued 1-cocycle.

(5) Let $T_i, i = 1,2$ be the induced representations corresponding to
the 1-cocycles $B_i, i = 1,2$, respectively. Then $T_1$ and $T_2$ are
unitary equivalent if and only if $B_1$ and $B_2$ are cohomologuous, i.e.
$$B_2(g,x) = C(x)^{-1}B_1(g,x)C(xg) \quad, $$ for some operator - valued
function $C(.)$ on $X$.
\end{rem}

\begin{defn} Let $T$ be a unitary representation of
$G$ on a Hilbert space ${\cal H }$, $X = H \setminus G$ and $C_0(X)$ the
*-algebra of continuous functions on $X$, vanishing at infinity, and $P$
a *-representation of $C_0(X)$ on the same Hilbert space ${\cal H}$. The
pair $(T,P)$ will be called a {\it unitary representation of the right
$G$-space $X$,} \index{representation!unitary - of G-space} iff $T$ and 
$P$ satisfy the so called {\it
quasi-invariance condition } $$T(g)P(f)T(g^{-1}) = P(R(g)f), \forall g
\in G, \forall f \in C_0(X) \quad,$$ where $(R(g)f)(x):= f(xg)$ is the
right {\it regular representation} \index{representation!regular -} of 
$G$. In this case we say that $T$
{\it can be extended } \index{representations!extended} to a unitary 
representation $(T,P)$ of the right $G$-space $X$.  \end{defn}

\begin{thm}[\bf Criterion of Inductibility]  Let T be a unitary
representation of $G$. Then there exist a closed subgroup $H$ and a
unitary representation $(\sigma,V)$ of $H$ such that $T = \Ind^G_H(\sigma
, V)$ if and only if $T$ can be extended to a unitary representation
$(T,P)$ of the right $G$-space $X = H \setminus G$.
\end{thm}

\begin{rem}

(1) The criterion of inductibility is equivalent to the following condition
on the {\it existence of systems of imprimitivity } \index{systems! of 
imprimitivity} on $X$: {\it There is
a one - to - one correspondence between the unitary representations
$(T,P)$ of the right $G$-space $X$ and the projection measures $\Delta$ on
$X$. }

Recall that a {\it projection measure } \index{measure!projection 
-} $\Delta$ on $X$ is by definition a
map from the $\sigma$-algebra ${\cal B} (X)$ of the Borel sets $E$ in $X$
to the projections in a Hilbert space, satisfying the following
conditions
\begin{enumerate}
\item[(a)]$\Delta(E_1 \cap E_2) = \Delta(E_1).\Delta(E_2)$,
for all Borel sets $E_1$ and $E_2$,

\item[(b)]$\Delta(\cup_{i=1}^{\infty} E_k) = \sum_{i=1}^{\infty}
\Delta(E_k)$, for all Borel sets $E_1,E_2,\dots $, such that $E_i \cap
E_j = \emptyset, i \ne j$,
v\item[(c)]$P(f) = \int_X f(x)\Delta(dx), \forall f \in C_0(X) $, and
\item[(d)]$T(g)\Delta(E)T(g^{-1}) = \Delta(Eg) \forall g \in G, \forall E
\in {\cal B} (X)$.
\end{enumerate}

The assertion above is clear in view of the following facts from
functional analysis: The representation $P$ of the *-algebra $C_0(X)$ is
just the quasi - invariant ( condition in Def. 1.5. ) integral,
corresponding to the quasi - invariant (condition (d)) projection ( condition (a) 
\& (b) ) measure $\Delta$, following the condition (c).

(2) Consider the convolution algebra ${\cal A} = C_0(G \times X) \cong
C_0(X \times G)$, consisting of all continuous functions with compact
support of type $\alpha = \alpha(x,g)$, and with the following
convolution product and involution: $$(\alpha_1 * \alpha_2)(x,g):=
\int_G \alpha_1(x,g_1^{-1})\alpha_2(xg_1^{-1},g_1g)d\mu_r(g_1), \forall x
\in X, g \in G \quad,$$ $$\alpha^*(x,g):= \overline{\alpha(xg,g^{-1})}
\Delta_G(g_1)^{-1} \quad,$$ where $\Delta_G$ is the {\it modular function
} \index{function!modular -} of the right invariant measure $d\mu_r$ on $G$.
\end{rem}

A *-representation $\phi$ of the convolution *-algebra ${\cal A }$ will be
called {\it consistent of type $T$ }, \index{representation!consistent - of 
type $T$} iff $$T(g_1)\phi(f)T(g_2) =
\phi(\tilde{f}) \quad, $$ where $$\tilde{f}(x,g):= f(xg_1,
g_1^{-1}gg_2^{-1}) \Delta_G (g_1)^{-1} \quad.$$ {\it There is a one - to
- one correspondence between the unitary representations of type $(T,P)$
of $C_0(X)$ and the consistent *-representation $\phi$ of the algebra
${\cal A} = C_0(X \times G)$ of type $T$.}

This remark is also clear from the point of view of the representation
theory: The algebra $C_0(X \times G)$ seems to be the group *-algebra for
the homogeneous space $X$. Thus the one - to - one correspondence between
the quasi - invariant projection measures on $X$ and the consistent
*-representations of $C_0(X \times G)$ is well-known.

{\sc Proof of the criterion of inductibility.

Necessity. }

Suppose that $T = \Ind^G_H(\sigma,V)$, realizing on the space
$L^2_{\mu_s}(X,V)$ of square - integrable sections of the induced bundle
${\cal E} _V$. For each element $f \in C_0(X)$, define $P(f)$ as the
operator of multiplication by $f$ in the space of sections
$L^2_{\mu_s}(X,V)$ of the representation $T$. It is easy to see that
$P(f)$ is a quasi-invariant in the sense of (1.5).

{\sc Sufficiency. } Suppose that our representation $T$ extends to a
unitary representation $(T,P)$ of the homogeneous space $X = H \setminus
G$. Following our remarks 1.7, there exist a consistent *-representation
$\phi$ of *-algebra ${\cal A } = C_0(X \times G)$ such that $$\phi(\alpha)
= \int_G P(\alpha(.,g))T(g)d\mu_r(g) \quad. $$ By restriction to
components, we can suppose that $\phi$ is a cyclic representation,
acting on $L_{\mu_s}^2(X,V)$ with {\it source, ( or cyclic vector ) } 
\index{vector!source -} \index{vector!cyclic -} $\xi$.

For every element $\beta \in C_0(G \times G)$ we can produce an element
$\alpha \in C_0(X \times G)$, such that $$\alpha(x,g_2) = \alpha(Hg_1,
g_2):= \int_H \beta(hg_1,g_2)d\mu_r(h) \quad.$$ Choose a positive
function $\alpha \geq 0$ on $G$ such that $$\int_H \rho(hg)d\mu_r(h)
\equiv 1 \quad,$$ we have $$\begin{array}{rl}(\phi(\alpha)\xi,\xi) &= \int_G
((\phi(\alpha)\xi)(g_1),\xi(g_1))_V \rho(g_1) d\mu_r(g_1) \cr
   &= \int_G \int_G \alpha(Hg_1,g_2)(\xi(g_1g_2),\xi(g_1))_V \rho(g_1)
d\mu_r(g_1)d\mu_r(g_2) \cr
   &= \int_H \int_G \int_G
\beta(hg_1,g_2)(\xi(g_1g_2),\xi(g_1))_Vd\mu_rd\mu_r(g_1)d\mu_r(g_2). 
\end{array}$$ Changing the variable $g_1 \rightsquigarrow h^{-1}g_1$, we have
$$(\phi(\alpha)\xi,\xi) = \int_G\int_G
\beta(g_1,g_2)(\xi(g_1g_2),\xi(g_1))_V d\mu_r(g_1)d\mu_r(g_2) \quad.$$
Taking $\beta(g_1,g_2) = \theta_1(g_1g_2)\overline{\theta_2(g_1)}
\Delta_G(g_1)$, with $\theta_i \in C_0(G)$, we have
$$(\phi(\alpha)\xi,\xi) = (\hat{\theta}_1, \hat{\theta}_2)_V \quad,$$
where, by definition, $$\hat{\theta}:= \int_G \theta(g)\xi(g)d\mu_r(g)
$$ and because $\xi$ is the cyclic vector, the set $\{ \hat{\theta};
\theta \in C_0(G) \}$ is everywhere dense in $V$.

In $C_0(G)$ we define $(\theta_1,\theta_2):= (\phi(\alpha)\xi,\xi)$ and
a unitary representation $\sigma$ of $H$ by $$(\sigma(h)\theta)(g):=
\Delta_H(h)^{-1/2}.\Delta(h)^{-1/2} \theta(h^{-1}g) \quad.$$ Consider the
induced representation $\Ind^G_H(\sigma,V)$ on the space ${\cal H}$ \enskip
of functions $F: G \rightarrow C_0(G)$ such that $$F(hg_1, g_2) =
\Delta_G^{-1}(h)F(g_1,h^{-1}g_2) \quad, $$ with the scalar product,
defined by $$(F_1,F_2) = \int (F_1(g,.)F_2(g,.))_{C_0(G)} \rho d\mu_r(g)
.$$ The dense subspace $$L_0 = \{ \phi(\alpha)\xi ; \alpha \in C_0(X
\times G) \}$$ in $L_{\mu_s}^2(X,V)$ can be isometrically mapped into
${\cal H}$ by $$ \tau: \phi(\alpha)\xi \in L_0 \mapsto F \in {\cal H }
\quad,$$ $$F(g_1,g_2) = \alpha(Hg_1,g^{-1}_1g_2)\Delta_G(g_1)^{-1} \quad
.$$ The map $\tau$ commutes with the action of $G$ and really is an
isometric linear operator with unitary closure. Thus we have $T \cong
\Ind^G_H(\sigma, \tilde{V})$, where $\tilde{V} = C_0(G)$ with the scalar
product $(\phi(\alpha)\xi,\xi)$ as above. The theorem is proved.

\subsection{The Mackey method of small subgroups }

The main application of the theory of induced representation is probably
the Mackey theory, which completely describes the set of all irreducible
representations of locally compact groups having proper closed normal
subgroups via the induction process.

\begin{thm}[\bf Mackey Method of Small Subgroups]
Let $G$ be a locally compact group, $N$ a closed normal subgroup,
$\hat{N}$ the dual object of $N$ (i.e. the set of unitary equivalence
classes of irreducible unitary representations), on which $G$ acts in
the natural way, $\forall (\sigma \in) \langle\sigma\rangle  \in \hat{N}, \forall g
\in G $,
$$(\sigma g)(n):= \sigma(gng^{-1}), \forall n \in N \quad.$$
Suppose that $G$-orbit space $\hat{N}/G$ is of class $T_0$. Let
$G_{\sigma}$ be the stabilizer of $\langle \sigma\rangle  \in \hat{N}$,which contains
$N$. Let $\tau$ be a representation such that $\langle \tau\rangle  \in
\hat{G}_{\sigma}$ and the restriction $\tau|_N$ is equivalent to a
multiple of $\sigma$. We can  therefore consider the induced
representation $\Ind_{G_{\sigma}}^G \tau$.

Under these assumptions, there is a one - to - one correspondence
between the dual object $\hat{G}$ of $G$ and the set of all these
classes of induced representations of type $\Ind_{G_{\sigma}}^G \tau$.
\end{thm}

\begin{pf}
It is well-known that $N$ is a type $I$ group. Let $T$ be an arbitrary
(unitary ) irreducible representation of $G$, in some Hilbert space
$L$. The the restriction $T|_N$ can be unitarily and uniquely decomposed
into a direct integral of representations
$$T|_N \simeq \int_{\hat{N}}^{\oplus} W_{\lambda} d\mu (\lambda) \quad,
\quad  L = \int_{\hat{N}}^{\oplus} L_{\lambda} d\mu(\lambda) \quad,$$
where $ W_{\lambda} \simeq U_{\lambda} \otimes S_{\lambda}, L_{\lambda}
\cong  V_{\lambda} \otimes {\Bbb C}^{n(\lambda)}$, $U_{\lambda} \in
\hat{N} $ and $S_{\lambda}$ is some trivial representation of dimension
$n(\lambda) = 1,2,3,...,\infty \quad $.

For every measurable bounded function $f$ on $\hat{N}$, we define
$P(f)$ as the diagonal operator of multiplication by $f(\lambda), \lambda \in \hat{N} $. It is easy to see that:

(1) $(T,P)$ is a unitary representation of $\hat{N}$,

(2) the measure $\mu$ is quasi-invariant with respect to the $G$-action, and

(3) if $T$ is irreducible, $\mu$ is an ergodic measure.

Therefore $\mu$ must be concentrated on a $G$-orbit in $\hat{N}$, say in
$X$. Then $X$ is a homogeneous $G$-space. We take a point $\langle
\sigma\rangle \in X$. Following the criterion of inductibility, there
exists some $\langle \tau\rangle \in \hat{G_{\sigma}}$, such that
 $\tau|_N \simeq \mult \sigma $ and finally $T \simeq \Ind_{G_{\sigma}}^G
\tau \quad $.
\end{pf}

\subsection{Projective representations and Mackey obstructions }

\begin{defn}
Let us denote by ${\cal H}$
a Hilbert space, ${\Bbb U}({\cal H}
)$ the group of unitary operators and ${\Bbb S}^1.Id = cent {\Bbb U}({\cal H}
)$. We define a  {\it projective unitary representation } 
\index{representation!unitary projective -}
$T$ to be a homomorphism $T: G \rightarrow {\Bbb U}({\cal H}
)/{\Bbb C}$, such that the map $G \times {\cal H}
\rightarrow {\cal H}
$ is continuous.
\end{defn}

\begin{rem} It is easy to see that a projective unitary representation $T$ can be
always lifted just to a continuous map, but not necessarily a homomorphism,$ T: G
\rightarrow {\Bbb U}({\cal H} )$, such that $$T(g_1)T(g_2) = c(g_1,g_2)T(g_1g_2)
\quad,$$ for some ${\Bbb S}^1$-valued function $c(.,.):=
T(g_1)T(g_2)T(g_1g_2)^{-1}$.  \end{rem}

From the properties of the homomorphism, it is easy to deduce that
$$c(g_1,g_2)c(g_1g_2,g_3) = c(g_1,g_2g_3)c(g_2,g_3) \quad.$$ Thus
$c(.,.)$ is a 2-cocycle in $H^2(G;{\Bbb T})$. This cocycle will be a
co-boundary if there is some function $b: G \rightarrow {\Bbb T}$, such
that $$c(g_1,g_2) = {b(g_1)b(g_2) \over b(g_1g_2)} \quad.$$ In this case
we have $$b(g_1)T(g_1).b(g_2)T(g_2) = b(g_1g_2)T(g_1g_2)$$ and therefore
$b(.)T(.)$ will be a unitary representation of $G$. In general case, the
cohomology class $\langle c(.,.)\rangle  \in H^2(G;{\Bbb T})$ is called the {\it
obstruction to lifting a projective representation to a unitary ( linear )
representation. }

Now we consider the situation of locally compact groups, having (proper)
closed normal subgroups, say $N$.

\begin{thm} Suppose that for every $h \in H$, the representation
$\sigma_h, \sigma_h(n) = \sigma(hnh^{-1})$, $\forall n \in N$, is
equivalent to $\sigma$. Then there exists a one - to - one correspondence
between the classes $\langle \tau\rangle  \in \hat{H}$, such that $\tau|_N \simeq \mult
\sigma $, and the set $\widehat{(H/N)}_{proj}$ of unitary equivalence
classes of irreducible projective representations of $H/N$.
\end{thm}

\begin{rem} If one is using the Mackey method of small subgroups $H =
G_{\sigma}$ the cocycles corresponding to the projective representations
from $\widehat{(H/N)}_{proj}$ are called the {\it Mackey obstructions }. 
\index{Mackey obstruction}
{\it If the Mackey obstructions vanish }, then we can obtain the dual
object $\hat{G}$ from the dual objects $\widehat{(G_{\sigma}/N)}$ in
killing the Mackey obstructions. And so we can use the induction process
on group dimension.  In the orbit method, this will be done by taking
either the so called ${\Bbb Z}/2{\Bbb Z}$-covering ( see the next sections
\S6 - \S8 ), or the so called $U(1)$-covering ( see Appendices A2. - A3. ). 
\end{rem}

{\sc Proof of the theorem. }

For every $k \in K:= H/N$, choose a representation $\rho(k) \in H$, for
example by a Borel section, such that every element $h \in H$ has a
unique decomposition $$ h = \rho(k).n, \forall k \in K, \forall n \in N
\quad.$$ Every $\langle \tau\rangle  \in \hat{H}$ such that $\tau|_N \simeq \mult \sigma
$, by definition, is of the form $$\tau(n) = I_{V_1} \otimes \sigma(n)$$
, acting on $V_1 \otimes V_2$. Because $\sigma \simeq \sigma_h, \forall
h \in H $, there exists a family of projective unitary operators $W(k), k
\in K $, such that $$\sigma(\rho(k)n\rho(k)^{-1}) =
W(k)\sigma(n)W(k)^{-1} \quad.$$ The operator $\tau(\rho(k))(I_{V_1}
\otimes W(k)^{-1})$ commutes with all the operators $\tau(n)$ and because
$\sigma$ is irreducible, $\tau(\rho(k))(I_{V_1} \otimes W(k)^{-1}) = S(k)
\otimes I_{V_2} $. Thus for $h = \rho(k)n$, we have $$\tau(h) = S(k)
 \otimes W(k) \sigma(n) \quad.$$ Now we verify that $S(.): K \rightarrow
{\Bbb U}(V_1)/{\Bbb T} $ is a projective representation. We see, on one
hand, that $\rho(k_1)\rho(k_2) = \rho(k_1k_2) \mod N $. Therefore,
$$\tau(\rho(k_1)\rho(k_2)) = S(k_1k_2) \otimes W(k_1k_2)\sigma(n)$$ for
some $ n \in N$. On the other hand, we have
$$\tau(\rho(k_1))\tau(\rho(k_2)) = \tau(\rho(k_1)\rho(k_2)) = S(k_1)S(k_2)
\otimes W(k_1)W(k_2) \quad.$$ Remark that $A \otimes B = A' \otimes B'$
if and only if there exists some constant $\lambda \in {\Bbb C}$ such that
$A' = \lambda A, B' = \lambda^{-1}B$. Thus from $$S(k_1k_2) \otimes
W(k_1k_2)\sigma(n) = S(k_1)S(k_2) \otimes W(k_1)W(k_2)$$ and from the fact
that $W$ is a projective family, we can conclude that $$W(k_1)W(k_2) =
\lambda(k_1,k_2)^{-1}W(k_1k_2)\sigma(n)$$ and $$S(k_1)S(k_2) =
\lambda(k_1,k_2)S(k_1k_2) \quad.$$ The last equation proves that $S(.)$
is a unitary projective representation of $K = H/N$. The irreducibility
and the inverse direction is easy.

\begin{rem} The generalized orbit method proposes a multidimensional
generalization of the idea of the Mackey method to the general situation.
The rest of this work is devoted to explain this idea.
\end{rem}

\section{Symplectic Manifolds with Flat Action of Lie Groups}

In this section we recall the definition of symplectic manifolds and then show that
locally, every symplectic manifold with a flat action of a Lie group can be
considered as some co-adjoint orbit.

\subsection{Flat action } 

\begin{defn} A {\it symplectic manifold } 
\index{manifold!symplectic -} $(M,\omega)$ is a real smooth
manifold $M$ jointed with a {\it symplectic structure } 
\index{structure!symplectic} $\omega$, 
i.e. a closed, non-degenerate differential 2-form $\omega$.
\end{defn}

It is easy to deduce from this definition the following:

\begin{cor}
Every symplectic manifold has even dimension.
\end{cor}

\begin{exam}
Every cotangent bundle can be transformed into a symplectic manifold.
\end{exam}

To show this it is enough to construct a closed non-degenerate differential 2-form
$\omega$. We shall construct the so called {\it Liouville form } 
\index{form!Liouville -} $\sigma$. Its
differential $B = d\sigma$ will be the desired symplectic structure.

Let us consider an arbitrary differentiable manifold $N$ and its cotangent bundle
$T^*N$. If $\tilde{U} \subset N$ is an arbitrary local coordinate chart with the
coordinate functions $q^1,q^2, \dots, q^k$, the cotangent space $T^*_mN,m \in N$,
has the linear coordinates $p_1,p_2,\dots,p_k$, dual to the basis ${\partial
\over \partial q^1},\dots,{\partial \over \partial q^k}$ of $T_mN$.

Let us denote by $p: T^*N \rightarrow N $ the natural projection. Then
$$p^{-1}\tilde{U} = U \approx \tilde{U} \times {\Bbb R}^k \hookrightarrow T^*N$$ is
a local coordinate chart of $T^*N$ with coordinates $q^1,q^2,\dots
,q^k,p_1,p_2,\dots,p_k$.We define the restriction $\sigma|_U$ of form $\sigma$ by
$$\sigma|_U:= \sum_{i=1}^k p_idq^i $$

Let $\tilde{V}$ be another coordinate chart in $N$, and $V:= p^{-1}\tilde{V}$ the
corresponding coordinate chart in $T^*N$ with the coordinate functions $\tilde{q}^1
, \tilde{q}^2, \dots \tilde{q}^k, \tilde{p}_1,\dots,\tilde{p}_k $. Then in the
intersection $U \cap V$ we have the relations 
$$\begin{array}{rl}\tilde{q}^i &= \tilde{q}^i(q^1,\dots q^k) \quad, i = 1,\dots, k
\cr 
 \tilde{p}_i &= \tilde{p}_i(q^1,\dots,q^k,p_1,\dots,p_k) \quad, i= 1,\dots,k
\end{array}$$ and $\tilde{p}_i, i = 1,\dots, k$, are linear functions of
$p_1,p_2,\dots,p_k$ of form $$\tilde{p}_i = {\sum_{j=1}^k} {\partial q^j \over
\partial\tilde{q}^i} p_j \quad.$$ Therefore, we have $$\begin{array}{rl}\sigma|_{V \cap U}
&= \sum_{i=1}^k \tilde{p}_id\tilde{q}^i \cr
    &= \sum_i (\sum_j {\partial q^j \over \partial \tilde{q}^i }p_j)(\sum_l
{\partial\tilde{q}^i \over \partial q^l} dq^l) \cr
  &= \sum_j \sum_l (\sum_i {\partial q^j \over \partial \tilde{q}^i }{\partial
\tilde{q}^i \over \partial q^l }) p_j dq^l \cr
  &= \sum_j \sum_l \delta_{jl}p_jdq^l \cr
  &= \sum_j p_jdq^j = \sigma|_{U \cap V}.\end{array} $$ Therefore $\{ \sigma|_U =
\sum_{i=1}^k p_idq^i \}$ defines a 1-form $\sigma \in \Omega^1(T^*N)$. This form is
the so called {\it Liouville form }. \index{form!Liouville -} Its 
differential $B = d\sigma$ has the
following local coordinate form $$B|_U = d\sigma|_U = \sum_{i=1}^k dp_i \wedge dq^i
\quad.$$ Therefore, it is a non-degenerate differential 2-form, $B \in
\Omega^2(T^*N)$ \quad. This example is in the foundation of the standard
Hamiltonism of classical mechanical systems.

\begin{defn}
A vector field $\xi \in \Vect(M)$ is called {\it Hamiltonian } 
\index{vector field!Hamiltonian -}
and denoted $\xi \in \Vect(M,\omega)$ iff holds one of the following equivalent 
conditions:
\begin{enumerate}
\item[(i)] The Lie derivative of $\omega$ along the field $\xi$ vanishes,
$$\xi\omega:= L_{\xi}\omega:= Lie_{\xi}\omega = 0 \quad,$$
\item[(ii)] $\imath(\xi)\omega$ is a closed 1-form.
\end{enumerate}
\end{defn}

A vector field $\xi \in \Vect(M)$ is called {\it strictly Hamiltonian } 
\index{vector field!strictly Hamiltonian -} and denoted
$\xi \in \Vect_0(M,\omega)$ if the 1-form $\imath(\xi)\omega$ is exact, i.e. there
exists the so called {\it generating function } 
\index{function!generating -} $f = f_{\xi}$, such that
$$\imath(\xi)\omega + df_\xi = 0 \quad.$$ In this case one says that $\xi = \xi_f$
is the (strictly) Hamiltonian vector field,{\it corresponding } to the function $f$
, or the {\it symplectic gradient } \index{gradient!symplectic -} of $f$.

\begin{rem} On the symplectic manifold $(M,\omega)$, there is a
one - 
to - one correspondence between the vector fields $\xi \in \Vect(M)$ and the 1-forms
$\imath(\xi)\omega \in \Omega^1(M)$. Therefore we have the following diagram of two
exact sequence of vector spaces 
$$\vbox{\halign{ #&#&#\hfill \cr
   &  &$0$   \cr
   &  &$\uparrow$   \cr
$0 \quad \rightarrow \quad {\Bbb R} \quad \rightarrow Ve$&$ct_0$&$(M,\omega)\quad \rightarrow \quad \Vect(M,\omega) \quad \rightarrow \quad H_{DR}^1(M,\omega) \quad \rightarrow 0$  \cr
   &  &$\uparrow$   \cr 
   &$C^{\infty}$&$(M,{\Bbb R})$   \cr
   &  &$\uparrow$   \cr
   &  &${\Bbb R}$   \cr 
   &  &$\uparrow$   \cr
   &  &$0$          \cr}}$$
\end{rem}

\begin{prop} The vertical and the horizontal sequences are in fact the exact
sequence of Lie algebras.  \end{prop} \begin{pf} First we prove that
$$[\Vect(M,\omega),\Vect(M,\omega)] \subseteq \Vect_0(M,\omega). $$ Really, let $\xi,
\eta \in \Vect(M,\omega)$. Consider the function defined by $$f = f_{[\xi, \eta]}
:= \imath(\xi)\imath(\eta)\omega = \omega(\xi, \eta) \quad.$$ From differential
geometry, we know that $$\imath([\xi,\eta]) = L_{\xi}\circ\imath(\eta) -
\imath(\eta)\circ L_{\xi}$$ and $$L_{\xi} = d\circ\imath(\xi) + \imath(\xi)\circ d.$$

 Then$$\begin{array}{rl} \imath([\xi,\eta])\omega &= L_{\xi}\imath(\eta)\omega -
\imath(\eta)L_{\xi}\omega \cr
  &= L_{\xi}\imath(\eta)\omega \cr
  &= d\imath(\xi)\imath(\eta)\omega + \imath(\xi)d\imath(\eta)\omega =
d\imath(\xi)\imath(\eta)\omega \cr
  &= -d f_{[\xi,\eta]}.\end{array} $$

Thus the quotient Lie algebra $\Vect(M,\omega)/Vetc_0(M,\omega)$ is commutative and
is isomorphic to the commutative Lie algebra $H_{DR}^1(M;{\Bbb R})$.

Recall that for two functions $f,g \in C^{\infty}(M,{\Bbb R})$, their {\it Poisson
brackets} \index{brackets!Poisson -} is defined as $$\{ f,g \}:= \xi_f g = - 
\xi_g f = \omega(\xi_f,\xi_g) =
f_{[\xi,\eta]} \quad.$$ Thus $$d \{ f,g \} = df_{[\xi_f,\xi_g]} =
\imath([\xi,\eta])\omega$$ and the vertical sequence is also an exact sequence of
Lie algebras. The proposition is proved.
\end{pf}

Let us consider now the smooth group action of a Lie group $G$ on $(M,\omega)$.The
Lie algebra ${\frakt g} = Lie G$ acts on $M$ by the infinitesimal action. More
precisely, for every $x \in M$ and $X \in {\frakt g} $, the one - parameter group
$\exp (tX)$ provides a smooth curve, passing through the point $x \in M$. Let thus
denote this curve by $\exp (tX) x $ and its tangent vector at $x$ by $\xi_X(x)$,
$$\xi_X(x):= {d \over dt}|_{t = 0} \exp (tX) x \quad.$$

Suppose that $\xi_X \in \Vect_0(M,\omega), \forall X \in {\frakt g}$ with generating
function $f_X$. The quantity $$c(X,Y):= \{ f_X,f_Y\} - f_{[X,Y]}$$ can be
considered as the {\it curvature of the group action } \index{group 
action! curvature of -} of $G$ on $M$.

\begin{defn}
The action of $G$ on $M$ is said to be {\it flat } \index{action!flat -}
iff the curvature $c(.,.)$ of the $G$-action vanishes, i.e.
$$ \{ f_X,f_Y\} - f_{[X,Y]}\equiv 0, \forall X,Y \in {\frakt g} \quad.$$
\end{defn}

\begin{rem} If the action of $G$ on $M$ is flat, the homomorphism ${\frakt g}
\rightarrow \Vect(M,\omega)$ can be lifted to a homomorphism ${\frakt g} \rightarrow
C^{\infty}(M,{\Bbb R})$, following the commutative diagram $$\vbox{\halign{
#&#&#&#&#&#\hfill \cr
   &  &$0$ &  &  &  \cr
   &  &$\uparrow$ & \cr $0 \rightarrow Ve$&$ct_0$&$(M,\omega)$&$
\rightarrow$&$\Vect$&$(M,\omega) \rightarrow H^1_{DR}(M,{\Bbb R}) \rightarrow 0$ \cr
   &  &$\uparrow $&$\nwarrow$&  &$\uparrow$   \cr
   &$C^{\infty}$&$(M,{\Bbb R})$&$\leftarrow$&  &${\frakt g}$ \cr 
   &            &$\uparrow$    &               &  &  \cr    
   &            &${\Bbb R}$    &               &  &  \cr
   &            &$\uparrow$    &               &  &  \cr
   &            &$0$             &               &  &  \cr}}$$
\end{rem}

\begin{exam} Every co-adjoint orbit ( or, following A.A. Kirillov, K-orbit ) is a
simplectic manifold with a flat action.

\end{exam}

Let us consider a connected and simply connected Lie group $G$ with Lie algebra
${\frakt g} = Lie G$ and the dual vector space ${\frakt g}^*:= Hom_{\Bbb R}({\frakt
g},{\Bbb R})$. With every element $g \in G$ we can associate a map $$A(g): G
\rightarrow G \quad,$$ $$A(g)x:= gxg^{-1} \quad,$$ fixing the identity element $e
\in G$ $$A(g)e = geg^{-1} = e \quad.$$ The differential of this map $$Ad(g):=
A(g)_*: {\frakt g} \cong T_eG \rightarrow {\frakt g} \cong T_eG $$ is called the
{\it adjoint representation } \index{representation!adjoint -} of $G$ in 
{\frakt g}, $$Ad: G \rightarrow \Aut{\frakt
g} \quad.$$ The corresponding {\it contragradient } representation 
\index{representation!contragradient} of $G$ in
${\frakt g}^*$ is called the {\it co-adjoint representation } 
\index{representation!co-adjoint} of $G$, $$\begin{array}{rl} K
:= coAd &: G \longrightarrow \Aut{\frakt g}^* \quad,\cr K(g):= (Ad(g^{-1}))^* &:
{\frakt g}^* \longrightarrow {\frakt g}^*.\end{array}$$ Under this co-adjoint action
of $G$ on ${\frakt g}^* $, the space ${\frakt g}^*$ is divided on
to the so called {\it co-adjoint } \index{orbit!co-adjoint -} 
( or, following A. A. Kirillov, the $K$- ) {\it orbits }
, ${\frakt g}^*/G:= {\cal O}
(G)$.

Remark that the {\it adjoint } action \index{action!adjoint -} of Lie 
algebra ${\frakt g}$ on itself,
$$(adX)(Y):= [X,Y] $$ is the differential of the {\it adjoint } action 
\index{action!adjoint -} of $G$ on {\frakt g}, $$ Ad(\exp X) = e^{adX} 
\quad.$$

Let us now consider a fixed $K$-orbit $\Omega \in {\cal O}(G)$. We fix a point $F
\in \Omega$ and consider the stabilizer $G_F$ at this point, and the bilinear form
$\langle F,[.,.]\rangle $ on {\frakt g}. It is easy to see that the kernel of this form is
coincided with the Lie algebra ${\frakt g}_F:= Lie G_F $ of the stabilizer $G_F$,
symbolically, $$Ker \langle F,[.,.]\rangle  = {\frakt g}_F \quad.$$ Therefore $\langle F,[.,.]\rangle $
induces a non-degenerate bilinear form $\omega_F(.,.)$ on the tangent at $F$ space
$T_F\Omega \cong {\frakt g/g}_F$. We write the action of $G$ on $\Omega = \Omega_F$
on the right. So it is natural to identify $\Omega$ with the coset space $G_F
\setminus G$, and the tangent space $T_F\Omega$ at $F$ with the quotient space
${\frakt g}_F \setminus {\frakt g} \cong {\frakt g/g}_F$. The orbit $\Omega$ is a
homogeneous right $G$-space and $\omega_F(.,.)$ is a non-degenerate skew form on the
tangent space $T_F\Omega = {\frakt g/g}_F$.

Remark that {\it the form $\omega_F: T_F\Omega \times T_F\Omega \rightarrow {\Bbb
R}$ is $AdG_F$-invariant }, i.e. $\forall \tilde{X} = X + {\frakt g}_F, \tilde{Y}
= Y + {\frakt g}_F $ and $ \forall g \in G_F$, $$\begin{array}{rl} \langle K(g)F,[X,Y]\rangle  &=
\omega_F(Adg^{-1}\tilde{X},Adg^{-1}\tilde{Y}) \cr
  &= (Adg)^*\omega_F(\tilde{X},\tilde{Y}) \quad,\end{array}$$ i.e.$$\langle F,[X,Y]\rangle  =
\omega_F(\tilde{X},\tilde{Y}) \quad.$$ Therefore $\omega_F$ can be extended to a
$G$-invariant 2-form $\omega_{\Omega}$ on $\Omega$.

We prove now that {\it this form $\omega_F$ is closed. } Really, $$\begin{array}{rl}
d\omega_F(\tilde{X},\tilde{Y},\tilde{Z}) &:=
\omega_F([\tilde{X},\tilde{Y}],\tilde{Z}) -
\omega_F([\tilde{X},\tilde{Z}],\tilde{Y}) +
\omega_F([\tilde{Y},\tilde{Z}],\tilde{X}) \cr
  &= \langle F,[[X,Y],Z] + [[Y,Z],X] +[[Z,X],Y]\rangle  \equiv 0,\end{array}$$ following the {\it
Jacobi identity } for Lie algebras. Thus $\omega_F$ is a $G$-invariant symplectic
structure on orbit $\Omega = \Omega_F$.

Finally, we show that {\it the $G$-action is flat. } 
\index{action!flat -} Really, every $X \in {\frakt
g}$ can be considered as a linear function $f_X:= \langle .,X\rangle $ on ${\frakt g}^*$, and
therefore a function on $\Omega$. Hence, $$df_X(F) = \langle F,[.,X]\rangle  =
-\imath(\xi_X)\omega_F$$ and $$\begin{array}{rl} f_{[X,Y]} &= \langle .,[X,Y]\rangle  = L_X\langle .,Y\rangle  \cr
                      &= L_Xf_Y = \{ f_X,f_Y\}.\end{array}$$ Thus, every $K$-orbit is a
$G$-homogeneous symplectic manifold with the flat co-adjoint action of $G$.

The space ${\frakt g}^*$ is decomposed into a disjoint union of the homogeneous
symplectic $G$-manifolds, with the flat co-adjoint $G$-action. Together they can be
considered as the symplectic leaves of a $G$-homogeneous {\it Poisson 
structure } \index{structure!Poisson -} on
${\frakt g}^*$. This idea goes back to the classical one from Sophus Lie.

\subsection{Classification} 

\begin{thm}[\bf Classification ] Every homogeneous symplectic manifold
$(M,\omega)$ with flat action of a Lie group $G$ is locally isomorphic to
an $K$-orbit of $G$ or a central extension of $G$ by ${\Bbb R}$.
\end{thm}

\begin{pf} Consider the map $\phi: M \rightarrow {\frakt g}^*$, defined by
$$\langle \phi(m),X\rangle := F_X(m) \quad,$$ where $F_X, X \in {\frakt g}$ are the generating
functions of $\xi_X, X \in {\frakt g} $. To verify that this map is $G$-equivalent
, in virtue of connectedness of $G$, it suffices to verify that $\phi$ commutes with
the action of elements of type $\exp{Y}, Y \in {\frakt g}$ on $M$ and on ${\frakt
g}^*$.  From the flatness of the $G$-action, $\{ F_X,F_Y \} = F_{[X,Y]}$, one
deduces that $\phi_*$ maps $\xi_Y$ on $M$ to $\eta_Y$ on ${\frakt g}^*$. This
suffices for the local equivariance and hence also for global equivariance, as
remarked before.  From this $G$-equivariance, $\phi(M)$ is an $K$-orbit in
${\frakt g}^*$.  It is easy to see that,if $X_1,X_2,\dots,X_{2k} \in {\frakt g}$
and the corresponding Hamiltonian fields $\xi_{X_1}, \dots, \xi_{X_{2k}}$ are
linear independent, then so are also the differentials of their generating
functions $df_{X_1},\dots, df_{X_{2k}}$. Therefore, $$\phi: M \twoheadrightarrow
\phi(M) = \Omega \subset {\frakt g}^*$$ is a local homeomorphism. This means that
$M$ is a covering of a $K$-orbit $\Omega \subset {\frakt g}^* $. If this $K$-orbit
$\Omega$ is simply connected, the covering is unique and $M$ is homeomorphic to
$\Omega$. In other case, $\pi_1(\Omega) \ne \{1\}$, $M$ is a connected covering
of $\Omega$, which is in one - to - one correspondence with subgroups $\Gamma$ of
$\pi_1(\Omega)$, and $M$ is homeomorphic to $\Gamma \setminus \tilde{\Omega}$,
where $\tilde{\Omega}$ is the universal covering of $\Omega$.

Let $\tilde{G}$ and $\tilde{M}$ are the corresponding universal coverings of $G$ and
$M$, respectively. Consider the Lie algebra ${\frakt g}_1$, generated by the
generating functions $f_X,X \in {\frakt g}$. it is naturally an extension of
{\frakt g} by ${\Bbb R}$,$$ 0 \rightarrow {\Bbb R} \rightarrow {\frakt g}_1
\rightarrow {\frakt g} \rightarrow 0 \quad.$$ We have a commutative diagram of Lie
algebras and homomorphisms $$\vbox{\halign{ #&#&#&#&#&# \hfill \cr
   &  &$0$    &    &    &   \cr
   &  &$\uparrow$ & & & \cr $0 \rightarrow Ve$&$ct_0$&$(\tilde{M},\omega)$
&$\buildrel \cong \over \rightarrow$ &$Ve$&$ct(\tilde{M},\omega) \rightarrow
H_{DR}^1(\tilde{M} ;{\Bbb R}) = 0$ \cr
  &  &$\uparrow$ &$\nwarrow$ &  &$\uparrow$  \cr
  &$C^{\infty}$&$(\tilde{M},{\Bbb R})$ &$\leftarrow$ &  &${\frakt g}_1$  \cr
  &  &$\uparrow$   &  &   &   \cr
  &  &${\Bbb R}$   &  &   &   \cr
  &  &$\uparrow$ &  &  &  \cr
  &  &$0$ & & & \cr }}.$$ Thus $\tilde{M}$ is a homogeneous symplectic manifold
with flat action of the simply connected Lie group $G_1$, corresponding to Lie
algebra ${\frakt g}_1$. The theorem is proved.
\end{pf}

\begin{rem} The classification theorem means that every Hamiltonian system with flat
action of a connected Lie group $G$ is locally isomorphic to a co-adjoint orbit of
$G$ or its central extension by ${\Bbb R}$. Roughly speaking, locally every flat
homogeneous Hamiltonian system is a co-adjoint orbit.
\end{rem}

\section{Prequantization}

In this section we define the so called procedure of quantization, the
rule of geometric quantization and show its application to the
representation theory. Our main intention is to do all in the
multidimensional context. In this situation, there arise s some
noncommutative summand in the expression of the corresponding curvature.

\subsection{Quantization procedure}

Let us consider a fixed co-adjoint orbit $\Omega \in {\cal O} (G)$, a
fixed point $F \in \Omega$, the stabilizer $G_F$ at the point $F$ and its
connected component of identity $(G_F)_0$. Recall that $\Omega \approx
G_F \setminus G$ is a right homogeneous $G$-space.

\begin{rem} $$\chi_F(\exp{(.))}:= \exp{({i \over \hbar}\langle F,.\rangle )}: (G_F)_0
\rightarrow {\Bbb T} = {\Bbb S}^1 \subset {\Bbb C}$$ {\it defines a
unitary character ( i.e. a 1-dimensional representation ) 
\index{character!unitary -} of $(G_F)_0$, }
where $\hbar:= {h \over 2\pi}$ is the so called {\it normed Planck
constant } \index{constant!normed Planck -}and $h$ is the {\it un-normed 
Planck constant}. \index{constant!un-normed Planck -} For the
mathematical theory exposed here, the concrete value of $\hbar$ does not
play any role. The only importance is that its value is nonzero. We use
it in the correspondence with the corresponding physics theories.
Normally, we can suppose that
 $\hbar = 1$.
\end{rem}

Really, we have $${\frakt g}_F = Lie G_F = Ker\langle F,[.,.]\rangle  \quad.$$ Then the linear functional ${i \over \hbar}\langle F,.\rangle $ is also a character ( i.e. 1-dimensional representation ) of ${\frakt g}_F$, i.e. $${i \over \hbar}F,\langle [{
\frakt g}_F,{\frakt g}_F]
\rangle \equiv 
0 \quad.$$

\begin{defn} The orbit $\Omega$ is said to be {\it admissible } ( resp.,
{\it integral }), \index{orbit!admissible integral -} iff there exists a 
unitary representation $\sigma$ of
the whole $G_F$,such that its restriction on $(G_F)_0$ is a multiple of
$\chi_F$,$$\sigma|_{(G_F)_0} \simeq \mult \enskip \chi_F \quad,$$ (resp.
, iff $\chi_F$ can be extended to a unitary character r of $G_F$ ). 
\end{defn}

\begin{rem} In virtue of the fact that $(G_F)_0$ is a normal closed
subgroup of $G_F$, the set of such $\sigma$ is in one - to - one
correspondence with the projective representations of the quotient group
$G_F/(G_F)_0 $, $$\{ \sigma \in \hat{G}_F ; \sigma |_{(G_F)_0 } \simeq
\mult \enskip \chi_F \} {\buildrel 1-1 \over \longleftrightarrow}
\widehat{(G_F/(G_F)_0)}_{proj} \quad. $$
\end{rem}

\begin{rem} It is easy to see that $G_F$ can be included in the symplectic
group of the tangent space ${\frakt g/g}_F$ of the orbit $\Omega_F$, $G_F
\hookrightarrow Sp({\frakt g/g}_F)$. The last group has the well-known
metaplectic two-fold covering $$1 \rightarrow
 {\Bbb Z}/2{\Bbb Z} \rightarrow Mp({\frakt g/g}_F) \rightarrow Sp({\frakt
g/g}_F) \rightarrow 1\quad.$$ Using the 5-homomorphism lemma, we can
construct the unique two-fold covering of the stabilizer $G_F$ such that
the following diagram is commutative $$\vbox{\halign{ #&#&#&#\hfill \cr
$1 \rightarrow {\Bbb Z}/$&$2{\Bbb Z} \rightarrow $&$G_F^{\frakt g}
\rightarrow $&$G_F \rightarrow 1$ \cr
  &$\downarrow$ &$\downarrow$ & $\downarrow$ \cr
$1 \rightarrow {\Bbb Z}/$&$2{\Bbb Z} \rightarrow \Mp$&$({\frakt g/g}_F) \rightarrow \Sp$&$({\frakt g/g}_F) \rightarrow 1$  \quad,\cr}}$$
on one hand. On the other hand, The symplectic group has also the well-known $U(1)$-covering $$ 1 \rightarrow U(1) \rightarrow Mp^c({\frakt g/g}_F) \rightarrow Sp({\frakt g/g}_F) \rightarrow 1\quad.$$  Also using the 5-homomorphism lemma, we can construct

 the unique $U(1)$-covering of the stabilizer $G_F$, such that the following diagram is commutative
$$\vbox{\halign{ #&#&#&# \hfill \cr
$1 \rightarrow U$&$(1) \rightarrow G$&$_F^{U(1)} \rightarrow G$&$_F \rightarrow 1$ \cr
   &$\downarrow$ &$\downarrow $&$\downarrow $  \cr
$1 \rightarrow U$&$(1) \rightarrow \Mp^c$&$({\frakt g/g}_F) \rightarrow \Sp$&$({\frakt g/g}_F) \rightarrow 1$ \quad. \cr}}$$
\end{rem}

Therefore there are two subsets $\widehat{G_F^{{\Bbb Z}/2{\Bbb Z}}}$ and
$\widehat{G_F^{U(1)}}$ of unitary representations in the projective dual
object $\widehat{(G_F/(G_F)_0)}_{proj}$.

In the future, avoiding the Mackey obstructions of the construction, we
shall use the two - fold covering in the Duflo's construction ( Sections
\S\S6-8 ) and the $U(1)$-covering in the later development ( Appendices
A1-A3 ). 

\begin{prop}  The orbit $\Omega$ is integral if and only if the
cohomology class $[\omega_{\Omega}]$ of the Kirillov's form
$\omega_{\Omega}$ is integral ; i.e. $$[\omega_{\Omega}] \in H^2(\Omega;
{\Bbb Z}).$$
\end{prop}

\begin{pf} 
Recall that we suppose always $G$
to be a connected and simply connected Lie group, and $G_F$ the closed
subgroup. Therefore,$\Omega \approx G_F \setminus G$ and we have the
long exact sequence of cohomology groups $$0 \rightarrow H^0(\Omega;{\Bbb
R}) \rightarrow H^0(G;{\Bbb R}) \rightarrow H^0(G_F;{\Bbb R})
\rightarrow$$ $$\rightarrow H^1(\Omega;{\Bbb R}) \rightarrow H^1(G;{\Bbb
R}) \rightarrow H^1(G_F;{\Bbb R}) \rightarrow$$ $$\rightarrow
H^2(\Omega;{\Bbb R}) \rightarrow H^2(G;{\Bbb R}) \rightarrow H^2(G_F;{\Bbb
R}) \rightarrow \dots.$$ We have $$H^0(G;{\Bbb R}) = 0 \quad,$$
$$H^1(G;{\Bbb R}) = 0 \quad.$$ Then $$H^2(\Omega;{\Bbb R}) \cong
H^1(G_F;{\Bbb R})$$ and $$\pi_1(\Omega) \cong \pi_0(G_F) = G/(G_F)_0 \quad
.$$

Let us denote by $$p:  G \longrightarrow \Omega \approx G_F \setminus G
$$ the natural projection. Then $$p^*\omega_{\Omega}.Id = d\alpha $$ for
some $\alpha \in H^1(G; {\Bbb Z})$. This form $\alpha$ is unique up to a
differential, say $df, f \in C^{\infty}(G,{\Bbb R})$ and $\alpha_0:=
\alpha|_{(G_F)_0}$ is closed. Really
, $ \alpha_0 = {i \over \hbar }\langle F,.\rangle $.

Now suppose that the orbit $(\Omega,\omega_{\Omega})$ is admissible.
Therefore, there exists $(\tilde{\sigma},\tilde{V}) \in \hat{G}_F$ such
that $$\tilde{\sigma}|_{(G_F)_0} = \mult \enskip \chi_F \quad.$$Let us
consider the action $\Ind_{G_F}^G(\tilde{\ sigma},\tilde{V})$ on the
sections of th associated bundle ${\cal E} _{\tilde{V}}(\Omega):= G
\times_{G_F} (\tilde{\sigma},\tilde{V}) $.
\end{pf}

\begin{defn}
 Under a {\it procedure of quantization }, \index{procedure of 
quantization} or some time a {\it rule of
quantization } \index{rule of quantization} we mean a correspondence 
associating to each function $f \in
C^{\infty}(\Omega,{\Bbb C})$ an Hermitian operator $\hat{f}$ ( which
become anti-auto-adjoint for real functions $f$ ) in a Hilbert space ${\cal
H}$, satisfying the so called {\it commutation relations } 
\index{commutation relation} what follow
$$\vbox{\halign{ #&#\quad &#\quad &#\hfill \cr $\{ \widehat{f,g}$&$\}$
&$=$ &${i \over \hbar}[\hat{f}_1,\hat{f}_2],\forall f_1,f_2 \in
C^{\infty}(\Omega,{\Bbb C})$ \cr
  &$\hat{1}$ &$=$ &$Id$ \quad.\cr}}$$ In this case the operators
$\hat{f}, f \in C^{\infty}(\Omega,{\Bbb C})$ will be called the {\it
quantized operators.}
\end{defn}

Recall that for the associated ( induced ) bundle ${\cal E}_V(\Omega)$ we
have to suppose that the principal bundle $$\vbox{\halign{ #&#\hfill \cr
$G_F \rightarrow$&$G$ \cr
                 &${\downarrow}\Gamma$  \cr
                 &$\Omega \approx G_F \setminus G$ \cr}}$$ is locally
trivial and that we fix a connection $\Gamma$ ( i.e. a trivialization ).
Then with every representation $(\tilde{\sigma},\tilde{V})$ of $G_F$ we
can associate an ( affine ) connection
 $\nabla$ on the induced bundle with the connection form $\alpha$.

\begin{defn} For every $f \in C^{\infty}(\Omega,\omega_{\Omega})$, denote
its hamiltonian vector field by $\xi_f$, $$ \imath(\xi_f)\omega + df = 0
\quad.$$ The operators $$\hat{f}:= f + {\hbar \over i}\nabla_{\xi_f} = f
+ {\hbar \over i}L_{\xi_f} + \alpha(vert(\xi _f))$$ are called the {\it
geometric quantization operators }, corresponding to $f \in
C^{\infty}(\Omega, {\Bbb C})$, where $f$ on the right hand side is the
operator of multiplication by $f$ and $vert(\xi_f)$ is the vertical
componenent of the vector field $\xi_f$, following the fixed connection
$\Gamma$ on the principal bundle. 
\end{defn}

\begin{thm}  The following three conditions are equivalent.
$$\xi\alpha(\eta) - \eta\alpha(\xi) - \alpha([\xi,\eta]) + {i \over
\hbar}[\alpha(\xi),\alpha(\eta)] = -\omega(\xi,\eta)Id,\forall \xi,\eta
\in \Vect(\Omega)\leqno(1)$$ $$\Curv(\nabla)(\xi,\eta) =
[\nabla_{\xi},\nabla_{\eta}] - \nabla_{[\xi,\eta]} = -{i \over
\hbar}\omega(\xi,\eta) = -{i \over \hbar}\omega(\xi,\eta)Id\leqno(2)$$ $$
\widehat{\{ f_1,f_2\}} = {i \over \hbar}[\hat{f_1},\hat{f_2}],\hat{1} =
Id\quad,\leqno(3)$$ i.e.$f \mapsto \hat{f}$ is a quantization procedure.
\end{thm}

$(1) \Longleftrightarrow (2)$ is just from differential
geometry of connections. From our definition 3.7 we have $$\nabla_{\xi}
:= L_{\xi} + {i \over \hbar}\alpha(vert(\xi))\quad.$$ For simplicity, we
denote $\alpha(vert(\xi))$ as $\alpha(\xi)$. Therefore we write
$$\nabla_{\xi} = L_{\xi} + {i \over \hbar}\alpha(\xi)\quad.$$

$$\begin{array}{rl}
 [\nabla_{\xi},\nabla_{\eta}] - \nabla{[\xi,\eta]} &= [L_{\xi} + {i
\over \hbar}\alpha(\xi),L_{\eta} + {i\over \hbar}\alpha(\eta)] -
L_{[\xi,\eta]} - {i \over \hbar}\alpha([\xi,\eta]) \cr
  &= [L_{\xi},L_{\eta}] + {i \over \hbar}[L_{\xi},\alpha(\eta)] - {i \over
\hbar}[L_{\eta},\alpha(\xi)] + ({i \over \hbar})^2
[\alpha(\xi),\alpha(\eta)] \cr
  &- L_{[\xi,\eta]} - {i \over \hbar}\alpha([\xi,\eta]) \cr
  &= \{[L_{\xi},L_{\eta}] - L_{[\xi,\eta]}\} + {i \over \hbar}\{
[L_{\xi},\alpha(\eta)] - [L_{\eta},\alpha(\xi)] \cr
  &\phantom{=}- \alpha([\xi,\eta]) + {i
\over \hbar}[\alpha(\xi),\alpha(\eta)]\} \end{array}$$ In virtue of Lie
derivation, the first summand vanishes. It rests only to prove the
following

\begin{lem}
$$ [L_{\xi},\alpha(\eta)] = \xi\alpha(\eta) \quad.$$
\end{lem}

{\it Really },applying the left hand side to a section $s$, we have
$$\begin{array}{rl} [L_{\xi},\alpha(\eta)]s &= L_{\xi}(\alpha(\eta)s) -
\alpha(\eta)L_{\xi}s \cr
   &= (L_{\xi}\alpha(\eta))s +\alpha(\eta)L_{\xi}s - \alpha(\eta)L_{\xi}s \cr
   &= \xi\alpha(\eta).s.\end{array}$$

$(2) \Longleftrightarrow (3)$ can be proved by a direct verification.
$$\begin{array}{rl} {i\over\hbar}[\hat{f}_1,\hat{f}_2] &= {i\over\hbar}[f_1 +
{\hbar\over i}\nabla_{\xi_{f_1}}, f_2 + {\hbar \over
i}\nabla_{\xi_{f_2}}] \cr
  &= {i\over \hbar }\{ {\hbar \over i }[f_1,\nabla_{\xi_{f_2}}] +
{\hbar\over i}[\nabla_{\xi_{f_1}},f_2] + ({\hbar\over
i})^2[\nabla_{\xi_{f_1}},\nabla_{\xi_{f_2}}] \} \cr
  &= [f_1,\nabla_{\xi_{f_2}}] + {\hbar \over i}\nabla_{\xi_{\{ f_1,f_2\}}}
- {\hbar \over i}\nabla_{\xi_{\{ f_1,f_2\}}} + [\nabla_{\xi_{f_1}},f_2] +
{\hbar \over i}[\nabla_{\xi_{f_1}},\nabla_{\xi_{f_2}}].\end{array}$$ From
differential geometry, we know that $$\xi_{\{ f_1,f_2\}} =
[\xi_{f_1},\xi_{f_2}] \quad.$$ We have therefore only to prove the
following

\begin{lem} $$[f_1,\nabla_{\xi_{f_2}}] = [\nabla_{\xi_{f_1}},f_2] =
-\xi_{f_2}(f_1) = \xi_{f_1}(f_2) = \omega(\xi_{f_1},\xi_{f_2}) = \{
f_1,f_2\}\quad. $$
\end{lem}

Really, Because $f_1,f_2$ are the scalar functions and
$\beta(f_1),\beta(f_2)$ are the operator valued functions,$\beta:=
{i\over\hbar}\alpha$,the connection form, we have then
$$[f_1,\beta(\xi_{f_2})] \equiv 0,$$ $$[f_1,\nabla_{\xi_{f_2}}] = [f_1,
L_{\xi{f_2}}] \quad.$$ Applying this to a section $s$, we have
$$\begin{array}{rl} [f_1,\nabla_{\xi_{f_2}}]s &= f_1(L_{\xi_{f_2}}) -
L_{\xi_{f_2}}(f_1s) \cr
  &= -L_{\xi_{f_2}}(f_1).s = \omega(\xi_{f_1},\xi_{f_2})s. \end{array}$$
The theorem is proved.
 
\subsection{Application.}

Let $X \in {\frakt g} = Lie G $ be an arbitrary element, and $\{
\exp{(tX)}\}_{t \in {\Bbb R}} \subset G$ the corresponding one-parameter
subgroup in $G$, and $\xi_X$ the corresponding strictly hamiltonian
vector field, with the generating function $f_ X:= f_{\xi_X}$ and the
corresponding Lie derivation $$L_X:= L_{\xi_X}\quad.$$

It is easy to see that $$[L_X,L_Y] = L_{[X,Y]} \quad,$$ $$L_Xf =
\{f_X,f\} $$ and $$c(X,Y) = \{ f_X,f_Y\} - f_{[X,Y]} \equiv 0 \quad,$$
in virtue of flatness of the $G$-action.

Applying the procedure of (pre)quantization, we define $$\Lambda(X):=
{i\over\hbar}\hat{f}_X = {i\over\hbar}f_X + \nabla_{\xi_X} \quad.$$ It is
easy to show that $$\begin{array}{rl} [\Lambda(X),\Lambda(Y)] &= \Lambda([X,Y]) +
c(X,Y) \cr
                                 &= \Lambda([X,Y]),\end{array}$$
following the flatness of the $G$-action.

{\bf Conclusion.} The procedure of geometric quantization provides a Lie
algebra representation $\Lambda(.)$ on the space of sections of the
induced bundle ${\cal E} _V(\Omega)$.

\begin{rem} There are some conditions, providing the corresponding Lie
group representations, namely the E. Nelson conditions: {\it If all the
operators $\Lambda(X),X \in {\frakt g} $ and the Casimir operators
$$\Lambda(C), C \in {\cal Z} ({\frakt g}) = cent U({\frakt g}) $$ have
operator closures, then $$\exp{\Lambda(X)},X \in {\frakt g}$$ provide a
unitaty representation of the universal covering group $\tilde{G}$ of the
Lie group $G$.}
\end{rem}

\section{Polarization}

In the previous section, we have discussed how to construct the so called
{\it quantum correspondence } \index{quantum correspondence} $$f \longmapsto 
\hat{f}, 
f \in C^{\infty}(\Omega,{\Bbb C}) \quad.$$ But we didn't discuss about the
arizing here Hilbert space ${\cal H}$, on which the quantized operators
$\hat{f},f \in C^{\infty}(\Omega,{\Bbb C})$ act. So we had the so called
{\it prequantization }. \index{prequantization} In general, from a 
symplectic manifold, one
constructs the $L^2$-space on `` coordinates '', but not on ``impulsion
''. The procedure of deleting from the coordinates of the symplectic
manifold the half-number of ``impulsion coordinates'' is, roughly
speaking, {\it polarization }. \index{polarization}

\subsection{Some ideas from physics }

Let us consider a general symplectic manifold $(M,\omega)$, with the
symplectic form $$\omega(x) = {\sum_{i,j = 1,i<j}^{2k}} a_{i,j}(x)dx^i
\wedge dx^j$$ in local coordinates $x = (x^1,x^2,\dots,x^{2k})$.

The so called {\it Darboux coordinates } \index{Darboux coordinates} are the 
local coordinates
$$(x^1,x^2,\dots,x^{2k}) = (q^1,q^2,\dots,q^k,p_1,\dots,p_k)$$ such that
in these local coordinates the symplectic form has the {\it canonical form
} \index{form!canonical -} $$\omega(p,q) = {\sum_{i=1}^k} dp_i \wedge dq^i 
\quad.$$ It is
well-known in Differential Geometry that such a system of coordinates
always exists. In these $(p,q)$ -coordinates the Poisson brackets are
just the {\it classical canonical commutation relations} \index{relation! 
canonical commutation -} 
$$\begin{array}{rl} \{ p_i,q^j \} &= 
\delta_{i,j}, \cr
            \{ p_i,p_j \} &= 0, \cr
            \{ q^i,q^j \} &= 0.\end{array}$$ Thus the so called
(classical) {\it polarization means a choice of a maximal commutative
subalgebra $P(U) \subset C^{\infty}(U)$ in a fixed Darboux coordinate
neighborhood $U \subset M$.}

There is a physical principle, asserting that the {\it local
transformations } \index{local transformations} (i.e. {\it symmetry} ) 
\index{symmetry} do not interchange {\it positions
} \index{position} ( i.e. $q$-coordinates ) with {\it impulsions } (i.e. 
$p$-coordinates). \index{impulsion}
This means that {\it the stabilizer $G_x$ \index{stabilizer} at every 
point $x \in M$
normalizes the commutative subalgebra $P(U)$. } This means also that the
semi - direct product ${\frakt h}:= P(U) \rtimes {\frakt g}_x $ is some
subalgebra of ${\frakt g} = Lie G$, where ${\frakt g}_x:= Lie G_x$.

In quantum mechanics, one considers the following model: One fixes a
(separable) Hilbert space ${\cal H} \cong L^2({\Bbb R}^k)$, consisting of
the complex - valued functions with square-integrable module, and with
the usual scalar product of functions $$\langle f_1,f_2\rangle := \int_{{\Bbb R}^k}
f_1(q^1,\dots,q^k).\overline{f_2(q^1,\dots,q^k)}dq^1\dots dq^k \quad.$$

 In this model, the {\it quantum states } \index{quantum state} $\langle f|$ 
are represented by the
normed vectors $f \in L^2({\Bbb R}^k), \parallel f \parallel = 1$, where
by definition, $$\parallel f \parallel^2:= \langle f,f\rangle \quad.$$ Some time the
physicists note $\langle f|$ for this normed vector $f$.  The {\it quantum
quantities } \index{quantum quantities} are represented by the {\it normal, 
i.e. $[A,A^*] \equiv 0$ \index{operator!normal -}
( perhaps, unbounded, but having the operator closure ) operators }.
Roughly speaking, a {\it quantum mechanical system } 
\index{system!quantum mechanical -} is a pair $({\cal H}
,{\frakt A})$ consisting of a Hilbert space and a *-algebra of normal
operators. Their {\it measurable values} \index{values!measurable -} are 
just their eigenvalues.
Therefore the commutator $[A,B]$ of two operators $A$ and $B$ has its
eigenvalues as the difference of two measuring order of physical
quantities, $AB$ and $BA$. Therefore every physical quantity $A$ and its
{\it adjoint } quantity \index{quantity!adjoint} $A^*$ are always {\it 
exactly commensurable } \index{quantities!commensurable} (
i.e. can be together at any state without errors measurable.  {\it
Quantization } means a process associating $$p_i \mapsto \hat{p}_i,i =
1,\dots,k \quad,$$ $$q^i \mapsto \hat{q^i}, i = 1,\dots,k \quad,$$ and
more generally to every classical quantity, i.e. a smooth function, $f
\in C^{\infty}({\Bbb R}^k)$ a quantum quantity $\hat{f} \in {\frakt A}$,
$$f \mapsto \hat{f}, $$ in such a way that $${\sqrt{-1} \over
\hbar}[\hat{p}_i,\hat{q^j}] = \widehat{\{ p_i,q^j\}} = \delta_{ij} Id
\quad,$$ where $\hbar:= {h \over 2\pi}, h$ is the Planck constant.

This means that $f \mapsto {\sqrt{-1}\over\hbar}\hat{f}$, is a
homomorphism from the Lie algebra of classical quantities with Poisson
brackets to the Lie algebra with operator brackets, associated with the
associate *-algebra of quantum quantities. In

 particular, the classical commutation relations must be mapped into the
corresponding quantum commutatiom relations.

The well-known {\it Stone-von Neumann theorem } asserts that {\it there
exists a unique, up to unitary equivalence, solution of the quatum
correspondence $$p_i \mapsto \hat{p}_i, q^i \mapsto \widehat{q^i},$$
satisfying the precedent commutation relations:

 $\widehat{q^i} = q^i$, multiplication by coordinate function $q^i$

$\hat{p}_i = {\hbar \over \sqrt{-1}} {\partial \over \partial q^i}$,
derivation by coordinate variable $p_i$. }

These operators act on the dense subspace ${\cal S} ({\Bbb R}^k) \subset
L^2({\Bbb R}^k)$, consisting of the Schwartz class functions, and have
normal operator closure on the Sobolev spaces $H^s({\Bbb R}^k)$, which
are also Hilbert spaces.

{\it This unique solution of the quantum correspondence $$p_i \mapsto
\hat{p}_i = {\hbar \over \sqrt{-1}}{\partial \over \partial q^i} \quad,$$
$$q^i \mapsto \widehat{q^j} = q^j \quad,$$ satisfying the quantum
commutation relations $$[\hat{p}_i,\widehat{q^i}] = {\hbar \over
\sqrt{-1}}\delta_{ij}Id \quad,$$ can be also obtained from our geometric
quantization procedure } $$f \mapsto \hat{f}:= f +
{\hbar\over\sqrt{-1}}\nabla_{\xi_f} \quad.$$ This explains the physical
meaning of our procedure of geometric quantization.The central problem of
constructing the Hilbert space ${\cal H} = L^2({\Bbb R}^k)$ of quantum
states is therefore reduced to the problem of delating the
$(q,p)$-coordinates, which means delating a maximal commutative
$G_x$-invariant subalgebra {\frakt h} od {\frakt g}, subordinating the
functional $F \in {\frakt g}^*$
 i.e. ${\frakt g}_x = Lie G_x \subset {\frakt h} \subset {\frakt g}$,
{\frakt h} is a maximal subalgebra of {\frakt g} such that $\langle F,[{\frakt
h,h}]\rangle  \equiv 0 \quad $, what means the so called {\it 
polarization }, \index{polarization}
following the terminology of physicists. We are going now to some exact
mathematical models of these physical ideas.

\subsection{$(F,\tilde{\sigma})$-polarizations and polarizations }

Let us now consider a co-Adjoint orbit $\Omega = G_F \setminus G$, and
$H_0$ the connected closed subgroup corresponding to {\frakt h}. The last
condition of polarization means that ${i\over\hbar}\langle F,.\rangle $ is a
1-dimensional representation of

 ${\frakt h}$ and therefore can be considered as the diferential of the
character $$\chi_F|_{(G_F)_0 \cap H_0} = \exp{({i\over\hbar}\langle F,.\rangle )} \quad
.$$

Suppose that the orbit $\Omega$ is admissible, i.e. there exists a
representation denoted $\tilde{\sigma}\chi_F \in \hat{G}_F$ such that its
restriction to the connected component $(G_F)_0$ is a multiple of the
character $\chi_F$, symbolically $$\tilde{\sigma}\chi_F|_{(G_F)_0} = \mult
\chi_F \quad.$$ We have therefore a pair $(G_F,\tilde{\sigma}\chi_F)$
consisting of a group $G_F$ and a representation $\tilde{\sigma}\chi_F$.
From the Lie algebra point of view, this representation can be extended
to representation of Lie subalgebra ${\frakt h \subset g }$ Thus
polarization proposes some extension of the pair
$(G_F,\tilde{\sigma}\chi_F)$ to some pair $(H,\sigma)$ consisting of $H:=
G_F \ltimes H_0$ and some representation $\sigma$ such that $\sigma|_{G_F}
= \mult \tilde{\sigma}\chi_F $.

\begin{rem} We consider the representations $\sigma$ of, the restrictions
of which to $(G_F)_0 \cap H_0$ are the multiples of $\chi_F$. Following
the theory of projective representations ( see \S1 ) they are in a one -to
- one correspondence with the projective representations of $H/(G_F)_0
\cap H_0 \cong (G_F)_0 \setminus G_F$. Therefore we can consider them as
linear representations, passing to the two-fold coverings of type
$G_F^{{\Bbb Z}/2{\Bbb Z}} = G_F^{\frakt g}$ or $U(1)$-coverings of type
$G_F^{U(1)}$.
\end{rem}

\begin{defn} A {\it real $(F,\tilde{\sigma})$-polarization } 
\index{polarization!real -}$(H,\sigma)$
is a maximal with respect to inclusion pair of group and irreducible
representation, such that

(1) $G_F \subseteq H \subseteq G$, where $H$ is a closed subgroup in $G$.

(2) The restriction of $\sigma$ to $G_F$ is a multiple of
$\tilde{\sigma}\chi_F$ ; $$\sigma|_{G_F} = \mult \tilde{\sigma}\chi_F \quad
. $$

(3) $Lie H = {\frakt h}$ is $\Ad G_F$-invariant and $\sigma|_{H_0} = \mult
\sigma_0, \sigma_0 \in \hat{H}_0 $ is $G_F$-fixed.
\end{defn}

\begin{prop}  If $(H,\sigma)$ is a $(F,\tilde{\sigma})$-polarization,
then $$\codim_GH = {1 \over 2}\dim{\Omega_F} \quad. $$
\end{prop}
\begin{pf} 
Because $$(G_F)_0 \cap H_0 \triangleleft H_0, \sigma_0 \in
\hat{H}_0 $$ and $$\sigma|_{(G_F) \cap H_0} = \mult \chi_F \quad,$$ there
is a one - to - one correspondence between these $\sigma_0 \in \hat{H}_0$
such that $\sigma|_{H_0} = \mult \sigma_0$ and the projective
representations of $((G_F)_0 \cap H_0) \setminus H_0$ or the linear
representations of the coverings ( ${\Bbb Z}/2{\Bbb Z}$ -coverings or
$U(1)$-coverings). Because $$ \sigma|_{(G_F)_0 \cap H_0} = \mult \chi_F
\quad,$$ we have $$ \sigma_0|_{(G_F)_0 \cap H_0} = \mult \chi_F \quad,$$
$$\sigma|_{H_0} = \mult \chi_F \quad.$$ Thus $$ \langle F,[{\frakt h,h}]\rangle  \equiv
0 \quad.$$ As a linear space, ${\frakt h}$ is a maximal isotropic
subspace of the Kirillov's form $$\omega_F(.,.):= \langle F,[.,.]\rangle  \quad.$$
Thus $$\begin{array}{rl} \codim_{\frakt g}{ \frakt h} &= {1\over 2 }rank
\omega_F \cr
    &= {1\over 2}(dim{\frakt g} - dim \Ker \omega_F) \cr
    &= {1\over 2}(dim{\frakt g} - dim {\frakt g}_F) \cr
    &= {1\over 2}dim \Omega_F.\end{array}$$
\end{pf}

\begin{defn} A {\it polarization of orbit $\Omega_F$ at point $F$ } 
\index{polarization!of orbit} is a
maximal with respect to inclusion pair $({\frakt h},\sigma)$, consisting
of a Lie subalgebra {\frakt h} and a representation $\sigma$ such that:

(1) ${\frakt g}_F \subseteq {\frakt h} \subseteq {\frakt g}$. Denote the
(connected ) analytic subgroup corresponding to {\frakt h},

(2) subgroup $H:= G_F \ltimes H_0 $ is closed in $G$,

(3) the subalgebra {\frakt h} is $\Ad G_F$-invariant,

(4) $\sigma|_{(G_F)_0} = \mult\chi_F $,

(5) $\sigma|_{H_0} = \mult \sigma_0, \sigma_0 \in \hat{H_0}$ and
$\sigma_0$ is $G_F$-fixed.
\end{defn}

\begin{rem} It is easy to see that from a
$(F,\tilde{\sigma})$-polarization $(H,\sigma)$ one can construct a
polarization $({\frakt h}, \sigma)$ and vice-versa.
\end{rem}

In oder to obtain some irreducible representations by the construction,
one consider a restricted condition:

{\bf L. Pukanszky condition 4.6. } \index{Pukanszky condition} {\it Let us 
denote $p: {\frakt 
g}^* \rightarrow {\frakt h}^*$ the dual map to the inclusion ${\frakt h}
{\buildrel i \over \hookrightarrow} {\frakt g}$. The affine subspace
$p^{-1}(pF) = F + {\frakt h}^{\perp}$ is contained in the orbit $\Omega_F$. }

\begin{defn}
If the Pukanszky is satisfied, we say that $F$ is {\it well-admissible }
\index{functional!well-polarizable -} and {\frakt h} is an {\it admissible 
polarization } \index{polarization!admissible -} at $F$.
\end{defn}

\begin{rem} If ${\frakt h}$ is an admissible polarization at $F$ then by
translations, $\Ad_g({\frakt h})$ is an admissible polarization at $K(g)F
\in \Omega_F$. Therefore we say that $\Omega$ is an {\it admissible 
orbit}. \index{orbit!admissible -} \end{rem}

Some more general notion of $(F,\tilde{\sigma})$-polarization and
polarization can be obtained in passing to the complex domain.

\subsection{Complex polarizations }

\begin{defn}{\it A complex $(F,\tilde{\sigma})$-polarization } 
\index{polarization!complex -} $({\frakt
p},H,\sigma)$ is a maximal quadruple satisfying the folowing conditions:

(a) {\frakt p} is a complex subalgebra of the complexified ${\frakt
g}_{\Bbb C} $, such that $({\frakt g}_F)_{\Bbb C} \subset {\frakt p}
\subset {\frakt g}_{\Bbb C}$.

(b) {\frakt p} is $\Ad G_F$-invariant.

(c) ${\frakt h}:= ({\frakt p} \cap \overline{\frakt p}) \cap {\frakt g} $
and ${\frakt m}:= ({\frakt p} + \overline{\frakt p}) \cap {\frakt g}$ are
real subalgebras.

(d) There exists the closed subgroups $H$ and $M$ such that $H \subset M
\subset G$, where $Lie H = {\frakt h}, Lie M = {\frakt m} $.

(e) $\sigma$ is an {\it irreducible } representation, $\langle \sigma\rangle  \in
\hat{H}$ such that $\sigma|_{G_F} \simeq \mult \tilde{\sigma}$.

(f) $\rho$ is a complex representation of ${\frakt p}$ such that
$\rho|_{\frakt h} = D\sigma$, the differential of $\sigma$.
\end{defn}

\begin{defn} A {\it complex polarization } \index{polarization!complex -} of 
$\Omega$ at $F$ is a triple
$({\frakt p}, \rho, \sigma_0)$, satisfying the following conditions:

(a) ${\frakt p}$ is a complex Lie subalgebra of ${\frakt g}_{\Bbb C}$,
such that $({\frakt g}_F)_{\Bbb C} \subseteq {\frakt p} \subseteq {\frakt
g}_{\Bbb C}$,

(b) ${\frakt p}$ is $\Ad G_F$-invariant,

(c) there exist the real Lie subalgebras {\frakt m} and {\frakt h} of
{\frakt g} such that $${\frakt p} + \overline{\frakt p} = {\frakt m}_{\Bbb
C}\quad,\quad {\frakt m}:= ({\frakt p} + \overline{\frakt p}) \cap
{\frakt g} \quad,$$ $${\frakt p} \cap \overline{\frakt p} = {\frakt
h}_{\Bbb C} \quad ; \quad {\frakt h}:= ({\frakt p} \cap \overline{\frakt
p}) \cap {\frakt g} = {\frakt p \cap g} \quad,$$

(d) the subgroups $M_0, H_0, M,H$ are closed, where $M_0$ and $H_0$ are
the analytic subgroups corresponding to the Lie algebras ${\frakt m}$ and
${\frakt h}$, respectively, and $M:= G_F \ltimes M_0, H:= G_F \ltimes
H_0$,

(e) $\sigma_0 \in \hat{H}_0 $ and the point $\langle \sigma_0\rangle $ in the dual
object $\hat{H}_0$ is fixed under the action of $G_F$ and $\sigma_0|_{G_F
\cap H_0} = \mult \tilde{\sigma}\chi_F|_{G_F \cap H_0} $,

(f) $\rho$ is a complex representation of ${\frakt p}$ such that
$\rho|_{\frakt h} \simeq D\sigma$, the differential of $\sigma$.
\end{defn}

\begin{thm}  There is a one - to - one correspondence between the
complex $(F,\tilde{\sigma})$-polarizations $({\frakt p},\rho,H,\sigma)$
and the complex polarizations $({\frakt p},\rho,\sigma_0)$.
\end{thm}
\begin{pf} 
It is clear that from a $(F,\tilde{\sigma})$-polarization
$({\frakt p},\rho,H,\sigma)$, one can construct easily a polarization
$({\frakt p},\rho,\sigma_0)$. Conversely, considering the surjection
$$G_F \times H_0 \longrightarrow H:= G_F \ltimes H_0$$
 with kernel $G_F \cap H_0$. From our assumptions, we have
$$\tilde{\sigma}\chi_F|_{G_F \cap H_0} \simeq \sigma_0|_{G_F \cap H_0 }
\quad.$$ Because $\sigma_0$ is $G_F$-fixed we can construct a
representation $\sigma$ of $H = G_F \ltimes H_0$, such th at
$$\sigma|_{G_F} \simeq \tilde{\sigma}\chi_F \quad, $$ $$\sigma|_{H_0}
\simeq \sigma_0 \quad.$$ It is easy now to see from the definition that
$({\frakt p},H,\rho,\sigma)$ is a $(F,\tilde{\sigma})$-polarization. The
theorem is proved.
\end{pf}

\subsection{Weak Lagrangian distributions }

\begin{rem}
We say that an integrable $G_F$-invariant tangent distribution $L$ on
 $\Omega \approx G_F \setminus G$ is {\it weak Lagrangian } 
\index{distribution!weak Lagrangian -} iff the
semidirect product of Lie algebras ${\frakt p}:= ({\frakt g}_F)_{\Bbb
C}\times \Gamma_G(L) $, where by $\Gamma_G(L)$ we denote the space of all
the $G$-invariant sections of the tangent distribution $L$, is a complex
Lie algebra and the representation $\tilde{\sigma}\chi_F$ of $G_F$ can be
extended to a representation of $H = G_F \ltimes H_0$ and to a complex
representation $\rho$ of ${\frakt p}$, satisfying the conditions in the
definition 4.10 above. 
\end{rem}

\begin{defn}
We say that this distribution is {\it closed } \index{distribution!closed 
-} iff the groups $M_0,H_0,M,H$ are closed, where by definition, $M_0$ and
$H_0$ are the analytic subgroups corresponding to the Lie algebras
$${\frakt m}:= ( {\frakt p} + \overline{\frakt p}) \cap {\frakt g} $$ and
$${\frakt h}:= ( {\frakt p} \cap \overline{\frakt p}) \cap {\frakt g}
\cong {\frakt p \cap g} \quad,$$ respectively, and $$M:= G_F \ltimes
M_0, H:= G_F \ltimes H_0 \quad.$$ Thus {\it there is a bijection
between the $(F,\tilde{\sigma})$-polarizations and the closed weak
Lagrangian $G$-invariant integrable tangent distributions. } 
\end{defn}

\subsection{Duflo data }

In general, the quotient group $(G_F \cap H_0) \setminus H_0$ is
noncommutative. We considered in our definition of polarizations
$\sigma_0 \in \hat{H}_0$, which is fixed under the adjoint action of
$G_F$ and the condition $$\sigma_0|_{H_0 \cap G_F} = \mult
\tilde{\sigma}\chi_F|_{H_0 \cap G_F} \quad.$$ The representations of this
type are in a one - to -one correspondence with the projective
representations of group $(H_0 \cap G_F) \setminus H_0$. In the
particular case, where we are taking not multi- but one-dimemsional
representations $\sigma_0$, we obtain the single corresponding projective
representation. As a projective represntation of Lie group, it can be
lifted to a linear representation in considering the coverings. Following
M. Duflo, we consider now the two - fold covering $H_0^{{\Bbb Z}/2{\Bbb
Z}}= H_0^{\frakt g}$ of $H_0$ $$1 \rightarrow {\Bbb Z}/2{\Bbb Z}
\rightarrow H_0^{\frakt g} \rightarrow H_0 \rightarrow 1 \quad.$$ We
denote the generator of the cyclic group ${\Bbb Z}/2{\Bbb Z}$ by
$\varepsilon$ and define the lifted character $\chi_F^{\frakt g}$ of
$H_0^{\frakt g}$ by the condition $$\chi_F^{\frakt g}(\varepsilon) = -1
\quad.$$ Then we consider the so called {\it odd representations } 
\index{representation!odd -}
$\sigma_0 \in \widehat{H_0^{\frakt g}}$ such that $$\sigma_0|_{G_F^{\frakt
g} \cap H_0^{\frakt g}} = \mult \tilde{\sigma}\chi_F^{\frakt g} \quad.$$
We have therefore $$\langle F,[{\frakt h,h}]\rangle  \equiv 0 \quad.$$ Thus passing to
the two-fold coverings $G_F^{\frakt g}, (G_F)_0^{\frakt g},H_0^{\frakt g}
, H^{\frakt g}$ of $G_F,(G_F)_0,H_0,H $, respectively, we can in this
particular situation consider the real polarization $({\frakt h},\sigma)$
simply as a maxim al possible with respect to inclusion Lie subalgebra
${\frakt h}$, satisfying the following conditions:

(a) ${\frakt g}_F \subseteq {\frakt h} \subseteq {\frakt g}, {\frakt h} $
is a Lie subalgebra. Let $H_0$ be the corresponding analytic subgroup.

(b) ${\frakt h} $ is $\Ad G_F$-invariant.
 
(c) $H_0$ and $H:= G_F \ltimes H_0$ are  closed subgroups in $G$.

(d) $\langle F,[{\frakt h,h}]\rangle  \equiv 0 $.

One can easily generalize this notion of polarization, passing into the
complex domain.

\begin{defn} A {\it ( complex ) polarization } 
\index{polarization!complex -} is a ( complex
) maximal with respect to inclusion Lie subalgebra {\frakt p} of ${\frakt
g}_{\Bbb C}$ satisfying the following conditions:
  
(a) $({\frakt g}_F)_{\Bbb C} \subseteq {\frakt p} \subseteq {\frakt
g}_{\Bbb C}$.

(b) {\frakt p} is $\Ad G_F$-invariant.

(c) ${\frakt h}:= ({\frakt p} \cap \overline{\frakt p}) \cap {\frakt g} =
{\frakt p \cap g}$ and $ m:= ({\frakt p} + \overline{\frakt p}) \cap
{\frakt g} $ are real Lie subalgebras of ${\frakt g}$.

(d) $H_0,M_0,H,M$ are closed subgroups in $G$, where $H_0$ and $M_0$ are
the analytic subgroups, corresponding to the real Lie subalgebras
${\frakt h,m}$ said above, $H:= G_F \ltimes H_0, M:= G_F \ltimes M_0$
\end{defn}
 
Then one considers a fixed representation $\tilde{\sigma}$ such that its
restriction to $(G_F)_0$ is a multiple of the character $\chi_F$ and lift
it to the odd character $\chi_F^{{\frakt g}_F}$ of the two - fold covering
$G_F^{{\frakt g}_F}$.

\begin{defn}
The pair $(F,\tilde{\sigma})$ is called a {\it Duflo's datum } 
\index{Duflo's data} 
and the complex Lie subalgebra ${\frakt p}$ is called a {\it polariszation }
in the Duflo's theory.
Now consider any representation $\sigma$, satisfying the conditions:

(e) $\langle \sigma\rangle  \in \widehat{H^{\frakt g}}, \sigma|_{(G_F)_0^{\frakt g}} =
\mult \tilde{\sigma}\chi_F^{\frakt g}$.

(f) there is a complex contination of the representation $D\sigma$ of the
real Lie algebra ${\frakt h} $ to a complex representation $\rho$ of the
complex Lie algebra ${\frakt p}$ such that the E. Nelson's conditions are
satisfied.
\end{defn}

We return then to our theory, exposed above.

\begin{rem} This particular case, $\dim{ \sigma} = 1 $
correspondence to the case of {\it Lagrangian } 
\index{distribution!Lagrangian} $G_F$-invariant integrable tangent 
distributions. \end{rem}

We finish this section by recalling a theorem of M. Vergne on existence of
polarizations. The reader could find more detail in the original works of
M. Vergne.

(1) Let $(V,\omega)$ be a symplectic vector space and suppose that $$0 =
V_0 \subset V_1 \subset V_2 \subset \dots \subset V_n = V \eqno({\frakt
S}) $$ is a chain of vector subspaces such that:$$\dim V_k = k, k =
0,1,\dots,n \quad.$$Define $$\omega_k:= \omega|_{V_k} $$ and $$
W({\frakt S},\omega):= {\sum_{k=1}^n} \Ker \omega_k \quad.$$ Then
$W({\frakt S},\omega)$ is a maximal isotropic subspace for $\omega$.

(2) If $V$ is a Lie algebra and $({\frakt S})$ is a chain of ideals $V_k
\triangleleft V$ and $\omega = \omega_F$ is the Kirillov's form for some
$F \in V^*$, then $W({\frakt S},\omega)$ is a maximal isotropic
subalgebra.

\begin{thm}[\bf Real version of M. Vergne's Theorem ] 
Let $G$ be a real exponential ( i.e. the exponential map is a 
diffeomorphism ) (
therfore, solvable ) Lie group, with Lie algebra ${\frakt g}:= Lie G,
F \in {\frakt g}^*$.  Then there exists a polarization ${\frakt h
\subseteq g} $ such that

(1) ${\frakt h} = W({\frakt S},\omega_F)$ for some chain of ideals
$({\frakt S})$,

(2) ${\frakt h}$ is $\Ad G_F$-invariant,

(3) $\langle F,[{\frakt h,h}]\rangle  \equiv 0 $, i.e. ${\frakt h}$ is a maximal
isotropic subalgebra of $\omega_F$,

(4) ${\frakt h}$ satisfies the Pukanszky condition: $$F + {\frakt
h}^{\perp} \subseteq \Omega_F \quad. $$ 
 
\end{thm}

\section{Bibliographical Remarks}
The first idea about a general construction of multi-dimensional
quantization procedure was appeared in 1979-1980 in \cite{diep5} and
\cite{diep12}. Later it was largely developed in many works
\cite{diep4}-\cite{diep6}, \cite{diep11}, \cite{duflo1}-\cite{duflo3}.

\chapter {Partially invariant holomorphly induced representations}
\section{Holomorphly Induced Representations. Lie Derivative}

From a $(F,\tilde{\sigma})$-polarization $({\frakt p},H,\rho,\sigma)$ we
can produce some representation of $G$ on the space of partially invariant
partially holomorphic sections of some induced bundle. After that we
shall show that the Lie derivative of
 this representation is just the representation of Lie algebra ${\frakt
g}$, arising from the procedure of geometric ( multidimensional )
quantization.

\subsection{Partially invariant holomorphly induced 
representations  }  

\begin{defn} Let $D$ be a closed ( not necessarily connected ) subgroup of
$G$, with Lie algebra ${\frakt d} := \Lie D$. Let $\tilde{\sigma}$ be a
fixed unitary ( not necessarily irreducible ) of $D$ in a Hilbert space
$\tilde{V}$. Then a $(D,\tilde{\sigma})$-{\it polarization } 
\index{polarization} $({\frakt
p},H;\rho,\sigma)$ will be any {\it maximal with respect to inclusion
quadruple }, satisfying the following conditions:

(a) ${\frakt p}$ is a complex subalgebra such that, ${\frakt d}_{\Bbb C}
\subset {\frakt p} \subset {\frakt g}_{\Bbb C}$.

(b) ${\frakt p}$ is $\Ad D$-invariant.

(c) $H$ is a closed subgroup of $G$ such that its Lie algebra is just $\Lie
H = {\frakt h} = {\frakt p \cap g }$.

(d) There exists a closed subgroup $M$ of $G$ with Lie algebra $$\Lie M :=
{\frakt m} = ({ \frakt p} + \overline{\frakt p}) \cap {\frakt g} \quad.
$$

(e) $\sigma$ is an irreducible representation of $H$, $\langle \sigma\rangle  \in
\hat{H}$, such that $\sigma|_D = \mult \tilde{\sigma} $.

(f) $\rho$ is a complex representation of ${\frakt p}$ such that all the
conditions of E. Nelson are satisfied and that the corresponding to
$\sigma$ Lie algebra representation $D\sigma, \rho|_{\frakt h} = D\sigma$.
\end{defn}

\begin{exam} Let us consider $D = G_F$, $\tilde{\sigma}\chi_F \in
\\Rep(G_F)$ and a $(F,\sigma)$-polarization $({\frakt p},\rho,H,\sigma)$ of
$\Omega$. it is easy to check that we have $({\frakt p},\rho,H,\sigma)$
as a $(D,\tilde{\sigma}\chi_F)$-polarization.
\end{exam}

Consider the principal bundle $$\begin{array}{rl} D \rightarrow &G \cr
   &\downarrow \cr
   &X = D \setminus G \quad.\end{array}$$ It is a locally trivial
principal bundle. Fix a {\it connection } on it ( i.e. a {\it
trivialization } ). Let us consider a $(D,\tilde{\sigma})$-polarization
$({\frakt p},H,\rho,\sigma)$ and a fixed connection $\Gamma$ on the
principal bundle $$\begin{array}{rl} H \rightarrow &G \cr
  &\downarrow \cr
  &H \setminus G.\end{array}$$ Consider the associated with respect to
the projection map $D \setminus G {\buildrel p \over \twoheadrightarrow} H
\setminus G$ principal bundle $$\begin{array}{rl} D \rightarrow &G \cr
  &\downarrow  \cr
  &D \setminus G.\end{array}$$ There is the so called associated
connection $\overline{\Gamma}$ on it such that the following diagram is
commutative $$\vbox{\halign{ #&#&#&#\hfill \cr $D \rightarrow$& $G $ &$H
\rightarrow$&$G$ \cr
               &${\downarrow}\overline{\Gamma}$ &  &${\downarrow}\Gamma$ \cr
 &$D \setminus G \qquad\quad {\buildrel p \over \twoheadrightarrow}$ &
&$H\setminus G$.\cr}}$$ We have a representation $\sigma$ of $H$ and
$\sigma|_D = \mult\tilde{ \sigma}$ of $D$. Therefore we have also the
corresponding affine connections $\nabla$ and $\overline{\nabla}$ on the
vector bundles ${\cal E} _V(H \setminus G)$ and $p^*{\cal E}_V(D \setminus
G)$, respectively, such that the following diagram is commutative
$$\vbox{\halign{ #&#&#&#&#&#\hfill \cr $V \rightarrow$ & &${\cal E} _V $
&$\tilde{V} \rightarrow$ & &$p^*{\cal E} _V$ \cr
  & &${\downarrow}\nabla$ & & &${\downarrow}\overline{\nabla}$ \cr
  &$H$&$\setminus G \qquad \quad {\buildrel p \over \leftarrow}$ &
&$D$&$\setminus G$. \cr}}$$ Therefore the corresponding spaces of sections
could be included as follows $$\Gamma({\cal E} _V(H\setminus G),\nabla)
\hookrightarrow \Gamma({\cal E} _V(D\setminus G),\overline{\nabla}) \quad
.$$ The image of this inclusion is just the space $$\{ s \in \Gamma({\cal
E} _V,\overline{\nabla}); s \enskip is \enskip left \enskip H-invariant
\}.$$ We consider also the next projection $$G \twoheadrightarrow
D\setminus G$$ and we have by analogy the diagram of triple vector bundles
with connection $$\vbox{\halign{ #&#&#&#&#&#\hfill \cr $p_1^*p_2^*$&${\cal
E} _V$ &$p_2^*$&${\cal E} _V$ & &${\cal E} _V$ \cr
  &${\downarrow}\overline{\overline{\nabla}}$ &  &${\downarrow}\overline{\nabla}$ &  &${\downarrow}\nabla$ \cr
  &$G {\buildrel p_1 \over \longrightarrow }$&$D$&$\setminus G {\buildrel
p_2 \over \longrightarrow}$&$H$&$\setminus G \quad.$ \cr}}$$ Remark that
the last bundle $$V \rightarrowtail p_1^*p_2^*{\cal E} _V
\twoheadrightarrow G$$ is trivial, and we can identify the sections of
this bundle with $\tilde{V}$-valued functions on $G$. It is easy to see
that $$\Gamma({\cal E} _V) \cong \{ f \in C(G,V) ; f(hx) = \sigma(h)f(x)
,\forall h \in H, \forall x \in G \},$$ $$\Gamma(p_2^*{\cal E} _V) \cong
\{ f \in C(G,V) ; f(hx) = \sigma|_D(h)f(x),\forall h \in D, \forall x
\in G \} \quad.$$

\begin{defn} Section $s \in \Gamma({\cal E} _V)$ is 
called 
{\it partially invariant and partially holomorphic } 
\index{section!partially invariant partially holomorphic -} iff its 
covariant
derivatives $\nabla_{\xi}, \xi \in \overline{\frakt p}$, along the
directions of $\overline{\frakt p}$ vanish.
\end{defn}

\begin{thm}  The space $\Gamma({\cal E} _{V,\rho,\sigma}(H \setminus
G))$ of partially invariant and partially holomorphic sections of the
induced bundle ${\cal E} _V(H \setminus G)$ is invariant subspace in
$\Gamma({\cal E} _V) $ and is isomorphic to the space
$C^{\infty}(G;{\frakt p},H,\rho,\sigma)$ of the $\tilde{V}$-valued
functions $f$ on $G$, satisfying the following equations $$f(hx) =
\sigma(h)f(x), \forall x \in G, \forall h \in H \quad, $$ $$[L_X +
\rho(\xi_X)]f
 = 0, \forall X \in \overline{\frakt p} \quad,$$ where, by definition,
$L_X := L_{\xi_X}$ is the Lie derivation along the vector field $\xi_X$,
$X \in {\frakt g}$.
\end{thm} 
\begin{pf} 
Recall the action $L_X$
:$$(L_Xf)(x) := {d \over dt}f(\exp{(-tX)} x)|_{t = 0} $$ and $$f(hx) =
\sigma(h)f(x)$$ are the left action. Therefore $C^{\infty}(G;{\frakt
p},H,\rho,\sigma)$ is a $G$-invariant subspace with respect to the right
translations.

By analogy, we have an inclusion $$\Gamma({\cal E} _{V,\rho,\sigma})
\hookrightarrow \Gamma({\cal E}_V) \quad $$ and it is a $G$-invariant
subspace. Looking at the diagram $$\vbox{\halign{ #&#&#&#\hfill \cr $V
\rightarrow G \times V \cong p^*$&${\cal E} _V $ &$V \rightarrow$&${\cal
E}_V$ \cr
  &$\downarrow$   &   &${\downarrow}\nabla$ \cr
  &$G \qquad\quad {\buildrel p \over \twoheadrightarrow}$ & &$H\setminus G
\quad,$ \cr}}$$ we have the natural isomorphism $$\Gamma({\cal E}
_{V,\rho,\sigma}(H\setminus G)) \cong C^{\infty}(G;{\frakt
p},H,\rho,\sigma) \quad.$$ The proof of the theorem is here achieved.
\end{pf}

\begin{rem} It is easy to see that the partially invariant and partially
holomorphic sections of the induced bundle ${\cal E} _V(H \setminus G)$
form a structural sheaf of the {\it partially invariant holomorphly
induced bundle }, \index{bundle!partially invariant holomorphly 
induced -} noted ${\cal E} _{V,\rho,\sigma}$ \end{rem}

\begin{defn} The natural action of $G$ on $\Gamma({\cal E}
_{V,\rho,\sigma})$, isomorphic to the right translations on
$C^{\infty}(G;{\frakt p},H,\rho,\sigma)$, is called the {\it partially
invariant and holomorphly induced representation } 
\index{representation!partially invariant holomorphly induced -} of 
$G$ and denoted by $\Ind(G;{\frakt p},H,\rho,\sigma)$.
\end{defn}

\begin{rem}
Our partially invariant and partially holomorphic sections are just the image of the partially holomorphic in the sense of R. B. Blattner sections with respect to the inclusion $$\Gamma({\cal E}
_{V,\rho,\sigma}(H \setminus G)) \hookrightarrow \Gamma({\cal E}
_{V,\rho|_{\frakt d},\sigma|_D}(D\setminus G))\quad.$$
\end{rem}

\subsection{Unitarization }

Let us denote by $\Delta_H$ and $\Delta_G$ the {\it modular functions } 
\index{function!modular -} (
for the Haar measures ) of $H$ and $G$, respectively. Then $$\delta =
\sqrt{\Delta_H / \Delta_G}$$ is a non-unitary character of $H$. We have
the so called {\it ${1 \over 2}$-density bundle } \index{bundle!${1\over 
2}$-density -} 
$${\cal M} ^{1 \over 2} := G \times_{D,\delta} {\Bbb C}$$ on $D\setminus 
G$. The bundle
$$ \tilde{\cal E} _{V,\rho,\sigma}(H\setminus G) := {\cal
E}_{V,\rho,\sigma}(H\setminus G) \otimes {\cal M}^{1 \over 2}$$ is a
$G$-bundle oner $X = D\setminus G$. Because $\sigma$ is unitary, for
every section $s \in \Gamma({\cal E} _{V,\rho,\sigma})$, $$\parallel s
\parallel^2_V \in \Gamma({\cal M}) \quad,$$ where by definition ${\cal M}
:= {\cal M} ^{1\over 2} \otimes {\cal M} ^{1\over 2}$ is the so called
{\it density bundle. } \index{bundle!density -} We define the scalar 
product on the vector space
$\Gamma(\tilde{\cal E}_{V,\rho,\sigma})$ by the formula $$(s_1,s_2) :=
\int_{H\setminus G} (s_1(x),s_2(x))_V dx \quad.$$ Let us denote the
completion of $\Gamma(\tilde{\cal E} _{V,\rho,\sigma}(H\setminus G))$ by
$L^2(\tilde{\cal E} _{V,\rho,\sigma})$. it is a Hilbert space, which can
be identified with the Hilbert space $L^2(G;{\frakt p},H,\rho,\sigma)$
,consisting of $V$-valued functions $f : G \rightarrow V $, satisfying
the equations $$f(hx) = \sigma(h)\sqrt{\Delta_H(h)/\Delta_G(h)} f(x),
\forall x \in G,\forall h \in H \quad,$$ $$(L_X + \rho(X) + {1\over
2}\tr _{{\frakt g}_{\Bbb C}/{\frakt h}} \ad X)f = 0,\forall X \in
\overline{\frakt p} \quad,$$ with the scalar product norm
$$\int_{H\setminus G} \parallel f(x)\parallel^2_V dx < \infty \quad.$$
From now on we mean $\Ind(G;{\frakt p},H,\rho,\sigma)$ this unitarized
representation of $G$ on $L^2(G;{\frakt p},H,\rho,\sigma)$ $ \cong
L^2(\tilde{\cal E}_{V,\rho,\sigma})$.

\subsection{Lie derivation }
 
Remember that $\Lie_{\xi}$ is a functor which can be applied to functions,
sections, differential forms, vector fields,... and also to
representations, and gives some infinitesimal ones. We apply here this
functor $\Lie_{\xi}$ to our induced representation to have an
infinitesimal representation of our Lie algebra ${\frakt g} = \Lie G$.

\begin{thm}
For every $X$ in ${\frakt g}$,
$$\Lie_X \Ind(G;{\frakt p},H,\rho,\sigma) = {i\over\hbar}\hat{f}_X \quad.$$
\end{thm}

\begin{rem} The quantization procedure gives us some Lie algebra
representation of type $$\Lambda : X \mapsto {i\over\hbar}\hat{f}_X \quad
.$$ This is ``local version '' of the theory. The construction of
partially invariant holomorphly induced representations gives us a
``global theory'' of group representations. The theorem asserts the
relationship between the local and global constructions.
\end{rem}

{\sc Proof of the Theorem. }

First of all, recall the construction of the Hilbert bundle
$\tilde{\cal E} _{V,\rho,\sigma}$. By definition, we have $$\tilde{\cal
E} _{V,\rho,\sigma} \cong G \times_{G_F} (\delta\sigma,V) = (G \times
V)/\sim \quad,$$ where by definition, $$(g,v) \sim (g',v')
\Longleftrightarrow g' = hg, v' = \delta(h)\sigma(h)v \quad,$$ for some
$h \in G_F$. The representation $\delta\sigma : G_F \rightarrow {\Bbb
U}(V) \times {\Bbb C} $ is unitary projective, $\delta\sigma(G_F) \subset
{\Bbb U}(V) \times {\Bbb C} $. The image of this representation is a
finite dimensional topological subgroup, and then by the positive
solution of the well-known Hilbert problem is a Lie subgroup, of
${\Bbb U}(V) \times {\Bbb C} $. Therefore, in virtue of the Stone
theorem, every one-parameter subgroup of ${\Bbb U}(V) \times {\Bbb C}$
has a skew-selfadjoint generator, say $A = -A^*$ and the 1-parameter
subgroup self is therefore of form $$u_t = \exp{({i\over\hbar}tA)}, t \in
{\Bbb R} \quad,$$ we can thus define the derivative of 1-parameter
subgroups.

Recall now the definition of Lie derivative of bundle homomorphisms,
following Y. Koschmann-Schwarzbach. Let us consider two $G$-bundles
${\cal E,E'}$ over $D\setminus G$ of type ${\cal E} _{V,\rho,\sigma}$.
Then $G$ acts also on the homomorphism bundle $$\End  V \rightarrowtail
\Hom ({\cal E,E'} )\twoheadrightarrow D\setminus G = M \quad.$$ Consider
$(u,u_M) \in \Hom ({\cal E,E'})$, $$\vbox{\halign{ #\quad&#\hfill \cr
${\cal E} \quad {\buildrel u \over \longrightarrow} $ &${\cal E} '$ \cr
$\downarrow $ &$\downarrow$ \cr $M \quad {\buildrel u_M \over
\longrightarrow} $ &$M$\quad. \cr}}$$ Then $G$-action on $\Hom ({\cal E,E'}
)$ is just $$g.(u,u_M) := (g.u.g^{-1},g.u_M.g^{-1}) \quad.$$Let $X\in
{\frakt g}$ and $g_t = \exp(tX) \subset G$ be the corresponding
1-parameter subgroup. We can define $$X.u := {d\over dt}(g_t.u)|_{t=0}
\quad, $$ which exists following the Stone theorem.

\begin{rem}
 (1) $X.u$ is a differential operator of degree 1, $$\vbox{\halign{
#&#&#&#&#&#\hfill\cr ${\cal E} $& &$\rightarrow$& &$u_M^*{\cal E} '$
&${\cal E} '$ \cr
  &$\searrow$&  &$\swarrow$&   &$\downarrow$ \cr
  &          &$M$& &${\buildrel u_M \over \longrightarrow}$ &$M\quad,$
\cr}}$$ $$X.u = X_{\varepsilon}.u-u.X_{\varepsilon'} \quad,$$ where
$X_{\varepsilon}$ is a differential operator in $\Gamma({\cal E} )$ and
$X_{\varepsilon'}$ is a differential operator in $\Gamma({\cal E}
')$.

(2) Consider the particular case, where $u_M = \Id  $. Then $$X.u \in
\Hom ({\cal E,E'} ) = {\cal E}^* \otimes {\cal E} ' \quad,$$is the Lie
derivative of $G$-action on ${\cal E}
 \otimes {\cal E}
'$.

(3) Consider the case where ${\cal E} = M$. Then $u$ is some section and
then $$X.u = \nabla_Xu$$ is just the covariant derivative of $u$ along the
vector field $\xi_X$.
\end{rem}

Return now to our situation of quantization. We have a
$(F,\tilde{\sigma})$-polarization $({\frakt p},H,\rho,\sigma)$ and we can
consider the inclusion $$\Gamma(\tilde{\cal E} _{V,\rho + {1\over 2}\tr
(\ad ),\delta\sigma}) \hookrightarrow \Gamma(\tilde{\cal E} _{V,\rho +
{1\over 2}\tr  (\ad ),\delta\sigma}) \quad,$$ $$p : G_F \setminus G
\twoheadrightarrow H \setminus G \quad.$$ Let us consider $$
\overline{\nabla} = p^*\nabla \quad.$$ We have $$\nabla_{\xi_X} =
L_{\xi_X} + {i\over \hbar}\alpha_1(\xi_X) \quad,$$ where
$$\alpha_1(\xi_X) := \rho(X) + {1\over 2}\tr _{{\frakt g}_{\Bbb
C}/{\frakt h}}\ad X$$ is the differential of the
 representation $\delta\sigma|_{G_F}$.

Consider the differential form $\beta$ on $M = G_F \setminus G$, defined
by the formula $$\langle \beta, \xi\rangle (y) := \langle y,\xi(F)\rangle , \forall y \in M \quad
.$$ The expression $\langle \beta,\xi_X\rangle $ can be considered also as the
generating function of vector field $\xi_X$, $$\langle \beta,\xi_X\rangle (y) =
\langle y,\xi_X(F)\rangle  = \langle y,X\rangle  = f_X(y) \quad.$$ Thus we have $$\begin{array}{rl}
\nabla_{\xi_X} &= L_{\xi_X} + {i\over \hbar}\alpha_1(\xi_X) \cr &=
L_{\xi_X} + {i\over\hbar}f_X + {i\over\hbar}(\alpha_1(\xi_X) -
\langle \beta,\xi_X\rangle ) \cr &= {i\over\hbar}\hat{f}_X.\end{array}$$ The theorem
is proved.

\begin{exam} The degenerate principal series 
representations of
semi-simple Lie groups can be obtained from our construction, as
$$\Ind_P^G((\tau \times \sigma)\otimes\nu), \sigma \in \hat{M}_{disc}, P
= MAU \quad,$$ where $G = KP$ is the well-known Gauss decomposition.
\end{exam}

\section{The Irreducible Representations of Nilpotent Lie Groups}

In this section we shall show the Duflo's construction of irreducible
representations as a particular case of our construction of partially
invariant and holomorphly induced representations. This Duflo's
construction proposes to pass to the two-fold covering of stabilizers.
It is natural to consider the so called Shale-Weil construction. In the
particular case of nilpotent Lie groups, where the stabilizers are
connected, the results give us the original Kirillov orbit method.

\subsection{Duflo construction}
 
Recall that our multidimensional quantization procedure starts with a
$(F,\tilde{\sigma})$-polarization $({\frakt p},H,\rho,\sigma_0)$,
satisfying the conditions (a)-(f) of the definition 4.10. This is
equivalent to giving a $(G_F,\tilde{\sigma}\chi_F)$ -polarization
$({\frakt p},H,\rho,\sigma)$ of the definition 5.1. In general case,
$\sigma_0$ is an arbitrary irreducible unitary representation of $H_0$
such that its restriction to $G_F \cap H_0$ is equivalent to the
restriction of $\tilde{\sigma}\chi_F$ to $G_F \cap H_0$, symbolically
$$\sigma_0|_{G_F \cap H_0} \simeq \tilde{\sigma}\chi_F |_{G_F \cap
H_0}\quad.$$ For example the degenerate principal series representations
of semi-simple Lie groups can be obtained in this way.

For the particular case, where $\sigma_0 = \chi_F $, one can
consider the two-fold coverings of type $$1 \rightarrow {\Bbb
Z}/2{\Bbb Z} \cong \{ 1,\varepsilon \} \rightarrow G_F^{\frakt g}
\rightarrow G_F \rightarrow 1 \quad,$$ and $$1 \rightarrow {\Bbb
Z}/2{\Bbb Z} \rightarrow H_0^{\frakt g} \rightarrow H_0
\rightarrow 1 \quad,$$ which is obtained by restriction from $$1
\rightarrow {\Bbb Z}/2{\Bbb Z}\rightarrow H^{\frakt g}\rightarrow
H := G_F \ltimes H_0 \rightarrow 1 \quad,$$ and some unitary
representation $\hat{\sigma}\in \widehat{H_0^{\frakt g}}$, such
that $$\hat{\sigma}_0|_{G_F^{\frakt g} \cap H_0^{\frakt g}}
\simeq \tilde{\sigma}\chi_F^{\frakt g}|_{G_F^{\frakt g} \cap
H_0^{{\frakt g}_F}} \quad, $$ where the {\it odd character } 
\index{character!odd -}
$\chi_F^{{\frakt g}_F}$ is the lifted character of the two-fold
covering $H^{\frakt g}$ from the character $\chi_F$ of $H$, such
that $$\chi_F^{\frakt g}(\varepsilon) = -1 \quad.$$

As it was remarked, considering only the particular case with the
lifted $\hat{\sigma}$ in place of un-lifted $\sigma$, we can restrict our
consideration to the Duflo's data, and polarization. Interesting to
remark that the resulting representation is in fact independent from a
particular polarization. It depends therefore only on the Duflo's data.

\begin{defn} A {\it Duflo Datum } \index{Duflo's data} is a pair 
$(F,\tilde{\sigma})$,
consisting of {\it admissible and well polarizable functional } $F \in
{\frakt g}^*$ and a some irreducible representation $\tilde{\sigma} \in
\hat{G}_F$ of the stabilizer $G_F$. Recall that admissibility of $F$
means existence such an irreducible representation, and the good
polarizability of $F$ means that there
 exists a maximal possible ( therefore of half dimension of the orbit ) (
complex, solvable following M. Duflo ) subalgebra ${\frakt p}$ in
${\frakt g}_{\Bbb C}$ such that the following conditions are satisfied : 

(a) ${\frakt g}_{\Bbb C} \subset {\frakt p} \subset {\frakt g}_{\Bbb C} $.

(b) ${\frakt p}$ is $\Ad  G_F$-invariant.

(c) ${\frakt h := p \cap g}$ and ${\frakt m} := ({\frakt p} +
\overline{\frakt p})\cap {\frakt g}$ are ( real ) Lie sub-algebras. Let us
denote by $H_0$ and $M_0$, the analytic subgroups of $G$ corresponding to
${\frakt h}$ and ${\frakt m}$, respectively,
 $H := G_F \ltimes H_0$, and $M := G_F \ltimes M_0$. Then $H_0,M_0, H,
M$ are closed.

(d) There exists $\hat{\sigma} \in \widehat{H^{\frakt g}}$, such
that $$\hat{\sigma}|_{(G_F)_0^{\frakt g}} = \mult \chi_F^{\frakt
g} \quad.$$

(e) $\langle F,[{\frakt p,p}]\rangle  \equiv 0 $.

(f) There exists a complex representation $\rho$ of the complex Lie
subalgebra ${\frakt p}$, such that its restriction to the real part
${\frakt h}$ is coincided with the character ${i\over\hbar}\langle F,.\rangle $ and the
E. Nelson conditions are satisfied.
\end{defn}

\begin{rem} The condition 6.1(e) is often referred as the {\it
subordinateness condition }.In the construction of irreducible unitary
representations one assumes always the Pukanszky condition ( see Theorem
4.17 ). \end{rem}

\begin{rem}
As it was remarked, the Duflo's data are defined as the pairs of
admissible well-polarizable functional $F$ and some irreducible
representation $\tilde{\sigma}$ of $G_F$, the restriction of which to
the connected component $(G_F)_0$ are multiples of the character
$\chi_F$. This character and representation, can be lifted to the
two-fold coverings $$1 \rightarrow {\Bbb Z}/2{\Bbb Z} \cong \{
1,\varepsilon\} \rightarrow G_F^{\frakt g} \rightarrow G_F
\rightarrow 1\quad,$$ and $$1 \rightarrow {\Bbb Z}/2{\Bbb Z} \cong \{
1,\varepsilon\} \rightarrow H^{\frakt g} \rightarrow H \rightarrow
1 \quad. $$ Therefore we can have also an equivalent version of
Duflo's data.
\end{rem}

\begin{defn}  
A {\it ( lifted ) Duflo datum } \index{Duflo's data}
is a pair $(F,\hat{\sigma})$ of an admissible and well-polarizable
functional $F \in {\frakt g}^*$ and an odd irreducible unitary
representation $\hat{\sigma} \in \hat{G}_F^{\frakt g}$, such that
$$\hat{\sigma}|_{(G_F)_0^{{\frakt g}_F}} \simeq \mult \chi_F^{\frakt
g} \quad.$$ \vskip 2truecm
\end{defn}

\subsection{Metaplectic Shale-Weil representation }

 Let us consider a symplectic vector space $(V,\omega)$. To this vector
space corresponds the {\it Heisenberg Lie algebra } \index{Lie 
algebra!Heisenberg -} ${\frakt h}_{2n}(V)$,
which is $V \oplus {\Bbb R}$ as a real vector space and with the {\it
commutation relations } : $$[v \oplus t, v' \oplus t'] = \omega(v,v')
\quad.$$ Let ${\Bbb H}_{2n}(V) := \exp{{\frakt h}_{2n}(V)}$ be the
corresponding simply connected Lie group. It is called the {\it
Heisenberg group }, corresponding to the symplectic space $(V,\omega)$.

Let $p_1,p_2,\dots, p_n, q^1,\dots,q^n $ be a symplectic ``
orthonormal'' basis of $V$ with respect to the symplectic form $\omega$,
i.e. $$\omega(p_i,q^j) = \delta_{ij}, \forall i,j = 1,2,\dots,n \quad
.$$ It is easy to see that $\cent {\frakt h}_{ 2n}(V) \cong {\Bbb R}$ and
$cent{\Bbb H}_{2n}(V) \cong {\Bbb R} $. Following the well-known {\bf
Stone-von Neumann Theorem : } \index{Theorem!Stone-von Neumann -} {\it There 
exists a unique, up to unitary
equivalence, unitary representation $T$ of ${\Bbb H}_{2n}(V)$ in the
Hilbert space $L^2({\Bbb R}^n)$, such that it restriction to the center
of ${\Bbb H}_{2n}(V)$ is a multiple of the character $\chi$, define d by
$$\chi(\exp{t}) = \exp{({i\over\hbar}t)} \Id  \quad.$$ Its differential is
the representation $DT$ of the Heisenberg algebra ${\frakt h}_{2n}(V)$,
defined by : 

$p_i \mapsto \hat{p}_i := {\hbar\over i}{\partial \over \partial q^i}$,
the $i^{th}$ partial derivation,

$q^i \mapsto \widehat{q^i} := q^i$, multiplication with the coordinate
function $q^i$,

$t \mapsto \hat{t} := {\sqrt{-1} \over \hbar}t$.}

These differential operators are defined on the dense subspace ${\cal S}
({\Bbb R}^n)$ consisting of the Schwartz class functions, in the
representation space ${\cal H} = L^2({\Bbb R}^n)$.

From the point of view of the orbit method, this representation $T$
corresponds to the parameter $t^* = 1^* \in {\frakt h}_{2n}(V)^*$, where
$\{ p_1^*,\dots,p_n^*,(q^1)^*,\dots,(q^1)^*,1^*\}$ is the dual basis of
${\frakt h}_{2n}(V)^*$ and the real p polarization $({\frakt
h},\chi_{1^*})$, generated by the $p$-coordinates $p_1,\dots,p_n$.

The symplectic group $\Sp(V)$ acts on $V$ and we can consider the
semi-direct product $\Sp(V) \ltimes {\Bbb H}_{2n}(V)$. We want to extend
the representation $T$ of the normal closed subgroup ${\Bbb H}_{2n}(V)$ to
some representation of $\Sp(V)\ltimes {\Bbb H}_{2n}(V)$ and also its
subgroups. This question is solved by the Shale-Weil theorem, what
follows. Recall that the fundamental group $\pi_1\Sp(V) = {\Bbb Z}$. Then
to the quotient group ${\Bbb Z}/2{\Bbb Z} \cong \{ 1,\varepsilon\} $,
there exists a two-fold covering $\Mp (V)$ of $\Sp(V)$, called the
metaplectic group, i.e. there is an exact sequence $$1 \rightarrow {\Bbb
Z}/2{\Bbb Z} \cong \{ 1,\varepsilon\} \rightarrow \Mp (V) \rightarrow \Sp(V)
\rightarrow 1\quad.$$

\begin{thm}[\bf Shale-Weil Theorem]
There exists a unique representation $S$ of $\Mp (V)$ in ${\cal H}
 = L^2({\Bbb R}^n)$ such that, for every $\hat{x} \in \Mp (V)$, with image
$x \in \Sp(V) $ and for each $n \in {\Bbb H}_{2n}(V) $,
$$S(\hat{x})T(n)S(\hat{x})^{-1} = T(x.n) \quad,$$ $$S(\varepsilon) = -\Id
\quad.$$
\end{thm}

\begin{rem} The two-fold covering $\Mp (V)$ of $\Sp(V)$ is called the {\it
metaplectic group } and the representation $S$ of $\Mp (V)$, extending the
representation $T$ of ${\Bbb H}_{2n}(V)$ is called the {\it metaplectic
representation } \index{representation!metaplectic} or the {\it Shale-Weil 
representation. } \index{representation!Shale-Weil -} \end{rem}

\begin{cor}
Every representation $\tilde{\sigma} \in \widehat{\Sp(V)}$, realized in the Hilbert space $\tilde{\cal H}
$ can be considered as a representation of the semi-direct product $\Sp(V) \ltimes {\Bbb H}_{2n}(V)$, which is trivial on the normal subgroup ${\Bbb H}_{2n}(V)$, or as a representation of $\Mp (V) \ltimes {\Bbb H}_{2n}(V)$, which is trivial on the normal

s

subgroup $({\Bbb Z}/2{\Bbb Z}) \ltimes {\Bbb H}_{2n}(V) \cong ({\Bbb Z}/2{\Bbb Z}) \times {\Bbb H}_{2n}(V)$. Therefore, $\tilde{\sigma}\otimes S.T$ is a representation of $\Mp (V) \ltimes {\Bbb H}_{2n}(V)$ such that its restriction to ${\Bbb H}_{2n}(V)$ i

s a multiple of $T$, $$\tilde{\sigma}\otimes S.T |_{{\Bbb H}_{2n}(V)} \simeq \mult T \quad.$$
Conversely,every representation of $\Sp(V) \ltimes {\Bbb H}_{2n}(V)$,
which is a multiple of $T$ on ${\Bbb H}_{2n}(V)$ can be considered as a
projective representation of $\Sp(V)$. This projective representation
can be transformed into
a linear ( unitary ) representation, if the Mackey obstruction cocycle is
trivial. Often one can do this for Lie groups, by going to the two -
fold coverings of type $$1 \rightarrow {\Bbb Z}/2{\Bbb Z} \rightarrow
\Mp (V) \rightarrow \Sp(V) \rightarrow 1 \quad, $$ or to the
$U(1)$-coverings, $U(1) \cong {\Bbb T} \cong {\Bbb S}^1$, of type $$1
\rightarrow U(1) \rightarrow \Mp ^c(V) \rightarrow \Sp(V) \rightarrow 1 \quad
. $$ The first version will be done in the Duflo's construction. The
second version will be done in more advanced studies.  \end{cor}

This mechanism of extension of representations is
now applied to our orbit method.

\subsection{Irreducible unitary representations of nilpotent Lie 
groups }  
 
   Let $G$ be a connected and simply connected nilpotent Lie group,
${\frakt g} = \Lie G $ its Lie algebra, ${\frakt g}^*$ the dual vector
space and $F \in {\frakt g}^*$ a fixed functional. In this situation,
$G_F = (G_F)_0$, and every functional $F$ is integral. Following Theorem
4.2, there exists real polarizations ${\frakt h} \subset {\frakt g}$, which is
$\Ad G_F$-invariant, maximal isotropic with respect to the Kirillov form
$\omega_F (.,.) = \langle F,[.,.]\rangle $, satisfying the Pukanszky conditions and
which is positive in the following sense : 

For every ( not necessarily real ) polarization ${\frakt p}, \omega_F
|_{{\frakt g/g}_F}$induces some Hermitian form $i\omega_F(v,\bar{v})$, $
v \in {\frakt l} := {\frakt p}/({\frakt g})_{\Bbb C}$. The last induces a
unique Hermitian mapping on {\frakt
 l}. The {\it positivity condition } \index{condition!positivity -} asserts 
that the negative inner 
index of $\omega_F$ vanishes, $$q({\frakt p}) := \# \{ strictly \enskip
negative \enskip eigenvalues \} = 0 \quad.$$ With this Hermitian form
$\omega_F$, we can consider $\chi_F$ as an 1-dimensional unitary
representation of $H = H_0 $, in this case, and we can induce it to
obtain $\Ind(G;{\frakt h},F)$.

\begin{thm}[\bf A. Kirillov, B. Kostant, L. Auslander etc.]

 (1) The representations $\Ind(G;{\frakt h},F)$ for real polarizations and
the representations $\Ind(G;{\frakt p},\rho,F)$ for complex polarizations
, satisfying the Pukanszky condition, are irreducible.

(2) These representations are independent from polarizations ${\frakt h}$,
from the fixed point $F$, but only from the orbit $\Omega = \Omega_F$.
We denote them therefore by $T_{\Omega}$.

(3) The correspondence $$\Omega \in {\cal O} (G) \mapsto T_{\Omega} \in
\hat{G} $$ is a  homeomorphism, $${\frakt g}^*/G = {\cal O} (G) \approx
\hat{G} \quad.$$
\end{thm}

\begin{rem} We omit a detailed proof of this theorem, remarking only that
to prove (1), one reduces the representations to the tensor products of
two extended representations of some Heisenberg groups corresponding to
the polarizations. To prove (2), one remarks that by changing the point
$F\in \Omega$, the polarizations are changed by conjugations, what does
not change the equivalent classes of representations. To prove (3), one
considers in more detail the topology of the orbit space ${\frakt g}^*/G$
and the well-known continuity of the induced representation functor
$\Ind_H^G$.
\end{rem}

\subsection{Irreducible representations of extensions of nilpotent 
Lie groups }

Using the Shale-Weil representations, we extend now the representations
of nilpotent Lie groups to more general situations.

Suppose that $G$ is a connected and simply connected group, $D$ a
subgroup of automorphisms of $G$ preserving $G_F$, i.e. $D$ is a subgroup
of diffeomorphisms of the orbit $\Omega = \Omega_F$, conserving the
symplectic structure, $$D {\buildrel i \over \hookrightarrow} \Sp({\frakt
p}/({\frakt g}_F)_{\Bbb C}) \quad.$$ Following the well-known five
homomorphisms diagram lemma, the following diagram $$\vbox{\halign{
#&#&#\hfill \cr $1\rightarrow {\Bbb Z}$&$/2{\Bbb Z}$ &$D\rightarrow 1$ \cr
    &$\downarrow$                     &${\downarrow}i$ \cr 
$1\rightarrow {\Bbb Z}$&$/2{\Bbb Z} \rightarrow \Mp ({\frakt p}/({\frakt g}_F)_{\Bbb C}){\buildrel p \over \rightarrow} \Sp({\frakt p}$&$/({\frakt g}_F)_{\Bbb C}) \rightarrow 1$  \cr }}$$
can be completed to a  commutative diagram follows 
$$\vbox{\halign{ #&#&#&#\hfill\cr
$1\rightarrow {\Bbb Z}$&$/2{\Bbb Z} \rightarrow $ &$D^{\frakt g} \rightarrow $ &$D\rightarrow 1$ \cr    
  &$\downarrow$ &$\downarrow$ &${\downarrow}i$ \cr $1 \rightarrow {\Bbb
Z}$&$/2{\Bbb Z} \rightarrow \Mp ({\frakt p}$&$/({\frakt g}_F)_{\Bbb C})
{\buildrel p \over \rightarrow} \Sp({\frakt p}$&$/({\frakt g})_{\Bbb C})
\rightarrow 1\quad.$ \cr}}$$ The group $D^{\frakt g}$ is just the fibered
product of $p$ and $i$, $$D^{\frakt g} := \{ (d,\hat{x}) \in D \times
\Mp ({\frakt p}/({\frakt g})_{\Bbb C}) \quad |\quad i(d) = p(\hat{x}) (
denoted \enskip x \enskip or \enskip d_{{\frakt p}/({\frakt g})_{\Bbb C}}
) \} \quad.$$ Consider a fixed element $d\in D$. Then, following
Auslander and Kostant, there exists a {\it positive } polarization 
\index{polarization!positive -}
{\frakt p} such that ${\frakt p} + \overline{\frakt p}$ is normalized by
$d$. The operator $$(S'_{\frakt p}(d)f)(x) := |\det d_{\frakt
g/p}|^{-1/2}f(d(x)) $$ does not depend on ${\frakt p}$. Then if $U$ is a
intertwining operator between the
 representations $\Ind(G;{\frakt p},F)$ and $\Ind(G;{\frakt p}',F')$, then
$$S'_{\frakt p}(d) = U^{-1}S'_{{\frakt p}'}U$$ and $$S'_{\frakt
p}(d)\Ind(G;{\frakt p},F)(x)S'_{\frakt p}(d^{-1}) = \Ind(G;{\frakt
p},F)(d(x)) \quad.$$ Thus $$d \in D \mapsto S'_{\frakt p}(d)/{unitary
\enskip equivalence} = S'(d) $$ is a projective representation of $D$.

The character $\delta = \det (.)^{1/2}$ is well defined on $\Sp({\frakt
p}/({\frakt g})_{\Bbb C})$, because the polarization is positive and is a
character. It can be lifted easily to a character $\delta^{\frakt g}$ of
the two-fold covering such that $$\delta^l{\frakt g}(\varepsilon) = -1
\quad.$$

\begin{thm}[ M. Duflo ]  For every $\hat{d} \in D^{\frakt g}$, the
operators $$S(\hat{d}) := \delta^{\frakt g}(\varepsilon)S'(d), \forall d
\in D $$ are unitary and $S(.)$ is a unitary representation of $D^{\frakt
g}$ such that $S(\varepsilon) = -\Id $, i.e. a so called odd 
representation.
\end{thm}

\begin{cor} For every so called odd representation $U$ of $D^{\frakt
g}$, i.e. $U(\varepsilon) = -\Id  $, the tensor product $U\otimes
ST_{\Omega}$ is a representation of $D \ltimes G$, $$ (U\otimes
ST_{\Omega})(d.x) := U(\hat{d}) \otimes S(\hat{d})T_{\Omega}(x ) \quad
.$$
\end{cor}

\begin{rem} We can not go to the two-fold covering $\Mp ({\frakt
p}/({\frakt g})_{\Bbb C})$ of $\Sp({\frakt p}/({\frakt g})_{\Bbb C})$, but
to the $U(1)$-covering $\Mp ^c({\frakt p}/({\frakt g})_{\Bbb C})$ to obtain
another version of the theory.
\end{rem}

\section {Representations of Connected Reductive Groups}

This section is devoted to expose the so called Duflo's construction.
It seems to be an analytic version  of the quantization procedure in
cases, where $\sigma_0 = \chi_F$. We start with some Duflo datum
$(F,\tau)$, then reduce to the case of connected reductive groups,
where we shall start our exposition.

First of all, we recall that in  the nilpotent case, the K-orbits are
affine planes in ${\frakt g}^*$. Therefore the Pukanszky condition : $F
+ {\frakt h}^\perp \subset \Omega_F$ is automatically satisfied. The
situation is not so in the reductive case.

\begin{lem}

(1) For a reductive Lie algebra ${\frakt g}$, $F \in {\frakt g}^*$ admit a
good polarization ( i.e. satisfied the Pukanszky condition ) if and only
if the Lie algebra ${\frakt g}_F$ is a Cartan subalgebra.

(2) Let ${\frakt a}$ be an ideal of ${\frakt g}$, which is contained in
${\frakt g}_F$, ${\frakt g}' := {\frakt g/a}$, $F' \in {{\frakt g}'}^*$
is the induced from $F$ form on ${\frakt g}'$. Then $F$ has a good
polarization if and only if such one has also
 $F'$. 

(3) Suppose ${\frakt n}$ to be an nilpotent ideal of ${\frakt g}, f :=
F|_{\frakt n}$, ${\frakt g}_1 := {\frakt g}_f \subset {\frakt g}$, $F_1
:= F|_{{\frakt g}_1}$. Then $F$ has a good polarization if and only if
$F_1 \in {\frakt g}_1^*$ has a good polarization.
\end{lem}

The proof of this lemma is not difficult. We refer the reader to M.
Duflo's works.

\begin{rem} The last assertion reduces the question about existence of
good polarizations to the reductive case by induction.
\end{rem}

Let us denote $\tilde{X}_G(F)$ the set of all the unitary representations
$\tau$ of the two-fold covering $G_F^{\frakt g}$ of $G_F$, which are
{\it odd } \index{representation!odd -} in the sense that 
$$\tau(\varepsilon) = -\Id  $$ and the
restrictions of which to $(G_F)_0^{\frakt g}$ are multiples of the odd
character $\chi_F^{\frakt g} = \chi_F \delta^F$. This set is in a one -
to-one correspondence with the set $X_G(F)$ of all representations
$\tilde{\sigma}$ of the stabilizer $G_F$ such that their restrictions
 to the connected component $(G_F)_0$ are multiples of the character
$\chi_F$.  Denote its subset of the irreducible representations by
$\tilde{X}_G^{irr}(F)$, $$\tilde{X}_G^{irr}(F) := \{ \tau \in
\widehat{G_F^{\frakt g}} ; \tau(\varepsilon) = -\Id, \tau|_{(G_F)_0} =
\mult \chi_F^{\frakt g} \}$$ and the corresponding subset of $X_G(F)$,
consisting of the irreducible ones by $X_G^{irr}(F)$. Recall that $F \in
{\frakt g}^*$ is said to be {\it admissible } 
\index{functional!admissible -} iff $X_G(F) \ne \emptyset$.
Recall also that $\delta^{\frakt g} = \delta^F$ is an odd character of
$G_F^{\frakt g}$, such that $\delta(\varepsilon) = -1$ and its
differential $D\delta^F \in \sqrt{-1}{\frakt g}_F^*$, $$\delta^F(\exp{X})
= {Sh((\ad X)/2) \over (\ad X)/2}.$$
  The multiplication by the odd character $\delta^F$ realizes the
bijections between the sets $X_G(F)$ and $\tilde{X}_G(F)$, and also
$X_G^{irr}(F)$ and $\tilde{X}_G^{irr}(F)$.

\begin{rem} Suppose that there exists a positive ( with respect to the
Hermitian form $i\omega_F(v,\bar{v})$ ), $\Ad G_F$-invariant Lagrangian
subspace ${\frakt l}$ in the symplectic vector space $({\frakt g/g}_F)_{\Bbb
C} \approx (T_F\Omega)_{\Bbb C}$, then such the odd character can be
easily constructed. It is realized in the single non trivial cohomology
group of the complex associated with the adjoint action of group on (
metric ) tangent bundle of the orbit $\Omega_F$ (see the works of M. Duflo
 for more details ).  This is the case, for example, for the solvable
or semi-simple Lie groups.
\end{rem}

\begin{defn} Denote by ${\cal AP} (G)$ the set of all {\it admissible and
well-polarizable functionals } \index{functional!well-polarizable -} $F$ in 
${\frakt g}^*$. For such a
functional, a {\it Duflo's datum } \index{Duflo's data} is a pair 
$(F,\tau)$, with $\tau \in
X_G(F)$. We define ${\cal X} (G) := \cup_{F\in {\cal AP} (G)} X_G(F)$.
\end{defn}

This definition is a particular case of Definition 4.15.
We recall here for the concrete usage for reductive groups case.

Let us consider a connected reductive group $G$ and an admissible well -
polarizable functional $F \in {\cal AP} (G)$. This means that $X_G(F) \ne
\emptyset $ and ${\frakt g}_F$ is a Cartan subalgebra, denoted ${\frakt h}$
. Let $H$ be the centralizer of ${\frakt h}$ in $G$ and $\Delta =
\Delta({\frakt g}_{\Bbb C}, {\frakt h}_{\Bbb C})$ be the associated root
system. For every root $\alpha \in \Delta$, denote the corresponding root
space by ${\frakt g}^{\alpha} \ne \{0\}$. It is easy to see that the
Cartan subalgebra ${\frakt h}$ can be decomposed into the sum $${\frakt h =
t \oplus a \oplus z } \quad, $$ of the center ${\frakt z}$ of ${\frakt g}$,
the split torus ${\frakt a}$, and the compact torus ${\frakt t}$. Therefore
$\sqrt{-1}{\frakt t + a}$ will be the real linear span of co-roots. Let us
denote by $M$ the centralizer ${\cal Z} _G(A)$ of the split torus $A$,
and ${\frakt m} := \Lie M$. Then $\Delta({\frakt m}_{\Bbb C},{\frakt
h}_{\Bbb C})$ is contained in $\Delta = \Delta({\frakt g}_{\Bbb C},{\frakt
h}_{\Bbb C})$ as the subset of purely imaginary roots.

If $\Sigma$ is a subset of $\Delta$, we denote $$\rho(\Sigma) = {1\over
2} \sum_{\alpha\in{\Sigma}} \alpha$$ the half-sum of roots in $\Sigma$.
With each $\alpha \in \Sigma$ one can associate a character
$$\chi_{\alpha} := \xi_{\alpha} := \exp{\langle \alpha,.\rangle }$$ of $H = G_{\alpha} $
, the stabilizer of the point $\alpha \in {\frakt h}^* \subset{\frakt
g}^*$. The last inclusion is possible in virtue of the root
 decomposition $${\frakt g} = {\frakt h} \oplus {\sum_{\alpha \in \Delta}}
{\frakt g}^{\alpha} \quad.$$ It is reasonable to recall that $ {\frakt h}
= {\frakt g}^0$, as eigenspace.

Recall that $\lambda \in {\frakt h}^*_{\Bbb C}$ is said to be ${\frakt
g}$-{\it regular } \index{functional!regular -} iff $\lambda(H_{\alpha}) 
\ne 0, \forall \alpha \in
\Delta$, where $H_{\alpha}$ is the co-root, corresponding to $\alpha$.
This notion can be applied also to ${\frakt m}$ and we can say for example
about ${\frakt m}$-regular elements

For every $\alpha \in \Delta$, consider $${\frakt s}^{\alpha}_{\Bbb C} :=
{\frakt g}^{\alpha} \oplus {\frakt g}^{-\alpha} \oplus {\Bbb C}H_{\alpha}
\quad.$$ For each $\alpha \in \Delta_{\frakt m} := \Delta({\frakt
m}_{\Bbb C}, {\frakt h}_{\Bbb C})$, $$ {\frakt s}^{\alpha}_{\Bbb C} =
({\frakt s}^{\alpha}_{\Bbb C} \cap {\frakt g})$$ and therefore $\alpha\in
\Delta_{\frakt m}$ is said
 to be {\it compact } \index{root!compact -} ( resp., {\it noncompact } ) 
\index{root!noncompact -} root, iff $${\frakt
s}^{\alpha} := {\frakt s}^{\alpha}_{\Bbb C} \cap {\frakt g} \cong su(2) $$
( resp., $\cong sl(2,{\Bbb R})$ ).  Denote $\Delta_{{\frakt m},c}$ (
resp., $\Delta_{{\frakt m},n}$ ) the set of all compact ( resp.,
noncompact ) roots.

Consider now a functional $\lambda \in {\frakt h}^*_{\Bbb C}$ such
that $\lambda$ is ${\frakt m}$-regular and $\lambda|_{\frakt t} \in
\sqrt{-1}{\frakt t}^*$. Denote $$\Delta^+_{\frakt m}(\lambda) := \{
\alpha\in \Delta^+_{\frakt m}; \lambda(H_{\alpha}) > 0 \}\quad,$$
$$\Delta^+_{{\frakt m},n} := \Delta^+_{\frakt m}(\lambda) \cap
\Delta_{{\frakt m},n} \quad,$$ $$\Delta^+_{{\frakt m},c}(\lambda) :=
\Delta^+_{\frakt m}(\lambda) \cap \Delta{{\frakt m},c} $$ and finally
, $$\delta^{\lambda} := \rho(\Delta^+_{{\frakt m},n}(\lambda)) -
\rho(\Delta^+_{{\frakt m},c}(\lambda)) \quad.$$

Let us recall some notations from Vogan : By $R(H)$ denote the set of all
the so called $M$-{\it regular unitary pseudo-characters } 
\index{pseudo-character!M-regular -}
$(\Lambda,\lambda)$ consisting of a ${\frakt m}$-regular functional
$\lambda\in\sqrt{-1}{\frakt h}$ and a unitary representation $\Lambda$ of
$H$ with differential $$D\Lambda = (\lambda + D\delta^{\lambda}) \Id  \quad
. $$ Remark that the last condition $$ D\Lambda = (\lambda
+D\delta^{\lambda} ) \Id  $$ is equivalent to the assertion that
$$\Lambda|_{H_0} = \mult \delta^{\lambda}\chi_{\lambda} \quad,$$ what
figures in the orbit method.

Denote by $R^{irr}(H)$ the subset of $R(H)$, consisting of the
irreducible pseudo-characters. For a fixed $\lambda\in \sqrt{-1}{\frakt
h}^*$, denote $$R(H,\lambda) := \{ (\Lambda,\lambda) \in R(H); D\Lambda =
(\lambda + \delta ^{\lambda})\Id  \} $$ and $$R^{irr}(H,\lambda) :=
R(H,\lambda) \cap R^{irr}(H) \quad.$$ The Weyl group $$W(G,H) := {\cal N}
_G(H)/{\cal Z} _G(H) = {\cal N} _G(H)/H$$ acts on both $R(H)$ and
$R^{irr}(H)$.

\subsection{Harish-Chandra construction of
$\pi(\Lambda,\lambda)$ }

Recall that $$\xi_{\alpha}(.) := \exp{\langle \alpha,.\rangle }$$ for each
$\alpha\in\Delta$. Let us denote $$F := \{ x\in H ; x \enskip
centralizers \enskip {\frakt m}\enskip and \enskip |\xi_{\alpha}| = 1,
\forall \alpha \in \Delta \} \quad.$$ Then $$H := FH_0 := F \ltimes
H_0\quad,$$ $$H \cap M_0 = H_0 $$ and $$FM_0 = F \ltimes M_0\quad.$$
Denote $${\frakt k_m} = ({\frakt h}_{\Bbb C} \oplus \bigoplus_
{\alpha\in\Delta_{{\frakt m},c}} {\frakt g}^{\alpha}) \cap {\frakt g} $$
and $K_{M_0}$ the corresponding analytic subgroup. Let us denote
$\pi^{M_0} (\lambda)$ the irreducible unitary representation of $M_0$,
which is square-integrable modulo the center of $M_0$, and which is
associated with $\lambda$. This representation, following Harish-Chandra
is characterized by the following condition.

{\it The restriction $\pi^{M_0}(\lambda)|_{K_{M_0}}$ contains the ( finite
dimensional ) irreducible unitary representation of $K_{M_0}$ with the
dominant weight $ \lambda +D\delta^{\lambda}$ with respect to
$\Delta^+_{{\frakt m},c}$, as a minimal $K_{M_0}$-type }

Now a representation $\pi^{FM_0}(\Lambda,\lambda)$ of $FM_0 := F\ltimes
M_0$ can be constructed as follows $$\pi^{FM_0}(\Lambda,\lambda)(y.x) :=
\Lambda(y) \otimes \pi^{M_0}(\lambda) (x), \forall x\in M_0,y\in F \quad
.$$Let $P = MN$ be a parabolic subgroup of $G$ with the Levi component $M$
and the unipotent radical $N$, $$M = FM_0 = FM_0 \quad,$$ $$P = MN =
(FM_0) \ltimes N \quad.$$ Define now $$\pi(\Lambda,\lambda) :=
Ind^G_{FM_0 \ltimes N}(\pi^{FM_0} \otimes \Id _N) \quad.$$ Recall that {\it
if $\Lambda$ is is irreducible and $\lambda$ is ${\frakt g}$-regular, the
representation $\pi(\Lambda,\lambda)$ is irreducible. }

\begin{thm} The representation of type $\pi(\Lambda,\lambda)$ can be also
obtained from the procedure of multidimensional quantization.  \end{thm}
\begin{pf} Let us consider $F\in {\frakt g}^*$ which is admissible, i.e.
$X_G(F) \ne \emptyset$ and well-polarizable, i.e. ${\frakt g}_F = {\frakt
h} $ is a Cartan subalgebra of ${\frakt g}$. Therefore, the stabilizer
$G_F = H$ is the corresponding Cartan subgroup
 ( because $G$ is connected ! ). By the assumption, $\lambda :=
\sqrt{-1}F|_{\frakt h}$ is a ${\frakt g}$-regular functional from ${\frakt
h}^*$. There is a positive $G_F$-invariant Lagrangian subspace ${\frakt
l} \subset ({\frakt g/g}_F)_{\Bbb C}$, $${\frakt l} :=
{\bigoplus_{\alpha\in\Delta^+ _{{\frakt m},n}(\lambda)}} {\frakt
g}^{\alpha} \oplus {\bigoplus_{\alpha\in\ \Delta^+_{{\frakt
m},c}}(\lambda) } {\frakt g}^{-\alpha} \oplus {\frakt n}_{\Bbb C} \quad
,$$ where ${\frakt n} := \Lie N$. In general, ${\frakt l}$ is not a
subalgebra. But in any case, there is a character $\delta^F$ of two -
fold covering $G_F^{\frakt g} = H^{\frakt g}$ with the differential
$D\delta^F$.

Now we have $$\Lambda\in R^{irr}(H,\lambda) = X_G^{irr}(F){\buildrel 1-1
\over \leftrightarrow} \tilde{X}^{irr}(F) \quad.$$ This means that there
exists a unique $\tau \in \hat{H} = \widehat{G_F}$, such that
$$\tau|_{H_0} = \mult(\chi_F\delta^F)$$ and

 $\Lambda = \tau\delta^F$, $$ \pi(\tau\delta^F,\lambda) = \Ind(G;{\frakt
l},\tau, F) \quad.$$ The theorem is proved.
\end{pf}

\subsection{( Possibly non connected ) reductive groups } 

Let us now consider the case, where $G$ is (non-)connected reductive with
Lie algebra ${\frakt g} = \Lie G$. Consider an admissible well-polarizable
$F\in {\cal AP} (G) $, i.e. $X_G(F) \ne \emptyset$, ${\frakt g}_F =
{\frakt h}$ is a Cartan subalgebra, $H := {\cal Z} _G({\frakt h})$. In
this case, we have $$H_0 \subseteq H \cap G_0 \subseteq H \subseteq G_F
\quad.$$ Let $\lambda = \sqrt{-1}F|_{\frakt h} $, $$\Delta^+ := \{
\alpha\in\Delta ; \Im \lambda(H_{\alpha}) > 0 \enskip or \enskip \Im
\lambda(H_{\alpha}) = 0 \enskip\& \enskip \Re (\lambda(H_{\alpha})) > 0
\},$$ $${\frakt n}^+ := \sum_{\alpha\in\Delta^+} {\frakt g}^{\alpha} \quad
.$$ Then $${\frakt b}^+ := {\frakt h}_{\Bbb C} \oplus {\frakt n}^+$$ is a
$G_F$-invariant polarization at $F$ and $${\frakt e}^+ := ({\frakt b}^+ +
\overline{{\frakt b}^+}) \cap {\frakt g}$$ is a parabolic subalgebra with
the reductive part ${\cal Z} _{\frakt g}(\gamma)$, where $\gamma :=
F|_{\frakt a}$, $ {\frakt a}$ is the split component, $\gamma$ can be
considered as a linear form on ${\frakt g}$, which is null on ${\frakt t
\oplus z \oplus [h,g]}$ and on the nilpotent part $(\sum_{\alpha\in\Delta
, \alpha(\gamma) > 0} {\frakt g}^{\alpha}) \cap {\frakt g}$.

Consider now any polarization ${\frakt b}$ at $F$, which is $G_F$-stable.
Let ${\frakt e} := ({\frakt b} \oplus \overline{\frakt b})\cap {\frakt g}$
be the parabolic subalgebra of ${\frakt g}$, with unipotent radical {\frakt
u} and with reductive part ${\frakt r}$, i.e. ${\frakt r} \oplus {\frakt u} = 
{\frakt e}$.
Let $R_0$ be the analytic subgroup with Lie algebra ${\frakt r}$ and $U$ its
unipotent radical, ${\frakt u} = \Lie U$. Define $R := G_F \ltimes R_0$
and $E := (G_F \ltimes R_0)\ltimes U$, $r = F|_{\frakt r}$. it is easy
to see that the stabilizers $G_F$ at $F$ in $G$ and $R_{\frakt r}$ at $r$
in $R$ are the same, $R_{\frakt r} = G_F$. Therefore $r$ has a good
polarization. Consider the two-fold coverings $$\vbox{\halign{
#&#&#\hfill\cr $1\rightarrow {\Bbb Z}$&$/2{\Bbb Z} \rightarrow G_F^{\frakt
g} \rightarrow $ &$G_F \rightarrow 1$ \cr
   &$\parallel $ &$\parallel $ \cr $1 \rightarrow {\Bbb Z}$&$/2{\Bbb Z}
\rightarrow R_r^{\frakt r} \rightarrow $ &$R_r \rightarrow 1 \quad.$
\cr}}$$ Therefore, we can define $\hat{x} \in G^{\frakt g}_F$ and
$\tilde{x}\in R_r^{\frakt r}$ for every $x\in G_F = R_r$ and $\delta^F$
and $\delta^r$ as in the connected case, which are not necessarily to be
a character, but a odd one. By the condition of admissibility, there
exists the well-known bijection $$\tau \in X_G^{irr}(F) {\buildrel 1-1
\over \leftrightarrow } \tau' \in X_R^{irr}(r)$$ and it holds
$$\delta^F(\hat{x})\tau(\hat{x}) = \delta^r(\tilde{x})\tau'(\tilde{x})
\quad.$$ The subspace ${\frakt c := b \cap r}_{\Bbb C}$ is a Lagrangian
subspace of $({\frakt r/r}_r)_{\Bbb C}$ which is {\it totally complex }, 
\index{subspace!totally complex -}
i.e. ${\frakt c} \cap \overline{\frakt c} = \{ 0\}$. Therefore, {\frakt
h} is also a {\it fundamental Cartan subalgebra } 
\index{subalgebra!Cartan -} of ${\frakt r}$ 
and $R_0 \cap H = H_0$, x$$\tilde{X}_{R_0}^{irr}(r) = \{ \chi_{\lambda} ;
\delta^r.\chi_{\lambda} \in char(H_0),\mbox{ with differential } \lambda
= D\delta^{\lambda} \},$$ which has only a unique element $\chi_{\lambda}$. 

As in the connected case, we have a representation $\pi^0$ of $R_0$, $$
\pi^0 = \pi(\chi_{\lambda}\delta^r,\lambda) = \Ind(R_0;{\frakt
b},\chi_{\lambda}, \lambda) $$ in the space ${\cal H} $. Recall that
${\frakt c := b \cap r}_{\Bbb C}$, with the unipotent radical ${\frakt v}
= \Rad_u{\frakt c}$, is a Lagrangian subspace of $({\frakt r/r}_r)_{\Bbb
C}$, ${\frakt c} \cap \overline{\frakt c} = \{ 0\}$, i.e. is totally
complex.

Recall the notation of root spaces $$V_{\mu} := \{ v \in V ; \exists n \in
{\Bbb N} ; (H-\lambda(H))^nv = 0, \forall H \in{\frakt h} \}\quad,$$
corresponding to the root $\mu \in {\frakt h}_{\Bbb C}^*$. Let
$q^r({\frakt c})$ be the number of the 
negative eigenvalues of $\omega|_{\frakt c \times c}$ and let us
 denote $$\rho({\frakt c}) := \rho(\Sigma) ,$$ where $\Sigma$ is the
set of roots of ${\frakt h}_{\Bbb C}$ in ${\frakt v} := \Rad_u{\frakt c}$.
The space ${\cal H} ^{\infty}$ of smooth vectors in ${\cal H}$ can be
considered as some ${\frakt h}_{\Bbb C}$-module and, following Vogan,
$$\dim{H_j({\frakt v};{\cal H} ^{\infty})_{\lambda + \rho({\frakt c})}} =
\cases 0 &\mbox{ if }j\ne q^r({\frakt c}) \cr
                                    1 &\mbox{ if }j = q^r({\frakt c}) 
\endcases$$ M. Duflo proved that there exists a unique representation $S$ of
the two-fold covering $R_r^{\frakt r} \cong G_F^{\frakt g}$ in the space
${\cal H}$ denote $$\rho({\frakt c}) := \rho(\Sigma) \quad,$$ where
$\Sigma$ is the set of roots of ${\frakt h}_{\Bbb C}$ in ${\frakt v} :=
\Rad_u{\frakt c}$.  The space ${\cal H}^{\infty}$ of smooth vectors in
${\cal H}$ can be considered as some ${\frakt h}_{\Bbb C}$-module and,
following Vogan, $$\dim{H_j({\frakt v};{\cal H}  ^{\infty})_{\lambda +
\rho({\frakt c})}} = \cases 0 &\mbox{ if } j\ne q^r({\frakt c}) \cr
 1 &\mbox{ otherwise}\endcases$$ of $\pi^0$, such that
$$S(\tilde{x})\pi^0(y)S(\tilde{x})^{-1} = \pi^0(xyx^{-1}), \forall y\in
R_0, \forall x\in G_F, \forall \tilde{x}\in R^{\frakt r}_r\quad.$$
Because $R^{\frakt r}_r \cong G_F^{\frakt g}$, $S$ can be considered as a
representation of $G_F^{\frakt g}$ in posing $S(\hat{x}) := S(\tilde{x})$
. Therefore, for every odd representation $\tau' \in X_R(r)$ of
$R_r^{\frakt r}$, the formula $$\pi(xy)
:= \tau'(\tilde{x}) \otimes S(\tilde{x})\pi^0(y),\forall \tilde{x}\in
R_r^{\frakt r}, \forall y \in R_0 \quad $$ provides a representation of
semi-direct product $R = R_r \ltimes R_0 = G_F \ltimes R_0$. This
representation can be considered as the 
representation $(\tau' \otimes S\pi^0)
\otimes \Id _U$ of $E$. Then $\Ind_E^G((\tau' \otimes S\pi^0)\otimes
\Id _U)$ is a representation of $G$, which is also denoted by $T_{F,\tau,{\frakt b}} := T_{F,\tau,{\frakt b}^+}$.

\begin{cor} The representations $T_{F,\tau,{\frakt b}} :=
\Ind((\tau'\otimes S\pi^0) \otimes \Id _U)$ are obtained from the
multidimensional quantization procedure.
\end{cor}

It is reasonable to remark that {\it various concrete polarizations give
us equivalent representations. } Thus the representation
$T_{F,\tau,{\frakt b}}$ is {\it independent from the polarizations. } We
denote it from now on simply as $T^G_{F,\tau}$. The upper index $G$
indicates the group, we are interesting on.

\subsection{The induction procedure for general ( separable ) Lie 
groups }  
 
Let us now consider the case of general ( separable ) Lie groups.
Consider a Duflo's datum $(F,\tau)$, i.e. $F\in {\cal AP} (G)$ is an
admissible and well-polarizable functional on ${\frakt g} = \Lie G$ and
$\tau\in X_G^{irr}(F) $.

Let $J$ be a subgroup of $G_F$, acting trivially on ${\frakt g/g}_F$.
Following the extension $$1 \rightarrow {\Bbb Z}/2{\Bbb Z} \cong \{
1,\varepsilon \} \rightarrow G_F^{\frakt g} \rightarrow G_F \rightarrow
1$$ we can for each $x\in G_F$ choose a preimage $\hat{x} = r^F(x) $ by a
map $$r^F : J \rightarrow G_F^{\frakt g} \quad,$$ such that
$$\delta^F(r^F(x)) \equiv 1\quad.$$ The map $r^F$ can be chosen to be a
homomorphism.

{\bf Assumption} {\it For each subgroup $C$ in the center $cent G$
of $G$, if $\tau \in X_G^{irr}(F)$ is such that $$\tau(r^F(x)) \equiv
\psi(x)\Id, \forall x \in C \quad,$$ for some character $\psi \in
\hat{C}$, then so is the representation $T_{F,\tau}$, $$T_{F,\tau}(x) =
\psi(x) \Id, \forall x \in C \quad.$$ }

Consider the case $\dim{\frakt g} = 0$. Then $G$ is a discrete group,
$F=0$ and $T_{F,\tau}(x) = \tau(r^F(x)), \forall x\in G$. The assumption
is satisfied.

Suppose now that $\dim{\frakt g} > 0 $ and that $T_{F_1,\tau_1}^{G_1}$ are
defined for every $G_1$ with $\dim{\frakt g}_1 < \dim{\frakt g}$ and
$(F_1,\tau_1)$ with $F_1 \in {\cal AP} (G)$ and $\tau_1\in X_{G_1}(F_1)$,
satisfying the assumption.

Consider the Duflo datum $(F,\tau) \in {\cal X} (G) := {\cup_{F\in {\cal
AP} (G)}} X_G(F)$. Let us denote by ${\frakt n}$ the nilpotent radical of
${\frakt g}$ and $N$ the corresponding analytic subgroup, which is the
nilpotent radical of $G$ ( i.e. a closed normal analytic nilpotent
subgroup ).

Denote $$f := F|_{\frakt n}\quad,$$ $${\frakt q} := kerF \cap {\frakt
n}_f \triangleleft {\frakt g}_f \quad,$$ $${\frakt g}_1 := {\frakt
g}_f/{\frakt q} \quad.$$

{\it $1^{st}$ case : } $\dim{\frakt g}_1 = \dim{\frakt g}$. Then ${\frakt
g}$ is reductive and the representations $T^G_{F,\tau}$ are defined as in
\S7.3.

{\it $2^{nd}$ case : } $\dim{\frakt g}_1 < \dim{\frakt g}$. Let us denote
${\frakt h := g}_f$, $h := F|_{\frakt h}$ and $H := G^{\frakt n}_f $.
Because $N$ is nilpotent, the stabilizers of K-orbits are connected, and
therefore the two-fold coverings are trivial, $$N^{ \frakt n}_f \cong
{\Bbb Z}/2{\Bbb Z} \times (N_f)_0 \cong \{1,\varepsilon\} \times N_f \quad
.$$ Because $F$ is admissible, so is also $f := F|_{\frakt n}$, i.e.
there exists a character $\chi$ of $N^{\frakt n}_f$ such that
$$\chi(\varepsilon) = 1, D\chi = {i\over\hbar}f|_{{\frakt n}_f} \quad.$$
 Denote $Q := (ker \chi)_0$, the connected component of identity. It is
a normal subgroup in $H = G_f^{\frakt n}$, with Lie algebra ${\frakt q}
:= \Lie Q \cong Ker D\chi \cap {\frakt n}_f \triangleleft {\frakt g}_f $.

Define $G_1 := G^{\frakt n}_f/Q$ and $p : H = G_f^{\frakt n} \rightarrow
G_1$, the canonical projection. It is easy to see that $(G_F^{\frakt
n})^{\frakt h} \subseteq H_h^{\frakt h}$. Really we have $$H_h^{\frakt h}
= (G_F^{\frakt n})^{\frakt h} \ltimes (N_f^{\frakt n})^{\frakt h} \quad
.$$ For some element $x\in G_F$, we denote
$\hat{x}$ a representative of the preimages of $x$ in $G_F^{\frakt g}$,

\"x a representative of the preimages of of $x$ in $G_F^{\frakt n}$,

$\tilde{\ddot{x}}$ a representative of the preimages of $x$ in 
$(G_F^{\frakt n})^{\frakt h} $.

Then there exists $\tau' \in X^{irr}_{G_F^{\frakt n}}$, acting on the
space ${\cal H}_{\tau}$, such that $$\tau'(\tilde{\ddot{x}})
\delta^h(\tilde{\ddot{x}})\delta^f(\ddot{x}) =
\tau(\hat{x})\delta^F(\hat{x}) \quad.$$ Because $N^{\frakt n}_f$ acts
trivially on ${\frakt h/n}_f$, there exists a section $$r^h : N^{\frakt
n}_f \rightarrow (N^{\frakt n}_f)^{\frakt h}$$ such that
$$\delta^h(r^h(\ddot{x})) \equiv 1 \quad.$$ Therefore, there exists
$\tau' \in X_H(h)$, such that $\tau'(y)$ is defined above for all $y \in
(G_F^{\frakt n})^{\frakt h}$ and $$\tau'(r^h(y)) = \chi(y) \Id, \forall y
\in N^{\frakt n}_f \quad.$$

It is easy to see that $$(G_F)_{F_1} = H^{\frakt h}_h/r^h(Q) $$ and there
exists $\tau_1 \in X_{G_1}(F_1)$, such that $\tau' = \tau_1\circ p$. By
induction, $T^{G_1}_{F_1,\tau_1}$ is well-defined.

Consider now $U := T^{G_1}_{F_1,\tau_1} \circ p$ and $C = N^{\frakt
n}_f/Q$.  Then by the assumption, $$U(y) = \chi(y) \Id, \forall y \in
N^{\frakt n}_f \quad.$$

Let us denote by $\tilde{N}$ the universal covering of $N$. Then, we
have $$ G_f \times \tilde{N} \twoheadrightarrow G_f \ltimes \tilde{N}
\quad.$$ The representation $U\otimes S_fT^N_f$ is trivial on the kernel
of this surjection ( i.e. the preimage of the identity element ).
Therefore it defines a representation, denoted $T^{G_1}_{F_1,\tau_1}
\otimes S_fT^N_f$ of $G_f.N$.

Finally for $(F,\tau) \in {\cal X} (G)$, we define $$T^G_{F,\tau} :=
\Ind^G_{G_f.N}(T^{G_1}_{F_1,\tau_1} \otimes S_fT^N_f) \quad.$$

\begin{rem} In the case where ${\frakt g}$ is reductive and
$\dim{\frakt g}_1 < \dim{\frakt g}$, we have therefore two construction
of $T^G_{F,\tau}$. But it is easy to see that they are equivalent.
\end{rem}

We finish this section by stating the functorial properties of the
construction of representation, referring the reader to the works of M.
Duflo for a more detailed exposition.

\begin{thm}[\bf M. Duflo]

(1) For every $(F,\tau) \in {\cal X}(G)$, the commuting ring
${\cal C}(T^G_{F,\tau})$ of intertwining operators of the representation
$T^G_{F,\tau}$ is isomorphic with the same one ${\cal C}(\tau)$ of the
representation $\tau$. In particular, if $\tau$ is irreducible, so is
also the induced representation $T^G_{F,\tau}$.

(2) For every automorphism $a \in AutG $, $$T^G_{a(F,\tau)} \simeq {^a}
T^G_{F,\tau} \quad.$$ In particular, the $G$-equivalent pairs $(F,\tau)
\sim_G (F',\tau')$ give the equivalent representations $T^G_{F,\tau}
\simeq T^G_{F',\tau'}$.  Therefore the representation $T^G_{F,\tau}$
depends only on the quasi-orbit $G.(F,\tau) \subseteq {\cal X} (G)$.

(3) If $F$ is not in the K-orbit of $F'$, then the representations
$T^G_{F,\tau}$ and $T^G_{F',\tau'}$ are disjoint.

(4) There is an isomorphism of the spaces of intertwining operators
$${\cal C} (T^G_{F,\tau},T^G_{F',\tau'}) \cong {\cal C}(\tau, \tau')
\quad.$$

(5) If $C$ is a subgroup of the center $centG$ of $G$ and $$\tau(r^F(x))
\equiv \psi(x) \Id, \forall x \in C $$ for some character $\psi \in
\hat{C}$, then
 $$T^G_{F,\tau}(x) \equiv \psi(x) \Id, \forall x \in C $$ i.e. the
assumption  is always satisfied.

(6) If ${\frakt n}$ is any ideal of ${\frakt g}$ and $N$ the corresponding
analytic subgroup in $G$, $(F,\tau) \in {\cal X} (G)$, $f = F|_{\frakt
n}$, etc.... as in the beginning of \S7.4., then $$T^G_{F,\tau} =
\Ind^G_{G_f.N}(T^{G_1}_{F_1,\tau_1} \otimes S_fT^N_f) \quad.$$

(7) The representations $T^G_{F,\tau}$ are {\rm normal } ( i.e. of type
GCR ) if and only if in this construction $\dim{\tau} < \infty$ and the
K-orbit $\Omega_f = G.F$ of $F$ in ${\frakt g}^*$ is {\rm locally closed }
, i.e. an intersection of closed and open subsets.  

\end{thm}

\section{Representations of
of almost algebraic Lie groups }

This section is devoted to the so called {\it almost algebraic Lie groups
}.  \index{Lie group!almost algebraic -} For this class of Lie groups, the 
construction of irreducible
representations can be reduced to the case of reductive groups. The set
of irreducible unitary representations is so big, enough to write out the
Plancher\`el formula, one of the most important part of the harmonic
$L^2$-analysis.
\subsection{Co-isotropic subalgebra s }

\begin{defn} The real separable Lie group $G$ is called {\it
almost algebraic }, if and only if there exists a discrete subgroup
$\Gamma$ of its center ${\cal Z} (G)$ and there exists an algebraic linear
group ${\frakt G}$ over ${\Bbb R}$, such that the quotient $G/\Gamma$ can
be included as an open subgroup in the group ${\frakt G}({\Bbb R})$ of
real points of ${\frakt G}$.
\end{defn}

\begin{defn} Let $F\in {\frakt g}^*, {\frakt g} = \Lie G$.
Subalgebra ${\frakt p \subset g}$ is called {\it co-isotropic } 
\index{subalgebra!co-isotropic -} iff
${\frakt p}^F \subset {\frakt p}$, where by by definition, $${\frakt
p}^F := \{ X\in{\frakt p} ; \langle F,[X,{\frakt p}]\rangle  \equiv 0 \} \quad.$$
\end{defn}

Let us define $p:= F|_{\frakt p}$. Then it is easy to see that ${\frakt
p}_p = {\frakt p}^F$. Let $(P_p)_0$ be the analytic subgroup of $G$,
corresponding to ${\frakt p}_p = {\frakt p}^F$. Then the subset
$(P_p)_0.F$ is open in $F + {\frakt p}^\perp$. Recall that ${\frakt p}$ {\it
satisfies the Pukanszky condition } iff $$(P_p)_0.F = F + {\frakt p}^\perp
\quad.$$ For a fixed $F$, let us denote $cos(F)$ the set of all
co-isotropic subalgebra s, satisfying the Pukanszky condition.

\begin{exam}

(1) For every $F\in{\frakt g}^*$, ${\frakt g} \in cos(F)$, i.e. $cos(F)
\ne \emptyset $.

(2) If ${\frakt p}$ is a real polarization of $F$, satisfying the Pukanszky
condition, then ${\frakt p}\in cos(F)$.

(3) If ${\frakt b}$ is a complex polarization, so that ${\frakt b} +
\overline{\frakt b}$ is a complexification of ${\frakt p} := ({\frakt b} +
\overline{\frakt b}) \cap {\frakt g}$, then ${\frakt p}$ is a co-isotropic
subalgebra and ${\frakt p}_p = {\frakt b \cap g}$.

Conversely, if ${\frakt p} := ({\frakt b} + \overline{\frakt b}) \cap
{\frakt g} \in cos(F)$, then ${\frakt b}$ is a complex polarization,
satisfying the Pukanszky condition.

(4) For a semi-simple Lie algebra ${\frakt g, p \subset g}$ is a
co-isotropic subalgebra iff ${\frakt p}$ is a parabolic subalgebra.
\end{exam}

\begin{prop}
 Let ${\frakt p} \in cos(F)$. Then $F$ has a
good polarization if and only if a such one has $p = F|_{\frakt p} $.
\end{prop}
\begin{pf}
(i) {\it Sufficiency. } Suppose that $p = F|_{\frakt p}$ has a
good polarization ${\frakt b \subseteq p}_{\Bbb C}$. Then ${\frakt b}$ is
also a good polarization at $F$ in ${\frakt g}_{\Bbb C}$.

(ii) {\it Necessity. } Suppose that ${\frakt b}$ is a good polarization at
$F$ in ${\frakt g}_{\Bbb C}$. We prove that there exists a good
polarization at $p$ in ${\frakt p}_{\Bbb C}$. This will be done by
induction on dimension.

If $\dim{\frakt g} = 0$, it is obvious.

Suppose that the assertion could be proved for the Lie algebras of
strictly less dimension. There are several cases (a)-(e) :

{\it Case (a) : There exists a non-trivial unipotent ideal ${\frakt a
\subset g}$, contained in $Ker F$. } In this case, the assertion is
trivially proved by considering the quotient algebra ${\frakt g/a}$,
which has the dimension, strictly less than $\dim{\frakt g}$.

Let us denote ${\frakt z}$ the unipotent part of the center $cent{\frakt
g}$. Due to the case (a) we can restrict our consideration to the case,
where $\dim{\frakt z} \leq 1$ and and $F|_{\frakt z} \ne 0 $ if ${\frakt
z} \ne \{0\}$. Let us denote ${
\frakt u} := \Rad_u({\frakt g})$ the unipotent radical of ${\frakt g}$.

{\it Case (b) : } ${\frakt u = z}$. In this case ${\frakt g = r \times
z}$, where ${\frakt r}$ is the reductive part and ${\frakt g}_F$ is a
Cartan subalgebra of ${\frakt g}$.  Therefore ${\frakt g}_F \subset {\frakt
p}$ will be a parabolic subalgebra, and ${\frakt p}_{\Bbb C}$ contains a
Borel subalgebra ${\frakt b \subset g}_{\Bbb C}$, which is a good
polarization at $p = F|_{\frakt p} $.

{\it Case (c) : Lie algebra ${\frakt g}$ contains a characteristic
commutative ideal ${\frakt a}$ such that ${\frakt z}$ $\subsetneq$
${\frakt a \subseteq u}$ and ${\frakt g \ne p \oplus a }$. }
  In this case, let us denote

$a := F|_{\frakt a}$,

${\frakt h := g}_a \ne {\frakt g}$,

${\frakt p}' := ({\frakt p \cap h}) + {\frakt a}$,

$h := F|_{\frakt h}$.

Then ${\frakt p}' \in cos(h)$. If ${\frakt b}$ is a good polarization at
$F$, then ${\frakt b}' := ({\frakt b \cap h}_{\Bbb C} + {\frakt a}_{\Bbb
C}$ will be a good polarization at $h$. By induction, ${\frakt p}'_{\Bbb
C}$ contains a good polarization,

 say ${\frakt b}''$ at $h$ and $${\frakt b}'' \subseteq {\frakt p}_{\Bbb
C} \oplus {\frakt a}_{\Bbb C} = {\frakt g}'_{\Bbb C} \quad.$$ Let $F' :=
F|_{\frakt g}'$. Then ${\frakt b}''$ is a good polarization at $F'$ in
${\frakt g}' \subsetneq {\frakt g}$. Therefore by the induction
hypothesis, there exists a good polarization, which is contained in
${\frakt p}$.

{\it Case (d) :${\frakt u}$ is a Heisenberg algebra with center ${\frakt z
\ne u}$. } Let us define $u := F|_{\frakt u}$, ${\frakt u := g}_u$, $r
:= F|_{\frakt r}$. Let ${\frakt p}' := {\frakt p \oplus u}$, ${\frakt k
:= p' \cap r}$, ${\frakt p'' := r \oplus u}$. Then ${\frakt k} \in
cos(r) $ and ${\frakt p}'' \in cos(F) $. Let ${\frakt b}$ be a good
polarization at $F$ in ${\frakt g}_{\Bbb C}$. Then there exists a good
polarization ${\frakt b}'$ at $F$ in ${\frakt g}_{\Bbb C}$, such that
${\frakt b' \cap u}_{\Bbb C}$ is a good polarization at $u$ in ${\frakt
u}_{\Bbb C}$ and ${\frakt b}' \cap {\frakt u}_{\Bbb C}$ is a good
polarization at $r$ in ${\frakt r}_{\Bbb C}$.

Suppose ${\frakt p \ne g}$. Apply the induction process hypothesis to
${\frakt p}' $ to find a good polarization at $p$ in ${\frakt p}_{\Bbb C}$
. 

Otherwise, ${\frakt g = p' = p + u}$. If ${\frakt p = g }$, it's
trivial.  If ${\frakt g \ne p}$, let us denote the unipotent radical of
${\frakt p}$ by ${\frakt v } = {^u}{\frakt p = p \cap u}$. Then, $${\frak
t} [u,v] \subseteq [u,u] \subseteq z \subset p.$$ Therefore
${\frakt v}$ is an ideal in ${\frakt g = p \oplus u}$. Let us denote
${\frakt w} := {\cal Z} _{\frakt u}({\frakt v})$ the centralizer of
${\frakt v}$ in ${\frakt u}$. Then ${\frakt v \cap w}$ is a commutative
proper ideal in ${\frakt z}$. This is the case because in the opposite case
, ${\frakt v \cap w = z}$, and, as ${\frakt p}$-modules, there exists a
${\frakt p}$-invariant subspace ${\frakt m}$ in ${\frakt w}$, such that
${\frakt m \oplus z = w}$. Let $F'\in {\frakt g}^*$ be such a functional
that $F'|_{\frakt p} = p$, $F'|_{\frakt m} = 0$. By the Pukanszky
condition, ${\frakt p}\in cos(F')$
and ${\frakt m} \subseteq ({\frakt p}')^{F'}$, what is a contradiction.

So we can apply the case (c) to the proper ideal ${\frakt v \cap w \ne z}$.

{\it Case (e) : ${\frakt g = p}$. }
It is so clear !

{\it Now let us show that, if $\dim{\frakt g} > 0$, one of the cases (a)
- (e) must be satisfied.

non (a) : } $\dim{\frakt z} \leq 1 $, $F|_{\frakt z} \ne 0$, ${\frakt u}
:= {^u}{\frakt g}$, the unipotent radical.

{\it non (b) : }
${\frakt u \ne z}$.

{\it non (c) : } There is no characteristic ideal ${\frakt a \triangleleft
g}$, i.e.  there is no such an ideal ${\frakt a}$, that ${\frakt z
\subsetneq a \subseteq u, g = p \oplus a}$.

{\it non (d) : }
${\frakt u}$ is not a Heisenberg algebra.

{\it non (e) : }
${\frakt g \ne p}$.

We show that it must be therefore some contradiction. Really, if {\frakt
u} isn't a Heisenberg algebra, there exists some commutative ideal,
chosen of minimal dimension, ${\frakt a \triangleleft g}$, such that
${\frakt z \subseteq a \subseteq u}$. We have ${\frakt [u,a] \subseteq z}$
, by minimality of $\dim{\frakt a}$, and ${\frakt p \cap a \subseteq z}$
. Let ${\frakt } := {^u}{\frakt p} $. Then ${\frakt v}$ is contained in
${\frakt u}$, because ${\frakt p \oplus u = g}$, and any ideal of {\frakt
p}, contained in ${\frakt v}$ is an ideal in ${\frakt g}$. Therefore by {\it
non (c) } there is no commutative ideal in ${\frakt v}$, which is different
from ${\frakt z}$. So ${\frakt v}$ is a Heisenberg algebra with center {\frakt
z} or it is ${\frakt z}$ itself.

Consider ${\frakt a}' := \cent {\frakt u}$, which is a commutative ideal
in ${\frakt g}$, and ${\frakt a' \cap v = z}$, ${\frakt a' \oplus v = u }$
. Replacing ${\frakt a}'$ by ${\frakt a}$, we assume that ${\frakt [a,a]}
\equiv 0 $. Suppose ${\frakt m}$ be
 a subspace of ${\frakt a}$, such that $${\frakt a = m \oplus z} \quad,$$
$${\frakt [g,m] \subseteq m} \quad.$$ Let $F' \in {\frakt g}^*$ be such a
functional that $F'|_{\frakt p} = p$, $F'|_{\frakt m} = 0$. By the
Pukanszky condition, ${\frakt p} \in
 cos(F')$. This is a contradiction to the fact that ${\frakt m}$ is a
non-zero subspace in ${\frakt p}^{F'}$. The proposition is then proved.
\end{pf}

\subsection{Irreducible representations }

Returning to the almost algebraic Lie group $(G,\Gamma,{\frakt G})$, let
us now consider a $G_F$-stable element ${\frakt p} \in cos(F)$ and
${\frakt P}^{\circ}$ the corresponding irreducible subgroup of ${\frakt
G}$, with Lie algebra ${\frakt p}$, and $P_0$ the analytic subgroup of $G$
with the same Lie
 algebra ${\frakt p}$, $P := G_F \ltimes P_0$, ${\frakt P}$ the algebraic
closure of $P/\Gamma$ in ${\frakt G}$. Therefore ${\frakt P}({\Bbb R})$ is
closed in ${\frakt G}({\Bbb R})$ and the group $P$ is open in the inverse
image of ${\frakt P}({\Bbb R})$ in $G$, and so, it is a closed subgroup
of $G$. Finally the triple $(P,\Gamma,{\frakt P})$ is an almost algebraic
group.

Let us denote ${\frakt v} := {^u}{\frakt p}_p$ the unipotent radical of
${\frakt p}_p$, $V({\Bbb R})$ the corresponding simply connected subgroup
of ${\frakt G}({\Bbb R})$ and therefore the corresponding analytic
subgroup $V$ of $G$ is closed and simply connected. From the Pukanszky
condition, one deduces that $P_p = G_F. V$ and that $V/(V \cap G_F)$ is
simply connected.

\begin{thm}

(i) $F$ is admissible if and only if $p = F|_{\frakt p}$ is admissible.

(ii) For every $\tau \in X_G(F)$ there exists a unique $\tau' \in X_G(p)$
, such that $$(\delta^F\tau)(x) = (\delta^p\tau')(x), \forall x\in G_F
\quad.$$

(iii) The correspondence $\tau \mapsto \tau'$ in (ii) provides bijections
$$X_G(F) {\buildrel 1-1 \over \leftrightarrow } X_P(p) \quad,$$
$$X_G^{irr}(F) {\buildrel 1-1 \over \leftrightarrow} X_P^{irr}(p) \quad
.$$  
\end{thm} 
\begin{pf} 
Some parts of this lemma were checked before.
The other ones are also not difficult to be checked.
\end{pf}

\begin{thm}[ The main result] 
Suppose that $\Gamma$ is finite. Let
$(F,\tau)\in {\cal X} (G)$, ${\frakt p} \in cos(F)$ is $G_F$-stable,$p
:= F|_{\frakt p}$ and $\tau'$ corresponds to $\tau$ in the bijection
$X_G(F) {\buildrel 1-1 \over \leftrightarrow}X_P(p)$. Then $(p,\tau') \in
{\cal X} (P)$ and $$T^G_{F,\tau} = \Ind^G_P(T^P_{p,\tau'}) \quad.$$
\end{thm}

\begin{rem}
 The unipotent radical ${\frakt n} := {^u}{\frakt g}$ of
${\frakt g}$ is $G$-invariant. Let $f := F|_{\frakt n}$ and ${\frakt p :=
g}_f \oplus {\frakt n}$, $\tilde{P} := G_f.N$, $P := H_F.P_0$, where
$P_0$ is the analytic subgroup corresponding to ${\frakt p}$. It is easy to
see that $${\frakt p}\in cos(F) \quad,$$ $$\tilde{P} = P_p \quad,$$
$$T^{\tilde{P}}_{p,\tau'} = \Ind^{\tilde{P}}_P(T^P_{p,\tau'}) \quad,$$ $$
T^P_{p,\tau'} = T^{G_1}_{F_1,\tau_1} \otimes S_fT^N_f$$ and $$T^G_{F,\tau}
= \Ind^G_{\tilde{P}} \Ind^{\tilde{P}}_P(T^P_{p,\tau'}) =
\Ind^G_P(T^P_{p,\tau'}) \quad.$$
\end{rem}

\begin{rem}
If ${\frakt p} \in {\frakt g} \in cos(F)$. then
$T^G_{F,\tau} = \Ind^G_{G_FG_0}(T^{G_FG_0}_{F,\tau}) $.
\end{rem}

{\sc Proof of the Theorem. } It is the same as in the proof of
Proposition 8.4, the induction on dimension of ${\frakt g}$ can be deduced
from the cases (a)-(e).

The assertions of the theorem are trivial in the start of induction,
where $\dim{\frakt g} = 0$.

The case, where $\dim{\frakt g} > 0 $ is covered by the cases (a) -
(e). For theses cases, the induction assumptions suppose that the
theorem is proved for every Lie algebra $\tilde{\frakt p}_1$ of
dimension $\dim{\tilde{
\frakt g}_1} < \dim{\frakt g}$.

{\it Case (a). }
This case is a particular case of a previous remark. Therefore
$$T^G_{F,\tau} = \Ind^G_{G_fA_0}(T^{G_fA_0}_{p,\tau'})\quad.$$

{\it Case (b). }
 We have ${\frakt z} := {^u}\cent {\frakt g} \cong {\frakt u} :=
{^u}{\frakt g}$. In this case, as in the proof of Proposition
(b), ${\frakt g}_F$ is a Cartan subalgebra and ${\frakt p}\in cos(F)$
is a parabolic subalgebra. Therefore ${\frakt p}_{\Bbb C}$ contains a
Borel subalgebra ${\frakt b}$, which is a good polarization at $p =
F|_{\frakt p}$ and $$T^G_{F,\tau} =
\Ind^G_{G_F.P_0}(T^{G_F.P_0}_{p,\tau'})\quad.$$ A more detailed proof
of this case is contained in Bouaziz' work [Bo1] about non connected
reductive groups.

{\it Case (c). } Let us denote ${\frakt a,h,p',g'}$ as in Proposition
8.4. Denote $H:= G_F.H_0$, $ P' := G_F.P_0$, $ G' = G_FG_0$. Therefore
with $\tau \in X_G(F)$ we can associate the representations
$T^H,T^{P'},T^{G'},T^P,T^G$ of $H,P',G',P,G$, respectively. We do

 not here precise the lower indices of these representations, they are
clear from the context. From the induction assumption, we have $$T^{G'}
= \Ind^{G'}_P(T^P) = \Ind^{G'}_{P'}(T^{P'}) \quad,$$ $$T^H =
\Ind_{P'}^H(T^{P'}) \quad,$$ $$T^G = \Ind^G_H(T^H)
=\Ind_H^G(\Ind_{P'}^H(T^{P'})) = \Ind^G_{P'}(T^{P'})\quad.$$

{\it Case (d). } Let us introduce ${\frakt p',p''}$ as in Proposition
8.4(d). Define $P' := G_F.P'_0$, $P'' := G_F.P''_0$. Let
$T^{P'},T^{P''},T^P,T^G$ be the corresponding representations in giving
$\tau\in X_G(F)$. Consider {\it the sub-case : } ${\frakt p' \ne g}$.
Then by induction assumption, we have $$T^{P'} = \Ind^{P'}_{P''}(T^{P''})
= \Ind^{P'}_P(T^P) \quad.$$ Therefore it's enough to prove that $T^G =
\Ind^G_{P''}(T^{P''}) $. Let us denote ${\frakt r} := {^u}{\frakt g}$,
the unipotent radical, $R := G^{\frakt u}_u$, ${\frakt k := p' \cap r}$
,$ k := F|_{\frakt k} $, $r := F|_{\frakt r}$, $K := (R_r)_0K_0 \subset
R$. Therefore, by Lemma 8.5(ii), there exist $(r,\tau'')\in X_R(r)$,
to provide $T^R$ and $(k,\tau''')\in X_K(k)$ to provide $T^K$ and $T^G =
\Ind^G_{G_u.U}(T^R \otimes S_u.T^U_u)$. By the induction hypothesis, we
have $T^R = \Ind^R_K(T^K)$.

Consider ${\frakt p''}_u$. It is easy to see that $(P''_u)^{\frakt u} =
K$, and $P'' = P''_u.U$,$ T^{P''} = T^K.S_u.T^U_u$. Therefore
$$\begin{array}{rl} T^G &= \Ind^G_{G_u.U}(\Ind_K^R(T^K) \otimes S_uT^U_u) \cr
  &= \Ind^G_{G_u.U}(\Ind^R_{P''}(T^K \otimes S_uT^U_u)) \cr
  &= \Ind^G_{P''}(T^{P''}).\end{array}$$
For the {\it sub-case : }
${\frakt p' = g \ne p}$, it is enough to apply the case (c).

{\it Case (e) }
, ${\frakt g = p}$, is trivial.

The theorem is therefore proved.

\section{The Trace Formula and The Plancher`el Formula }
 
This section is devoted to the two central theorems of harmonic
analysis in the class of square integrable functions with integrable
module : Trace formula and the Plancher\`el formula. We propose here
only a quick review of the theory without technical proofs. It seems
to be more clear for our ``Introduction to the theory''.

\subsection{\bf 9.1. Trace formula }

Recall first of all a result of M. Duflo about his construction : {\it
The representation $T^G_{F,\tau}$ is of class $GCR$ if and only if
$\dim{\tau} < \infty$ and $\Omega_F$ is locally closed. }

Recall that the representation $T^G_{F,\tau}$ is traceable if for every
smooth function on $G$ with compact support $\alpha\in
C^{\infty}_0(G)$, the operator $$T^G_{F,\alpha}(\alpha) := \int_G
T^G_{F,\tau}(x) \alpha(x)d\mu_r(x) $$ is traceable and there exists a
distribution, noted $trT_{F,\tau}(\alpha) \in {\cal D}
(G) :=
[C^{\infty}_0(G)]' $, such that $$\langle \tr T^G_{F,\tau },\alpha\rangle  := \int_G
\tr T^G_{F,\tau}(x)\alpha(x) d\mu_r(x) := \tr T^G_{F,\tau}(\alpha)\quad.$$
For each $X\in{\frakt g}$, we introduce $$J(x)^{1/2} :=
|\det({\Sh((\ad X)/2) \over (\ad X)/2}))|^{1/2} \quad,$$ $\beta_{\Omega} :=
(1\pi)^{-d}(d!)^{-1} |\omega_{\Omega} \wedge \dots \wedge
\omega_{\Omega}|$, with $d = {1\over 2}\dim{\Omega}$, and
$\omega_{\Omega}$ is the Kirillov symplectic form on $\Omega =
\Omega_F$.

Therefore $\beta_{\Omega} $ is a positive Borelian measure on ${\frakt
g}^*$, concentrated on $\Omega = \Omega_F \subset {\frakt g}^*$.

\begin{thm}[\bf The character formula ]
Suppose that
$T_{F,\tau}$ is a normal representation and that ${\frakt g}_F$ is
nilpotent. Then $T_{F,\tau}$ is traceable if and only if the measure
$\beta_{\Omega}$ is {\rm tempered }. In this case there exists a neighborhood
of $0$ in ${\frakt g}$, in which the following formula holds
$$J(X)^{1/2}\tr T_{F,\tau}(\exp{X}) = \dim{\tau} \int_{{\frakt g}^*}
\exp{(\sqrt{-1}\langle F,X\rangle )} d\beta_{\Omega}(F)$$ in sense of distributions,
i.e. for every $\alpha\in C^{\infty}_0({\frakt g})$, with compact
support in the indicate neighborhood of $0$ in ${\frakt g}$,
$$\begin{array}{rl} \langle J^{1/2}(.) trT_{F,\tau}(\exp{(.)}),\alpha(.)\rangle
   &= \dim{\tau} \int_{\frakt g}\int_{{\frakt g}^*} \exp{(\sqrt{-1}\langle F,X\rangle )}
\alpha(X) d \beta_{\Omega}(F) dX \cr
   &= \dim{\tau} \int_{{\frakt g}^*} \hat{\alpha}(F) d\beta_{\Omega}(F)
,\end{array}$$ where $$\hat{\alpha}(F) := \int_{\frakt g} \alpha(X)
\exp{(\sqrt{-1}\langle F,X\rangle )} dX$$ is the {\rm Fourier transform }
of $\alpha$.
\end{thm}

This theorem is the so called {\it universal character formula }, 
\index{universal character formula}
proposed firstly by A. A. Kirillov in the more simple case, where
$\dim{\tau} \equiv 1$.

Suppose now, for simplicity, $G$ to be {\it connected}. Let us
denote $U({\frakt g}_{\Bbb C})$ the universal enveloping algebra,
${\cal Z}
({\frakt g}_{\Bbb C}) := \cent U({\frakt g}_{\Bbb C})$,
$S({\frakt g}_{\Bbb C})$ the symmetric algebra, on which $G$ acts and
$I({\frakt g}_{\Bbb C}) := S({\frakt g}_{\Bbb C})^G$ the set of
$G$-invariants. It is well-known that there exists an isomorphism
$$a : {\cal Z}
({\frakt g}_{\Bbb C}) {\buildrel \simeq \over
\longrightarrow} I({\frakt g}_{\Bbb C}) \quad.$$

\begin{thm}[\bf The infinitesimal character formula ]
For every $u\in {\cal Z}
({\frakt g}_{\Bbb C})$, the operators
$T_{F,\tau}(u)$, acting on $C^{\infty}$-vectors, is scalar,
$$T_{F,\tau}(u) = a(u)(\sqrt{-1}F) \Id \quad.$$
\end{thm}

\begin{rem}
If the trace formula holds for all the orbits of
maximal dimension, then the Laplace operators are just the operators
of multiplication by the
$G$-invariants $a(u), u\in {\cal Z}
({\frakt g})$. It is true also
the inverse assertion. If the Laplace operators $T_{F,\tau}(u)$, $u
\in {\cal Z}
({\frakt g})$ are the operators of multiplication by
invariants $a(u)$, on the union of all the orbit of maximal dimension
, then for the orbit of this type, the both sides of the trace
formula are the distributions satisfying the same system of elliptic
differential equations. Therefore they are coincided.
This is the
relation between the universal character formula and the infinitesimal
character formula.
\end{rem}

\subsection{Plancher\`el formula for unimodular groups }

Let us consider an almost algebraic Lie group $(G,\Gamma,{\frakt G})$
, with Lie algebra ${\frakt g} := \Lie G$. Denote ${\frakt u} :=
{^u}{\frakt g}$ the unipotent radical, ${\frakt j}$ a subalgebra such
that ${\frakt j} \subseteq \cent{\frakt g} \subseteq {\frakt u}$, $J$
the analytic subgroup ( which is closed and central ). Then $\Gamma
J$ is also closed and central, because its image in $G/\Gamma$ is
closed and central, $\gamma \cap J \subseteq \{ 1\}$.

Let us consider $\Xi\in char(\Gamma J)$, such that its differential
$d\Xi =\sqrt{-1}\xi \in \sqrt{-1}{\frakt j}^*$. Consider $G/(ker\Xi)_0$
if necessary, we can suppose that $\xi$ is injective : $\dim{\frakt j}
\leq 1$ and that, if $\dim{\frakt j}  = 1$, $ \xi \ne 0 $.

Let us denote  $$\hat{G}_{\Xi} := \{ T \in \hat{G} ; T|_{\Gamma J} 
\simeq \mult \Xi \}$$ and $${\frakt g}^*_{\xi} := \{ F \in {\frakt g}^* ; 
F|_{\frakt j} = \xi \} \quad.$$ We introduce the set of all the so called
$\Xi$-admissible functionals $$X_G(F,\Xi) := \{ \tau \in X_G(F) ; 
\tau(r^F(x)) \equiv \Xi(x)\Id  \} $$ and its subset $X_G^{irr}(F,\Xi)$
of the irreducible ones. It is easy to see that the representation
$T_{F,\tau}\in \hat{G}_{\Xi}$ if and only if $\tau \in X^{irr}_G(F,\Xi)$.

Consider  the set  of smooth   $\Xi$-functions  with compact  support
modulo $\Gamma J$, $$C^{\infty}_0(\Gamma  J  \setminus  G,\Xi) := \{ \phi
\in C^{\infty}_0 (\Gamma \setminus G) ; \phi(yx)  = \Xi(y)\phi(x),
\forall y\in  \Gamma  J, \forall  x \in G \} \quad.$$ Let us denote
$dX$ the Lebesgue measure on ${\frakt j \setminus g}$ and $dx$ the
right-invariant measure on $\Gamma \setminus G$, $L^2(\Gamma J
\setminus G, \Xi)$ the completion of $C^{\infty}_0(\Gamma J\setminus G,
\Xi)$, on which there is a natural right regular representation of $G$.

From now on, suppose that $G$ is {\it unimodular }.
\index{Lie group!unimodular -}
Let us recall that if $\pi$ is a unitary representation of $G$, and
$\phi \in L^2(\Gamma J \setminus G, \Xi)$, then $$\pi(\phi) :=
\int_{\Gamma J \setminus G} \phi(x)\pi(x) dx \quad,$$ which is
well-defined.

With the assumption about uni-modularity of $G$, there exists a unique
so called {\it Plancher\`el measure } \index{measure!Plancher\`el -}
$\mu$ on $\hat{G}_{\Xi} := \{\pi \in \hat{G} ; \pi|_{\Gamma J} = \mult 
\Xi\}$, such
that for every $\phi\in L^2(\Gamma J \setminus G,\Xi)$, the so called
{it Plancher\`el Formula }
$$\phi(1) = \int_{\hat{G}_{\Xi}} tr(\phi) d\mu(\pi)$$ holds.

\begin{rem}  If the index $(G_F : \Gamma (G_F)_0)$ is finite,
the set $X^{irr}_G(F,\Xi)$ is also finite and the Plancher\`el formula
for $\phi\in L^2((\Gamma (G_F)_0)^{\frakt g} \setminus G_F^{\frakt g}
,\Xi)$ is $$ \phi(1) = {1 \over \#\Gamma(G_F)_0 \setminus G_F}
\sum_{\tau \in X^{irr}_G(F,\Xi)} \dim{\tau}\enskip tr\tau(\phi)$$ for $\phi
\in C^{\infty}_0((\Gamma (G_F)_0)^{\frakt g} \setminus G_F^{\frakt
g},\Xi)$.
\end{rem}

Let us consider a {\it regular } \index{functional!regular -}
$F \in {\frakt g}^*$, i.e. ${\frakt g}_F$ is commutative and $\dim{\frakt
g}_F$ is minimal. If its reductive factor ${\frakt s}_F$ in ${\frakt
g}_F$ is of maximal possible dimension,
we say that $F$ is {\it strongly regular } \index{functional!strongly 
regular -} and note $F\in{\frakt g}^*_{st}$.

It is not hard to prove that ${\frakt g}^*_{st}$ is an {\it Zariski open 
set } \index{set!Zariski open -} in ${\frakt g}^*$ and there exists a {\it 
finite subset }
$\{{\frakt s}_1,{\frakt }_2,\dots, {\frakt s}_N\} \subseteq \{ {\frakt s}_F,
F\in {\frakt g}^*_{st} \}$ of representatives of conjugacy classes.

Consider now ${\cal U}
 := {\frakt g}_{st}^* \cap {\frakt g}^*_{\xi}$  and a $G$-invariant closed subset ${\cal P \subseteq U}
$, $${\cal P}
 := \{ strongly \enskip regular \enskip \& \enskip \Xi-admissible \enskip 
F \in {\frakt g}_{\xi}^* \}\quad.$$
It is not hard to see that, if $F\in {\cal P}
$, then $F$ has a good polarization. Denote $${\cal Y = Y}
_G(\Xi) := \{ (F,\tau) \in X^{irr}_G(F,\Xi) ; F \in {\cal P}
\} \subseteq {\cal X}
(G) \quad.$$ The Duflo's construction of representations gives us a map
$${\cal Y}
/G \longrightarrow \hat{G}_{\Xi}\quad,$$ $$ (F,\tau) \mapsto T^G_{F,\tau}$$
and a natural projection $${\cal Y}
/G \twoheadrightarrow {\cal P}
/G \quad,$$ $$(F,\tau) \mapsto F \quad.$$

\begin{thm}[\bf Plancher\`el Formula ]
If $G$ is a unimodular almost algebraic Lie group, there exists a
$G$-invariant function $\zeta$ on ${\cal Y}
$, with values in $(0,\infty)$, such that for every $\phi \in
C^{\infty}_0(\Gamma J \setminus G,\Xi)$, $$\psi_{\phi} := \sum_{\tau\in
X^{irr}_G(F,\Xi)} {\dim{\tau}.\zeta(F,\tau) \over \#\Gamma (G_F)_0 \setminus 
G_F } \tr T_{F,\tau}(\phi) $$ is a $G$-invariant $p(F)dF$-measurable function
on ${\cal P}
$ and $$\psi(1) = \int_{{\cal P}
/G} \psi_{\phi}(F)dm(F) \quad,$$ where $dm(F)$ is the corresponding to
$p(F)dF$ quotient measure on ${\cal P}
/G$ and $p(F)$ some function on ${\cal P}
$.
\end{thm}

This theorem is the highest point of the Duflo theory.

\section{Bibliographical Remarks}
The main idea about partially invariant holomorphly induced
representations are due to the author\cite{diep4}-\cite{diep5}. The
analytic version of this construction of representations were due to
{\sc M. Duflo}\cite{duflo1}-\cite{duflo2}.

\chapter{Reduction, Modification and Superversion}
\section{Reduction to the Semi-simple or Reductive Cases}

The aim of this appendix is to suggest a reduction of the procedure
of multidimensional quantization to the indicate cases. The geometric
construction is based on some ideas of M. Duflo about reduction. By
using a new comprehension of polarizations, we can construct the
representations of $G$, starting from the solvable or unipotent
co-isotropic tangent distributions. We shall modify also the
construction of partially invariant and holomorphly induced
representations $\Ind(G;\tilde{L},B,\sigma_0)$.

\subsection{Co-isotropic tangent distributions }

As usually let us denote by ${\frakt g}$ the Lie algebra of $G$ and by
${\frakt g}^*$ its dual vector space.  The group $G$ acts on ${\frakt
g}^*$ by the co-adjoint representation, some time K-representation.
Let $F$ be an arbitrary point in an K-orbit, say $\Omega = \Omega_F$
and $G_F$ the stabilizer at this point. Denote, as usually ${\frakt g}_F$
its Lie algebra, ${\frakt r}_F$ the radical of ${\frakt g}_F$ and
$R_F$ the corresponding analytic subgroup in $G$. Let $S_F$ be the semi-simple
component of $G_F$ in its Cartan-Levi-Mal'tsev's decomposition
$G_F = S_F \ltimes R_F$.

Recall that we have had suppose that on the principle bundle $$G_F
\rightarrowtail G \twoheadrightarrow \Omega \approx G_F\setminus G$$ a
connection $\Gamma$ ( or the same, a trivialization ) was fixed. We
have therefore an induced connection on the principal bundle
$$\vbox{\halign{ #&#&#\hfill\cr $S_F \rightarrow$&$R_F$&$\setminus G$ \cr
& &$\downarrow$ \cr &$G_F$&$\setminus G \approx \Omega$ \cr}}.$$ This
means in particular that we obtain a fixed decomposition of the tangent
bundle into the corresponding horizontal and vertical parts
$$T(R_F\setminus G) = T^h(R_F \setminus G) \oplus T^v(R_F\setminus G)
.$$ In particular the Kirillov's symplectic form $\omega_{\Omega}$ of
the K-orbit $\Omega_F$ induces a non-degenerate closed $G$-invariant
2-form $\tilde{\omega}_{\Omega}$ on $T^h(R_F\setminus G)$, defined by
the formula $$\tilde{\omega}_{\Omega}(\tilde{X},\tilde{Y}) = \omega_{\Omega}
(k_*\tilde{X},k_*\tilde{Y}),$$ where $k: R_F \setminus G
\twoheadrightarrow G_F\setminus G$ is the natural projection of the 
principal bundle $$\vbox{\halign{ #&#&#\hfill \cr $S_F \rightarrow$&$
R_F$&$\setminus G$ \cr & &${\downarrow}k$ \cr &$G_F$&$\setminus G$
\cr}}$$ and $k_*\tilde{X},k_*\tilde{Y}$ are the lifted following the
connection vector fields, which are just sections of the tangent
distribution $T^h(R_F\setminus G)$.

\begin{defn}
A smooth tangent distribution $\tilde{L} \subseteq T(R_F 
\setminus G) $ is called a {\it solvable co-isotropic distribution } 
\index{distribution!solvable co-isotropic -} iff:

i) $\tilde{L}$ is integrable and $G$-invariant,

ii) $\tilde{L}$ is $\Ad G_F$-invariant,

iii) $\tilde{L}$ is horizontal, i.e. $\tilde{L} \subseteq T^h(R_F\setminus
G)$,

iv) $\tilde{L}$ is co-isotropic at $f \in R_F\setminus G$ such that
$k(f) = F$ with respect to the form $\tilde{\omega}_{\Omega}$, i.e.
$$(\tilde{L}_f)^f\subseteq \tilde{L}_f,$$ where $(\tilde{L}_f)^f$ is
the set of all elements $\tilde{X} \in T^h(R_F \setminus G)$ such that
$$\tilde{\omega}_{\Omega}(f) (\tilde{X},\tilde{Y}) \equiv 0,\forall
\tilde{Y} \in \tilde{L}_f\quad.$$
\end{defn}

It follows from the definition that
if $\tilde{L}$ is co-isotropic at one point $f\in R_F\setminus G$, then
so is it at all the other points of $R_F\setminus G$.

\begin{thm} 
There is a one-to-one correspondence
between solvable co-isotropic distributions and $\Ad G_F$-invariant
co-isotropic Lie sub-algebras of ${\frakt g}$.
\end{thm}
\begin{pf} 
Let $\tilde{L}\subseteq T(R_F\setminus G)$ be a solvable
co-isotropic distribution. according to the Frobenius theorem,
$\tilde{L}$ is a sub-algebra of the tangent subspace $T^h_f(R_F\setminus
G)$ and then $L_F = k_*\tilde{L}_f$ is a sub-algebra of $T_F\Omega
\approx {\frakt g/g}_F$. Thus the inverse image ${\frakt b}$ of the
sub-algebra ${\frakt l} = L_F$ under the natural projection $$p:
{\frakt g} \rightarrow {\frakt g/g}_F \cong T_F\Omega $$ is an $\Ad
G_F$-invariant sub-algebra of ${\frakt g}$.

Denote by ${\frakt b}^F$ the orthogonal complement of ${\frakt b}$ in
${\frakt g}$ with respect to the Kirillov form $\omega_F$. Now we
verify the co-isotropic property of ${\frakt b}$. Indeed, let $X\in
{\frakt b}^F$. We have therefore $$\omega_F(X,Y) = 0, \forall Y\in
{\frakt b} = p^{-1}({\frakt l}),$$ then
$$\omega_{\Omega}(F)(\bar{X},\bar{Y}) = 0, \forall \bar{Y} \in {\frakt
l} = L_F; \bar{X} \in T_F\Omega.$$
 This means that $$\omega_{\Omega}(F)(k_*\tilde{X},k_*\tilde{Y}) = 0,
\forall \tilde{Y} \in \tilde{L}_f, \tilde{X}\in T^h_f(R_F\setminus G).$$
Therefore, $\tilde{X} \in (\tilde{L}_f)^f \subseteq \tilde{L}_f$,
since $ \tilde{L}$ is co-isotropic. It follows that $\bar{X} =
k_*\tilde{X} \in L_F = {\frakt l}$. Thus $X\in {\frakt b} =
p^{-1}({\frakt l})$ and ${\frakt b}^F \subseteq {\frakt b}$, i.e.
${\frakt b}$ is a co-isotropic sub-algebra.

Suppose now that ${\frakt b} \subseteq {\frakt g}$ is an $\Ad
G_F$-invariant co-isotropic sub-algebra. We define a smooth
distribution $L \subseteq T\Omega$ by the formula $$L_F:= p({\frakt
b}),$$ $$L_{F'}:= K(g)L_F; \forall F' = K(g)F \in \Omega.$$ By
definition $L$ is integrable, $G$-invariant and $\Ad G_F$-invariant.
Let us denote by $\tilde{L}$ the horizontal lift of $L$ into the
tangent bundle $T(R_F\setminus G)$, i.e. $\tilde{L} \subseteq
T^h(R_F\setminus G)$ and $k_*(\tilde{L}) = L$. Then $\tilde{L}$ is
$G$-invariant, integrable, horizontal and $\Ad G_F$-invariant.  Now
verify that $\tilde{L}$ is co-isotropic. Indeed, for every
$\tilde{X} \in (\tilde{L}_f)^f$, we have
$$\tilde{\omega}_{\Omega}(F)(\tilde{X},\tilde{Y}) = 0, \forall
\tilde{Y} \in \tilde{L}_f.$$ By definition, this means that
$$\omega_{\Omega}(F)(k_*\tilde{X},k_*\tilde{Y}) \equiv 0, \forall
k_*\tilde{Y} \in L_F $$ or $$\langle F,[X,Y]\rangle  \equiv 0, \forall Y \in
{\frakt b}; \bar{Y} = k_*\tilde{Y}.$$ Hence, $X \in {\frakt b}^F
\subseteq {\frakt b} $, since ${\frakt b}$ is co-isotropic. Thus
$\bar{X} \in L_F = k_*\tilde{L}_f$. The last means that $$\tilde{X}
\in k_*^{-1}(L_F) \cap T^h_f(R_F \setminus G) = \tilde{L}_f.$$ Then
we have $(\tilde{L}_f)^f \subseteq \tilde{L}_f$. The theorem is
therefore proved.
\end{pf}

We can do the same for the unipotent radical in place of the solvable
radical. For this aim suppose that $G$ is an algebraic Lie group.
Let $F \in {\frakt g}^*$. Denote by $U_F$ the unipotent radical of
$G_F$ and by ${\frakt u}_F$ its Lie algebra. Let $Q_F$ be the
reductive component of $G_F$ in the Cartan-Levi's decomposition
$G_F = Q_F \ltimes U_F$. Using the $Q_F$-principal bundle
$$\vbox{\halign{ #&#&#\hfill \cr $Q_F \rightarrow$&$U_F$&$\setminus G$ \cr
& &$\downarrow$ \cr &$G_F$&$\setminus G \approx \Omega$,\cr}}$$ we can
construct a non-degenerate closed $G$-invariant 2-form
${^u}\tilde{\omega}_{\Omega}$ on the horizontal component
$T^h(U_F\setminus G)$ of $T(U_F \setminus G) = T^h(U_F 
\setminus G) = T^h(U_F \setminus G) \oplus T^v(U_F \setminus G)$

\begin{defn} A smooth distribution $L \subseteq T(U_F \setminus G)$ is
called a {\it unipotent co-isotropic distribution } \index{distribution! 
unipotent co-isotropic -} if it is integrable,
$G$-invariant, $\Ad G_F$-invariant, horizontal and co-isotropic at $f$
with respect to ${^u}\tilde\omega_{\Omega} $.
\end{defn}

The following theorem can be proved by the similar arguments.

\begin{thm}
There is a one-to-one correspondence between the unipotent co-isotropic 
distributions and the $\Ad G_F$invariant co-isotropic Lie sub-algebras of the
Lie algebra ${\frakt g}$.
\end{thm}

\subsection{$(\sigma,\chi_F)$-polarizations }

\begin{defn} 
A point $F$ of ${\frakt g}^*$ is said to be
{\it r-admissible } \index{functional!r-admissible -} 
( the prefix $r-$ is for $solvable \enskip radical$ ) iff there
exists a character ( i.e.  a one -dimensional representation )
$\chi_F$ of $R_F$ such that its differential $D\chi_F$ is the
restriction ${\sqrt{-1}\over \hbar}\langle F,.\rangle |_{{\frakt r}_F}$.
\end{defn}

Denote by ${^r}X^{irr}_G(F)$ the set of all equivalent classes of
irreducible unitary representations of $G_F$ such that the restriction
of each of them to $R_F$ is a multiple of the character $\chi_F$. Then
there is one-to-one correspondence between ${^r}X^{irr}_G(F)$and the set
of all equivalent classes of irreducible projective representations of
the semi-simple Lie group $R_F \setminus G_F$.

We can do the same for the unipotent radical ${\frakt u}_F$ of ${\frakt
g}_F$.

\begin{defn} 
Let $G$ be an algebraic Lie group. We say
that a functional $F\in {\frakt g}^*$ is {\it $u$-admissible } 
\index{functional!u-admissible -} iff there
exists a character, the differential of which is equal to the
restriction of ${\sqrt{-1}\over \hbar}\langle F,.\rangle $ to the unipotent radical
${\frakt u}_F$.
\end{defn}

Note that the Lie algebra ${\frakt u}_F$ of the
unipotent radical $U_F$ is unipotent. So from the fact that in this
case the exponential map is a diffeomorphism, it follows that there
exists such a character, say $\chi_F$ of $U_F$ with the differential
equal to the restriction ${\sqrt{-1}\over \hbar}\langle F,.\rangle |_{{\frakt u}_F}$
. Hence we can say that every point $F\in {\frakt g}^*$ is
$u$-admissible.

Denote by ${^u}X^{irr}_G(F)$ the set of all equivalent classes of irreducible
unitary representation of $G_F$ such that the restriction of each of them 
to $U_F$ is a multiple of the character $\chi_F$. Following the general theory of projective representations, there exists a one-to-one correspondence
between ${^u}X^{irr}_G(F)$ and  the set of equivalent classes or irreducible projective unitary representations of the reductive group $U_F\setminus G_F$.

Since the reductive component $Q_F$ of $G_F$ has only a trivial
covering, we can identify ${^u}X^{irr}_G(F)$ with the subset of the
set of all equivalent classes of irreducible unitary representations
of $Q_F$.

With the aim to find out the irreducible unitary representations, as
usually, we consider an important generalization of the notion of
co-isotropic distribution by going over to the complex domain. This
means that we define the co-isotropic distribution $\tilde{L}$ in such
a way that $\tilde{L}_f$ is a complex subspace of the complexified
horizontal part of the tangent bundle $T^h_f(R_F\setminus G)_{\Bbb C}$
. Then the Theorems 1.1, 1.2 rest true also for these complexes
versions.

Let $\tilde{L} \subseteq T^h_f(R_F\setminus G)_{\Bbb C}$ be a
co-isotropic distribution such that $\tilde{L} \cap
\overline{\tilde{L}}$ and $\tilde{L} + \overline{\tilde{L}}$ are the
complexifications of some real distributions.  In this case, the
corresponding, following Theorem A1.2 complex sub-algebra ${\frakt p}
\subseteq {\frakt g}_{\Bbb C}$. satisfies the condition: ${\frakt p}^F
\cap \overline{\frakt p}^F$ and ${\frakt p}^F + \overline{\frakt p}^F$
are the complexifications of the real Lie algebras ${\frakt h}^F:=
{\frakt p}^F \cap {\frakt g}$ and ${\frakt m}^F:= ({\frakt p}^F +
\overline{\frakt p}^F) \cap {\frakt g}$. Denote by $H^F$ and $M^F$
the corresponding analytic subgroups in $G$.

Similarly, we can construct also the unipotent co-isotropic
distribution $\tilde{L} \subseteq T^h(U_F\setminus G)_{\Bbb C}$. Now
suppose that the sub-algebra ${\frakt h}:= {\frakt p} \cap {\frakt g}$
is algebraic co-isotropic.

\begin{defn}
A solvable ( resp., unipotent ) co-isotropic distribution $\tilde{L}$
is called {\it closed } \index{distribution!solvable 
(unipotent) co-isotropic -}
, iff all the subgroups $H^F$, $M^F$, and the semi-direct products
$H:= G_F \ltimes H^F$, $M:= G_F \ltimes M^F$ are closed in $G$.
\end{defn}

\begin{defn} Let $\tilde{\sigma}$ be some fixed
irreducible unitary representation of $G_F$ in some Hilbert space
$\tilde{V}$ such that its restriction to the radical $R_F$ ( resp.,
unipotent radical $U_F$ ) is a multiple of the character $\chi_F$.
The triple $(\tilde{L},\rho,\sigma_0)$ is called a {\it
$(\tilde{\sigma},\chi_F)$-solvable }
( resp., {\it $(\tilde{\sigma},\chi_F)$-unipotent }
) polarization, and $\tilde{L}$ is called {\it weakly Lagrangian
distribution } \index{distribution!weakly Lagrangian -}
, iff:

(1) $\sigma_0$ is an irreducible representation of the group $H^F$
in a Hilbert space $V'$, such that: (a) The restrictions are equal,
$\sigma_0|_{G_F \cap H^F} = \tilde{L}|_{G_F \cap H^F}$, (b) The point
$\sigma_0$ in the dual $\hat{H^F}$ is fixed under the natural action of 
$G_F$.

(2) $\rho$ is a representation of the complex Lie sub-algebra ${\frakt
p}^F$ in $V'$, which satisfies E. Nelson's conditions and
$\rho|_{{\frakt h}^F} = D\sigma_0$, the differential of $\sigma_0$.
\end{defn}

\begin{prop}
Suppose that $F \in {\frakt g}^*$ is r-admissible, $\tilde{L}$ is
 closed and $(\tilde{L},\rho,\sigma_0)$ is either a $(\tilde{\sigma},
\chi_F)$-solvable or $(\tilde{L},\chi_F)$-unipotent polarization. Then

(1) There exists a structure of mixed manifold of type $(k,l)$ on the
space $H^F \setminus G$, where $k = \dim{G} - \dim{M}$, $l = {1\over 2}
(\dim{M} - \dim{H^F})$.

(2) There exists a unique irreducible representation $\sigma$ of the 
subgroup $H:= G_F \ltimes H^F$ such that $$\sigma|_{G_F} = \tilde{\sigma}
, \sigma|_{H^F} = \sigma_0$$ and $$\rho|_{{\frakt p}^F} = D\sigma_0 \quad.
$$  
\end{prop}
\begin{pf}
(1). The assertion follows from Theorem 1 in Kirillov's book.

(2). Note that ${\frakt p}^F$ is invariant under the action $\Ad$ of
$G_F$ and $G_F$ acts naturally on the dual $\widehat{H^F}$ of the
group $H^F$.  From the assumptions, $\sigma_0$ is fixed under this
action of $G_F$. {\it The formula $$(x,b) {\buildrel \tau \over \mapsto}
\tilde{\sigma}(x) \sigma_0(b); \forall x\in G_F, \forall b \in H^F
,$$ define a representation of the product $G_F \times H^F$ in the
space $V = \tilde{V} \otimes V'$.}

Indeed, since $\sigma_0$ is fixed under the natural action of
$G_F$ in $\widehat{H^F}$, we have on one hand that
$$ \begin{array}{rl} \tau(x,b)\tau(x',b') &= \tilde{\sigma}(x)\sigma_0(b) 
                                      \tilde{\sigma}(x')\sigma_0(b') \cr 
   &= \tilde{\sigma}(xx')[\tilde{\sigma}(x')^{-1}\sigma_0(b)
      \tilde{\sigma}(x')]\sigma_0(b') \cr 
   &= \tilde{\sigma}(xx')[\sigma_0(b)\sigma_0(b')] \cr 
   &= \tilde{\sigma}(xx')\sigma_0(bb').\end{array}$$
On the other hand, by definition, we have
$$\tau((x,b).(x',b')) =\tau(xx',bb') = \tilde{\sigma}(xx')\sigma_0(bb').$$
It is clear that the representation $\tau$ is trivial on the kernel (, which
is by definition the inverse image of the identity element ) of the
surjection $$G_F \times H^F \longrightarrow H = G_F\ltimes H^F,$$
$$(x,b) \mapsto x.b.$$ Thus, there exists a unique representation
of the semi-direct product $H = G_F\ltimes H^F$. We denote this
representation by $\sigma$. Obviously, $\sigma$ is an irreducible
representation and $\sigma|_{G_F} = \tilde{\sigma}$, $\sigma|_{H^F}
= \sigma_0$. The proposition is proved.
\end{pf}

\subsection{Induced representations obtained from the solvable or 
unipotent polarizations } 

Suppose that $\sigma$ is the representation obtained from Proposition
1.1  Denote by $\overline{\cal E}
_{V,\sigma}:= G
\times_{H,\sigma} V$ the smooth $G$-bundle associated with the
representation $\sigma$ of $H = G_F \ltimes H^F$. Similarly, let us
consider the $G$-bundles ${\cal E} 
_{V,\sigma}:= G
\times_{G_F,\sigma|_{G_F}} V$ and $\tilde{\cal E} 
_{V,\sigma}:=
G \times_{R_F,\sigma|_{R_F}} V$ for solvable radical, or $\tilde{\cal E}
_{V,\sigma}:= G
\times_{U_F, \sigma|_{U_F}}$  for unipotent radical.

Recall that to obtain the unitary representations, we apply as usually the
construction of unitarization, considering the non-unitary character $\delta
= \sqrt{\Delta_H/\Delta_G}$ of $H$, the half-density bundle ${\cal M}
^{1/2}:= G \times_{R_F, \delta|_{R_F}} {\Bbb C}$ and finally the tensor
product $\tilde{\cal E}
_{V,\delta\sigma}:= \tilde{\cal E}
_{V,\sigma} \otimes {\cal M}
^{1/2}$. The last is called the {\it unitarized induced bundle. } 
\index{bundle!unitarized induced -}

According to the construction, the unitarized induced $G$-bundle
$\tilde{\cal E}
_{V.\delta\sigma}$
can be identified with the set of pairs $(g,v)\in G \times V$ factorized 
by the 
following equivalence relation: $(g,v) \sim (g',v')$ iff there exists
an element $h\in R_F$, such that $g' = hg$ and $v' = \delta(h)\sigma(h)v$.
Then we have a natural isomorphism of vector spaces $$\Gamma(\tilde{\cal E}
_{V,\delta\sigma}) \cong C^{\infty}(G;R_F,\delta\sigma),$$ $$\tilde{s} \mapsto
f_{\tilde{s}},$$ where $C^{\infty}(G;R_F,\delta\sigma)$ is the space of the
$V$-valued smooth functions on $G$ satisfying the following equations $$f(hg)
= \delta(h)\sigma(h)f(g), \forall h\in R_F, \forall g \in G.$$
Similarly, we can obtain the unitarized induced $G$-bundles ${\cal E}
_{V,\delta\sigma}$ and $\overline{\cal E}
_{V,\delta\sigma}$.

Recall that a section $\tilde{s} \in \Gamma(\tilde{\cal
E}_{V,\delta\sigma})$ is said to be $S_F$-{\it equivariant } 
\index{section!equivariant -} iff the
corresponding $V$-valued function $f_{\tilde{s}}$ satisfies the following
equations: $$f_{\tilde{s}}(hg) = \delta(h)\sigma(h)f_{\tilde{s}}(g),
\forall h \in R_F, \forall g \in G.$$ Denote by
$\Gamma_{S_F}(\tilde{\cal E} _{V, \delta\sigma})$ and
$\Gamma_{S_F.H^F}(\tilde{\cal E} _{V,\delta\sigma})$ the vector spaces of
$V$-valued $S_F$-equivariant or $S_F.H^F$-equivariant sections on the
unitarized induced $G$-bundle $\tilde{\cal E} _{V,\delta\sigma}$,
respectively.

\begin{prop}  There exist isomorphisms of vector spaces
$$\Gamma_{S_F}(\tilde{\cal E} _{V,\delta\sigma}) \cong \Gamma({\cal
E}_{V,\delta\sigma})$$ and $$\Gamma_{S_F.H^F}(\tilde{\cal
E}_{V,\delta\sigma}) \cong \Gamma(\overline{\cal E} _{V,\delta\sigma}).$$
\end{prop}

These assertions follow directly from the definition of $S_F$-equivariant
section and the construction of unitarized induced $G$-bundles
$\tilde{\cal E} _{V,\delta\sigma}$ and ${\cal E} _{V,\delta\sigma}$.

Fixing a connection $\overline{\Gamma}$ on the principal bundle
$$\vbox{\halign{ #&#&#\hfill\cr $B \rightarrow$& &$G$ \cr & &$\downarrow$
\cr &$B$&$\setminus G,$ \cr}}$$ we obtain the corresponding connection
$\overline{\nabla}$ on the unitarized induced $G$-bundle $\tilde{\cal E}
_{V,\delta\sigma}$. Using the natural projection $\pi: G_F \setminus G
\rightarrow H\setminus G$ and the projection $k: R_F \setminus G
\rightarrow G_F \setminus G$, we obtain the affine connection
$\tilde{\nabla}$ on the unitarized induced $G$-bundle $\tilde{\cal E}
_{V,\delta\sigma}$. We have the following diagram $$\vbox{\halign{
#&#&#&#&#&#&#&#&#&#&#&#&#&#\hfill\cr $(\overline{\cal E}
_{V}$&$_{,\delta\sigma}$&$;\overline{\nabla})$& &$({\cal E}
$&$_{V,\delta\sigma}$&$;\nabla)$& &$(\tilde{\cal E}
$&$_{V,\delta\sigma}$&$;\tilde{\nabla})$& &$(G \times V$&$;d)$ \cr
  & &$\downarrow $ & & & &$\downarrow $ & & & &$\downarrow $ & &
&$\downarrow $ \cr &$B$&$\setminus G $&${\buildrel \pi\over\leftarrow}$&
&$G_F$&$\setminus G $ &${\buildrel k\over\leftarrow}$& &$R_F$&$\setminus G
$ &$\leftarrow$ & &$G \quad.$ \cr }}$$

\begin{defn} A section $\tilde{s} \in \Gamma_{S_F.H^F}(\tilde{\cal E}
_{V,\delta\sigma})$ is called {\it partially invariant and partially
holomorphic }, \index{section!partially invariant holomorphic -} iff 
the corresponding function $f_{\tilde{s}}$ satisfies
the following equations $$[L_X + \rho(X) + D\delta(X)]f \equiv 0, \forall
X \in \overline{{\frakt p}^F}.$$
\end{defn}

Denote by ${\cal H}
:= L^2(\tilde{\cal E}_{V,\delta\sigma})
:= L^2_V(G;\tilde{L},H=G_F\ltimes H^F,\rho,\sigma_0)$ the completion of
the space of all the partially invariant and partially holomorphic
square-integrable sections of the unitarized induced $G$-bundle
$\tilde{\cal E} _{V,\delta\sigma}$. The natural unitary $G$-action on
${\cal H}$ is called the {\it partially invariant and holomorphly
induced representation } \index{representation!partially invariant 
holomorphly induced -} of $G$ and denote by 
$\Ind(G;\tilde{L},H,\rho,\sigma_0)$.

\begin{rem} Suppose that the Lie group $G$ is algebraic and $\sigma$ is
the unitary representation obtained in Proposition A1.9. corresponding to
a $(\tilde{\sigma},\chi_F)$-unipotent polarization
$(\tilde{L},\rho,\sigma)$. We can then construct also the unitarized
induced $G$-bundles $\overline{\cal E} _{V,\delta\sigma}:= G
\times_{H,\delta\sigma} V$, ${\cal E} _{V,\delta\sigma}:= G
\times_{G_F,(\delta\sigma)|_{G_F}} V$ and $\tilde{\cal E}_{V,\delta\sigma}
:= G \times_{U_F,(\delta\sigma)|_{U_F}} V$. We also denote the resulted
representation by $\Ind(G;\tilde{L},\rho,\sigma_0)$.
\end{rem}

\begin{thm} 
The representation $\Ind(G;\tilde{L}, \rho,\sigma_0)$ of
the Lie group $G$ in the space ${\cal H} $ is equivalent to the
representation of this group by the right translations in the space
$L^2(G;\tilde{L},H,\rho,\sigma_0)$ of the $V$-valued square-integrable
functions on $G$, satisfying the following equation in the sense of
distributions $$f(hg) = \delta(h)\sigma(h)f(g); \forall h \in
H:=G_F\ltimes H^F, \forall g \in G, \leqno{(1)} $$ $$[L_X + \rho(X) +
D\delta(X)]f = 0, \forall X \in \overline{\frakt p}^F, \leqno{(2)} $$
where $L_X$ is the Lie derivation along the vector field $\xi_X$ on $G$
corresponding to $X$.
\end{thm}
\begin{pf}
According to the definition, the partially invariant and partially holomorphic
sections $\tilde{s}$ are identified with the $V$-valued functions
$f_{\tilde{s}}$ on $G$ satisfying the equations (1) and (2). Then the
action of $g\in G$ on a section $\tilde{s}$ is identified with the 
action by right translation of the function $f_{\tilde{s}}$.
\end{pf}

\subsection{Unitary representations arising in the reduction of 
the multidimensional quantization procedure } 

As the model of the quantum mechanical system, we choose the space ${\cal
H} = L^2(\tilde{\cal E} _{V,\delta\sigma}) = L^2(G;\tilde{L},H=G_F\ltimes
H^F,\rho,\sigma_0)$. The ( reduction of the ) multidimensional
quantization procedure proposes the same quantization correspondence
$$\vbox{\halign{ #&#&#\hfill\cr $\hat{(.)}: C^{\infty}$&$(\Omega)$&$
\rightarrow {\cal L(H)},$ \cr
                      &$f$ &$\mapsto \hat{f}:= f +
{\hbar\over\sqrt{-1}}\tilde{\nabla}_{\tilde{\xi}_f} \quad,$\cr}}$$ where
$\tilde{\nabla}_{\tilde{\xi}_f}$ is the covariant derivation associated
with the connection $\tilde{\nabla}$ on the $G$-bundle $\tilde{\cal E}
_{V,\delta\sigma}$. We recall that it is defined by the following formula
$$\tilde{\nabla}_{\tilde{\xi}_f} = L_{\tilde{\xi}_f} +
{\sqrt{-1}\over\hbar}\alpha(\vert \tilde{\xi}_f),$$ where
${\sqrt{-1}\over\hbar}\alpha$ is the connection form of $\tilde{\nabla}$,
$L_{\tilde{\xi}_f}$ is the Lie derivation along $\tilde{\xi}_f$, which is
the horizontal lift of the strictly Hamiltonian field $\xi_f$ and $\vert
\tilde{\xi}_f$ is the vertical component of $\tilde{\xi}_f$.

We have also the following result 

\begin{thm} The following three conditions are equivalent: \item{(i)} The
application $f \mapsto \hat{f}$ is a procedure of quantization,
\item{(ii)} $\curv \tilde{\nabla}(\tilde{\xi},\tilde{\eta}) =
-{\sqrt{-1}\over\hbar}\tilde{\omega}_{\Omega}(\tilde{\xi},\tilde{\eta})Id
$, \item{(iii)}$d_{\tilde{\nabla}}\alpha(\tilde{\xi},\tilde{\eta}) =
-\tilde{\omega}_{\Omega}(\tilde{\xi},\tilde{\eta})Id$, where
$\tilde{\xi}$, $\tilde{\eta}$ are the horizontal lifts of the strictly
Hamiltonian fields $\xi$, $\eta$ and $\tilde{\omega}_{\Omega}$ is the
lifted on $\Omega$ 2-form.
\end{thm}

We obtain therefore a representation $\Lambda$ of the Lie algebra
${\frakt g}$ by the (perhaps unbounded) normal ( i.e. self -
commutative and admitting operator closure ) operators,$$\Lambda:
{\frakt g} \rightarrow {\cal L(H)},$$ $$X \mapsto \Lambda(X) =
{i\over\hbar}\hat{f}_X,$$ where $X\in{\frakt g}$ and $f_X\in
C^{\infty}(\Omega)$ is the generating function of the Hamiltonian
field $\xi_X$ corresponding to $X$.

If $G$ is connected and simply connected, we obtain a unitary
representation $T$ of $G$ defined by the following formula
$$T(\exp{X}):= \exp{(\Lambda(X))}; X\in {\frakt g}.$$ We say that
it is the {\it representation arising from the reduction of the
procedure of multidimensional quantization. } 
\index{representation!arising from quantization}

\section{Multidimensional Quantization and U(1)-Covering}

As it was remarked before, to avoid the Mackey's obstructions, M.
Duflo lifted every things to the ${\Bbb Z}/2{\Bbb Z}$-covering by using the
metaplectic group $\Mp( {\frakt g/g}_F)$. By using the technique of
{\sc P.L. Robinson and J. H. Rawnsley} \cite{robinsonrawnsley}, we shall lift all things to
the $U(1)$-covering via the $\Mp^c$-structure,i.e. the group extensions
of type $$1\rightarrow U(1) \rightarrow \Mp^c({\frakt g/g}_F)
\rightarrow \Sp({\frakt g/g}_F) \rightarrow 1$$ in place of the
metaplectic structure $\Mp$. Our purpose is to eliminate the Mackey's
obstructions of the arising projective representations to obtain
linear ones for our induction procedure. As in the ${\Bbb Z}/2{\Bbb
Z}$ case, for the Bargman-Segal model we lift the character
$\chi_F$ of the connected component of stabilizer $G_F$ to an fundamental
character $\chi_F^{U(1)}$ of $G_F^{U(1)}$, modify the notion of positive
polarization and then construct the induced representation. A
reduction of this construction is just proposed to construct
representations, starting from the so called semi-simple or unipotent.
This appendix is a revised exposition of the works \cite{vui1},\cite{vui2}.

\subsection{Positive polarizations.}

As usually, let $G$ be a connected and simply connected Lie group,
${\frakt g} = \Lie G$ its Lie algebra and ${\frakt g}^*$ the dual
vector space. The group $G$ acts on ${\frakt g}$ by the adjoint
representation $\Ad$ and on ${\frakt g}^*$ by the co-adjoint
representation $K(.):= \Ad((.)^{-1})^*$. Let $F\in {\frakt g}^*$ be an
arbitrary point in a $K$-orbit $\Omega = \Omega_F$, $G_F$ the
stabilizer of this point, ${\frakt g}_F:= \Lie G_F$ the Lie algebra
of this stabilizer, $\omega_F(X,Y):= \langle F,[X,Y]\rangle $, which has its
kernel $\Ker \omega_F = {\frakt g}_F $, $\tilde{\omega}_F$ the
corresponding symplectic form on the quotient space ${\frakt g/g}_F$ and
also on its complexification $({\frakt g/g}_F)_{\Bbb C}$ and finally
$\omega_{\Omega}$ the Kirillov 2-form on $\Omega$. Denote by $\Sp({\frakt g/g}_F
,\tilde{\omega}_F)$ the symplectic group of this symplectic vector space.

Recall
\begin{defn} 
A complex sub-algebra ${\frakt p}$ of ${\frakt
g}_{\Bbb C}$ is called a {\it positive polarization } 
\index{polarization!positive -} iff:
\begin{enumerate}
\item $({\frakt g}_F)_{\Bbb C} \subseteq {\frakt p} \subseteq {\frakt
g}_{\Bbb C}$ and it is invariant under the cation $\Ad$ of the
stabilizer $G_F$,

\item the subspace ${\frakt l}:= {\frakt p}/({\frakt g}_F)_{\Bbb C}$ of the
symplectic vector space $(T_F\Omega)_{\Bbb C} \cong ({\frakt
g/g}_F)_{\Bbb C} \cong {\frakt g}_{\Bbb C}/({\frakt g}_F)_{\Bbb C}$
satisfies the following conditions:

$(\alpha) \qquad\qquad 
\codim_{\Bbb C}{\frakt p} = {1\over 2}\dim_{\Bbb R}{\Omega} $,
 
$(\beta) \qquad\qquad
\tilde{\omega}_F(X,Y) \equiv 0, \forall X,Y \in {\frakt l}$,

$(\gamma)\qquad\qquad \sqrt{-1}\tilde{\omega}_F(X,\bar{X}) \geq 0,
\forall X\in {\frakt l}  $, where $\bar{X}$ is the
conjugation of $X$ in the complex space $({\frakt g/g}_F)_{\Bbb C}$.

\item ${\frakt p}$ is invariant in the sense that $\forall g \in G$,
$\Ad(g^{-1})({\frakt p}/({\frakt g}_{F'})_{\Bbb C})$ satisfies $(\gamma)$
in the space ${\frakt g}_{\Bbb C} / ({\frakt g}_{F'})_{\Bbb C}$, where
$F':= K(g)F$.
\end{enumerate}

We say that ${\frakt p}$ is {\it strictly positive } 
\index{polarization!strictly positive -} if the
inequality $(\gamma)$ is strict for each nonzero element $X \in
{\frakt l}$.
\end{defn}

Let us recall the standard notations ${\frakt h}:= {\frakt p \cap g = p} \cap
\overline{\frakt p \cap g}$, ${\frakt m}:=({\frakt p} + \overline{\frakt p})
\cap {\frakt g}$, $H:= G_F \ltimes H_0$ and $M:= G_F \ltimes M_0$, where
$H_0$ and $M_0$ are the analytic subgroups corresponding to the Lie
sub-algebras ${\frakt h}$ and ${\frakt m}$, respectively.

Suppose that our orbit $\Omega$ is admissible, i.e. there exists a ( perhaps
projective ) representation $U$, which is a multiple of the character $\chi_F$
in restricting to the connected component of the stabilizer 
, i.e. $$U(\exp{X}) = \chi_F(\exp{X})
Id = \exp{({i\over \hbar}\langle F,X\rangle )} Id,$$ for all $X \in {\frakt g}_F$.

We shall now, following {\sc R. L. Robinson and J. H. Rawnsley} construct
the vector space of vacuum states on which the subgroup $H$ acts. For
this reason, let us consider a positive polarization ${\frakt l}:=
{\frakt p}/({\frakt g}_F)_{\Bbb C}$ in the complexified symplectic
space $(({\frakt g/g}_F)_{\Bbb C}, \tilde{\omega}_F)$. Denote by
${\frakt d:= h/g}_F$, we have ${\frakt l} \cap \overline{\frakt l} =
{\frakt d}_{\Bbb C}$, ${\frakt d}^{\perp} = {\frakt m/g}_F =
(({\frakt p} + \overline{\frakt p}) \cap {\frakt g})/{\frakt g}_F$, where by
definition, $${\frakt d}^{\perp}:= \{ \tilde{X} \in {\frakt g/g}_F;
\tilde{\omega}(\tilde{X},\tilde{Y}) \equiv 0, \forall \tilde{Y} \in
{\frakt d} \}.$$  It is therefore true that ${\frakt d \subseteq
d^{\perp}}$, i.e. ${\frakt d}$ is an co-isotropic subspace of the
symplectic space $({\frakt g/g}_F,\tilde{\omega}_F)$. By using the
co-isotropic reduction then $\tilde{\omega}_F$ descends to a symplectic
structure on the space ${\frakt d}^{\perp}/{\frakt d}$, denoted by
$\tilde{\omega}_{F,{\frakt d}}$.We have $${\frakt d}^{\perp}/{\frakt
d} \cong ({\frakt m/g}_F)/({\frakt h/g}_F) \cong {\frakt m/h},$$ thus
there exists the symplectic structure on the space ${\frakt m/h}$ with
the strictly positive polarization $\Gamma_{\frakt d}$ of $({\frakt
d}^{\perp}/{\frakt d},\tilde{\omega}_{F,{\frakt d}})$ in the sense of
{\sc P. L. Robinson and J. H. Rawnsley} (see \cite{robinsonrawnsley},\S3] ) as
follows.

Consider the canonical projection $$\pi_{\frakt d}: {\frakt
d}^{\perp}_{\Bbb C} \rightarrow ({\frakt d}^{\perp}/{\frakt d})_{\Bbb
C},$$ $$\Gamma_{\frakt d} = \pi_{\frakt d}(L) = {\frakt l/d_{\Bbb C}}
= ({\frakt p}/({\frakt g}_F)_{\Bbb C}) /({\frakt h}_{\Bbb C}/({\frakt
g}_F)_{\Bbb C}) \cong {\frakt p/h}_{\Bbb C}.$$

\begin{prop}
The subgroup $H$ acts on the Bargman space ${\Bbb H}({\frakt m/h})$
and preserves the one dimensional subspace of the so called {\rm vacuum
states} 
$\{ {\cal E}
v'({\frakt m/h})\}^{{\frakt p/h}_{\Bbb C}}$.
\end{prop}
\begin{pf}
Since ${\frakt m/h}$ is a symplectic space, we have a rigged
Hilbert space in the sense of I. M. Gel'fand as follows $${\cal E}
({\frakt m}_{\Bbb C} / {\frakt h}_{\Bbb C} ) \subseteq {\Bbb H}({\frakt
m}_{\Bbb C} / {\frakt h}_{\Bbb C}) \subseteq {\cal E} '({\frakt m}_{\Bbb
C} / {\frakt h}_{\Bbb C})$$ In this case we have the fixed point set
$${\cal E} '({\frakt m}_{\Bbb C}/{\frakt h}_{\Bbb C})^{{\frakt p/h}_{\Bbb
C}} = \{ f\in {\cal E} '({\frakt m}_{\Bbb C} / {\frakt h}_{\Bbb C}) |
\dot{W}(v)f \equiv 0, \forall v\in {\frakt p/h}_{\Bbb C} \},$$ where
$\dot{W}: {\frakt m}_{\Bbb C} / {\frakt h}_{\Bbb C} \rightarrow {\cal E}
'({\frakt m}_{\Bbb C} / {\frakt h}_{\Bbb C})$, defined by $$(\dot{W}(v_1
+ \sqrt{-1} v_2)f)(z):= -df_z(v_1 + Jv_2) + {1\over
2\hbar}\langle z,v_1-Jv_2\rangle f(z),$$ $\langle .,.\rangle $ is the scalar product in ${\frakt
m}_{\Bbb C}/{\frakt h}_{\Bbb C}$ and $J$ denote the multiplication by $i =
\sqrt{-1}$ on ${\frakt m/h}$.

Since $\Gamma_{\frakt d} = {\frakt p/h}_{\Bbb C}$ is strictly positive
polarization, ${\cal E}
'({\frakt m}_{\Bbb C}/{\frakt h}_{\Bbb
C})^{{\frakt p/h}_{\Bbb C}}$ is a complex line with the basis vector
$f_{\frakt d}$ ( see \cite{robinsonrawnsley},\S4 ). The action of $H$ on ${\cal E}
'({\frakt m}_{\Bbb C}/{\frakt h}_{\Bbb C})^{{\frakt p/h}_{\Bbb C}}$
given by $$(U(h)f_{\frakt d})(v) = U(h) f_{\frakt d}(v), \forall v\in
{\frakt m}_{\Bbb C}/{\frakt h}_{\Bbb C}.$$ In the neighborhood of the
identity of subgroup $H$ we have $$\chi_F(\exp{X}) f_{\frakt d} =
\exp{({i\over\hbar}\langle F,X\rangle )}f_{\frakt d}.$$ We may therefore regard the
action of $H$ on ${\Bbb H}({\frakt m}_{\Bbb C}/{\frakt h}_{\Bbb C})$
given by the formula $$\exp{X} \mapsto
\exp{({i\over\hbar}\langle F,X\rangle )}f_{\frakt d}Id.$$ Therefore $H$ conserves the
one-dimensional space ${\cal E} 
'({\frakt m/h})^{{\frakt p/h}_{\Bbb
C}}$. The proposition is proved.
\end{pf}

\begin{prop} $\chi_F(X)$ can be considered as a 
unitary integral operator with $L^2$-kernel $u_X(z,w)$, $$u_X(z,w) =
\exp{({1\over \hbar}[\sqrt{-1}\langle F,X\rangle  + {1\over 2}\langle z,w\rangle  - {1\over 4}\langle w,w\rangle ])}
,$$ where $z,w \in {\frakt m}_{\Bbb C}/{\frakt h}_{\Bbb C}$.
\end{prop}
\begin{pf}
It is enough to remark that the identity operator $Id$ in ${\Bbb
H}({\frakt m}_{\Bbb C} /{\frakt h}_{\Bbb C})$ has its $L^2$-kernel
$$I(z,w) = \exp{({1\over 2\hbar}\langle z,w\rangle  -{1\over 4\hbar}\langle w,w\rangle )}.$$
\end{pf}

\begin{rem}

(1) If $(\tilde{\sigma}, \tilde{V})$ is an irreducible unitary
 representation of $G_F$, then the tensor product $\tilde{\sigma}
\otimes \chi_F$ is a representation of $G_F$ in the Hilbert space of
$\tilde{V}$- valued holomorphic functions on ${\frakt m/h}$,
square-integrable with respect to the Gaussian measure. For every $X\in
{\frakt g}_F$ the representation operator $\tilde{\sigma}(\exp{X}).
\chi_F(\exp{X})$  has its $L^2$-kernel $u(.,.)$,
$$u(z,w) = \tilde{\sigma}(\exp{X}) \exp{\{ {\sqrt{-1}\over \hbar}\langle F,X\rangle  +
{1\over 2\hbar}\langle z,w\rangle  - {1\over 4\hbar}\langle w,w\rangle  \}}.$$

(2) From the principal $H$-bundle $H\rightarrowtail M \twoheadrightarrow
H\setminus M$ and the representation $U$ of $H:= G_F \ltimes H_0$ in the
one-dimensional space ${\cal E} '({\frakt m}_{\Bbb C}/{\frakt h}_{\Bbb
C})^{{\frakt p/h}_{\Bbb C}}$ with the basis vector $f_{\frakt d}$ and the
scalar product $$(f_{\frakt d}, f_{\frakt d}) = \int_{\frakt m/h}
|f_{\frakt d}|^2 d\mu,$$ where $\mu$ is the Gaussian measure on ${\frakt
m/h}$ associated with the density function $\theta$, $$\theta(z):=
(2\pi\hbar)^{-m}\exp{(-{|z|^2 \over 2\hbar})},\forall z\in {\frakt m/h};
m:= \dim_{\Bbb C}{\frakt m/h},$$ we can construct the induced unitary
representation $\Ind^M_H$ in the space ${\Bbb H}_{\tilde{V}}( {\frakt
m}_{\Bbb C}/{\frakt h}_{\Bbb C})^{{\frakt p/h}_{\Bbb C}}$ \enskip ( see
(\cite{robinsonrawnsley},\S13) for more detail ).
\end{rem}

\subsection{Lifted characters }

Recall that if ${\frakt g}_F \ne {\frakt g}$, the symplectic group
$\Sp({\frakt g/g}_F)$ is non trivial and has the well known
$U(1)$-covering by the so called the $\Mp^c$-structure group
$\Mp ^c({\frakt g/g}_F)$, which is an extension $$1 \rightarrow U(1) \rightarrow
\Mp^c({\frakt g/g}_F) {\buildrel \pi \over
\rightarrow} \Sp({\frakt g/g}_F) \rightarrow 1 $$
(see \cite{robinsonrawnsley},\S2 for more detail ). Each element $U$ of this group
$\Mp^c({\frakt g/g}_F)$ can be presented by two parameters of type
$U = (\lambda, g)$, such that $g \in \Sp({\frakt g/g}_F)$,$ \lambda
\in {\Bbb C}$; $|\lambda^2 \det C_g| =1$; $C_g:= {1\over 2}(g - \sqrt{-1}
\circ g \circ \sqrt{-1})$. In this parameterization, the surjection $\pi$
of the extension can be precised as the projection on the second factor
$$\pi: \Mp^c({\frakt g/g}_F) \rightarrow \Sp({\frakt g/g}_F),$$
$$\pi(\lambda,g) = g $$
with kernel $\Ker \pi = U(1)$

It is clear from the definition of the stabilizer $G_F$ of the point
$F$ in the $K$-orbit $\Omega_F$ under the co-adjoint action that the following
assertion is true.

\begin{prop}
There exists a natural 
homomorphism $j: G_F \rightarrow \Sp({\frakt g/g}_F)$.
\end{prop}

\begin{prop}
There exists an extension $G_F^{U(1)}$ of $G_F$ with help of $U(1)$ such 
that the following diagram is commutative $$\vbox{\halign{ #&#&#&#&#&#&#\cr
$1\rightarrow $ &$U$&$(1) \rightarrow $ &$\quad G$&${_F^{U(1)}} {\buildrel 
\pi_j\over \rightarrow} $&$\enskip G$&$_F \rightarrow 1 $ \cr
  &  &${\downarrow}Id $ &  &${\downarrow}k $ &  &${\downarrow}j $ \cr
$1 \rightarrow $&$U$&$(1) \rightarrow $ &$\Mp^c$&$({\frakt g/g}_F) {\buildrel \pi
\over \rightarrow} $&$\Sp$&$({\frakt g/g}_F) \rightarrow 1. $\cr}}$$
\end{prop}
\begin{pf}
We have the following diagram $$\vbox{\halign{ #&#&#&#&# \cr
$1 \rightarrow $&$U$&$(1) $                        &$G$&$_F \rightarrow 1$ \cr
                &   &${\downarrow}Id $             &   &  ${\downarrow}j $ \cr
$1 \rightarrow $&$U$&$(1) \rightarrow \Mp^c({\frakt g/g}_F) {\buildrel \pi
\over \rightarrow} $ &$S$&$({\frakt g/g}_F) \rightarrow 1.$ \cr}}$$
Then from the well-known five-homomorphism lemma, we can construct
$$G_F^{U(1)}:= \{ (U,g) | \pi(U) = j(g) = \widetilde{\Ad g^{-1}}\}, $$
where by definition, $U = (\lambda,f)\in \Mp^c({\frakt g/g}_F)$ with $f\in \Sp(
{\frakt g/g}_F)$, $\lambda \in {\Bbb C}$; $$|\lambda^2 \det{\{ {1\over 2}
(f - \sqrt{-1}\circ f \circ \sqrt{-1})\}}| = 1,$$ such that the desired 
diagram is commutative.  So we have $\pi(U) = \pi(\lambda,f):= f = j(g) =
\widetilde{\Ad g^{-1}}$ and
the general form of the elements of the $U(1)$-covering $G_F^{U(1)}$ 
is $((\lambda,\widetilde{\Ad g^{-1}}),g)$.
\end{pf}

Now we consider an $K$-orbit $\Omega = \Omega_F$ passing through $F\in
{\frakt g}^*$. We do not assume that the orbit is integral. From the
exact sequence of Lie group $$1 \rightarrow U(1) \rightarrow
G_F^{U(1)} \rightarrow G_F \rightarrow 1$$ we have the corresponding
{\it split } \index{sequence!split -}
exact sequence of Lie algebras $$0 \rightarrow {\frakt u}(1) \rightarrow 
{\frakt g}_F = \Lie G_F \rightarrow {\frakt g}_F \rightarrow 0$$
Therefore ${\frakt g}_F^{U(1)} \cong {\frakt g}_F \ltimes {\frakt u}(1)$.

\begin{defn}
The $K$-orbit $\Omega_F$ is called $U(1)$-{\it admissible } 
\index{orbit!admissible -}
iff there exists a unitary character $$\chi_{F,k}^{U(1)}: G_F^{U(1)}
\rightarrow {\Bbb S}^1 = {\Bbb T}^1 \subset {\Bbb C}^{\times}, $$
such that its differential is $$D\chi_{F,k}^{U(1)}(X,\varphi) = {\sqrt{-1}
\over \hbar}(\langle F,X\rangle  + k\varphi),\forall (X,\varphi) \in {\frakt g}_F^{U(1)}
= {\frakt g}_F \ltimes {\frakt u}(1),$$ for some foxed $k \in {\Bbb Z}$.
The character $\chi_F^{U(1)}:=
\chi_{F,1}^{U(1)}$ is called {\it fundamental }. 
\index{character!fundamental -} \end{defn}

It is clear that if $k = 0 $, the $K$-orbit is integral;
some $K$-orbit could not be integral, but
$U(1)$-admissible.

\begin{prop}
In a neighborhood of identity of $G_F^{U(1)}$, we have $$\chi_F^{U(1)}
((\lambda,\widetilde{\Ad g^{-1}}),g) = \exp{\{ {\sqrt{-1}\over \hbar}(\langle F,X\rangle
+ \varphi ) \}},$$ where $\varphi \in {\Bbb R}$, $0 \leq \varphi < 2\pi\hbar$
such that $$\lambda^2 \det{(\widetilde{\Ad g^{-1}} - \sqrt{-1} \circ
\widetilde{\Ad g^{-1}} \circ \sqrt{-1})}
 = \exp{({\sqrt{-1} \over \hbar }\varphi)}.$$ As an operator, acting
on the Hilbert space ${\Bbb H}_{\tilde{V}}({\frakt m}_{\Bbb C}/{\frakt
h}_{\Bbb C} )$, it has the following integral $L^2$-kernel $$u(z,w) =
\exp{\{ {\sqrt{-1} \over \hbar}(\langle F,X\rangle  + \varphi) + {\sqrt{-1}\over 2\hbar}
\langle z,w\rangle  - {1\over 4\hbar}\langle w,w\rangle  \}},$$ for all $z,w \in {\frakt m}_{\Bbb C}/
{\frakt h}_{\Bbb C}.$
\end{prop}
\begin{pf}
This proposition is easy to check directly. Applying the operator ${d \over
dt}|_{t=0}$ to the expression $\chi_F^{U(1)}((\lambda, \widetilde{\Ad \exp{
(-tX)}}),\exp{(-tX)})$ in combining with the previous proposition, we
have these formulae.
\end{pf}

It is easily now to construct the $U(1)$-covering $H^{U(1)} = G_F^{U(1)}
\ltimes H_0
$ for the
polarization group $H:= G_F \ltimes H_0$, starting from the
homomorphism $$pr_1: H:= G_F \ltimes H_0 \rightarrow G_F$$ and the
$U(1)$- covering $$1 \rightarrow U(1) \rightarrow G_F^{U(1)}
{\buildrel \pi_j \over \rightarrow} G_F \rightarrow 1,$$ such that
the following diagram is commutative $$\vbox{\halign{ #&#&#&#&#&#&#\cr
$1 \rightarrow $&$U$&$(1) \rightarrow $&$G^{U(1)}_F$&$\ltimes H_0
\rightarrow $&$G_F$&$ \ltimes H_0 \rightarrow 1$ \cr
   &  &${\downarrow}Id $ &  &${\downarrow}\pr_1 $ &  &${\downarrow}\pr_1 $ \cr
$1 \rightarrow $&$U$&$(1) \rightarrow $&$G$&$^{U(1)}_F {\buildrel \pi_j 
\over\rightarrow }$&$G$&$_F \rightarrow$ \cr}}.$$
Remark that the Lie algebra of $H^{U(1)}$ is just ${\frakt h}^{U(1)} =
{\frakt h} \ltimes {\frakt u}(1)$.

Now assume that the $K$-orbit $\Omega_F$ is $U(1)$-admissible,
${\frakt p}$ a positive polarization in ${\frakt g}_{\Bbb C}$.

\begin{defn} 
A unitary representation $\tilde{\sigma}$ of
$G_F$ in a separable Hilbert space $\tilde{V}$, such that its
restriction to the connected component $(G_F)_0^{U(1)}$ is a multiple
of the  character $\chi_F^{ U(1)}$ is said to be {\it
fundamental. }

Fix one of this type representation $\tilde{\sigma}$. The triple $({\frakt p},\rho,\sigma_0)$ is called a ( positive )
$(\tilde{\sigma},\chi_F^{U(1)})$-{\it polarization }
iff:

(1) $({\frakt g}_F)_{\Bbb C} \subseteq {\frakt p} \subseteq {\frakt
g}_{\Bbb C}$,

(2) ${\frakt p}$ is $\Ad G_F$-invariant,

(3) ${\frakt h}:= ({\frakt p} \cap \overline{\frakt p}) \cap {\frakt g}
= {\frakt p \cap g}$ and ${\frakt m}:= ({\frakt p} + \overline{\frakt p}
) \cap {\frakt g}$ are  real Lie sub-algebras of ${\frakt g}$,

(4) $H_0$, $M_0$, $H$, $M$ are closed subgroups in $G$, where $H_0$
and $M_0$ are the analytic subgroups, corresponding to the real Lie algebras
${\frakt h,m}$, said above, $H:= G_F \ltimes H_0$, $M:= G_F \ltimes M_0$,

(5) $\langle F,[{\frakt h,h}]\rangle  \equiv 0 $, $\codim_{\frakt g} {\frakt h} =
{1\over 2} \dim \Omega_F$ and 
${\frakt p}/({\frakt g}_F)_{\Bbb C}$ is a positive polarization of the
symplectic space ${\frakt g}_{\Bbb C}/({\frakt g}_F)_{\Bbb C}$,

(6) $\sigma_0$ is an irreducible unitary representation of the $U(1)$-covering
$(H_0)^{U(1)}$ in a Hilbert space $V$ such that $$\sigma_0|_{G_F^{U(1)} \cap 
(H_0)^{U(1)}}  = \tilde{\sigma}|_{G_F^{U(1)} \cap (H_0)^{U(1)}},$$

(7) $\rho$ is a representation of the complex Lie algebra ${\frakt p \times u}
(1)$
in $V$ such that its restriction $$
\rho|_{{\frakt h \times u}(1)} = D\sigma_0|_{{\frakt h \times u}(1)}.$$
\end{defn}

As in the previous consideration, the following result is easily to be
checked.

\begin{prop} 
Let $\Omega_F$ to be $U(1)$-admissible
. The set $\tilde{X}_G(F)$ of fundamental representations is then
non-empty. Let $\tilde{\sigma}$ be an irreducible fundamental
representation and let $({\frakt p},\rho,\sigma_0)$ is a
$(\tilde{\sigma},\chi_F^{U(1)})$-polarization. Then there exists a
unique unitary irreducible representation $\sigma$ of the
$U(1)$-covering such that $$\sigma|_{G_F^{U(1)}} = \tilde{\sigma}\circ
\pi_j \chi_F^{U(1)}$$ and $$D\sigma = \rho|_{{\frakt h \times u}(1)}.$$
\end{prop}

\subsection{Induced representations }

 We introduce the following convention: For the trivial functional $F_0
= 0 \in {\frakt g}^*$, the stabilizer is therefore the whole group
$G_{F_0} = G$. By definition, $\Mp^c({\frakt g/g}_{F_0},\omega_{F_0}) =
U(1)$, we put $G^{U(1)}:= G \times U(1)$

\begin{prop}  There are isomorphisms of the fibered products $$G
\times_H H^{U(1)} \cong G \times_{G_F} G_F^{U(1)} \cong G^{U(1)} $$ in the
category of principal bundles.
\end{prop} 
\begin{pf} 
By fixing a
connection $\overline{\Gamma}$ on the principal bundle $H \rightarrowtail
G \twoheadrightarrow H\setminus G$ we have the projection map
$pr_{\overline{\Gamma}}: G \rightarrow H$. Therefore by fixing a
trivialization, we can construct the commutative diagram $$\vbox{\halign{
#&#&#&#\cr $1 \rightarrow U$&$(1) \rightarrow$&$G^{U(1)} {\buildrel
\pi_j\over\rightarrow} $&$G \rightarrow 1$ \cr
   &${\downarrow}Id $ &${\downarrow}\widetilde{\pr}_{\overline{\Gamma}} $ &
${\downarrow}\pr_{\overline{\Gamma}}$ \cr $1 \rightarrow U$&$(1)
\rightarrow $&$H^{U(1)} \rightarrow $&$H \rightarrow 1$ \cr}}$$. The rest
of the proposition is trivial.
\end{pf}

\begin{rem}
From two previous isomorphisms, we have a principal $H^{U(1)}$-bundle
over $H\setminus G$ and a principal $G_F^{U(1)}$-bundle over
$G_F\setminus G$ and a projection between two bases $$\vbox{\halign{
#&#&#&#\cr $H^{U(1)} \rightarrow$&$G^{U(1)} $ &$G_F^{U(1)} \rightarrow
$&$G^{U(1)} $\cr &$\downarrow $ & &$\downarrow $ \cr &$H\setminus G
{\buildrel \pi\over\leftarrow} $ & &$G_F\setminus G,$ \cr }}$$ where
$\pi$ is the natural projection from $G_F\setminus G$ onto $H\setminus
G$.
\end{rem}

It is therefore natural to consider the bundle ${\cal E}_{V,\sigma,\rho}
:= G^{U(1)} \times_{H^{U(1)},\sigma} V$ over $H\setminus G$,associated
with the representation $\sigma$ of $H^{U(1)}$ from Proposition A2.10.
The inverse image of this bundle $\pi^*{\cal E} _{V,\sigma,\rho}$ is a
vector bundle over the orbit $\Omega = \Omega_F = G_F \setminus G$. It is
easy to see that in the category of vector bundles, $\pi^*{\cal
E}_{V,\sigma,\rho}$ and $G^{U(1)} \times_{G_F^{U(1)},\sigma|_{G_F^{U(1)}}}
V$ are equivalent. Denote by $\Gamma(\pi^*{\cal E} _{V,\sigma,rho})$ the
space of smooth sections of the bundle $\pi^*{\cal E} _{V,\rho,\sigma}$.
The connection $\overline{\Gamma}$ on the principal $H^{U(1)}$-bundle
$H^{U(1)} \rightarrowtail G^{U(1)} \twoheadrightarrow H\setminus G$
induces a connection $\Gamma$ on the principal $G_F^{U(1)}$-bundle
$G_F^{U(1)} \rightarrowtail G^{U(1)} \twoheadrightarrow \Omega \approx
G_F\setminus G$, we then obtain an affine connection $\nabla^{\Gamma}$ on
the associated bundle $\pi^*{\cal E} _{V,\sigma,\rho} = G^{U(1)}
\times_{G_F^{U(1)},\sigma|_{G_F^{U(1)}}} V$.

Recall that there is a one-to-one correspondence between the
 $(\tilde{\sigma},\chi_F^{U(1)})$-polarizations of type
 $({\frakt p},\rho,\sigma_0)$ and the integrable $\Ad G_F$-invariant closed
weakly Lagrangian tangent distributions ${\frakt l}$. We put by
definition $$\Gamma_{\frakt l}(\pi^*{\cal E}_{V,\rho,\sigma}):= \{ s \in
\Gamma(\pi^*{\cal E} _{V,\sigma,\rho}) | \nabla_{\xi}^{\Gamma} s \equiv 0
, \forall \xi \in {\frakt l} \}.$$ The natural $G$-action on this
subspace of sections is denoted by $\Ind(G;{\frakt p},\rho,\sigma_0)$ and
will be called also as the {\it partially invariant and holomorphly
induced representation }. \index{representation!partially invariant 
holomorphly induced -} 

\begin{rem}
Now we show that with the technique developed in \cite{robinsonrawnsley} one can obtain
also these induced representations.
\end{rem}

 First we construct a canonical
representation of the metaplectic group $\Mp^c({\frakt g/g}_F)$.

 For every $v\in {\frakt g/g}_F$, consider the unitary operator $W(v)$
defined by $$(W(v)f)(z) = \exp{\{ {1\over 4\hbar}\langle 2z,w\rangle \}} f(z-v),\forall
f\in {\Bbb H}({\frakt g}_{\Bbb C}/({\frakt g}_F)_{\Bbb C}), \forall z\in
{\frakt g}_{\Bbb C}/({\frakt g}_F)_{\Bbb C}.$$ The map $$ W: {\frakt
g/g}_F \longrightarrow \Aut  {\Bbb H}({\frakt g}_{\Bbb C}/({\frakt
g}_F)_{\Bbb C})$$ is an irreducible projective unitary representation of
the vector group ${\frakt g/g}_F$ with multiplier $\exp{({1\over
\sqrt{-1}\hbar}\tilde{\omega}_F)}$. This is therefore the canonical
projective representation of the symplectic group $\Sp({\frakt g/g}_F)$ and
gives the corresponding unitary representation $\mu$ of the
$U(1)$-covering group $\Mp^c({\frakt g/g}_F)$, called the {\it (complex)
metaplectic representation.} \index{representation!(complex) metaplectic -}

Using the rigged structure $${\cal E} ({\frakt g}_{\Bbb C}/({\frakt
g}_F)_{\Bbb C}) \subset {\Bbb H}({\frakt g}_{\Bbb C}/({\frakt g}_F)_{\Bbb
C}) \subset {\cal E} '({\frakt g}_{\Bbb C}/({\frakt g}_F)_{\Bbb C})$$ we
can compute the corresponding representation $\dot{W}$ of $({\frakt
g}/({\frakt g}_F)_{\Bbb C}$ on the space ${\cal E} '({\frakt g}_{\Bbb
C}/({\frakt g}_F)_{\Bbb C})$ as follows $$(\dot{W}(v_1 +
\sqrt{-1}v_2)f)(z) = -df_z(v_1 + Jv_2) + {1\over 2\hbar} \langle z,v_1 -
Jv_2\rangle f(z), \forall z,v_1,v_2\in {\frakt g/g}_F, $$ where $f\in {\cal E}
'({\frakt g}_{\Bbb C}/({\frakt g}_F)_{\Bbb C})$.

For convenience in what follows, we write ${\frakt P}$ for the principal
$G_F^{U(1)}$-bundle $G_F^{U(1)} \rightarrowtail G^{U(1)}
\twoheadrightarrow G_F\setminus G$.  By using the homomorphism $k:
G_F^{U(1)} \rightarrow \Mp^c({\frakt g/g}_F)$ and the metaplectic
representation $$\mu: \Mp^c({\frakt g/g}_F) \rightarrow \Aut  {\cal
E}'({\frakt g}_{\Bbb C}/({\frakt g})_{\Bbb C}),$$ we have the associated
with the representation $\mu\circ k$ vector bundle with the typical fiber
${\cal E} '({\frakt g}_{\Bbb C}/({\frakt g}_F)_{\Bbb C} )$. For a
fundamental representation $(\sigma,V)$ of $G_F^{U(1)}$ the tensor product
$\sigma \otimes (\mu\circ k)$ defines a representation of $G_F^{U(1)}$ in
the space $V \otimes {\cal E} '({\frakt g}_{\Bbb C}/({\frakt g}_F)_{\Bbb
C})$. Denote by ${\cal E} '_{V,\sigma,\rho}({\frakt P})$ the vector
bundle associated with the representation $\sigma \otimes (\mu\circ k)$ of
the typical fiber $G_F^{U(1)}$ of the principal bundle ${\frakt P}$.

\begin{prop}  For each $F'\in \Omega_F$ there exists a canonical
linear map $$\dot{W}_{F'}: (T_{F'}\Omega)_{\Bbb C} \longrightarrow \End
({\cal E} '_{V,\sigma,\rho}({\frakt P}))_{F'},$$ such that
$$[\dot{W}_{F'}(\tilde{X}), \dot{W}_{F'}(\tilde{Y})] =
-{\sqrt{-1}\over\hbar}\tilde{\omega}_{F'}(\tilde{X}, \tilde{Y}), \forall
\tilde{X},\tilde{Y} \in (T_{F'}\Omega)_{\Bbb C}.$$
\end{prop} 
\begin{pf}
It enough to prove the commutation relation at one point $F$ because at
the other points it is transported in changing by conjugations, what do
not change the first. Because the representation $\sigma$ is fundamental
, its restriction to $(G_F)_0^{U(1)}$ is just a multiple of the
fundamental character $\chi_F^{U(1)}$. The commutation relation at $F$ is
therefore deduced from a direct computation.
\end{pf}

\begin{rem} We assume that our $K$-orbit is of dimension $\dim{\Omega}_F =
2m $. Recall (see \cite{robinsonrawnsley}) that the line bundle ${\cal
K}^{\tilde{L}}:= \wedge^m(\tilde{L}^{\perp})$ is a complex line sub-bundle
of the bundle $\wedge^m(T\Omega_{\Bbb C})^*$ with the basis vector section
$K_{\tilde{L}}$.
\end{rem}
  
Put $${\cal E} '_{V,\sigma,\rho}({\frakt P})^{\tilde{L}}_{F'}:= \{ f\in
{\cal E} '_{V,\sigma,\rho}({\frakt P})_{F'} | \dot{W}_{F'}(X)f \equiv 0
,\forall \tilde{X}\in \tilde{L}_{F'} \}. $$ Then ${\cal
E}'_{V,\sigma,\rho}({\frakt P})^{\tilde{L}}$ is a tensor product of the
complex line bundle ${\cal E} '({\frakt P})^{\tilde{L}}$ with the bundle
${\cal E} _{V,\sigma,\rho}({\frakt P})$, associated with the fundamental
representation $\sigma$ of the fiber $G_F^{U(1)}$ of ${\frakt P}$. It is
not hard to prove the following result.

\begin{prop}  There is a canonical isomorphism of the complex bundles
$${\cal E} '_{V,\sigma,\rho}({\frakt P})^{\tilde{L}} \otimes {\cal E}
'({\frakt P})^{\tilde{L}} \otimes {\cal K} ^{\tilde{L}} {\buildrel
\cong\over \rightarrow} {\cal E}_{V,\sigma,\rho}({\frakt P}). $$
\end{prop}

\begin{rem}
In our multidimensional situation, by
putting
$${\cal Q} ({\frakt P})^{\tilde{L}}:= {\cal E} '_{V,\sigma,\rho}({\frakt
P})^{\tilde{L}} \otimes {\cal K} ^{\tilde{L}}$$ we have then $${\cal Q}
({\frakt P})^{\tilde{L}} = {\cal E} _{V,\sigma,\rho}({\frakt P}) \otimes
{\cal E} '({\frakt P})^{\tilde{L}} \otimes {\cal K} ^{\tilde{L}}.$$ Hence
, $$({\cal Q} ({\frakt P})^{\tilde{L}})^2 = {\cal
E}_{V,\sigma,\rho}({\frakt P}) \otimes {\cal E} _{V,\sigma,\rho}({\frakt
P}) \otimes {\cal E} '({\frakt P})^{\tilde{L}} \otimes {\cal E} '({\frakt
P})^{\tilde{L}} \otimes {\cal K} ^{\tilde{L}} \otimes {\cal K}
^{\tilde{L}}.$$ Taking Proposition A2.16 into account, we get
$$\begin{array}{rl} ({\cal Q} ({\frakt P})^{\tilde{L}})^2 &= {\cal
E}_{V,\sigma,\rho}({\frakt P}) \otimes {\cal E} '_{V,\sigma,\rho}({\frakt
P})^{\tilde{L}} \otimes {\cal E} '_{V,\sigma,\rho}({\frakt P})^{\tilde{L}}
\otimes {\cal K} ^{\tilde{L}} \otimes {\cal K} ^{\tilde{L}} \cr
  &= {\cal E} _{V,\sigma,\rho}({\frakt P})^{\otimes 2} \otimes {\cal K}
^{\tilde{L}}. \end{array}$$
\end{rem}

As in the above section, on ${\cal E} _{V,\sigma,\rho}({\frakt P})$ there
is a connection $\nabla^{\Gamma}$ and $\nabla^{\Gamma} \otimes I + I
\otimes \Lie$ provides a connection on $[{\cal E} _{V,\sigma,\rho}({\frakt
P})]^{\otimes 2}$, and therefore defines uniquely
 a connection $\nabla^{\tilde{L}}$ in ${\cal Q}
({\frakt P})^{\tilde{L}}$. Denote by $\Gamma_{\tilde{L}}({\cal Q}
({\frakt P})^{\tilde{L}})$ the space of all the sections of ${\cal Q}
({\frakt P})^{\tilde{L}}$, the covariant derivatives of which vanish
along the directions of $\tilde{L}$, $$ \nabla_{\xi}^{\tilde{L}}s \equiv 0
, \forall \xi \in \tilde{L}.$$ We consider its subspace of sections
,say $s$ with compact support square module $\langle s,s\rangle _{\tilde{L}}$.Then in
this subspace, we can consider the scalar product of type
$\langle s,t\rangle _{\tilde{L}}$. The completion of this subspace with respect to this
scalar product is a Hilbert space, denoted ${\cal H} _{\tilde{L}}$.

We summarize the obtained result as follows 

\begin{thm} 
With any
$(\tilde{\sigma},\chi_F^{U(1)})$-polarization $({\frakt p},\rho,\sigma_0)$
, there exists a natural unitary representation of $G$ in ${\cal
H}_{\tilde{L}}$, denoted by $\Ind(G;{\frakt p},\rho,\sigma_0)$.
\end{thm}

\subsection{Multidimensional quantization }

As a model of the quantum system, we choose the Hilbert space ${\cal
H}_{\tilde{L}}$ as in the previous subsection. We shall use the bundle
${\cal Q} ({\frakt P})^{\tilde{L}} = {\cal E} _{V,\sigma,\rho}({\frakt P})
\otimes {\cal E} '({\frakt P})^{\tilde{L}} \otimes {\cal K} ^{\tilde{L}}$
to construct the procedure of quantization, $$\vbox{\halign{ #&#\cr
$\tilde{V} \otimes {\cal E} '({\frakt g}_{\Bbb C}/({\frakt g}_F)_{\Bbb
C})^{\tilde{L}} \otimes {\cal K} ^{\tilde{L}} \rightarrow {\cal
E}_{V,\sigma,\rho}({\frakt P}) $&$\otimes {\cal E} '({\frakt
P})^{\tilde{L}} \otimes{\cal K} ^{\tilde{L}}$ \cr
  &${\downarrow}\nabla^{\tilde{L}} $ \cr
  &$\Omega \approx G_F\setminus G.$ \cr}}$$ We define the quantization
procedure on this $U(1)$-covered situation as follows $$\vbox{\halign{
#&#&# \cr $\hat{.}: C^{\infty}$&$(\Omega)$&$\quad \rightarrow {\cal L}
({\cal H}_{\tilde{L}})$ \cr
  &$f $&$ \rightarrow \hat{f}:= f + {\hbar\over
\sqrt{-1}}\nabla_{\xi_f}^{\tilde{L}},$ \cr}}$$ where ${\cal L} ({\cal
H}_{\tilde{L}})$ is the space of ( perhaps unbounded ) Hermitian operators
,
 admitting operator closure on ${\cal H}_{\tilde{L}}$ and
$\nabla_{\xi_f}^{\tilde{L}}$ is the covariant derivation associated with
the connection $\nabla^{\tilde{L}}$ on the bundle ${\cal Q} ({\frakt
P})^{\tilde{L}}$. Recall that $$\nabla_{\xi_f}^{\tilde{L}} = L_{\xi_f} +
{\sqrt{-1} \over \hbar}\alpha(\vert(\xi_f)).$$ The first summand is the
well-known Lie derivation and the second summand is the connection form.

By a similar argument as in previous section, we have the following
result 

\begin{prop}
The following three conditions are equivalent:

(i) $d_{\nabla^{\tilde{L}}} \alpha(\xi,\eta) = -{\sqrt{-1}\over
\hbar}\omega_{\Omega} (\xi,\eta) Id $. 

(ii) $\curv(\nabla^{\tilde{L}})(\xi,\eta) = -{\sqrt{-1}\over\hbar}
\omega_{\Omega} (\xi,\eta) Id $. 

(iii) The correspondence $f \mapsto \hat{f}$ is a procedure of quantization.
\end{prop}

\begin{rem} Having this procedure of quantization, we
obtain also the corresponding representation $\Lambda$ of the Lie algebra
${\frakt g}$ $$\vbox{\halign{ #&#&#\cr
 $\Lambda:$ &${\frakt g}$ &$ \rightarrow {\cal L}
({\cal H}_{\tilde{L}})$ \cr
  &$X$ &$ \mapsto \Lambda(X) = {\sqrt{-1}\over \hbar}\hat{f}_X, $ \cr}}$$
where $X\in{\frakt g}$ and $f_X \in C^{\infty}(\Omega)$ is the generating
function of the corresponding Hamiltonian field $\xi_X$ corresponding to
$X$. If $G$ is connected and simply connected, we obtain the
corresponding representation $T$ of $G$, $$T(\exp{X}) =
\exp{(\Lambda(X))}, \forall X \in {\frakt g}.$$
\end{rem}

It can be also prove that $$ \Lie_X \Ind(G;{\frakt P}, \sigma, \rho) =
\Lambda(X), \forall X \in {\frakt g}. $$

\subsection{U(1)-covering of the radicals and the semi-simple or 
reductive data }

We shall now reduce the procedure of multidimensional quantization on
$U(1)$-coverings of the stabilizers to the one with the semi-simple or
reductive data, following an idea of {\sc M. Duflo} and another one of {\sc P. L.
Robinson and J. H. Rawnsley}.

As usually, let us denote ${\frakt r}_F$ the solvable radical and
${\frakt u}_F$ the unipotent radical of the Lie algebra ${\frakt g}_F$ of
stabilizer $G_F$, respectively. Denote by ${\cal AP}^{U(1)}(G)$ the set
of all $F$ in ${\frakt g}^*$, which are $U(1)$-admissible and positive
well-polarizable, $X^{U(1)}_{irr}(F)$ the set of all equivalent classes
of {\it fundamental } irreducible unitary representations of $G_F^{U(1)}$
such that the restriction of each of them to $(G_F)_0^{U(1)}$ is a
multiple of the fundamental character $\chi_F^{U(1)}$. Then the elements
$\tilde{\sigma}\in X^{U(1)}_{irr}(F)$ are in a one-to-one
correspondence with the projective representations of the discrete group
$(G_F)_0\setminus G_F$. Recall that the set of the $U(1)$-data is defined
as $${\cal X}^{U(1)}(G) = \{ (F,\tilde{\sigma}) | F\in {\cal
AP}^{U(1)}(G),\tilde{\sigma} \in X^{U(1)}_{irr}(F) \}.$$

For every $U(1)$-data $(F,\tilde{\sigma})\in {\cal X}^{U(1)}(G)$, we can
choose a $(\tilde{\sigma},\chi_F^{U(1)})$-polarization $({\frakt
p},\rho,\sigma)$. Then we can construct the corresponding
 representation $\sigma$ of the $U(1)$-covering $H^{U(1)}:= G_F^{U(1)}
\ltimes H_0$ and finally induce to obtain a unitary representation of $G$.
We reduce here the consideration to the case of semi-simple or reductive 
groups. For this aim, let us denote $S_F$ the semi-simple component of
$G_F$ in its Cartan-Levi-Mal'tsev's decomposition., $G_F = S_F \ltimes R_F$.
Due to the fact that the principal $S_F$-bundle $S_F \rightarrowtail R_F 
\setminus G {\buildrel k \over \longrightarrow} G_F \setminus G$, is locally
trivial, there exists a connection ( i.e. a trivialization ) on it. Therefore
the tangent vector fields on the base $\Omega \approx G_F\setminus G$ can be 
uniquely lifted to the corresponding horizontal vector fields on $R_F\setminus 
G$. Then the Kirillov 2-form $\omega_{\Omega}$ of our orbit $\Omega \approx
G_F\setminus G$ can be lifted a non-degenerate closed $G$-invariant 2-form 
$\tilde{\omega}_{\Omega}$ on the horizontal component $T^h(R_F\setminus G)$
defined by the formula $$\tilde{\omega}_{\Omega}(f)(\tilde{X},\tilde{Y}):=
\omega_{\Omega}(F)(k_*\tilde{X},k_*\tilde{Y}),$$
for all the horizontally lifted vector fields $\tilde{X}$,$\tilde{Y}$, where $f\in
R_F\setminus G$, $F = k(f) \in \Omega$.

The symplectic group $\Sp(T^h_f(R_F\setminus G),\tilde{\omega}_{\Omega}(f))$
has also a $U(1)$-covering $\Mp^c(T^h_f(R_F\setminus G),\tilde{\omega}_{\Omega}
(f))$, $$1 \rightarrow U(1) \rightarrow \Mp^c(T^h_f(R_F\setminus G)) {\buildrel
k\over \longrightarrow} \Sp(T^h_f(R_F\setminus G)) \rightarrow 1.$$

\begin{prop}
There is a natural group homomorphism $j$ from the solvable radical
$R_F$ and a natural group homomorphism 
from the unipotent radical $U_F$ of 
$G_F$ to $\Sp(T^h_f(R_F\setminus G))$, and therefore a metaplectic
representation of $R_F$ and one of $U_F$.
\end{prop}
\begin{pf}
To prove it is enough to remark that due to lifting we have
 $$k_*:T^h_f(R_F\setminus G) {\buildrel \approx \over \longrightarrow } T_{k(f)}\Omega \approx {\frakt
g/g}_F $$ and every element $g\in R_F \subseteq G_F$ acts on the tangent space 
$T_{k(f)}\Omega \approx {\frakt g/g}_F$ by the Adjoint action $K(g) = 
\widetilde{\Ad g^{-1}}$.
The same argument is true for the unipotent radical $U_F$.
\end{pf}

\begin{cor}
There are natural $U(1)$-coverings $R_F^{U(1)}$ and $U_F^{U(1)}$ of $R_F$
and $U_F$,respectively.
\end{cor}
\begin{pf}
It's enough to use the well-known five-homomorphisms lemma
from homological algebra, $$R_F^{U(1)}:= \{ ((\lambda,\phi),g) |
\delta(\lambda,\phi) = j(g):= k_*^{-1} \widetilde{\Ad g^{-1}} k_*\}
,$$ where $\phi \in \Sp(T^h_f(R_F\setminus G))$, $\lambda \in {\Bbb
C}$; $$ |\lambda^2 \det{\{ {1\over 2}(\phi - \sqrt{-1} \circ \phi
\circ \sqrt{-1} )\}}| = 1.$$ Therefore, every element of $R_F^{U(1)}$
is of type $((\lambda,\widetilde{\Ad g^{-1}}),g)$.
\end{pf}

Return now to the situation of the $K$-orbit $\Omega$. We have a short
exact sequence of Lie groups $$1 \rightarrow U(1) \rightarrow R_F^{U(1)}
 \rightarrow R_F \rightarrow 1$$, and therefore a short exact sequence of
Lie algebras $$0 \rightarrow {\frakt u}(1) \rightarrow {\frakt r}_F^{U(1)}
:= \Lie R_F^{U(1)} \rightarrow {\frakt r}_F \rightarrow 0.$$ Thus as
vector space, ${\frakt r}_F^{U(1)} \cong {\frakt r}_F \oplus {\frakt
u}(1) $. Moreover in (\cite{robinsonrawnsley},\S5) it was shown that it is the direct sum of
Lie algebras.

\begin{defn}
 A point $F\in {\frakt g}^*$ is called
$(r,U(1))$- (for solvable radical case, or $(u,U(1))$- for unipotent
radical case ){\it admissible } \index{point!admissible -}
 iff there exists a unitary character
$$\chi_{F,k}^{U(1)}:R_F^{U(1)} \rightarrow {\Bbb S}^1 = {\Bbb T}^1
\subset {\Bbb C}^{\times},$$ with the differential
$$D\chi_{F,k}^{U(1)}(X,\varphi) = {\sqrt{-1}\over\hbar}(\langle F,X\rangle  +
k\varphi),$$ where $(X,\varphi)\in {\frakt r}_F^{U(1)}$, $k\in {\Bbb
Z}$ is some fixed integral number. If $k = 0$,
$(r,U(1))$-admissibility is coincided with the previously defined
$r$-admissibility and $(u,U(1))$-admissibility is coincided with
$u$-admissibility. If $k=1$, we say that the lifted character
$\chi_F^{U(1)} = \chi_{F,1}^{U(1)}$ is {\it fundamental.  } 
\index{character!lifted fundamental -}  \end{defn}

Denote by ${^r}X^{U(1)}_{irr}(F)$ ( resp., ${^u}X^{U(1)}_{irr}(F)$ )
the set of all equivalent classes of the so called fundamental
irreducible unitary representations of $G_F^{U(1)}$ such that the
restriction of each of them to $R_F^{U(1)}$ ( resp., $U_F^{U(1)}$ ) is
a multiple of the fundamental character $\chi_F^{U(1)}$.

\begin{defn}
A smooth complex tangent distribution $\tilde{L} \subset T(R_F\setminus 
G)_{\Bbb C}$ is called a {\it positive solvable (tangent) distribution } 
\index{distribution!positive solvable} iff:

(i) $\tilde{L}$ is an integrable and $G$-invariant sub-bundle of
$T^h(R_F\setminus G)_{\Bbb C}$,

(ii) $\tilde{L}$ is invariant under the action $\Ad$ of $G_F$,

(iii) For each $f \in R_F \setminus G$, the fiber $\tilde{L}_f$
is a positive polarization of the symplectic vector space 
$(T^h(R_F\setminus G)_{\Bbb C},\tilde{\omega}_{\Omega}(f))$, i.e.
$$\dim{\tilde{L}_f} = {1\over 2}\dim{T^h_f(R_F\setminus G)},\leqno{(\alpha)} $$
$$\tilde{\omega}_{\Omega}(f)(\tilde{X},\tilde{Y}) \equiv 0, \forall
\tilde{X},\tilde{Y} \in T^h_f(R_F\setminus G)_{\Bbb C}, \leqno{(\beta)} $$
$$\sqrt{-1}\tilde{\omega}_{\Omega}(f)(\tilde{X},\overline{\tilde{X}}) \geq 0,
\forall \tilde{X} \in T^h_f(R_F\setminus G)_{\Bbb C}. \leqno{(\gamma)}$$
We say that the distribution $\tilde{L}$ is {\it strictly positive } 
\index{distribution!strictly positive -}
iff the inequality $(\gamma)$ is strict for every non-zero $\tilde{X}\in
\tilde{L}_f$.
\end{defn}

It is easily to see that {\it if $\tilde{L}$ is a positive solvable
distribution, then the inverse image ${\frakt p}:= p^{-1}(L_F)$,
where $L_F:= k_*\tilde{L}_f$, for $F:= k(f)$, under the natural
projection $$p: {\frakt g}_{\Bbb C} \longrightarrow {\frakt g}_{\Bbb
C}/({\frakt g}_F)_{\Bbb C}$$ is a positive ( solvable ) polarization
in ${\frakt g}_{\Bbb C}$ at $F$. }

Suppose that $\tilde{L} \subset T^h(R_F\setminus G)_{\Bbb C}$ is a
positive solvable distribution such that $\tilde{L} \cap
\overline{\tilde{L}}$ and $\tilde{L} + \overline{\tilde{L}}$ are the
complexifications of some real distributions. In this case, the
corresponding complex sub-algebra ${\frakt p}:=
p^{-1}(k_*\tilde{L}_f)$ satisfies the following conditions: ${\frakt
p} \cap \overline{\frakt p}$ and ${\frakt p} + \overline{\frakt p}$
are the complexifications of the real Lie sub-algebras ${\frakt h:= p
\cap g}$ and ${\frakt m}:= (p + \overline{
\frakt p}) \cap {\frakt g}$. Denote by $H_0$ and $M_0$ the corresponding
analytic subgroups.

\begin{defn}
A positive (solvable) distribution $\tilde{L}$ is called {\it closed
} \index{distribution!closed -} iff the subgroups $H_0$, $M_0$ and the 
semi-direct product 
$H:= R_F \ltimes H_0$ and $M:= R_F \ltimes M_0$ are closed in $G$.
\end{defn}

{\it In what follows we assume that $\tilde{L}$ is closed. }
It is also easy to prove the following.

\begin{prop}
In a neighborhood of identity of $R_F^{U(1)}$, we have $$\chi_F^{U(1)}
((\lambda,\widetilde{\Ad g^{-1}}),g) = \exp{({\sqrt{-1}\over \hbar}(\langle F,X\rangle
+ \varphi)),}$$ where $\varphi \in {\Bbb R}$, $0 \leq \varphi < 2\pi\hbar$
such that $$\lambda^2 \det{(\widetilde{\Ad g^{-1}} - \sqrt{-1} \circ
\widetilde{\Ad g^{-1}} \circ \sqrt{-1})} = \exp{({\sqrt{-1}\over\hbar}\varphi)}
.$$ As an operator, acting on the Hilbert space ${\Bbb H}_{\tilde{V}}({\frakt
 m}_{\Bbb C}/{\frakt h}_{\Bbb C})$, it has the following integral $L^2$-
kernel $$u(z,w) = \exp{\{ {\sqrt{-1}\over\hbar}(\langle F,X\rangle  + \varphi) +
{\sqrt{-1} \over 2\hbar}\langle z,w\rangle  - {1\over 4\hbar}\langle w,w\rangle  \}},\forall z,w\in
{\frakt m}_{\Bbb C}/{\frakt h}_{\Bbb C}.$$ 
\end{prop}

According to the definitions, $H_0$ is a normal subgroup in $H$ and
there is the adjoint action of $R_F$ on $H_0$. Moreover, we have the
projection epimorphism $$\pi_j: R_F^{U(1)} \rightarrow R_F,$$ Thus
$R_F^{U(1)}$ acts on $H_0$. Hence we can define the semi-direct product
$R_F^{U(1)}\ltimes H_0$ and the following commutative diagram $$\vbox{\halign{
#&#&#&#\cr  $1 \rightarrow U$&$(1) \rightarrow R_F^{U(1)}$&$\ltimes H_0
\rightarrow R_F$&$\ltimes H_0 \rightarrow 1$ \cr  &${\downarrow}Id $ 
&${\downarrow}\pr_1 $ &${\downarrow}\pr_1 $ \cr $1 \rightarrow U$&$(1) 
\rightarrow\quad $&$R_F^{U(1)} \rightarrow $ &$R_F \rightarrow 1,$ \cr}}$$
where by $pr_1$ be denote the projection map on the first factor.
Then $H^{U(1)}:= R_F^{U(1)}\ltimes H_0$ is the $U(1)$-covering of $H
:= R_F \ltimes H_0$. It has ${\frakt h \oplus u}(1)$ the corresponding
Lie algebra. Denote as usually by $H_0^{U(1)}$ the inverse image of $H_0$
in $H^{U(1)}$ under the $U(1)$-covering projection map.

Let us say that a  fixed irreducible unitary representation $\tilde{
\sigma}$ of
$G_F^{U(1)}$ in a separable Hilbert space $\tilde{V}$ is {\it fundamental 
} \index{representation!fundamental -}
iff its restriction to $R_F^{U(1)}$ is a multiple of the fundamental character 
$\chi_F^{U(1)}$.

\begin{defn}
The triple $(\tilde{L},\rho,\sigma_0)$ is said to be a {\it solvable $(\tilde{
\sigma},\chi_F^{U(1)})$-polarization } \index{polarization!solvable -}
and $\tilde{L}$ is a {\it weakly Lagrangian distribution }
iff 

(i) $\sigma_0$ is an irreducible representation of the subgroup $H_0^{U(1)}$
in a Hilbert space $V = \tilde{V} \otimes V'$ such that $$\sigma_0|_{R_F^{U(1)}
\cap H_0^{U(1)}} = \tilde{\sigma},$$

(ii) $\rho$ is a representation of the complex Lie algebra ${\frakt p
\oplus u}(1)$ in $V$ such that $$\rho|_{{\frakt h \oplus u}(1)} = 
D\sigma_0, $$ the corresponding infinitesimal representation of the
Lie algebra.
\end{defn}

It is easy to prove also an analogue of the Proposition A2.10 as follows 

\begin{prop} 
If $\Omega_F$ is $(r,U(1))$-admissible, $(\tilde{L},\rho,\sigma_0)$ 
is a$(\tilde{\sigma},\chi_F^{U(1)})$-solvable polarization, then there
exists a unique irreducible unitary representation $\sigma$ of $H^{U(1)}$ 
in $V$ such that $$\sigma|_{R_F^{U(1)}} = \tilde{\sigma} $$ and $$D\sigma
= \rho|_{{\frakt h \oplus u}(1)}.$$
\end{prop}

It is reasonable to remark that the same can be done for the unipotent
radical in place of the solvable radical. We have an analogous
reduction of the procedure of quantization on the $U(1)$-coverings.

\subsection{Induction from semi-simple data }

Let $(F,\tilde{\sigma})$ is a semi-simple $U(1)$-data, i.e. $F$ is a
$(r, U(1))$-admissible, therefore there exists the fundamental
character $\chi_F^{ U(1)}$ of the $U(1)$-covering $R_F^{U(1)}$ of
$R_F$, and $\tilde{\sigma}$ is a fundamental representation of
$G_F^{U(1)}$ in the sense that its restriction to the $U(1)$-covering 
$R_F^{U(1)}$ of the radical $R_F$ of $G_F$ is a multiple of the fundamental
character $\chi_F^{U(1)}$.

Recall that for the element $F_0 = 0 \in {\frakt g}^*$, the
stabilizer is the whole group, $G_F = G$. It is reasonable to put
$\Mp^c({\frakt g/g}_{F_0},\tilde{\omega}_{F_0}) = U(1)$ and $G^{U(1)}
 = G \times U(1) $.  Fixing a connection $\overline{\Gamma}$ on $H
\rightarrowtail G \twoheadrightarrow H\setminus G $ we have the corresponding 
ones $\Gamma$, $\tilde{\Gamma}$ on $G_F \rightarrowtail G
\twoheadrightarrow G_F\setminus G$ and $R_F \rightarrowtail G
\twoheadrightarrow R_F\setminus G$, respectively. Passing to the
$U(1)$-covering of the leaves i.e. the structural group are lifted to
the corresponding $U(1)$-covering, we have the following isomorphisms
between the total spaces of principal bundles $$G^{U(1)} \cong
G_{\Gamma}^{ U(1)} \cong G_{\overline{\Gamma}}^{U(1)} \cong
G_{\tilde{\Gamma}}^{ U(1)},$$  we have a principal $H^{U(1)}$-bundle
on $H\setminus G$, a principal $G_F^{U(1)}$-bundle on the orbit
$\Omega$, a principal $R_F^{U(1)}$-bundle on on $R_F \setminus G$
and two homomorphisms between them as follows $$\vbox{\halign{
#&#&#&#&#&#&#&#&#&#&#\cr $H^{U(1)} \rightarrow
$&$G$&$_{\overline{\Gamma}}^{U(1)}$ & &$G_F^{U(1)} \rightarrow $
&$G$&$_{\Gamma}^{U(1)}$ & &$R_F^{U(1)} \rightarrow $
&$G$&$_{\tilde{\Gamma}}^{U(1)} $ \cr & &$\downarrow $ & & &
&$\downarrow $ & & & &$\downarrow $ \cr &$B$&$\setminus G$ &$
{\buildrel \pi \over \leftarrow}$ & &$G_F$&$\setminus G$ &$ {\buildrel
k\over\leftarrow}$ & &$R_F$&$\setminus G,$ \cr}}$$ where $\pi$ and $k$
are the natural projections.

With the fundamental representation $\sigma$ of $H^{U(1)}$ we can
construct the induced bundles $\overline{\cal E}_{V,\sigma,\rho}$, ${\cal
E}_{V,\sigma|_{G_F^{U(1)}},\rho|_.}$, $\tilde{\cal
E}_{V,\sigma|_{R_F^{U(1)}},\rho|_.}$ over the bases $H\setminus G$, $G_F
\setminus G$, $R_F\setminus G$, respectively. In the category of smooth
 vector bundles $k^* \pi^* \overline{\cal E} _{V,\sigma,\rho}$, $k^*{\cal
E} _{V,\sigma|_{G_F^{U(1)}},\rho|_.}$ and $\tilde{\cal
E}_{V,\sigma|_{R_F^{U(1)}},\rho|.}$ are equivalent. The fixed connection
$\overline{\Gamma}$ induces a connection $\Gamma$ on the principal bundle
$G_F^{U(1)} \rightarrow G^{U(1)} \rightarrow \Omega \approx G_F\setminus
G$, and a ( affine ) connection $\nabla^{\Gamma}$ on the associated
bundle ${\cal E} _{V,\sigma|_{G_F^{U(1)}},\rho|_.}$. By analogy, we
obtain also an affine connection $\widetilde{\nabla^{\Gamma}}$ on
$\tilde{\cal E} _{V,\sigma|_{R_F^{U(1)}},\rho|.}$. We put $${\cal
S}_{\tilde{L},S_F}(\tilde{\cal E} _{V,\sigma|_{R_F^{U(1)}},\rho|.}):= \{
s\in {\cal S} _{S_F}(\tilde{\cal E} _{V,\sigma|_{R_F^{U(1)}},\rho|.}) |
\widetilde{\nabla^{\Gamma}}_{\tilde{\xi}}s
 \equiv 0, \forall \tilde{\xi}\in \tilde{L} \},$$ where ${\cal
S}_{S_F}(\tilde{\cal E} _{V,\sigma|_{R_F^{U(1)}},\rho|.})$ is the vector
space of $S_F$-equivariant sections of the bundle $\tilde{\cal
E}_{V,\sigma|_{R_F^{U(1)}},\rho|.}$. The natural action of $G$ on this
space of sections is denoted by $\Ind(G;\tilde{L},\rho,\sigma_0)$ and is
called the {\it (reduction of) the partially invariant and holomorphly
induced representation. } It is easily to unitarize, and do, this
representation to the unitary one on the completion ${\cal H}_{\tilde{L}}$
of ${\cal S} _{S_F}(\tilde{\cal E} _{V,\sigma|_{R_F^{U(1)}},\rho|.})$

As above, for our convenience in what follows we write ${\frakt P}$
for the principal bundle $R_F^{U(1)} \rightarrowtail G^{U(1)}
\twoheadrightarrow R_F \setminus G$. By using the well-defined
homomorphism $$l: R_F^{U(1)} \longrightarrow  \Mp^c(T^h_f(R_F\setminus G))$$
and the complex metaplectic representation $$\mu: \Mp^c(T^h_f(R_F\setminus G))
\longrightarrow \End {\cal E}
'(T^h_f(R_F\setminus G)),$$ we have the bundle associated with the principal
bundle ${\frakt P}$ via the homomorphism $\mu\circ l$.
Consider now the tensor product representation $\sigma. \chi_F^{U(1)}$
of the fiber $R_F^{U(1)}$. Denote by ${\cal E}
'_{V,\sigma|_{R_F^{U(1)}},\rho|.}({\frakt p})$ the corresponding associated 
bundle.

\begin{thm} With every solvable
$(\tilde{\sigma},\chi_F)$-polarization $(\tilde{L}, \rho,\sigma_0)$,
there exists a natural unitary representation of $G$ in the Hilbert space
${\cal H}_{\tilde{L}}$.
\end{thm}
\begin{pf} 
According to Proposition A2.14, and since $k_*$ is the
linear lifting isomorphism, there exists for each $f'\in R_F\setminus
G$ a canonical linear map $$W_{f'}: T^h_{f'}(R_F\setminus G)_{\Bbb C}
\longrightarrow \End({\cal E}
'_{V,\sigma|_{R_F^{U(1)}},\rho|.}(p)_{f'}),$$ such that
$$[W_{f'}(\tilde{X}),W_{f'}(\tilde{Y})] = -{\sqrt{-1}\over\hbar} \tilde{
\omega}_{\Omega}(f')(\tilde{X},\tilde{Y}), \forall \tilde{X},\tilde{Y}\in
T^h_{f'}(R_F\setminus G)_{\Bbb C}.$$
We assume that $\dim_{\Bbb R}{T^h_f(R_F\setminus G)} = 2m$, the top exterior
power ${\cal K}
^{\tilde{L}}$ is a complex line bundle in $\wedge^mT^h_f(FF\setminus G)$
 with the basis vector $K_{\tilde{L}}$. Putting $${\cal E}
'_{V,\sigma|_{R_F^{U(1)}},\rho|.}({\frakt P})^{\tilde{L}}_{f'}:= \{
\varphi \in {\cal E} '_{V,\sigma|_{R_F^{U(1)}},\rho|.}({\frakt P})_{f'} |
W_{f'}(\tilde{X})\varphi = 0,\forall \tilde{X}\in\tilde{L}_{f'} \},$$
then $${\cal E} '_{V,\sigma|_{R_F^{U(1)}},\rho|.}({\frakt P})^{\tilde{L}}
:= \bigcup_{f'} {\cal E} '_{V,\sigma|_{R_F^{U(1)}},\rho|.}({\frakt
P})^{\tilde{L}}_{f'}$$ is the tensor product of the complex line bundle
${\cal E} '({\frakt P})^{\tilde{L}}$ with the associated bundle ${\cal
E}_{V,\sigma|_{R_F^{U(1)}},\rho|.}({\frakt P})$. As in Proposition A2.16
, there is a canonical isomorphism of complex bundles $${\cal E}
'_{V,\sigma|_{R_F^{U(1)}},\rho|.}({\frakt P})^{\tilde{L}} \otimes {\cal E}
'({\frakt P})^{\tilde{L}} \otimes {\cal K} ^{\tilde{L}} \quad{\buildrel
\cong\over \rightarrow}\quad {\cal
E}_{V,\sigma|_{R_F^{U(1)}},\rho|.}({\frakt P}) \quad.$$

Putting ${\cal Q} ({\frakt P})^{\tilde{L}}:= {\cal
E}'_{V,\sigma|_{R_F^{U(1)}},\rho|.}({\frakt P})^{\tilde{L}} \otimes {\cal
K} ^{\tilde{L}}$, we have as in the previous subsection $$[{\cal Q}
({\frakt P})^{\tilde{L}}]^{\otimes 2} = [{\cal
E}_{V,\sigma|_{R_F^{U(1)}},\rho|.}({\frakt P})]^{\otimes 2} \otimes {\cal
K} ^{\tilde{L}}.$$ We have also the affine connection
$\tilde{\nabla}^{\Gamma}$ on the associated bundle ${\cal
E}_{V,\sigma|_{R_F^{U(1)}},\rho|.}({\frakt P})$, the affine connection
$\tilde{ \nabla}^{\tilde{L}}$ on the associated bundle ${\cal Q} ({\frakt
P})^{\tilde{L}})$. Denote by ${\cal S}_{\tilde{L},S_F}({\cal Q} ({\frakt
P})^{\tilde{L}})$ the space of the $S_F$-equivariant sections of ${\cal Q}
({\frakt P})^{\tilde{L}}$ such that $$\tilde{\nabla}^{\tilde{L}}s \equiv 0
,\forall \xi\in\tilde{L}.$$ Define $${\cal H}_{\tilde{L}}:= \{ s\in
{\cal S}_{\tilde{L},S_F}({\cal Q} ({\frakt P})^{\tilde{L}}) |
\langle s,s\rangle _{\tilde{L}} \in C_0^{\infty}(G_F\setminus G) \}$$ and denote also by
${\cal H}_{\tilde{L}}$ the completion with respect to this norm. The
natural unitary representation of $G$ in the Hilbert space is the required
representation.
\end{pf}

\subsection{A reduction of the multidimensional quantization 
procedure on U(1)-covering }

We use the associated bundle $${\cal Q} ({\frakt P})^{\tilde{L}} = {\cal
E}_{V,\sigma|_{R_F^{U(1)}},\rho|.}({\frakt P}) \otimes {\cal E} '({\frakt
P})^{\tilde{L}} \otimes {\cal K} ^{\tilde{L}}$$ to realize our reduction of
the quantization procedure on $U(1)$-covering. Define $$\vbox{\halign{
#&#&#&#\cr $\hat{.}: C^{\infty}$&$(\Omega)$&$ \rightarrow $&${\cal L}
({\cal H} _{\tilde{L}}),$ \cr
   &$ f $&$ \mapsto $&$\hat{f} = f + {\hbar\over\sqrt{-1}}
\tilde{\nabla}^{ \tilde{L}}_{\tilde{\xi}_f}\quad,$ \cr}}$$ where ${\cal
L} ({\cal H}_{\tilde{L}})$ is the space of all ( perhaps unbounded )
Hermitian operators admitting operator closure on the space ${\cal
H}_{\tilde{L}}$ and $\tilde{ \nabla}_{\tilde{\xi}_f}^{\tilde{L}}$ is the
covariant derivation associated with the connection
$\tilde{\nabla}^{\tilde{L}}$ on the $G$-bundle ${\cal Q} ({\frakt
P})^{\tilde{L}}$. Recall that by definition, $$\tilde{\nabla}^{
\tilde{L}}_{\tilde{\xi}_f} = L_{\tilde{\xi}_f} + {\sqrt{-1}\over\hbar}
\alpha(\tilde{\xi}_f),$$ where ${\sqrt{-1}\over\hbar}\alpha$ is the form
of connection and $L_{\tilde{\xi}_f}$ is the Lie derivation along $\tilde{
\xi}_f$ which is the horizontal lift of the strictly Hamiltonian vector
field $\xi_f$ corresponding to $f$.

It is easy also to check the following results 

\begin{prop}
The following three conditions are equivalent:

(i) $d_{\tilde{\nabla}^{\tilde{L}}} \alpha(\tilde{\xi},\tilde{\eta}) = 
-\tilde{\omega}_{\Omega}(\tilde{\xi},\tilde{\eta}) Id $.

(ii) $\curv \tilde{\nabla}^{\tilde{L}}(\tilde{\xi},\tilde{\eta}) =
-{\sqrt{-1}\over\hbar}\tilde{\omega}_{\Omega}(\tilde{\xi},\tilde{\eta}) Id $.
v
(iii) The application $f \mapsto \hat{f}$ is a procedure of multidimensional
 quantization.
\end{prop}

\begin{thm}
If one of the previous condition holds, with every solvable  $(\tilde{
\sigma},\chi_F^{U(1)})$-polarization $(\tilde{L},\rho,\sigma_0)$ there is
 a natural representation, denoted $\Ind(G;\tilde{L},\rho,\sigma_0)$ of $G$
in the Hilbert space of quantum states ${\cal H}
_{\tilde{L}}$, which is the completion of the space of of partially
invariant partially holomorphic $S_F$-invariant sections of the induced 
bundle $\tilde{\cal E}
_{V,\sigma|_{R_F^{U(1)}},\rho|.}({\frakt P})$ with respect to the indicated 
scalar product $\langle .,.\rangle _{\tilde{L}}$ of sections. The Lie derivation of this representation is just the representation of Lie algebra ${\frakt g}$ associated with the reduction of the reduction of quantization procedure on $U(1)$-covering, $$
\Lie_X \Ind(G;
\tilde{L},\rho,\sigma_0) = \Lambda(X):= {\sqrt{-1}\over\hbar}
\hat{f}_X,\forall X\in {\frakt g}.$$
\end{thm}

\section{Globalization  over U(1)-Coverings }

In the previous sections, we have realized representations in section
spaces of the induced bundles. It is well-known that this is not
enough in some situations, say the discrete series representations of
semi-simple Lie groups is better to be realized in $L^2$-cohomologies.
W. Schmid and J.A. Wolf \cite{schmidwolf} proposed an algebraic model of the
geometric multidimensional quantization to realize the discrete series
representations of semi-simple Lie groups, using the ${\Bbb Z}/2{\Bbb
Z}$-covering of the stabilizers of $K$-orbits. Using the
$U(1)$-covering, we can also give an algebraic version of the
multidimensional quantization procedure. In this section we expose a
revised version of the work \cite{dong1},\cite{dong2}.

\subsection{Classical constructions and three geometric complexes }

Let $G$ be a connected linear semi-simple Lie group, ${\frakt g} =
\Lie G$ its Lie algebra and ${\frakt g}^*$ its dual space. Recall that
$G$ acts on ${\frakt g}$ via the adjoint representation $\Ad$, and on
${\frakt g}^*$ by the co-adjoint ( i.e. contragradient ) representation
. Let $F\in {\frakt g}^*$ and $G_F$ be the stabilizer of this point.
Suppose that ${\frakt h =g}_F$ is a {\it Cartan sub-algebra }. 
. This is the case, for example for the discrete series representations.
In this case the functional $F\in {\frakt g}^*$ is admissible and
well-polarizable. Then the stabilizer group $H = G_F$ at this point
$F$ is a Cartan subgroup of $G$. It is reasonable to suppose that it is
$U(1)$-admissible, i.e. there exists a unitary character, said to be
fundamental $\chi_F^{U(1)}$, such that $$D\chi_F^{U(1)} = {\sqrt{-1}\over
\hbar}(\langle F,X\rangle  + \varphi), \forall (X,\varphi) \in {\frakt h\oplus u}(1)_{
\Bbb C}.$$

Let {\frakt b} be a closed positive polarization in ${\frakt g}_{\Bbb
C}$. It is well-known that ${\frakt b}$ is therefore a Borel
sub-algebra of ${\frakt g}_{\Bbb C}$, containing the Cartan sub-algebra
${\frakt h}$. Let $\tilde{\sigma}$ be some fixed fundamental
irreducible unitary representation of $H^{U(1)}=G_F^{U(1)}$ such that
its restriction to $H_0^{U(1)}$ is a multiple of the fundamental character
$\chi_F^{U(1)}$. We can then consider the associated bundle $${\cal E}
_{\tilde{V},\tilde{\sigma}}:= G_{\Gamma}^{U(1)} \times_{H^{U(1)},\tilde{
\sigma}} \tilde{V} \twoheadrightarrow H\setminus G \approx \Omega_F.$$

Suppose that $\dim_{\Bbb C}{\Omega_F} = m$. Recall that in the orbit
$\Omega_F$ there is a natural complex structure related with the
polarization ${\frakt b}$.Let ${\cal C} ^q({\cal
E}_{\tilde{V},\tilde{\sigma}})$ denote the sheaf of differential forms of
type $(o,q)$ on $\Omega_F \approx H\setminus G$ with coefficients in the
induced bundle ${\cal E} _{\tilde{V},\tilde{\sigma}}$. Each differential
form of this type is a section of the bundle ${\cal
E}_{\tilde{V},\tilde{\sigma}} \otimes \wedge^q {\cal N} ^*$ over the base
$H\setminus G$, where by definition ${\cal N} \twoheadrightarrow
H\setminus G$ is the homogeneous vector bundle with the fiber ${\frakt n
\cong b/h}$ and ${\cal N} ^*$ is its dual. Denote by ${\cal O} ({\cal
E}_{\tilde{V},\tilde{\sigma}})$ the sheaf of germs of partially
holomorphic $C^{\infty}$ sections of ${\cal E}_{\tilde{V},\tilde{\sigma}}$
that are annihilated by the action of ${\frakt n
 \oplus u}(1)_{\Bbb C}$. Then we have a sequence of sheaves $$0
\rightarrow {\cal O}_{\frakt n}({\cal E}_{\tilde{V},\tilde{\sigma}})
{\buildrel i\over \rightarrow} {\cal C} ^0({\cal E}
_{\tilde{V},\tilde{\sigma}}) {\buildrel \overline{\partial}_{{\cal E}
_{\tilde{V},\tilde{\sigma}}} \over \longrightarrow} \dots {\buildrel
\overline{\partial}_{{\cal E} _{\tilde{V},\tilde{\sigma}}}\over
\longrightarrow} {\cal C} ^m({\cal E} _{\tilde{V},\tilde{\sigma}})
\rightarrow 0,$$ where the map $i$ is induced by the inclusion of the
space of partially invariant and partially holomorphic sections of\ ${\cal
E} _{\tilde{V},\tilde{\sigma}}$ into the space of smooth sections, and
the maps $\overline{\partial}_{{\cal E} _{\tilde{V},\tilde{\sigma}}}$ are
induced by the usual differential operator, mapping a $(0,q)$-form into a
$(0,q+1)$-form.
Therefore we have the corresponding sequence of global section spaces $$0
\rightarrow C^{\infty}(H \setminus G;{\cal O}_{\frakt n}({\cal
E}_{\tilde{V},\tilde{\sigma}})) \rightarrow C^{\infty}(H\setminus G; {\cal
C} ^0({\cal E} _{\tilde{V},\tilde{\sigma}})) \rightarrow\dots $$
$$\dots \rightarrow C^{\infty}(H\setminus G; {\cal C} ^m({\cal
E}_{\tilde{V},\tilde{\sigma}})) \rightarrow 0 $$
 and this sequence of abelian groups form a cochain complex, denoted by
$(C^{\infty}(H\setminus G;{\cal E}_{\tilde{V}.\tilde{\sigma}} \otimes
\wedge^* {\cal N} ^*),\overline{\partial}_{{\cal
E}_{\tilde{V},\tilde{\sigma}}})$. For the cohomology groups of this
cochain complex it is easy to prove the following analogue of the
well-known Dolbault's theorem

\begin{prop} {\it There are canonical isomorphisms of between this type
cohomology groups $$H^p(C^{\infty}(H\setminus G;{\cal
E}_{\tilde{V},\tilde{\sigma}} \otimes \wedge^* {\cal N} ^*)) \cong
H^p(H\setminus G;{\cal O} _{\frakt n}({\cal
E}_{\tilde{V},\tilde{\sigma}})), \forall p \geq 0,$$ where the right
handside is the sheaf cohomology group of the space $H\setminus G$ of
degree $p$ with coefficients in ${\cal O} _{\frakt n}({\cal
E}_{\tilde{V},\tilde{\sigma}})$. }
\end{prop}

We refer to this complex as {\it the first. } \index{complex!the first -}

Let us now consider {\it the second complex. } \index{complex!the 
second -} It is easy to see that the
differential $\overline{\partial}_{{\cal E}_{\tilde{V},\tilde{\sigma}}}$
can be extended to the {\it hyperfunctions } section \index{hypersection} of 
the sheaves. So
we have the second complex $(C^{-\omega}(H \setminus G;{\cal E}
_{\tilde{V},\tilde{\sigma}} \otimes \wedge^* {\cal N}
^*),\overline{\partial}_{{\cal E} _{\tilde{V},\tilde{\sigma}}}) $. Let us
denote by $X$ the {\it flag variety } of the Borel sub-algebras of 
${\frakt
g}_{\Bbb C}$. It is a well-defined complex manifold. The $G$-orbit
$G.{\frakt b}$ of ${\frakt b}$ is a $G$-invariant analytic submanifold of
the complex manifold $X$. Therefore $S$ has also the structure of $CR$
manifold ( i.e. Cauchy-Riemann structure ). Since the subgroup $H=G_F$
normalizes ${\frakt b}$, we have a natural $G$-invariant fibration
$$H\setminus G \twoheadrightarrow S:= G.{\frakt b} \approx B\setminus
G\subseteq X.$$

Recall that there is a unique fundamental irreducible representation
$(\sigma,V)$ of the $U(1)$-covering $B^{U(1)} = H^{U(1)} \ltimes B_0$ of
$B:= H \ltimes B_0$ such that $\sigma|_{H^{U(1)}} = \tilde{\sigma}$ Then
the bundle ${\cal E}_{\tilde{V},\tilde{\sigma}} \twoheadrightarrow
H\setminus G$ can be considered as the push up of the bundle ${\cal
E}_{V,\sigma} \twoheadrightarrow S \approx B\setminus G$, or equivalently
the first one push down to the section one. we obtain as in
\cite{sulankewintgen} the {\it
Cauchy-Riemann complex } $(C^{-\omega}(S;{\cal E} _{V,\sigma} \otimes
\wedge^* {\cal N} _S), \overline{\partial}_S)$, where ${\cal N} _S:=
T^{0,1}(S)$, the second component of the decomposition of the tangent
bundle into holomorphic and antiholomorphic parts, with respect to the
complex structure defined by ${\frakt n}/({\frakt n}\cap \overline{\frakt
n})$

Denote by $X^{U(1)}$ the flag variety of $U(1)$-invariant Borel
sub-algebras of ${\frakt g \oplus u}(1)_{\Bbb C}$. We obtain the natural
projection $$\pi_X: X^{U(1)} \rightarrow X.$$ By using the Gauss
decomposition $G = K.B$, for some maximal compact subgroup of $G$, we
have $B\setminus G \cong B^{U(1)} \setminus K.B^{U(1)}$. Note also that
$K.B^{U(1)}$ acts on the flag variety $X^{U(1)}$. Let us denote
$$S^{U(1)} = (K.B^{U(1)}).({\frakt b \oplus u}(1)_{\Bbb C}) \approx
B^{U(1)} \setminus K.B^{U(1)}$$ the orbit passing through ${\frakt b
\oplus u}(1)_{ \Bbb C}$ in $X^{U(1)}$. Then there exists a diffeomorphism
of $S^{U(1)}$ onto $S$. By projection $\pi_X: X^{U(1)}
\twoheadrightarrow X$, we can make the homogeneous induced bundle
$$\pi_X^*{\cal E}_{V,\sigma} \twoheadrightarrow S^{U(1)},$$ and we have
the complex $$(C^{ -\omega}(S^{U(1)};\pi^*_X{\cal E} _{V,\sigma} \otimes
\wedge^* {\cal N} ^*_{S^{U(1)}}), \overline{\partial}_{S^{U(1)}}),$$
where by definition, ${\cal N} _{S^{U(1)}}:= \pi^*_X {\cal N} _S$, and
$\overline{\partial}_{S^{U(1)}}$ is induced from the $CR$ operator
$\overline{\partial}_S$.

The bundle $\pi^*_X{\cal E}_{V,\sigma} \otimes \wedge^p {\cal N}
^*_{S^{U(1)}} \longrightarrow S^{U(1)} \approx B\setminus G$ pull back to
the trivial bundle on $G$, so we have an isomorphism of the complexes
$$(C^{-\omega}(S^{U(1)};\pi^*_X {\cal E} _{V,\sigma} \otimes \wedge^*
{\cal N} ^*_{D^{U(1)}}), \overline{\partial}_{S^{U(1)}}) \cong $$
$$ \cong (\{
C^{-\omega}(G) \otimes V \otimes \wedge^* ({\frakt n}/{\frakt n} \cap
\overline{\frakt n} )^*\}^{{\frakt n}/{\frakt n}\cap{\frakt n},B^{U(1)}}
,\partial_{{\frakt n, n} \cap \overline{\frakt n}}),$$ for relative
$B^{U(1)}$-equivariant cohomology of the pair $({\frakt n}/{\frakt n}\cap
\overline{\frakt n})$ with the hyperfunctions coefficients.

A section $\tilde{s}\in C^{-\omega}(S^{U(1)};\pi^*_X{\cal E}_{V,\sigma}
\otimes \wedge^p{\cal N} _{S^{U(1)}})$ is said to be $H^{U(1)}$-{\it
invariant } iff $$\tilde{s}(hx) = \tilde{\sigma}(h)\tilde{s}(x),\forall
h\in H^{U(1)},\forall x\in S^{U(1)} \cong B\setminus G.$$ Let us denote
by $C^{-\omega}_{ H^{U(1)}}(S^{U(1)};\pi^*_X{\cal E} _{V,\sigma} \otimes
\wedge^p{\cal N} _{S^{U(1)}})$ the space of $H^{U(1)}$-invariant partially
holomorphic $C^{ \infty}$ sections of $C^{-\omega}(S^{U(1)};\pi^*_X{\cal
E} _{V,\sigma}\otimes \wedge^p{\cal N} _{S^{U(1)}})$, we have a canonical
isomorphism between the vector spaces
$$C^{\infty}_{H^{U(1)}}(S^{U(1)};\pi^*_X{\cal E} _{V,\sigma} \otimes
\wedge^p{\cal N} _{S^{U(1)}}) \cong C^{-\omega}(S;{\cal
E}_{V,\sigma}\otimes \wedge^p{\cal N} _S^*).$$

We see that the fibration $H\setminus G \rightarrow S^{U(1)}$ has
Euclidean space fibers. By applying the well-known Poincar\'e Lemma to
those fibers, we see that the inclusion of the complex $(\{
C^{-\omega}(G)\otimes V\otimes \wedge^*({\frakt n/n} \cap \overline{\frakt
n})^*\}^{{\frakt n}\cap\overline{\frakt n},B^{U(1)}}, \partial_{{\frakt
n,n}\cap \overline{\frakt n}})$ into the complex $(C^{-\omega}(H\setminus
G;{\cal E}_{V,\sigma} \otimes \wedge^*{\cal N} ^*),\partial_{{\cal
E}_{V,\sigma}})$ induces an isomorphism of the corresponding cohomology
theories The following result is therefore proved

\begin{prop} There are canonical isomorphisms between the cohomology
theories  $$H^p(C^{-\omega}(H\setminus G;{\cal E} _{V,\sigma} \otimes
\wedge^.{\cal N} ^*)) \cong H^p(C^{-\omega}_{H^{U(1)}}(S^{U(1)};\pi^*{\cal
E} _{V,\sigma} \otimes\wedge^.{\cal N} ^*_{S^{U(1)}}))$$ $$ \cong H^p(\{
C^{-\omega}(G) \otimes \wedge^.({\frakt n/n}\cap \overline{\frakt
n})^*\}^{{\frakt n}\cap\overline{\frakt n},B^{U(1)}}).$$
\end{prop}

Now we consider {\it the third complex. } \index{complex!the third -} Let us 
denote by 
$\tilde{S}$ the germ of $S$ in $X$. Then the bundle ${\cal E} _{V,\sigma}
\twoheadrightarrow S$ has a unique holomorphic ${\frakt g}$-equivalent
extension $\tilde{\cal E} _{V,\sigma} \twoheadrightarrow \tilde{S} \subset
X$ and we obtain {\it an analogue of the Dolbeault complex } 
\index{complex!analogue of the Dolbeault -}
$(C^{-\omega}(\tilde{S};\tilde{\cal E} _{V,\sigma}\otimes
T^{0,1}_X),\overline{\partial})$, with coefficients that are
hyperfunctions on $\tilde{S}$ with support in $S$.

Similarly, we have also the vector bundle $\pi^*_X\tilde{\cal E}
_{V,\sigma} \rightarrow \tilde{S}^{U(1)} \subset X^{U(1)}$ over the
$U(1)$-covering and then obtain the complex $(C^{-\omega}(\tilde{
S}^{U(1)};\pi^*_X\tilde{\cal E} _{V,\sigma} \otimes
T^{0,1}_{X^{U(1)}}),\overline{\partial}_{U(1)})$. By using the canonical
isomorphism $$C^{-\omega}_{H^{U(1)}}(\tilde{S}^{U(1)}; \pi^*_X\tilde{\cal
E} _{V,\sigma} \otimes \wedge^.T^{0,1}_{X^{U(1)}}) \cong
C^{-\omega}(S;\tilde{ \cal E} _{V,\sigma} \otimes T^{0,1}_X),$$ we havec
also the following result about cohomology of the third complex

\begin{prop}  There is a natural isomorphism of cohomology groups
$$H^p(C^{-\omega}_{ H^{U(1)}}(\tilde{S}^{U(1)};\pi^*_X\tilde{\cal E}
_{V,\sigma} \otimes \wedge^* T^{0,1}_{X^{U(1)}})) \cong
H^p_S(\tilde{S};{\cal O} (\tilde{\cal E} _{V,\sigma})),$$ where the right
hand side is the well-known local cohomology along $S$. 
\end{prop}

\subsection{Isomorphisms of cohomologies  }

We fix a basic datum $(H,{\frakt b},\tilde{\sigma})$ as in the
previous subsection. Denote by $Y$ the variety of ordered Cartan
sub-algebras. As homogeneous $G_{\Bbb C}$-space, we have $Y \approx
H_{\Bbb C}\setminus G_{\Bbb C}$, where $G_{\Bbb C}$ is the adjoint
group of ${\frakt g}_{\Bbb C}$, and $H_{\Bbb C}$ is the connected
subgroup corresponding to the Lie algebra ${\frakt h}_{\Bbb C}$.
Since $H_{\Bbb C}$ normalizes ${\frakt b}$, there is a natural
projection $$ p: Y \longrightarrow X$$ with fiber $p^{-1}({\frakt b})
= \exp{\frakt n}$. Let $S_Y = G.{\frakt h} \subset Y$ be the
$G$-orbit of the base point in $Y$, we have $$p: S_Y \longrightarrow
S$$ with fibers of type $\exp{({\frakt n\cap g})} = \exp{({\frakt n}\cap 
\overline{\frakt n}\cap {\frakt g})}$. Remark that by definition ${\frakt n}
= {\frakt b/h} \hookrightarrow {\frakt g}$ for some fixed inclusion,
related with the connection on the principal bundle. Then $S_Y$ is a
real form of the complex manifold $Y$ and $u:= \codim_{\Bbb R}(S)$ is
the dimension of fibers of the projection $p: S_Y \rightarrow S$.

Pushing up to the $U(1)$-covering of the map $p: Y\rightarrow X$, we
have the following commutative diagram
$$\vbox{\halign{ #&#&#\cr $Y^{U(1)}$ &${\buildrel p^{U(1)}\over \longrightarrow
}$ &$X^{U(1)} $ \cr ${\downarrow}\pi_X$ &$ \circlearrowleft $ &${\downarrow}
\pi_X $ \cr $Y $ &${\buildrel p \over \longrightarrow}$ &$X,$ \cr}}$$
where $\pi_X$ and $\pi_Y$ are the natural projections. Let us denote by $
T_{Y|X}$ the complexified relative tangent bundle of the fibration $ p: Y
\rightarrow X$ and $T_{Y|X} = T^{0,1}_{Y|X} \oplus T^{0,1}_{Y|X}$ its 
decomposition into the subbundles of holomorphic and antiholomorphic 
directions. We have also an $G_{\Bbb C}$-invariant isomorphisms $$p^*T_X
\oplus T_{Y|X} \cong T_Y,$$ which is compatible with the comple structure
and the Lie bracket.

We obtain also the complex $$(C^{-\omega}(S_Y;p^*\tilde{\cal
E}_{V,\sigma}\otimes (T^{1,0}_{Y|X})^*),\overline{\partial}).$$ Let
$S_Y^{U(1)} = (K.B^{U(1)}).({\frakt h \oplus u}(1)_{\Bbb C}) \subset
Y^{U(1)}$ be the orbit passing through ${\frakt h \oplus u}(1)_{\Bbb C}$.
We see that $S_Y^{U(1)}\approx S_Y$. By using the previous diagram, we
have the complex $$(C^{-\omega}(S_Y^{U(1)};(p^{U(1)})^*\pi_X^*\tilde{\cal
E} _{V,\sigma}\otimes \wedge^*
((T_{Y|X}^{1,0})^{U(1)})^*),\overline{\partial}_{ U(1)}),$$ where by
definition, $(T^{1,0}_{Y|X})^{U(1)} = (\pi_X)_*T^{1,0}_{Y| X}$, and
$\overline{\partial}_{U(1)} = (\pi_X)_*\partial$.

We see that $T_{Y|X}$ is a bundle with typical fiber ${\frakt n = b \oplus
u}(1)_{\Bbb C}/{\frakt h \oplus u}(1)_{\Bbb C}$ and $(\pi_X)_*\tilde{\cal
E}_{V,\sigma}$ is also a bundle with $H^{U(1)}$-module $V$ as the typical
fiber. It is not
 hard to see that the complex
$$(C^{-\omega}_{H^{U(1)}}(S^{U(1)};(p^{U(1)})^*( \pi_X)^*\tilde{\cal
E}_{V,\sigma} \otimes \wedge^*
((T^{1,0}_{Y|X})^{U(1)})^*),\overline{\partial}_{U(1)})$$ coincides with
the complex $$(C^{-\omega}(H\setminus G;{\cal E} _{V,\sigma}\otimes
\wedge^*{\cal n} ^*),\overline{\partial}_{{\cal E} _{V,\sigma}}).$$ We
obtain therefore

\begin{prop} There are canonical isomorphisms between algebraic
$G$-modules $$\begin{array}{rl} H^p(C^{-\omega}(H\setminus G;{\cal E} 
_{V,\sigma} \otimes \wedge^* {\cal N} ^*)) &\cong
H^p(C^{-\omega}_{H^{U(1)}}(S^{U(1)};\pi_V^*{\cal E} _{V,\sigma}\otimes
\wedge^*{\cal N} ^*_{S^{U(1)}})) \\   &\cong H^{p+u}_S(\tilde{S};{\cal O}
(\tilde{\cal E} _{V,\sigma})).\end{array} $$ \end{prop}

We fix a Cartan involution $\theta$ of $G$ such that $\theta K = K$.
Then the cartan subgroup $H = G_F$ can be decomposed into the direct
product $ H = T \times A$, such that the corresponding Lie algebra
decomposition ${\frakt h = t \oplus a}$ is the root decomposition of
${\frakt h}$ into the $\pm 1$-eigenspaces of the corresponding
endomorphism $\theta|_{\frakt h}$. Put $A:= \exp{({\frakt a \cap
g})}$. Consider the orbit $S = G,{\frakt b}\subset X$, where ${\frakt b}$
is as usually a fixed Borel sub-algebra, containing ${\frakt h}$.

{\it We can suppose for instance that ${\frakt b = b}_{max}$ is maximally
real }. This condition will be removed by the so called {\it change of
polarizations }. \index{change of polarizations} We have then $S_{max} = 
G.{\frakt b}_{max}$, ${\frakt h} \subset {\frakt b}_{max}$.

Consider a parabolic subgroup $P = MA.N_H:= (M \times A)\ltimes N_H$,
where $\theta M = M $, i.e. $M$ is a maximal compact semi-simple subgroup
of $P$ and ${\frakt b = b}_{max} \subset {\frakt p}:= \Lie P$. The
fibrations $S \twoheadrightarrow S_{max}$ and $S_{max} \twoheadrightarrow
P\setminus G$ induce a fibration $S \twoheadrightarrow P\setminus G$.
Since $S^{U(1)} \approx S$, we obtain therefore a fiberation $S^{U(1)}
\twoheadrightarrow P\setminus G$. Let ${\cal C} ^{-\omega}_{P\setminus
G}(S^{U(1)})$ be the sheaf of germs of hyperfunctions on $S^{U(1)}$, that
are $C^{\infty}$ along the fibers of $S^{U(1)} \twoheadrightarrow
P\setminus G$. Then ${\cal C} ^{-\omega}_{P\setminus G}(S^{U(1)})$
defines a complex of sheaves ${\cal C} ^{-\omega}_{P\setminus
G}(S^{U(1)};\pi_X^*{\cal E}_{V,\sigma} \otimes \wedge^p{\cal
N}^*_{S^{U(1)}})$ of germs of $H^{U(1)}$-equivariant sections of the
bundle $\pi^*_X{\cal E} _{V,\sigma}\otimes \wedge^p {\cal N} ^*_{S^{U(1)}}
\twoheadrightarrow S^{U(1)}$ with coefficients in $C^{-\omega}_{P\setminus
G}(S^{U(1)})$. Taking global sections, we have a subcomplex
$(C^{-\omega}_{P\setminus G}(S^{U(1)};\pi_X^*{\cal E} _{V,\sigma}\otimes
\wedge^*{\cal N} ^*_{S^{U(1)}}),\overline{\partial}_{S^{U(1)}})$ of
$(C^{-\omega}_{H^{U(1)}}(S^{U(1)};\pi^*_X{\cal E} _{V,\sigma} \otimes
\wedge^* {\cal N} ^*_{S^{U(1)}}), \overline{\partial}_{S^{U(1)}})$.

\begin{prop} The inclusion $$(C^{-\omega}_{P\setminus
G}(S^{U(1)};\pi_X^*{\cal E} _{V,\sigma}\otimes \wedge^*{\cal N}
^*_{S^{U(1)}}),\overline{\partial}_{S^{U(1)}}) \hookrightarrow
(C^{-\omega}_{H^{U(1)}}(S^{U(1)};\pi^*_X{\cal E} _{V,\sigma} \otimes
\wedge^* {\cal N} ^*_{S^{U(1)}}), \overline{\partial}_{S^{U(1)}})$$
induces isomorphisms of cohomology groups.
\end{prop} 
\begin{pf}
Applying the usual Dolbeault Lemma and the standard argument on
hyperfunctions, we see that the sheaves ${\cal C} ^{-\omega}_{P\setminus
G}(S^{U(1)};\pi_X^*{\cal E} _{V,\sigma}\otimes \wedge^*{\cal N}
^*_{S^{U(1)}})$ and ${\cal C} ^{-\omega}_{H^{U(1)}}(S^{U(1)};\pi^*_X{\cal
E} _{V,\sigma} \otimes \wedge^* {\cal N} ^*_{S^{U(1)}})$ are soft and the
inclusion of ${\cal C} ^{-\omega}_{P\setminus G}(S^{U(1)};\pi_X^*{\cal
E}_{V,\sigma}\otimes \wedge^*{\cal N} ^*_{S^{U(1)}})$ into ${\cal C}
^{-\omega}_{H^{U(1)}}(S^{U(1)};\pi^*_X{\cal E} _{V,\sigma} \otimes
\wedge^* {\cal N} ^*_{S^{U(1)}})$ induces isomorphisms of cohomology
sheaves. 

On the other hand, it follows easily that the inclusion of sheaves induces
an isomorphism of hyperfunction coefficient cohomology. Since both complexes
consist of soft sheaves, the hyperfunction coefficient cohomology is just
 the cohomology of the associated complexes of global sections. The
proposition is therefore proved.
\end{pf}

Remark that the theory of hyperfunctions with values in a reflexive Banach
space is developed exactly in the same way as the one for complex valued
hyperfunctions. By a similar argument as (\cite{schmidwolf},\S7) we obtain the
following result.

\begin{prop}  The vector spaces $C^{-\omega}_{P\setminus
G}(S^{U(1)};\pi_X^*{\cal E} _{V,\sigma} \otimes \wedge^p{\cal
N}_{S^{U(1)}})$ have natural Fr\'echet topologies. On those topologies
$\overline{\partial}_{S^{U(1)}}$ is continuous and the natural actions of
$G$ are Fr\'echet representations. 
\end{prop}

\subsection{Maximal real polarizations and change of 
polarizations }

Recall now some notions from \cite{robinsonrawnsley}: 

\begin{defn} An admissible Fr\'echet $G$-module is said {\it to have
property (MG) } \index{property (MG)}iff it is the maximal globalization of 
its underlying
Harish-Chandra module. {\it A complex of Fr\'echet $G$-modules has
property (MG) } iff its differential $d$ has closed range and each
cohomology group $H^p(C^., dB)$ is admissible, of finite length as
$G$-module, and has property (MG).
\end{defn}

Given a basic datum $(H,{\frakt b},\tilde{\sigma})$, we say that {\it the
corresponding homogeneous vector bundle ${\cal E} _{V,\sigma}
\twoheadrightarrow S^{U(1)}$ has the property (MG) } iff the associated
partially smooth Cauchy-Riemann complex $(C^{-\omega}_{P\setminus G}(
S^{U(1)};\pi^*_X{\cal E} _{V,\sigma} \otimes
\wedge^*{\cal N}_{S^{U(1)}}),\overline{\partial}_{S^{U(1)}})$ has property
(MG).

Let us denote $H^p(S^{U(1)};{\cal E} _{V,\sigma}):=
H^p(C^{-\omega}_{H^{U(1)}}(S^{U(1)};\pi_X^*{\cal E} _{V,\sigma}\otimes
\wedge^*{\cal N} ^*_{S^{U(1)}}))$. Proposition A3.5 shows that
$H^p(S^{U(1)};{\cal E} _{V,\sigma})$ is calculated by a Fr\'echet complex
, then we can consider the Fr\'echet subcomplex $H^p(S^{U(1)};{\cal
E}_{V,\sigma})_{(K)}$ of $K$-finite forms in this Fr\'echet complex. In
particular, we can define morphisms $$H^p(S^{U(1)};{\cal
E}_{V,\sigma})_{(K)} \longrightarrow A^p(G,H,{\frakt b},\tilde{\sigma})
,$$ where $$A^p(G,H,{\frakt b},\tilde{\sigma}) \cong
H^p(C^{for}(H\setminus G; {\cal E} _{V,\sigma}\otimes \wedge^*{\cal N}
^*_{S^{U(1)}})_{(K)})$$ are the well-known Harish-Chandra modules for $G$
.  Then $H^p(S^{U(1)};{\cal E} _{V,\sigma})$ {\it will be the
globalization \index{globalization} of $A^p(G,H,{\frakt b},\tilde{ 
\sigma})$ if the homomorphism
$$H^p(S^{U(1)};{\cal E}_{V,\sigma})_{(K)} \longrightarrow A^p(G,H,{\frakt
b},\tilde{\sigma}),$$ are really isomorphisms. }

Recall another notion from \cite{robinsonrawnsley}: 

\begin{defn} The bundle ${\cal E} _{V,\sigma} \twoheadrightarrow S^{U(1)}$
is said {\it to have property (Z) } \index{property (Z)} if the homomorphisms
$$H^p(S^{U(1)};{\cal E} _{V,\sigma})_{(K)} \longrightarrow A^p(G,H,{\frakt
b},\tilde{\sigma}),$$ are isomorphisms, i.e. we have a globalization.
\end{defn}

Note that the functional $f\in {\frakt h}^* $ can be identified with the
functional $F\in ({\frakt h \oplus u}(1)_{\Bbb C})^*$ such that
$F|_{{\frakt u}(1)_{\Bbb C}} \equiv 0 $, i.e. we have an inclusion
${\frakt h}^* \hookrightarrow ({\frakt h \oplus u}(1)_{\Bbb C})^*$.

{\bf Condition A3.9. }
{\it There exists a positive root system $\Phi^+$ and a number $C > 0$, such
 that: If the bundle ${\cal E} _{V,\sigma} \twoheadrightarrow S^{U(1)}$
is irreducible, $\lambda = D\chi_F^{U(1)}|_{\frakt h} \in {\frakt h}^*$,
$\lambda_{\Bbb R}:= \lambda|_{{\frakt h}_{\Bbb R}} $, the restriction of
$\lambda$ to the real form ${\frakt h}_{\Bbb C}$ on which roots take real
values, and $\langle \lambda_{\Bbb R},\alpha\rangle  > C$, for all $ \alpha \in
\Phi^+$, then the bundle ${\cal E} _{V,\sigma} \twoheadrightarrow
S^{U(1)}$ has both properties (MG) and (Z).}

\begin{prop}
Fix a pair $(H,{\frakt b})$. If the condition A3.9 holds, then for
 every datum $(H, {\frakt b},\tilde{\sigma})$, the associated bundle ${
\cal E} _{V,\sigma} \twoheadrightarrow S^{U(1)}$ has both properties (MG)
and (Z). \end{prop}

\begin{pf} We can always reduce the argumentation to the case where the
fiber bundle is irreducible because if the assertion of the property is
not true, it's must be for some irreducible component. Choose $r_0$ as
in \cite{robinsonrawnsley}, we see that $\lambda$ satisfies the condition A3.9: Suppose
$\lambda_0\in{\frakt h}^*_{\Bbb R}$; $\langle \lambda_0,\alpha\rangle  > 0,\forall
\alpha\in\Phi^+$. Thus the bundle ${\cal E} _{V,\sigma}
\twoheadrightarrow S^{U(1)}$ has both properties (MG) and (Z). Fix such a
$\lambda_0$ and put $$s_1:= \sup\{r>0 | \parallel \lambda_{\Bbb R} -
\lambda_0\parallel < r \enskip implies\enskip (MG) \enskip and \enskip (Z)
\enskip for \enskip {\cal E} _{V,\sigma} \twoheadrightarrow S^{U(1)} \}
.$$ We see that $s_1 \geq r_0$. This number $s_1$ must be infinite,
because in other case, we could choose $s_2 >s_1$ with the same
properties, see \cite{robinsonrawnsley} for more detail.  \end{pf}

Recall the notation $S:= G.{\frakt b}\subset X$, and $u:= \codim_{\Bbb
R}{ S}$.

\begin{prop} The properties (MG) and (Z) satisfy for every maximally real
polarization ${\frakt b}$.  \end{prop} 
\begin{pf}

{\it Property (Z). }

 Recal that ${\frakt b \subset p}$, for some cuspidal parabolic subgroup
$P = (M \times A) \ltimes N_H$, $M \subseteq K$, $H = T \times A$, with
$T = H \cap K$, $A = \exp{({\frakt a \cap g})} $. Then $S({U(1)} \cong
(H.N_H) \setminus G$ and $S^{U(1)}$ is a fibration over $P\setminus G$
with holomorphic fibers $T\setminus M$. Let us suppose, what we can
always do, that the fiber bundle ${\cal E} _{V,\sigma} \twoheadrightarrow
S^{U(1)}$ is irreducible, $\lambda = D\chi_F^{U(1)}|_{\frakt h} \in
{\frakt h}^*$ and $\chi_T^{U(1)}:= \chi_F^{U(1)}|_{T^{U(1)}}$, where
$T^{U(1)}$ is the $U(1)$-covering, i.e. the inverse image of $T$ in
$H^{U(1)}$. We see that $D\chi_T^{U(1)} |_{\frakt t} = r|_{\frakt t}$
Suppose that the scalar curvature of the bundle ${\cal
E}_{V,\sigma}|_{T\setminus M}$ is sufficiently negative, i.e. less than
some
 negative number. Then it is not hard to check the isomorphism
$$H^p(T\setminus M;{\cal E} _{V,\sigma}|_{T\setminus M})_{(K\cap M)}
{\buildrel \cong \over \longrightarrow } A^p(M;T,{\frakt b\cap
m},\chi_T^{U(1)}).$$ These $({\frakt m},K\cap M)$-modules are non zero
just for $$p = \dim_{\Bbb C}{T\setminus (K\cap M)}.$$

Let $Z^p$ and $B^p$ denote the corresponding spaces of closed and exact,
 respectively $(K\cap M)$-finite ${\cal E} _{V,\sigma}$-valued
$(0,p)$-forms om $T\setminus M$, and ${^\circ}Z^p$ and ${^\circ}B^p$
denote the corresponding spaces with ``smooth'' replaced by ``formal power
series'' for the coefficients. It is also not hard to check by using the
Taylor series decomposition, that as $({\frakt m},K\cap M)$-modules,
$$B^p\setminus Z^p \cong {^\circ}B^p \setminus {^\circ}Z^p.$$

Applying the Poincar\'e Lemma to the fibers of the fibration $N_H
\rightarrowtail H\setminus G \twoheadrightarrow (H.N_H)\setminus G \approx
S^{U(1)}$ we see that $A^p(G,H,{\frakt b},\chi_F^{U(1)})$ can be computed
from the complex of left $K$-finite, righ$K\cap M$-invariant functions
from $K$ to the Zuckerman complex for $T\setminus M$, i.e. it is induced
from $A^p(M,T,{\frakt b\cap m},\chi_T^{U(1)})$. Thus applying the functor
$\Ind$ to both the sides of $$H^p(T\setminus M;{\cal
E}_{V,\sigma}|_{T\setminus M})_{(K\cap M)} {\buildrel \cong \over
\longrightarrow } A^p(M;T,{\frakt b\cap m},\chi_T^{U(1)}).$$ we have
isomorphism $$H^p(S^{U(1)};{\cal E} _{V,\sigma})_{(K)} {\buildrel \cong
\over \longrightarrow} A^p(G,H,{\frakt b},\chi_F^{U(1)}).$$

{\it Property (MG). } We can also here suppose that the fiber bundle
${\cal E} _{V,\sigma} \twoheadrightarrow S^{U(1)}$ is irreducible,
$\lambda = D\chi_F^{U(1)}|_{\frakt h} \in {\frakt h}^*$. Let $\lambda =
\nu + \sqrt{-1}\sigma$, $\nu \in \sqrt{-1}({\frakt t \cap g})^*$ is in
the interior of the negative Weyl chamber of the root system $\Phi({\frakt
m,t})$. We see that $$H^p(S^{U(1)};{\cal E} _{V,\sigma}) =
H^p(C^{-\omega}_{P\setminus G}(S^{U(1)};\pi_X^*{\cal E}_{V,\sigma}\otimes
\wedge^*{\cal N} _{S^{U(1)}}))$$ vanishe except in degree $p_0 =
\dim_{\Bbb C}{(T\setminus K\cap M)}$ and in this dimension,
$$\Ind^G_{MAN_H}(\eta \otimes e^{\sqrt{-1}\sigma} ):=
H^{p_0}(S^{U(1)};{\cal E} _{V,\sigma}).$$ The induced module
$H^{p_0}(S^{U(1)};{\cal E} _{V,\sigma})$ has finite length because $\eta$
is irreducible and its restriction to $A$ is a multiple of $
\chi_{\nu}^{U(1)}$. As in (\cite{robinsonrawnsley}, Lemma 9.8) we see that this induced
module $H^{p_0}(S^{U(1)};{\cal E} _{V,\sigma})$ satisfies the property
(MG), and the operator $\overline{\partial}_{S^{U(1)}}$ has closed range.
In particular it inherits a Fr\'echet topology from the space
$C^{-\omega}_{P\setminus G}(S^{U(1)};{\cal
E}_{V,\sigma}\otimes\wedge^{p_0}{\cal N}^*_{S^{U(1)}})$. This completes
the proof of the proposition.  \end{pf}

It is therefore proved for every maximally real polarization the
following result 

\begin{thm}
For any maximally real polarization ${\frakt b}$ and any basic datum 
 $(H,{\frakt b},\chi_F^{U(1)})$, there are topological isomorphisms
between Fr\'echet $G$-modules 
$$\begin{array}{rl} H^p(C^{-\omega}(H \setminus G;{\cal
E}_{V,\sigma} \otimes\wedge^*{\cal N} ^*)) &\cong
H^p(C^{-\omega}_{H^{U(1)}}(S^{U(1)};\pi^*_X{\cal
E}_{V,\sigma}\otimes\wedge^*{\cal N}^*_{S^{U(1)}}))\\   & \cong
H^{p+u}(\tilde{S};{\cal O} (\tilde{{\cal E}}_{V,\sigma})),\end{array}$$ 
which are
canonically and topologically isomorphic to the action of $G$ on the
maximal globalization of $A^p(G,H,{\frakt b},\chi_F^{U(1)})$.
\end{thm}

Now we consider arbitrary polarization, not only the maximally real. The
final result is

\begin{thm}
Fix the Cartan sub-algebra $H$ and consider an arbitrary, not necessarily
maximally real polarization ${\frakt b}$. Then for a basic datum $(H,{\frakt b},\chi_F^{U(1)})$, the associated bundle ${\cal E}
_{V,\sigma} \twoheadrightarrow S^{U(1)}$ has both the properties (MG) and (Z).
In other words, the theorem A3.12 holds for  arbitrary basic data $(H,{\frakt
 b},\chi_F^{U(1)})$.
\end{thm}
\begin{pf} 
Suppose that $H = G_F$ is fixed, ${\frakt b \subset g}_{\Bbb
C}$ is a polarization such that ${\frakt h \subset b}$ and that ${\frakt
b}$ isn't maximally real.  It is easy to see that there exists a complex
simple root $\alpha$ such that $\overline{\alpha} \not\in\Phi^+$. Denote
$\Phi^+_0:= s_{\alpha}\Phi^+$, ${\frakt b}_0:= s_{\alpha}{\frakt b}$
and $S_0:= G.{\frakt b}_0$.

Given $\gamma \in \Phi({\frakt g}_{\Bbb C},{\frakt h})$, we can view
$\gamma$ as an element of $({\frakt h \oplus u}(1)_{\Bbb C})^*$. Since
${\frakt h}$ is the Cartan sub-algebra of ${\frakt g}_{\Bbb C}$, we obtain
a representation $\chi_{\gamma}:= e^{\gamma}: H^{U(1)} \rightarrow {\Bbb
C}^{\times}$.  Then we have a vector bundle ${\cal L} _{\gamma}^{U(1)}
\twoheadrightarrow S_0^{U(1)}$ and a vector bundle ${\cal
L}_{\gamma}^{U(1)} \twoheadrightarrow S^{U(1)}$. Applying Lemma 10.6 from
\cite{robinsonrawnsley}, we can obtain $G$-equivariant morphisms of
complexes
$$C^{-\omega}_{H^{ U(1)}}(S_0^{U(1)};\pi_X^*{\cal E} _{V,\sigma} \otimes
\wedge^p{\cal N} ^*_{S_0^{U(1)}}) \rightarrow
C^{-\omega}_{H^{U(1)}}(S^{U(1)};\pi_X^*{\cal E} _{V,\sigma}\otimes {\cal
L} _{-\alpha}^{U(1)} \otimes \wedge^{p+1}{\cal N}^*_{S^{U(1)}})$$ and this
morphism restricts to a morphism of subcomplexes $$C^{-\omega}_{P\setminus
G}(S_0^{U(1)};\pi^*_X{\cal E} _{V,\sigma}\otimes \wedge^p{\cal N}
^*_{S_0^{U(1)}}) \rightarrow C^{-\omega}_{P\setminus
G}(S^{U(1)};\pi^*_X{\cal E} _{V,\sigma} \otimes {\cal L} _{\gamma}^{U(1)}
\otimes \wedge^{p+1}{\cal N} ^*_{S^{U(1)}}).$$

Let $C^{-\omega}_{S_0^{U(1)}}(S^{U(1)};\pi^*_X{\cal E} _{V,\sigma}\otimes
{\cal L}_{-\alpha}^{U(1)} \otimes \wedge^*{\cal N} ^*_{S^{U(1)}})$ be the
subcomplex of the complex $C^{-\omega}(S^{U(1)};\pi^*_X{\cal E}_{V,\sigma}
\otimes {\cal L} _{-\alpha}^{U(1)}\otimes \wedge^*{\cal N} ^*_{S^{U(1)}})$
, consisting of forms $\omega$ such that $\overline{\partial
}_{S^{U(1)}}\omega$ vanish on $(0,1)$ vectors tangent to the fibers of
$S^{ U(1)} \twoheadrightarrow S_0^{U(1)}$. Applying the Dolbeault Lemma,
we see
 that the inclusion $$C^{-\omega}_{S_0^{U(1)}}(S^{U(1)};\pi^*_X{\cal
E}_{V,\sigma}\otimes{\cal L} _{-\alpha}^{U(1)}\otimes \wedge^*{\cal N}
^*_{S^{U(1)}}) \hookrightarrow C^{-\omega}(S^{U(1)};\pi^*_X{\cal
E}_{V,\sigma}\otimes {\cal L} _{-\alpha}^{U(1)}\otimes \wedge^*{\cal
N}^*_{S^{U(1)}})$$ induces isomorphisms on cohomology.

On the other hand we have a morphism of complexes
$$C^{-\omega}(S_0^{U(1)};\pi^*_X{\cal E} _{V,\sigma}\otimes \wedge^p{\cal
N} ^*_{S^{U(1)}}) \rightarrow
C^{-\omega}_{S_0^{U(1)}}(S^{U(1)};\pi^*_X{\cal E} _{V,\sigma}\otimes{\cal
L} _{-\alpha}^{U(1)} \otimes \wedge^{p+1}{\cal N} ^*_{S^{U(1)}}).$$

Let $${\frakt b}_{\alpha} = {\frakt b \oplus g}_{\alpha} = {\frakt
b}_0 \oplus {\frakt g}_{-\alpha}.$$ Denote by $X_{\alpha}$ the flag
manifold of parabolic sub-algebras of ${\frakt g}_{\Bbb C}$ which are 
$Int({\frakt g}_{\Bbb C})$-conjugate to ${\frakt b}_{\alpha}$ and consider the
 orbit $S_{\alpha}:= G.{\frakt b}_{\alpha} \subset X$. The natural
projection $p_{\alpha}: X \twoheadrightarrow X_{\alpha}$ is holomorphic
and there is a $U(1)$-covering homomorphism $p_{\alpha}^{U(1)}: X^{U(1)}
\rightarrow X_{\alpha}^{U(1)}$ such that $$p_{\alpha} \circ \pi_X =
\pi_{X_{\alpha}}
 \circ p_{\alpha}^{U(1)},$$ where $$\pi_{X_{\alpha}}: X_{\alpha}^{U(1)}
\rightarrow X_{\alpha}$$ is the natural projection.

Let $U_{\alpha} \subset S_{\alpha}$ be an open subset, whose
$\overline{U_{\alpha}}$ is compact and has an $X_{\alpha}$-open
neighborhood over which $p_{\alpha}: X \rightarrow X_{\alpha}$ is
holomorphly trivial. Let $$U_0^{U(1)}:= S_0^{U(1)} \cap
(p_{\alpha}{U(1)})^{-1} \pi^{-1}_{X_{\alpha}}(U_{\alpha})$$ and $$U^{U(1)}
:= S^{U(1)} \cap (p_{\alpha}^{U(1)})^{-1}\pi_{X_{\alpha}}^{-1}
(U_{\alpha}),$$ we see that $$C^{-\omega}(S_0^{U(1)};\pi^*_X{\cal E}
_{V,\sigma}\otimes \wedge^p{\cal N} ^*_{S^{U(1)}}) \rightarrow
C^{-\omega}_{S_0^{U(1)}}(S^{U(1)};\pi^*_X{\cal E} _{V,\sigma}\otimes{\cal
L} _{-\alpha}^{U(1)} \otimes \wedge^{p+1}{\cal N} ^*_{S^{U(1)}}) $$
localizes to maps $$C^{-\omega}_{H^{U(1)}}(U_0^{U(1)};\pi^*_X{\cal E}
_{V,\sigma}\otimes \wedge^p{\cal N} ^*_{S_0^{U(1)}}) \rightarrow
C^{-\omega}_{H^{U(1)},S_0^{U(1)}}(U^{U(1)};\pi^*_X{\cal
E}_{V,\sigma}\otimes{\cal L} _{-\alpha}^{U(1)} \otimes \wedge^{p+1}{\cal
N} ^*_{S^{U(1)}}).$$

Let $Cl(U_0^{U(1)})^{\sim}$, the closure, and $Bd(U_0^{U(1)})^{\sim}$,
the boundary, denote germs of neighborhood of
 $Cl(U_0^{U(1)})$ and $Bd(U_0^{U(1)})$ in $S^{U(1)} \cup S_0^{U(1)}$.
From the theory of hyperfunctions, it is not hard to see that
$$C^{-\omega}( U_0^{U(1)};\pi_X^*{\cal E} _{V,\sigma}\otimes \wedge^*{\cal
N} ^*_{S_0^{U(1)}}) = {C^{\omega}(Cl(U_0^{U(1)});{\cal F}
^{U(1)}\otimes\wedge^{c-p}{\cal N} ^*_{S_0^{U(1)}})' \over
C^{\omega}(Bd(U_0^{U(1)});{\cal F} ^{U(1)}\otimes\wedge^{c-p}{\cal N}
^*_{S_0^{U(1)}})' }, $$ where $c:= \dim_{CR}{S_0}$, the complex
dimension with respect to the complex structure given by the indicate
CR-structure, and ${\cal F} ^{U(1)}:= {\cal E} _{V,\sigma} \otimes {\cal
L} _{-2\rho + 2\alpha}^{U(1)}$.

Similarly, we have $$C^{-\omega}_{S_0^{U(1)}}(U^{U(1)};\pi^*_X{\cal
E}_{V,\sigma}\otimes{\cal L} _{-\alpha}^{U(1)} \otimes\wedge^{p+1}{\cal N}
^*_{S^{U(1)}}) = {
C^{\omega}_{S_0^{U(1)}}(Cl(U_0^{U(1)})^{\sim};\tilde{\cal F} ^{U(1)}
\otimes\wedge^{c-p}{\cal N} ^*_{S^{U(1)}})' \over
C^{\omega}_{S_0^{U(1)}}(Bd(U_0^{U(1)})^{\sim};\tilde{\cal F} ^{U(1)}
\otimes \wedge^{c-p}{\cal N} ^*_{S^{U(1)}})'}, $$ where $c+1 =
\dim_{CR}{S_{\alpha}}$, ${\cal F} ^{U(1)} = (\pi_X^*{\cal E} _{V,\sigma}
\otimes {\cal L} _{-\alpha}^{U(1)}) \otimes {\cal L} _{-2\rho +
\alpha}^{U(1)}$. Thus by a similar argument, we obtain the dual
statement, as follows: The restriction maps
$$C^{\omega}_{H^{U(1)},S_0^{U(1)}} (Cl(U_0^{U(1)})^{\sim};\tilde{\cal F}
^{U(1)} \otimes\wedge^*{\cal N} ^*_{S^{U(1)}}) \rightarrow
C^{\omega}_{H^{U(1)}}(Cl(U_0^{U(1)});{\cal F} ^{U(1)}\otimes\wedge^*{\cal
N} ^*_{S_0^{U(1)}}),$$
$$C^{\omega}_{H^{U(1)},S_0^{U(1)}}(Bd(U_0^{U(1)})^{\sim};\tilde{\cal F}
^{U(1)} \otimes \wedge^*{\cal N} ^*_{S^{U(1)}}) \rightarrow
C^{\omega}_{H^{U(1)}}(Bd(U_0^{U(1)});{\cal F} ^{U(1)} \otimes
\wedge^*{\cal N} ^*_{S_0^{U(1)}})$$
 induce isomorphisms in cohomology.

We know that these restriction maps are continuous and surjective, and
are dual via $$C^{-\omega}( U_0^{U(1)};\pi_X^*{\cal E} _{V,\sigma}\otimes
\wedge^*{\cal N} ^*_{S_0^{U(1)}}) = {C^{\omega}(Cl(U_0^{U(1)});{\cal F}
^{U(1)}\otimes\wedge^{c-p}{\cal N} ^*_{S_0^{U(1)}})' \over
C^{\omega}(Bd(U_0^{U(1)});{\cal F} ^{U(1)}\otimes\wedge^{c-p}{\cal N}
^*_{S_0^{U(1)}})' }, $$ and
$$C^{-\omega}_{S_0^{U(1)}}(U^{U(1)};\pi^*_X{\cal E}_{V,\sigma}\otimes{\cal
L} _{-\alpha}^{U(1)} \otimes\wedge^{p+1}{\cal N} ^*_{S^{U(1)}}) = {
C^{\omega}_{S_0^{U(1)}}(Cl(U_0^{U(1)})^{\sim};\tilde{\cal F} ^{U(1)}
\otimes\wedge^{c-p}{\cal N} ^*_{S^{U(1)}})' \over
C^{\omega}_{S_0^{U(1)}}(Bd(U_0^{U(1)})^{\sim};\tilde{\cal F} ^{U(1)}
\otimes \wedge^{c-p}{\cal N} ^*_{S^{U(1)}})'}, $$ to the maps
$$C^{-\omega}_{H^{U(1)}}(U_0^{U(1)};\pi^*_X{\cal E} _{V,\sigma}\otimes
\wedge^p{\cal N} ^*_{S_0^{U(1)}}) \rightarrow
C^{-\omega}_{H^{U(1)},S_0^{U(1)}}(U^{U(1)};\pi^*_X{\cal
E}_{V,\sigma}\otimes{\cal L} _{-\alpha}^{U(1)} \otimes \wedge^{p+1}{\cal
N} ^*_{S^{U(1)}}).$$

Thus we obtain the following statement: {\it Suppose that
$\tilde{\sigma}\in X^{U(1)}_{irr}$, $\lambda = D\chi_F^{U(1)}|_{\frakt h}
\in {\frakt h}^*$, and suppose that $2{\langle \lambda + \rho - \alpha, \alpha\rangle
\over \langle \alpha,\alpha\rangle }$ is not a positive integer. Then
$$C^{-\omega}_{H^{ U(1)}}(S_0^{U(1)};\pi_X^*{\cal E} _{V,\sigma} \otimes
\wedge^p{\cal N} ^*_{S_0^{U(1)}}) \rightarrow
C^{-\omega}_{H^{U(1)}}(S^{U(1)};\pi_X^*{\cal E} _{V,\sigma}\otimes {\cal
L} _{-\alpha}^{U(1)} \otimes \wedge^{p+1}{\cal N}^*_{S^{U(1)}})$$ induces
an isomorphism of the corresponding cohomology groups. }

One rests only to apply Proposition A3.11 and we may assume by induction
on $\dim{S^{U(1)}} -\dim{S^{U(1)}_{max}}$ that every ${\cal E} _{V,\sigma}
\twoheadrightarrow S_0^{U(1)}$ admits both the properties (MG) and (Z).
Since the cohomologies and the maps that occur in theorem A3.12 all are
compatible with coherent continuation, we may assume that $2\langle \lambda +
\rho - \alpha, \alpha\rangle  / \langle \alpha,\alpha\rangle $ is not a positive integer,
where $\tilde{\sigma}$ is irreducible, $\lambda = D\chi_F^{U(1)}|_{\frakt
h} \in {\frakt h}^*$. From the last assertion {\it Suppose that
$\tilde{\sigma}\in X^{U(1)}_{irr}$, $\lambda = D\chi_F^{U(1)}|_{\frakt h}
\in {\frakt h}^*$, and suppose that $2{\langle \lambda + \rho - \alpha, \alpha\rangle
\over \langle \alpha,\alpha\rangle }$ is not a positive integer. Then
$$C^{-\omega}_{H^{ U(1)}}(S_0^{U(1)};\pi_X^*{\cal E} _{V,\sigma} \otimes
\wedge^p{\cal N} ^*_{S_0^{U(1)}}) \rightarrow
C^{-\omega}_{H^{U(1)}}(S^{U(1)};\pi_X^*{\cal E} _{V,\sigma}\otimes {\cal
L} _{-\alpha}^{U(1)} \otimes \wedge^{p+1}{\cal N}^*_{S^{U(1)}})$$ induces
an isomorphism of the corresponding cohomolohy groups }, one follows that
$$C^{-\omega}_{P\setminus G}(S^{U(1)};\pi^*_X{\cal E} _{V,\sigma} \otimes
{\cal L} _{-\alpha}^{U(1)} \otimes \wedge^{p+1}{\cal N} ^*_{S^{U(1)}})$$
has both the properties (MG) and (Z). This completes the proof of the
theorem.
\end{pf}

\section{Quantization of Mechanical Systems with Supersymetry}

{\it `` Graded Lie algebras have recently become a topic of interest
in Physics in the context of ``supersymmetry'' relating particles of
different statistics '' } 
, as it is pointed out in the survey \cite{corvinneemansternberg}
. In the physical systems where the Bose-Einstein particles and
Fermi-Dirac particles interact together symmetry must be replaced by
supersymmetry. Lie superalgebras and Lie supergroups are therefore
important mathematical tools of physics and they must be studied
seriously.

B. Kostant has developed in his work \cite{kostant2} the representation theory
of Lie supergroups along the line of the Kirillov-Kostant orbit
method for the ordinary Lie groups. His theory is founded in
differential geometry and uses symplectic structures, Hamiltonian
formalism, integrality condition, line bundles with connection and
prequantization for the supergroup case. But polarizations and
unitarity are lacking by himself in general. Perhaps the main
difficulty is the fact that the Lie's theorem is not true in in the
general case of solvable Lie superalgebras: They have also
irreducible finite-dimensional representations rather than the
(one-dimensional) characters. This means that we must work with
multidimensional induced bundles in place of line bundles.

Using the new notion of polarization, exposed in the previous chapters,
one can develop also a theory of quantization for Lie supergroups, see
\cite{diep22}.

\subsection{Hilbert superbundles with connection }

\begin{defn} As usual, let us denote by ${\Bbb C}$ the ground field of
complex numbers, by ${\Bbb Z}/2{\Bbb Z}$ the residue field consisting of
two elements ${\bar 0}$ and ${\bar 1}$. Recalll that a {\it (complex)
vector superspace } \index{superspace} $V$ is by definition a ${\Bbb 
Z}/2{\Bbb Z}$-graded
vector space $V = V_{\bar 0} \oplus V_{\bar 1}$. The elements of $V_{\bar
0}$ are called {\it even }, \index{element!even -} i.e. their graded degree 
is 0, $|x| = 0$;
those of $V_{\bar 1}$ are called {\it odd }, \index{element!odd -} i.e. 
their graded degree is 1, $|x| = 1$.  \end{defn}

Throughout what follows if for some $X$ expression $|x|$ occurs, then it
is assumed that $X$ is homogeneous, and that the expression extends to
other elements by linearity.

Suppose that the vector superspace $V = V_{\bar 0} \oplus V_{\bar 1}$
admits a complex sesquilinear form $b: V \times V \rightarrow {\Bbb C}$,
a so called {\it scalar product } \index{product!scalar -} which is linear 
relative to the 
first variable, and:  \begin{itemize} \item {\it superhermitian }, 
\index{superhermitian} i.e. for all homogeneous $x$, $y$ in $V$, $$b(x,y) =
(-1)^{|x|.|y|}\overline{b(y,x)}, $$ \item {\it consistent }, 
\index{product!consistent -} i.e. for
all $x$ and $y$ of different graded degrees, $$b(x,y) \equiv 0,$$ and
\item {\it nondegenerate }, \index{product!nondegenerate -} i.e. if $b(x,y) 
\equiv 0, \forall y\in V$ then $x = 0$.  \end{itemize}

It is easy to deduce the following 

\begin{cor} Let $(V,b)$ be a vector superspace with scalar product.
\begin{enumerate} \item{(i)} The restriction of $b$ on $V_{\bar 0}$ is a
scalar product, and the restriction of $b$ to $V_{\bar 1}$ is a
nondegenerate skew-symmetric form. \item{(ii)} The correspondence $z
\mapsto b(.,z)$ establishes a monomorphism $V \hookrightarrow V^*$.
\end{enumerate} \end{cor}

\begin{rem} In the category of superspaces there is a special {\it functor
of changing the graduation degrees $\Pi$, } \index{functor!of changing 
graduation} $$(\Pi V)_{\bar 0}:= V_{\bar 1}, (\Pi V)_{\bar 1}:= V_{\bar 0}.$$
\end{rem}

{\it Suppose that there is also a scalar product $b^{\Pi}$ on $\Pi V$
and that $(V,b)$ and $(\Pi V, b^{\Pi})$ are the Hilbert spaces. In
this case we say that $(V,b,b^{\Pi})$ is a Hilbert superspace. }
It must be remarked that in work \cite{diep22} it was not true that the
Hilbert superspace structure can be deduced from a single scalar
product $b(.,)$. Indeed it must be defined by a pair, consisting
of two scalar product $b(.,)$ and $b^{\Pi}(.,.)$,with respect to which
the space $V$ is complete.

\begin{defn}
For an operator $u \in \Aut (V,b,b^{\Pi})$, we define
$u^{\Pi}:= \Pi \circ u \circ \Pi \in \Aut ( b,b^{\Pi})$. We have therefore
a commutative diagram $$\vbox{\halign{ #&#&#&#&#\cr 
  &$V$ &$ {\buildrel u \over \longrightarrow} $ &  &$V$ \cr
  &${\downarrow}\Pi$ & & &${\downarrow}\Pi$ \cr $\Pi$&$V$ &$ {\buildrel
u^{\Pi} \over \longrightarrow} $&$\Pi$&$V.$ \cr }}$$ \end{defn}

\begin{cor} $u$ is symplectic in $V_{\bar 1}$ if and only if $u^{\Pi}$
is unitary on the Hilbert space $W:= \Pi V_{\bar 1}$.  \end{cor}

Let us consider a Hilbert superspace $(V,b,b^{\Pi})$. Denote by $q$ and
$q^{\Pi}$ the corresponding quadratic forms on $V$ and $\Pi V$
respectively.

\begin{cor}
\begin{enumerate}
 \item[(i)] $q$ is nondegenerate on $V_{\bar 0}$, but identically 0 on
$V_{\bar 1}$.  \item[(ii)] $q^{\Pi}$ is nondegenerate on $\Pi V_{\bar 1}$
, but identically 0 on $\Pi V_{\bar 0}$ \end{enumerate} \end{cor}

\begin{defn} We now define the norm of supervectors on $V$ by the formula
$$\parallel v \parallel = \sqrt{q(x) + q^{\Pi}(x)}, $$ forall $x\in V
\cong \Pi V$ as complex vector spaces, and as usually the associated norm
of operators, $$\parallel A \parallel
 = \sup_{\parallel x \parallel \leq 1}{\parallel Ax \parallel }. $$
\end{defn}

Now it is easy to see that the vector superspace $V = V_{\bar 0}
\oplus V_{\bar 1}$ with scalar product $b$ on $V$ and $b^{\Pi}$ on
$\Pi V = (\Pi V)_{\bar 0} \oplus (\Pi V)_{\bar 1}:= V_{\bar 1} \oplus
V_{\bar 0} \cong V$ ( as complex vector spaces ) becomes a Hilbert
superspace iff it is complete with respect to this norm of
supervectors.

By an {\it (even) unitary-symplectic operator }
$u\in \USp (V)$ we mean an even automorphism $u \in \Aut (V,b,b^{\Pi})_{\bar 0}$
, i.e. $u V_{\bar i} \subseteq V_{\bar i}, i = {\bar 0},{\bar 1}$.

\begin{cor} 
The set of all unitary-symplectic operators forms a group denoted
also $\USp (V)$. On this group the strong topology is equivalent to the
weak topology.
\end{cor} 

The proof is the same as in the classical case by using the well-known
Cauchy-Bounjakowski-Schwarz inequality for scalar product.

\begin{lem}
The topological group $\USp (V)$ is isomorphic to the direct product of
the unitary groups of $V_{\bar 0}$ and $(\Pi V)_{\bar 1}$.
\end{lem}
\begin{pf}
By definition we have $b(ux,uy) = b(x,y)$ and $b^{\Pi}(ux,uy) = b^{\Pi}(x,y)$,
i.e. $$ u = \pmatrix u|_{V_{\bar 0}} & 0 \cr
                      0 & u|_{V_{\bar 1}}\endpmatrix 
         \cong \pmatrix u|_{V_{\bar 0}} & 0 \cr
                         0 & u|_{(\Pi V)_{\bar 0}}. \endpmatrix$$
We have therefore $$\begin{array}{rl} \Aut(V,b, b^{\Pi})_{\bar 0} &= \Aut(V_{\bar 0},b)
 \times \Aut(V_{\bar 1}, b)\cr   &\cong \Aut(V_{\bar 0},b) \times \Aut(\Pi
V_{\bar 1}, b^{\Pi}).\end{array}$$
\end{pf}

\begin{defn}
A linear superoperator $A: V \rightarrow V$ is said to be {\it
antisupersymmetric }
iff for every homogeneous element $x,y \in V$, $$b(Ax,y) = -(-1)^{|A|.|x|}b(x,
Ay), $$ $$b^{\Pi}(Ax,y) = -(-1)^{|A|.|x|}b^{\Pi}(x,Ay).$$
\end{defn}

\begin{lem}
Antisupersymmetric superoperators form a Lie superalgebra, the even
part of
which is a Lie algebra consisting of all pairs of antisymmmetric (unbounded)
operator on the direct product of Hilbert spaces $V_{\bar 0} \times \Pi V_{
\bar 1}$.
\end{lem}
\begin{pf}
It is enough to verify that the vector superspace of antisupersymmetric 
superoperators is stable under the brackets $$[A,B]:= AB -(-1)^{|A|.|B|}BA
.$$ But this is clear from the definition; see also (\cite{kac},\S5.3.4(b)).
\end{pf}

From the well-known Stone theorem for Hilbert superspaces, it is easy to
deduce the corresponding superversion.

\begin{thm}[\bf Stone Theorem ] Every one (real) parameter continuous
subgroup $\{u(t)\}_{t\in {\Bbb R}}$ of even unitary-symplectic
superoperators on a Hilbert superspace $V$ admits a generator $\sqrt{-1}A$
, which is antisymmetric continuous (perhaps, unbounded) superoperator.
\end{thm}

\begin{defn} Let $(G,A)$ be a Lie supergroup and let ${\frakt g}:= \Lie
(G,A)$ be the Lie superalgebra, $U({\frakt g})$ the enveloping
superalgebra of ${\frakt g}$, ${\Bbb R}(G)$ the function algebra on $G$
and $A(G)^*:= E(G,{\frakt g}) = {\Bbb R}(G) \circledS U({\frakt g})$ the
Lie-Hopf coalgebra, where by $\circledS$ we denote the so called {\it
smash product }, see also \cite{kostant2}.  Let $(V,b,b^{\Pi})$ be a
Hilbert superspace, and let $\End V$ be the superalgebra of continuous
superoperators on $V$. By a {\it smooth unitary representation } of
$(G,A)$ on $V$, we mean a homomorphism of superalgebras $$r: A(G)^*
\longrightarrow \End V,$$ such that:  \begin{enumerate} \item[(1)] The
restriction $r|_{G}$ is a continuous representation of the Lie group $G$
in the group of even unitary-symplectic automorphisms $$r|_{G}: G
\longrightarrow \USp (V).$$ \item[(2)] Each vector $v\in V$ is smooth,
i.e. the map $G \rightarrow V$, $g \mapsto r(g)v$ is of class
$C^{\infty}$.  \item[(3)] $r(X)v = {d \over dt}(r(\exp{tX})v)|_{t=0},
\forall v\in V, \forall X \in {\frakt g}_0 $. \end{enumerate} \end{defn}

The rest of this section is in large similar to the same of (\cite{kostant2},\S4), where
the theory is developed purely for line superbundles. Our need is
essentially in the multidimensional case. So we are trying to modify these
results and notions for the multidimensional situation.

{\bf Notation} Let $(X,A)$ be a supermanifold and let $U \subseteq X$ be
an open set, $\Der  A(U)$the Lie superalgebra of all superderivations in
$A(U)$, which is also a $A(U)$-module, $T(U)$ the tensor algebra of $\Der
A(U)$ over $A(U)$, $$T(U) = \bigoplus_{b=0}^{\infty} T^b(U) =
\bigoplus_{b=0}^{\infty} \underbrace{\Der  A(U) \otimes_{A(U)} \dots
\otimes_{A(U)} \Der  A(U)}_{b\rm\;times},$$ which is ${\Bbb Z} \oplus {\Bbb
Z}/2{\Bbb Z}$-(bi)graded.

Now let $J(U)$ be the two-side ${\Bbb Z} \oplus {\Bbb Z}/2{\Bbb
Z}$-graded ideal in $T(U)$ generated by the elements in $T^2(U)$ of
the form $\xi \otimes \eta + (-1)^{|\xi|.|\eta|}\eta \otimes \xi$,
where $\xi,\eta \in \Der  A(U)$ are homogeneous. Denote also $J^b(U) =
T^b(U) \cap J(U)$.

Let us denote by $E$ a projective (locally free) $A$-module sheaf on
$(X,A)$. Then for every $b \in {\Bbb N}$, $Hom_{A(U)}(T^b(U),E(U))$
can be considered as the set of all bilinear maps on $\Der  A(U)$ with
values in $E(U)$ which satisfy the conditions $$\langle \xi_1,\dots,f
\xi_l,\dots, \xi_b|\beta\rangle  = (-1)^{|f|\sum_{i=1}^l|\xi_i|} f\langle \xi_1,
\dots,\xi_b|\beta\rangle .$$

Now let $\Omega^b(U;E)$ the set of all $\beta \in
\Hom_{A(U)}(T^b(U),E(U))$ which vanish on $J^b(U)$. It is easy to see
that the elements $\beta$ in $ \Omega^b(U;E)$ are characterized by the
additional condition $$\langle \xi_1,\dots,\xi_j,\xi_{j+1},\dots,\xi_b|\beta\rangle
= (-1)^{1+|\xi_j|.|\xi_{j+1}|}\langle \xi_1,\dots,\xi_{j+1},\xi_j,\dots,\xi_b|
\beta\rangle .$$

\begin{cor}
$\Omega^b(U,E) = \Omega^b(U,A) \otimes E(U)$.
\end{cor}
\begin{pf}
It is enough to mension that for the locally free sheaf $E(U)$, $$
\Hom_{A(U)}(T^b(U),E(U)) \cong \Hom_{A(U)}(T^b(U),A(U)) \otimes_{A(U)}
E(U).$$
\end{pf}

\begin{defn}
Let us denote by $\End E$ the projective (locally free) $A$-module
sheaf of all endomorphisms of the superbundle $E$. Then it is easy to
see that $\Omega^b(U,E)$ and $\Hom_{A(U)}(T^b(U),E(U))$ are $\End
E(U)$-modules, $$\langle \xi_1,\dots,\xi_b|\beta f\rangle  =
\langle \xi_1,\dots,\xi_b|\beta\rangle f, \forall f \in \End E(U) $$ and
$$\langle \xi_1,\dots,f\xi_l,\dots,\xi_b|\beta\rangle  = (-1)^{|f|\sum_{i=l}^b
|\xi_i|}\langle \xi_1,\dots,\xi_b|f\beta\rangle ,\forall f \in A(U).$$ Also
$\Omega^b(U,E)$ is ${\Bbb Z}/2{\Bbb Z}$-graded:
$$\langle \xi_1,\dots,\xi_b|\beta\rangle  \in E(U)_k, $$ $$ k:= |\beta| +
\sum_{i=1}^b |\xi_i|.$$ $$\Omega^0(U,E):= E(U), $$ $$\Omega(U,E) =
\oplus_{b=0}^{\infty} \Omega^b(U,E).$$
\end{defn}

Let $\beta \in \Omega^b(U,E)$, $\alpha \in \Omega^a(U,E')$,
$\langle .,.\rangle $be the pairing of $E$ and $E'$ with values in $F$, then one
defines $\Omega^{a+b}(U,F)$ in an usual way; for example, $\langle .,.\rangle  =
\Hom(E\otimes E',F)$ is an important case.

Remark that if $V\subseteq U$ is an open subset, one has a
restriction map $$\rho_{U,V}: \Omega(U,E) \rightarrow \Omega(V,E),
$$ such that if $\xi_i \in \Der  A(U)$, $\beta \in \Omega^b(U,E)$ then
$\rho_{U,V}\beta \in \Omega^b(V,E)$ is characterized by $$\langle \rho_{U,V}
\xi_1,\dots,\rho_{U,V}\xi_b|\rho_{U,V}\beta\rangle  =
\rho_{U,V}\langle \xi_1,\dots,\xi_b|\beta\rangle .$$ It is clear that the
correspondence $U \mapsto \Omega(U,E)$ defines the sheaf of
differential superforms with values in a projective (locally free)
$A$-module sheaf $E$.

Assume that $\dim{(X,A)} = (m,n)$. An open set $U$ is called $A$-{\it
parallelizable } 
if there exist $\eta_l \in \Der  A(U)$, $l=1,\dots,m+n$, such that
$\eta_l \in \Der  A(U)_0 $, if $l\leq m$ and $\eta_i \in \Der  A(U)_1$ if
$l > m$, and such that every $\xi \in \Der  A(U)$ can be uniquely
decomposed as $$ \xi = \sum_{l=1}^{m+n} f_l\eta_l,$$ where $f_l \in
A(U) $.

Remark that every $A$-coordinate neighborhood is parallelizable.

Now, if the open $U$ is $A$-parallelizable, one defines $\alpha_l
\in \Omega^1(U,E)$ by the condition that $$\langle \xi|\alpha_l\rangle  = f_l,$$ if
$$ \xi = \sum_{l=1}^{m+n} f_l\eta_l \in \Der  A(U).$$ Thus
$$\langle \eta_k,\alpha_l\rangle  = \delta_{kl}.Id_U.$$ Put $$\beta_l = \alpha_l,
l \leq m $$ and $$\gamma_l = \alpha_{l+m}, 1 \leq l \leq n, $$ then
$\beta_l \in \Omega^1(U,A)_0 $ and $\gamma_l \in \Omega^1(U,A)_1$ and
$$\beta_l \alpha_i = -\alpha_i\beta_l,$$ $$\gamma_i \gamma_j =
\gamma_j \gamma_i.$$ We introduce the usual multiindex notation ( see
also (\cite{kostant2}, \S4.2) ), $$\beta_{\mu}\gamma^{\nu}:= \beta_{\mu_1}\dots
\beta_{\mu_k}\gamma_1^{\nu_1}\dots \gamma_n^{\nu_n}. $$ we see that
in the $A$-parallelizable open set $U$, $\Omega(U,E)$ is a free
$A(U)$-module such that every differential superform $\omega$ can be
written in form $$\omega = \sum_{\mu,\nu} \beta_{\mu}\gamma^{\nu}
f_{\mu \nu}, f_{\mu\nu} \in E(U).$$

Now assume that $U$ is an $A$-splitting coordinate neighborhood with
an coordinate system $\{r_i,s_j\}$, $ i = 1,\dots,m$, $j =
1,\dots,n$. Then $\Der  A(U)$ is a free $A(U)$-module with the
corresponding basis $\{ {\partial \over \partial r_i},{\partial\over
\partial s_j }\}$. We can choose just $\{ \beta_i,\gamma_j \} = \{
dr_i,ds_j \} $. So we have $$dr_{\mu}ds^{\nu} \in
\Omega^{k(\mu)+|\nu|}(U,A) \hookrightarrow \Omega^{k(\mu)
+|\nu|}(U,E)$$ and every $\beta \in \Omega(U,E)$ can be written as
$$\beta = \sum{\mu,\nu} dr_{\mu}ds^{\nu} f_{\mu\nu},$$ where $f_{\mu\nu}
\in E(U)$.

Recall that, as in the line super-bundle case, we can also construct
the map $$\sigma^*: \Omega(Y,E') \longrightarrow \Omega(X,E),$$
which is associated to the super-bundle sheaf morphism $$\sigma: E
\longrightarrow E'.$$

Finally, by a {\it super-bundle } \index{super-bundle}
we means a projective ( locally free ) $A$-module sheaf $E$, such that
there exists a covering by opens which are {\it principal } 
\index{super-bundle!principal -}
for $E$ in the sense that: {\it $E(U)$ is a free $A(U)$-module with a
basal system of even generators $t_i \in E(U)_0$, $i =
1,\dots,rank_{A(U)}E(U)$.}

\begin{defn} Let $E$ be a vector super-bundle over the super-manifold
$(X,A)$. By a {\it connection } $\nabla$ on $E$ we mean a covariant
super-derivation such that for any open set $U \subseteq X$ and any vector
super-field $\xi \in \Der  A(U)$, one has a linear map $$\nabla_{\xi}: E(U)
\rightarrow E(U),$$ where $|\nabla_{\xi}| = |\xi|$, which is compatible
with the restriction maps to smaller open sets and is such that

(1)\qquad $\nabla_{\xi}(ft) = (\xi f)t + (-1)^{|f|.|\xi|}f\nabla_{\xi}(t)$
, for $f \in A(U)$ and $t \in E(U)$, and

(2)\qquad the map $$\vbox{\halign{ #&#&#&#\cr $\Der $&$A(U)$ & $
\longrightarrow $ & $\End E(U),$ \cr &$\xi $ &$ \mapsto $ & $\nabla_\xi $ \cr
}}$$ is $A(U)$- linear.  \end{defn}

The complexification of $\Der  A(U)$ may be taken to be the complex Lie
super-algebra $\Der  A_{\Bbb C}(U)$ of super-fields of vectors. By
linearity we may take $\xi$ and $f$ in $\Der  A_{\Bbb C}(U)$ and
$A_{\Bbb C}(U)$, respectively and $\Omega_{\Bbb C}(U,E)$ will denote
the complexification of $\Omega(U,E)$.

Now assume that $(E,\nabla)$ is a complexified vector super-bundle with
connection, $U \subseteq X$ is principal for $E$ and $t_i$, $i=
1,\dots,rk_{A}E$. Then for every $\xi\in \Der  A_{\Bbb C}(U)$, there
exists $g = g(\xi)\in \End E_{\Bbb C}(U)$, such that $\nabla_{\xi}t =
g(\xi)t$. The correspondence $\xi \mapsto g(\xi)$ defines an $A_{\Bbb
C}(U)$ -linear map $\Der  A_{\Bbb C}(U) \rightarrow \End E(U)$.
Therefore, there exists a unique element $\alpha(t) \in
\Omega^1_{\Bbb C}(U,\End E)$ such that $$\nabla_{\xi}t = {\sqrt{-1}\over\hbar}
\langle \xi|\alpha(t)\rangle t,$$ for all $\xi \in \Der  A_{\Bbb C}(U)$ and
$|\alpha(t)| = |\alpha| = 0 $, i.e. $|\langle \xi|\alpha(t)\rangle | = |\xi| $.
Now if $s_i \in E(U)_0$, then clearly $s_i,
i=1,\dots,rank_{A(U)}E(U)$ is a basal system of generators iff $s =
tf$ for some $f\in \Aut  E(U)$. In this case we have

\begin{lem}
$\alpha(s) - f \alpha(t) t^{-1} = {\hbar\over\sqrt{-1}}df.f^{-1}$.
\end{lem}
\begin{pf}
Because $s=ft$, we have by definition, $$\begin{array}{rl} \nabla_{\xi}s &=
{\sqrt{-1}\over\hbar}\langle \xi|\alpha(s)\rangle s = {\sqrt{-1}\over\hbar}\langle \xi|\alpha(ft)\rangle
 ft \cr
  &= \nabla_{\xi}(ft) = \xi(f)t + (-1)^{|\xi||f|}f\nabla_{\xi}t \cr
  &= {\sqrt{-1}\over\hbar}\langle \xi|{\hbar\over\sqrt{-1}}df\rangle t + (-1)^{|\xi|.|f|} f.
{\sqrt{-1}\over\hbar}\langle \xi|\alpha(t)\rangle t \cr
  &= {\sqrt{-1}\over\hbar}\langle \xi|{\hbar\over\sqrt{-1}}df. f^{-1}\rangle ft + {\sqrt{-1}
\over\hbar}\langle \xi|f\alpha(t)\rangle f^{-1}.ft.\end{array}$$
From this the lemma is deduced.
\end{pf}

{\bf Notation.}
The classical formulas for the differential of a vector valued
function holds also in the super-context. However we must be careful
to use the right or the left $A(U)$-module structure in
$\Omega^1(U,E)$, for example if $f\in E(U)$ is arbitrary, one has
$$df = \sum_{i=1}^m dr_i {\partial f \over \partial r_i} +
\sum_{j=1}^n ds_j {\partial f\over \partial s_j}.$$ So we have in an
arbitrary open set $U\subseteq X$ the map $$d: \Omega^0(U,E)
\rightarrow \Omega^1(U,E)$$ has zero ${\Bbb Z}/2{\Bbb Z}$-graded degree.

Remark that $\Omega(U,E)$ is ${\Bbb Z} \oplus {\Bbb Z}/2{\Bbb Z}$-graded
and $\End \Omega(U,E)$ is also ${\Bbb Z} \oplus {\Bbb Z}/2{\Bbb Z}$-graded
.  Thus an element $u \in \End \Omega(u,E)$ is of bidegree $(c,j)$ iff $$u(
\Omega^b(U,E)_i) \subseteq \Omega^{b+c}(U,E)_{i+j},$$ for any $(b,i) \in
{ \Bbb Z} \oplus {\Bbb Z}/2{\Bbb Z}$. We will say about $u$ as a
derivation of bidegree $(c,j)$ iff for each $\alpha \in \Omega^b(U,E)_i$
and $\beta\in \Omega^{b'}(U,E')_j$, and if $E$ and $E'$ are paired with
values in $F$, we have $$u(\alpha \beta) = u(\alpha)\beta + (-1)^{bc +
ij}\alpha u(\beta).$$

Remark that if a derivation is defined, it has also functorial property.

Let us now denote by $\nabla$ the affine connection of our super-bundle $E$
, and by $\alpha = \alpha_{\nabla}$ its connection form. Let $\beta \in
\Omega^b(U,E)$, then we define the differential $d_{\nabla}$, and the
inner product $\imath(\xi)$ and Lie derivative $\theta(\xi)$ for any
vector super-field $\xi \in \Der  A(U,E)$ by the following formulae
$$\begin{array}{rl} \langle \xi,\dots,\xi_{b+1}|d_{\nabla}\beta\rangle  &=
\sum_{i=1}^{b+1}
(-1)^{i-1+j_{i-1}|\xi_i|}\nabla_{\xi_i}\langle \xi_1,\dots,\check{\xi}_i,
\dots,\xi_{b+1}|\beta\rangle  \cr &\quad + \sum_{k\langle l} (-1)^{d_{k,l}}\langle [\xi_k,
\xi_l],\xi_1,\dots,\check{\xi}_k,\dots,\check{\xi}_l,\dots,\xi_{b+1}
|\beta\rangle , \end{array}$$ where $$j_i:= \sum_{k=1}^i |\xi_k|, d_{k,l}:=
|\xi_k|j_{k-1} + |\xi_l|j_{l-1} + |\xi_k||\xi_l| + k + l, $$
$$\langle \xi_1,\dots,\xi_{b-1}|\imath(\xi)\beta\rangle  = (-1)^{|\xi|\sum_{i=1}^b
|\xi_i|} \langle \xi,\xi_1,\dots,\xi_{b-1}|\beta\rangle ,$$ and finally we define
 $$\theta(\xi) = \theta_{\nabla}(\xi):= d_{\nabla}\circ \imath(\xi) +
\imath(\xi)\circ d_{\nabla}.$$

It is easy to see that $d_{\nabla}$, $\imath(\xi)$, and
$\theta_{\nabla}$ are the derivations of degrees $(1,0)$,
$(-1,|\xi|)$, $(0,|\xi|)$, respectively. As in the ordinary case,
we have the superb-racket relations $$\imath(\xi)\imath(\eta) +
(-1)^{|\xi||\eta|}\imath(\eta)\imath(\xi) = [\imath(\xi),\imath(\eta)]
= 0, \leqno{(1)} $$ $$\theta(\xi)\imath(\eta) -
(-1)^{|\xi|.|\eta|}\imath(\eta)\theta(\xi) =
[\theta(\xi),\imath(\eta)] = \imath([\xi,\eta]), \leqno{(2)}$$
$$\theta(\xi)\theta(\eta)-(-1)^{|\xi|.|\eta|}\theta(\eta)\theta(\xi) =
[\theta(\xi),\theta(\eta)] = \theta([\xi,\eta]), \leqno{(3)} $$ and
the relationship between contraction and Lie derivation $$\begin{array}{rl}
\xi\langle \xi_1, \dots,\xi_b|\beta\rangle  &= \sum_{i=1}^b
(-1)^{|\xi|\sum_{k=1}^{i-1} |\xi_k|}\langle \xi_1,
\dots,[\xi,\xi_i],\dots,\xi_b|\beta\rangle  \cr &\quad
+(-1)^{|\xi|\sum_{k=1}^b |\xi_k|} \langle \xi_1,\dots,\xi_b|\theta(\xi)\beta\rangle
. \end{array}$$ So we have a de Rham complex of global sections of $E$,
$$\dots \rightarrow \omega^b(X,E) {\buildrel d_{\nabla} \over
\longrightarrow} \Omega^{b+1}(X,E) \rightarrow \dots
\leqno{\Omega^*(X,E):}$$ if the connection $\nabla$ is flat, i.e. for all
$\xi,\eta \in \Der  A(U)$, $$\curv(\nabla)(\xi,\eta) =
[\nabla_{\xi},\nabla_{\eta}] - \nabla_{[\xi,\eta]} \equiv 0. $$

Let $p\in U \subseteq X$, $T_p(X,A)$, the tangent super-space at $p$
, $\Omega^b_E(p)$ be the linear super-space of all $E(U)$-valued
$b$-linear forms $z$ on $T_p(X,A) = T_p(X) \oplus T_p(X,A)_1$ such
that $$\langle v_1,\dots,v_j,v_{j+1},\dots,v_b|z\rangle  = (-1)^{1 + |v_j||v_{j+1}|}
\langle v_1,\dots,v_{j+1},v_j,\dots,v_b|z\rangle .$$ Note that $\Omega^b_E(p)$ is
${\Bbb Z}/2{\Bbb Z}$-graded such that if $z$ is homogeneous, the
expression $\langle v_1,\dots,v_b|z\rangle $ vanishes unless the case $$|z| =
\sum_{i=1}^b |v_i|.$$ So $z|_{T_p(X)}$ is an $E(U)_0$-valued form on
$T_p(X)$ and $z|_{T_p(X,A)_1}$ is a symplectic $b$-form on
$T_p(X,A)_1$. We define $$\Omega^0_E(p):= E(U)_0,$$ $$\Omega_E(p)
:= \bigoplus_{b=0}^{\infty} \Omega^b_E(p).$$ We observe that the map
$A(U) \rightarrow C^{\infty}(U)
; f \mapsto \tilde{f}$ extends to a homomorphism
$\Omega(U,E) \rightarrow \Omega_E(U); \beta \mapsto \tilde{\beta} $,
$$\langle \xi_1,\dots,\xi_b|\beta\rangle ^{\sim} =
\langle \tilde{\xi}_1,\dots,\tilde{\xi}_b|\tilde{\beta}\rangle .$$

Now let $\Omega_E(X) \rightarrow \Omega_{E_0}(X) = \Hom(T_p(X),E(U)_0)$
be the restriction map from the complex $\Omega_E(X)$ to the ordinary
$E(U)_0$-valued de Rham complex, then we have a commutative diagram
$$\vbox{\halign{ #&#&#&#&# \cr
$\Omega($&$X,E)$&$\enskip\longrightarrow$&$\Omega_{E_0}$&$(X)$ \cr
&${\sim}\searrow $ & &${\nearrow}res $ & \cr & &$\Omega_E(X)$& &.
\cr}}$$ Remark that $$k: \Omega(X,E) \longrightarrow \Omega_{E_0}(X)$$
commute with $d_{\nabla}$ and suppose that the connection is flat, we
have then really a commutative diagram of complexes and complex
homomorphisms $$\vbox{\halign{ #&#&#&#&# \cr
$(\Omega(X,E)$&$,d_{\nabla})$ &$ \quad{\buildrel k \over \longrightarrow} $
&$(\Omega_{E_0}(X)$&$,d_{\nabla})$ \cr &$ \searrow $ & &$ \nearrow $ &
\cr & &$(\Omega_E(X),d_{\nabla})$ & &. \cr}}$$

\begin{lem}[\bf Poincar\'e Lemma]
Suppose that $f \in E(U) = \Omega^0(U,E)$. Then $d_{\nabla}f =
0$ in the connected open set $U$ iff $f = \lambda.1_U$, where
$\lambda$ is a constant super-function with a single value in $\End
E(pt)$.  If, in addition, $U$ is contractible $A$-coordinate
neighborhood and $\beta \in \Omega^b(U,E)$, $d_{\nabla}\beta = 0 $,
then there exists $\omega \in \Omega^{b-1}(U,E)$ such that $\beta =
d_{\nabla}\omega$. \end{lem}
\begin{pf} 
By assumption, our connection is flat. Then locally we can
consider its trivial form of connection. So, $d_{\nabla} = d$, $$
df = \sum_{i=1}^m dr_i {\partial f \over \partial r_i} + \sum_{j=1}^n
ds_j {\partial f \over \partial s_j}$$ and because $$f = \sum_{\mu}
f_{\mu}s_{\mu} $$ for every scalar super-function the first assertion
is trivial.

Locally, the map $$(\Omega(U,E),d_{\nabla}) {\buildrel k \over
\longrightarrow} (\Omega_{E_0}(X),d_{\nabla})$$ is an complex
isomorphism. So for a contractible open $U$, the acyclicity of the
usual de Rham complex deduces the acyclicity of our graded de Rham
complex, the lemma is therefore proved.
\end{pf}

\begin{thm} 
There is a commutative diagram of complexes and their
isomorphisms  
$$\vbox{\halign{ #&#&#&#\cr $\Coh (\Omega(X$&$,E),d_{\nabla})$&$
\qquad\quad {\buildrel \bar{k} \over \longrightarrow } $
&$\Coh (\Omega_{E_0}(X),d_{\nabla})$ \cr &${\cong}\searrow $&
&${\nearrow}\cong $ \cr & &$H(X;\End E(X)_0)$&
.\cr}}$$
\end{thm}
\begin{pf}
For the complex $\Omega(X,A)$, B. Kostant (\cite{kostant2},\S4.7) has
constructed a flasque resolution of the constant sheaf. Our complex
is its tensor product with $\End
E(pt)$. Thus we have a flasque resolution by our complex for the
constant sheaf.
\end{pf}

\begin{defn} 
Let $(E,\nabla)$ be a vector super-bundle
with connection form ${\sqrt{-1}\over\hbar}\alpha =
{\sqrt{-1}\over\hbar}\alpha_{\nabla}$ on super-manifold $(X,A)$. Then
there exists a unique differential 2-super-form $\omega\in
\Omega^2(X,\End E)$, such that $-\omega = d_{\nabla}\alpha $, i.e.
$$-\langle \xi,\eta|\omega\rangle  = \langle \eta|\xi\alpha\rangle  -
(-1)^{|\xi|.|\eta|}\langle \xi|\eta\alpha\rangle  - \langle [\xi,\eta]|\alpha\rangle + $$ $$
+{\sqrt{-1}\over \hbar}[\langle \xi|\alpha\rangle ,\langle \eta|\alpha\rangle ].$$ This
2-super-form is called the {\it curvature form } 
of the connection $\nabla$.
\end{defn}

As in the classical case, it is easy to deduce the following result

\begin{prop}
$$\langle \xi,\eta|\curv(E,\nabla)\rangle  = [\nabla_{\xi},\nabla_{\eta}] -
\nabla_{[\xi,\eta]} = -{\sqrt{-1}\over \hbar}\langle \xi,\eta|\omega\rangle .$$
\end{prop}

\begin{defn} 
Let $(E,\nabla)$ be a
super-bundle with connection over $(X,A)$ and $\{ (U_i,t_{ij})\}$,
$i\in I$, $j= 1,\dots,\rank _AE$, is a local system for $E$. We denote
$t_i = (t_{i1},\dots,t_{i,\rank _AE})$ and $c_{ij} \in \End E(U)$ the
transition functions defined by $t_jc_{ij} = t_j$, and we will then
refer to the set $(c_{ij},\alpha_i)$; $$\alpha_j -
c_{ij}\alpha_ic_{ij}^{-1} = {\hbar\over\sqrt{-1}} dc_{ij}.c^{-1}_{ij}
$$ as {\it local datum } \index{local datum} 
for $(E,\nabla)$.
\end{defn}

If $(c_{ij},\alpha'_i)$ is an another local data of some vector
super-bundle $(E',\nabla')$, then $(E,\nabla)$ is equivalent to
$(E',\nabla')$ if and only if there exist $\lambda_i \in Iso(E,E')$
such that $$\lambda_ic_{ij}\lambda_j^{-1} = c_{ij}' $$ and $$\alpha_i'
- \lambda_i\alpha_i\lambda_i^{-1} = {\hbar\over\sqrt{-1}}
d\lambda_i.\lambda_i^{-1}. $$ Since every vector super-bundle with
connection admits at least one local datum with respect to a
contractible covering, it follows that the notion of curvature is an
equivalence invariant and hence $\curv[(E,\nabla)]:=
[\curv(E,\nabla)]$ is well defined. Note that the set ${\cal L}
_c(X,A)$ of all equivalent classes of vector super-bundles with connection 
has the structure of an Abelian group. $$[(E,\nabla)] = [(E',\nabla')] +
[(E'',\nabla'')],$$ iff $$c_{ij} = c_{ij}'c_{ij}'',$$ $$\alpha_i =
\alpha'_i + \alpha''_i,$$ and $$\curv[(E,\nabla)] = \curv[(E',\nabla')] +
\curv[(E'',\nabla'')].$$

Now for any closed 2-super-form $\omega \in \Omega^2_c(X,\End E(U)_0)$, let
${\cal L} _{\omega}(X,A)$ be the set of all classes $[(E,\nabla)] \in
{\cal L} _c(X,A)$ such that $\omega = \curv[(E,\nabla)]$, $${\cal
L}_c(X,A) = \cup_{\omega} {\cal L} _{\omega}(X,A) $$ is a disjoint union
over the set of all closed 2-super-forms $\omega \in \Omega^2(X,\End
E(U)_0).$

Now given a closed 2-super-form $\omega \in \Omega^2(X,\End E(U))$ the
question is to decide whether ${\cal L} _{\omega}(X,A)$ is nonempty.
Observe that one has the same answer as in the ungraded case.

If $\Aut  E(pt)$ is the automorphism group of fiber transformations. Then
the cohomology group $H^1(X;\Aut  E(X)_0)$ operates on ${\cal L} _c(X,A)$.

Let $\{ U_i\}_{i\in I}$ be the contractible covering of $X$ and assume
that $(E,\nabla)$ is a vector super-bundle with connection over $(X,A)$.
Let $(c_{ij},\alpha_i)$ be the corresponding local datum for $(E,\nabla)$.
Let $z_{ij}$ be the cocycle for the constant sheaf $\Aut  E(pt)_0$.

\begin{lem}  The cohomology group $H^1(X,\Aut  E(X)_0)$ operates as
follows: $$[z_{ij}].[(E,\nabla)] = [(E',\nabla')], $$ where
$(E',\nabla')$ has the local datum $(c_{ij},z_{ij},\alpha'_i)$ with
respect to the covering $ \{ U_i\}_{i\in I}$, $\alpha'_i =
z_{ij}\alpha_iz_{ij}^{-1}$.
\end{lem} 
\begin{pf} 
It is easy to see that
$(c_{ij}z_{ij},z_{ij}\alpha_iz_{ij}^{-1})$ is a local datum of some
super-bundle with connection $(E',\nabla')$. We must prove that
$\curv(E',\nabla') = \omega$. We see that $(z_{ij},0)$ is also a local
datum, then there exists a flat super-bundle with connection
$(E_z,\nabla_0)$ such that $(z_{ij},0)$ is its local datum. By the
Abelian group structure on ${\cal L}_c(X,A)$, we have $$\begin{array}{rl}
\curv(E',\nabla') &= \curv(E,\nabla) + \curv(E_z,\nabla_0) \cr &= \omega + 0
\cr &= \omega.\end{array}$$
\end{pf}

\begin{cor}
${\cal L}_{\omega}(X,A) \cong H^1(X;\Aut  E(X)_0).$
\end{cor}

\begin{rem} Let us denote by $\exp$ the exponential map $$f\in \End E(U)_0
\mapsto \exp{({\sqrt{-1}\over\hbar}f)} \in \Aut  E(U), $$ where as usually
$$\exp{({\sqrt{-1}\over\hbar}f)}:= \sum_{n=0}^{\infty}
{({\sqrt{-1}\over\hbar}f)^n \over n!}. $$ The exponential series
converges absolutely on the operator norm topology, as usually.
\end{rem}

Remark that elements of the form $I + \End E^1(U)_0$, the unipotents,
have the unique logarithms in $\End E^1(U)_0$, hence there is an
isomorphism $$\End E^1(U)_0 \cong I_U + \End E^1(U)_0.$$

Let us denote by $\Gamma$ the kernel of the exponential map, one has
an exact sequence of group sheaves $$0 \rightarrow (\End E)_0 {\buildrel
exp \over \longrightarrow} \Aut  E_0 \rightarrow 1. $$ So by the long
exact cohomology sequence we have 

\begin{cor}  For a fixed super-manifold $(X,A)$ there is a natural
isomorphism $${\cal L} _c(X,A) \cong H^2(X;\Gamma).$$  \end{cor}

\begin{rem} By the well-known Kuiper theorem, $\Aut  E(U)$ is homotopically
trivial in the infinite-dimensional case. We are interesting in the
multidimensional quantization theory, to the classification of all
infinite-dimensional Hilbert super-bundles, associated with principal
super-bundles with finite dimensional Lie super-groups as structural groups
.  \end{rem}

\begin{lem} Let $(X,A)$ be a super-manifold and let $\omega\in
\Omega^2_{Bbb C}( X,\End E)_0$ be a closed 2-form. Then ${\cal L}
_{\omega}(X,A)$ is nonempty if and only if the cohomology class $[\omega]$
belong to the cohomology group $H^2(X;\Gamma) \hookrightarrow H^2(X,\Aut
E(X)_0) $.
\end{lem} 
\begin{pf} 
Assume that the class $[\omega]$ is
$\Gamma$-valued. Let $\{U_i\}_{i\in I}$ is a contractible covering of $X$
. By the Poincar\'e lemma, there exists $\alpha_i\in \Omega^1(U_i,E)$
such that $d\alpha_i = \omega|_{U_i}$. Hence in the intersection $U_i
\cap U_j$, $d(\alpha_i - \alpha_j) \equiv 0 $. Thus there exists $f_{ij}
\in \Aut  E(U_i \cap U_j)$ such that $\alpha_i - \alpha_j = df_{ij}$. Then
, in the intersection $U_i \cap U_j \cap U_k$, $d(f_{ij} + f_{jk} -
f_{ik}) \equiv 0 $. So there is some $z_{ijk}\in \Aut  E(X)_0 $ such that
$$f_{ij} + f_{jk} - f_{ik} = z_{ijk}Id_{U_i \cap U_j \cap U_k}. $$
Because $\omega$ is $\Gamma$-valued, we can choose $z_{ijk} \in \Gamma$.
So, $$c_{ij}:= \exp{(f_{ij})} \in \Aut  E(U)_0.$$ It is easy to see that
$( c_{ij},\alpha_i)$ is a local datum of some super-bundle with connection
, say $(E,\nabla)$. Clearly, $\omega = \curv(E,\nabla)$. This means that
$ {\cal L} _{\omega}(X,A)$ is nonempty.

The inverse assertion is proved in the same way as in the classical case.
\end{pf}

\subsection{Quantization super-operators }

{\bf Notation}
Let $(M,A;\omega)$ be a symplectic super-manifold of dimension $(m,n)$.
recall from 4.4 the map $k: (\Omega(X,E),d_{\nabla}) \rightarrow
(\Omega_{E_0}(X),d_{\nabla})$. Then $(X,k\omega)$ is a
symplectic manifold and the dimension $m$ must be even, $m = 2m_0$
and there is at every point a Darboux coordinate system $(p_i,q_i)$
in coordinate neighborhood $U$. In other hand, because the restriction
of $\omega$ to the odd part of the tangent super-spaces is a non-degenerate
symmetric form, there exists ( see also \cite{kostant2},\S5.3) a Morse canonical
 system of coordinates $s_j,j= 1,\dots,n$. Together, there exists at
 every point a so called {\it local $A$-Darboux coordinate system }
$(\{(p_i,q_i)\}_{i=1,\dots,m_0},\{s_j\}_{j=1,\dots,n})$, where $p_i,
q_i \in A(U)_0$, and $s_j\in A(U)_1$, and in which $$\omega =
\sum_{i=1}^{m_0} dp_i \wedge dq_i + \sum_{i=1}^n {\varepsilon_j \over 2}
(ds_j)^2,$$ where $\varepsilon_j = \pm 1$.

It is easy to see that in this $A$-Darboux coordinate system, the Hamiltonian
vector super-field $\xi_f$ corresponding to a super-function $f$; $\imath(\xi_f
)\omega = df$ now becomes $$\xi_f = \sum_{i=1}^{m_0} ({\partial f \over 
\partial q_i}{\partial \over \partial p_i} - {\partial f \over \partial p_i}
{\partial \over \partial q_i}) + \sum_{j=1}^n \varepsilon_j{\partial f \over 
\partial s_j}{\partial \over \partial s_j}.$$ In particular, we have $$
\xi_{p_k} = -{\partial \over \partial q_k}, \xi_{q_k} = {\partial \over
\partial p_k}, $$  $$\xi_{s_j} = \varepsilon_j{\partial \over \partial s_j}.
$$ In such a $A$-Darboux coordinate system the Poisson brackets have the form 
$$\begin{array}{rl} f,g \in A(U) \mapsto \{f,g\} &= \xi_f(g) = -(-1)^{|f|.|g|}
\xi_g(f) = \langle \xi_f,\xi_g|\omega\rangle  \cr  &= \sum_{i=1}^{m_0}({\partial f\over
\partial q_i}{\partial g \over \partial p_i} - {\partial f \over \partial 
p_i}{\partial g \over \partial q_i}) + \sum_{j=1}^n (-1)^{|f|} \varepsilon_j
{\partial f\over \partial s_j}{\partial g\over \partial s_j}.  \end{array}$$
In particular we have the Poisson brackets between the generators as the 
canonical super-commutator relations $$\{p_i,p_k\} =\{q_i,q_l\} =\{s_j,p_k\}
= \{s_j,q_l\} = 0,$$ $$\{q_k,p_l\} = \delta_{kl}I_U,$$ $$\{s_i,s_j\} =
\varepsilon_j\delta_{ij}I_U.$$

Now we consider the quantization problem. As in the classical case,
  by definition, a {\it quantization procedure }
is a correspondence associating to each super-function $f$ a super-symmetric 
super-operator $\hat{f}$, which is anti-auto-adjoint super-operator if $f$ is
a real super-function,
in some fixed Hilbert super-space, such that $$\vbox{\halign{ #&# \cr
$\widehat{\{f_1,f_2\}} \quad$&$= \quad {\sqrt{-1}\over \hbar}[\hat{f_1},
\hat{f_2}]$ \cr $\quad\hat{1} \quad$&$= \quad Id,$ \cr}}$$ where as usually
$\hbar = {h\over 2\pi}$ is the normed Planck's constant.

\begin{defn} 
Let us denote by $(E,\nabla)$ a vector super-bundle with connection and
Hilbert super-space fibers, such that the connection conserves the
Hilbert structure. This means that the corresponding parallel
transpose conserves the scalar product. Let us denote also by
${\sqrt{-1}\over\hbar}\alpha(.)$ the connection super-form of $\nabla$
. Then the values of the super-form $\alpha$ are the anti-super-symmetric
super-operators for any complex vector super-fields and are
anti-auto-adjoint for real super-fields.

We define now for each super-function $f\in A(M)$ the corresponding 
{\it quantized super-operator }
$\hat{f}$, $$\hat{f} = f + {\hbar \over \sqrt{-1}}\nabla_{\xi_f} =
f + {\hbar\over\sqrt{-1}}\theta_{\xi_f} + \alpha(\vert(\xi_f)).$$
\end{defn}

The same as in the classical case, we have also the following
quantization condition.

\begin{thm}
The following three conditions are equivalent.

(1) The super-operator-valued differential 2-super-form $-\omega$ is the
differential of $\alpha$ with respect to $d_{\nabla}$, i.e.
$$\begin{array}{rl} \langle \xi,\eta|d_{\nabla}\alpha\rangle  &= \langle \eta|\xi\alpha\rangle  -
(-1)^{|\xi|.|\alpha|}\langle \xi| \eta\alpha\rangle  -\langle [\xi,\eta]|\alpha\rangle   +
\cr &+{\sqrt{-1}\over \hbar}[\langle \xi|\alpha\rangle ,\langle \eta| \alpha\rangle ] \cr
&=-\langle\xi,\eta|\omega\rangle Id, \forall \xi,\eta \in \Der  A(M). \end{array}$$

(2) The curvature of $\nabla$ is symplectic, more precisely $$\langle \xi,\eta|
\curv(E,\nabla)\rangle := [\nabla_{\xi},\nabla_{\eta}] - \nabla_{[\xi,\eta]} =
-{\sqrt{-1}\over \hbar}\langle \xi,\eta|\omega\rangle .$$

(3) The correspondence $f \mapsto \hat{f}$ is a quantization
correspondence.
\end{thm}
\begin{pf}
$(1) \Longleftrightarrow (2)$. As in the classical case, we have
$$\nabla_{\xi} = \theta_{\xi} + {\sqrt{-1}\over\hbar}\langle \xi|\alpha\rangle $$ and
therefore $$\begin{array}{rl}
[\nabla_{\xi},\nabla_{\eta}]-\nabla_{[\xi,\eta]} &= [\theta(\xi)
+{\sqrt{-1}\over\hbar}\langle \xi|\alpha\rangle ,\theta(\eta)+{\sqrt{-1}
\over\hbar}\langle \eta|\alpha\rangle ] - \theta([\xi,\eta])\\  & - {\sqrt{-1}\over\hbar}
\langle [\xi,\eta]|\alpha\rangle  \cr &=[\theta(\xi),\theta(\eta)]
+{\sqrt{-1}\over\hbar} [\theta(\xi),\langle \eta|\alpha\rangle ] -
{\sqrt{-1}\over\hbar}(-1)^{|\xi|.|\eta|}[\theta( \eta),\langle \xi|\alpha\rangle ] \cr
&+ ({\sqrt{-1}\over\hbar})^2[\langle \xi|\alpha\rangle ,\langle \eta|\alpha\rangle ] -
\theta([\xi,\eta]) - {\sqrt{-1}\over\hbar}\langle [\xi,\eta]|\alpha\rangle  \cr &=
[\theta(\xi),\theta(\eta)] - \theta([\xi,\eta]) + {\sqrt{-1}\over \hbar}
\{
[\theta(\xi),\langle \eta|\alpha\rangle ]-\\
        &-(-1)^{|\xi|.|\eta|}[\theta(\eta),\langle \xi|\alpha\rangle ]
\quad -\langle [\xi,\eta]|\alpha\rangle  +
{\sqrt{-1}\over\hbar}[\langle \xi|\alpha\rangle ,\langle \eta| \alpha\rangle ] \}.\end{array}$$ In
virtue of the relation $$[\theta(\xi),\theta(\eta)] = \theta([\xi, \eta])
$$ the rest of the proof of this part is
\end{pf}

\begin{lem}
$$[\theta(\xi),\langle \eta|\alpha\rangle ] = \langle \eta|\xi\alpha\rangle .$$
\end{lem}
\begin{pf} 
For any super-section $f\in E(U)$ we have $$\begin{array}{rl}
[\theta(\xi),\langle \eta|\alpha\rangle ]f &= \theta(\xi)(\langle \eta|\alpha\rangle f) +
(-1)^{|\xi|.|\eta|}\langle \eta|\alpha\rangle \theta(\xi)f-\\  & - (-1)^{|\xi|.|\eta|}
\langle \eta|\alpha\rangle \theta(\xi)f \cr &= \theta(\xi)\langle \eta|\alpha\rangle f.\end{array}$$
\end{pf}

$(2) \Longleftrightarrow (3)$. We have directly from the definition
$$\begin{array}{rl} {\sqrt{-1}\over\hbar} [\hat{f}_1,\hat{f}_2] &=
{\sqrt{-1}\over\hbar}[f_1 + {\hbar\over\sqrt{-1}} \nabla_{\xi_{f_1}},f_2 +
{\hbar\over\sqrt{-1}}\nabla_{\xi_{f_2}}] \cr &= {\sqrt{-1}\over\hbar}\{
{\hbar\over\sqrt{-1}}[f_1,\nabla_{\xi_{f_2}}] +
{\hbar\over\sqrt{-1}}[\nabla_{xi_{f_1}},f_2] + ({\hbar\over\sqrt{-1}})^2[
\nabla_{xi_{f_1}},\nabla_{\xi_{f_2}}] \} \cr &=[f_1,\nabla_{\xi_{f_2}}] +
{\hbar\over\sqrt{-1}}\nabla_{\xi_{\{f_1,f_2\}}} -
{\hbar\over\sqrt{-1}}\nabla_{\xi_{\{f_1,f_2\}}} + [\nabla_{\xi_{f_1}},f_2]
+\\  &+ {\hbar\over\sqrt{-1}}[\nabla_{\xi_{f_1}},\nabla_{\xi_{f_2}}].
\end{array}$$ It is easy to see that $$\xi_{\{f_1,f_2\}} =
[\xi_{f_1}.\xi_{f_2}].$$ So the
 proof of the theorem shall be completed by proving the following

\begin{lem} $$[f_1,\nabla_{\xi_{f_2}}] = [\nabla_{\xi_{f_1}},f_2] =
\theta(\xi_{f_1})f_2
=-(-1)^{|\xi_{f_1}|.|\xi_{f_2}|}\theta(\xi_{f_2})(f_1) =$$ $$=
\langle \xi_{f_1},\xi_{f_2}| \omega\rangle  = \{f_1,f_2\}.$$
\end{lem} 
\begin{pf} 
The multiplication by a super-function is super-commuting with the
multiplication by any super-operator-valued super-function. So we have
$$[f_1,\nabla_{\xi_{f_2}}] = [f_1,\theta(\xi_{f_2})].$$ For every section
$s \in E(U)$, one has $$\begin{array}{rl} [f_1,\theta(\xi_{f_2})]s &=
f_1(\theta(\xi_{f_2})s) -(-1)^{|f_1|.|f_2|}\theta(\xi_{f_2})(f_1s) \cr &=
f_1(\theta(\xi_{f_2})s)-(-1)^{|\xi_{f_2}|.|f_1|}(\theta(\xi_{f_2})f_1)s -
(-1)^{2|\xi_{f_2}.|f_1|}f_1(\theta(\xi_{f_2})s) \cr &=
-(-1)^{|\xi_{f_2}|.|f_1|}(\theta(\xi_{f_2})f_1)s \cr &=
\langle \xi_{f_1},\xi_{f_2}|\omega\rangle s. \end{array}$$ By analogy, we have also $$
[\nabla_{\xi_{f_1}},f_2] = \theta(\xi_{f_1})f_2 =
\langle \xi_{f_1},\xi_{f_2}|\omega\rangle .$$
\end{pf}

\begin{rem} Here we talk about a prequantization procedure, because as
usually the covariant derivation $\nabla_{\xi}, \xi \in \Der  A(U)$
operates as a element of $\End E(U)$ for every open $U$. Only with the
help of some super-polarization, we can construct the corresponding Hilbert
super-space of the quantum states of the system under consideration.  We
shall do this as super-polarization in what follows. \end{rem}

\begin{rem} Remember that the expression $$d_{\nabla}\alpha = d\alpha +
{\hbar \over 2\sqrt{-1}}\alpha \wedge \alpha$$ is non linear with respect
to $\alpha$. Thus in the condition (1) of the theorem we have a nonlinear
differential equation $$-\langle \xi,\eta|d_{\nabla}\alpha\rangle  = \langle \xi,\eta|\omega\rangle
Id $$ and in the condition (2) of the theorem, we have an equivalent
nonlinear differential equation $$\curv(E,\nabla) = -{\sqrt{-1}\over
\hbar}\langle .,.|\omega\rangle id. $$ In the temporary physics there are some sorts of
this equation, related with the  curvature of some connections. An
interesting example of this type is the well-known Yang-Mills equation
$${^*}\curv(E,\nabla) = \pm \curv(E,\nabla),$$ where the star sign is the
Hodge star. By Radon-Penrose transform this equation becomes the
Cauchy-Riemann condition for the corresponding bundle in the space of
light rays of the compactified Minkowski space.  \end{rem}

In the work \cite{diep22} it was posed the following 

{\bf Conjecture. } Our quantization conditions could be also transformed
by the Radon-Penrose transformation to some algebraic condition of the
corresponding super-bundle.  \vskip 2truecm

We discuss now about application of this (pre)quantization procedure to
the representation theory of Lie super-groups.

{\bf Notations } Recall that if $G$ is any group and $K$ is a fixed
ground field (${\Bbb R}$ or ${\Bbb C}$), the group algebra $K(G)$ is a
commutative Hopf algebra with antipode over $K$, $$\Delta: K(G)
\rightarrow K(G) \otimes K(G),$$ so that $$\Delta(g) = g \otimes g, s(g)
= g^{-1}, 1_E(g) = 1,\forall g\in G.$$

Now assume that ${\frakt g}$ is a Lie super-algebra and one has a
representation $$\Pi: G \longrightarrow \Aut  {\frakt g}.$$ Then $\Pi$
extends uniquely to a representation of $G$ by automorphisms of the
universal enveloping super-algebra $U({\frakt g})$ Now the {\it smash
product $K(G) \circledS U({\frakt g})$ \index{product!smash -}with 
respect to $\Pi$, } or 
simply {\it smash product } if $\Pi$ is understood, is by definition a
co-commutative Hopf super-algebra with antipode such that:

(1) as a vector super-space, it is the graded tensor product $K(G) \hat{
\otimes} U({\frakt g})$,

(2) as an algebra $K(G)$ and $U({\frakt g})$ are sub-algebras with the
relations $$gug^{-1} = \Pi(g)u, \forall g\in G, \forall u\in U({\frakt
g}),$$

(3) with respect to the diagonal map $\Delta$ the elements of $G$ are {\it
group-like } \index{element!group-like -} and the elements of ${\frakt g}$ 
are {\it primitive }\index{element!primitive -} and finally

(4) one has $s(g) = g^{-1}$, $s(X) = -X$ and $1_E(g) = 1$, $1_E(X) = 0$,
$\forall g\in G$, $\forall X \in {\frakt g}$.

Conversely, Let $E$ be any commutative Hopf super-algebra with antipode
over an algebraically closed field $K$ of characteristic zero, $G$ the
group of all group-like elements in $E$ and ${\frakt g}$ be the Lie
super-algebra of all primitive elements in $E$. Then there is a
representation say $\Pi: G \rightarrow \Aut  {\frakt g}$ such that
$$gXg^{-1} = \Pi(g)X,\forall g\in G, \forall X \in U({\frakt g})$$ and
$E$ is isomorphic to the smash product of these parts with respect to the
representation $\Pi$, $E \cong K(G) \circledS U({\frakt g})$.

{\bf Notations} Now assume that $G$ is a group, ${\frakt g}$ is a Lie
super-algebra over the ground field $K$ and $E$ is the Hopf super-algebra
with respect to some representation $\Pi: G \rightarrow \Aut  {\frakt g}$.
Following B. Kostant \cite{kostant2},\S3.4, we will say that $E$ has the structure
of a {\it Lie-Hopf super-algebra } \index{Lie-Hopf super-algebra} iff:

(1) $G$ has a structure of ( not necessarily connected ) Lie group,

(2) ${\frakt g = g}_{\bar 0} + {\frakt g}_{\bar 1}$ is a finite
dimensional Lie super-algebra such that ${\frakt g}_{\bar 0} \cong \Lie G$
is the Lie algebra of $G$,

(3) $Ad_{\frakt g}$ is defined on the identity component $G_0$ of $G$ and
 $\Pi|_{G_0} \cong Ad_{\frakt g}$.

Denote in this case $E:= E(G,{\frakt g})$ this Lie-Hopf super-algebra.

If $E(G,{\frakt g})$ is a Lie-Hopf super-algebra, denote $E(G,{\frakt g}_{
\bar 0})$ the Lie-Hopf algebra obtained by replacing ${\frakt g}$ by its
even part ${\frakt g}_0$. As Hopf algebras one knows that $$E(G,{\frakt
g}_0) \cong C^{\infty}(G)^*,$$ where $C^{\infty}(G)^*$ is the set of
distributions with finite support.

The Lie-Hopf algebras form a category in which a morphism $$E(G,{\frakt g})
\rightarrow E(H,{\frakt h})$$ is a morphism of Hopf super-algebras such that
 the restriction to the even parts is a Hopf algebra morphism
$$E(G,{\frakt g}_{\bar 0}) \rightarrow E(H,{\frakt h}_{\bar 0}),$$
induced by a Lie group
 morphism $G \rightarrow H$.

{\bf Notations} Let us consider a super-manifold $(Y,B)$ with the sheaf
$B(Y)$ of super-functions on $(Y,B)$. Consider the full dual $B(Y)'$ of
$B(Y)$. One has certainly an injection $$0 \rightarrow B(Y)' \otimes
B(Y)' \rightarrow (B(Y) \otimes B(Y))'.$$ One has also the diagonal map
$$\Delta: B(Y)' \rightarrow (B(Y) \otimes B(Y))'$$ defined by the
relation $$v(f \otimes g) = v(fg), \forall f,g \in B(Y),$$ and $v\in
B(Y)'$.

We consider the subspace $B(Y)^*$ defined as the super-subspace of $B(Y)'$
, consisting of all $v\in B(Y)'$, which vanishes on some ideal of finite
co-dimension of $B(Y)$. One knows that if $v\in B(Y)'$, $\Delta v \in
B(Y)^* \otimes B(Y)^*$, iff $v\in B(Y)^*$. So there is a natural
morphism $$\Delta: B(Y)^* \rightarrow B(Y)^* \otimes B(Y)^*,$$ which
induces on $B(Y)^*$ a co-commutative super-algebra structure.

Recall that if $$\Delta: C \rightarrow C \otimes C $$ is a super-algebra
co-multiplication, an element $\delta \in C$ is called {\it group-like } 
\index{element!group-like -}  iff it is non zero even element and $$\Delta 
\delta = \delta \otimes
\delta.$$ An element $v$ is called {\it primitive } 
\index{element!primitive -} with respect to a
group-like element $\delta$ iff $$\Delta v = \delta \otimes v + v \otimes
\delta.$$ Remark that $B(Y)^*$ is just the set of all distributions with
finite support on $(Y,B)$.

Recall also that a morphism $$\tau: B(Y)^* \longrightarrow C(Z)^*$$ of
super-co-algebras is said to be {\it smooth } \index{supermorphism!smooth -} 
iff it is an induced morphism
$\tau = \sigma_*$ for some morphism of super-manifolds $$\sigma: (Y,B)
\rightarrow (Z,C).$$

Now let $(G,A)$ be a super-manifold of dimension $(m,n)$, $$\Delta:
A(G)^* \rightarrow A(G)^* \otimes A(G)^* $$ be the diagonal map with
respect to which $A(G)^*$ is cocommmutative super-algebra. The co unit is
given by the identity element $1_G \in A(G), 1_G(v) = v(1_G), \forall
v\in A(G)^*$.

Recall that $(G,A)$ has the structure of a Lie super-group if $A(G)$ has
also the structure of an algebra such that:

(1) $(A(G)^*,1_G,\Delta)$ is a Hopf super-algebra with antipode $s$,

(2) the map $$A(G)^* \otimes A(G)^* \longrightarrow A(G)^*$$ given by the
multiplication and the map $$s: A(G)^* \longrightarrow A(G)^*$$ given by
the antipode are {\it smooth. }

It is well-known ( see \cite{kostant2},\S3.5 ) that if $(G,A)$ is a Lie group,
then $G$ is a Lie group with respect to the underlying manifold structure
, with Lie algebra $\Lie G \cong {\frakt g}_{\bar 0}$, the even part of
the Lie super-algebra of primitive elements of $A(G)^*$ and $A(G)^* =
E(G,{\frakt g})$ with respect to the representation $$\Pi: G
\longrightarrow \Aut  {\frakt g},$$
 $$\Pi(g)x {\buildrel \rm def \over =} gxg^{-1},\forall x\in {\frakt g},
\forall g\in G $$ and $\Pi|_{G_0} \cong \Ad_{\frakt g}$, the restriction
to the connected component of identity.

\begin{defn} Let $(G,A)$ be a Lie super-group and let $(Y,B)$ be a
super-manifold. We will say that $(G,A)$ {\it acts }
 on $(Y,B)$ or $(Y,B)$ is a $(G,A)$-{\it space } iff the following map is
smooth $$\Delta: A(G)^* \times B(Y)^* \rightarrow B(Y)^*,$$ $$\Delta(vw)
= \sum_{i,j} (-1)^{|w'_j|.|u''_i|} u'_iw'_j \otimes u''_iw''_j,$$ if
$$\Delta(u) = \sum_i u'_i \otimes u''_i,$$ and $$ \Delta(w) = \sum_j w'_j
\otimes w''_j,$$ for each $u\in A(G)^*$ and $w\in B(Y)^*$.  \end{defn}

In this case $B(Y)^*$ becomes an $A(G)^*$-module. By duality, the
commutative Lie super-algebra $B(Y)$ becomes also some $A(G)^*$-module:
$$\langle w,u.f\rangle  = (-1)^{|u|.|w|}\langle s(u)w,f\rangle .$$ Observe that if $f,g \in B(Y)$ and
$$\Delta(u) = \sum_i u'_i \otimes u''_i,$$ then $$u.(fg) = \sum_i (-1)^{
|f|.|u''_i|} (u'_i.f)(u''_i.g).$$

Recall that a Lie super-subgroup $(H,B)$ in $(G,A)$ will be called {\it
closed } if $H$ is a closed subgroup in $G$. Let $$\vbox{\halign{
#&#&#&#\cr $\rho:$&$G$&$ \rightarrow $&$H\setminus G,$ \cr &$g$&$
\mapsto $&$Hg $ \cr}}$$ be the coset projection map. Put $V =
\rho^{-1}(U)$ if $U$ is some open subset in $H\setminus G$. Then
$(V,A(V))$ is some $(H,B)$-super-space and the restriction map $$\rho:
A(G) \rightarrow A(V)$$ is some $(H,B)$-module map. Now put $$(B\setminus
A)(U) = \{ f\in A(V); R^V_Wf = (-1)^{|w|.|f|}fs(w) = \langle w,1_H\rangle f \}.$$ It is
easy to see that $(B\setminus A)(U)$ is a commutative super-algebra
containing in $A(V)$ and the correspondence $$U \mapsto (B\setminus
A)(U)$$ form a sheaf of commutative super-algebras. As it is pointed out
in (\cite{kostant2},\S3.9), the sheaf $B\setminus A$ on $H\setminus G$
together with the homomorphism $$(B\setminus A)(U) \rightarrow
C^{\infty}(U)$$ defines
 a super-manifold structure $(H\setminus G,B\setminus A)$ of dimension $(m -
m',n - n')$, if $\dim{(G,A)} = (m,n)$ and $\dim{(H,B)} = (m',n')$.
Furthermore, we have also the local triviality of the projection map:
{\it For the sufficiently small open sets $U$, one has an isomorphism
$$\theta: (U \times H, (B\setminus A) \times B) {\buildrel \cong \over
\longrightarrow} (G,A).$$ In other words, we have also a principal
super-bundle associated with each closed Lie super-subgroup. }

It is not hard to show ( see also \cite{kostant2},\S3.10.3 ) that {\it
if $(H,B)$ is a closed Lie super-subgroup of a Lie super-group $(G,A)$, then
with respect to the action of $(G,A)$ on $(H\setminus G,B\setminus A)$,
$(H\setminus G,B\setminus A)$ is a homogeneous $(G,A)$-super-space. Conversely
, if $(X',A')$ is a homogeneous $(G,A)$-super-space, then $(X',A') \cong
(H\setminus G,B\setminus A)$, where $(H,B)$ is the stabilizer of a point
,say $p \in X'$. }

\begin{defn} Suppose that our Lie super-group $(G,A)$ acts on a symplectic
super-manifold $(M,B;\omega)$ by a representation $$\Pi(.): A(G)^*
\rightarrow \End B(M),$$ such that its restriction to the Lie super-algebra
${\frakt g}$ is a representation of the Lie super-algebra ${\frakt g}$ by
the canonical transformations $$X \in {\frakt g} \mapsto \xi_X \in \Ham
B(M) \subseteq \Der  B(M).$$ Denote by $L_X$ the Lie derivation along the
vector super-field $\xi_X, $ i.e. $L_X:= \theta(\xi_X)$.

We have also a natural exact sequence of Lie super-algebras $$0 \rightarrow
{\Bbb R} 1_M \rightarrow B(M) \rightarrow \Ham B(M) \rightarrow 0.$$ Hence
, for each $X \in {\frakt g}$, there exists a super-function, the so
called {\it generating function }, $f_X \in B(M)$; $$\imath(\xi_X)\omega
+ df_X = 0.$$ By the calculus on super-manifolds, we have also
$$[L_X,L_Y] = L_{[X,Y]} $$ and $$L_Xf =\{f_X, f\}.$$ \end{defn}

Now suppose that $f_X$ depends linearly on $X$, we have then a 2-cocycle
of Lie super-algebra $$c(X,Y) = \{f_X,f_Y\} - f_{[X,Y]}.$$ By the
quantization conditions, we have $$[\Lambda(X),\Lambda(Y)] =
\Lambda([X,Y]) + c(X,Y),$$ where $$\Lambda(X):=
{\sqrt{-1}\over\hbar}\hat{f}_X = {\sqrt{-1}\over\hbar} f_X +
\nabla_{\xi_X}.$$

\begin{defn} We will say that the action of $(G,A)$ by the canonical
transformations on the symplectic super-manifold $(M,B;\omega)$ is {\it
flat } iff the 2-cocycle $c(.,.)$ is zero.  \end{defn}

Remark that if the $(G,A)$-action on $(M,B)$ is flat, the Lie
super-algebra homomorphism $${\frakt g} \rightarrow \Ham_{loc} (M,B); X
\mapsto \xi_X$$ can be lifted to a super-algebra homomorphism ${\frakt g}
\rightarrow B(Y)$.  So we return to the Kirillov's notion of {\it
strictly Hamiltonian action } in the Lie group situation.

In this case, we have a Lie super-algebra representation $${\frakt g}
\rightarrow End B(Y); X \mapsto \{f_X,.\},$$ which is an {\it
integrable Poisson representation } in the Kostant's sense 
(\cite{kostant2},\S5.4).

In the Lie super-group flat action case, we have a representation of our
Lie super-algebra ${\frakt g}$ by super-functions, the so called {\it
classical ( physical ) quantities }, and some representation $\Lambda(.)$
of our Lie super-algebra by the quantum quantities $\hat{f_X}$, the
antisymmetric super-operators. If the E. Nelson conditions are satisfied,
we will have a Lie super-group representation $$\exp{X}\in G \mapsto
\exp{({\sqrt{-1}\over\hbar}\hat{f}_X)} $$ of the universal covering group
$(\tilde{G},A)$ of $(G,A)$.

\begin{defn} By a {\it mechanical system with super-symmetry } we mean a
symplectic super-manifold together with a flat homogeneous action of some
Lie super-group ( say also, of symmetry ).  \end{defn}

So starting from a mechanical system with super-symmetry we can obtain some
representations, i.e. the corresponding quantum systems, by using the
quantization procedure also as in ungraded cases.  \vskip 2truecm

\subsection{Super-polarizations and induced representations }

We see that to construct the quantum system, we must take in account a
prequantization procedure and a Hilbert super-bundle with connection
$(E,\nabla) \in {\cal L} _{\omega}(M,B)$ for the quantum model, the
sections of which are the quantum states ( with internal super-symmetry ).
The question is: Whether this set is nonempty ? As usually, the integral
condition guarantees a positive solution of this question.

\begin{defn} Let ${\frakt g}$ be a finite dimensional real Lie
super-algebra, $(G,A)$ the corresponding simply connected Lie super-group,
$(M,B)$ an arbitrary homogeneous $(G,A)$-super-space with a flat action of
$(G,A)$, $m\in M$ a fixed point and $(G_m,A_m)$ the isotropy subgroup at
this point, and finally ${\frakt g}_m:= \Lie (G_m,A_m) \subseteq {\frakt
g}$ the corresponding Lie super-algebra.  The point $m$ can be considered
as a super-functional $\delta_m \in A(G)^*$ on {\frakt g}. Therefore we
can define $$\langle m,X\rangle  = D\chi_m({\sqrt{-1}\over\hbar}f_X):=
{\sqrt{-1}\over\hbar}f_X(m).$$ In the flat action case there is some
character $\chi_m$ of the connected component of identity $((G_m)_0,A_m)$
. Then as usually, The point $m$ is said to be {\it admissible ( resp.,
integral ) } iff there is an extension of this super-character $\chi_m$
from the connected component of identity $((G_m)_0,A_m)$ to an irreducible
representation $\tilde{\sigma}$ (resp., to the same character ) of the
whole stabilizer group $(G_m,A_m)$.  It is not hard to check an
super-analogous result \end{defn}

\begin{prop} If $(M,B)$ is a homogeneous flat action $(G,A)$-super-space,
then the following conditions are equivalent:

(1) Some point $m\in M$ is admissible ( resp.,integral).

(2) Every point of $M$ is admissible ( resp., integral).

(3) The set ${\cal L}_{\omega}(M,B)$ is nonempty.
\end{prop}

Suppose now $(M,B)$ is admissible, i.e. there exists a representation
$\tilde{\sigma}$ of $(G_m,A_m)$ in some Hilbert super-space $\tilde{V}$.
Taking a fixed connection on the principal super-bundle $$(G_m,A_m)
\rightarrowtail (G,A) \twoheadrightarrow (M,B)$$, we obtain some affine
connection on the associated super-bundle. As in the classical cases the
set of these possible quantum associated bundles can be parameterized by
the first cohomology group $H^1(M;\Aut (\tilde{V})_0)$. Therefore we have
also in the super-case a result analogous to the same one in the classical
ungraded case.

Let us now consider a symplectic flat action  $(G,A)$-manifold $,
(\Omega,C; \omega_{\Omega} \in {\cal O}
(M,B):= (M,B;\omega)/(G,A)$ some $(G,A)$- orbit. It can be also
proved an analogue of the classification theorem of homogeneous flat
symplectic manifolds.  So that locally we can consider our symplectic
super-manifold as some Kirillov's co-adjoint orbit. Let us fix a point
$x$ in $\Omega$, its stabilizer $(G_x,A_x)$, the connected component
of identity $((G_x)_0,A_x)$, a fixed representation $(\tilde{\sigma},
\tilde{V})$ of $(G_x,A_x)$,etc....

Let $T(\Omega,C)$ be the complexified tangent super-bundle of the orbit
$(\Omega,C;\omega_{\Omega})$. Consider some complex tangent
distribution $L \subseteq T(\Omega,C)$, which is $(G,A)$-invariant
and integrable, i.e.  there is some sub-super-bundle, the tangent
bundle of which is just $L$.  It is not hard also to prove an
super-version of the {\it Frobenius Theorem: The tangent distribution
L is integrable iff the set of all global sections form a Lie
super-algebra with respect to the superb-rackets of super-sections. }
It is also easy to check that {\it The Lie super-algebra of invariant
global sections of an integrable invariant distribution is isomorphic
to the quotient of some complex Lie super-subalgebra by the complexified
Lie super-algebra of the stabilizer. }
It is therefore easy to
construct the complex super-subalgebra ${\frakt p} \subseteq {\frakt
g}_{\Bbb C}$ such that $L_x \cong {\frakt p}/({\frakt g}_x)_{\Bbb C}$ 
and $$L_x \oplus \overline{L}_x \cong ({\frakt p} + \overline{\frakt p})/
({\frakt g}_x)_{\Bbb C}.$$

It is clear that the invariance guarantees the inclusion $[{\frakt g}_x,
L_x] \subseteq L_x$.

\begin{defn} We will say that the distribution $L$
is {\it closed } iff the connected super-subgroup $(H_0,\tilde{F})$,
$(H_0,\tilde{I})$ corresponding to the Lie super-subalgebra ${\frakt h:=
p \cap g}$ and ${\frakt m}:= ({\frakt p} + \overline{\frakt p}) \cap
{\frakt g} $, and the super-subgroup $(H,F):= (G_x,A_x) \ltimes
(H_0,\tilde{F})$, $(M,I):= (G_x,A_x) \ltimes (M_0,\tilde{I})$ are closed.
\end{defn}

\begin{defn} Fix a representation $(\tilde{\sigma},\tilde{V})$ of the
stabilizer $(G_x,A_x)$, which is a multiple of $\chi_x$ in restricting to
the connected component $((G_x)_0,A_x)$ of the admissible orbit $(\Omega
,C;\omega_{\Omega})$.  We say that $(L,\rho,\sigma_0)$ is a
$(\chi_x,\tilde{\sigma})$-{\it polarization and L is weakly Lagrangian }
iff: (a) $\sigma_0$ is an irreducible unitary representation of the
super-group $(H_0,F)$ in some Hilbert super-space $V$ such that:
\begin{enumerate} \item[(1)] The restriction $$\sigma_0|_{(G_x,A_x)\cap
(H_0,\tilde{F})} \simeq mult\tilde{\sigma}|_{(G_x,A_x) \cap
(H_0,\tilde{F})},$$ \item[(2)] The point $\sigma_0$ is fixed under the
contragradient action of the Lie super-group $(G_x,A_x)$ in the dual
$\widehat{(H_0,\tilde{F})}$, \end{enumerate} (b) $\rho$ is a
representation of the complex Lie super-algebra ${\frakt p}$ in the Hilbert
super-space $V$ such that its restriction to the real part {\frakt h} is
the Lie derivative of the representation $\sigma_0$ and the E. Nelson
conditions are supposed to be satisfied.
\end{defn}

The notion of polarization can be also formulated in terms of Lie
super-subalgebra. That is the so called $(\tilde{\sigma},x)$-polarization.

\begin{defn} 
We say that $({\frakt p},\rho,\sigma_0)$ is a $(\tilde{\sigma},
x)$-polarization iff:

(a) ${\frakt p}$ is some complex Lie super-subalgebra of ${\frakt g}_{\Bbb
C}$, containing ${\frakt g}_x$,

(b) The super-subalgebra {\frakt p} is $\Ad_{{\frakt g}_{\Bbb
C}}(G_x)$-invariant
,

(c) The vector super-space ${\frakt p} + \overline{\frakt p}$ is the
complexification of some real super-algebra ${\frakt m} = ({\frakt p} +
\overline{\frakt p}) \cap {\frakt g}$,

(d) The super-subgroups $(M_0,\tilde{I}) $, $(H_0,\tilde{F})$,
$(M,I)$, $( H,F)$ are closed, where $(M_0,\tilde{I})$ (resp.,
$(H_0,\tilde{F})$ ) is the connected super-subgroup of $(G,A)$ with the
Lie super-subalgebra ${\frakt m}$ ( resp., ${\frakt h = p \cap g}$ )
and $(M,I):= (G_x,A_x) \ltimes (M_0,\tilde{I})$, $(H,F):= (G_x,A_x)
\ltimes (H_0,\tilde{F})$,

(e) $\sigma_0$ is an irreducible unitary representation of the super-group
$(H_0,\tilde{F})$ such that \begin{enumerate} \item[(1)] The restriction
$\sigma_0|_{(G_x,A_x) \cap (H_0,\tilde{F})}$ is a multiple of the
restriction $\tilde{\sigma}|_{(G_x,A_x) \cap (H_0,\tilde{F})}$, where by
definition $$\chi_x(\exp{X}) = \exp{({\sqrt{ -1}\over\hbar}\langle x,X\rangle )},$$ and
by the definition $$\tilde{\sigma}|_{((G_x)_0, A_x)} = mult \chi_x $$ and
\item[(2)] The point $\sigma_0$ is fixed under $(G_x,A_x)$-action on the
dual $\widehat{(H_0,\tilde{F})}$,
\end{enumerate} 

(f) $\rho$ is some representation of the complex Lie super-algebra ${\frakt p}$
in Hilbert super-space $V$ such that its restriction to the real part ${\frakt 
h = p \cap g}$ is equivalent to the derivative $D\sigma_0$ of the 
representation $\sigma_0$.
\end{defn}

\begin{thm}
$(L,\rho,\sigma_0)$ is a $(\chi_x,\tilde{\sigma})$-polarization iff
$({\frakt p},H,\rho,\sigma_0)$ is a $(\tilde{\sigma},x)$-polarization.
\end{thm}
\begin{pf} 
If $(L,\rho,\sigma_0)$ is a
$(\chi_x,\tilde{\sigma})$-polarization, $L$ is weakly Lagrangian,
invariant, integrable tangent super-distribution of the tangent
super-bundle $T(\Omega,A_x\setminus A)$. In the notations before
we have remarked that we could reconstruct the Lie super-subalgebra {\frakt p}
 form the distribution $L$, ${\frakt p} = ({\frakt g}_x)_{\Bbb C} \ltimes
L_x$. It is easy to verify that we have in this case a $(\tilde{\sigma},
x)$-polarization. Conversely, it is easy to reconstruct a
$(\chi_x,\tilde{ \sigma})$-polarization $(L,\rho,\sigma_0)$ from some
$(\tilde{\sigma}, x)$-polarization $({\frakt p},H,\rho,\sigma_0)$.
\end{pf}

With the same reason as in the ordinary Lie group case, we can easily
prove

\begin{cor}
Suppose that $(\Omega,C;\omega_{\Omega})$ is an admissible ( oder
integral ) orbit of a mechanical system with super-symmetry
$(M,B;G,A,\omega)$, $(L,\rho,\sigma_0)$ a
$(\chi_x,\tilde{\sigma})$-polarization, where as above
$\tilde{\sigma}$ is a representation of $(G_x,A_x)$, the restriction
of which to the connected component $((G_x)_0,A_x)$ is a multiple of
the character $\chi_x$. Then

(1) The homogeneous super-space $\Omega$ admits the structure of a mixed
manifold of type $(k,l,m)$ in the sense of \cite{diep4},\cite{diep13}.

(2) There exists a unique irreducible unitary representation $\sigma$
of the super-group $(H,F)$ such that its restriction to the stable
super-subgroup $(G_x,A_x)$ is a multiple of the representation
$\tilde{\sigma}$ and its derivative is the restriction of the representation
$\rho$ to the real part ${\frakt h = p \cap g}$.
\end{cor}

{\bf Construction } ( see also (\cite{kostant2}, \S6.1) for the line
super-bundle case ).  Assume that $(L,\rho,\sigma_0)$ is a
$(\chi_x,\tilde{\sigma})$-polarization of our orbit $\Omega$. From the
preceding theorem, we can obtain some representation $\sigma$ of the
polarization (closed) super-subgroup $(H,F)$ in some Hilbert super-space $V$
. Let us denote by $\tau_G$ the natural projection $(G,A)
\twoheadrightarrow (H\setminus G,F\setminus A)$. Let $U\subseteq
H\setminus G$ be an open set and let $V = \tau^{-1}_G(U) \subseteq G$.
one has therefore $A(V)^*B(H)^* \subseteq A(V)^*$.

Denote by $A(V,\sigma)$ the super-space of all $V$-valued super-space
$f\in E_V(U)$ such that $$\langle wv,f\rangle  = \sigma(w)\langle v,f\rangle ,\forall v\in A(V)^*
, w\in B(H)^*. $$ It is easy to see that $(F\setminus A)(U)$ can be
embedded into $A(V) $ as a sub-super-space of super-functions such that
$$\langle wv,g\rangle  = \langle w,1_H\rangle \langle v,g\rangle , \forall v\in A(V)^*, w\in B(H)^*.$$ So, the
correspondence $$U \subseteq H\setminus G \mapsto E_V(U) =
E_V(\tau^{-1}_G( U),\sigma)$$ is a sheaf on $H\setminus G$. Let
$\hat{\sigma}$ be the element of $E_V(\tau^{-1}_G(U))$ such that for
any $w\in B(H)^*$, $$\langle w,\hat{\sigma}\rangle  = \sigma(w)t_0,$$ and for a
fixed $t_0$ in $E_V(\tau^{-1}_G(U)) $. Because the representation is
irreducible, we see that the open $U$ is principal for the sheaf $E_V$.
Thus we have some sheaf.

It is easy to see also that $E_V(G,A,\sigma)$ is a closed
sub-super-space of $V$-valued super-functions on $(G,A)$ and it is stable
under the action of our super-group $(G,A)$ on the right. In
particular, the elements of the Lie super-algebra {\frakt g} acts via
the Hamiltonian super-fields on the right. Hence for each open set $U$
in $H\setminus G$ the super-subspace $E_V(\tau^{-1}_G(U),\sigma,\rho)$
consisting of the sections the covariant derivative of which vanishes
along the vector super-fields from $\overline{\frakt p}$ for the polarization
${\frakt p}$, form a sub-sheaf, which gives us also an invariant closed
super-subspace of global sections of our quantum super-bundle ${\cal E}
_{V,\rho,\sigma}$. We refer to this invariant super-subspace of global sections
as the {\it induced representation $\Ind(G,A;{\frakt p},H,\rho,\sigma_0)$. }

\begin{cor} The natural representation called {\rm partially invariant
holomorphly induced representation } and denoted by
$\Ind(G,A;L,x,\rho,\sigma_0)$ of the Lie super-group $(G,A)$ in the
super-space of so called {\rm partially invariant and partially holomorphic
sections} of the induced super-bundle ${\cal E} _{V,\rho,\sigma_0}$ is
equivalent to the natural right regular representation by right
translations on the super-space $A(G,A;L,x,\rho,\sigma_0)$ of smooth
super-functions on $G$ with values in $V$ such that $$\langle wv,f\rangle  =
\sigma(w)\langle v,f\rangle , \forall w\in B(H)^*,\forall v \in A(V)^*,$$
$$\nabla_{\xi_X}f \equiv 0, \forall X \in \overline{\frakt p},$$ where
by definition, $$\nabla_{\xi_X} = \theta(\xi_X) + \langle X,\rho\rangle .$$ \end{cor}

Finally by the same computation as in the ordinary Lie group case, we can
compute the Lie derivative of the introduced partially invariant and
holomorphly induced representation.

\begin{thm}
The Lie derivative of the partially invariant and holomorphly induced
representation $\Ind(G,A;{\frakt p},(H,F),\rho,\sigma_0) $ is equivalent to the
Lie super-algebra representation $$X \in {\frakt g} \mapsto {\sqrt{-1}\over
\hbar}\hat{f}_X $$ of the Lie super-algebra ${\frakt g} = \Lie (G,A)$  via the
multidimensional quantization procedure.
\end{thm}

\section{Bibliographical Remarks}

The main ideas of these reductions, modifications and super-version were
due to the author of this book. The author posed problems for former
Ph.D. students Tran Vui and Tran Dao Dong. Tran Vui has completed and presented
there Dissertations in 1994 and Tran Dao Dong - in 1995.
In a revised version, this chapter reproduces the works \cite{dongvui},
\cite{dong1}, \cite{dong2}, \cite{vui1}-\cite{vui3} and \cite{diep22}.

\chapter{Index of Type I C*-algebras}

As said in the introduction, the problem of description of irreducible
unitary representations of groups is a preliminary step for using the
general Gel'fand transform to study the structure of
$L^1(G)$ and $C^*(G)$. In this chapter we introduce a general
construction of indices of type I C*-algebras. we shall prove that if  $\pi$
is an irreducible *-representation of  a  C*-algebra $A$ of type I then an
element $a\in  A$ is transformed under  $\pi$ into a compact operator
if and only if  it belongs into the intersection of kernels of the
irreducible unitary representations  $\pi'\ne \pi$ of  $A$, which belong
to the closure of the single point set $\{ \pi \}$. This result give us
a possibility to construct canonical composition series with CCR
quotients, by which we shall find out a complete system of invariants of
C*-algebras of type I. At the last section we shall show some concrete
applications to group C*-algebras.

\section{Compact Type Ideals in Type I C*-Algebras}

Let $A$ to be a C*-algebra of type I with identity element, we shall add
it formally if the algebra has no one, $\hat{A}$  its dual object with
the usual Jacobson Dixmier-Glimme-Sakai topology, see {\sc A.
Kirillov}\cite{kirillov}, {\sc J. Dixmier}\cite{dixmier}.
For simplicity we shall identify the equivalence classes of irreducible
*-representations  with their representatives.

A C*-algebra $A$ is of type I if and only if the image of $A$ under any
non-degenerate irreducible *-representations contains the ideal ${\cal
K}({\cal H})$ of compact operators in the separable Hilbert space ${\cal
H}= {\cal H}_\pi$, in which the representation is realized.

\begin{defn} The preimage $${\cal K}_\pi = \pi^{-1}({\cal
K}{\cal H})$$ of the ideal ${\cal K}({\cal H})$ in  $A$ is called an
{\it ideal of compact type}, associated with the representation $\pi$.
\end{defn}

From this definition, we can deduce directly that the compact type ideals
are always closed and two-sided in $A$, see {\sc J.
Dixmier}\cite{dixmier}.

If the representation $\pi$ is not of CCR, following the definition, the
ideal  ${\cal K}_\pi$ can not coincided with the whole C*-algebra $A$.
In other side, the C*-algebra $A$ is of type I iff ${\cal K}_\pi \ne 0$.
Hence, we have an proper ideal.
One pose a natural question is how to describe these ideals in the
general case.

We recall a well-known result.
\begin{thm}[\bf Rosenberg \cite{rosenberg}- J.  M. G. Fell \cite{fell}]
Irreducible *-representation $\pi$ of C*-algebra of type I $A$ is  CCR if
and only if the one-point set  $\{\pi \}$ is closed in the topology of
the dual object $\hat{A}$ cur a $A$.
\end{thm}

For a subset $X$ in the dual object $\hat{A}$ of a C*-algebra $A$, we
denote by $\bar{X}$ the closure of $X$.

\begin{thm}
Let $a$ to be an element of $A$,   $\pi$ an irreducible non-degenerate
*-representation of $A$. Then the operator $\pi(a)$ id a compact
operator if and only if $\pi'(a)   =   0$  for all $\pi'\in  \{\pi\}$
and   $\pi'\ne \pi$.\end{thm}
\begin{pf}
As remarked above, ${\cal K}_\pi$ is a two-sided ideal of a type I, or
in other words, GCR C*-algebra $A$.  The restriction of $\pi$ on ${\cal
K}_\pi$ is a nontrivial *-representation, because $\pi({\cal  K}) =
{\cal K}({\cal H})$. Following {\sc J. Dixmier}(\cite{dixmier},
\S3.2.1),   $\pi\vert_{\cal K} \in \hat{\cal
K}$ is an irreducible  CCR *-representation of ${\cal K}$.  Following
the theorem of {\sc A. Rosenberg - J. M. G. Fell},   the one-point set
$\{ \pi \}$ is   closed  in the topology of $\hat{\cal K}$, which is a
closed set in $\hat{A}$.

Assume that $\pi'\ne  \pi$ belongs to the closure $\overline{\{\pi\}}$
of the one-point set $\{\pi\}$ in the topology of $\hat{A}$. Following
{\sc J. Dixmier}(\cite{dixmier},Prop.3.2.1), if      $\pi'\vert_{\cal  K}  \ne  0$,   then    $\pi'\in \hat{\cal K}$ and
$\pi'$ is a weak limit of $\pi$. This is impossible, as said above 
. Thus, we have $\pi'\vert_{\cal K}  = 0$. In other words, if $\pi(a)$
is a compact operator, i.e.  $\pi(a)\in  {\cal  K}$,  then
$\pi'(a)=0$ for all        $\pi'\in\overline{\{\pi\}} \setminus \{\pi\}$.

Conversely, suppose that  $\pi'(a)   =   0$  for all $\pi'\in
\overline{\{\pi\}}\setminus  \{\pi\}$. We must prove that $\pi(a)$ is a
compact operator.

Suppose that
$$B = \bigcap_{\pi' \in \overline{\{\pi\}}\atop \pi'\ne \pi} \Ker \pi'.$$
Then $B$ is an two-sided closed ideal in $A$.  Thus, the dual object
$\hat{B}$ of $B$ is included in the dual object $\hat{A}$ of $A$ as an
open subset. More exactly, $$\hat{B} = \hat{A}   \setminus
\overline{(\overline{\{\pi\}}\setminus \{\pi\})}.$$

Really, $\tilde{\pi}\in  \overline{\overline{\{\pi\}}\setminus  \{\pi\}}$
means that
$$\Ker \tilde{\pi}  \supseteq  \bigcap_{\pi'\in  \overline{\{\pi\}}\setminus
\{\pi\}} \Ker\pi'.$$
Hence,       $\tilde{\pi}\vert_B   \equiv   0,   \forall   \tilde{\pi}\in
\overline{\overline{\{\pi\}}\setminus \{\pi\}}$. In other words, we have
$$B = \bigcap_{\pi'\in\overline{\overline{\{\pi\}}\setminus
\{\pi\}}}\Ker\pi'.$$
Following {\sc J. Dixmier}(\cite{dixmier},\S3.2.2), we have
$$\hat{B} = \hat{A}   \setminus   \overline{\overline{\{\pi\}}\setminus
\{\pi\}}.$$
If      $\pi\vert_B = 0$, then nothing is need to check. If
$\pi\vert_B \ne 0$, we have           $\pi\in \hat{B}$, as above.
Because $\hat{B}$ is dense in $\hat{A}$, and following the definition of
$B$, we conclude that $\{\pi\}$ is closed in $\hat{B}$. This means that
$\pi$ is a CCR *-representation of $A$, in other words, if $pi'(a) = 0$
with every $\pi' \in \overline{\{\pi\}} \setminus \{\pi\}$, we have $\pi(a) \in
{\cal K}({\cal H})$.
\end{pf}

\begin{cor}
The compact type ideal ${\cal K}_\pi$ can be described as
$${\cal K}_\pi  = \{  a\in A\quad  ; \quad  \pi'(a) = 0,   \forall \pi' \in
\overline{\{\pi\}}\setminus \{\pi\} \}.$$
\end{cor}

\begin{rem}
We have therefore a solution of our problem of description of compact
type ideal. We also remark that the set $\{\pi\}$ is closed if and only
if $\overline{\{\pi\}} \setminus \{\pi\}  $ is an empty set. Therefore,
our theorem contain the Theorem of A. Rosenberg  -  J. M. G. Fell as a
particular case.
\end{rem}

\section{Canonical Composition Series}

We construct now  canonical composition series for each "good"
C*-algebra, which satisfy the following definition.

\begin{defn} We say that element $\pi\in \hat{A}$ is a {\it boundary point
of degree 0 },   if it  belongs to the closure of no other point of
$\hat{A}$. We denote by $X_0$ the set of all boundary points of degree
0. The point    $\pi\in \hat{A}\setminus X_0$  is said to have  {\it
boundary degree 1}, if it belongs to the closure of no other point of
$\hat{A}\setminus X_0$ with the induced topology. We denote by $X_1$ the
set of all points of boundary degree 1. The point $\pi\in \hat{A}$ is
said to have boundary degree 2, if it belongs to the closure of no other
point of $\hat{A}\setminus X_0 \cup X_1$ with induced topology, etc....

We continue this procedure until infinity, if $\hat{A} \setminus (X_0 \cup
X_1 \cup... )$ is not empty.

Suppose  that $\rho$ is the first transfinite cardinal number and we
suppose that for all preceded cardinal,
$\tilde{\rho} \leq  \rho$,   the sets $X_{\tilde{\rho}}$ were already
constructed. We pose $X_\rho   =   \varnothing$   and consider the set
$\hat{A}   \setminus
\bigcup_{\tilde{\rho}\leq \rho}  X_{\tilde{\rho}}$  with the induced
topology. If this set is nonempty, we continue this procedure,.... We
obtain thus a sequence of sets $X_0,  X_1,  X_2,\dots,  X_\alpha$, where
$\alpha$ is a fixed cardinal number
$$\hat{A} = \bigcup_{\rho \leq \alpha} X_\rho.$$

We say that a C*-algebra has {\it boundary property} iff
$X_{\tilde{\rho}}$ is a nonempty open set in $\hat{A} \setminus
\bigcup_{\rho' < \tilde{\rho}} X_{\rho'}$ in every case when
$\tilde{\rho}$ is a non-limit cardinal, $0 \leq \tilde{\rho} < \alpha$ and
$$\hat{A} \setminus \bigcup_{\rho' < \tilde{\rho}} X_{\rho'}\ne
\varnothing. $$ \end{defn}

We see that the set of GCR C*-algebras with boundary property is rather
wide. The most of C*-algebras $C^*(G)$  of locally compact groups in our
works and in other works of {\sc J. Rosenberg} have this property. 
We pose $${\cal E}_\rho:= \bigcap_{\pi\in A_\rho} \ker\pi,$$ where
$$A_\rho:= \hat{A} \setminus \bigcap_{\rho' < \rho} X_{\rho'},  \text{
with} \rho > 0, \text{ and } A_0:= \hat{A}.$$

\begin{thm}
Suppose that  $A$ is a GCR C*-algebra with boundary property. Then:
\begin{enumerate}
\item[a)] All ${\cal E}_0,   {\cal E}_1,\dots,{\cal E}_\alpha$ are
two-sided closed ideals.
\item[b)] For each cardinal $\rho$, $0 \leq \rho < \alpha$,
$$\widehat{\cal E}_\rho = \bigcap_{\tilde{\rho} < \rho} X_{\tilde{\rho}}.$$
\item[c)] We have a tower of ideals
$$\{ 0\} \subset {\cal E}_0 \subset {\cal E}_1  \subset \dots \subset {\cal
E}_\alpha = A.$$
\item[d)] For each cardinal $\rho$,  $0 \leq \rho < \alpha$,  the
quotient C*-algebra ${\cal E}_{\rho+1}/{\cal E}_\rho$ is    CCR.
\end{enumerate}
Hence for each C*-algebra $A$ with boundary property, we have a
canonical composition series.
\end{thm}

Before to prove the theorem, we prove the following preliminary lemma.

\begin{lem}
Let $T$ be a topological space, $U$ an open subset of $T$,  $F=
T\setminus U$, $O$ is an open set in the induced topology of $F$. Then,
$U \coprod O$ is an open set in $T$.
\end{lem}
\begin{pf}
Really, let $x$ to be an arbitrary point in $U  \coprod O$. Then either
$x\in U$, or       $x\in O$. In the first case, there exists a
neighborhood
$U_1$ of     $x$ in   the topology of       $T$,  such that   $U_1
\subseteq U \subseteq U\coprod O$. In the second case, there exists a
neighborhood $U_1$ of     $x$
in the induced topology of $F$, such that $O_1  \subseteq O \subseteq U\coprod
O$. Following definition of induced topology, there exists an open $V$ in
$T$, such that  $V \cap F = O_1$. Then, $V =  V \cap T = (V  \cap  U) \coprod
(V \cap F) = (V \cap U) \coprod  O_1 \subseteq U \coprod O$.  Thus, in
any case, we can always choose an open neighborhood of
$x$ in $U \coprod O$. This means that $U \coprod O$ is open in $T$. 
 \end{pf}

\begin{cor}
Let $A$ to be a GCR C*-algebra with boundary property. Then each set
$A_\rho = \hat{A} \setminus  \bigcup_{\rho'<\rho}  X_\rho'$  is closed.
\end{cor}
\begin{pf}
Really, because $A$ is a GCR C*-algebra with boundary property, the sets
$X_\rho$
, $0 \leq \rho < \alpha$  are nonempty open sets, if $\rho$ is not a
limit cardinal number. Thus, $A_1 =  \hat{A} \setminus X_0$  is closed.
Because $X_1$  is open in $A_1$, following the above lemma, $X_0
\cup X_1$ is open in $\hat{A}$, and  $A_2  =  \hat{A}  \setminus
X_0 \cup X_1$ is closed, etc....

For each finite number $n\in  {\Bbb  N}$, the set $X_0  \cup  \dots \cup
X_n$ is open in $\hat{A}$. For the first limit cardinals $\omega$,
$$\bigcup_{\rho   <   \omega}   X_\rho   =  X_0  \cup  \dots  \cup  X_n\cup
\dots $$
is open in $A$, then 
$$A_\omega  = \hat{A}  \setminus \bigcup_{\rho <
\omega} X_\rho$$ 
is closed. Continue the procedure infinitely time,we see that all the
sets $A_\rho$ are closed in $\hat{A}$.
\end{pf}

Finally, we return to Proof of the theorem.

(a) Because $A$ is a GCR C*-algebra with boundary property then each set
$X_\rho$ is open in the closed set $$A_\rho =
\bigcup_{\tilde{\rho}<\rho} X_{\tilde{\rho}}.$$
Following the corollary, all the set $A_\rho$ are closed. Thus,
following definition, all the set $${\cal E}_\rho =  \bigcap_{\pi\in
A_\rho} \Ker \pi$$
are closed two-sided ideals and ${\cal E}_0 = \{  0 \}$,   ${\cal E}_\alpha = A$.

(b) Because $\cup_{\tilde{\rho}}  X_{\tilde{\rho}}$ is open in
$\hat{A}$, we have naturally
$$\hat{\cal E}_\rho = \bigcup_{\tilde{\rho} < \rho} X_{\tilde{\rho}}$$

(c) Because $A_0 = A$,  $$A_\rho  =  \hat{A}  \setminus
\bigcup_{\tilde{\rho} <
\rho} X_{\tilde{\rho}},$$ we have a tower of inclusions
$$A_0 = \hat{A} \supseteq A_1 \supseteq \dots \supseteq A_\alpha = \varnothing.
$$

(d) Following definition, we have $$\widehat{({\cal E}_{\rho + 1}/{\cal
E}_\rho)}  = \widehat{\cal E}_{\rho +
1} \setminus \widehat{\cal E}_\rho = X_\rho.$$

Finally, we need only to show that each irreducible *-representation $\pi\in
\widehat{\cal E}_{\rho + 1} / {\cal E}_\rho$  is    CCR. Indeed, $$\pi\in
\widehat{({\cal E}_{\rho + 1}  / {\cal E}_\rho)} = X_\rho  \subseteq A_\rho
= X_\rho \coprod A_{\rho  +1}$$  and moreover $$\overline{\{\pi  \}}
\setminus \{\pi\} \subseteq A_{\rho + 1}.$$ Because
$${\cal E}_{\rho + 1} = \bigcap_{\pi'\in A_{\rho + 1}} \Ker \pi',$$
$\{\pi\}$ is closed in the induced topology of $X_\rho$. Following the
theorem of compact type ideals, the representation  $\pi$ is CCR. 

\section{Index of type I C*-Algebras}

We construct now the canonical composition representation for GCR
C*-algebras with boundary property.

Let $\alpha$ to be a fixed cardinal number, $C_1,  C_2,  \dots C_\alpha$
a system of C*-algebras such that $C_\rho = \{ 0 \}$, if $\rho$ is a
limit cardinal number. We denote
$\widetilde{\Ext_0}(C_1,C_2,\dots,C_\alpha)$ the
set of all increasing towers of C*-algebras of type I $$\{ 0\} \subset
E_) \subset E_1 \subset \dots  \subset E_\alpha = A,$$ such that
\begin{enumerate}
\item[a)] all the C*-algebras $E_\rho$ are closed two-sided ideals in $A$,
\item[b)] if $\rho$ is a limit cardinal, $E_\rho$ coincided with the
closure of the ideal $\bigcup_{\tilde{\rho} < \rho} E_{\tilde{\rho}}$ in    $A$,
\item[c)] for each non-limit cardinal number $E_{\rho+1}/E_\rho
\cong C_{\rho + 1}, 0 \leq \rho < \alpha.$
\end{enumerate}

In the set    $\widetilde{\Ext_0}(C_1,C_2,\dots,C_\alpha)$ we define an
equivalence relation.

\begin{defn}
Two towers of ideals $\alpha$
$$\{ 0\} \subset  E_0  \subset  E_1  \subset  \dots  \subset  E_\alpha  = A
\leqno{({\frakt A})}$$ and
$$\{ 0 \} F_0 \subset   F_1   \subset   \dots   \subset   F_\alpha   =   B
\leqno{({\frakt B})} $$
are called {\it  equivalent} if there exists an isomorphism $\theta:
A \to B$  inducing *-isomorphisms between the corresponding ideals
$$\theta_i = \theta|_{E_i} \;:\; E_i \to F_i.$$
\end{defn}

We denote by $\Ext_0(C_1,C_2,\dots,C_\alpha)$ the set of all equivalent
classes of towers of ideals in C*-algebras of type I.
The main problem is how to find out the invariants characterizing the
equivalence classes from $\Ext_0(C_1,C_2,\dots,C_\alpha)$.

First we recall a little bits about  {\sc  R. C. Busby} invariant
\cite{busby}.
Let $D$ to be an arbitrary C*-algebra, the pair $(T',T'')$ of two maps
$$T',  T'': D
\longrightarrow D$$ satisfying the properties $$xT'(y) \equiv T''(x) y,
\forall x, y \in D$$
is called a {\it double multiplier}. \index{multiplier!double -} We 
denote by $M(D)$ the set 
of all double   multipliers of $D$. Equip to $M(D)$ natural operations and
the norm of operators, $M(D)$ become a C*-algebra, called the {\it
multiplier C*-algebra}. Remark that $M(D)$ contains  $D$ as an ideal and
$M(D)$ is the smallest C*-algebra, which contains $D$ as a proper closed
ideal. We denote $O(D) =  M(D)/D$ the  quotient C*-algebra and denote
$\Hom_0(C,O(D))$ the set of all  *-homomorphism from a C*-algebra  $C$
to $O(D)$.

In the set    $\Hom_0(C,O(D))$, we introduce an equivalence relation.
Two elements $\gamma_1,   \gamma_2$ of $\Hom(C,O(D))$ are called
{\it equivalent } \index{extensions!equivalent} if there exists 
*-isomorphisms $\omega \;:\;  D \to D$ and   $\sigma
\;: \; C \to C$ such that
$$\gamma_2 = \hat{\omega}\gamma_1 \sigma^{-1},$$
where following definition, $\hat{\omega} \;:  \; O(D) \to O(D)$ is
*-isomorphism, associated with $\omega$. We denote by $\Hom(C,O(D))$ the
set of all these equivalence classes.

Assume $\gamma \in \Hom_0(C,O(D))$ is a fixed *-homomorphism. We
consider the fibered product
$$E_\gamma:= \{ (m,c) \in M(D) \times C \;: \; \pi(m) = \gamma(c)\},$$
where $\pi$ is the canonical projection from  $M(D)$
onto the quotient C*-algebra  $O(D) = M(D)/D$. Following {\sc R. C.
Busby}\cite{busby},   $E_\gamma$
is a C*-algebra. Assume $\mu \;:   \;   D \hookrightarrow M(D)$ is
the natural inclusion. We denote $E^0$ the "universal" extension
$$\CD 0 @>>> D @>\mu >> M(D) @>\pi >> O(D) @>>> 0\endCD \leqno{E^0}$$

We can complete the commutative diagram
$$\CD
@.     @.        @.        C         @.\\
@.     @.        @.        @VV\gamma V @. \\
0 @>>> D @>\mu >> M(D) @>\pi>> O(D) @>>> 0
\endCD$$
by the homomorphisms $f:  D \to E_\gamma$ following the formula $f(d):= (\mu(d),0),
\forall d \in D$ and   $g:   E_\gamma \to C$ following the formula  $g(m,c):  = c,
\forall d\in D, m\in M(D)$
in order to have the following commutative diagram
$$\CD
0 @>>> D @>f>> E_\gamma @>g>> C @>>> 0\\
@.    @VVId V   @VVpr_1 V       @VV\gamma V @. \\
0 @>>> D @>\mu>> M(D) @>\pi>> O(D) @>>> 0.
\endCD$$

The obtained extension of C*-algebras
$$\CD 0 @>>> D @>>> E_\gamma @>>> C @>>> 0 \endCD$$
is denoted by  $E^0_\gamma$ and its equivalence class by $[E^0_\gamma]$.
{\sc R. C. Busby}\cite{busby} has had proved the following result.

\begin{thm}[\bf R. C. Busby]
The map
$$[\gamma]\in \Hom(C,O(D)) \mapsto [E^0_\gamma]  \in \Ext(D,C)$$ is a
bijection.
\end{thm}
We remark that {\sc R.  C. Busby} denoted our set $\Ext(D,C)$ as  $\Ext(C,D)$.

We consider now the set  $\Ext(C_1,C_2,\dots,C_\alpha)$ of equivalence
classes of towers of ideals of C*-algebras of type
$$\{ 0\} \subset E_0  \subset  E_1  \subset  \dots  \subset  E_\alpha  =  A
\leqno{({\frakt A})}$$
We have naturally the extensions
$$\CD 0 @>>> E_1 \cong C_1 @>>> E_2 @>>> E_2/E_1 \cong C_2 @>>> 0. \endCD$$
There exists an element  $[\gamma_1] \in \Hom(C_2, O(C_1))$, defining
the equivalence class of this extension, i.e. the equivalence class of
the tower  
$$\CD \{ 0\}= E_0 \subset E_1 \cong C_1 \subset E_2.\endCD$$
Thus up to equivalence, we can replace $E_2$ by  $E^0_{\gamma_1}$.

Then following the extension
$$\CD 0 @>>> E^0_{\gamma_1} \cong E_2 @>>> E_3 @>>> E_3/E_2 \cong C_3 @>>> 0.
\endCD$$
As above, there exists an element $[\gamma_2]   \in
\Hom(C_3,O(E^0_{\gamma_1}))$,  defining the equivalence class of the
extension and is defined uniquely the equivalence class of the tower 
$$\{ 0\} = E_0 \subset E_1 \subset E_2 \subset E_3,$$
etc....

If  $\rho$ is a limit cardinal, then following the conditions of the
definition of boundary property.
$$E_\rho = \overline{\bigcup_{\tilde{\rho} < \rho} E_{\tilde{\rho}}}$$
and
$E_\rho$ is uniquely defined by the system of invariants $\{
[\gamma_{\tilde{\rho}}]\}_{\tilde{\rho < \rho}}$.

Then the equivalence class of the tower of ideals
$$\{ 0\} = E_0  \subset  E_1  \subset  E_2  \subset  \dots  E_\rho  \subset
E_{\rho +1}$$
is uniquely defined by the invariant
$$[\gamma_\rho]  \in  \Hom(C_{\rho+1},   O(\overline{\bigcup_{\tilde{\rho} <
\rho}E^0_{\tilde{\rho}}})). $$

We have proved that:
With each increasing tower of ideals of C*-algebras of type
$$\{ 0\} = E_0 \subset E_1 \subset  E_2 \subset \dots \subset E_\alpha  = A
\leqno{({\frakt A})}$$
in the set        $\Ext_0(C_1, C_2, \dots, C_\alpha)$,  there exists an
unique system of invariants, uniquely defining the equivalence class of
the tower ${\frakt A}$:
$$ \begin{array}{rl} [\gamma_1] &\in \Hom(C_2,O(C_1)),\cr
             [\gamma_2] &\in \Hom(C_3,O(E^0_{\gamma_1})),\cr
             \ldots\ldots  & \ldots\ldots\ldots\ldots\ldots\ldots\ldots\ldots\ldots \cr
             [\gamma_\rho]   &\in   \Hom(C_{\rho+1},
O(\overline{\bigcup_{\tilde{\rho}<   \rho} E^0_{\gamma_{\tilde{\rho}}}}  )),
\cr 
                \ldots\ldots\ &\ldots\ldots\ldots\ldots\ldots\ldots\ldots\ldots\ldots \end{array}$$

\begin{defn}
Assume that the C*-algebra $A$ can be decomposed into an increasing
tower of ideals with length  $\alpha$,
$$\{ 0\} = E_0 \subset E_1 \subset E_2 \subset \dots \subset E_\alpha  = A.
\leqno{({\frakt A})}$$
Then the system of invariants is called the {\it index } \index{index of 
C*-algebra} of C*-algebra $A$ and is denoted by $\Ind_\alpha A$.
\end{defn}

We have therefore proved the following results 

\begin{thm}
The system of invariants $\Ind_\alpha A = \{ [\gamma_\rho]\}_{0 \leq
\rho < \alpha}$ uniquely defines,up to equivalence the structure of GCR
C*-algebra  $A$, if it admits a   composition series of ideals of type
${\frakt A}$.
\end{thm}

\section{Compactness Criteria for Group C*-Algebras}

Applying the K-homology to studying the structure of C*-algebras, we
introduced the method of topological invariants. In the previous chapter,
this method was constructed in most general for the class of C*-algebras
of type I with boundary property, including the C*-algebras of groups. In
this theory the ideals of compact type plays a central role. We described
theses ideals in the language of the topology of dual objects. This
topology is just the Zariski topology for primitive ideals and in general
is complicate to describe. In this last section of this chapter, we shall
describe these ideals in a more geometrical language. For the
representations induced from irreducible unitary representations of
closed normal subgroups, we describe the compactness by sing the partial
Fourier-Gel'fand transforms.

\subsection{Compactness Criteria}

Let $G$ be a countable at infinity locally compact group, $G_1$ a closed
normal subgroup of  $G$, $X= G_1 \setminus G$ the quotient group of $G$
by   $G_1$. Choose representatives of the coset classes from
$X=G_1\backslash G$ in such a way that we have a Borel section  $s:  X \to G$
corresponding to the right invariant Haar
measure on $G$. We have a decomposition $g =
bs(a),$ where          $b = b(g)$ and $a = a(g)$, for all        $g$
in        $G$. As a set,
 $G$ can be identify with $G_1 \times
X$ and the group multiplication la is defined by the formulae
$$(b,a).(b',a') =
(b[\alpha(a)b']\beta(a,a'),aa'),$$ where the map $\alpha: X
\to \Aut G_1$ is defined by the formula
$$\alpha(a)(b,e) = (e,a)(b,e)(e,a)^{-1},$$ and the map        $\beta: X
\times X \to G_1$ is defined by the formula
c $$(e,a)(e,a') = (\beta(a,a'), aa').$$ We can  therefore
choose the right (quasi-)invariant Haar measure  $dg$, $dg_1$, $dx$ on      $G$, $G_1$, $X$, resp. such
that $$d_rg:= dg = dg_1dx, \quad d_lg =
(\Delta_G/\Delta_{G_1})(g_1)d_lg_1d_lx,$$
 if  $g = g_1.s(x).$
Then
$$d(xg)/dx = (\Delta_{G_1} / \Delta_G)(h),$$
if               $$gs(x)g = hs(x).$$

Let                $T = \Ind_{G_1}^G\xi$  be an irreducible unitary
representation of $G$, induced in th sense of  {\sc G. Mackey} from a
CCR presentation $\xi$ of a closed normal subgroup
$G_1$. Then it is realized in the Hilbert space  $L^2(X,dx)$ of
square-integrable functions with values in the space of the
representation $\xi$, following the formula
$$T(b,a)f(x) = (\Delta_{G_1}/\Delta_G)^{1\over
2}(h)\xi(h)f(xa),$$ where          $h = [\alpha(x)b]\beta(a,x)$ and
$\Delta_{G_1}$ and $\Delta_G$ are the modular functions of Haar
measures on $G_1$ and   $G$, respectively.

We denote by $\overline{S}$ the closure of the set  $S$ in the topology
of the dual $\hat{G}$ of group  $G$. We shall use the same letters to
denote  representations and their corresponding equivalence classes of
C*-algebras. For example, we use notation $$T(\varphi) = \int_X\int_{G_1}
T(b,a)\varphi(b,a)dbda,$$ for every       $\varphi \in L^1(G)$.

We define the partial Fourier - Gel'fand transform in variable $b$ of
functions            $\varphi(b,a)$
from    $L^1(G,dbda)$, by the formula            $$\tilde{\varphi}(x,a) =
\int_{G_1}(\xi.(\Delta_{G_1}/\Delta_G)^{1\over 2})(\alpha(x)b)
\varphi(b,a) db.$$ When $x$ runs over the set       $X$, the
representations
$\alpha^*(x)\xi$ run over an orbit 
$O = \alpha^*(X)\xi$ in the dual object  $\hat{G}_1$ of
$G_1$, where          $$(\alpha^*(x)\xi)(b):= \xi(\alpha(x)b).$$

It is easy to see that the functions $\Vert
\tilde{\varphi}(.,a)\Vert$ and   $$\Vert \int
\tilde{\varphi}(.,a)(\xi(\Delta_{G_1}/\Delta_G)^{1\over
2}(\beta(.,a))R(a)da\Vert$$, where          $R(.)$ is the regular representation
of the group $X$,
are lower semi-continuous when the quasi-orbit  $O$ is separable at
infinity and they vanish at this infinity of the quasi-orbit $O$. 

When $\alpha^*(x)\xi$ tend to the boundary of the quasi-orbit $O$ or to
the infinity of the quasi-orbit $O$, we have the restriction on
$\partial O \cup \infty$, which we shall denote by
$\tilde{\varphi}(., a)\vert_{\partial O\cup \infty}$, etc....

Let $\zeta\in \partial O$ and   $\omega = G.\omega$ be the quasi-orbit
of $\zeta$, $\tau$ an irreducible representation of the stabilizer
$G_\zeta$, under the action of $G$, the restriction of which onto $G_1$
is a multiple of $\zeta$.  We denote
$T^{\omega,\tau} = \Ind_{G_\zeta}^G(\tau)$, which is irreducible
following the Mackey theory of small subgroups.

We suppose that the quasi-orbit $O$ is separable. We denote by $\tr$ the
trace function. In particular, we have $$tr\tilde{\varphi}(x,a) = tr
\int(\xi(\Delta_{G_1}/\Delta_G)^{1\over
2}(\alpha(x)b)\varphi(b,a) db.$$

\begin{thm} For every $\varphi\in L^1(G)$, the following conditions are
equivalent:
\begin{enumerate}      \item[(1)] $T(\varphi)$ is a compact operator.
            \item{(2)} $$S(\varphi) = 0, \forall S \in \overline{T}
            \mbox{
and} S \ne T. $$ \item{(3)} $$T^{\omega, \tau}(\varphi) = 0, \mbox{
for every quasi-orbit } \omega \subseteq \partial O, $$
and for every representation $\tau$ the restriction of which onto
$G_1$ is a multiple of $\zeta$,
$\tau\vert_{G_1} = \mult  \zeta.$ \item{(4)} $\Ind_{G_1}^{G_\zeta}(\varphi)
= 0, \forall \zeta\in \partial O$. \item{(5)} $\tilde{\varphi}(.,
a)\vert_{\partial O} = 0, \text{ almost everywhere w.r.t.}a. $ \item{(6)}
$\int_X \tilde{\varphi}(.,a)(\xi(\Delta_{G_1}/\Delta_G)^{1\over 2}
(\beta(.,a))R(a) da\vert_{\partial O} = 0.$ \item{(7)} $tr \int_X
\tilde{\varphi}(., a)(\xi(\Delta_{G_1}/\Delta_G)^{1\over
2}(\beta(.,a))R(a) da \vert_{\partial O} = 0, $ if the last trace exists
and $\varphi$ is a positive function in $L^1(G)$.  \end{enumerate}
\end{thm}

\begin{cor}  The compactness criteria are independent from the functions
$\alpha(.)$ and $\beta(.,. )$ \end{cor}

{\sc Proof of the Theorem}
$(1) \Longrightarrow (2)$
was proved in the previous chapter.

$(2) \Longrightarrow (3)$ can b proved directly by using the Mackey
theory of representations induced from small subgroups and the results of
{\sc J. M. G. Fell}[46] on continuity of the induction functor: we
denote th property "weak containment" by ``$\in$". If
$\zeta\in \partial O$, $\zeta$ is weak contained in
$\overline{\{\xi\}}$, and thus         $\Ind_{G_1}^G\zeta$ is weak
contained in  $\Ind_{G_1}^G\xi$. Because $\tau\vert_{G_1} = \mult  \zeta$,
we have $\tau$ weak contain in
$\Ind^{G_\zeta}_{G_1}\zeta$ and $T^{\omega,\tau} = \Ind^G_{G_\zeta}\tau
\in \Ind^G_{G_\zeta} \Ind^{G_\zeta}_{G_1} \zeta = \Ind^G_{G_1}\zeta.$
From this, we have
$$T^{\omega,\tau} \in \Ind^G_{G_1}\zeta \in
\overline{\Ind^G_{G_1}\xi} = \overline{\{ T\}}.$$

$(3) \Longrightarrow (4)$. This assertion is clear,because
$\Ind^{G_\zeta}_{G_1} \zeta\vert_{G_1} = \mult  \zeta$. Hence,
$$\Ind_{G_1}^{G_\zeta}\zeta = \int^\oplus \tau d\mu(\tau),$$ where the
direct integral is taken over the set of all irreducible unitary
representations of
$G_\zeta$ the restriction of which onto $G_1$ is a multiple of
$\zeta$, and $d\mu(\tau)$ the fixed measure on it,
and finally, $$\Ind^G_{G_1}\zeta = \int^\oplus T^{\omega,\tau}
d\mu(\tau).$$

$(4) \Longrightarrow (5)$ is proved from the well-known properties of
regular representations. In particular, there exists an approximative
unity $$\CD \Vert
\varphi \* f_n - \varphi \Vert @>> n\to \infty > 0,\endCD$$ where $\{
f_n\}$ is a sequence of Dirac $\delta$-like functions.

$(5) \Longrightarrow (6)$ is clear.

$(6) \Leftrightarrow (7)$ is proved easily that: For positive operators,
the trace is vanishing iff nd only the operator is zero.

$(6) \Longrightarrow (1)$ can be proved with the help of two lemmas what
follow.

\begin{lem}[\bf Particular case]
Assume $\varphi(b,a) = \chi(b)\psi(a),$ where $\chi\in
L^1(G_1)$, $\tilde{\xi}\vert_{\partial O} = 0$; $\tilde{\xi}(x):=
\int_{G_1} \xi(\Delta_{G_1}/\Delta_G)^{1\over 2}(\alpha(x)b)\chi(b) db$,
$\psi\in L^1(X)$ and   $\Vert \psi\Vert_{L^1} \ne 0.$ Then $T(\phi)$ is
a compact operator.
\end{lem}
\begin{pf}
Because $X$ is countable at infinity, we can choose a sequence
of included one-in-another compact $$K_n \Subset K_{n+1} \Subset \dots
\Subset O,$$ such that $\cup_n K_n = O$ and continuous functions
$theta_n$ such that
$0 \leq \theta_n(x) \leq 1$ and   $$\theta_n(x) = \cases 1
&\mbox{if }       x\in K_n\cr
                       0 &\mbox{if } x\not\in K_{n+1}\cr\endcases $$
Because $\tilde{\chi}(x)$ is continuous in    $x$, and in virtue of
hypothesis,
   $$\tilde{\chi}\vert_{\partial O} = 0,$$ the operator   $T(\varphi)$
can be approximate by a sequence of operators  $A_n$, $$A_nf(x):=
\theta_n(x)\tilde{\chi}(x)\int_X f(xa)
(\xi(\Delta_{G_1}/\Delta_G)^{1/2}(\beta(x,a))\psi(a)da,$$ which are
{\it compact}.

Indeed, we have
\begin{enumerate}
\item[(a)]
$$\begin{array}{rl} \Vert   T(\varphi)   -   A_n   \Vert   &=   \Vert
(1-\theta_n)\tilde{\chi}\int_X   R(a) (\xi (\Delta_{G_1} / \Delta_G)^{1/2}
(\beta(x,  a))\psi(a) da\cr
  &\leq \Vert (1-\theta_n)\tilde{\chi}\Vert.
\Vert\int_XR(a)(\Delta_{G_1}/\Delta_G)^{1/2}(\beta(.,a))\cr
  &\psi(a)da\Vert, \end{array}$$ where $(1-\theta_n)\tilde{\chi}$ is the
  operator of multiplication by a function.

Because     $\tilde{\chi}\vert_{\partial O} = 0$ and
$(1-\theta_n)\tilde{\chi}\vert_{\partial O} = 0$, we can choose $n$ big
enough that $$\Vert
(1-\theta_n)\tilde{\chi}\Vert.\Vert \int_X
R(a)(\xi(\Delta_{G_1}/\Delta_G)^{1/2}(\beta(.,a))\psi(a) da\Vert$$ is
small enough.

\item[(b)]
Operator $A_n$ has an operator-valued kernel
$K_n(x,a)$, $$A_nf(x) = \int_X K_n(x,t) f(t) dt,$$ with    $$K_n(x,a) =
\theta_n(x)\tilde{\chi}(x)(\xi(\Delta_{G_1}/\Delta_G)^{1/2}(\beta(x,
a))\psi(x^{-1}a)\Delta_X(x),$$ where $\Delta_X$ is the modular function
of the Haar measure on $X$. Because the integral
$$\int_X dx \int_X \Vert K_n(x,a)\Vert da $$ converges absolutely,
following Fubini Theorem, the operators
$A_n$ are all compact ones.
\end{enumerate}
The lemma is therefore proved.
\end{pf}

\begin{lem}[\bf General case]
Let $\varphi$ be an element of $L^1(G)$ such that
$\varphi(., a)\vert_{\partial O} = 0$, Then the operator 
$T(\varphi)$ is compact. 
 \end{lem}
 \begin{pf}
 Because     $\varphi \in L^1(G_1
\times X)$, we can construct an approximative sequence
$$\CD \Vert \varphi_n - \varphi \Vert_{L^1} @>> n \to \infty >
0,\endCD$$ where  $$\varphi_n(b,a) = \sum_{k=1}^{N_n}
\chi_k(b)\psi_k(a)$$ and $\chi_k \in L^1(G_1)$, $\psi_k\in L^1(X)$,
$\Vert \psi_k\Vert_{L^1} \ne 0$. Because     $$\int_X
\tilde{\varphi}(.,a)(\xi(\Delta_{G_1}/\Delta_G)^{1/2}(\beta(.,
a))R(a)da\vert_{{\partial O}\cup \infty} = 0,$$ we have
$$\begin{array}{rl} \sum_{k=1}^{N_n}\tilde{\chi}_k(.)\int_X \psi_k(a)
(\xi(\Delta_{G_1}/\Delta_G)^{1/2}(\beta(.,a))R(a) da\Vert_{{\partial
O}\cup \infty}\cr = \Vert \int_X\tilde{\varphi}_n(.,
a)(\xi(\Delta_{G_1}/\Delta_G)^{1/2}(\beta(.,a))R(a)da - \cr -
\int_X\tilde{\varphi}(.,a)(\xi(\Delta_{G_1}/\Delta_G)^{1/2}(\beta(.,
a))R(a)da \Vert_{{\partial O}\cup \infty}\cr \leq \int_X \int_{G_1} \vert
\varphi_n(b,a) - \varphi(b,a)\vert dbda.\end{array}$$
\end{pf}

Following the previous Lemma, $T(\varphi_n)$ can be approximated by a
sequence of operators
$A_{nm}$, $$A_{nm}f(x) = \theta_m(x) \int_X \tilde{\varphi}(x,
a)(\xi(\Delta_{G_1}/\Delta_G)^{1/2})(\beta(x,a))f(xa) da,$$ such that
$$\Vert T(\varphi_n) - A_{nm}\Vert = $$ $$\CD
 \Vert  (1-\theta_m)(.)\sum_{k=1}^{N_m}\tilde{\chi}(.)\int_X\psi_k(a)
(\xi(\Delta_{G_1}/\Delta_G)^{1/2})(\beta(.,a))R(a)   da\Vert
   @>>{n  \to  \infty \atop m \to \infty}> 0.\endCD$$

We have hence, $$\CD \Vert T(\varphi) - A_{nm}\Vert \leq
\Vert T(\varphi) - T(\varphi_n)\Vert + \Vert T(\varphi_n) - A_{nm}\Vert
@>>{n \to \infty \atop m \to\infty}> 0.\endCD$$ Following the previous
lemma, all the operators  $A_{nm}$ are compact operators. Thus the
operator $T(\varphi)$ is a compact operator and the proof of the theorem
is achieved.

\section{Application to Lie Group Representations}

We apply now the compactness criteria to description of compact type
ideals, associated with the representations, obtained from the
multidimensional orbit method.

\subsection{The case of solvable Lie groups}
For type I connected and simply connected solvable Lie groups, all the
irreducible unitary representations can be obtained from the orbit
method, see {\sc L. Auslander - B. Kostant}\cite{auslanderkostant}.

\begin{thm} Every irreducible unitary representation of a type I
connected and simply connected solvable or nilpotent Lie group can be
obtained as a representation, induced from an irreducible unitary
representation of a closed normal subgroup.
And hence, the compactness criteria are applicable.
 \end{thm}
 \begin{pf}
Indeed, for type I connected and simply connected solvable or nilpotent
Lie groups, every irreducible unitary representation can be obtained
from the orbit method, i.e. can be expressed in the form $T =
T^{\Omega_F,\chi_F}$ with some $F$ in
${\frakt g}^*$: $$T = T^{\Omega_F,\chi_F} = \Ind(G;{\frakt
p},H,\rho,\sigma).$$ Following the analysis of the structure of induced
representations in the previous chapter,
$$T = \Ind(G,{\frakt p},H,\rho,\sigma) \simeq \Ind_A^G\xi,$$ where
$\xi = \xi_1 \otimes \xi_2$, $A = M.N$ is a semi-direct product, where
$A$ is a normal subgroup in $G$, $$\xi_1 =
T^{\Omega_f,\chi_f} = \Ind(N;{\frakt p}_1,\rho,\chi_f),$$
$$\begin{array}{rl} \xi_2
&= T^{\Omega_l,\chi_l} \cr
                  &=   \Ind(M;{\frakt  p}_2,\rho, \chi_l)\cr 
                  &=  (\xi_2)_+\circ \pi \cr
                  &\simeq \Ind(M_+;  ({\frakt p}_2)_+, \rho_+,
\chi_{l_+})\circ \pi\end{array}$$ are irreducible unitary representations
of the nilpotent Lie groups $N$, $M_+$. Thus,
$\xi = \xi_1 \otimes \xi_2$ is a CCR representation
of the closed normal subgroup $A = M.N$.
 \end{pf}

\subsection{Generic Representations of Reductive Lie Groups}

We recall now the detailed{\sc R.  L. Lipsman}analysis for irreducible
representations of semi-simple and reductive Lie groups.  We shall show
that all the generic irreducible presentations can be
obtained as representations, repeatedly induced from irreducible
representation of
some appropriate normal subgroups, and therefore the compactness
criteria can be repeated applicable on each induction steps, see
\cite{diep7}.

Let $(\Phi,\tau) \in {\cal B}$ be a Duflo data, i.e.  $\Phi\in {\frakt g}^*$
is admissible in the sense of Duflo, integral and well-polarizable
and $\tau$ an irreducible unitary representation of two-fold covering of
$G_\Phi$, such that its restriction to
$G_\Phi$ is a multiple of $\chi_\Phi$. We denote the corresponding
representation, obtained by using the orbit method, by $\pi(\Phi,\tau)$.
Let $N$ be the unipotent radical of $G$, ${\frakt n} = Lie N$ is its Lie
algebra $\theta = \Phi\vert_{\frakt
n}$ is the restriction of $\Phi$ onto ${\frakt n}$.  We denote
$G^1 = G_\theta.N$, ${\frakt g}^1 = Lie G^1$, and
$\Phi^1 = \Phi\vert_{{\frakt g}^1}$.  It is easy to see that   $\Phi^1 \in
{\cal A}{\cal P}(G^1)$, and there exists a representation      $\tau^1\in
{\cal X}^i_{G^1}(\Phi^1)$ uniquely defined by $\tau$ such that $\pi(\Phi,\tau) =
\Ind^G_{G^1}\pi_{G^1}(\Phi^1,\tau^1)$. 

Repeat the procedure with $(\Phi^1,  \tau^1) \in {\cal B}(G^1)$,
suppose that $N^1$ is the unipotent radical of    $G^1$,...,   we
have a finite sequence, say $r$ steps, 
$$G^r = (G^r)_{\theta^r}.N^r,\quad (\Phi^r,\tau^r) \in {\cal B}(G^r).$$
Following the rule of step induction, we have
$$\begin{array}{rl} (\Phi,\tau)   &=
\Ind^G_{G^1}\Ind^{G^1}_{G^2}...\Ind^{G^{r-1}}_{G^r   }\pi_{G^r}(\Phi^r,
\tau^r)\cr  
            &= \Ind^G_{G^r}\pi_{G^r}(\Phi^r,\tau^r).\end{array}$$
            Suppose that
$\gamma = \pi_{G^r}(\theta^r)$ is the representation of nilpotent Lie
group, defined by {\sc A.  Kirillov}. 
We choose the  reductive factor$S$ following the  Levi decomposition of
$G^r$, in $(G^r)_{\theta^r}$. We have a semi-direct product
 $G^r = S.N^r$.  Because the action of $S$  fix  $\theta^r$,
then $\gamma$ can be canonically extended into a projective
representation $\tilde{\sigma}$ of $S$ in the space of the
representation
$\gamma$.  It is well-known that there exists an irreducible unitary
representation (can be projective), uniquely defined by
$\omega$ of     $S$ such that
$$\pi_{G^r}(\Phi^r,\tau^r) = (\omega \times \tilde{\gamma}) \times \gamma
. $$ Passing to the two-fold-covering, we Cain always consider $\omega$
as a unitary representation of $S$.  If      $\zeta = \Phi\vert_{{\frakt s}=
Lie S}$ then there exists a $\nu\in {\cal X}_S(\zeta)$, uniquely defined
by $\tau$ such that $\omega = \pi_S(\xi,\nu).$

\begin{thm} Every generic irreducible unitary representation of a
reductive or semi-simple Lie group can be obtained as a representation,
repeatedly induced  from irreducible representations of closed normal
subgroups. And hence in each induction step, we can apply the
compactness criteria.
\end{thm}
\begin{pf}
From the construction, as analyzed above, it is easy to see that
$G^i$ are  closed normal subgroup in  $G^{i-1}$. 

{\sc R. L. Lipsman}\cite{lipsman2} has had seen that one can restrict to
consider only the square-integrable  $\pi_S(\xi,\nu)$
In this case  $\pi_{G^r}(\Phi^r, \tau^r)$ is  CCR, then following the
Gel'fand - Piateskij-Schapiro theorem on     CCR-property, on each
induction step, we can apply the compactness criteria.
\end{pf}

\section{Bibliographical Remarks}

The main idea of this chapter appeared in \cite{diep1} and \cite{diep2}. 
The author started from a concrete example of the affine transformation 
group and then generalized into the main situation of C*-algebras of type 
I, which admit some boundary properties, guaranting existence of a 
canonical composition series. 
For the group C*-algebras, the author funded out some geometric criteria. 
This results were published in \cite{diep7}.

\chapter{Invariant Index of Group C*-Algebras}
This chapter is devoted to studying harmonic $L^1$-analysis on Lie
groups.  It is more complicate than the problem of studying harmonic
$L^2$-analysis exposed in the previous section. We must
consider the so called the
Fourier - Gel'fand transformation for functions  
in place of the well-defined Fourier transformation.
We shall see that
there are some analytical difficulties which require introducing K-theory 
in order to find  topological invariants for function algebras.
Below, we introduce the problem, define the construction of indices of
group $C^*$-algebras and  
illustrate the construction with examples.

\section{The Structure of Group C*-Algebras }

Let $G$ be a locally compact group and $dg$ the right- ( left- ) invariant
Haar measure. In the previous section it was shown that the constructed
irreducible representations are enough to decompose the ( right ) regular
representation of $G$ in $L^2(G,dg)$ into a sum of the so called discrete
series part $\sum^{\oplus}$ and the continuous series part
$\int^{\oplus}$,
$L^2(G,dg) = \sum^{\oplus} \dots \oplus \int^{\oplus} \dots d\mu(.)$.
It is more complicate to study the class $L^1$ function algebra
$L^1(G,dg)$ with the well - known convolution product $\varphi *
\psi$, $$(\varphi * \psi)(x):= \int_G \varphi(y)\psi(y^{-1}x)dy $$
and with involution $\varphi \mapsto \varphi^*$, $$ \varphi^*(x):=
\overline{\varphi(x^{-1})} \quad.$$ Remember that with the $L^1$-norm
$\parallel.\parallel_{L^1}$, $$\parallel \varphi \parallel_{L^1}:=
\int_G |\varphi(x)|dx $$ $L^1(G,dg)$ is a Banach ( involutive )
algebra. But its norm $\parallel. \parallel_{L^1}$ is not regular ;
in general, $$\parallel a^**a\parallel_{L^1} \leq \parallel
a\parallel^2_{L^1} \quad.$$ We try to introduce therefore some
regular norm. Recall the well - known Fourier - Gel'fand
transformation $$\varphi \in L^1 \mapsto
\hat{\varphi};\hat{\varphi}(\pi) = \pi(\varphi):= \int_G
\pi(x)\varphi(x) dx, \forall \pi \in \hat{G}.$$

We define the new $C^*$-norm $\parallel.\parallel_{C^*(G)}$ by
$$\parallel \varphi \parallel_{C^*(G)}:= \sup_{\pi\in
\hat{G}}{\parallel \pi(\varphi)\parallel}. $$ This norm is {\it
regular } 
, i.e. $$\parallel \varphi^* * \varphi\parallel_{C^*(G)} =
\parallel \varphi\parallel^2_{C^*(G)}.$$ The {\it group C*-algebra } 
\index{C*-algebra!group -}
$C^*(G)$ is defined as the completion of $L^1(G,dg)$ with respect to
this regular norm $\parallel. \parallel_{C^*(G)}$.

If $G$ is commutative, it returns to the Fourier -Gel'fand transform
and $C^*(G) \cong {\Bbb C}(\hat{G})$. With the Plancher\`el theorem
is related the well - known {\it Pontrijagin duality } 
\index{duality!Pontrijagin -} . If $G$ is compact 
but possibly non-commutative, $$C^*(G) \cong
\prod_{i=1}^{\infty} \Mat_{n_i}({\Bbb C}).$$ In particular, $C^*({\Bbb S}^1)
 \cong {\Bbb C}_0({\Bbb Z})$, what is just the Fourier analysis
theory of class $L^1$ functions. The Plancher\`el type theorem on
this class is the so called {\it Tanaka - Krein duality } 
\index{duality!Tanaka-Krein -}

The situation is rather complicate in cases, where $G$ is non-compact
and noncommutative. Some spectacular results were obtained in 1962 by
J. M. G. Fell for $C^*(\SL_2({\Bbb C}))$ and some times later by others
for $C^*(\SL_2({\Bbb R}))$, $C^*(\widetilde{\SL_2({\Bbb R})})$,....
by the same method of analytic description of the ( noncommutative )
Fourier - Gel'fand transforms.

The question about the structure of the group C*-algebra $C^*(G)$ was
open for example for the simplest non-trivial solvable Lie group of
affine transformations of the real straight line ${\Bbb R}$,
$$\Aff {\Bbb R}:= \{ g = (a,b): {\Bbb R} \rightarrow {\Bbb
R}; x \mapsto ax+b, \forall x \in {\Bbb R} ; a \ne 0 ; a,b \in {\Bbb
R}\} $$

For $G = \Aff {\Bbb R}$ it was well - known a detailed description of
the dual object $\hat{G}$ already in 1945 by I. M. Gel'fand and M. A.
Naimark: $\hat{G}$ consists of one infinite dimensional
representation $T$, which is itself dense in the whole $\hat{G}$.
All the other representations are the one-dimensional representations
( more precisely, the characters ) $U^{\varepsilon}_{\lambda},
\varepsilon = 0,1, \lambda \in {\Bbb R}$,
$$U^{\varepsilon}_{\lambda}(a,b):=
(\sgn (a))^{\varepsilon}.|a|^{\sqrt{-1}\lambda}.$$ The irreducible
unitary representation $T$ can be realized in the space $L^2({\Bbb
R}^*,{dx \over |x|}) $, where ${\Bbb R}^*:= {\Bbb R} \setminus (0)$
, by the formula $$ (T(a,b)f)(x) = \exp{(\sqrt{-1}bx)}f(ax).$$ Now if
we consider the connected component of identity $$G_0:= (\Aff {\Bbb
R})_0:= \{ (a,b)\in \Aff {\Bbb R} ; a>0 \},$$ there are two infinite dimensional
representations, the closure of them contains one-parameter family
$U_{\lambda}, \lambda \in {\Bbb R}$ of one-dimensional
representations ( more precisely, the characters ),
$$U_{\lambda}(a,b) = a^{\sqrt{-1}\lambda}.$$ These representations
$T_{\pm}, U_{\lambda};\lambda\in{\Bbb R}$ exhaust the dual object
$\hat{G}_0$ and are in a one-to-one correspondence with the co-adjoint
orbits $${\frakt g}^*/G = \{ \Omega_{\pm}, b^* \in {\Bbb R}\}.$$

However, the structure of the group C*-algebra $C^*(\Aff {\Bbb R})$ and
$C^*((\Aff {\Bbb R})_0)$ was unknown until 1975, when the BDF-K-functor
${\cal E}
xt_*$ was used to calculate the 
topological invariant $Index C^*(G)$ for these groups. This topological
invariant $Index C^*(G)$ defines the topological equivalence class of
$C^*(G)$ up to some isomorphisms. So it was computed that $$Index C^*
(\Aff {\Bbb R}) = (1,1) \in {\cal E}
xt({\Bbb S}^1 \vee {\Bbb S}^1) \cong {\Bbb Z} \oplus {\Bbb Z},$$
$$Index C^*((\Aff {\Bbb R})_0) = (1,1) \in {\cal E}
xt({\Bbb C}({\Bbb S}^1),{\cal K \oplus K}
) \cong {\Bbb Z} \oplus {\Bbb Z} \quad.$$

\section{Construction of $Index C^*(G)$}

Let us denote by G a connected and simply connected Lie group,${\frakt g}
= \Lie  G$ its Lie algebra, ${\frakt g}^* = \Hom_{\Bbb R}({\frakt g},{\Bbb
R})$ the dual vector space, ${\cal O = O} (G)$ the space of all the
co-adjoint orbits of G in {\frakt g}*. This space is a disjoint union of
subspaces of co-adjoint orbits of fixed dimension, i.e.  $${\cal O} =
\amalg_{0 \leq 2n \leq \dim G}{\cal O}_{2n},$$ $${\cal O}_{2n}:= \{
\Omega \in {\cal O} ; \dim \Omega = 2n \}.$$ We define $$V_{2n}:=
\cup_{\dim \Omega = 2n}\Omega.$$ Then it is easy to see that $V_{2n}$ is
the set of points of a fixed rank of the Poisson structure bilinear
function $$\{X,Y\}(F) = \langle F,[X,Y]\rangle ,$$
 hence it is a foliation, at least for $V_{2n},2n=max$.

\begin{prop}
The foliation $V_{2n}$
can be obtained by the associated action of ${\Bbb R}^{2n}$ on $V_{2n}$ via 
2n times repeated action of ${\Bbb R}$.
\end{prop} 
\begin{pf}
Indeed, fixing any basis $X_1,X_2,\dots,X_{2n}$ of the tangent space ${\frakt g} /{\frakt g}_F$ of $\Omega$  at the point $F \in \Omega$, we can define an action ${\Bbb R}^{2n} \curvearrowright V_{2n}$ as
 $$({\Bbb R} \curvearrowright ({\Bbb R} \curvearrowright (\dots {\Bbb R} \curvearrowright V_{2n})))$$
by 
$$(t_1,t_2,\dots,t_{2n}) \longmapsto \exp(t_1X_1)\dots\exp(t_{2n}X_{2n}) F.$$
Thus we have the Hamiltonian vector fields 
$$\xi_k:= {d \over dt} |_{t=0}\exp(t_kX_k)F, k = 1,2,\dots,2n$$
and the linear span 
$$F_{2n} = \{\xi_1,\xi_2,\dots,\xi_{2n}\}$$ provides a tangent distribution.
The proposition is proved.
\end{pf}

\begin{thm}
$(V_{2n},F_{2n})$ is a measurable foliation in the sense of A. Connes.
\end{thm}
\begin{pf}
Let us denote by $f_k$ the generating function of the
Hamiltonian vector field $\xi_k$, i.e. $$df_k + \imath(\xi_k)\omega_F = 0
,$$ where $\omega_F$ is the symplectic structure on co-adjoint orbit
$\Omega_F$. It is well-known that in every symplectic manifold $f_k$ is
uniquely ( up to an additive constant ) defined by its Hamiltonian vector
field $\xi_k$. We have $$df_1 \wedge df_2 \wedge \dots \wedge df_{2n} =
(-1)^{2n}\imath(\xi_1)\omega_F\wedge \dots \wedge
\imath(\xi_{2n})\omega_F$$ $$= (\Pfaffian \enskip \omega_F)^2 \times Volume
\enskip element \ne 0$$ The theorem is proved.
\end{pf}

\begin{cor} 
{\it The Connes C*-algebra $C^*(V_{2n},F_{2n}), 0 \leq {2n} \leq
\dim G$ are well defined.}
\end{cor}

  Now we assume that the orbit method gives us a complete list of
irreducible representations of $G$, $$\pi_{\Omega_F,\sigma} =
\Ind(G,\Omega_F,\sigma,{\frakt p}), \sigma \in {\cal X}_G(F),$$ the
finite set of Duflo's data.

    Suppose that $${\cal O} = \cup_{i=1}^k{\cal O}_{2n_i}$$ is the
decomposition of the orbit space on a stratification of orbits of
dimensions $2n_i$, where $n_1 > n_2 \dots > n_k > 0$

We include $C^*(V_{2n_1},F_{2n_1})$ into $C^*(G)$. It is well known
that the Connes C*-algebra of foliation can be included in the algebra
of pseudo-differential operators of degree 0 as an ideal. This algebra
of pseudo-differential operators of degree 0 is included in C*(G).

We define $$J_1 = {\bigcap_{\Omega_F \in {\cal O}
(G) \setminus {\cal O}
_{2n_1}}} \Ker  \pi_{\Omega_F,\sigma},$$
and $$ {A_1} = {C^*(G)/J_1}.$$
Then $$C^*(G)/C^*(V_{2n_1},F_{2n_1}) \cong A_1$$
and we have 

$$\vbox{\halign{ #\quad &#\quad &#\quad &#\hfill\cr
$ 0 \rightarrow $ & $J_1 \rightarrow$ & $C^*(G) \rightarrow$ & $A_1 \rightarrow 0$ \cr
          \hfill & $\enskip\downarrow$ & $\quad \downarrow Id$ &
$\quad\downarrow$ \cr $0 \rightarrow$ & $C^*(V_{2n_1},F_{2n_1})
\rightarrow$ & $C^*(G) \rightarrow$ & $C^*(G)/C^*(V_{2n_1},F_{2n_1})
\rightarrow 0$ \cr}}$$

Hence $J_1 \simeq C^*(V_{2n_1},F_{2n_1})$ and we have 
$$O \rightarrow C^*(V_{2n_1},F_{2n_1}) \rightarrow C^*(G) \rightarrow A_1 \rightarrow 0.$$
Repeating the procedure in replacing $$C^*(G),C^*(V_{2n_1},F_{2n_1}),A_1,J_1 
\enskip by \enskip A_1,C^*(V_{2n_1},F_{2n_1}),A_2,J_2,$$we have 
$$0 \rightarrow C^*(V_{2n_2},F_{2n_2}) \rightarrow A_1 \rightarrow A_2 \rightarrow 0$$
etc....

So it is proved the following result.

\begin{thm}
The group C*-algebra C*(G)  can be included in a finite sequence of extensions
$$(\gamma_1): \qquad 0 \rightarrow C^*(V_{2n_1},F_{2n_1}) \rightarrow C^*(G) \rightarrow A_1 \rightarrow 0$$
$$(\gamma_2): \quad\qquad 0 \rightarrow C^*(V_{2n_2},F_{2n_2}) \rightarrow A_1 \rightarrow A_2 \rightarrow 0$$
$$\ldots\dots\dots \dots \dots$$
$$(\gamma_k):\qquad 0 \rightarrow C^*(V_{2n_k},F_{2n_k}) \rightarrow A_{k-1} \rightarrow A_k \rightarrow 0$$
where $\widehat{A_k} \simeq \Char(G)$
\end{thm}

\begin{cor} $Index C^*(G)$ is reduced to the system $Index
C^*(V_{2n_i},F_{2n_i}), i = 1,2,\dots, k$ by the invariants $$[\gamma_i]
\in KK(A_i,C^*(V_{2n_i},F_{2n_i})), i= 1,2,\dots,k.$$
\end{cor}

 {\it Ideally, all these invariants $[\gamma_i]$ could be computed
step-by-step from $[\gamma_k]$ to $[\gamma_1]$.}

\section{Reduction of the Indices}

Let us consider $C^*(V_{2n_i},F_{2n_i})$ for a fixed i. We introduce the following assumptions which were considered by Kasparov in nilpotent cases:

{\bf Assumption} $(A_1)$ \quad There exists $k \in {\Bbb Z}, 0 < k \leq
2n_i$ such that $$V_{gen}:= V_{2n_i} \setminus (\Lie  \Gamma)^\perp$$ has
its C*- algebra $$C^*(V_{gen},F|_{V_{gen}}) \cong C({\cal O}_{gen}^\sim)
\otimes {\cal K} (H),$$ where $$\Gamma:= {\Bbb R}^k \hookrightarrow
{\Bbb R}^{2n_i} \hookrightarrow G,$$ $$\Lie  \Gamma = {\Bbb R}^k
\hookrightarrow {\frakt g}/{\frakt g}_{F_i},(\Lie  \Gamma)^\perp \subset
{\frakt g}^* \cap V_{2n_i}.$$

\begin{exam}
 If $V_{gen}$ is a principal bundle, or the space ${\cal O}_{gen} =
V_{gen}/G$ is a Hausdorff space, then $C^*(V_{gen},F|_{V_{gen}}) \simeq
C({\cal O}_{gen}^\sim) \otimes {\cal K} (H)$
\end{exam}

 It is easy to see that if the condition $(A_1)$ holds
,$C^*(V_{2n_i},F_{2n_i})$ is an extension of $C^*(V_{2n_i} \setminus
V_{gen},F_{2n_i}|_.)$ by $C({\cal O}_{gen}^\sim) \otimes {\cal K} (H)$
,where ${\cal O}_{gen}^\sim = \{ \pi_{\Omega_F,\sigma}; \Omega_F \in {\cal
O}_{gen},\sigma \in {\cal X}_G(F)$. If $k = 2n_i, ({\Bbb
R}^{2n_i})^\perp = \{ O \}$, $V_{2n_i} = V_{gen}$, we have
$$C^*(V_{2n_i},F_{2n_i}) \simeq C({\cal O}_{2n_i}^\sim \otimes {\cal K}
(H).$$ If $k = k_1 < 2n_i$, then ${\Bbb R}^{2n_i-k_1}$ acts on $V_{2n_i}
\setminus V_{gen}$ and we suppose that a similar assumption $(A_2)$ holds

$(A_2)$ There exists $k_2,0 < k_2 \leq 2n_i-k_1$ such that $$(V_{2n_i}
\setminus V_{gen})_{gen}:= (V_{2n_i} \setminus V_{gen}) \setminus ({\Bbb
R}^{k_2})^\perp$$ has its C*-algebra $$C^*((V_{2n_i} \setminus
V_{gen})_{gen},F_{2n_i}|.) \simeq C(({\cal O}_{2n_i} \setminus {\cal O}
_{gen})_{gen})^\sim \otimes {\cal K} (H).$$ As above, if $k_2 = 2n_i -
k_1$, $C^*(V_{2n_i} \setminus V_{gen},F_{2n_i}|.) \simeq C(({\cal O}
_{2n_i} \setminus {\cal O}_{gen})_{gen}^\sim) \otimes {\cal K} (H)$. In
other case we repeat the procedure and go to assumption $(A_3)$, etc.... 

The procedure must be finished after a finite number of steps, say in
m-th step,$$C^*((\dots(V_{2n_i} \setminus V_{gen}) \setminus (V_{2n_i}
\setminus V_{gen})_{gen} \setminus \dots,F_{2n_i}|_.) \simeq C((\dots
({\cal O}_{2n_i} \setminus {\cal O}_{gen}) \setminus \dots )^{\sim})
\otimes {\cal K} (H).$$ Thus we have the following result.

\begin{thm}
If all the arising assumptions $(A_1),(A_2),\dots$ hold, the
C*-algebra $C^*(V_{2n_i},F_{2n_i})$can be included in a finite sequence of
extensions $$ 0 \rightarrow C({\cal O}_{gen}^\sim) \otimes {\cal K} (H)
\rightarrow C^*(V_{2n_i},F_{2n_i}) \rightarrow C^*(V_{2n_i} \setminus
V_{gen},F_{2n_i}|_.) \rightarrow 0$$ $$ 0 \rightarrow C(({\cal O} _{2n_i}
\setminus {\cal O} _{gen})^\sim) \otimes {\cal K} (H) \rightarrow
C^*(V_{2n_i} \setminus V_{gen},F_{2n_i}) \rightarrow C^*(\dots)
\rightarrow 0$$ $$\ldots \ldots \ldots$$ $$ 0 \rightarrow C((\dots({\cal
O} _{2n_i} \setminus {\cal O} _{gen}) \setminus ({\cal O} _{2n_i}
\setminus {\cal O} _{gen}))_{gen}\dots{^\sim}) \otimes {\cal K} (H)
\rightarrow C^*(\dots) \rightarrow C^*(\dots) \otimes {\cal K} (H)
\rightarrow 0$$
\end{thm}

\section{General Remarks on Computation of Indices}

We see that the general computation procedure of Index C*(G) is reduced to
the case of short exact sequences of type $$(\gamma): \qquad 0
\rightarrow C(Y) \otimes {\cal K} (H) \rightarrow {\cal E} \rightarrow
C(X) \otimes {\cal K} (H) \rightarrow 0$$ $$ [\gamma] = Index {\cal E} \in
KK(X,Y).$$ The group $KK_i(X,Y)$ can be mapped onto $$\oplus_{j \in {\Bbb
Z}/(2)}\Hom_{\Bbb Z}(K^j(X),K^{i+j}(Y))$$ with kernel $$\oplus_{j \in {\Bbb
Z}/(2)}\Ext _{\Bbb Z}^1(K^j(X),K^{i+j+1}(Y))$$ by the well known cap-product
, see [K2]. So $[\gamma] = (\delta_0,\delta_1)$ $$\delta_0 \in \Hom_{\Bbb
Z}(K^0(X),K^1(Y)) = \Ext _0(X) \wedge K^1(Y)$$ $$\delta_1 \in \Hom_{\Bbb
Z}(K^1(X),K^0(Y)) = \Ext _1(X) \wedge K_0(Y).$$ Suppose $e_1,e_2,\dots,e_n
\in \pi^1(X)$ to be generators and $\phi_1,\phi_2,\dots,\phi_n \in {\cal
E}$ \quad the corresponding Fredholm operators, $T_1,T_2,\dots,T_n$ the
Fredholm operators, representing the generators of $K^1(Y) =
Index[Y,\Fred]$. We have therefore $$[\delta_0] = \sum_jc_{ij}\enskip
Index T_j,where$$ $$\delta_0 = (c_{ij}) \in \Mat _{\rank  K_0(X) \times \rank
K^1(Y)}({\Bbb Z}).$$ In the same way $\delta_1$ can be computed.

\section{Bibliographical Remarks}
This general ideas are due to the author of this book. It was firstly 
appeared in \cite{diep29}; then published in \cite{diep11}.


\begin{bibliography}{lbl}

\end{bibliography}
\printindex

\thebibliography{ Bibliogrphy}
\bibitem[A]{atiyah}
{\sc  Atiyah, M. F.} {\it K-theory}, Benjamin, New York 1976.

\bibitem[AK]{auslanderkostant} {\sc Auslander, L. and Kostant, B.}, {\it
Polarizations and unitary representations of solvable Lie groups } ,
Invent. Math., {\bf 14}(1971), No 4, 255-354. 

\bibitem[AS1]{atiyahsinger1} {\sc Atiyah, M. F. and Singer}, I. M., {\it
The index of elliptic operators.I } , Ann. Math., {\bf 87}(1968), 484-530. 

\bibitem[AS2]{atiyahsinger2} {\sc Atiyah, M. F. and Singer, I. M.}, {\it
The index of elliptic operators.III } , Ann. Math., {\bf 87}(1968),
546-604. 

\bibitem[AS3]{atiyahsinger3} {\sc Atiyah, M. F. and Singer, I. M.}, {\it
The index of elliptic operators. IV } , Ann. Math., {\bf 93}(1971),
119-138. 

\bibitem[AS4]{atiyahsinger4}
{\sc Atiyah, M. F. and Singer, I. M.}, {\it The index of elliptic operators.V }
, Ann. Math., {\bf 93}(1971), 139-149.

\bibitem[Ba]{bak}
{\sc Bak, A.}, {\it K-theory of forms }
, Princeton Univ. Press, Princeton, 1981.

\bibitem[BC]{baumconnes} {\sc Baum, P. and Connes, A.}, {\it leafwise
homotopy equivalence and rational Pontrijagin classes } , Preprint IHES,
Sept. 1983. 

\bibitem[BCD]{bcd} {\sc Bernat, P., Conze.N., Duflo, M. et coll.}, {\it
Repr\'esentations de groupes de Lie r\'esolubles } , Paris, Dunot, 1972. 

\bibitem[Bla1]{blattner1}
{\sc Blattner, P.}, {\it On induced representations.I }
, Amer. J. Math., {\bf 83}(1961), No. 1, 79-98.

\bibitem[Bla2]{blattner2}
{\sc Blattner, P.}, {\it On induced representations.II }
, Amer. J. Math., {\bf 83}(1961), No. 3, 499-512.


\bibitem[BDF1]{bdf1}
{\sc Brown, L. G., Douglas, R. G. and Fillmore, P. A.}, {\it Extensions
of C*-algebras }
, Ann. Math., Ser. 2, {\bf 105}(1977), No. 2, 265-324.

\bibitem[BDF2]{bdf2}
{\sc Brown, L. G., Douglas, R. G. and Fillmore, P. A.}, {\it Extensions
of C*-algebras, operators with compact self-commutators and K-homology }
, Bull. A.M.S., {\bf 79}(1973), No. 5, 973-978.

\bibitem[BM]{boyermartin} 
{\sc Boyer, R.  and Martin, R.},  {\it The group C*-algebra of the De
Sitter group.}

\bibitem[Bo1]{bouaziz}
{\sc Bouaziz, A.}, {\it Sur les representations des groupes de Lie
r\'eductifs non connexes }, preprint 1983.

\bibitem[Bou1]{bourbaki1}
{\sc Bourbaki,N.}, {\it Groupes et alg\`ebres de Lie, ch. I-III}, Hermann,
Paris 1971,1972.

\bibitem[Bou2]{boubaki2}
{\sc Bourbaki, N.}, {\it Vari\'et\'es diff\'erentielles et analytiques }
, Hermann, Paris, 1971.

\bibitem[Br]{brown}
{\sc Brown, L. G.}, {\it Extension and the structure of C*-algebras }
, Inst. Naz. Alta Matematica, Symp. Math., {\bf 20}(1976), 539-566.


\bibitem[Bu]{busby}
{\sc Busby, R. C.}, {\it Double centralizers and extensions of C*-algebras }
, Trans. A.M.S., {\bf 132}(1968), 79-99.

\bibitem[C1]{connes1} 
{\sc Connes, A.}, {\it A survey of foliations and operator algebras, }
Proc. Symp. Pure Math., {\bf 38}(1982), 521--628.


\bibitem[C2]{connes2} 
{\sc Connes, A.}, {\it C*-alg\`ebres et g\'eometrie diff\'erentielle, }
C. R. Acad. Sci. Paris, S\'erie A, {\bf 280}(1980), 599--604.


\bibitem[C3]{connes3} {\sc Connes, A.}, {\it An analogue of the Thom
isomorphism for crossed products of C*-algebras by an action of ${\Bbb R}$
} , Adv. Math., 38(1980), Adv. in Math. {\bf 39}(1981), 31-55. 

\bibitem[C4]{connes4} {\sc Connes, A.}, {\it Sur la th\'eorie
noncommutative de l'integration } , in {\it Alg\`ebres d'operateurs } ,
Lecture Notes in Math., No. 725, pp. 19 -43, Springer - Verlag, Berlin -
Heidelberg - New York, 1979. 

\bibitem[CNS]{corvinneemansternberg} {\sc Corvin, L., Nee'man,Y., and
Sternberg, S.}, {\it Graded Lie algebras in mathematics and in physics
(Bose-Fermi symmetry) } , Rev. Mod. Phys., {\bf 47}(1975), 573-604. 

\bibitem[CS]{connesskandalis} {\sc Connes, A. and Skandalis, G.}, {\it The
longitudinal index theorem for foliations } , Preprint IHES 1982. 

\bibitem[Cu]{cuntz}
{\sc Cuntz, J.}, {\it A survey of some aspects of noncommutative geometry, }
Preprint, Math. Inst. Uni Heidelberg, {\bf No 35}
(1992).

\bibitem[CQ1]{cuntzquillen1}
{\sc Cuntz, J. and Quillen, D.},
{\it Algebra extensions and nonsingularity, }
Preprint, Math. Inst. Uni Heidelberg, {\bf No 34}(1992).

\bibitem[CQ2]{cuntzquillen2}
{\sc Cuntz, J. and Quillen, D.},
{\it Operators on noncommutative differential forms and cyclic homology, }
Preprint, Math. Inst. Uni Heidelberg, {\bf No 38}(1992).

\bibitem[CQ3]{cuntzquillen3}
{\sc Cuntz, J. and Quillen, D.},
Preprint in preparation, 1993.

\bibitem[De]{delaroche} 
{\sc Delaroche, C.}, {\it Extensions des C*-alg\`ebres, }
Bull. Soc. Math. France, {\bf 29}(1972). 
           
\bibitem[D1]{diep1} 
{\sc Do Ngoc Diep}, {\it The structure of the group C*-algebra
of the group of affine transformations of the straight line, }
Funkt. Anal. i Priloz., {\bf 9}
(1975), No 1, 63-- 64.

\bibitem[D2]{diep2} 
{\sc Do Ngoc Diep}, {\it The structure of C*-algebras of type I, }
Vestnik Moskov. Uni., {\bf 1978 },  No 2, 81 --87.

\bibitem[D3]{diep3} 
{\sc Do Ngoc Diep}, {\it Applications of the homological K-functor
${\cal E}xt$ to studying the structure of the C*-algebras of some 
solvable Lie groups, }
Ph. D. Thesis, Moskov. Uni. {\bf 1977.}

\bibitem[D4]{diep4} 
{\sc Do Ngoc Diep}, {\it Quelques aspects topologiques en analyse
harmonique, }
Preprint IHES/M/83/22; Acta Math. Vietnam., {\bf 8}
(1983), No 2, 35--131.

\bibitem[D5]{diep5} 
{\sc Do Ngoc Diep}, {\it Construction des repr\'esentations
unitaires par les K-orbites et quantification, }
C. R. Acad. Paris, S\'erie A, {\bf 291}
(1980), 295--298.

\bibitem[D6]{diep6} 
{\sc Do Ngoc Diep}, {\it Quantification des syst\`emes
hamiltoniens \`a l'action plate d'un groupe de Lie, }
C. R. Acad. Sci. Paris, S\'erie I, {\bf 295}
(1982), 345--348.

\bibitem[D7]{diep7} 
{\sc Do Ngoc Diep}, {\it Id\'eaux de type compact associ\'es
aux repr\'esentations irr\'eductibles induites par des rep\'esentations 
liminaires de sous-groupes invariants, }
C. R. Acad. Sci. Paris, S\'erie I, {\bf 294}
(1982), 189--192.

\bibitem[D8]{diep8} 
{\sc Do Ngoc Diep}, {\it Construction et reduction of the K-theory
 invariant Index C*(G) of group C*-algebras, }
 Sonderforschungbereich 343 ``Diskrete Strukturen in der Mathematik", Uni Bielefeld, {\bf 92-015}  (1992),
 I.1-I.10.

\bibitem[D9]{diep9}
{\sc Do Ngoc Diep}, {\it Discrete series for loop groups I, }
Sonderforschungbereich 343 ``Diskrete Strukturen in der Mathematik", Uni Bielefeld {\bf 92-015}
(1992),IV.1-IV.16.

\bibitem[D10]{diep10}
{\sc Do Ngoc Diep}, {\it Multidimensional quantization and Fourier integral
operators, }
Forschungsgruppe ``Nichtkommutative Geometrie und Topologie", Math,
Inst. Uni Heidelberg, {\bf No 52}
(1992).

\bibitem[D11]{diep11}
{\sc Do Ngoc Diep}, {\it A survey of noncommutative geometry methods for group
algebras}, J. of Lie Theory, {\bf 3}(1993), 149-176.


\bibitem[D12]{diep12}
{\sc Do Ngoc Diep}, {\it Multidimensional quantization.I The general construction }
, Acta Math. Vietnam., {\bf 5}(1980), No 2, 42-55.

\bibitem[D13]{diep13}
{\sc Do Ngoc Diep}, {\it Methods of algebraic topology in harmonic analysis }
, Announcement of research results  ( in Vietnamese ), Hanoi, 1980.

\bibitem[D14]{diep14}
{\sc Do Ngoc Diep}, {\it Functor of projective limit in Banach categories }
, J. Math. ( in Vietnamese ), {\bf 9}(1981) No. 1, 16-20.

\bibitem[D15]{diep15}
{\sc Do Ngoc Diep}, {\it Multidimensional quantization II. The covariant derivation }
, Acta Math. Vietnam., {\bf 7}(1982), No. 1, 87-93.

\bibitem[D16]{diep16}
{\sc Do Ngoc Diep}, {\it Quantification multidimensionnelle III. Applications : Sur les repr\'esentations irr\'eductibles de groupes de diff\'eomorphismes }
, Acta Math. Vietnam., {\bf 8}(1983), No. 1, 59-72.

\bibitem[D17]{diep17}
{\sc Do Ngoc Diep}, {\it Geometric quantization }
, J. Math. ( in Vietnamese ), {\bf 11}(1983), No 3, 1-4.

\bibitem[Do18]{diep18}
{\sc Do Ngoc Diep}, {\it C*-complexes de Fredholm I. }
, Preprint IHES/M/83/53 ; Acta Math. Vietnam., {\bf 9}(1984), No. 1, 121-130.

\bibitem[D19]{diep19}
{\sc Do ngoc Diep}, {\it C*-complexes de Fredholm II. }
, Preprint IHES/M/83/64 ; Acta Math. Vietnam., {\bf 9}(1984), No. 2, 193-199.

\bibitem[D20]{diep20}
{\sc Do Ngoc Diep}, {\it The Borel-Serre compactification for Langlands type discrete groups }
, Preprint Series, No. 28,Inst. Math. Hanoi, 1986.

\bibitem[D21]{diep21}
{\sc Do Ngoc Diep}, {\it Multidimensional quantization and the generic representations }
, Preprint Series, Inst. Math. Hanoi, No. 29, 1986.

\bibitem[D22]{diep22} {\sc Do Ngoc Diep}, {\it Multidimensional
quantization of the Hamiltonian systems with supersymmetry } , Preprint
series, Inst. Math. Hanoi, No. 31, 1986. 

\bibitem[D23]{diep23}
{\sc Do Ngoc Diep}, {\it Multidimensional quantization and
Fourier integral operators } , Preprint, Inst. Math. Hanoi, 1986.

\bibitem[D24]{diep24} {\sc Do Ngoc Diep}, {\it On the Langlands type
discrete groups I. The Borel-Serre compactification } , Acta Math.
Vietnam., {\bf 12}(1987), No. 1, 41-54. 

\bibitem[D25]{diep25} {\sc Do Ngoc Diep}, {\it Multidimensional
quantization IV. The generic representations } , Acta Math. Vietnam., {\bf
13}(1988), 67--72. 

\bibitem[D26]{diep26} {\sc Do Ngoc Diep}, {\it Multidimensional
quantization V. The mechanical systems with supersymmetry } , Acta Math.
Vietnam., {\bf 15}(1990), No 1, 11--40. 

\bibitem[Do27]{diep27} {\sc Do Ngoc Diep}, {\it On the Langlands type
discrete groups II. The theory of Eisenstein series } , Acta Math.
Vietnam., {\bf 16}(1991), no 1, 77-90. 

\bibitem[D28]{diep28} {\sc Do Ngoc Diep}, {\it On the Langlands type
discrete groups III. The continuous cohomology } , Sonderforschungbereich
343 ``Diskrete Strukturen in der Mathematik", {\bf 92-015}, Uni Bielefeld. 

\bibitem[D29]{diep29}
{\sc Do Ngoc Diep}, {\it Construction and reduction of
the K-theory invariant Index C*(G) of group C*-algebras }
, Preprint 1992.

\bibitem[D30]{diep30} {\sc Do Ngoc Diep}, {\it Noncommutative geometry
methods for group algebras}, {\bf D. Sci. Dissertation}, Institute of
Math., NCST of Vietnam, Hanoi, 1995, 147pp. 

\bibitem[D31]{diep31} {\sc Do Ngoc Diep}, {\it Multidimensional
quantization and degenerate principal series}, Vietnam J. of Math., {\bf
23}(1995), 127-132. 

\bibitem[D32]{diep32}
{\sc Do Ngoc Diep}, {\it Vanishing theorem for representations with
regular lowest weight of loop groups}, Forschungsgruppe ``
Nichtkommutative Geometrie und Topologie", Uni Heidelberg, {\bf
75}(1993), 1-21.

\bibitem[DT]{diepthu} {\sc Do Ngoc Diep and Nguyen Van Thu}, {\it Homotopy
invariance of entire current periodic cyclic homology}, Vietnam J. of Math.
(to appear). 

\bibitem[DVS]{diepvietson}
{\sc Do Ngoc Diep, Ho Huu Viet and Vuong Manh Son},
{\it Sur la structure des C*-alg\`ebres d'une classe de groupes de Lie }
Preprint Series, Inst. Math. Hanoi, No. 7(1981) ; Acta math. Vietnam.
{\bf 8}(1983), No. 2, 90-125.

\bibitem[Di]{dixmier} 
{\sc Dixmier, J.}, {\it Les C*-alg\`ebres et leurs repr\'esentations,  }Paris,
Gauthier-Verlag, 1969.

\bibitem[Do1]{dong1} {\sc Tran Dao Dong}, {\ On globalization over
$U(1)$-covering of Zuckermann $({\frakt g} , K)$-modules, } Tap Chi Toan
Hoc (J. of Math.), {\bf 19} (1991), No 1, 60-72. 

\bibitem[Do2]{dong2} {\sc Tran Dao Dong}, {\it A reduction of the
globalization and $U(1)$-covering, } Preprint ICTP, Trieste, Italy, 1993. 


\bibitem[Du1]{duflo1}
{\sc Duflo, M.}, {\it Th\'eorie de Mackey pour les groupes de Lie alg\'ebriques}
, Acta Math. {\bf 149}
(1982), 153--213.

\bibitem[Du2]{duflo2}
{\sc Duflo, M.}, {\it On the trace formula for almost algebraic real
Lie groups }
, Proc. Conf. "Lie group representations III", Lecture Notes in Math.
, No 1077(1984), 101-165.

\bibitem[Du3]{duflo3}
{\sc Duflo, M.}, {\it Sur les extensions des repr\'esentations
irr\'eductibles des groupes de Lie nilpotents }
, Ann. Scient. Ec. Norm. Sup., 5(1972), 71-120.

\bibitem[DV]{dongvui}
{\sc Tran Dao Dong and Tran Vui}, {\it On the procedure of multidimensional
quantization, }
Acta Math. Vietnam., {\bf 14}
(1989), 19--30.


\bibitem[Em]{Emch}
{\sc Emch, G. G.}, {\it Algebraic methods in statistical physics and quantum field theory }
, Wiley - Interscience, Inc., New York, London, Sydney, Toronto, 1972.

\bibitem[Ev]{evans}
{\sc Evans, B.}, {\it C*-bundles and compact transformation
groups, }
Memoirs A.M.S., {\bf No 269}
(1982).

\bibitem[Fa]{Fack}
{\sc Fack, T.}, {\it K-th\'eorie Bivariante de Kasparov }
, S\'eminaire N. Bourbaki, exp. 605, F\'evrier, 1983.

\bibitem[Fe1]{fell}
{\sc Fell, J. M. G.}, {\it The structure of algebras of operator
fields, }
Acta Math., {\bf 106}
(1961), 233--280.

\bibitem[Fe2]{Fell2}
{\sc Fell, J. M. G.}, {\it A new prove that nilpotent groups are CCR }
, Proc. A.M.S., {\bf 13}(1962), 93-99.

\bibitem[Fe3]{Fell3}
{\sc Fell, J. M. G.}, {\it Weak containment and induced representations of groups I}
, Canad. J. Math., {\bf 14}(1962), No. 2, 237-268 ; II, Trans. A.M.S.,
{\bf 110}(1964), No. 3, 424-447.

\bibitem[Ge]{gel'fand}
{\sc Gel'fand, I. M.}, {\it On the elliptic equations }
, Uspechi Math. Nauk ( in Russian ), {\bf 15(113)}(1960), No. 3,.

\bibitem[GN]{gelfandnaimark}
{\sc Gel'fand, I. M. and Naimark, M.A.,} {\it Unitary representations of the group of affine transformations of the straight line, }
Dokl. Akad. Nauk SSSR, {\bf 55}
(1947), No 7, 571--574.
 
\bibitem[GoR]{gootmanrosenberg}
{\sc Gootman, E. C. and Rosenberg, J.}, {\it The structure of crossed product C*-algebras : A proof of the generalized Effros-Hahn conjecture }
, Invent. Math., {\bf 52}(1979), 283-298.

\bibitem[Gr]{green} 
{\sc Green, P.}, {\it C*-algebras of transformation groups
with smooth orbit space, }
Pacific J. Math., {\bf 72}
(1977), No 1.

\bibitem[Gu]{gutt}
{\sc Gutt, S.}, {\it Deformation quantization, }
Lectures at the Workshop on Representation Theory of Lie Groups,
International Centre for Theoretical Physics, Trieste, Italy, March 15 - April 2, 1993.
 
\bibitem[HW]{hiltomwiley}
{\sc Hilton, P. and Wiley, S.}, {\it Homology theory }
, Cambridge, 1960.


\bibitem[JT]{jensenthomsen}
{\sc Jensen, K. and Thomsen, K.},
{\it Elements of KK-theory, }
Birkh\"ause, Boston-Basel-Berlin, 1991.

\bibitem[Ka]{kac}
{\sc Kac, V.}, {\it Lie super-algebras }
, Adv. Math. {\bf 26}(1977), No. 1, 8-96.
 
\bibitem[KaS]{kaminkerschoket}
{\sc Kaminker, J. and Schoket, C.}, {\it Steenrod homology and operator algebras }
, Bull. A.M.S., {\bf 81}(1975), No 2, 431-434.

\bibitem[Kar1]{karoubi}
{\sc Karoubi, M.}, {\it K-theory: An introduction, }
Grundlehren der Math. Wissenschaften, {\bf No 226 }
, Springer-Verlag, Berlin-Heidelberg-New York, 1978.
 
\bibitem[Kar2]{karoubi2}
{\sc Karoubi, M.}, {\it Connections, courbures et classes caract\'eristiques en K-th\'eorie alg\'ebrique }
, Canadian Math. Soc. Conf. Proc., {\bf Vol. 2}, part 1, 1982, 19-27.

\bibitem[Kas1]{kasparov1}
{\sc Kasparov, G. G.}, {\it Topological invariants of elliptic operators I. K-homologies }
, Izvestija AN SSSR, Ser, mat., {\bf 39}(1975), No. 4, 796-838.

\bibitem[Kas2]{kasparov2}
{\sc Kasparov, G. G.}, {\it K-functor in the theory of C*-algebra extensions }
, Funkt. Anal. Priloz., {\bf 13}(1979), No 4, 73-74.

\bibitem[Kas3]{kasparov3}
{\sc Kasparov, G. G.}, {\it The operator K- functor and extensions of C*-algebras, }
Math. USSR Izvestija, {\bf 16}
(1981), No 3, 513--572.

\bibitem[Kas4]{kasparov4}
{\sc Kasparov, G. G.}, {\it The operator K-theory and its applications, }
Itogi Nauki Tekh., Ser. Soverem. Probl. Math., {\bf 27}
(1985), 3--31.

\bibitem[Kas5]{kasparov5}
{\sc Kasparov, G. G.}, {\it K-theory, group C*-algebras and higher signatures }
( Konspectus ), Parts I, II, Preprints, Chernogolovka, 1981.

\bibitem[Ki]{kirillov} 
{\sc Kirillov, A. A.}, {\it Elements of the theory of representations, }
Springer-Verlag, Berlin-Heidelberg-New York, 1976. 

\bibitem[KM]{kraljevicmilicic}
{\sc Kraljevic H. and Milicic, D.}, {\it The C*-algebra of the universal
covering group of $SL(2,{\Bbb R})$, }
Glasnik Math., Ser. III, {\bf 1972}
, 35--48.

\bibitem[KN]{kobayashinomizu}
{\sc Kobayashi, S. and Nomizu, K.}, {\it Foundations of differential
geometry }
, New York - London, 1963.

\bibitem[KoNS]{korwinneemansternberg}
{\sc Korwin, L., Ne'eman, Y. and Sternberg, S.}, {\it Graded Lie
algebras in mathematics and in physics ( Bose-Fermi symmetry ) }
, Rev. Mod. Phys., {\bf 47}(1975), 573-604.

\bibitem[KoS]{koschmanschwarzbach}
{\sc Kosman-Schwarzbach, Y.}, {\it D\'eriv\'ee de Lie de morphismes de
fibr\'es }
, Publ. Univ. Paris VII, t.3, G\'eometrie diff\'erentielle.

\bibitem[Kos1]{kostant1}
{\sc Kostant, B.}, {\it Quantization and unitary representations, Part
I. Prequantization }
, Lecture Notes in Modern analysis and applications I, Lecture Notes in
Math., No. 170, pp 87-208, Springer-Verlag, Berlin-Heidelberg-New York, 1978.

\bibitem[Kos2]{kostant2}
{\sc Kostant, B.}, {\it Graded manifolds, graded Lie theory and representations }
, Lecture Notes in Math., No. 570, pp. 177-306, Springer -Verlag, Berlin-Heidelberg-New York.

\bibitem[Li1]{lipsman1}
{\sc Lipsman, L.}, {\it Group representations }
, Lecture Notes in Math., No. 338, Springer - Verlag, Berlin - Heidelberg - New York, 1974.

\bibitem[Li2]{lipsman2}
{\sc Lipsman, L.}, {\it Generic representations are induced from square-integrable representations }
, Trans. A.M.S., {\bf 285}(1984), No. 2, 845-854.

\bibitem[Mac]{mackey}
{\sc Mackey, M. G.}, {\it Infinite dimensional Lie group representations}, Coll. Lecture given
at Stillwater, Oklahoma, August 29 - September 01, 1961 (Sixty-six Summer 
Meeting of the AMS).

\bibitem[Man]{manin}
{\sc Manin, Yu. I.}, {\it New dimensions in geometry }
, Uspekhi Mat. Nauk, {\bf 39}(1984), No 6(240), 47-73.

\bibitem[Mi1]{milicic1}
{\sc Milicic, D.}, {\it Topological representations of the group C*-algebra of
$SL(2,{\Bbb R})$, }
Glasnik Math., Ser. III, {\bf 1971}
,231--246.

\bibitem[Mi2]{milicic2}
{\sc Milicic, D.}, {\it Algebraic ${\cal D}$-modules and representations of
semi-simple Lie groups, }
Lectures at the Workshop on Representation Theory of Lie Groups, International Centre for Theoretical Physics, Trieste, Italy, March 15 - April 2, 1993.

\bibitem[Mis1]{mischenko1}
{\sc Mischenko, A. S.}, {\it Homotopy invariants of non-simply connected manifolds I. Rational invariants }
, Izvestija AN SSSR ( in Russian ), {\bf 4}(1970), 506-519.

\bibitem[Mis2]{mischenko2}
{\sc Mischenko, A. S.}, {\it Infinite dimensional representations of discrete groups and higher signatures }
, Izvestija AN SSSR ( in Russian ), {\bf 8}(1974), 85-111.

\bibitem[MisF1]{Mischenkofomenko1}
{\sc Mischenko, A. S. and Fomenko, A. T.}, {\it Generalized Liouville method of integration of Hamiltonian systems }
, Funkt. Anal. Priloz., {\bf 12}(1978), No. 2, 46-56.

\bibitem[MisF2]{MischenkoFomenko2}
{\sc Mischenko, A. S. and Fomenko, A. T.}, {\it Index of elliptic operators over C*-algebras }
, Izvestija AN SSSR ( in Russian ), {\bf 15}(1980), 87-112.

\bibitem[MisS]{MischenkoSolovev}
{\sc Mischenko, A. S. and Solov'ev, I. P.}, {\it On infinite dimensional representations of fundamental groups and formulae of Hizerbruch type }
, Dokl. AN SSSR ( in Russian ), {\bf 18}(1977), No. 3, 767-771.

\bibitem[N1]{nghiem1}
{\sc Nghiem Xuan Hai}, {\it Alg\`ebres de Heisenberg et g\'eom\'etrie symplectique de alg\`ebres de Lie }
, Publ. Math. Orsay, No. 78-08, Univ. Paris XI, 1978.

\bibitem[N2]{nghiem2}
{\sc Nghiem Xuan Hai}, {\it Une variante de la conjecture de Gel'fand - Kirillov et la transformation de Fourier - Plancher\`el }
, C. R. Acad. Sci. Paris, 293(1981), S\'erie I, 381-384.


\bibitem[N3]{nghiem3}
{\sc Nghiem Xuan Hai}, {\it La tranformation de Fourier - Plancher\`el analytique des groupes de Lie I. Alg\`ebres de Weyl et op\'erateurs diff\'erentiels }
, Publ. Math. Orsay, No. 81T22, Univ. Paris XI, 1981.

\bibitem[N4]{nghiem4}
{\sc Nghiem Xuan Hai}, {\it La transformation de Fourier - Plancher\`el analytique des groupes de Lie II. Les groupes nilpotents }
, Publ. Math. Orsay, No. 81T23, Univ. Paris XI, 1982.

\bibitem[N5]{nghiem5}
{\sc Nghiem Xuan Hai}, {\it La transformation de Fourier - Plancher\`el analytique des groupes de Lie r\'esolubles }
, Pr\'epublication, No. 80T39, Univ. Paris XI, 1982.

\bibitem[Per]{perdrizet}
{\sc Perdrizet, F.}, {\it Topologie et trace sur les C*-alg\`ebres, }
Bull. Soc. Math. France, {\bf 99}
(1971), 193--239.

\bibitem[Puk]{pukanszky}
{\sc Pukanszky, L.}, {\it Unitary representations of solvable Lie groups }
, Ann. Sci. \'Ecole Norm. Sup., S\'erie A, {\bf 4}(1971), No 4, 457-608.

\bibitem[Pus]{puschnigg}
{\sc Puschnigg, M.}, {\it Asymptotic cyclic cohomology, }
Dissertation, 1993, Univ. Heidelberg.


\bibitem[RoR]{robinsonrawnsley}
{\sc Robinson, P. L. and Rawnsley, J. H.}, {\it The metaplectic representation, $Mp^c$ structures and geometric quantization }
, Mem. A. M. S., {\bf Vol. 81}, No 410, 1989.


\bibitem[Ros]{rosenberg}
{\sc Rosenberg, A.}, {\it The number of irreducible representations of simple rings with no minimal ideals }
, Amer. J. Math., {\bf 75}(1953), 523-530.

\bibitem[Ros1]{rosenberg1}
{\sc Rosenberg, J.}, {\it The C*-algebras of some real and p-adic solvable groups, }
Pacific J. Math., {\bf 65}
(1976), No 1, 175--192. 

\bibitem[Ros2]{rosenberg2}
{\sc Rosenberg, J.}, {\it K-theory of group C*-algebras, foliation C*-algebras
and crossed products, }
Contemp. Math., {\bf 1988}
, 251--301.

\bibitem[Ros3]{rosenberg3}
{\sc Rosenberg, J.}, {\it Homological invariants of extensions of C*-algebras }
, Proc. Symp. Pure Math., {\bf Vol. 38}, pp. 35-75, A.M.S., Providence, R.I..

\bibitem[Ros4]{rosenberg4}
{\sc Rosenberg, J.}, {\it Realization of square-integrable representations of unimodular Lie groups in $L^2$-cohomology spaces }
, Trans. A.M.S., {\bf 261}(1980), No. 1, 1-31.

\bibitem[Ros5]{rosenberg5}
{\sc Rosenberg, J.}, {\it C*-algebras, positive scalar curvature and the Novikov conjecture }
, Appendix to {\it M. Gromov  and H. B. Lawson, Positive scalar curvature and the Dirac operators on complete Riemannian manifolds }
,.

\bibitem[RosS]{rosenbergschochet}
{\sc Rosenberg, J.and Schochet, C.}, {\it The classification of extensions of
C*-algebras, }
Bull. A.M.S., New series, {\bf 4}
(1981), 105-110.

\bibitem[Sa]{saito}
{\sc Saito, M.}, {\it Sur certains groupes de Lie r\'esolubles,} Sci. Papers of
The College of General Education, Univer. of Tokyo, {\bf 7}(1957)1-11, 
157-168.

\bibitem[ScW]{schmidwolf}
{\sc Schmid, W. and Wolf, J. A.},{\it Geometric quantization and deri
derived functor modules for semi-simple Lie groups }
, MSRI Berkeley, California, J. Funct. Anal.,.


\bibitem[SoV]{sonviet}
{\sc Vuong Manh Son and Ho Huu Viet}, {\it Sur la structure des C*-alg\`ebres d'une
classe de groupes de Lie, }
J. Operator Theory, {\bf 11}
(1984),77--90.

\bibitem[SuW]{sulankewintgen}
{\sc Sulanke, M. and Wintgen, P.}, {\it Differentialgeometrie und Faserb\"undel, }
Berlin. Vel. Deutscher Verlag der Wissenschaften.


\bibitem[Tam]{tamura}
{\sc Tamura, I.}, {\it Topology of foliations}, Mir, Moscow
1979 (Russian translation).

\bibitem[Tay]{taylor}
{\sc Taylor, J. L.}, {\it Banach algebras and Topology},
in "Algebra in Analysis", J. William Ed., Acad. Press, New York, 1975.

\bibitem[Tor]{torpe}
{\sc Torpe, A. M.}, {\it K-theory for the leaf space of foliations by
Reeb components, }
J. Funct. Anal., {\bf 61}
(1985), 15--71. 

\bibitem[Vi1]{viet1}
{\sc Ho Huu Viet}, {\it Sur la structure des C*-alg\`ebres d'une classe
de groupes de Lie r\'esolubles de dimension 3, }
Acta Math. Vietnam., {\bf 11}
(1986), No 1, 86--91.

\bibitem[Vi2]{viet2}
{\sc Ho Huu Viet}, {\it Application of the sequence of repeated extensions
to studying the structure of type I C*-algebras, }
Preprint Series, Inst. Math. Hanoi, {\bf No 87.24}
(1987), (in Russian).

\bibitem[Vo]{voiculescu} 
{\sc Voiculescu, D.}, {\it Remarks on the singular extensions in C*-algebra
of the Heisenberg group, }
J. Operator Theory. 

\bibitem[Vu1]{vu} {\sc Le Anh Vu}, {\it The foliation formed by the
K-orbits of maximal dimension of the real diamond group } , Tap chi Toan
hoc ( J. Math. ), {\bf 15}(1987), No 3, 1-10. 


\bibitem[Vu2]{vu1}
{\sc Le Anh Vu}, {\it On the structure of the C*-algebra of foliation formed
by the K-orbits of maximal dimension of the real diamond group, }
J. Operator Theory, {\bf 24}
(1990), No 2, 227--238.

\bibitem[Vu3]{vu2}
{\sc Le Anh Vu}, {\it On the foliation formed by the generic K-orbits of the
MD4-groups, }
Acta Math. Vietnam., {\bf 15}
(1990), No 2, 35--55. 

\bibitem[Vui1]{vui1}
{\sc Tran Vui}, {\it Multidimensional quantization and $U(1)$-covering, }
Acta Math. Vietnam., {\bf 16}
(1991), 103--119.

\bibitem[Vui2]{vui2}
{\sc Tran Vui}, {\it A reduction of the multidimensional quantization and
$U(1)$-covering, }
Tap Chi Toan Hoc (J. Math.), {\bf 19}
(1992), 1--12.

\bibitem[Vui3]{vui3}
{\sc Tran Vui}, {\it Geometric construction of unitary representations of Lie
groups via multidimensional quantization, }
Preprint ICTP, Trieste, Italy, 1993.

\bibitem[Wan1]{wang1} {\sc Wang, X.}, {\it Les C*-alg\`ebres d'une classe
de groupes de Lie r\'esolubles, } C. R. Acad. Sci. Paris, S\'erie I, {\bf
306} (1988), 765--767. 

\bibitem[Wan2]{wang2}
{\sc Wang, X.}, {\it On the C*-algebras of a family of solvable Lie groups
and foliations, }
Ph. D. Thesis, Univ. Maryland,  {\bf 1985.}

\bibitem[Wass]{wasserman} 
{\sc Wasserman, A}, 
{\it Une d\'emonstration de
la conjecture de Connes-Kasparov pour les groupes de Lie lin\'eaires
convese r\'eductifs, } Preprint IHES, 1983. 

\bibitem[Win]{winkelkemper}
{\sc  Winkelkemper, H. E.}, {\it The graph of foliation}, Preprint

\bibitem[Wo]{wolf}
{\sc  Wolf, J. A.}, {\it Admissible representations and
the geometry of flag manifold, }
Lectures at the Workshop on Lie  Representation Theory of Lie Groups, 
International Centre for Theoretical Physics, Trieste, Italy, 
March 15 - April 2, 1993; 
to appear in Proceedings of the A.M.S. Summer Conference 
``Penrose Transform and Analytic Cohomology''. 



\endthebibliography
\end{document}